\documentclass{amsbook}
\usepackage{times}
\usepackage{amsthm, amssymb, amsmath, a4wide}
\usepackage{bbm}
\usepackage{bm} 
\usepackage[inline]{enumitem}
\usepackage{graphicx}
\usepackage{csquotes} 
\usepackage{mathrsfs} 
\usepackage{dsfont}
\usepackage{exercise}

\usepackage{subcaption} 

\usepackage{varioref}
\usepackage[bookmarks=false]{hyperref}
\usepackage{cleveref}

\usepackage{xcolor}
\hypersetup{
    colorlinks,
    linkcolor={red!50!black},
    citecolor={blue!50!black},
    urlcolor={blue!80!black}
}

\usepackage{tikz}
\usetikzlibrary{calc}
\usetikzlibrary{matrix}
\usetikzlibrary{cd} 
\usetikzlibrary{arrows}
\usepackage{pgf}
\usepackage{pgfplots}
\pgfplotsset{compat=1.14}
\usetikzlibrary{positioning}

\makeatletter

\newcommand{\Rmnum}[1]{\expandafter\@slowromancap\romannumeral #1@}
\makeatother

\DeclareMathOperator\aff{aff}

\DeclareMathOperator\cl{Cl}
\DeclareMathOperator\conv{conv}

\DeclareMathOperator\gr{gr}
 
\DeclareMathOperator\ord{ord}

\DeclareMathOperator\rank{Rank}
\DeclareMathOperator\sing{Sing} 

\DeclareMathOperator\supp{Supp}
\DeclareMathOperator\vol{Vol} 

\newcommand{\scrA}{\ensuremath{\mathcal{A}}}
\newcommand{\scrB}{\ensuremath{\mathcal{B}}}
\newcommand{\scrC}{\ensuremath{\mathcal{C}}}
\newcommand{\scrD}{\ensuremath{\mathcal{D}}}
\newcommand{\scrE}{\ensuremath{\mathcal{E}}}
\newcommand{\scrF}{\ensuremath{\mathcal{F}}}
\newcommand{\scrG}{\ensuremath{\mathcal{G}}}

\newcommand{\scrI}{\ensuremath{\mathcal{I}}}
\newcommand{\scrJ}{\ensuremath{\mathcal{J}}}
\newcommand{\scrK}{\ensuremath{\mathcal{K}}}
\newcommand{\scrL}{\ensuremath{\mathcal{L}}}
\newcommand{\scrM}{\ensuremath{\mathcal{M}}}
\newcommand{\scrN}{\ensuremath{\mathcal{N}}}

\newcommand{\scrP}{\ensuremath{\mathcal{P}}}
\newcommand{\scrQ}{\ensuremath{\mathcal{Q}}}
\newcommand{\scrR}{\ensuremath{\mathcal{R}}}
\newcommand{\scrS}{\ensuremath{\mathcal{S}}}
\newcommand{\scrT}{\ensuremath{\mathcal{T}}}
\newcommand{\scrU}{\ensuremath{\mathcal{U}}}
\newcommand{\scrV}{\ensuremath{\mathcal{V}}}
\newcommand{\scrW}{\ensuremath{\mathcal{W}}}

\newcommand{\scrY}{\ensuremath{\mathcal{Y}}}
\newcommand{\scrZ}{\ensuremath{\mathcal{Z}}}

\newcommand{\cc}{\ensuremath{\mathbb{C}}}

\newcommand{\kk}{\ensuremath{\mathbb{K}}}

\newcommand{\pp}{\ensuremath{\mathbb{P}}}
\newcommand{\qq}{\ensuremath{\mathbb{Q}}}
\newcommand{\rr}{\ensuremath{\mathbb{R}}}

\newcommand{\zz}{\ensuremath{\mathbb{Z}}}

\newcommand{\aaa}{\ensuremath{\mathfrak{a}}}
\newcommand{\bbb}{\ensuremath{\mathfrak{b}}}
\newcommand{\ccc}{\ensuremath{\mathfrak{c}}}

\newcommand{\mmm}{\ensuremath{\mathfrak{m}}}
\newcommand{\nnn}{\ensuremath{\mathfrak{n}}}

\newcommand{\jjj}{\ensuremath{\mathfrak{j}}}

\newcommand{\ppp}{\ensuremath{\mathfrak{p}}}
\newcommand{\qqq}{\ensuremath{\mathfrak{q}}}
\newcommand{\rrr}{\ensuremath{\mathfrak{r}}}

\newcommand{\sheaf}{\ensuremath{\mathcal{O}}}
\newcommand{\ntorus}{(\cc^*)^n}
\newcommand{\nktorus}{(\kk^*)^n}
\newcommand{\ktorus}{\kk^*}

\newcommand{\im}{\ensuremath{\Rightarrow}}

\newcommand{\dsum}{\ensuremath{\bigoplus}}
\newcommand{\into}{\ensuremath{\hookrightarrow}}

\newcommand{\D}[2]{\frac{\partial #2}{\partial #1}}

\newtheorem{thm}{Theorem}[chapter]
\newtheorem*{thm*}{Theorem}
\newtheorem{claim}[thm]{Claim}
\newtheorem*{claim*}{Claim}
\newtheorem{cor}[thm]{Corollary}
\newtheorem{lemma}[thm]{Lemma}
\newtheorem*{lemma*}{Lemma}
\newtheorem{problem}[thm]{Problem}
\newtheorem*{problem*}{Problem}
\newtheorem*{conjecture*}{Conjecture}
\newtheorem{prop}[thm]{Proposition}
\newtheorem*{prop*}{Proposition}

\newtheorem*{propinition*}{Proposition-Definition}

\newtheorem*{thminition*}{Theorem-Definition}

\newtheorem{prolemma}{Lemma}[thm]
\newtheorem{proclaim}[prolemma]{Claim}

\theoremstyle{definition}

\newtheorem*{constrinition*}{Construction-Definition}

\newtheorem*{convention*}{Convention}
\newtheorem{defn}[thm]{Definition}
\newtheorem*{defn*}{Definition}

\newtheorem*{definotation*}{Definition-Notation}
\newtheorem{example}[thm]{Example}
\newtheorem*{example*}{Example}

\newtheorem{exercise}{Exercise}[chapter]
\newtheorem{exercise*}[exercise]{Exercise$^*$}
\newtheorem*{fact*}{Fact}

\newtheorem*{facts*}{Facts}

\newtheorem*{bold-note*}{Note}
\newtheorem{bold-question}[thm]{Question}
\newtheorem*{bold-question*}{Question}
\newtheorem{rem}[thm]{Remark}

\newtheorem*{reminition*}{Remark-Definition}

\newtheorem{remexample*}{Remark-Example}

\newtheorem*{remtation*}{Remark-Notation}

\newtheorem*{remuestion*}{Remark-Question}

\newtheorem*{remvention*}{Remark-Convention}

\theoremstyle{remark}
\newtheorem*{rem*}{Remark}
\newtheorem*{note*}{Note}
\newtheorem*{notation*}{Notation}
\newtheorem*{question*}{Question}
\newtheorem{question}[thm]{Question}
\newtheorem*{questions*}{Questions}

\author{Pinaki Mondal}
\title{How many zeroes?\\ {\small Counting solutions of systems of polynomials via toric geometry at infinity\\ Draft III
}}

\DeclareRobustCommand{\SkipTocEntry}[5]{} 

\def\*#1{\bm{#1}}

\newcommand{\Ai}{\scrA^I}

\newcommand{\Aibprime}{\scrA^{I_{B'}}}
\newcommand{\Aj}{\scrA^J}
\DeclareMathOperator\area{Area}
\newcommand{\Atildei}{\scrA^{\tilde I}}

\newcommand{\badj}{\scrB}

\newcommand{\Bi}{\scrB^I}
\DeclareMathOperator\bl{Bl}
\newcommand{\blnukn}{\bl_\nu(\kk^n)}
\newcommand{\blnuprimekn}{\bl_{\nu'}(\kk^n)}
\newcommand{\blnuprimeknn}[1]{\bl_{\nu'}(\kk^{#1})}
\DeclareMathOperator\character{char}

\newcommand{\Ci}{\cc^I}

\DeclareMathOperator\cone{cone}

\newcommand{\corbit}[1]{V_{#1}}
\DeclareMathOperator\degsep{\deg_{sep}}
\DeclareMathOperator\deginsep{\deg_{insep}}

\newcommand{\emptyiA}{\scrE^I_\scrA}
\newcommand{\Enu}{E_\nu}
\newcommand{\Enuprime}{E_{\nu'}}
\DeclareMathOperator\fund{fund}
\newcommand{\ggamma}{\boldsymbol{\bar \gamma}}
\newcommand{\goodI}{\mscrG}
\DeclareMathOperator{\graph}{gr}
\DeclareMathOperator{\grassman}{Gr}
\newcommand{\hatVzeroprime}{\hat \scrV'_{\origin}}
\newcommand{\hot}{\text{h.o.t.}}

\newcommand{\id}{\mathds{1}}
\DeclareMathOperator\In{In}
\newcommand\Innub{\In_{\nu_B}}
\newcommand{\Inalphapsinu}{\In'_{\alpha_\nu, \psi_\nu}}

\newcommand{\Inbetapsinu}{\In'_{\beta_\nu, \psi_\nu}}

\newcommand{\ki}{k^I}
\newcommand{\Ki}{\kk^I}
\newcommand{\kii}[1]{k^{#1}}
\newcommand{\Kib}{\Kii{I_B}}
\newcommand{\Kii}[1]{\kk^{#1}}
\newcommand{\Kiprime}{\Kii{I'}}
\newcommand{\Kitilde}{\Kii{\tilde I}}
\renewcommand{\kk}{\mathbbm{k}}
\newcommand{\Kns}{\kk^n_\mscrS}

\newcommand{\kstari}{(k^*)^I}
\newcommand{\Kstari}{(\kk^*)^I}
\newcommand{\Kstarib}{\Kstarii{I_B}}
\newcommand{\Kstariprime}{\Kstarii{I'}}
\newcommand{\Kstaritilde}{\Kstarii{\tilde I}}
\newcommand{\kstarii}[1]{(k^*)^{#1}}
\newcommand{\Kstarii}[1]{(\kk^*)^{#1}}

\DeclareMathOperator\lcm{lcm}
\DeclareMathOperator\ld{ld}
\newcommand{\local}[2]{\sheaf_{#2,#1}}  
\newcommand{\hatlocal}[2]{\hat {\sheaf}_{#2,#1}}

\newcommand{\mat}[1]{[#1]}

\newcommand{\mscrD}{\mathscr{D}}
\newcommand{\mscrE}{\mathscr{E}}

\newcommand{\mscrG}{\mathscr{G}}
\newcommand{\mscrI}{\mathscr{I}}
\newcommand{\mscrIA}{\mathscr{I}_\scrA}

\newcommand{\mscrS}{\mathscr{S}}

\newcommand{\mscrT}{\mathscr{T}}

\newcommand{\mult}[2]{[#1, \ldots, #2]}

\newcommand{\multA}[1]{\multp{\scrA_1}{\scrA_n}{#1}}
\newcommand{\multAiso}[1]{\multpiso{\scrA_1}{\scrA_n}{#1}}
\newcommand{\multAprimeiso}[1]{\multpiso{\scrA'_1}{\scrA'_n}{#1}}

\newcommand{\multAkniso}{\multAiso{\kk^n}}

\newcommand{\multAntorusiso}{\multAiso{\nktorus}}
\newcommand{\multAzero}{\multA{\origin}}

\newcommand{\multB}[1]{\multp{\scrB_1}{\scrB_n}{#1}}
\newcommand{\multBiso}[1]{\multpiso{\scrB_1}{\scrB_n}{#1}}
\newcommand{\multBzero}{\multB{\origin}}

\newcommand{\multGammazero}{\multzero{\Gamma_1}{\Gamma_n}}
\newcommand{\multGammazerostar}{\multzero{\Gamma_1}{\Gamma_n}^*}

\newcommand{\multf}[1]{\multp{f_1}{f_n}{#1}}
\newcommand{\multfiso}[1]{\multpiso{f_1}{f_n}{#1}}

\newcommand{\multfkniso}{\multfiso{\kk^n}}

\newcommand{\multfntorus}{\multf{\nktorus}}
\newcommand{\multfntorusiso}{\multfiso{\nktorus}}
\newcommand{\multfnctorusiso}{\multfiso{(\cc^*)^n}}
\newcommand{\multfzero}{\multf{\origin}}

\newcommand{\multg}[1]{\multp{g_1}{g_n}{#1}}
\newcommand{\multgiso}[1]{\multpiso{g_1}{g_n}{#1}}

\newcommand{\multgkniso}{\multgiso{\kk^n}}

\newcommand{\multgntorusiso}{\multgiso{\nktorus}}
\newcommand{\multgnctorusiso}{\multgiso{(\cc^*)^n}}
\newcommand{\multgzero}{\multg{\origin}}

\newcommand{\multkniso}[2]{\multpiso{#1}{#2}{\kk^n}}
\newcommand{\multntorus}[2]{\multp{#1}{#2}{\nktorus}}
\newcommand{\multntorusiso}[2]{\multpiso{#1}{#2}{\nktorus}}
\newcommand{\multnntorus}[3]{\multp{#1}{#2}{\nktoruss{#3}}}

\newcommand{\multp}[3]{\mult{#1}{#2}_{#3}}
\newcommand{\multpiso}[3]{\mult{#1}{#2}_{#3}^{iso}}
\newcommand{\multpstar}[3]{\mult{#1}{#2}^*_{#3}}
\newcommand{\multpnodots}[2]{[#1]_{#2}}
\newcommand{\multpisonodots}[2]{[#1]^{iso}_{#2}}
\newcommand{\multpstarnodots}[2]{[#1]^*_{#2}}

\newcommand{\multpzeronodots}[1]{[#1]_{\origin}}

\newcommand{\multzero}[2]{\multp{#1}{#2}{\origin}}
\newcommand{\multzerostar}[2]{\multpstar{#1}{#2}{\origin}}

\newcommand{\multPiso}[1]{\multpiso{\scrP_1}{\scrP_n}{#1}}

\newcommand{\multPkniso}{\multPiso{\kk^n}}

\newcommand{\multPntorusiso}{\multPiso{\nktorus}}

\DeclareMathOperator\mv{MV}
\DeclareMathOperator\nd{ND}
\newcommand{\nctoruss}[1]{(\cc^*)^{#1}}
\newcommand{\nktoruss}[1]{(\kk^*)^{#1}}
\newcommand{\norm}[1]{|| #1 ||}
\DeclareMathOperator\np{NP}
\newcommand{\Nstrong}{\check \scrN}
\newcommand{\Onu}{\orbit{\nu}}
\newcommand{\Onuprime}{\orbit{\nu'}}
\newcommand{\orbit}[1]{O_{#1}}
\newcommand{\origin}{0}

\newcommand{\partialxone}[1]{\partial #1/\partial x_1}
\newcommand{\partialxtwo}[1]{\partial #1/\partial x_2}
\newcommand{\partialxi}[1]{\partial #1/\partial x_i}
\newcommand{\partialxifrac}[1]{\frac{\partial #1}{\partial x_i}}
\newcommand{\partialxii}[2]{\partial #1/\partial x_{#2}}
\newcommand{\partialxiifrac}[2]{\frac{\partial #1}{\partial x_{#2}}}
\newcommand{\partialxj}[1]{\partial #1/\partial x_j}
\newcommand{\partialxk}[1]{\partial #1/\partial x_k}
\newcommand{\partialxn}[1]{\partial #1/\partial x_n}

\DeclareMathOperator{\perm}{perm}

\newcommand{\qplus}{\qq_{> 0}}
\newcommand{\qpluss}[1]{\qplus^{#1}}
\newcommand{\qzero}{\qq_{\geq 0}}

\DeclareMathOperator\relint{relint}
\newcommand{\ri}{\rr^I}
\newcommand{\rii}[1]{\rr^{#1}}
\newcommand{\riprime}{\rr^{I'}}
\newcommand{\rj}{\rr^J}
\newcommand{\rnnuperp}{\rr^n_{\nu^\perp}}
\newcommand{\rnstar}{(\rr^n)^*}
\newcommand{\rnnstar}[1]{(\rr^{#1})^*}
\newcommand{\rplus}{\rr_{> 0}}
\newcommand{\rpluss}[1]{\rplus^{#1}}
\newcommand{\rzero}{\rr_{\geq  0}}
\newcommand{\rzeroo}[1]{\rzero^{#1}}

\newcommand{\sepclosure}[1]{\overline{#1}_{sep}}
\DeclareMathOperator{\sign}{sign}

\newcommand{\talphanu}{T_{\alpha_\nu}}
\newcommand{\talphanuu}[1]{T_{#1}}
\newcommand{\tangent}[2]{T_{#2}(#1)}
\newcommand{\tAone}{\mscrT_{\scrA,1}}
\newcommand{\tAAone}[1]{\mscrT_{#1,1}}
\newcommand{\tAprimeone}{\mscrT_{\scrA',1}}
\newcommand{\thp}{T^H(\scrP)}
\newcommand{\thprimep}{T^{H'}(\scrP)}

\newcommand{\tiA}{T^I_{\scrA}}
\newcommand{\tiAprime}{T^I_{\scrA'}}
\newcommand{\tiprimeA}{T^{I'}_\scrA}
\newcommand{\titildeA}{T^{\tilde I}_\scrA}

\newcommand{\tiiAA}[2]{T^{#1}_{#2}}
\newcommand{\tjA}{T^J_{\scrA}}
\newcommand{\touch}{\mscrT}

\newcommand{\touchUA}{\touch(U,\scrA)}
\newcommand{\touchUAone}{\touch_1(U,\scrA)}
\newcommand{\touchUAprimeone}{\touch_1(U,\scrA')}

\newcommand{\touchUUAAone}[2]{\touch_1(#1, #2)}
\DeclareMathOperator{\trd}{tr.deg.}

\newcommand{\ua}{U_\scrA}

\newcommand{\Vhatprime}{\hat \scrV'}
\newcommand{\Vhatprimeinfty}{\hat \scrV'_\infty}
\newcommand{\Vhatprimess}[1]{\hat \scrV'_{#1}}
\newcommand{\Vi}{\scrV^I}
\newcommand{\Vitilde}{\scrV^{\tilde I}}
\newcommand{\Vinfty}{\scrV_\infty}
\newcommand{\Viprimeinfty}{\scrV'^I_\infty}
\newcommand{\Viprimess}[1]{\scrV'^I_{#1}}
\newcommand{\Vprimeinfty}{\scrV'_\infty}
\newcommand{\Vs}{\scrV_\mscrS}

\newcommand{\Vprimes}{\scrV'_\mscrS}
\newcommand{\Vprimess}[1]{\scrV'_{#1}}
\newcommand{\vonep}{V_\scrP^1}

\newcommand{\Vzero}{\scrV_{\origin}}
\newcommand{\Vzeroprime}{\scrV'_{\origin}}

\newcommand{\volsub}[1]{V^{-}_{#1}}
\newcommand{\volsubk}{\volsub{k}}
\newcommand{\volsubn}{\volsub{n}}

\newcommand{\woutlog}{without loss of generality}
\newcommand{\Woutlog}{Without loss of generality}
\newcommand{\wpn}{\pp^n(\omega)}
\newcommand{\wpnprime}{\pp^n(\omega')}
\newcommand{\wwpn}[1]{\pp^n(#1)}
\newcommand{\wwpnn}[2]{\pp^{#2}(#1)}
\newcommand{\xa}{\xpp{\scrA}}
\newcommand{\xaa}[1]{\xpp{#1}}
\newcommand{\xb}{\xpp{\scrB}}
\newcommand{\xbprime}{X_{\scrB'}}

\newcommand{\xonep}{\xp^1}

\newcommand{\xzerob}{\xb^0}
\newcommand{\xzeroa}{\xa^0}
\newcommand{\xzeroaa}[1]{\xpp{#1}^0}
\newcommand{\xp}{X_\scrP}
\newcommand{\xpp}[1]{X_{#1}}
\newcommand{\xpinfty}{X_{\scrP, \infty}}
\newcommand{\xq}{X_\scrQ}
\newcommand{\xzerop}{\xp^0}
\newcommand{\xzeroq}{\xq^0}
\newcommand{\znnuperp}{\zz^n_{\nu^\perp}}
\newcommand{\znnuzero}{\zz^n_{\nu \geq 0}}
\newcommand{\znstar}{(\zz^n)^*}

\newcommand{\znzero}{\zzeroo{n}}

\newcommand{\znzeropreceqalpha}{[\znzero]_{\preceq \alpha}}

\newcommand{\zplus}{\zz_{> 0}}
\newcommand{\zpluss}[1]{\zplus^{#1}}
\newcommand{\zb}{Z_\scrB}
\newcommand{\ziA}{\scrT'^I_\scrA}

\newcommand{\zzero}{\zz_{\geq 0}}
\newcommand{\zzeroo}[1]{\zzero^{#1}}

\newlist{defnlist}{enumerate}{3}
\setlist[defnlist,1]{label=(\alph*)}
\setlist[defnlist,2]{label=(\arabic*), ref=(\alph{defnlisti}.\arabic*)}
\setlist[defnlist,3]{label=(\roman*), ref=(\alph{defnlisti}.\arabic{defnlistii}.\roman*)}

\newlist{prooflist}{enumerate}{3}
\setlist[prooflist,1]{label=(\roman*)}
\setlist[prooflist,2]{label=(\arabic*), ref=(\roman{prooflisti}.\arabic*)}
\setlist[prooflist,3]{label=(\alph*), ref=(\roman{prooflisti}.\arabic{prooflistii}.\alph*)}

\crefname{claim}{claim}{claims}
\crefname{cor}{corollary}{corollaries}
\crefname{defn}{definition}{definitions}
\crefname{exercise}{exercise}{exercises}
\crefname{fact}{fact}{facts}
\crefname{observation}{observation}{observations}
\crefname{problem}{problem}{problems}
\crefname{proclaim}{claim}{claims}
\crefname{prop}{proposition}{propositions}
\crefname{thm}{theorem}{theorems}

\makeatletter
\newtheorem*{rep@theorem}{\rep@title}
\newcommand{\newreptheorem}[2]{%
	\newenvironment{rep#1}[1]{%
		\def\rep@title{#2 \ref{##1}}%
	 	\begin{rep@theorem}
 	}%
 	{\end{rep@theorem}}
}
\makeatother

\newreptheorem{thm}{Theorem}
\newreptheorem{lemma}{Lemma}

\makeatletter
\renewcommand*{\p@section}{\thechapter.} 
\renewcommand*{\p@subsection}{\thechapter.} 
\renewcommand*{\p@subsubsection}{\thechapter.} 
\makeatother
\renewcommand{\thechapter}{\Roman{chapter}}

\makeindex
\begin{document}

\let\cleardoublepage\clearpage

\maketitle
\setcounter{tocdepth}{1}






\pagenumbering{gobble}
\addtocontents{toc}{\SkipTocEntry}
\chapter*{Preface}
In this book we describe an approach through {\em toric geometry} to the following problem: ``estimate the number (counted with appropriate multiplicity) of isolated solutions of $n$ polynomial equations in $n$ variables over an algebraically closed field $\kk$.'' The outcome of this approach is the number of solutions for ``generic'' systems in terms of their {\em Newton polytopes}, and an explicit characterization of what makes a system ``generic.'' The pioneering work in this field was done in the 1970s by Kushnirenko, Bernstein and Khovanskii, who completely solved the problem of counting solutions of generic systems on the ``torus'' $(\kk \setminus \{0\})^n$. In the context of our problem, however, the natural domain of solutions is not the torus, but the affine space $\kk^n$. There were a number of works on extending Bernstein's theorem to the case of affine space, and recently it has been completely resolved, the final steps having been carried out by the author. \\

The aim of this book is to present these results in a coherent way. We start from the beginning, namely Bernstein's beautiful theorem which expresses the number of solutions of generic systems on the torus in terms of the {\em mixed volume} of their Newton polytopes. We give complete proofs, over arbitrary algebraically closed fields, of Bernstein's theorem, its recent extension to the affine space, and some other related applications including generalizations of Kushnienko's results on {\em Milnor numbers} of hypersurface singularities which in 1970s served as a precursor to the development of toric geometry. Our proofs of all these results share several key ideas, and are accessible to someone equipped with the knowledge of basic algebraic geometry. 
This book can serve as a companion to introductory courses on algebraic geometry or toric varieties. While it does {\em not} provide a comprehensive introduction to algebraic geometry, it does develop the relevant parts of the subject from the beginning (modulo some explicitly stated basic results) with lots of examples and exercises, and can be used as a quick introduction to basic algebraic geometry. We hope the readers who take that undertaking will be rewarded by a deep understanding of the affine B\'ezout problem.

\subsection*{Acknowledgements}
It was Pierre Milman who wanted me to write a book; it would not have been possible without his constant encouragement and support - with mathematics, and all sorts of things beyond mathematics - throughout these years. Even though the scope of the final version is considerably limited compared to his vision, I offer it as a first step. The encouragement from Eriko Hironaka worked as a catalyst during a critical period when the project was stuck. I sincerely thank Jan Stevens who sent numerous corrections after reading one of the earlier drafts. Najma Ahmad, Kinjal Dasbiswas, Naren Hoovinakatte, and especially, Jonathan Korman read parts of earlier drafts and gave important suggestions. Thanks are also due to the referees and editors, especially Keith Taylor, whose suggestions significantly improved the exposition. Over the last few years the work on this book took a great portion of my time owed to my friends and family, especially my mother Purnima Mondal and brother Protim Mondol. The application of points at infinity to Chickens' Road Crossing problem is due to Shatabdi Sarker; its presentation given in this book is due to Tanzil Rashid.

\tableofcontents

\cleardoublepage
\pagenumbering{arabic}

\chapter{Introduction} \label{introduction}
\def\picfontsize{\small}
\newlength{\subwidth}
\def\colorP{blue}
\def\colorL{red}
\def\opazero{0.5}
\def\colorzero{green}
\def\colorone{black}
\def\colortwo{orange}
\def\colordot{red}
\def\nsamples{103}

\section{The problem and the results}
This book is about the problem of computing the number of solutions of systems of polynomials, or equivalently, the number of points of intersection of the sets of zeroes of polynomials. In this section we formulate the precise version of the problem we are going to study and give an informal description of the results. One natural observation that simplifies the problem is that {\em intersection multiplicity} should be taken into account, e.g.\ even though a tangent line intersects a parabola at only one point, it should be counted with multiplicity two (see \cref{fig:parabola-tangent}).

\begin{figure}[h]
\begin{center}
\begin{tikzpicture}[
    thick,
    dot/.style = {
      draw,
      fill,
      circle,
      inner sep = 0pt,
      minimum size = 3pt
    },
    xscale=1, yscale=0.5]
  \def\x{1.8}
  \draw [\colorP] (-\x,\x^2) parabola [bend pos=0.5] bend +(0,-\x^2) +(2*\x,0);
  \draw (-1.5,0) -- (2, 0);
   \foreach \x in {0.4, 1.4, 2.6} {
   		\draw (-1.5,\x*0.9/4.4) -- (2, \x); 
	}
	\draw (0,0) node[dot, label = {below:O}] {};
	\def\x{2}
	\def\y{3}
	\def\txtwdt{6cm}
	\node [below right, text width= \txtwdt,align=justify] at (\x,\y) {
 		\small
			As secants approach the tangent at $O$ more and more closely, both of the two points of intersection move arbitrarily close to $O$.
		};
\end{tikzpicture}
\caption{A tangent line intersects a parabola at a point with multiplicity two} \label{fig:parabola-tangent}
\end{center}
\end{figure}

The geometric intuition for intersection multiplicity is the ``principle of continuity,'' the principle that continuous perturbations of systems result in continuous changes of associated metrics or invariants\footnote{``Consider an arbitrary figure in general position \ldots\ Is it not obvious that if \ldots\ one begins to change the initial figure by insensible steps, or applies to some parts of the figure an arbitrary continuous motion, then is it not obvious that the properties and relations established for the initial system remain applicable to subsequent states of this system provided that one is mindful of particular changes, when, say, certain magnitudes vanish, change direction or sign, and so on---changes which one can always anticipate a priori on the basis of reliable rules.'' -- J.\ V.\ Poncelet, the foremost exponent of the principle of continuity, in the introduction of {\em Trait\'e des propri\'et\'es projectives des figures} (1822), as cited in \cite{rosenfeld}.}. Since the number of points of intersection is a {\em discrete} invariant of a system, it follows that it must not change under a continuous perturbation. However, over real numbers points of intersection may disappear upon an infinitesimal deformation (see \cref{fig:vanishing-intersection}). On the other hand, this problem disappears if one also counts ``imaginary'' solutions (this is why the intersection theory over complex numbers, or, more generally, an algebraically closed field, is easier than the intersection theory over real numbers). In this book we will consider polynomial systems defined over algebraically closed fields\footnote{\textellipsis\ which Poncelet probably would not have approved of, given his attitude towards consideration of complex solutions; see \cite[Section 4.2]{gray} for a most interesting account of this history.}.

\begin{figure}[h]
\begin{center}
\begin{tikzpicture}[
    thick,
    dot/.style = {
      draw,
      fill,
      circle,
      inner sep = 0pt,
      minimum size = 3pt
    },
    xscale=1, yscale=0.5]
  \coordinate (O) at (0,0);
  \def\x{1.8}
  \draw [\colorP] (-\x,\x^2) parabola [bend pos=0.5] bend +(0,-\x^2) +(2*\x,0);

  \def\x{1.8}
  \def\yroot{1.1}
  \draw(-\x,\yroot^2) -- (\x, \yroot^2);
  \draw(-\x,0) -- (\x, 0);
  \draw(-\x,-\yroot^2) -- (\x, -\yroot^2);
  \draw (0,0) node[dot, label = {above:{\small $(0,0)$}}] {} ;
	\def\x{1.8}
   \node [\colorP, right] at (-\x+0.1,\x^2-0.1) {
   		\small
   		$y = x^2$
   	};
	\node [right] at (\x,\yroot^2) {
 		\small
		$y = \epsilon$ intersects the parabola at two real points with multiplicity one.
		};
	\node [right] at (\x,0) {
 		\small
		$y = 0$ intersects the parabola at one real point with multiplicity two.
		};
	\node [right] at (\x,-\yroot^2) {
 		\small
		$y = -\epsilon$ intersects the parabola at two imaginary points with multiplicity one.
		};
\end{tikzpicture}
\caption{Disappearance of real points of intersection} \label{fig:vanishing-intersection}
\end{center}
\end{figure}


If there are infinitely many solutions of a system of polynomials, then the solution set has positive dimensional components, and assigning multiplicity to these components is trickier; we bypass this problem in this book and consider only the number of {\em isolated}\footnote{A point is {\em isolated} in a set $S$ if it is open in $S$.} solutions. This implies in particular we do not consider ``underdetermined systems,''\footnote{A system is {\em underdetermined} or {\em overdetermined} depending on whether the number of equations is smaller or greater than the number of variables. \label{uo-determined}} since an underdetermined system over an algebraically closed field can only have either positive dimensional or empty sets of solutions. We also ignore ``overdetermined systems''\textsuperscript{\ref{uo-determined}} because of the relative difficulty in assigning multiplicities. The final form of the subject of this book is thus the following:

\begin{problem} [Affine B\'ezout problem] \label{counting-problem}
Given $n$ polynomials in $n$ variables over an algebraically closed field $\kk$, give a sharp estimate of the number of its isolated solutions counted with appropriate multiplicity, and determine the conditions under which it is exact.
\end{problem}

For $n = 1$, the fundamental theorem of algebra gives a complete answer: a polynomial of degree $d$ has precisely $d$ zeroes counted with multiplicity. For $n \geq 2$, there is a problem: points of intersection may run off to infinity (see \cref{fig:infinite-intersection}).

\begin{center}
\begin{figure}[h]
\begin{tikzpicture}[
    thick,
    dot/.style = {
      draw,
      fill,
      circle,
      inner sep = 0pt,
      minimum size = 3pt
    },
    xscale =0.66,
    yscale=0.33]

  \def\x{2.7}
  \draw [\colorP] (-\x,\x^2) parabola [bend pos=0.5] bend +(0,-\x^2) +(2*\x,0);
   \node [\colorP, right] at (-\x+0.1,\x^2-0.1) {
   		\small
   		$y = x^2$
   	};
   	\def\ybottom{-0.9}
   	\draw (0, \ybottom) -- (0, \x^2);
    \def\ymiddle{0.5}
    \def\xright{2.7}
   \foreach \x in {1.2,1.75, 2.6} {
   		\draw ({\x*(\ybottom - \ymiddle)/(\x^2 - \ymiddle)}, \ybottom)
  		   		-- (\xright,{\xright*(\x^2 - \ymiddle)/\x + \ymiddle});
	}
	\draw (0,0) node[dot, label = {below right:O}] {};
	\def\x{2.7}
	\def\y{2.7*2.7}
	\def\txtwdt{6.3cm}
	\node [below right, text width= \txtwdt,align=justify] at (\x,\y) {
 		\small
			As secants approach the vertical line at $O$ more and more closely, one of the points of intersection approaches $O$ and the other goes to infinity.
		};
\end{tikzpicture}
\caption{A  vertical line intersects the parabola at one point with multiplicity one}
 \label{fig:infinite-intersection}
\end{figure}
\end{center}

Any reasonable approach to \cref{counting-problem} therefore must take into account ``intersections at infinity.'' A theorem named after E.\ B\'ezout (1730--1783) is the most basic result that does it satisfactorily.



\begin{thm} [B\'ezout's theorem, affine version] \label{bezout}
The number of isolated solutions in $\kk^n$ of $n$ polynomials in $n$ variables is at most the product of their degrees. Moreover, this bound is exact if and only if the only common solution of the \index{Leading form}{\em leading forms}\footnote{The {\em leading form} of a polynomial is the sum of its monomial terms with the highest degree; e.g.\ if $f = 2x^3 + 7x^2y - 9y^2 + 7xy -x + 1$, then its degree is $3$ and the leading form is $2x^3 + 7x^2y$.} of the polynomials is the origin.
\end{thm}

\begin{example} \label{bezout-counterexample}
Consider the system in \cref{fig:infinite-intersection} consisting of the parabola $y - x^2 = 0$ and a line $ax + by + c = 0$. The B\'ezout bound is $2 \times 1 = 2$, and the leading forms are $-x^2$ and $ax + by$. As long as $b \neq 0$, the only solution to $-x^2 = ax + by = 0$ is $(0,0)$, so that the bound is exact. However, if $b = 0$,  i.e.\ the line is vertical, then any point of the form $(0, k)$, $k \in \kk$, is a common solution of the leading forms. Consequently the B\'ezout bound overestimates the number of solutions in this case, as illustrated in \cref{fig:infinite-intersection}.
\end{example}

From the perspective of {\em projective geometry}, the B\'ezout bound is the number of intersections of polynomial hypersurfaces in the {\em projective space} $\pp^n$, which is a {\em compactification} of the affine space $\kk^n$ formed by adjoining a ``hyperplane at infinity.'' Therefore the B\'ezout bound is exact if and only if the hypersurfaces do not intersect at any point at infinity on $\pp^n$. 
However, as Gauss famously remarked,\footnote{Discussing his friend H.\ Schumacher's purported proof of the parallel postulate, Gauss wrote to him (as cited in \cite{waterhouse}), ``I protest first of all against the use of an infinite quantity as a completed one, which is never permissible in mathematics. The infinite is only a fa\c{c}on de parler, where one is really speaking of limits to which certain ratios come as close as one likes while others are allowed to grow without restriction.'' } infinity is the limit of some process, and curves which approach arbitrarily close to each other in one process may grow apart in another. A natural class of compactifications of $\kk^n$ containing the projective space is that of {\em weighted projective spaces}. Given an $n$-tuple $\omega =(\omega_1, \ldots, \omega_n)$ of positive integers, the corresponding {\em weighted rational curve} $C^\omega_a$ through a point $a = (a_1, \ldots, a_n) \in \kk^n$ is the curve parametrized by the map $t \mapsto (a_1t^{\omega_1}, \ldots, a_nt^{\omega_n})$. In the same way that in the projective space straight lines with different slopes are separated at infinity, in the weighted projective space $\pp^n(1, \omega)$ the curves $C^\omega_a$ corresponding to distinct $a$ are separated at infinity. See \cref{fig:at-infinity} for an example with $\omega = (1, 2)$, in which case $\{C^\omega_a\}_a$ is the family of parabolas $\{{a_1}^2y - {a_2}x^2 = 0\}$. The ``weight'' of a monomial $x_1^{\alpha_1}x_2^{\alpha_2} \cdots x_n^{\alpha_n}$ corresponding to $\omega$ is $\omega_1\alpha_1 + \cdots + \omega_n\alpha_n$. If $f$ is a polynomial, then the corresponding {\em weighted degree} $\omega(f)$ of $f$ is the maximum of the weights of all the monomials appearing in $f$.  The {\em leading weighted homogeneous form} of $f$ is the sum of all monomials (with respective coefficients) of $f$ with the highest weight. Computing intersection numbers on $\pp^n(1, \omega)$ leads to the ``weighted B\'ezout theorem,'' of which the original theorem of B\'ezout (\cref{bezout}) is a special case (corresponding to $\omega = (1, \ldots, 1)$).

\begin{center}
\begin{figure}[h]

\def\ymin{-1.5}
\def\ymax{3}
\def\ymaxx{4}
\def\paray{3.5}
\def\tminsq{-1.6}
\def\tmaxsq{1.6}
\def\tx{0}
\def\ty{-0.5}
\def\N{2}

\def\cmin{0.75}
\def\cmax{1.25}

\setlength{\subwidth}{0.3\textwidth}
\def\scalefactor{0.5}

\begin{subfigure}[b]{\subwidth}
\begin{tikzpicture}[scale=\scalefactor]

\def\plotdiff{0.5}
\def\adiff{0.5}
\def\tmaxtoadiff{3}
\def\mmin{0.6}
\def\mmax{1.2}
\pgfmathsetmacro\b{\ymax/2}
\pgfmathsetmacro{\a}{\tmaxsq+\plotdiff+\adiff}

\pgfmathsetmacro{\xmin}{\tminsq}
\pgfmathsetmacro{\xmax}{\a + \b/\mmin}

\begin{axis}[
xmin = \xmin, xmax=\xmax, ymin = \ymin, ymax= \ymaxx,
axis equal=true, axis equal image=true, hide axis
]
\foreach [evaluate= \s as \colorfraction using (\s+1)*100/(\N+1)] \s in {0,...,\N}{
	\pgfmathsetmacro{\c}{\cmin + \s*(\cmax-\cmin)/\N};
	\pgfmathsetmacro{\m}{\mmin + \s*(\mmax-\mmin)/\N};
	\pgfmathsetmacro{\tmin}{\a - \b/\m};
	\pgfmathsetmacro{\tmax}{\a + \b/\m};
	\edef\temp{
		\noexpand
		\addplot[\colorP!\colorfraction, thick, domain=\tminsq:\tmaxsq, samples=\nsamples] ({x} ,{\c*x^2});
		\noexpand
		\addplot[\colorL!\colorfraction, thick, domain=\tmin:\tmax,samples=2](x,{\m*(x-\a) + \b});
		}
	\temp
}
\draw (axis cs:0,\paray) node [\colorP] {parabolas};
\draw (axis cs:4.2,\paray) node [\colorL] {lines};
\end{axis}
\end{tikzpicture}
\caption{$\kk^2$}
\label{fig:families-k^2}
\end{subfigure}
\begin{subfigure}[b]{\subwidth}
\begin{tikzpicture}[scale=\scalefactor]

\def\tminsq{-1.7}
\def\tmaxsq{1.7}

\def\plotdiff{1}
\def\adiff{2}
\def\tmaxtoadiff{3}
\def\mmin{0.75}
\def\mmax{1.25}
\pgfmathsetmacro\b{\ymax/2}
\pgfmathsetmacro{\a}{\tmaxsq+\plotdiff+\adiff}

\pgfmathsetmacro{\xmin}{\tminsq}
\pgfmathsetmacro{\xmax}{\a + \b/\mmin}

\def\initialshift{3}

\begin{axis}[
xmin = \xmin, xmax = \xmax, ymin = \ymin, ymax = \ymaxx,
axis equal image = true,
axis y line = none,
axis x line = middle,
axis line style = {-, \colorone},
xlabel = {line at infinity},
xlabel style = {at={(axis cs:\xmax,0)}, \colorone, anchor = north},
ticks = none
]

\pgfmathsetmacro{\tmin}{0};
\pgfmathsetmacro{\tmax}{\ymax};
\pgfmathsetmacro\b{0.07}
\def\mmin{0.0185}
\def\mmax{0.02}
\pgfmathsetmacro\du{1/\mmax - \plotdiff}

\foreach [evaluate= \s as \colorfraction using (\s+1)*100/(\N+1)] \s in {0,...,\N}{
	\pgfmathsetmacro{\c}{\cmin + \s*(\cmax-\cmin)/\N};
	\pgfmathsetmacro{\m}{\mmin + \s*(\mmax-\mmin)/\N};
	\edef\temp{
		\noexpand
		\addplot[\colorP!\colorfraction, thick, domain=\tminsq:\tmaxsq, samples=\nsamples] ({x} ,{\c*x^2});
		\noexpand
		\addplot[\colorL!\colorfraction, thick, domain=\tmin:\tmax,samples=2]({((\m*\a - \b)*x+1)/\m - \du},x);
		}
	\temp
}
\draw (axis cs:0,\paray) node [\colorP] {parabolas};
\draw (axis cs:5.4,\paray) node[\colorL]  {lines};
\end{axis}

\end{tikzpicture}
\caption{$\pp^2$}
\label{fig:families-P^2}
\end{subfigure}
\begin{subfigure}[b]{\subwidth}
\begin{tikzpicture}[scale=\scalefactor]

\def\plotdiff{1}
\def\adiff{2}
\def\tmaxtoadiff{3}
\def\mmin{0.75}
\def\mmax{1.25}
\pgfmathsetmacro\b{8}
\pgfmathsetmacro{\a}{6}

\pgfmathsetmacro{\xmin}{-3}
\pgfmathsetmacro{\xmax}{3}

\begin{axis}[
xmin = \xmin, xmax = \xmax, ymin = \ymin, ymax = \ymaxx,
axis equal image = true,
axis y line = none,
axis x line = middle,
axis line style = {-, \colortwo},
xlabel = {line at infinity},
xlabel style = {at={(axis cs:\xmax,0)}, anchor = north,\colortwo},
ticks = none
]

\pgfmathsetmacro\cmin{\xmin + 0.5}
\pgfmathsetmacro\cmax{\xmin+ 2}
\pgfmathsetmacro\du{0}
\foreach [evaluate= \s as \colorfraction using (\s+1)*100/(\N+1)] \s in {0,...,\N}{
	\pgfmathsetmacro{\c}{\cmin + \s*(\cmax-\cmin)/\N};
	\pgfmathsetmacro{\m}{\mmin + \s*(\mmax-\mmin)/\N};
	\edef\temp{
		\noexpand
		\addplot[\colorP!\colorfraction, thick, domain=\ymin:\ymax, samples=2] (\c, x);
		\noexpand
		\addplot[\colorL!\colorfraction, thick,domain=-2:2,samples=\nsamples]({\m*x + (\b - \m*\a)*x^2 + \du},x);
		}
	\temp
}
\draw (axis cs:-1.75,\paray) node [\colorP] {parabolas};
\draw (axis cs:2,1.9) node [\colorL] {lines};
\end{axis}

\end{tikzpicture}
\caption{$\pp^2(1,1,2)$}
\label{fig:families-WP^2}
\end{subfigure}

\caption{$\pp^2$ separates lines, but not parabolas, at infinity, whereas $\pp^2(1,1,2)$ separates parabolas, but not lines, at infinity}
\label{fig:at-infinity}
\end{figure}
\end{center}


\begin{thm}[Weighted B\'ezout theorem for positive weights]\label{wt-bezout}
Let $\omega$ be a weighted degree on the ring of polynomials with positive weights $\omega_i$ for $x_i$, $i = 1, \ldots, n$. Then the number of isolated solutions of polynomials $f_1, \ldots, f_n$ on $\kk^n$ is bounded above by $(\prod_j \omega(f_j))/(\prod_j \omega_j)$. This bound is exact if and only if the leading weighted homogeneous forms of $f_1, \ldots, f_n$ have no common solution other than the origin.
\end{thm}

\begin{example} \label{wt-example}
Let $\omega = (1,2)$, $f = y - x^2$ and $g = ax+ c$, $a \neq 0$. Then $\omega(f) = 2$, $\omega(g) = 1$, and the leading weighted homogeneous forms of $f$ and $g$ are respectively $y-x^2$ and $ax$ . The only solution to the leading weighted homogeneous forms of $f$ and $g$ with respect to $\omega$ is $(0,0)$, so \cref{wt-bezout} implies that the number of solutions of $f = g = 0$ is precisely the weighted B\'ezout bound $(\omega(f)\omega(g)/(\omega(x)\omega(y)) = (2\times 1)/(1 \times 2) = 1$, as we saw in \cref{fig:infinite-intersection}.
\end{example}

\begin{center}
\begin{figure}[h]
\def\xmin{-0.5}
\def\xmax{3.5}
\def\ymin{-.5}
\def\ymax{2.5}
\def\tx{0}
\def\ty{5}
\def\tw{5cm}
\setlength{\subwidth}{0.25\textwidth}
\def\scalefactor{0.6}

\tikzstyle{dot} = [\colordot, circle, minimum size=4pt, inner sep = 0pt, fill]

\begin{subfigure}[b]{\subwidth}
\begin{tikzpicture}[scale=\scalefactor]
\draw [gray,  line width=0pt] (\xmin, \ymin) grid (\xmax,\ymax);
\draw [<->] (0, \ymax) |- (\xmax, 0);

\draw[thick, fill=\colorzero, opacity=\opazero ] (0,0) -- (0,2) --  (2,1) -- (3,0) -- cycle;
\draw[thick, \colorone] (3,0) -- (2,1);
\draw[thick, \colortwo ] (0,2) -- (2,1);

\coordinate (oneperp) at (3.5,1.5);
\draw [thick, \colorone, ->] (2.5,0.5) -- (oneperp);
\node [anchor= south, \colorone] at (oneperp) {\picfontsize $\omega = (1,1)$};

\coordinate (twoperp) at (2,3.5);
\draw [thick, \colortwo, ->] (1,1.5) -- (twoperp);
\node [anchor= south, \colortwo] at (twoperp) {\picfontsize $\omega = (1,2)$};

\end{tikzpicture}
\caption{$\scrP$}\label{fig:P}
\end{subfigure}	
\begin{subfigure}[b]{\subwidth}
\begin{tikzpicture}[scale=\scalefactor]

\def\mmin{0.9}
\def\mmax{1.2}
\pgfmathsetmacro\b{8}
\pgfmathsetmacro{\a}{6}
\pgfmathsetmacro{\xmin}{-3}
\pgfmathsetmacro{\xmax}{2}
\def\ymin{-1.5}
\def\ymax{3}
\def\ymaxx{4}
\def\paray{3.5}
\def\N{2}

\begin{axis}[
xmin = \xmin, xmax = \xmax, ymin = \ymin, ymax = \ymaxx,
axis equal image = true,
axis y line = none,
axis x line = middle,
x axis line style = {-, \colortwo},
xlabel = {$c_{(1,2),\infty}$},
xlabel style = {at={(axis cs:\xmax,0)}, anchor = west,\colortwo},
ticks = none
]
\addplot[\colorone, thick, domain=\ymin:\ymax, samples=2] (0, x);
\draw (axis cs:0,\paray) node {$c_{(1,1),\infty}$};

\pgfmathsetmacro\cmin{\xmin + 0.5}
\pgfmathsetmacro\cmax{\xmin+ 2}
\pgfmathsetmacro\du{0}
\foreach [evaluate= \s as \colorfraction using (\s+1)*100/(\N+1)] \s in {0,...,\N}{
	\pgfmathsetmacro{\c}{\cmin + \s*(\cmax-\cmin)/\N};
	\pgfmathsetmacro{\m}{\mmin + \s*(\mmax-\mmin)/\N};
	\edef\temp{
		\noexpand
		\addplot[\colorP!\colorfraction, thick, domain=\ymin:\ymax, samples=2] (\c, x);
		\noexpand
		\addplot[\colorL!\colorfraction, thick,domain=-2:-0.1,samples=\nsamples]({(1- \m*x)/ ((\b - \m*\a)*x^2) + \du},x);
		\noexpand
		\addplot[\colorL!\colorfraction, thick,domain=0.1:2,samples=\nsamples]({(1- \m*x)/ ((\b - \m*\a)*x^2) + \du},x);
		}
	\temp
}
\draw (axis cs:-1.75,\paray) node [\colorP] {parabolas};
\draw (axis cs:1.2,0.9) node [\colorL] {lines};
\end{axis}

\end{tikzpicture}
\caption{A coordinate chart near infinity on $\xp$} \label{fig:at-xp-infinity}
\end{subfigure}	

\caption{Parabolas and lines near curves at infinity on $\xp$}  \label{fig:at-toric-infinity}

\end{figure}
\end{center}

The main class of compactifications considered in this book are {\em toric varieties} associated to {\em convex integral polytopes}\footnote{A {\em convex integral polytope} in $\rr^n$ is the convex hull of finitely many points in $\rr^n$ with integer coordinates.}. If $\scrP$ is an $n$ dimensional convex integral polytope in $\rr^n$, then the outer normal to each of its $(n-1)$-dimensional faces determines (up to a constant of proportionality) a weighted degree, and in the corresponding toric variety $\xp$, weighted rational curves corresponding to each of these weights are separated. See \cref{fig:at-toric-infinity} for an example of a toric variety in which {\em both} parabolas and lines are separated at infinity. It has two curves at infinity (with respect to $\kk^2$) corresponding to the two edges of $\scrP$ which are not along the axes; we denote these curves by $c_{\omega, \infty}$, where $\omega$ is the corresponding weight. Each $c_{\omega, \infty}$ separates the family of weighted rational curves corresponding to $\omega$. Computing intersection numbers of hypersurfaces on toric varieties yields a beautiful result of D.\ Bernstein, which we now describe.

\begin{figure}[h]
\begin{center}
\begin{tikzpicture}[
	dot/.style = {
      draw,
      fill,
      circle,
      inner sep = 0pt,
      minimum size = 3pt,
    }, scale=0.6
    ]
	\begin{scope}[shift={(0,0)}]
		\draw [gray,  line width=0pt] (-0.5,-0.5) grid (3.5,3.5);
		\draw [<->] (0,3.5) node (yaxis) [above] {$y$}
       	 |- (3.5,0) node (xaxis) [right] {$x$};
       	
       	\draw[ultra thick, fill=green, opacity=\opazero ] (0,0) -- (3,0) -- (2,1) -- (0,2) -- cycle;
       	\draw (0,0) node[dot, red] {};
        \draw (1,0) node[dot, red] {};
        \draw (3,0) node[dot, red] {};
        \draw (2,1) node[dot, red] {};
        \draw (0,2) node[dot, red] {};       	
        \draw (-0.5,-0.5) node [below right, text width= 4.5cm] {
			\small
			$\np(f)$
			\newline
		    $f = 1-x + 3x^3 +4x^2y - 7y^2 $
		};    	

	\end{scope}
	
	\begin{scope}[shift={(10,0)}]
		\draw [gray,  line width=0pt] (-0.5,-0.5) grid (3.5,3.5);
		\draw [<->] (0,3.5) node (yaxis) [above] {$y$}
     	 |- (3.5,0) node (xaxis) [right] {$x$};
       	
       	\draw[ultra thick, fill=green, opacity=\opazero] (0,0) -- (1,0) -- (2,1) -- (1,2) -- (0,1) -- cycle;
       	\draw (0,0) node[dot, red] {};
        \draw (1,0) node[dot, red] {};
        \draw (2,1) node[dot, red] {};
        \draw (1,2) node[dot, red] {};       	
        \draw (0,1) node[dot, red] {};
        \draw (1,1) node[dot, red] {};

        \draw (-0.5,-0.5) node [below right, text width= 4.5cm] {
			\small
			$\np(g)$ \newline
			$g =2 + x - y + xy + x^2y + xy^2$
		};    	
	\end{scope}
\end{tikzpicture}
\caption{Some Newton polytopes in dimension $2$}  \label{fig:np}
\end{center}
\end{figure}

The {\em Newton polytope} of a polynomial is the convex hull of all the exponents that appear in its expression, see \cref{fig:np}. V.\ I.\ Arnold noticed sometime in 1960s or 1970s that invariants of ``generic'' systems of polynomials tend not to depend on precise values of the coefficients of their monomials, but only on the combinatorial relations of the exponents of these monomials. The study of this phenomenon was a recurring topic at his seminars at Moscow University. While working on Arnold's question on determination of the {\em Milnor number}\footnote{The {\em Milnor number} is an invariant of a singularity, see \cref{milnor-defn}.} at the origin of a generic polynomial, A.\ Kushnirenko discovered that if all polynomials have the same Newton polytope, then for generic systems the number of isolated solutions which do not belong to any coordinate hyperplane has a strikingly simple expression: it is simply $n!$ times the volume of this polytope! D.\ Bernstein soon figured out how to remove the restriction on Newton polytopes (about 130 years before this F.\ Minding \cite{minding} discovered a special case of Bernstein's theorem in dimension two\footnote{A.\ Khovanskii gives a summary of Minding's approach in \cite[Section 27.3]{burago-zalgaller}; an English translation of \cite{minding} by D.\ Cox and J.\ M.\ Rojas appears in \cite{topics-modeling}.}).

\begin{thm} \label{bkk-thm-intro}
Let $N$ be the number (counted with appropriate multiplicities) of the isolated zeroes of polynomials $f_1, \ldots, f_n$ on $\nktorus := \kk^n \setminus \bigcup_i \{x_i = 0\}$.
\begin{enumerate}
\item {\em Kushnirenko \cite{kush-poly-milnor}:} If each $f_j$ has the same Newton polytope $\scrP$,  then $N \leq n! \vol(\scrP)$. If $\vol(\scrP)$ is nonzero, then the bound is exact if and only if the following condition holds:
\begin{align}
\parbox{.6\textwidth}{%
for each nontrivial weighted degree $\omega$, the corresponding leading forms of $f_1, \ldots, f_n$ do not have any common zero on $\nktorus$.}  \tag{$*$} \label{bkk-cond-intro}
\end{align}
\item {\em Bernstein \cite{bern}:} In general $N$ is bounded above by the {\em mixed volume}\footnote{The {\em mixed volume} is the canonical multilinear extension (as a functional on convex bodies) of the volume to $n$-tuples of convex bodies in $\rr^n$, see \cref{mixed-section} for a precise description.} of the Newton polytopes of $f_j$. If the mixed volume is nonzero, then the bound is exact if and only if \eqref{bkk-cond-intro} holds.
\end{enumerate}
\end{thm}


\begin{example}
If the Newton polytope of each polynomial contains the origin, then \cref{bkk-thm-intro} in fact gives an upper bound on the number of isolated solutions on $\kk^n$ and it is in general better than the bounds from \cref{bezout,wt-bezout}. For example, using the fact that mixed volume of two planar bodies $\scrP$ and $\scrQ$ is simply $\area(\scrP + \scrQ) - \area(\scrP) - \area(\scrQ)$ (\cref{mixed-example}), we see that Bernstein's bound for the number of solutions of $f = g = 0$ (where $f, g$ are as in \cref{fig:np}) is the area of the region shaded in blue in \cref{fig:np-sum}, which is equal to $8$. B\'ezout bound, on the other hand is $3 \times 3 = 9$; it is not hard to show that the $9$ is also the best possible weighed B\'ezout bound.
\end{example}

\begin{figure}[h]
\begin{center}
\begin{tikzpicture}[
	dot/.style = {
      draw,
      fill,
      circle,
      inner sep = 0pt,
      minimum size = 3pt
    }, scale=0.6
    ]
    \def\shiftone{4}
    \def\shifttwo{5}
    \draw (\shiftone,1.5) node {$+$};    	
    \draw (\shiftone+\shifttwo,1.5) node {$=$};    	
	\begin{scope}[shift={(0,0)}]
		\draw [gray,  line width=0pt] (-0.5,-0.5) grid (3.5,3.5);
       	
       	\draw[ultra thick, fill=green, opacity=\opazero] (0,0) -- (3,0) -- (2,1) -- (0,2) -- cycle;
		\draw (-0.5,-0.5) node [below right] {
     			\small
     			$\np(f)$
     		};    	
	\end{scope}
	
	\begin{scope}[shift={(\shifttwo,0)}]
		\draw [gray,  line width=0pt] (-0.5,-0.5) grid (3.5,3.5);
       	
       \draw[ultra thick, fill=green, opacity=\opazero] (0,0) -- (1,0) -- (2,1) -- (1,2) -- (0,1) -- cycle;
       \draw (-0.5,-0.5) node [below right] {
     			\small
     			$\np(g)$
     		};    	
	\end{scope}
	
	\begin{scope}[shift={(2*\shifttwo,0)}]
			\draw [gray,  line width=0pt] (-0.5,-0.5) grid (5.5,4.5);
	       	
	       	\draw [ultra thick] (0,0) -- (4,0) -- (5,1) -- (3,3) -- (1,4) -- (0,3) -- cycle;
	       	\draw[fill=blue, opacity=\opazero] (0,2) -- (0,3) -- (1,4) -- (3,3) -- (4,2) -- (3,1) -- (3,0) -- (2,1) -- cycle;
       		\draw[dashed, ultra thick, fill=green, opacity=\opazero ] (0,0) -- (3,0) -- (2,1) -- (0,2) -- cycle;
	       	\begin{scope}[shift={(3,0)}]
       			\draw[dashed, ultra thick, fill=green, opacity=\opazero ] (0,0) -- (1,0) -- (2,1) -- (1,2) -- (0,1) -- cycle;
	       	\end{scope}
	        \draw (-0.5,-0.5) node [below right] {
       			\small
       			$\np(f) + \np(g)$
       		};    	
		\end{scope}
		
\end{tikzpicture}
\caption{Minkowski sum of Newton polytopes of $f$ and $g$}  \label{fig:np-sum}
\end{center}
\end{figure}

The natural domain of solutions of systems of polynomials over a field $\kk$ is however not the torus $\nktorus$, but the affine space $\kk^n$. There are at least two different ways to extend Bernstein's formula to $\kk^n$. The approach motivated by the {\em polynomial homotopy} method for solving polynomial systems is as follows: given polynomials $f_1, \ldots, f_n$, one starts with a deformed system $f_1 = c_1, \ldots, f_n = c_n$ with nonzero $c_j$. For generic $f_1, \ldots, f_n$ all solutions of the deformed system are in fact on the torus, and their number is given by Bernstein's theorem. Then one counts how many of these solutions approach isolated solutions of $f_1, \ldots, f_n$ as each $c_j \to 0$. This approach is taken in \cite{khovanus,hurmfels-polyhedra,li-wang-bkk,rojas-wang,rojas-toric}. In particular, B.\ Huber and B.\ Sturmfels \cite{hurmfels-polyhedra} found the general formula through this approach; however they proved it in a special case, and only in characteristic zero. J.\ M.\ Rojas \cite{rojas-toric} observed that Huber and Sturmfels' formula works over all characteristics. The other approach is closer to Bernstein's original proof of his theorem: here one computes the number of ``branches'' of the curve defined by $f_2 = \cdots = f_n = 0$ and then the sum of the order of $f_1$ along these branches. General formulae through this approach were obtained by A.\ Khovanskii [unpublished]\footnote{Khovanskii described his result to the author at the {\em Askoldfest} in 2017. \label{khovanote}} and the author \cite{toricstein}. This formula requires knowing the {\em intersection multiplicity} at the origin of generic systems of polynomials. As an illustration we now state the weighted  B\'ezout formula for weighted degrees with possibly {\em negative} weights\footnote{For simplicity here we do not allow zero weights; see \cref{wt-bezout-21} for the statement without this restriction.}. Let $\omega$ be a weighted degree with nonzero weights $\omega_1, \ldots, \omega_n$. If $I_- := \{i: \omega_i < 0\}$, then the ``general weighted B\'ezout bound'' for the number of isolated zeroes of $f_1, \ldots, f_n$ is
\begin{align}
		\sum_{I \subseteq I_-} (-1)^{|I_-| - |I|}
		\cfrac{
			\prod_j \left(\max\{\omega(f_j),0\} + \sum_{i \in I_-}|\omega_i|\deg_{x_i}(f_j) \right)
		}{
			\prod_i |\omega_i|
		}
\end{align}
(\cref{wt-bezout-21}). Note that this reduces to the weighted B\'ezout bound from \cref{wt-bezout} in the case that each $\omega_i$ is positive, i.e.\ $I_- = \emptyset$. This bound is exact for generic $f_1, \ldots, f_n$, provided $\omega(f_j)$ is {\em nonnegative} for each $j$. In the general case, define 
\begin{align*}
\scrP_\omega(f)
	:= \{
		\alpha = (\alpha_1, \ldots, \alpha_n) \in \rr^n: 
			\alpha_i \geq 0\ \text{for each}\ i,\ 
			\langle \omega, \alpha \rangle \leq \omega(f),\
			\alpha_k \leq \deg_{x_k}(f)\ \text{for each}\ j \in I_-
	\}
\end{align*}
\begin{center}
\begin{figure}[h]

\def\shiftx{7.5}
\def\colorzero{green}
\def\colorone{red}
\def\colortwo{blue}
\def\colorthree{yellow}
\def\colorg{orange}
\def\opazero{0.5}
\def\viewx{120}
\def\viewy{15}
\def\titlex{2}
\def\titley{-1}
\def\scale{0.6}

\begin{subfigure}[b]{0.3\textwidth}
\begin{tikzpicture}[scale=\scale]
\pgfplotsset{every axis title/.append style={at={(0,-0.2)}}, view={\viewx}{\viewy}, axis lines=middle, enlargelimits={upper}}

\begin{axis}
\addplot3[ thick, draw, fill=\colortwo,opacity=\opazero] coordinates{(1,0,0) (3,2,0) (0,2,3) (0,0,1)};
\addplot3[ thick, draw, fill=\colorone,opacity=\opazero] coordinates{(0,2,0) (3,2,0) (0,2,3)};
\end{axis}
\end{tikzpicture}
\caption{$\omega = (1,-1,1)$, $\omega(f) = 1$}  \label{fig:pomega0-pnp}
\end{subfigure}\hspace{0.15\textwidth}
\begin{subfigure}[b]{0.3\textwidth}
\begin{tikzpicture}[scale=\scale]
\pgfplotsset{every axis title/.append style={at={(0,-0.2)}}, view={\viewx}{\viewy}, axis lines=middle, enlargelimits={upper}}

\begin{axis}
\addplot3[fill=white,opacity=0] coordinates{(1,0,0) (3,2,0) (0,2,3) (0,0,1)};
\addplot3[ thick, draw, fill=\colortwo,opacity=\opazero] coordinates{(0,1,0) (1,2,0) (3,2,2) (1,0,2) (0,0,1)};
\addplot3[ thick, draw, fill=\colorone,opacity=\opazero] coordinates{(3,2,2) (1,2,0) (0,2,0) (0,2,2)};
\addplot3[ thick, draw, fill=\colorg,opacity=\opazero] coordinates{(3,2,2) (1,0,2) (0,0,2) (0,2,2)};
%
\end{axis}
\end{tikzpicture}
\caption{$\omega = (1,-1,-1)$, $\omega(f) = -1$}  \label{fig:pomega0-pnn}
\end{subfigure}
\caption{$\scrP_\omega(f)$ for $f = x_2^2 + x_2x_3^2 + x_1x_2x_3$}  \label{fig:pomega0}
\end{figure}
\end{center}
(see \cref{fig:pomega0}). Given $I \subseteq \{1, \ldots, n\}$, let $\ri$ be the $|I|$-dimensional coordinate subspace of $\rr^n$ spanned by all $x_i$, $i \in I$. Then there is a collection $\touch$ of subsets of $\{1, \ldots, n\}$ such that for each $I \in \touch$, the number of distinct $j$ such that $\scrP_\omega(f_j)$ touches $\ri$ is precisely $|I|$, and the number of isolated zeroes of $f_1, \ldots, f_n$ is bounded by
\begin{align}
	\sum_{I \in \touch}
		\mv\left( \scrP_\omega(f_{j_1}) \cap \ri, \ldots, \scrP_\omega(f_{j_{|I|}}) \cap \ri \right)
		\times
		\multzero{\pi_{I'}(\scrP_\omega(f_{j'_1})}{\pi_{I'}(\scrP_\omega(f_{j'_{n-|I|}}))}
\label{general-wt-formula-00}
\end{align} 
(see \cref{wt-bezout-3} for the precise statement), where
\begin{itemize}
\item $j_1, \ldots, j_{|I|}$ (respectively, $j'_1, \ldots, j'_{n-|I|}$) is the collection of indices $j$ such that $\scrP_\omega(f_j)$ touches (respectively, does not touch) $\ri$;
\item $\mv(\cdot, \ldots, \cdot)$ is the mixed volume;
\item $I' := \{1, \ldots, n\} \setminus I$ is the complement of $I$, and $\pi_{I'}$ is the natural projection onto the coordinate subspace of $\rr^n$ spanned by all $x_{i'}$, $i' \in I'$, and
\item $\multzero{\cdot}{\cdot}$ is the intersection multiplicity at the origin of systems of generic polynomials with given Newton polytopes. 
\end{itemize}
The general formula for generic number of solutions on the affine space is no more difficult; it is of the same type as \eqref{general-wt-formula-00}, i.e.\ it is a sum of products of mixed volumes and generic intersection multiplicities at the origin (see \cref{extended-bkk-bound-0}). However, to use it one needs to compute the generic intersection multiplicity at the origin (i.e.\ the second factor in the summands of \eqref{general-wt-formula-00}). In the special case that each polynomial is ``convenient,''\footnote{A polynomial or power series is {\em convenient} if for each $j$, there is $m_j \geq 0$ such that the coefficient of $x_j^{m_j}$ is nonzero.} a formula for generic intersection multiplicity was given by L.\ Ajzenberg and A.\ Yuzhakov \cite{aizenberg-yuzhakov}; a Bernstein-Kushnirenko type ``non-degeneracy'' condition, i.e.\ the condition for the bound being exact, was also known for convenient systems (see e.g.\ \cite[Theorem 5]{esterov}). In the general case Rojas \cite{rojas-toric} gave a formula via Huber and Sturmfels' polynomial homotopy method. The non-degeneracy condition for the general case was established by the author in \cite{toricstein}.\\

As hinted above, the formula for the generic number of solutions on $\kk^n$ is straightforward once one has the formula for generic intersection multiplicity at the origin. Sufficient criteria under which the bound is exact can also be obtained easily by adapting the Bernstein-Kushnirenko non-degeneracy condition \eqref{bkk-cond-intro}; such criteria were given by several authors including Khovanskii \cite{khovanus}, Rojas \cite{rojas-toric}. Precise non-degeneracy conditions, i.e.\ which are {\em necessary and sufficient} for the bound to be exact on $\kk^n$, are however more subtle than \eqref{bkk-cond-intro}; consider e.g.\ the problem of characterizing non-degenerate systems on $\kk^3$ of the form
\begin{align*}
	f_1 & = a_1 + b_1x_1x_2 + c_1x_2x_3 + d_1x_3x_1       \\
	f_2 & = a_2 + b_2x_1x_2+ c_2x_2x_3  + d_2x_3x_1       \\
	f_3 & = x_3(a_3 + b_3x_1x_2 + c_3x_2x_3  + d_3x_3x_1)
\end{align*}
where $a_j, b_j, c_j, d_j \in \ktorus$ (this system is discussed in \cref{ex-aff-non-deg-1}). If all $a_j, b_j, c_j, d_j$ are generic, then it is straightforward to check directly that all common zeroes of $f_1, f_2, f_3$ on $\kk^3$ are isolated and they appear on $\nktoruss{3}$. Consequently, Bernstein's theorem implies that the number of solutions is the mixed volume of the Newton polytopes of the $f_j$, which equals $2$. Now if $a_1 = a_2$, $b_1 = b_2$, and the remaining coefficients are generic, then \eqref{bkk-cond-intro} continues to be true, so that Bernstein's theorem applies and number of solutions on $\nktoruss{3}$ is still $2$; in particular, the system continues to be non-degenerate on $\kk^3$. However, in this case the set of common zeroes of $f_1, f_2, f_3$ on $\kk^3$ also has a {\em positive} dimensional component, namely the curve $\{x_3 = a_1 + b_1x_1x_2 = 0\}$. This situation never arises in the case of Bernstein's theorem; indeed, existence of a positive dimensional component makes a system violate \eqref{bkk-cond-intro} and its straightforward adaptations. Unlike the Bernstein-Kushnirenko non-degeneracy criterion, the correct non-degeneracy criterion for $\kk^n$ needs to accommodate existence of positive dimensional components - it has to be able to differentiate between the cases when such a component leads to a loss of isolated solutions and when it does not; 
such a criterion was given by the author in \cite{toricstein}. \\

We mentioned above that the pioneering work of Kushnirenko on counting solutions of polynomial systems was motivated by his work on Milnor numbers of hypersurface singularities. In \cite{kush-poly-milnor} he gave a beautiful formula for a lower bound on the Milnor number, and showed that the bound is achieved by {\em Newton non-degenerate} singularities if either the characteristic is zero or if the polynomial is convenient. It was however clear from the beginning that Newton non-degeneracy is not necessary for the formula to hold, and it also does not imply ``finite determinacy.''\footnote{i.e.\ it does not ensure that the singularity at the origin is isolated.} C.\ T.\ C.\ Wall \cite{wall} introduced another notion of non-degeneracy which implies finite determinacy and which also guarantees that the Milnor number can be computed by Kushnirenko's formula.
S.\ Brzostowski and G.\ Oleksik \cite{brzostowski-oleksik} found the combinatorial condition which under Newton non-degeneracy is equivalent to finite determinacy. The Milnor number of a hypersurface at the origin is same as the intersection multiplicity at the origin of the partial derivatives of the defining polynomial (or power series). The non-degeneracy condition for intersection multiplicity therefore gives a natural starting point to study Milnor numbers. This condition generalizes both Newton non-degeneracy (for isolated singularities) and Wall's non-degeneracy condition; the author showed in \cite{toricstein} that in positive characteristic this condition is sufficient, and in zero characteristic it is both necessary and sufficient, for the Milnor number to be generic. \\

The purpose of this book is to give a unified exposition of the results described above. In addition to Bernstein's theorem (over arbitrary algebraically closed fields), classical results proved in this book include weighted homogeneous and multi-homogeneous versions of B\'ezout's theorem; complete proofs (or even, statements) of these results are otherwise hard to find. We followed Bernstein's original proof for establishing the non-degeneracy conditions of his theorem; in particular we present his simple and ingenious trick to construct a curve of solutions that runs off to infinity in the case that the non-degeneracy condition \eqref{bkk-cond-intro} is not satisfied\footnote{The bound from Bernstein's theorem and the sufficiency of \eqref{bkk-cond-intro} for the bound can be established without much difficulty (and in a very elegant way) using the general machinery of intersection theory (see e.g.\ \cite[Section 5.4]{fultoric}). However, we do not know of any proof of the {\em necessity} of \eqref{bkk-cond-intro} using this approach which does not involve an adaptation of Bernstein's trick; in all probability it would be much more difficult otherwise, since establishing positivity of excess intersections is in general a hard problem. Bernstein's trick is a nontrivial example of an elementary argument faring better than a formidable machinery.}. This book is the first part of a series of works on a constructive approach to compactifications of affine varieties started in the author's PhD thesis \cite{pisis}, for which the affine B\'ezout problem served as a motivation. Based on the results of this book, in the next part we give a solution to the general version of the affine B\'ezout problem, i.e.\ give a recipe to compute the precise number (counted with multiplicity) of solutions of any given system of $n$ polynomials in $n$ variables. The algorithm is inductive; it consists of finitely many steps, and at each step a non-degeneracy criterion determines if the correct number has been computed. The estimate and non-degeneracy criterion for the number of solutions on $\kk^n$ from \cref{affine-chapter} of this book serve as the initial step of that algorithm. 

\section{Prerequisites}
We tried to ensure that this book is accessible to someone with the mathematical maturity and algebra background of a second year mathematics graduate student. In the ideal case a reader would be familiar with the properties of algebraic varieties discussed in \cref{var-chapter}, so that (s)he could start with toric varieties in \cref{toric-part} and only refer to results from \cref{pretoric-part} if necessary. However, \cref{pretoric-part} is self contained (modulo the dependencies explicitly stated in \cref{noproof-section,rcurve-section} and some commutative algebra results stated in \cref{algebra-section,proof-scheme-section}) - with proper guidance it can be used as the material for a first course in algebraic geometry. One possible strategy for such a course would be to cover the chapters on algebraic varieties (\cref{var-chapter}), toric varieties (\cref{toric-intro}), Bernstein-Kushnirenko theorem (\cref{bkk-chapter}) and (weighted) B\'ezout's theorem (\cref{bezout-chapter}). The chapters on intersection multiplicity (\cref{mult-chapter}) and polytopes (\cref{appolytopes}) are included for completion - in a first course the required results from these chapters can simply be explained, perhaps via examples and/or pictures, instead of working out the details of the proofs. In particular, the proofs (and exercises) given in \cref{appolytopes} (polytopes) are elementary and a student should not have much difficulty in following them. The most sophisticated part of \cref{mult-chapter} (intersection multiplicity) is the concept of a ``closed subscheme'' of a variety and the fact that it can be locally defined by ideals determined by regular functions; the other results are basic facts about intersection multiplicity of $n$ regular functions at a nonsingular point $a$ of an $n$ dimensional variety (e.g.\ that they can be defined via the ``order'' at $a$ of one of the functions along the curve defined by the other functions) and relevant properties of the ``order'' function at a point on a (possibly non-reduced) curve. While the proofs use somewhat complicated algebra, the statements are intuitive, at least if one has some familiarity with basic properties of (complex) analytic functions.

\section{Organization}
\Cref{pretoric-part} and the first chapter of \cref{toric-part} have been designed as parts of a textbook, with many exercises and examples. The goal was to develop efficiently (and in an elementary way) the theory needed to prove the results in the subsequent chapters. These latter chapters are more like those of a monograph; there are no exercises, but they do contain a number of examples. We now give a short description of each chapter. \\

In \cref{var-chapter} we develop the required theory of algebraic varieties. We tried to stress the geometric point of view where possible. A number of results have been developed through exercises; ample hints have been provided to ensure that no single step of any exercise is very difficult. %
In \cref{mult-chapter} we describe basic properties of intersection multiplicity (of $n$ regular functions at a nonsingular point of an $n$ dimensional variety), in particular how it can be computed using curves. After giving simple examples to illustrate that a satisfactory treatment of intersection multiplicity would need to incorporate non-reduced rings, we give a short introduction to ``closed subschemes of a variety''\footnote{We decided to omit definitions of general sheaves and schemes since we do not use these notions anywhere in this book. On the other hand, once one really understands the special cases of ``sheaves of ideals'' and ``closed subschemes of a variety,'' which are discussed in \cref{mult-chapter}, the leap to the general notions will be natural.}. A number of examples presented in \cref{var-chapter,mult-chapter} were taken from answers to the question {\em Algebraic geometry examples} \cite{mathoverflow-algeom-examples} posed by R.\ Borcherds on {\em MathOverflow}. \Cref{appolytopes} is a compilation (with complete proof) of the properties of convex polyhedra which, together with the results of \cref{var-chapter}, constitute the foundation on which we introduce toric varieties in \cref{toric-intro}. In \cref{toric-intro} we mainly discuss those properties of toric varieties which are required for the results in the subsequent chapters. In \cref{bkk-chapter} we prove Bernstein's theorem and present some of its basic applications to convex geometry. In \cref{bezout-chapter} we apply Bernstein's theorem to prove the weighted homogeneous and multi-homogeneous versions of B\'ezout's theorem. \Cref{multiplicity-chapter} contains the results on the generic bound and non-degeneracy conditions for intersection multiplicity at the origin, which we use in \cref{affine-chapter} to compute the generic bounds and non-degeneracy conditions for the number of solutions of polynomial systems on $\kk^n$. It turns out that one can as easily replace $\kk^n$ by an arbitrary Zariski open subset of $\kk^n$ - the results of \cref{affine-chapter} are derived in this greater generality. In \cref{affine-chapter} we also use the main results to derive generalizations of weighted homogeneous and multi-homogeneous versions of B\'ezout's theorem applicable to weighted degrees with possibly zero or negative weights. In \cref{milnor-chapter} we apply the results from \cref{multiplicity-chapter} to the study of Milnor numbers; in particular, we derive and generalize classical results of Kushnirenko on Milnor numbers. \Cref{bkk-chapter,affine-chapter,milnor-chapter} end with selections of open problems (mostly combinatorial in nature).

\newpage
\thispagestyle{plain}
\begin{center}
{\bf Chapter Dependencies}
\vspace{1cm}

\newcommand{\chapterwidth}{0.24\textwidth}

\begin{tikzcd} [row sep = 1.2cm, column sep = 1cm]
&
\parbox{\chapterwidth}{\centering\Cref{appolytopes}\\Polytopes}
\arrow{d}
&
\\
\parbox{\chapterwidth}{\centering \Cref{var-chapter}\\Varieties}
\arrow[d]
\arrow{r}
&
\parbox{\chapterwidth}{\centering \Cref{toric-intro}\\Toric varieties}
\arrow{d}
\\
\parbox{\chapterwidth}{\centering \Cref{mult-chapter}\\Intersection multiplicity}
\arrow{r}
&
\parbox{\chapterwidth}{\centering \Cref{bkk-chapter}\\Generic number of zeroes on $\nktorus$ (Bernstein-Kushnirenko theorem)}
\arrow{d}
\arrow{r}
&
\parbox{\chapterwidth}{\centering \Cref{bezout-chapter}\\(Weighted) B\'ezout's theorem}
\\
&
\parbox{\chapterwidth}{\centering \Cref{multiplicity-chapter}\\ Generic intersection multiplicity}
\arrow{d}
\arrow{r}
&
\parbox{\chapterwidth}{\centering \Cref{milnor-chapter}\\ Generic Milnor number}
\\
&
\parbox{\chapterwidth}{\centering \Cref{affine-chapter}\\ Generic number of zeroes on $\kk^n$}
&
\end{tikzcd}
\label{dependent-label}
\end{center} 

\newpage
\thispagestyle{plain}

\mbox{}
\vspace{3cm}

\begin{center}
\begin{figure}[h]
\includegraphics[width=0.9\textwidth]{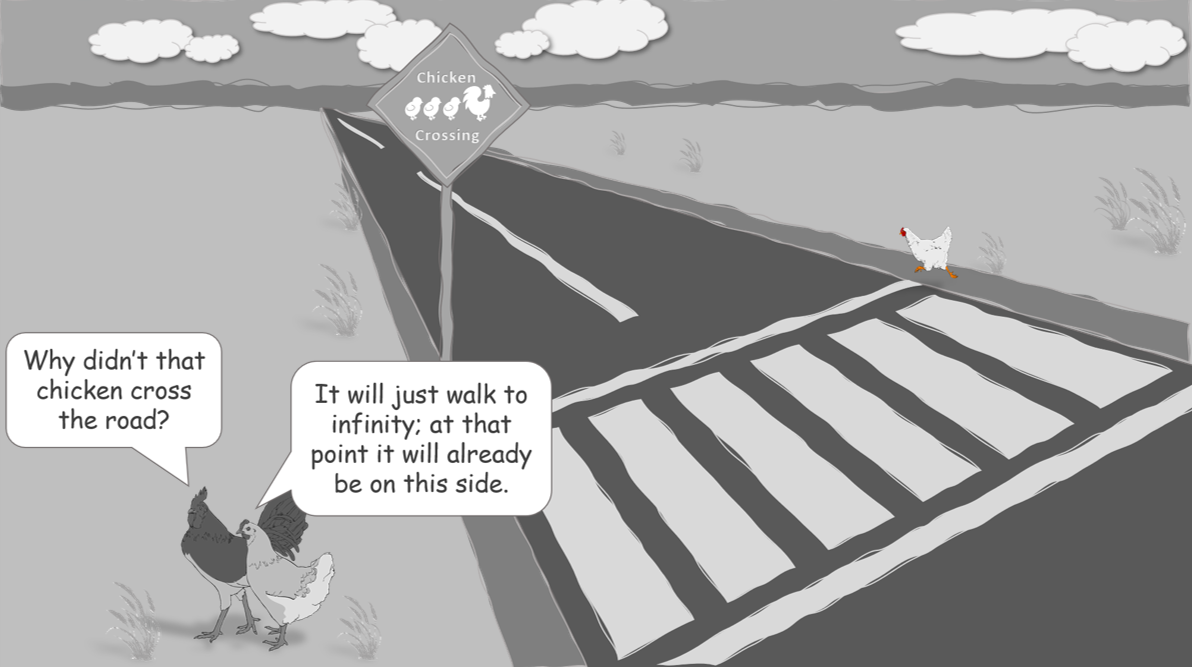}
\caption*{On the projective space chickens have more than one way of crossing roads}
\label{fig:chicken}
\end{figure}
\end{center}

\chapter{A brief history of points at infinity in geometry} \label{brief-chapter}
In this chapter we give a brief historical overview of the concept of points at infinity in geometry and the subsequent introduction of homogeneous coordinates on projective spaces.

\section{Points at infinity}
Points at infinity seem to have first cropped up in Johannes Kepler's work on conics in {\em Ad Vitellionem paralipomena quibus astronomiae pars optica traditur}\footnote{``Literally `Things omitted by' (or `Supplements to') `Witelo with which the optical part of astronomy is concerned'. \textelp{} Witelo's {\em Perspectiva}, probably written in the 1270s, appeared in several new editions in the sixteenth century, and seems to have been the standard textbook on Optics'' \cite[pp.\ 221-222]{field-gray}.} (1604). It is in this text that Kepler introduces the term {\em focus}\footnote{{\em Focus} is the Latin word for hearth. ``Since light was reflected to the focus, \ldots\ the focus of the mirror was the position in which one would place the material one wished to burn'' \cite[p.\ 222]{field-gray}.} to denote each of the (unique) pair of points inside a conic such that the rays from any point on the conic make equal angles to the tangent at that point. For a circle the foci coincide at the center, and they separate as the circle deforms into an ellipse. As one continues to deform the ellipse so that in the end it turns into a parabola, Kepler concludes that ``In the Parabola one focus \ldots\ is inside the conic section, the other to be imagined either inside or outside, lying on the axis at an infinite distance from the first, so that if we draw the straight line \ldots\ from this blind focus to any point \ldots\ on the conic section, the line will be parallel to the axis \ldots'' \cite[pp.\ 186-187]{field-gray}, see \cref{fig:foci}.

\begin{figure}[h]
\begin{center}

\begin{tikzpicture}[scale=0.33]

\def\colorzero{green}
\def\colorone{red}
\def\colortwo{blue}

\def\a{2}
\def\t{0.78}
\def\anglezero{15}
\def\angleone{75}
\def\ytwo{6}
\def\ctangentfactor{0.5}
\def\etangentfactor{3}
\def\ptangentfactor{3}

\def\xmax{15}
\def\ypar{9}
\def\ymax{9}

\pgfmathsetmacro\rc{2*\a};
\pgfmathsetmacro\rex{2*\a/\t^2};
\pgfmathsetmacro\rey{2*\a/\t};
\pgfmathsetmacro\reshift{sqrt(1-\t^2)};
\pgfmathsetmacro\cnormslope{-cot(\anglezero)}
\pgfmathsetmacro\enormslope{-cot(\angleone)*\t}
\pgfmathsetmacro\pnormslope{2*\a/\ytwo}

\tikzstyle{dot} = [circle, minimum size=4pt, inner sep = 0pt, fill]
\tikzstyle{smalldot} = [circle, minimum size=0pt, inner sep = 0pt, fill]

\draw[\colorzero, thick] (\rc, 0) circle (\rc);
\draw[\colorone, thick] (\rex, 0) ellipse ({\rex} and {\rey});
\draw[\colortwo, thick] plot[smooth, domain=-\ypar:\ypar, variable=\y] ({(\y)^2/(4*\a)}, \y);

\node[dot, \colorzero] (c) at (2*\a,0) {};
\node[dot, \colorone] (e1) at ({\rex*(1- \reshift)},0) {};
\node[dot, \colorone] (e2) at ({\rex*(1 + \reshift)},0) {};
\node[dot, \colortwo] (p) at (\a,0) {};

\node[smalldot] (cp) at ($(c) +(\anglezero:\rc)$) {};
\node[smalldot] (ep) at ($(\rex,0) +(\angleone:{\rex} and {\rey})$) {};
\node[smalldot] (pp) at ({(\ytwo)^2/(4*\a)}, \ytwo) {};

\draw[thick, dashed] (c) -- (cp);
\draw[thick, dashed] (e1) -- (ep) -- (e2);
\draw[thick, dashed] (p) -- (pp) -- (\xmax,\ytwo);

\draw[thick] ( $(cp) + \ctangentfactor*(1, \cnormslope)$) -- ( $(cp) - \ctangentfactor*(1, \cnormslope)$)  ;
\draw[thick] ( $(ep) + \etangentfactor*(1, \enormslope)$) -- ( $(ep) - \etangentfactor*(1, \enormslope)$)  ;
\draw[thick] ( $(pp) + \ptangentfactor*(1, \pnormslope)$) -- ( $(pp) - \ptangentfactor*(1, \pnormslope)$)  ;

\draw (0,-\ymax)  -- (0, \ymax);
\draw (0,0) -- (\xmax,0);

\end{tikzpicture}
\caption{Foci of families of conics}
\label{fig:foci}
\end{center}
\end{figure}

Another inspiration, albeit indirect, of points at infinity is the theory of linear perspective. Application of perspectives were already present in early fourteenth century paintings from Italy \cite[Chapter I]{andersen-geometry-art}, the earliest surviving written account of geometric construction of perspective being Leon Battista Alberti's {\em De Pictura} (1435). By the seventeenth century there were numerous treatises on perspective. In 1639 Girard Desargues, who had worked as a military engineer and written on perspective, circulated fifty copies of his {\em Brouillon project d'une atteinte aux evenmens des rencontres du cone avec un plan} (``Rough draft of an essay on the results of taking plane sections of a cone''). At the very beginning of {\em Brouillon project} Desargues introduced the notion that parallel lines intersect at a point at infinity and parallel planes intersect at a line at infinity; constructing essentially the projective plane $\pp^2(\rr)$ and the three dimensional projective space $\pp^3(\rr)$ over $\rr$. He made extensive use of the lines and planes at infinity to give a unified treatment of families of parallel lines and families of lines through a common point. The subject of projective geometry was born in {\em Brouillon project}.

\section{Homogeneous coordinates}
The birth however, went practically unnoticed. Desargues's manuscript was thought to have been lost and it did not inspire much new work (other than Blaise Pascal's {\em Essay pour les coniques} (1640) which contains Pascal's theorem on conics). Projective geometry was revived in the nineteenth century largely due to Jean-Victor Poncelet, who fought in Napoleon's army in the battle of Krasnoi in November, 1812, and then was a prisoner of war in Saratov till Napoleon's defeat in mid 1814. In the prison ``he occupied himself summarising all he knew of the mathematical sciences in notebooks that he then distributed to his fellow prisoners who wanted to finish an education disrupted by the incessant military campaigns'' \cite[p.\ 13]{gray}. In the process he discovered, and upon his return to France, championed, the unifying aspect of projective geometry (as opposed to the ``analytic geometry'' of Ren\'e Descartes). A fundamental tool of this new geometry was the duality between points and lines on the plane. Initially applied by Charles Julien Brianchon and Poncelet to conics, the duality principle was extended to all planar curves by Joseph Diaz Gergonne\footnote{``\ldots\ it is one thing to realise that dualising a figure is a good way to obtain new theorems, which is what Poncelet did, and quite another thing to claim that points and lines are interchangeable concepts which must logically be treated on a par. This was the view that Gergonne put forward in 1825. Interpreted in such generality, Gergonne's principle of duality is one of the most profound and simple ideas to have enriched geometry since the time of the Greeks \ldots'' \cite[p.\ 55]{gray}.}. All the details of the duality principle however were not clear, e.g.\ even though the principle suggests that dualising twice one should get back the original curve, it was soon discovered that dualising a curve of degree higher than two results in a curve of degree higher than that of the original curve. Poncelet had some ideas about resolving this paradox by taking into account the effects of cusps and double points on a curve, but his ideas were not very precise. The resolution came through the algebraic treatment of projective geometry by August M\"obius in {\em Der Barycentrische Calc\"ul} (1827). M\"obius observed that weights $w_0,w_1$ placed at the ends of a (weightless) rod uniquely determines a point $P$ on the rod, namely their \index{Barycenter}{\em Barycenter}, i.e.\ the center of gravity; the ordered pair $[w_0:w_1]$ (we write it in this way to distinguish from the Cartesian coordinates of $P$) are the \index{Barycenter!Barycentric coordinate}{\em Barycentric coordinates} of $P$. It is straightforward to work out the relation between the Cartesian and barycentric coordinates, e.g.\ if we identify the rod with the closed interval $[a,b]$ on the real line, then the barycentric coordinsates $[w_0:w_1]$ of $x \in [a,b]$ satisfies:
\begin{align*}
x &= \frac{aw_0 + bw_1}{w_0 + w_1}
\end{align*}
This formula can be readily extended to allow for $w_0$ and $w_1$ to be zero or negative. It follows that each point on the real line has barycentric coordinates, see \cref{fig:barycentric-1}. It is also clear that the barycentric coordinates are {\em homogeneous}, i.e.\ $[w_0\lambda: w_1 \lambda]$ denote the same point as $[w_0: w_1]$ for every nonzero $\lambda \in \rr$. Finally, note that if $w_0 + w_1 = 0$, then $[w_0: w_1]$ does not correspond to any point on the line; M\"obius defined it as a point lying at infinity.

\begin{center}
\begin{figure}[h]
\def\scalefactor{2.5}
\def\xmin{-2}
\def\xmax{3}
\def\ly{0.2}
\def\picfontsize{\small}
\tikzstyle{dot} = [red, circle, minimum size=4pt, inner sep = 0pt, fill]

\begin{tikzpicture}[scale=\scalefactor]
\draw [<->] (\xmin,0) -- (\xmax,0);
\draw[blue,  ultra thick] (0,0) -- (1,0);

\node[dot] (-1) at (-1,0) {};
\node[below = \ly of -1] {\picfontsize $-1$};
\node[above = \ly of -1] {\picfontsize $[2:-1]$};

\node[dot] (origin) at (0,0) {};
\node[below = \ly of origin] {\picfontsize $0$};
\node[above = \ly of origin] {\picfontsize $[1:0]$};

\node[dot] (half) at (0.5,0) {};
\node[below = \ly of half] {\picfontsize $0.5$};
\node[above = \ly of half] {\picfontsize $[0.5:0.5]$};

\node[dot] (1) at (1,0) {};
\node[below = \ly of 1] {\picfontsize $1$};
\node[above = \ly of 1] {\picfontsize $[0:1]$};

\node[dot] (2) at (2,0) {};
\node[below = \ly of 2] {\picfontsize $2$};
\node[above = \ly of 2] {\picfontsize $[-1:2]$};

\end{tikzpicture}
\caption{Barycentric coordinates in dimension one with respect to the interval $[0,1]$}
\label{fig:barycentric-1}
\end{figure}
\end{center}

In dimension two one starts with a triangle $\Delta$; assume for convenience that the vertices of $\Delta$ are the points with Cartesian coordinates $(0,0)$, $(1, 0)$ and $(0,1)$. Then the Cartesian coordinates $(x,y)$ and the barycentric coordinates $[w_0:w_1:w_2]$ of a point $P$ on the plane with respect to $\Delta$ are related as follows:
\begin{align*}
x = \frac{w_1}{w_0+w_1+w_2}, \quad y = \frac{w_2}{w_0 + w_1 + w_2}
\end{align*}

\begin{center}
\begin{figure}[h]
\def\scalefactor{1.5}
\def\xmin{-2}
\def\xmax{2}
\def\ymin{-1}
\def\ymax{2}
\def\ly{0.2}
\def\opazero{0.5}
\def\picfontsize{\small}
\def\colorzero{green}
\def\colorone{blue}

\tikzstyle{dot} = [red, circle, minimum size=4pt, inner sep = 0pt, fill]

\begin{tikzpicture}[scale=\scalefactor]
\draw [gray,  line width=0pt] (\xmin, \ymin) grid (\xmax,\ymax);
%

\foreach \x in {\xmin, ..., \xmax}
	\node[anchor=north] at (\x,\ymin) {\x};
\foreach \y in {\ymin, ..., \ymax}
	\node[anchor=east] at (\xmin, \y) {\y};

\draw [<->] (0, \ymax) |- (\xmax, 0);
\draw[ultra thick, \colorone, fill=\colorzero, opacity=\opazero ] (0,0) --  (0,1) -- (1,0) -- cycle;

\node[dot] (origin) at (0,0) {};
\node[anchor=north] at (origin) {\picfontsize $[1:0:0]$};

\node[dot] (a) at (1,0) {};
\node[anchor=north] at (a) {\picfontsize $[0: 1:0]$};

\node[dot] (b) at (0,1) {};
\node[anchor=south] at (b) {\picfontsize $[0:0:1]$};

\node[dot] (c) at (1,1.5) {};
\node[anchor=south] at (c) {\picfontsize $[-1.5:1.5:1]$};

\node[dot] (d) at (-1,1.5) {};
\node[anchor=south] at (d) {\picfontsize $[0.5:-1:1.5]$};

\node[dot] (e) at (-1,-0.5) {};
\node[anchor=north] at (e) {\picfontsize $[2.5:-1:-0.5]$};

\node[dot] (f) at (1,-0.5) {};
\node[anchor=north] at (f) {\picfontsize $[0.5:1:-0.5]$};

\end{tikzpicture}
\caption{Barycentric coordinates in dimension two with respect to $\Delta$}
\label{fig:barycentric-2}
\end{figure}
\end{center}

As in the case of the real line, the barycentric coordinates of the points at infinity are $[w_0:w_1:w_2]$ with $w_0 + w_1 + w_2 = 0$. M\"obius observed that many computations with barycentric coordinates become simpler upon a change of coordinates of the form
\begin{align*}
[w_0:w_1:w_2] \mapsto [w_0+w_1+w_2: w_1: w_2]
\end{align*}
These new coordinates are nowadays usually denoted as {\em homogeneous coordinates}. In particular, the equation of the line $ax + by + c = 0$ changes in the barycentric coordinates to $aw_1 + bw_2 + c(w_0 + w_1 + w_2) = 0$, and this in turn becomes $aw_1 + bw_2 + cw_0 = 0$ in homogeneous coordinates. And in homogeneous coordinates the points at infinity are described by $w_0 = 0$. As he was finishing {\em Der Barycentrische Calc\"ul}, M\"obius heard of the duality between points and lines studied by the French geometers, and noticed that the homogeneous coordinates gives a natural algebraic approach to duality, namely the line $aw_1 + bw_2 + cw_0 = 0$ corresponds simply to the point with homogeneous coordinates $[a:b:c]$, and vice versa. This automatically ensured that concurrent lines go to collinear points under duality and that dualizing twice one gets back to the original curve. Julius Pl\"ucker, possibly independently of M\"obius, gave an analogous theory of homogeneous coordinates in 1830, and later used it to completely resolve the duality paradox. The homogeneous coordinates were soon extended to higher dimensions, which opened the door to algebraic study of higher dimensional projective spaces.

\section{Projective space} \label{projective-intro}
Take an arbitrary field $\kk$ and fix coordinates $(x_0, \ldots, x_n)$ on $\kk^{n+1}$, $n \geq 0$. The \index{Projective!space}{\em $n$-dimensional projective space $\pp^n$ over $\kk$} is the set of lines (with respect to $(x_0, \ldots, x_n)$) in $\kk^{n+1}$ through the origin. Every point $(a_0, \ldots, a_n) \in \kk^{n+1}\setminus \{0\}$ determines a unique line through the origin which passes through it; the \index{Homogeneous!coordinate}{\em homogeneous coordinate} of this line is $[a_0: \cdots : a_n]$. For each $j=0, \ldots, n$, let $U_j := \{[a_0: \cdots: a_n]: a_j \neq 0\}$. The map 
\begin{align*}
(a_1, \ldots, a_n) \mapsto [a_1: \cdots: a_{j-1}: 1: a_j: \cdots : a_n]
\end{align*}
gives a one-to-one correspondence between $\kk^n$ and $U_j$. In the case that $\kk = \rr$ or $\cc$, one can use this correspondence to induce a topology on $U_j$ (by declaring a subset of $U_j$ to be open if and only if its pre-image in $\rr^n$ or $\cc^n$ is open). It is straightforward to check that these topologies are compatible (i.e.\ they induce the same topology on their intersections), and accordingly turns $\pp^n = \bigcup_{j=0}^n U_j$ into a manifold. For a general $\kk$, the usual topology put on $\kk^n$ is the \index{Zariski!topology}{\em Zariski topology}, in which the closed subsets are zero-sets of systems of polynomials, and the identification of the $U_j$ with $\kk^n$ is used to give $\pp^n$ the structure of an {\em algebraic variety over $\kk$}. We review algebraic varieties and Zariski topology in \cref{var-chapter}. The complement $H_j$ of $U_j$ in $\pp^n$ is the set of lines (through the origin) which lie on the $j$-th coordinate hyperplane, so that if we identify $\kk^n$ with $U_0$, then the set of points at infinity, i.e.\ the complement of $\kk^n$ in $\pp^n$, is precisely $H_0 := \{[0: a_1: \cdots: a_n] \}$, and the homogeneous coordinates on $\pp^n$ are precisely those introduced by M\"obius. Note that $H_0$ is naturally isomorphic to $\pp^{n-1}$. 

\begin{prop} \label{prop:p^n-line-closure}
The closure in $\pp^n$ of each straight line on $\kk^n$ has a unique point at infinity. Two coplanar lines intersect at a common point at infinity if and only if they are parallel.
\end{prop}

\begin{proof}
Here we treat the case that $\kk = \cc$ and the topology on $\pp^n$ is that induced from the Euclidean topology on $\cc^n \cong \rr^{2n}$; see \cref{exercise:projective-closure-example} for the case of general $\kk$ and Zariski topology on $\pp^n$. Let $L = \{(a_1, \ldots, a_n) + t(b_1, \ldots, b_n): t \in \cc\}$ be a line on $\cc^n \cong U_0$ (note that this means $(b_1, \ldots, b_n) \neq (0, \ldots, 0)$). In homogeneous coordinates  
\begin{align*}
L = \{[1: a_1 + tb_1 : \cdots :a_n + tb_n]: t \in \cc\}
\end{align*}
Therefore the set of all points on $\cc^{n+1}$ which correspond to points on $L$ is 
\begin{align*}
L' = \{(s, sa_1 + stb_1, \ldots, sa_n + stb_n): s, t \in \cc, s \neq 0\}
\end{align*}
Since $(b_1, \ldots, b_n) \neq 0$, it is straightforward to check that the closure $\bar L'$ of $L'$ in $\cc^{n+1}$ is the plane spanned by $(1, a_1, \ldots, a_n)$ and $(0, b_1, \ldots, b_n)$ (see \cref{fig:parallel-intersection}). The points at infinity on the closure of $L$ in $\pp^n$ correspond to the points $(x_0, \ldots, x_n)$ on $\bar L' \setminus \{0\}$ with $x_0 = 0$, i.e.\ the set of points $(0, \lambda b_1, \ldots, \lambda b_n)$, $\lambda \in \cc \setminus \{0\}$. Since all these points correspond to the single point $[0: b_1: \cdots : b_n]$ on $\pp^n$, this proves both assertions of the proposition. 
\end{proof}

\def\colorzero{blue}
\def\opazero{0.5}
\def\viewx{120}
\def\viewy{30}
\def\scalefactor{1}
\def\colorzero{green}
\def\colorone{blue}
\def\colortwo{red}

\def\height{2}
\def\length{4}
\def\slope{1.5}
\def\lengthzero{8}
\pgfmathsetmacro{\startpoint}{\slope * \height}
\pgfmathsetmacro{\endpoint}{\startpoint + \length}
\pgfmathsetmacro{\endpointone}{\slope + \length}

\begin{figure}[h]
\begin{center}
\begin{tikzpicture}[scale=\scalefactor]
\pgfplotsset{every axis/.append style = {view={\viewx}{\viewy}, axis lines=middle, enlargelimits={upper}}}

\begin{axis}[
	 xlabel = {\picfontsize $x_1$},
	 ylabel = {\picfontsize $x_2$},
	 zlabel = {\picfontsize $x_0$}
]
	\addplot3 [fill=\colorzero, opacity=\opazero, thick] coordinates{(0,0,1) (0,\lengthzero,1) (\lengthzero,\lengthzero,1) (\lengthzero,0,1)};
	\addplot3 [fill=\colorone, opacity=\opazero, thick] coordinates{(0,0,0) (\startpoint,0,\height) (\endpoint,\length,\height) (\length,\length,0)};
	\addplot3 [fill=\colortwo, opacity=\opazero, thick] coordinates{(0,0,0) (0,\startpoint,\height) (\length,\endpoint,\height) (\length,\length,0)};
	\addplot3 [thick] coordinates{(\slope,0,1) (\endpointone,\length,1)};
	\addplot3 [thick] coordinates{(0,\slope,1) (\length,\endpointone,1)};
	\addplot3 [thick, dashed] coordinates{(0,0,0) (\length,\length,0)};
\end{axis}

\end{tikzpicture}
\end{center}
\caption{$U_0$ can be identified with the hyperplane $x_0 = 1$ on $\kk^{n+1}$. Every line on $U_0$ corresponds to a plane on $\kk^{n+1}$. Parallel lines on $U_0$ intersect at a point at infinity on $\pp^n$ since the corresponding planes on $\kk^{n+1}$ intersect along a line through the origin on $x_0 = 0$.} \label{fig:parallel-intersection}
\end{figure}

\Cref{prop:p^n-line-closure} shows that the projective space incorporates the intuition from the theory of perspectives that two parallel lines intersect at a point at infinity. The connection of projective spaces with the affine B\'ezout problem comes from the following property (which you will see in \cref{example:principal-closure}): if $\kk$ is algebraically closed, then for each polynomial $f \in \kk[x_1, \ldots, x_n]$, there is a correspondence between the following sets:
\begin{align}
\parbox{0.3\textwidth}{
points at infinity on the closure in $\pp^n$ of $\{x \in \kk^n: f(x) = 0\}$
}
\quad
\longleftrightarrow
\quad
\{x \in \kk^n\setminus\{\origin\}: \ld(f)(x) = 0\}
\label{bezout-infinity-correspondence}
\end{align}
where $\ld(f)$ is the leading form of $f$. This correspondence provides the geometric explanation for the condition from B\'ezout's theorem (\cref{bezout}) under which the B\'ezout estimate for the number of solutions is exact. A proof of B\'ezout's theorem usually consists of showing that
\begin{itemize}
\item given $n$ polynomials $f_1, \ldots, f_n \in \kk[x_1, \ldots, x_n]$, the number of isolated points (counted with appropriate multiplicity) on the intersection of the closures in $\pp^n$ of $\{f_i = 0\}$ is at most the product of the degrees of the $f_i$, and
\item if $\kk$ is algebraically closed, then this bound is attained with points inside $\kk^n$ if and only if there is no ``intersection at infinity'' on $\pp^n$ i.e.\ 
\begin{align*}
\bigcap_i \overline{\{f_i = 0\}} = \bigcap_i\{f_i = 0\} 
 \end{align*}
\end{itemize}
Most of the mathematics of this book takes place on {\em toric varieties} (introduced in \cref{toric-intro}), a class of algebraic varieties of which the projective space is a special case. We will study natural analogues of B\'ezout's theorem on different classes of toric varieties, and B\'ezout's theorem will fall out as a special case (\cref{bezout-1}).  

\part{Preliminaries} \label{pretoric-part}
\chapter{Quasiprojective varieties over algebraically closed fields} \label{var-chapter}
In this chapter we give a quick introduction to algebraic varieties, focusing mainly on the properties used in \cref{toric-part,postoric-part}. \Cref{algebra-section} includes a discussion of the relevant concepts and results from commutative algebra including Hilbert's basis theorem (\cref{thm:Hilbert-basis}) and Nullstellensatz (\cref{thm:Hilbert-nulls}) which are part of the foundation of modern algebraic geometry. There are a few results including Krull's ``principal ideal theorem'' (\cref{thm:principal-ideal}) that we use without proof - these are listed in \cref{noproof-section}. Throughout this book $\kk$ denotes an algebraically closed field. We start the discussion with {\em affine varieties} over $\kk$. General {\em quasiprojective varieties} are studied from \cref{quasi-section} onward. Unless explicitly stated otherwise, in this book a {\em variety} will mean a quasiprojective variety over $\kk$.


\section{Affine varieties} \label{affection}
An \index{Affine!space}{\em affine space} over $\kk$ is simply the set $\kk^n$ of $n$-tuples of elements from $\kk$ for some $n \geq 0$. Fix a system of coordinates $(x_1, \ldots, x_n)$ on $\kk^n$, i.e.\ $\kk^n = \{(x_1, \ldots, x_n): x_1, \ldots, x_n \in \kk\}$. A \index{Subvariety!of an affine space}{\em subvariety} of $\kk^n$ is the set of zeroes of a collection of polynomials in $(x_1, \ldots, x_n)$\footnote{Traditionally in the definition of a subvariety it was common to include also the requirement that a subvariety $V$ should be {\em irreducible}, i.e.\ whenever $V = V_1 \cup V_2$, where $V_1, V_2$ are sets of zeroes of systems of polynomials, then either $V = V_1$ or $V = V_2$; however, our definition is also widely used now.}. An \index{Affine!variety}\index{Variety!affine}{\em affine variety} is simply the subvariety of an affine space. A \index{Hypersurface}{\em hypersurface} of $\kk^n$ is the set of zeroes of a single polynomial.

\begin{center}
\begin{figure}[htb]
\begin{subfigure}[b]{0.2\textwidth}
\includegraphics[width=25mm]{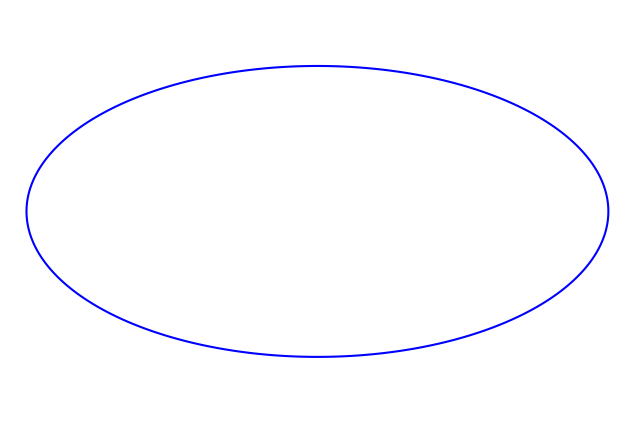}
\caption{$x^2 + 4y^2= 9$}
\label{fig:ellipse}
\end{subfigure}\hspace{0.1\textwidth}
\begin{subfigure}[b]{0.2\textwidth}
\includegraphics[width=3cm]{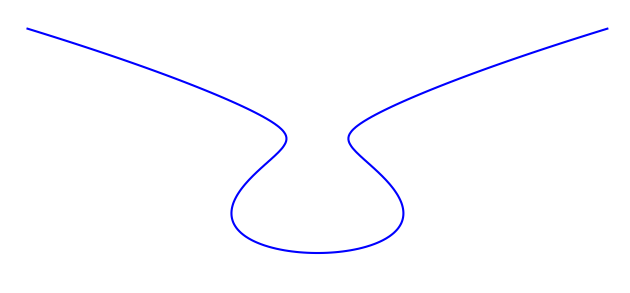}
\caption{$x^2 = y^3 - 4y + 4$}
\label{fig:elliptic-curve}
\end{subfigure}\hspace{0.1\textwidth}
\begin{subfigure}[b]{0.36\textwidth}
\includegraphics[width=2cm]{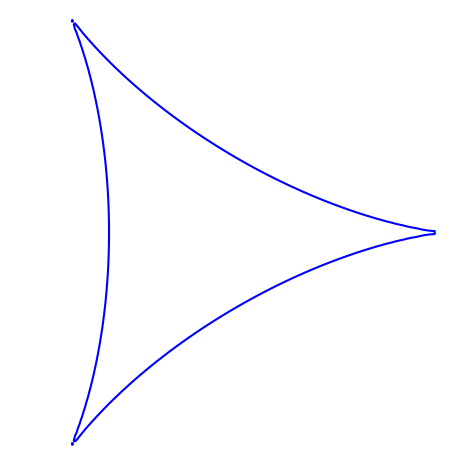}
\caption{$(x^2 + y^2 + 9)^2  - 8(x^3 - 3xy^2) = 108$}
\label{fig:deltoid}
\end{subfigure}
\caption{Real points of some curves on $\cc^2$}
\label{fig:planar curves}
\end{figure}
\end{center}

\begin{example} \label{example:subvarieties-aff1}
Since every nonzero polynomial in a single variable has only finitely many zeroes, every proper subvariety of $\kk$ consists of finitely many points. Hypersurfaces of $\kk^2$ are special cases of (affine) {\em algebraic curves}\footnote{A {\em curve} is a variety of {\em dimension} one. We discuss dimensions in \cref{dimension-section}.}. See \cref{fig:planar curves} for pictures of real points of some curves on $\cc^2$. The curve in \cref{fig:elliptic-curve} is an ``elliptic curve'' and the one in \cref{fig:deltoid} is a ``deltoid'' which was investigated by L.\ Euler in 1745 in relation to a problem in optics; we refer to \cite[Chapter I]{bries-horst-curves} for many pictures of algebraic curves along with their history.
\end{example}

Given a subset $\qqq$ of $\kk[x_1, \ldots, x_n]$, we write $V(\qqq)$ for the subvariety determined by $\qqq$, i.e.\ $V(\qqq) = \{a \in \kk^n : f(a) = 0$ for all $f \in \qqq\}$. Hilbert's Basis Theorem (\cref{thm:Hilbert-basis}) implies that the ideal $\langle \qqq \rangle$ of $\kk[x_1, \ldots, x_n]$ generated by $\qqq$ is in fact generated by {\em finitely many} polynomials $f_1, \ldots, f_k$. But then it is immediate to check that $V(\qqq)$ is precisely the set of common zeroes of $f_1, \ldots, f_k$, i.e.\ $V(\qqq) = V(f_1, \ldots, f_k)$. This proves the fundamental fact that every affine variety is determined by {\em finitely many} polynomials:

\begin{prop} \label{prop:affinite-basis}
Every subvariety of $\kk^n$ is the set of common zeroes of finitely many polynomials. \qed
\end{prop}

It is straightforward to check that there is a unique topology on $\kk^n$ for which the closed sets are precisely the subvarieties on $\kk^n$ (\cref{exercise:zariski-topology-aff}); it is called the \index{Zariski!topology}{\em Zariski topology}. In this book $\kk^n$ (and its subsets) will always be assumed to be equipped with the Zariski topology; in particular, by ``closed'' or ``open'' subsets we will mean Zariski closed or Zariski open subsets. If $X$ is a subset of $\kk^n$, we write $I(X)$ for the set of all polynomials $f \in \kk[x_1, \ldots, x_n]$ such that $f(a) = 0$ for all $a \in X$. It is straightforward to check that $I(X)$ is a {\em radical ideal} of $\kk[x_1, \ldots, x_n]$ (\cref{exercise:radicalIX}).

\begin{example} \label{example:Isingleton}
If $X$ is a singleton consisting of a single point $a = (a_1, \ldots, a_n) \in \kk^n$, then $I(X)$ is the ideal $\mmm_a$ of $\kk[x_1, \ldots, x_n]$ generated by $x_1 - a_1, \ldots, x_n - a_n$. Indeed, it is clear that $\mmm_a \subseteq I(X)$. Since $I(X)$ is a proper ideal of $\kk[x_1, \ldots, x_n]$ and $\mmm_a$ is a {\em maximal} ideal of $\kk[x_1, \ldots, x_n]$ (since $\kk[x_1, \ldots, x_n]/ \mmm_a \cong \kk$), it follows that $I(X) = \mmm_a$.
\end{example}

\begin{example} \label{example:Ikn}
If $X = \kk^n$, then $I(X) = 0$. Indeed, since an algebraically closed field is infinite, no nonzero polynomial in $(x_1, \ldots, x_n)$ vanishes at all points of $\kk^n$ (\cref{exercise:pol-zero-on-Kn}), so that $I(\kk^n) = 0$.
\end{example}

Some basic properties of the operators $V(\cdot)$ and $I(\cdot)$ are presented in \crefrange{exercise:V-containment}{exercise:zariski-closure-aff} - the reader is urged to go over them. Given an ideal $\qqq$ of $\kk[x_1, \ldots, x_n]$, it can be seen directly from the definitions that $V(\sqrt{\qqq}) = V(\qqq)$ and $I(V(\qqq)) \supseteq \sqrt{\qqq}$ (e.g.\ see \cref{exercise:power-V,exercise:I-containment}); Hilbert's Nullstellensatz (\cref{thm:Hilbert-nulls}) implies that the latter containment is in fact an equality:
\begin{align}
I(V(\qqq)) &= \sqrt{\qqq} \label{eqn:nullstellensatz}
\end{align}

The Nullstellensatz sets up the basic correspondence between algebra and geometry underlining algebraic geometry over algebraically closed fields:

\begin{thm} \label{thm:affine-nullspondence}
There is a one-to-one correspondence between subvarieties of $\kk^n$ and radical ideals of $\kk[x_1, \ldots, x_n]$ given by $I(\cdot)$ and $V(\cdot)$. Given a subvariety $X$ of $\kk^n$ and a radical ideal $\qqq$ of $\kk[x_1, \ldots, x_n]$, one has $V(I(X)) = X$ and $I(V(\qqq)) = \qqq$.
\end{thm}

\begin{proof}
This follows immediately from identity \eqref{eqn:nullstellensatz} and \cref{exercise:radicalIX,exercise:zariski-closure-aff}.
\end{proof}

In many ways the Zariski topology is distinctly different from the Euclidean topology on $\rr$ and $\cc$, in part due to ``fewer'' closed sets. We have already seen that the only proper closed subsets of $\kk$ are finite subsets. See \cref{exercise:Zariski-quirk:Hausdorff,exercise:Zariski-quirk:compact} for some other quirks of Zariski topology.

\subsection{Exercises}
\begin{exercise} \label{exercise:V-containment}
 If $\qqq_1 \subseteq \qqq_2 \subseteq \kk[x_1, \ldots, x_n]$, show that $V(\qqq_1) \supseteq V(\qqq_2)$.
\end{exercise}

\begin{exercise} \label{exercise:V-union-intersection}
Given subsets $\qqq_1, \qqq_2$ of $\kk[x_1, \ldots, x_n]$, show that
\begin{enumerate}
\item $V(\qqq_1) \cap V(\qqq_2) = V(\qqq_1 \cup \qqq_2) = V(\qqq_1 + \qqq_2)$.
\item $V(\qqq_1) \cup V(\qqq_2) = V(\qqq_1 \cap \qqq_2) = V(\qqq_1\qqq_2)$, where $\qqq_1\qqq_2 = \{f_1f_2: f_j \in \qqq_j,\ j = 1, 2\}$. [Hint: the inclusions $V(\qqq_1) \cup V(\qqq_2) \subseteq V(\qqq_1 \cap \qqq_2) \subseteq V(\qqq_1\qqq_2)$ follows from \cref{exercise:V-containment}, so that it suffices to show $V(\qqq_1\qqq_2) \subseteq V(\qqq_1) \cup V(\qqq_2)$. For the latter containment pick $a \in V(\qqq_1\qqq_2) \setminus V(\qqq_1)$. There is $f \in \qqq_1$ such that $f(a) \neq 0$. Consider the products $fg$ with $g \in \qqq_2$ to show that $a \in V(\qqq_2)$.]
\end{enumerate}
\end{exercise}

\begin{exercise} \label{exercise:zariski-topology-aff}
Show that the collection of subvarieties on an affine space satisfies the axioms of a topology, namely that it contains the empty set and the affine space itself, and it is closed under finite unions and arbitrary intersections. [Hint: use \cref{exercise:V-union-intersection}.]
\end{exercise}

\begin{exercise} \label{exercise:separate-points-aff}
Given finitely many distinct points $a_1, \ldots, a_N \in \kk^n$, construct a polynomial $f \in \kk[x_1, \ldots, x_n]$ which is nonzero at $a_1$ but zero at $a_j$ for each $j \neq 1$. Conclude that there is a Zariski open neighborhood of $a_1$ in $\kk^n$ which does not contain any $a_j$ for $j \neq 1$. 
\end{exercise}

\begin{exercise} \label{exercise:power-V}
Given an ideal $\qqq$ of $\kk[x_1, \ldots, x_n]$ and $m \in \zzero$, let $\qqq^m$ be the ideal generated by all $f_1\cdots f_m$ for $f_1, \ldots, f_m \in \qqq$. Show that $V(\qqq^m) = V(\qqq)$ for each $m \geq 1$. Conversely, show that $V(\qqq) = V(\sqrt{\qqq})$, where $\sqrt{\qqq}$ is the {\em radical} of $\qqq$.
\end{exercise}

\begin{exercise} \label{exercise:basic-open-cover-Kn}
Given any $f \in \kk[x_1, \ldots, x_n]$, show that the set $\kk^n \setminus V(f)$ is {\em Zariski open} (i.e.\ open with respect to the Zariski topology). Deduce that every Zariski open subset of $\kk^n$ has a open covering by finitely many subsets of the form $\kk^n \setminus V(f)$. [Hint: for finiteness you need to use Hilbert's basis theorem (\cref{thm:Hilbert-basis}).]
\end{exercise}

\begin{exercise} \label{exercise:radicalIX}
Show that $I(X)$ is a radical ideal of $\kk[x_1, \ldots, x_n]$ for each $X \subseteq \kk^n$.
\end{exercise}

\begin{exercise} \label{exercise:IX-examples}
Compute $I(X)$ for $X \subseteq \kk^n$ in the following cases:
\begin{enumerate}
\item $n = 1$, $X$ consists of two distinct points in $\kk$.
\item $n = 2$, $X$ is the $x$-axis.
\item $n = 2$, $X$ is the union of $x$ and $y$-axes.
\item $n = 2$, $X$ is the union of $x$-axis and the point $(1, 0)$.
\item $n = 1$, $\kk = \cc$, $X = \zz = \{0, \pm 1, \pm 2, \ldots\}$.
\end{enumerate}
\end{exercise}

\begin{exercise} \label{exercise:I-containment}
If $X_1 \subseteq X_2 \subseteq \kk^n$, show that $I(X_1) \supseteq I(X_2)$.
\end{exercise}

\begin{exercise} \label{exercise:zariski-closure-aff}
Given $X \subseteq \kk^n$, the closure of $X$ in $\kk^n$ under the Zariski topology is the subvariety $\bar X := V(I(X))$ (we say that $\bar X$ is the \index{Zariski!closure}{\em Zariski closure} of $X$ in $\kk^n$). [Hint: it suffices to show that every subvariety $V$ of $\kk^n$ containing $X$ also contains $V(I(X))$. If $V = V(\qqq)$, then show that $\qqq \subseteq I(X)$.]
\end{exercise}


\begin{exercise}\label{exercise:pol-zero-on-Kn}
Let $k$ be a field and $x_1, \ldots, x_n$, $n \geq 1$, be indeterminates over $k$.
\begin{enumerate}
\item If $k$ is infinite, then show that for each nonzero polynomial $f \in k[x_1, \ldots, x_n]$, there is $a \in k^n$ such that $f(a) \neq 0$. [Hint: prove it for $n = 1$. In the general case, after renumbering the $x_j$ if necessary, $f$ can be expressed as a nonzero polynomial in $x_n$ with coefficients in $k[x_1, \ldots, x_{n-1}]$. Apply induction and reduce to case $n = 1$.]
\item If $k$ is finite, then show that there is a nonzero polynomial $f \in k[x_1, \ldots, x_n]$ such that $f(a) = 0$ for all $a \in k^n$.
\end{enumerate}
\end{exercise}

\begin{exercise}\label{exercise:Zariski-quirk:Hausdorff}
\begin{enumerate}
\item Show that any pair of nonempty open subsets of $\kk^n$ has a nonempty intersection. [Hint: due to \cref{exercise:basic-open-cover-Kn} it suffices to show that $V(f_1f_2) \neq \emptyset$ for non-constant polynomials $f_1, f_2$. Now use \cref{exercise:pol-zero-on-Kn}.]
\item Deduce that $\kk^n$ is {\em not} Hausdorff for any $n \geq 1$.
\end{enumerate}
\end{exercise}

\begin{exercise} \label{exercise:Zariski-quirk:compact}
Recall that a topological space is {\em compact} if each of its open covers has a finite subcover. Show that every affine variety is compact\footnote{In many algebraic geometry and commutative algebra texts this property is defined as {\em quasicompactness}, and ``compactness'' is reserved for spaces which are both quasicompact and Hausdorff.}. [Hint: due to \cref{exercise:V-union-intersection} it suffices to prove that ``if $\qqq$ is the sum of a collection $\{\qqq_i\}_{i \in \scrI}$ of ideals of $\kk[x_1, \ldots, x_n]$, then $\qqq$ is actually the sum of finitely many of the $\qqq_i$.'' Now use the fact that $\kk[x_1, \ldots, x_n]$ is Noetherian.]
\end{exercise}

\begin{exercise} \label{exercise:nonempty-hypersurface-aff}
This exercise illustrates a fundamental property of algebraically closed fields.
\begin{enumerate}
\item Show that every non-constant polynomial over $\kk$ in $(x_1, \ldots, x_n)$ vanishes at some point of $\kk^n$. [Hint: proceed by induction on $n$. Treat polynomials in $(x_1, \ldots, x_n)$ over $\kk$ as polynomials in $x_n$ over $\kk[x_1, \ldots, x_{n-1}]$, and use the inductive hypothesis to reduce to the case of $n = 1$.]
\item Show by examples that the preceding statement may be false if $\kk$ is not algebraically closed.
\end{enumerate}
\end{exercise}

\section{(Ir)reducibility} \label{irreducible-section-aff}
A topological space $X$ is called \index{Reducible topological space}{\em reducible} (respectively, \index{Irreducible!topological space}{\em irreducible}) if it can (respectively, can not) be represented as the union of two proper closed subsets. The following is a compilation of a few basic properties of irreducible sets - all these follow directly from the definition of irreducible sets; their verification is left as an exercise.

\begin{prop} \label{prop:irreducible-properties}
Let $X$ be a topological space.
\begin{enumerate}
\item \label{assertion:irreducible} $X$ is irreducible if and only if any two of its nonempty open subsets has a nonempty intersection. In particular, if $X$ is irreducible, then no proper subset of $X$ can be both open and closed in $X$. 
\item \label{assertion:irreducible-union} Let $V_1, \ldots, V_k$ be closed subsets of $X$. If $V$ is a closed irreducible subset of $X$ such that $V \subseteq \bigcup_j V_j$, then $V \subseteq V_j$ for some $j$.
\item \label{assertion:irreducible-dense} If $W$ is a dense subset of $X$, then $W$ is irreducible if and only if $X$ is irreducible.
\item \label{assertion:irreducible-open} If $X$ is irreducible, then every nonempty open subset of $X$ is irreducible and dense in $X$. 
\item \label{assertion:irreducible-continuous} If $X$ is irreducible, then the image of $X$ under a continuous map is irreducible. \qed
\end{enumerate}
\end{prop}

The relevance of irreducibility in algebraic geometry comes from the observation that under the basic correspondence between subvarieties of $\kk^n$ and radical ideals of $\kk[x_1, \ldots, x_n]$ given in \cref{thm:affine-nullspondence}, irreducible subvarieties correspond to {\em prime} ideals:

\begin{prop} \label{prop:irred=prime}
A subvariety $X$ of $\kk^n$ is irreducible if and only if $I(X)$ is a prime ideal of $\kk[x_1, \ldots, x_n]$.
\end{prop}

\begin{proof}
If $I(X)$ is not prime, then there are $f_1, f_2 \not\in I(X)$ such that $f_1f_2 \in I(X)$. Then $X \subseteq V(f_1f_2) = V(f_1) \cup V(f_2)$ (\cref{exercise:V-union-intersection}), but $X \not\subseteq V(f_j)$ for any $j$. It follows that $X$ is not irreducible (\cref{prop:irreducible-properties}). On the other hand, if $X$ is not irreducible, then there are $\qqq_1, \qqq_2 \subseteq \kk[x_1, \ldots, x_n]$ such that $X \subseteq V(\qqq_1) \cup V(\qqq_2)$, but $X \not\subseteq V(\qqq_j)$ for any $j$. Then we can pick $f_j \in \qqq_j$ which does not vanish everywhere on $X$. But $f_1f_2$ vanishes on $V(\qqq_1) \cup V(\qqq_2) \supseteq X$, so that $f_1f_2 \in I(X)$ even though neither $f_1$ nor $f_2$ is in $I(X)$. It follows that $I(X)$ is not prime, as required.
\end{proof}

\begin{example} \label{example:kn-irreducible}
A finite set of points in $\kk^n$ is irreducible if and only if it consists of only one point. Since the intersection of any pair of nonempty open subsets of $\kk^n$ is nonempty (\cref{exercise:Zariski-quirk:Hausdorff}), \cref{prop:irreducible-properties} implies that $\kk^n$ is irreducible. Note that $I(\kk^n)$, being the zero ideal (\cref{example:Ikn}), is prime in the polynomial ring. The concept of irreducibility is simple but powerful; see \cref{exercise:cayley-hamilton} for a simple proof of the {\em Cayley-Hamilton theorem} on matrices using the irreducibility of $\kk^n$.
\end{example}

An \index{Irreducible!component!of a topological space}{\em irreducible component} of $X$ is a closed irreducible subset which is not properly contained in any other irreducible subset of $X$. Each point of a finite subset $S$ of $\kk^n$ is an irreducible component of $S$. Consider e.g.\ the case that $X = \{a_1, \ldots, a_k\} \subset \kk$. Then $I(X)$ is the ideal of $\kk[x]$ generated by $(x - a_1)(x - a_2)\cdots (x - a_k)$. Note that
\begin{align*}
I(X) = \bigcap_j \langle x - a_j \rangle
\end{align*}
and the ideals $\langle x - a_j \rangle$ are precisely the ideals of the irreducible components of $X$; this is a manifestation of the following general property of affine varieties:

\begin{thm} \label{thm:irred-components-affine}
Let $X$ be a subvariety of $\kk^n$. Then there is a unique {\em minimal representation}\footnote{A representation $\qqq = \bigcap_{i=1}^k \qqq_i$ is {\em minimal} if $\qqq_j \not\supseteq \bigcup_{i \neq j} \qqq_i$ for any $j$.}
\begin{align}
I(X) = \bigcap_j \ppp_j \label{eq:I(X)-primary-decomp}
\end{align}
of $I(X)$ as the intersection of finitely many prime ideals. The operators $I(\cdot)$ and $V(\cdot)$ induce a one-to-one correspondence between irreducible components of $X$ and the prime ideals $\ppp_j$ that appear in \eqref{eq:I(X)-primary-decomp}. In particular, $X$ has {\em finitely many} irreducible components.
\end{thm}

\begin{proof}
At first we prove the following: every subvariety of $\kk^n$ is the union of finitely many closed irreducible subsets. Indeed, otherwise there is a subvariety $X_0$ of $\kk^n$ which is not the union of finitely many closed irreducible subsets. In particular $X_0$ is reducible, and it can be expressed as the union of proper closed subsets $Y_1, Y_2$. At least one of the $Y_j$ must also have the property that it is not the union of finitely many closed irreducible subsets; denote it by $X_1$. By the same arguments $X_1$ has a proper closed subset $X_2$ which is not the union of finitely many closed irreducible subsets, and continuing in this way we can construct an infinite chain $X \supsetneq X_1 \supsetneq X_2 \supsetneq \cdots$ of subvarieties of $\kk^n$. But then there is an infinite strictly ascending chain of ideals $I(X_0) \subsetneq I(X_1) \subsetneq I(X_2) \subsetneq \cdots$ of ideals of $\kk[x_1, \ldots, x_n]$, which violates the Noetherianity of polynomial rings (\cref{thm:Hilbert-basis}). This proves the claim. Now pick closed irreducible  subsets $X_1, \ldots, X_k$ of $X$ such that
\begin{align}
X = \bigcup_j X_j \label{eqn:Xunion}
\end{align}
Discarding some $X_j$ from the union if necessary, we may assume that the presentation in \eqref{eqn:Xunion} is {\em minimal}, i.e.\ $X_j \not\subseteq \bigcup_{i \neq j} X_i$ for any $j$. It is then straightforward to see (e.g.\ using assertion \eqref{assertion:irreducible-union} of \cref{prop:irreducible-properties}) that $X_j$ are precisely the irreducible components of $X$. This in particular proves the last assertion of \cref{thm:irred-components-affine}. \\

Let $\ppp_j := I(X_j)$, $j = 1, \ldots, k$. Then $\ppp_j$ are prime (\cref{prop:irred=prime}) and $X = V(\bigcap_j \ppp_j)$ (\cref{exercise:V-union-intersection}). Since $\bigcap \ppp_j$ is radical (\cref{exercise:radical-intersection=radical}), the Nullstellensatz (\cref{thm:Hilbert-nulls}) implies that $I(X) = \bigcap_j \ppp_j$. The minimality of this representation of $I(X)$ follows from the minimality of the representation in \eqref{eqn:Xunion}. It is straightforward to see that any such minimal representation is unique; it is left as an exercise (\cref{exercise:prime-containment}).
\end{proof}

\begin{example} \label{example:affine-hypersurface}
Consider the hypersurface $V(f)$ of $\kk^n$ determined by $f \in \kk[x_1, \ldots, x_n]$. Recall that $\kk[x_1, \ldots, x_n]$ is a {\em unique factorization domain}. If $f_1, \ldots, f_k$ are the distinct {\em irreducible factors} of $f$, then the irreducible components of $V(f)$ are precisely $V(f_j)$, $j = 1, \ldots, k$. Moreover, the ideals generated by the $f_j$ are prime and the ideal generated by $\prod_j f_j$ is radical, so that $I(V(f)) = \langle \prod_j f_j \rangle = \bigcap_j \langle f_j \rangle = \bigcap_j I(V(f_j))$.
\end{example}

\subsection{Exercises}
\begin{exercise}\label{exercise:irreducible-properties}
Prove \cref{prop:irreducible-properties}.
\end{exercise}

\begin{exercise} \label{exercise:radical-intersection=radical}
Show that intersections of any collection of radical ideals in a ring is also a radical ideal.
\end{exercise}

\begin{exercise}\label{exercise:prime-containment}
If $\ppp$ is a prime ideal of a ring $R$ containing the intersection of finitely many ideals $\qqq_1, \ldots, \qqq_k$, then show that $\ppp \supseteq \qqq_j$ for some $j$. Conclude that an ideal of $R$ can have (up to reordering) at most one minimal presentation as the intersection of finitely many prime ideals.
\end{exercise}

\begin{exercise}\label{exercise:cayley-hamilton}
Given an $n \times n$ matrix $A$ over a field $F$, its {\em characteristic polynomial} $p(\lambda)$ is the polynomial $\det(A - \lambda \id_n) \in F[\lambda]$, where $\det$ denotes the determinant, $\id_n$ is the $n\times n$ identity matrix and $\lambda$ is an indeterminate over $F$. The {\em Cayley-Hamilton theorem} states that $p(A)$ is the zero matrix. In this exercise following a suggestion on {\em MathOverflow} \cite{mathoverflow-algeom-examples-cayham} we outline a proof of the Cayley-Hamilton theorem using the following fact from linear algebra: if $p(\lambda)$ has $n$ distinct roots in $F$, then $A$ is {\em diagonalizable}\footnote{Recall that a matrix $A$ is {\em diagonalizable} if $A = PDP^{-1}$ for an invertible matrix $P$ and a diagonal matrix $D$.}.
\begin{enumerate}
\item Show that to prove Cayley-Hamilton theorem it suffices to assume that $F$ is algebraically closed. In all steps below assume $F$ is algebraically closed.
\item Show that Cayley-Hamilton theorem is true for diagonalizable matrices.
\item Identify the space of $n \times n$ matrices over $F$ with the affine space $F^{n^2}$. Show that to prove Cayley-Hamilton theorem it suffices to show that the set of diagonalizable matrices is Zariski dense in $F^{n^2}$. [Hint: $I(F^{n^2})$ is the zero ideal.]
\item In \cref{example:distinct-roots} we will see that the following is true:
\begin{align}
\parbox{0.8\textwidth}{
Given an algebraically closed field $\kk$ and a positive integer $n$, there is a nonempty Zariski open subset $U$ of $\kk^{n+1}$ such that for each $(c_0, \ldots, c_n) \in U$, the polynomial $c_0\lambda^n + c_1 \lambda^{n-1} + \cdots + c_n$ has $n$ distinct roots in $\kk$.
} \label{property:discriminant}
\end{align}
Use this fact to show that the set of diagonalizable matrices contains a dense Zariski open subset of $F^{n^2}$. [Hint: use the irreducibility of the affine space, \cref{exercise:basic-open-cover-Kn} and assertion \eqref{assertion:irreducible-open} of \cref{prop:irreducible-properties}.]
\item It turns out that there is a unique irreducible polynomial $\Delta$ in $\kk[c_0, \ldots, c_n]$ such that $c_0\lambda^n + c_1 \lambda^{n-1} + \cdots + c_n$ has multiple roots if and only if $\Delta(c_0, \ldots, c_n) = 0$; it is called the {\em discriminant}. For $n = 2$ use the quadratic formula to explicitly compute the discriminant. Discriminants are discussed in many introductory algebra or algebraic geometry books, e.g.\ \cite[Section II.2]{griffiths}.
\end{enumerate}
\end{exercise}

\section{Regular functions, coordinate rings and morphisms of affine varieties} \label{regular-section-aff}
Let $X$ be a subvariety of $\kk^n$ with coordinates $(x_1, \ldots, x_n)$. A \index{Regular function!on a variety}{\em regular function} on $X$ is a function $\phi: X \to \kk$ which is the restriction of a polynomial in $(x_1, \ldots, x_n)$ over $\kk$. The set of regular functions on $X$, equipped with the natural $\kk$-algebra structure, is called the \index{Coordinate ring}{\em coordinate ring} of $X$, and denoted as $\kk[X]$. The restriction map induces a natural surjective homomorphism $\kk[x_1, \ldots, x_n] \to \kk[X]$. Since $f|_X \equiv 0$ if and only if $f \in I(X)$, it follows that
\begin{align}
\kk[X] \cong \kk[x_1, \ldots, x_n]/I(X) \label{eqn:K[X]-affine}
\end{align}

\begin{example} \label{example:K[X]-singleton}
If $X$ is a singleton, then $\kk[X] \cong \kk$ (this follows from \cref{example:Isingleton} and identity \eqref{eqn:K[X]-affine}). On the other extreme, if $X = \kk^n$, then $\kk[X] \cong \kk[x_1, \ldots, x_n]$. If $H$ is the {\em hyperbola} $V(xy - 1) \subseteq \kk^2$, then $\kk[H] \cong \kk[x,y]/\langle xy - 1 \rangle \cong \kk[x, x^{-1}]$.
\end{example}

\begin{example} \label{example:K[X]-component}
All the coordinate rings computed in \cref{example:K[X]-singleton} are {\em integral domains}. In fact it turns out that $\kk[X]$ is an integral domain if and only if $X$ is irreducible (\cref{exercise:prime=integral}). If $X$ is the union of the $x$ and $y$-axes in $\kk^2$, then $I(X) = \langle xy \rangle$ (\cref{example:affine-hypersurface}), so that $\kk[X] \cong \kk[x,y]/\langle xy \rangle$, which in particular is not an integral domain. If $X$ is the union of the hyperbola $H := V(xy - 1) \subseteq \kk^2$ from \cref{example:K[X]-singleton} and the $x$-axis (which we denote by $A$) on $\kk^2$, then
\begin{align}
\kk[X]
	&\cong \kk[x,y]/\langle y(xy - 1) \rangle
	\cong \kk[x,y] /\langle y \rangle \times \kk[x,y] / \langle xy - 1 \rangle
	\cong \kk[A] \times \kk[H]
\label{eqn:K[X]-hyperbolaxunion}
\end{align}
(\cref{exercise:K[X]-hyperbolaxunion}). More generallly, you will prove in \cref{exercise:K[X]-emptyinterunion} that the coordinate ring of a pairwise disjoint union of affine varieties is isomorphic to the product of their coordinate rings.
\end{example}

A map $\phi: X \to Y$ between affine varieties is called a \index{Morphism}{\em morphism} if for every regular function $h$ on $Y$, the pullback $h \circ \phi$ is a regular function on $X$. An \index{Isomorphism of varieties}{\em isomorphism} is a bijective morphism whose inverse is also a morphism; the notation $X \cong Y$ is a shorthand for the statement that $X$ and $Y$ are {\em isomorphic}, i.e.\ there is an isomorphism between $X$ and $Y$. It is straightforward to check that a morphism $\phi: X \to Y$ induces a $\kk$-algebra homomorphism $\phi^*:\kk[Y] \to \kk[X]$ given by $h \mapsto h \circ \phi$, and that $\phi$ is an isomorphism if and only if $\phi^*$ is a $\kk$-algebra isomorphism between $\kk[X]$ and $\kk[Y]$ (\cref{exercise:affine-pullback}).

\begin{example}
Projections constitute a basic source of morphisms. Let $X$ be the parabola $y = x^2$ in $\kk^2$. Note that $\kk[X] = \kk[x,y]/\langle y - x^2 \rangle \cong \kk[x]$, so that $X \cong \kk$. Indeed, the projection onto $x$-axis realizes this isomorphism, but the projection onto $y$-axis induces a two-to-one morphism $X \to \kk$ (\cref{exercise:parabola-to-k}). Note that the map $(x,y) \mapsto (x, y - x^2)$ is an {\em automorphism}\footnote{An {\em automorphism} of $\kk^n$ is an isomorphism from $\kk^n$ to itself.} of $\kk^2$ which maps $X$ onto the $x$-axis. This is a special case of the celebrated {\em line embedding theorem} of S.\ S.\ Abhyankar and T.\ T.\ Moh \cite{abhya-moh-line} which states that ``if $X$ is a subvariety of $\kk^2$ isomorphic to $\kk$ then there is an automorphism of $\kk^2$ which maps $X$ onto the $x$-axis.'' The following analogous statement remains a conjecture (named after Abhyankar and A.\ Sathaye, who was a student of Abhyankar) for $n \geq 2$: ``if $X$ is a subvariety of $\kk^{n+1}$ isomorphic to $\kk^n$ then there is an automorphism of $\kk^{n+1}$ which maps $X$ onto the the hyperplane $x_n = 0$.''
\end{example}

\begin{example} \label{example:k-to-cusp-cubic}
The image of the morphism $\kk \to \kk^2$ given by $t \mapsto (t^2, t^3)$ is the variety $X = V(x^3 - y^2)$. Considered as a morphism from $\kk$ to $X$, this induces a bijection between the points of $\kk$ and $X$, but it is {\em not} an isomorphism (\cref{exercise:k-to-cusp-cubic}).
\end{example}

\begin{example} \label{example:Kn-minus-V(f)}
Let $f$ be a nonzero polynomial in $(x_1, \ldots, x_n)$ and $X$ be the hypersurface of $\kk^{n+1}$ defined by $x_{n+1}f = 1$. Then $\kk[X] \cong \kk[x_1, \ldots, x_n, 1/f]$. Let $\pi: X \to \kk^n$ be the projection onto the first $n$-coordinates. The image of $\pi$ is the proper open subset $U_f := \kk^n \setminus V(f)$ of $\kk^n$. It is clear that $\pi$ is one-to-one; in fact $\pi$ induces a {\em homeomorphism} with respect to Zariski topology (\cref{exercise:Kn-minus-V(f)-homeo}). We will see in \cref{example:affine-complement} that $\pi$ is actually an {\em isomorphism of quasiprojective varieties.}
\end{example}

\begin{example} \label{example:constructible-0}
The morphism $\sigma: \kk^2 \to \kk^2$ given by $(x,y) \mapsto (x, xy)$ maps the whole $y$-axis to the origin, but it is one-to-one at every point not on the $y$-axis. The image of $\sigma$, which is the union of the origin and all points {\em not} on the $y$-axis, is neither open nor closed in $\kk^2$ (\cref{exercise:sigma-homeo}). However, it is the union of a closed subset and an open subset, which is a ``constructible set''; in \cref{constructible-section} we discuss constructible sets and prove Chevalley's theorem that the image of every morphism is a constructible set.
\end{example}

\subsection{Exercises}

\begin{exercise} \label{exercise:prime=integral}
Show that an affine variety is irreducible if and only if its coordinate ring is an integral domain. [Hint: use \cref{prop:irred=prime} and identity \eqref{eqn:K[X]-affine}.]
\end{exercise}

\begin{exercise} \label{exercise:K[X]-hyperbolaxunion}
Consider the ``diagonal'' map $\delta: \kk[x,y] \mapsto \kk[x,y] /\langle y \rangle \times \kk[x,y] / \langle xy - 1 \rangle$ which maps $f \mapsto (f + \langle y \rangle, f + \langle xy - 1 \rangle)$.
\begin{enumerate}
\item Show that $\delta$ is a surjective $\kk$-algebra homomorphism and $\ker(\delta) = \langle y(xy-1) \rangle$.
\item Verify the $\kk$-algebra isomorphisms presented in \eqref{eqn:K[X]-hyperbolaxunion}. [Hint: use \cref{example:affine-hypersurface} and the preceding assertion.]
\item Find a polynomial $f \in \kk[x,y]$ such that the $f|_A = x|_A$ and $f|_H = y|_H$ (where $A$ is the $x$-axis and $H$ is the hyperbola $V(xy-1)$ in $\kk^2$).
\end{enumerate}
\end{exercise}

\begin{exercise} \label{exercise:K[X]-emptyinterunion}
Let $X := \bigcup_{j=1}^k X_j$, where $X_j$ are subvarieties of $\kk^n$ with pairwise empty intersection. Write $\qqq_j := I(X_j)$ and let $\delta: \kk[x_1,\ldots, x_n] \to \prod_j \kk[x_1, \ldots, x_n]/\qqq_j$ be the ``diagonal'' map that sends $f \mapsto \prod_j(f + \qqq_j)$. Show that
\begin{enumerate}
\item $I(X) = \cap_j \qqq_j$.
\item $\ker(\delta) = \cap_j \qqq_j$.
\item $\qqq_i$ and $\qqq_j$ are \index{Coprime ideal}{\em coprime}\footnote{Two ideals $\ppp,\qqq$ of a ring $R$ are {\em coprime} if $\ppp + \qqq = R$.} if $i \neq j$. [Hint: you need to use that $X_i \cap X_j = \emptyset$.]
\item Deduce that $\delta$ is surjective. [Hint: it suffices to show there is $f \in \kk[x_1, \ldots, x_n]$ such that $\delta(f) = (1, 0, \ldots, 0)$. Due to the preceding assertion for each $i > 1$, there is $f_i \in \qqq_i$ and $g_i \in \qqq_1$ such that $f_i + g_i = 1$. Take $f := \prod_{i>1} f_i$.]
\item Deduce that $\kk[X] \cong \prod_j \kk[X_j]$.
\end{enumerate}
\end{exercise}

\begin{exercise}\label{exercise:affine-pullback}
Show that a morphism $\phi: X \to Y$ between affine varieties induces a $\kk$-algebra homomorphism $\phi^*:\kk[Y] \to \kk[X]$ via pullback. Show that $\phi$ is an isomorphism if and only if $\phi^*$ is a $\kk$-algebra isomorphism.
\end{exercise}

\begin{exercise}\label{exercise:affine-isomeo}
Show that an isomorphism of affine varieties induces a homeomorphism between them.
\end{exercise}

\begin{exercise}\label{exercise:parabola-to-k}
Let $X$ be the parabola $y = x^2$ in $\kk^2$. Let $\pi_x$ (respectively, $\pi_x$) be the morphism from $X$ to $\kk$ induced by the projection onto $x$-axis (respectively, $y$-axis). Compute the pullback maps $\kk[t] \to \kk[X]$ (where $t$ is the coordinate on $\kk$) induced by $\pi_x$ and $\pi_y$. Deduce that $\pi_x$ is an isomorphism, but $\pi_y$ is not. Compute $\pi_x^{-1}:\kk \to X$.
\end{exercise}

\begin{exercise} \label{exercise:k-to-cusp-cubic}
Prove the statements from \cref{example:k-to-cusp-cubic}. [Hint: compute $\kk[X]$.]
\end{exercise}

\begin{exercise} \label{exercise:algebra-to-affine-variety}
Let $R$ be a finitely generated $\kk$-algebra which is \index{Reduced!algebra}{\em reduced}\footnote{A ring is {\em reduced} if it does not have any nonzero nilpotent elements.}. Show that there is an affine variety $X$ such that $\kk[X] \cong R$. [Hint: Construct a surjective $\kk$-algebra homomorphism from a polynomial ring over $\kk$ to $R$.] \Cref{exercise:affine-pullback} shows that $X$ is unique up to an isomorphism. 
\end{exercise}

\begin{exercise} \label{exercise:affine-pullback-stronger}
The correspondence between $\kk$-algebra homomorphisms of finitely generated $\kk$-algebras and morphisms of affine varieties can be pushed a little bit further than \cref{exercise:affine-pullback}. Given affine varieties $X,Y$ and any $\kk$-algebra homomorphism $\Phi: \kk[Y] \to \kk[X]$, show that there is a morphism $\phi: X \to Y$ such that $\Phi = \phi^*$.  This, together with \cref{exercise:affine-pullback,exercise:algebra-to-affine-variety}, shows that the categories of affine varieties and finitely generated reduced $\kk$-algebras are equivalent.
\end{exercise}

\begin{exercise} \label{exercise:morphism-zariski-continuous}
\begin{enumerate}
\item Show that every open subset of an affine variety $X$ has a finite covering by open subsets of the form $X \setminus V(f)$ for regular functions $f$ on $X$. [Hint: use \cref{exercise:basic-open-cover-Kn}.]
\item Show that a morphism $\phi:Y \to X$ of affine varieties is a continuous map with respect to the Zariski topology.
\end{enumerate}
\end{exercise}

\begin{exercise} \label{exercise:Kn-minus-V(f)-homeo}
Let $f \in \kk[x_1, \ldots, x_n]$. Show that the projection onto the first $n$-coordinates of $\kk^{n+1}$ induces a homeomorphism between the subvariety $V(x_{n+1}f - 1)$ of $\kk^{n+1}$ and the open subset $\kk^n \setminus V(f)$ of $\kk^n$.
\end{exercise}

\begin{exercise} \label{exercise:sigma-homeo}
Consider the morphism $\sigma:\kk^2 \to \kk^2$ given by $(x,y) \mapsto (x, xy)$. Show that
\begin{enumerate}
\item $\sigma$ induces a homeomorphism from $\kk^2\setminus V(x)$ to itself.
\item Image of $\sigma$ is $\{(0,0)\} \cup (\kk^2\setminus V(x))$. 
\item Image of $\sigma$ is not closed in $\kk^2$. [Hint: $\kk^2$ is irreducible. Use \cref{prop:irreducible-properties}.] 
\item Image of $\sigma$ is not open in $\kk^2$. [Hint: the $y$-axis is isomorphic to $\kk$, so any proper open set of the $y$-axis is infinite.]
\end{enumerate}
\end{exercise}

\section{Quasiprojective varieties} \label{quasi-section}
Let $(x_0, \ldots, x_n)$ be a system of coordinates on $\kk^{n+1}$, $n \geq 0$. Consider the relation $\sim$ on $\kk^{n+1} \setminus \{0\}$, $n \geq 0$, defined as follows: if $a, b \in \kk^{n+1} \setminus \{0\}$ with coordinates (with respect to $(x_0, \ldots, x_n)$) respectively $(a_0, \ldots, a_n)$ and $(b_0, \ldots, b_n)$, then $a \sim b$ if and only if there is $\lambda \in \kk\setminus \{0\}$ such that $a_j = \lambda b_j$ for each $j = 0, \ldots, n$. Note that $a \sim b$ if and only if $a$ and $b$ belong to the same line through the origin on $\kk^{n+1}$. This immediately implies that $\sim$ is an equivalence relation; the \index{Projective!space}{\em $n$-dimensional projective space} $\pp^n(\kk)$, or simply $\pp^n$, is the set of the equivalence classes of $\sim$; in other words $\pp^n$ is the set of lines through the origin on $\kk^{n+1}$. We denote the equivalence class of $\sim$ containing $a$ by $[a_0: \cdots : a_n]$ and say that $[a_0: \cdots: a_n]$ is the \index{Homogeneous!coordinate}{\em homogeneous coordinate} of the point of $\pp^n$ determined by $a$. Let $f$ be a \index{Homogeneous!polynomial}{\em homogeneous}\footnote{Recall that a polynomial is {\em homogeneous} if each of its monomials has the same degree.} polynomial of degree $d$ in $(x_0, \ldots, x_n)$, then it is straightforward to check that
\begin{align}
f(\lambda a_0, \ldots, \lambda a_n) = \lambda^d f(a_0, \ldots, a_n) \label{eqn:homogeneous-vanishing}
\end{align}
for each $\lambda, a_0, \ldots, a_n \in \kk$. It follows that given $a = (a_0, \cdots, a_n) \in \kk^{n+1}$, $f(a_0, \ldots, a_n) = 0$ if and only if $f(b_0, \ldots, b_n) = 0$ for all $(b_0, \ldots, b_n)$ in the equivalence class of $a$. Consequently the set $V(f) := \{[a_0: \cdots: a_n]: f(a_0, \ldots, a_n) = 0\}$ of zeroes of $f$ is a well defined subset of $\pp^n$; we say that $V(f)$ is the \index{Hypersurface}{\em hypersurface} of $\pp^n$ determined by $f$. There is a unique topology on $\pp^n$ in which the basic closed subsets are intersections of hypersurfaces of $\pp^n$ (\cref{exercise:zariski-topology}); it is called the \index{Zariski!topology}{\em Zariski topology}. A \index{Variety!projective}\index{Projective!variety}{\em projective variety} is a Zariski closed subset of the projective space (equipped with the topology induced from the Zariski topology of $\pp^n$). A \index{Quasiprojective!variety}\index{Variety!quasiprojective}{\em quasiprojective variety} is a Zariski open subset of a projective variety (also equipped with the Zariski topology induced from $\pp^n$). Given a quasiprojective variety $X$, a \index{Quasiprojective!subset}{\em quasiprojective subset} of $X$ is a subset which is also a quasiprojective variety, and a \index{Subvariety}{\em subvariety} of $X$ is a Zariski closed subset of $X$.

\begin{example} \label{example:basic-clopen}
For each $j = 0, \ldots, n$, $V(x_j) \subseteq \pp^n$ is the set of all $[a_0: \cdots: a_n]$ such that $a_j = 0$. It is straightforward to check that the projection onto the coordinates excluding the $j$-th coordinate induces a homeomorphism between $V(x_j)$ and $\pp^{n-1}$ (\cref{exercise:Pnminus1-homeo}). The complementary Zariski open subset $U_j := \pp^n \setminus V(x_j)$ has a one-to-one correspondence with $\kk^n$ via the map
\begin{align}
\kk^n \ni (x_1, \ldots, x_n) \mapsto [x_1: \cdots:x_j: 1: x_{j+1}: \cdots : x_n] \in U_j \label{map:k^n-to-U_j}
\end{align}
We leave it as an exercise (\cref{exercise:U_j-to-k^n}) to check that the inverse to the above map is given by
\begin{align}
U_j \ni [x_0: \cdots: x_n] \mapsto
	(\frac{x_0}{x_j}, \ldots, \frac{x_{j-1}}{x_j}, \frac{x_{j+1}}{x_j}, \ldots, \frac{x_n}{x_j}) \in \kk^n \label{map:U_j-to-k^n}
\end{align}
It is clear that $\pp^n = \bigcup_j U_j$. We say that the $U_j$ are \index{Basic open subsets of $\pp^n$}{\em basic open subsets} of $\pp^n$.
\end{example}

\begin{example} \label{example:closed-in-P1}
If $f$ is a homogeneous polynomial of degree $d$ in $(x_0, x_1)$ then $f$ can be expressed as $\prod_{i=1}^d (a_i x_0 - b_ix_1)$ for $(a_i, b_i) \in \kk^2 \setminus \{0\}$ (\cref{exercise:dehomogenize-x0x1}), so that $V(f) = \{[b_i:a_i]: i = 1, \ldots, d\} \subseteq \pp^1$. It follows that all proper subvarieties of $\pp^1$ are finite sets. This is a special case of the general fact that all proper subvarieties of a ``curve'' are finite sets (\cref{example:curve-in-curve}).
\end{example}

\begin{example} \label{example:not-clopen}
The subet $(\pp^n \setminus V(x_0)) \cup V(x_0, x_1)$ of $\pp^n$ is neither closed nor open for $n \geq 2$ (\cref{exercise:not-clopen}). This is the projective analogue of the neither-closed-nor-open subset of $\kk^2$ from \cref{example:constructible-0}. 
\end{example}

Let $X = V(f_i: i \in \scrI)$ be the subvariety of $\pp^n$ defined by a collection $\{f_i\}_{i \in \scrI}$ of homogeneous polynomials in $(x_0, \ldots, x_n)$. The \index{Cone!over a projective variety}{\em cone $C(X)$ over $X$} is the affine subvariety of $\kk^{n+1}$ determined by $\{f_i\}_{i \in \scrI}$ \Cref{prop:XCX} below shows that $C(X)$ is a union of lines through the origin (i.e.\ $C(X)$ is an actual ``cone''), and these lines are in one-to-one correspondence with points on $X$.

\begin{example}
The real points of the cone over $X := V(x_0^2 - x_1^2 - x_2^2) \subseteq \pp^2(\cc)$ is the ``circular double cone'' pictured in \cref{fig:cone-over-circle}. Different cross sections of $C(X)$ yields different representations of points on $X$, e.g. the real points of $X$ are in one-to-one correspondence with
\begin{itemize}
\item a circle (intersection with $x_0 = 1$),
\item a hyperbola (intersection with $x_1 = 1$) and two ``points at infinity'' (namely $[1:0:1], [1:0:-1]$),
\item a parabola (intersection with $x_0 + x_1 = 1$) and one ``point at infinity'' (namely $[1:-1:0]$).
\end{itemize}
This is a manifestation of the equivalence of ``compactifications'' of conic sections under ``projective transformations.''
\end{example}

\begin{figure}[htb]
\newcommand\planecolor{green}
\newcommand\conecolor{blue}
\newcommand\intersectcolor{magenta}
\tikzstyle{sdot} = [yellow, circle, minimum size = 3pt, inner sep = 0pt, fill]
\pgfmathsetmacro\scalefactor{0.6}
\pgfmathsetmacro\axiscalefactor{1}
\pgfmathsetmacro\samplenum{5}
\pgfmathsetmacro\sampley{45}
\pgfmathsetmacro\opacity{0.3}
\pgfmathsetmacro\fillopacity{0.3}

\begin{subfigure}[b]{0.3\textwidth}
\begin{tikzpicture}[scale = \scalefactor]
\pgfmathsetmacro\lengthzero{1.35}
\pgfmathsetmacro\b{0.75}
\pgfmathsetmacro\c{1}

\pgfmathsetmacro\conerpluslimit{1}
\pgfmathsetmacro\conerminuslimit{1}

\begin{scope}
\begin{axis}[
    axis lines = middle,
    axis equal,
    view/h = 45, view/v = 15,
    scale = \axiscalefactor,
    xticklabels = {,,},
    yticklabels = {,,},
    zticklabels = {,,},
    x label style={anchor = north},
    y label style={anchor = north},
    z label style={anchor = south},
    xlabel = {\picfontsize $x_1$},
    ylabel = {\picfontsize $x_2$},
    zlabel = {\picfontsize $x_0$}
]
\addplot3 [
    surf,
    color = \conecolor,
    opacity = \opacity,
    fill opacity = \fillopacity,
    faceted color = \conecolor,
    samples = \samplenum,
    samples y = \sampley,
    domain = 0 : \conerpluslimit,
    y domain = 0 : 360
] (
	{x*cos(y)*\c},
    {x*sin(y)*\c},
   	{x}
);

\addplot3 [
    surf,
    color = \conecolor,
    opacity = \opacity,
    fill opacity = \fillopacity,
    faceted color = \conecolor,
    samples = \samplenum,
    samples y = \sampley,
    domain = 0 : \conerminuslimit,
    y domain = 0 : 360
] (
	{x*cos(y)*\c},
    {x*sin(y)*\c},
   	{-x}
);

\addplot3 [fill=\planecolor, opacity=\opacity] coordinates{(-\lengthzero,-\lengthzero,\b) (\lengthzero,-\lengthzero,\b) (\lengthzero,\lengthzero,\b) (-\lengthzero,\lengthzero,\b)};

\addplot3 [
    samples = \sampley,
    samples y = 0,
    domain = 0 : 360,
    \intersectcolor, thick
] (
	{\b*\c*cos(x)},
	{\b*\c*sin(x)},
    {\b}
);


\end{axis}
\end{scope}

\end{tikzpicture}

\end{subfigure}
\begin{subfigure}[b]{0.3\textwidth}
\begin{tikzpicture}[scale = \scalefactor]
\pgfmathsetmacro\lengthzero{1.35}
\pgfmathsetmacro\b{0.3}
\pgfmathsetmacro\c{1}

\pgfmathsetmacro\conerpluslimit{1}
\pgfmathsetmacro\conerminuslimit{1}

\begin{scope}
\begin{axis}[
    axis lines = middle,
    axis equal,
    view/h = 45, view/v = 15,
    scale = \axiscalefactor,
    xticklabels = {,,},
    yticklabels = {,,},
    zticklabels = {,,},
    x label style={anchor = north},
    y label style={anchor = north},
    z label style={anchor = south},
    xlabel = {\picfontsize $x_1$},
    ylabel = {\picfontsize $x_2$},
    zlabel = {\picfontsize $x_0$}
]
\addplot3 [
    surf,
    color = \conecolor,
    opacity = \opacity,
    fill opacity = \fillopacity,
    faceted color = \conecolor,
    samples = \samplenum,
    samples y = \sampley,
    domain = 0 : \conerpluslimit,
    y domain = 0 : 360
] (
	{x*cos(y)*\c},
    {x*sin(y)*\c},
   	{x}
);

\addplot3 [
    surf,
    color = \conecolor,
    opacity = \opacity,
    fill opacity = \fillopacity,
    faceted color = \conecolor,
    samples = \samplenum,
    samples y = \sampley,
    domain = 0 : \conerminuslimit,
    y domain = 0 : 360
] (
	{x*cos(y)*\c},
    {x*sin(y)*\c},
   	{-x}
);

\addplot3 [fill=\planecolor, opacity=\opacity] coordinates{(\b, -\lengthzero,-\conerminuslimit) (\b,\lengthzero,-\conerminuslimit) (\b,\lengthzero,\conerpluslimit) (\b,-\lengthzero,\conerpluslimit)};

\addplot3 [
    samples = \sampley,
    samples y = 0,
    domain = {\b/\c} : \conerpluslimit,
    \intersectcolor, thick
] (
	{\b},
	{sqrt((\c*x)^2 - \b^2)},
    {x}
);
\addplot3 [
    samples = \sampley,
    samples y = 0,
    domain = {\b/\c} : \conerpluslimit,
    \intersectcolor, thick
] (
	{\b},
	{-sqrt((\c*x)^2 - \b^2)},
    {x}
);
\addplot3 [
    samples = \sampley,
    samples y = 0,
    domain = {\b/\c} : \conerminuslimit,
    \intersectcolor, thick
] (
	{\b},
	{sqrt((\c*x)^2 - \b^2)},
    {-x}
);
\addplot3 [
    samples = \sampley,
    samples y = 0,
    domain = {\b/\c} : \conerminuslimit,
    \intersectcolor, thick
] (
	{\b},
	{-sqrt((\c*x)^2 - \b^2)},
    {-x}
);


\end{axis}
\end{scope}

\end{tikzpicture}


\end{subfigure}\hspace{-0.05\textwidth}
\begin{subfigure}[b]{0.3\textwidth}
\begin{tikzpicture}[scale = \scalefactor]
\pgfmathsetmacro\lengthzero{1.35}
\pgfmathsetmacro\b{0.6}
\pgfmathsetmacro\c{1}

\pgfmathsetmacro\conerpluslimit{1}
\pgfmathsetmacro\conerminuslimit{1}
\pgfmathsetmacro\offset{0.3}
\pgfmathsetmacro\conerminuslimitoff{\conerminuslimit - \offset}

\pgfmathsetmacro\xtwolimit{sqrt(2*\b*\c*\conerpluslimit -\b^2)}

\begin{scope}
\begin{axis}[
    axis lines = middle,
    axis equal,
    view/h = 45, view/v = 15,
    scale = \axiscalefactor,
    xticklabels = {,,},
    yticklabels = {,,},
    zticklabels = {,,},
    x label style={anchor = north},
    y label style={anchor = north},
    z label style={anchor = south},
    xlabel = {\picfontsize $x_1$},
    ylabel = {\picfontsize $x_2$},
    zlabel = {\picfontsize $x_0$}
]
\addplot3 [
    surf,
    color = \conecolor,
    opacity = \opacity,
    fill opacity = \fillopacity,
    faceted color = \conecolor,
    samples = \samplenum,
    samples y = \sampley,
    domain = 0 : \conerpluslimit,
    y domain = 0 : 360
] (
	{x*cos(y)*\c},
    {x*sin(y)*\c},
   	{x}
);

\addplot3 [
    surf,
    color = \conecolor,
    opacity = \opacity,
    fill opacity = \fillopacity,
    faceted color = \conecolor,
    samples = \samplenum,
    samples y = \sampley,
    domain = 0 : \conerminuslimit,
    y domain = 0 : 360
] (
	{x*cos(y)*\c},
    {x*sin(y)*\c},
   	{-x}
);

\addplot3 [fill=\planecolor, opacity=\opacity] coordinates{({\c*\conerminuslimitoff +\b}, -\lengthzero,-\conerminuslimitoff) ({\c*\conerminuslimitoff + \b},\lengthzero,-\conerminuslimitoff) ({-\c*\conerpluslimit + \b},\lengthzero,\conerpluslimit) ({-\c*\conerpluslimit + \b},-\lengthzero,\conerpluslimit)};

\addplot3 [
    samples = \sampley,
    samples y = 0,
    domain = 0 : \xtwolimit,
    \intersectcolor, thick
] (
	{\b/2 - x^2/(2*\b)},
	{x},
    {(x^2 + \b^2)/(2*\b*\c)}
);
\addplot3 [
    samples = \sampley,
    samples y = 0,
    domain = 0 : \xtwolimit,
    \intersectcolor, thick
] (
	{\b/2 - x^2/(2*\b)},
	{-x},
    {(x^2 + \b^2)/(2*\b*\c)}
);


\end{axis}
\end{scope}

\end{tikzpicture}


\end{subfigure}
\caption{Real points and cross sections of $C(X)$ for $X := V(x_0^2 - x_1^2 - x_2^2) \subset \pp^2$}  \label{fig:cone-over-circle}

\end{figure}

\begin{prop} \label{prop:XCX}
For every $a = [a_0: \cdots :a_n] \in X$, every point of the line in $\kk^{n+1}$ determined by $a$ is in $C(X)$. Conversely, for every $(a_0, \ldots, a_n) \in C(X)\setminus\{0\}$, the corresponding point $[a_0: \cdots : a_n]$ in $\pp^n$ is in $X$.
\end{prop}

\begin{proof}
This is a straightforward consequence of the homogeneity of the $f_i$ and identity \eqref{eqn:homogeneous-vanishing} above.
\end{proof}

Recall that the affine variety $C(X)$ determines an ideal $I(C(X))$ of $\kk[x_0, \ldots, x_n]$ consisting of all polynomials that vanish on $C(X)$. The same construction can be used to define an ideal associated to $X$, provided we consider only {\em homogeneous} generators. A \index{Homogeneous!ideal}{\em homogeneous ideal} of $\kk[x_0, \ldots, x_n]$ is an ideal generated by homogeneous polynomials. The \index{Homogeneous!ideal!of a projective variety}{\em homogeneous ideal of $X$}, denoted by $I(X)$, is the ideal of $\kk[x_0, \ldots, x_n]$ generated by {\em homogeneous} polynomials that vanish on all points of $X$.

\begin{prop} \label{prop:IXICX}
$I(X) = I(C(X))$. In particular, there are {\em finitely many} homogeneous polynomials $f_1, \ldots, f_k$ in $(x_0, \ldots, x_n)$ such that $X = V(f_1, \ldots, f_k) \subseteq \pp^n$.
\end{prop}

\begin{proof}
Recall that a \index{Homogeneous!component!of a polynomial}{\em homogeneous component} of $f = \sum_\alpha c_\alpha x^{\alpha} \in \kk[x_0, \ldots, x_n]$ (where $x^\alpha$ is a shorthand for $x_0^{\alpha_0} \cdots x_n^{\alpha_n}$) is a polynomial of the form $\sum_{|\alpha| = d} c_\alpha x^\alpha$ for some $d \geq 0$, where $|\alpha| := \alpha_0 + \cdots + \alpha_n$. Since $I(X)$ is generated by homogeneous polynomials, it has the following property (\cref{exercise:homog-ideal-homog-gen}):
\begin{align}
\parbox{0.5\textwidth}{
a polynomial $f$ in $(x_0, \ldots, x_n)$ is in $I(X)$ if and only if all homogeneous components of $f$ are in $I(X)$.
} \label{property:IXhomogeneous}
\end{align}
Property \eqref{property:IXhomogeneous} coupled with Hilbert's Basis theorem (\cref{thm:Hilbert-basis}) implies that $I(X)$ is generated by finitely many homogeneous polynomials $f_1, \ldots, f_k$; it is straightforward to check that $X = V(f_1, \ldots, f_k)$. It remains to prove that $I(X) = I(C(X))$. The containment $I(X) \subseteq I(C(X))$ follows immediately from \cref{prop:XCX} and identity \eqref{eqn:homogeneous-vanishing}. For the opposite containment pick $a = [a_0: \ldots : a_n] \in \pp^n$. Let $L(a) := \{(ta_0, \ldots, ta_n): t \in \kk\}$ be the ``line in $\kk^{n+1}$ represented by $a$.''

\begin{proclaim} \label{line-claim}
A polynomial vanishes on $L(a)$ if and only if each of its homogeneous components vanishes at $a$.
\end{proclaim}

\begin{proof}
Pick $f  \in \kk[x_0, \ldots, x_n]$ such that $f|_{L(a)} \equiv 0$. Write $f= \sum_{d = 0}^e f_d$, where each $f_d$ is homogeneous of degree $d$. Then $f(ta_0, \ldots, ta_n) = \sum_{d=0}^e t^df_d(a)$ is a polynomial in $t$ with infinitely many zeroes, so it must be identically zero. Therefore $f_d(a) = 0$ for each $d$, as required.
\end{proof}

\Cref{line-claim,prop:XCX} imply that $I(C(X)) \subseteq I(X)$, which completes the proof.
\end{proof}

There is a homogeneous ideal of $\kk[x_0, \ldots, x_n]$ which is {\em not} the homogeneous ideal of any nonempty subvariety of $\pp^n$: it is the ideal $\mmm_+$ generated by $x_0, \ldots, x_n$ [why?]. It is sometimes called the \index{Irrelevant ideal}{\em irrelevant} ideal. Note that {\em all} proper homogeneous ideals of $\kk[x_0, \ldots, x_n]$ are contained in $\mmm_+$. The correspondence between affine varieties and radical ideals described in \cref{thm:affine-nullspondence} has a projective counterpart, also due to the Nullstellensatz, provided one excludes $\mmm_+$ - this is described in \cref{exercise:projective-nullspondence}.

\subsection{Exercises}
\begin{exercise} \label{exercise:zariski-topology}
Show that the collection of Zariski closed sets of $\pp^n$ satisfies the axioms of a topology. [Hint: mimic the solution to \cref{exercise:zariski-topology-aff}.]
\end{exercise}

\begin{exercise} \label{exercise:basic-open-cover-Pn}
Show that every Zariski open subset of $\pp^n$ has an open cover by subsets of the form $\pp^n \setminus V(f)$ for homogeneous polynomials $f \in \kk[x_0, \ldots, x_n]$.
\end{exercise}

\begin{exercise} \label{exercise:Pnminus1-homeo}
Show that the map from the subvariety $V(x_0)$ of $\pp^n$ to $\pp^{n-1}$ defined by $[0:x_1: \cdots : x_n] \mapsto [x_1: \cdots : x_n])$ induces a homeomorphism with respect to the Zariski topology.
\end{exercise}

\begin{exercise} \label{exercise:U_j-to-k^n}
Show that the map from $U_j$ to $\kk^n$ from \eqref{map:U_j-to-k^n} is well defined and it is the inverse to the map from $\kk^n$ to $U_j$ given in \eqref{map:k^n-to-U_j}.
\end{exercise}

\begin{exercise} \label{exercise:dehomogenize-x0x1}
Let $f$ be a homogeneous polynomial in $(x_0, x_1)$ of degree $d$.
\begin{enumerate}
\item Show that $f/x_0^d$ is a polynomial in $x_1/x_0$ over $\kk$.
\item Since $\kk$ is algebraically closed, deduce that $f$ can be expressed as a product of {\em linear} homogeneous polynomials.
\end{enumerate}
\end{exercise}

\begin{exercise} \label{exercise:homog-ideal-homog-gen}
Let $I$ be a homogeneous ideal of $\kk[x_0, \ldots, x_n]$. 
\begin{enumerate}
\item Show that for every $f \in I$, all homogeneous components of $f$ are also in $I$. [Hint: Given $f \in I$, express it as a sum of products of polynomials and homogeneous generators of $I$. Equate homogeneous components of both sides of the equation.]
\item Deduce that $I$ is prime if and only if for all {\em homogeneous} polynomials $f_1, f_2 \in \kk[x_0, \ldots, x_n]$, $f_1f_2 \in I$ if and only if either $f_1$ or $f_2$ is in $I$.
\item Show that the radical $\sqrt I$ of $I$ is also homogeneous. [Hint: assume by contradiction that $\sqrt I$ is not homogeneous. Then there is $f \in \sqrt I$ such that {\em no} homogeneous component of $f$ is in $\sqrt I$. Pick $k \geq 1$ such that $f^k \in I$. If $f_d$ is the homogeneous component of $f$ of degree $d = \deg(f)$, then $f_d^k \in I$, so that $f_d \in \sqrt{I}$, which is a contradiction.]
\end{enumerate}
\end{exercise}

\begin{exercise} \label{exercise:projective-nullspondence}
Let $\mmm_+$ be the maximal ideal of $\kk[x_0, \ldots, x_n]$ generated by $x_0, \ldots, x_n$. Use the Nullstellensatz (\cref{thm:Hilbert-nulls}) to prove that there is a one-to-one correspondence between nonempty subvarieties of $\pp^n$ and radical homogeneous ideals properly contained in $\mmm_+$. [Hint: use \cref{exercise:homog-ideal-homog-gen}.]
\end{exercise}

\begin{exercise} \label{exercise:irreducible-projective}
Let $X$ be a subvariety of $\pp^n$. Show that $X$ is irreducible if and only if $I(X)$ is a prime ideal of $\kk[x_0, \ldots, x_n]$. [Hint: mimic the proof of \cref{prop:irred=prime}. Use \cref{exercise:homog-ideal-homog-gen}.]
\end{exercise}

\begin{exercise} \label{exercise:Pn-irreducible}
Show that $I(\pp^n) = 0$. Deduce that $\pp^n$ is irreducible.
\end{exercise}

\begin{exercise} \label{exercise:not-clopen}
Show that $(\pp^n \setminus V(x_0)) \cup V(x_0, x_1)$ is  neither open nor closed in $\pp^n$ if $n \geq 2$. [Hint: $\pp^n$ and $V(x_0)$ are irreducible. Use \cref{prop:irreducible-properties}.]
\end{exercise}

\begin{exercise} \label{exercise:basic-open-cover-qproj}
Show that every Zariski open subset of a quasiprojective variety $X \subseteq \pp^n$ has a {\em finite} open cover by subsets of the form $X\setminus V(f)$ for homogeneous polynomials $f \in \kk[x_0, \ldots, x_n]$. [Hint: use \cref{prop:IXICX}.]
\end{exercise}

\begin{exercise} \label{exercise:zariski-closure}
Given $X \subseteq \pp^n$, show that the closure of $X$ in $\pp^n$ under the Zariski topology is $\bar X := V(I(X)) = \{[x_0: \cdots :x_n] : f(x_0, \ldots, x_n) = 0$ for each homogeneous $f \in I(X)\}$; we say that $\bar X$ is the \index{Zariski!closure}{\em Zariski closure} of $X$ in $\pp^n$.
\end{exercise}

\begin{exercise} \label{exercise:nonempty-hypersurface-proj}
let $f$ be a homogeneous polynomial in $(x_0, \ldots, x_n)$, $n \geq 1$, and let $U_j = \pp^n \setminus V(x_j)$ be a basic open subset of $\pp^n$ from \cref{example:basic-clopen}.
\begin{enumerate}
\item \label{nonempty-hypersurface-proj:basic} Show that $V(f) \cap U_j = \emptyset$ if and only if $f = cx_j^m$ for some $m \geq 0$ and $c \in \kk\setminus \{0\}$. [Hint: let $m = \deg(f)$. Then $f/x_j^m$ induces a well defined map from $U_j$ to $\kk$ which can be expressed as a polynomial in $(x_i/x_j)_{i \neq j}$. Use \cref{exercise:nonempty-hypersurface-aff} and the one-to-one correspondence between $U_j$ and $\kk^n$ given by \eqref{map:U_j-to-k^n}.]
\item Deduce that if $f$ is not a constant, then the hypersurface $V(f) \subset \pp^n$ is nonempty. \item Show by examples that the preceding statement may be false if $\kk$ is not algebraically closed.
\end{enumerate}
\end{exercise}


\begin{exercise} \label{exercise:complement-hypersurface}
Let $S_1, S_2$ be finite subsets of $\pp^n$ such that $S_1 \cap S_2 = \emptyset$. Show that there is a hypersurface $X$ of $\pp^n$ containing $S_1$ but not containing any point of $S_2$. [Hint: for each $a \in S_1$, one can choose $b_0, \ldots, b_n \in \kk$ such that $b_0x_0 + \cdots + b_nx_n$ vanishes at $a$ but not at any point of $S_2$.]
\end{exercise}

\section{Regular functions} \label{quasi-regular-section}
If $f,g$ are homogeneous polynomials in $(x_0, \ldots, x_n)$ of the same degree, then $f/g$ is a well defined function on the Zariski open subset $\pp^n\setminus V(g)$ of $\pp^n$; we say that $f/g$ is a {\em rational function} on $\pp^n$ which is {\em regular} on $\pp^n \setminus V(g)$. In general, a \index{Regular function!on a variety}{\em regular function} on a quasiprojective variety $X$ is a function $\phi: X \to \kk$ which can be ``locally represented by rational functions,'' i.e.\ for each $x \in X$, there is an open neighborhood $U$ of $x$ in $X$ such that $\phi|_{U}$ is the restriction to $U$ of a rational function $f/g$ such that $g$ does {\em not} vanish at any point of $U$. The set of regular functions on $X$ has the natural structure of a $\kk$-algebra; we denote it by $\kk[X]$.

\begin{example} \label{example:basic-regular}
For each $i,j$, $x_i/x_j$ is a regular function on the basic open set $U_j := \pp^n \setminus V(x_j)$ of $\pp^n$. It follows that all polynomials in $(x_0/x_j, \ldots, x_n/x_j)$ are regular functions on $U_j$. We will shortly see (in \cref{prop:regular-equivalence} below) that these are in fact {\em all} regular functions on $U_j$.
\end{example}

\begin{example} \label{example:ruled-0}
Let $X := V(x_0x_3 - x_1x_2) \setminus V(x_1, x_3) \subseteq \pp^3$. Note that $X$ is the union of open subsets $X_1 := X \setminus V(x_1)$ and $X_3 := X \setminus V(x_3)$. Since $x_0/x_1 = x_2/x_3$ on $X_1 \cap X_3$, it follows that the function $f: X \to \kk$ defined by $x_0/x_1$ on $X_1$ and by $x_2/x_3$ on $X_3$ is a regular function on $X$. In \cref{exercise:ruled-0-regular} you will prove that $\kk[X] = \kk[f]$.
\end{example}

Given $J \subseteq \kk[X]$, we denote by $V(J)$ the ``subvariety of $X$ determined by $J$,'' i.e.\ the set of points on $X$ on which each $f \in J$ vanishes; \cref{prop:regular-zeroes} below implies that $V(J)$ is indeed a subvariety of $X$. \Cref{prop:regular-zeroes} can be verified directly from the definitions; we leave its proof as an exercise.

\begin{prop} \label{prop:regular-zeroes}
If $f$ is a regular function on a quasiprojective variety $X$, then $V(f) := \{x \in X: f(x) = 0\}$ is Zariski closed in $X$. If in addition $X$ is irreducible and $f$ is zero on a nonempty open subset of $X$, then it is zero everywhere on $X$. \qed
\end{prop}

Consider the basic open subsets $U_j$ of $\pp^n$, $j = 0, \ldots, n$. Recall (from \cref{example:basic-clopen}) that the map
\begin{align*}
U_j \ni [x_0: \cdots: x_n] \mapsto
	(\frac{x_0}{x_j}, \ldots, \frac{x_{j-1}}{x_j}, \frac{x_{j+1}}{x_j}, \ldots, \frac{x_n}{x_j}) \in \kk^n
\end{align*}
induces a bijection between $U_j$ and $\kk^n$. Denote the coordinates on $\kk^n$ by $(u_1, \ldots, u_n)$ so that the above map from $U_j$ to $\kk^n$ is given by
\begin{align}
u_i &=
	\begin{cases}
	x_{i-1}/x_j & \text{if}\ 1 \leq i \leq j \\
	x_{i+1}/x_j & \text{if}\ j < i \leq n.
	\end{cases}
	\label{map:U_j-to-k^n-1}
\end{align}
Recall that $\kk^n$ already comes with a Zariski topology and a set of regular functions; we next show that these objects are compatible with the corresponding objects on $U_j$ arising from its structure as a quasiprojective variety.

\begin{prop} \label{prop:Zariski-equivalence}
The Zariski topology on $U_j$ induced from $\pp^n$ is the same as the Zariski topology on $U_j$ induced from its identification with the affine space $\kk^n$ with coordinates $(u_1, \ldots, u_n)$ via the map given by \eqref{map:U_j-to-k^n-1}.
\end{prop}

\begin{proof}
A subset $V$ of $U_j$ which is closed with respect to the Zariski topology induced from $\kk^n$ is the set of zeroes of a collection of polynomials in $(x_i/x_j)_{i \neq j}$. Since $x_i/x_j$ are regular functions on $U_j$, \cref{prop:regular-zeroes} implies $V$ is closed with respect to the Zariski topology induced from $\pp^n$. Conversely, if a subset $V'$ of $U_j$ is closed with respect to the Zariski topology induced from $\pp^n$, \cref{exercise:basic-Zariski-equivalence} implies that $V'$ is also closed with respect to the Zariski topology induced from $\kk^n$, as required.
\end{proof}

\Cref{prop:Zariski-equivalence} in particular implies that ``subvarieties'' of $U_j$ are the same regardless of whether we identify $U_j$ with the affine space $\kk^n$ or an open subset of $\pp^n$. In \cref{prop:regular-equivalence} we will show ``regular functions'' also remain the same regardless of the consideration. Prior to that we give a description of regular functions on quasiprojective subsets of $U_j$ in terms of rational functions\footnote{Recall that classically a ``rational function'' means quotients of polynomials.} in $(u_1, \ldots, u_n)$.

\begin{prop}\label{prop:rational-on-Uj}
Let $X$ be a quasiprojective subset of $U_j$. Given a map $\phi : X \to \kk$, the following are equivalent:
\begin{enumerate}
\item $\phi$ is a regular function on $X$,
\item for each $x \in X$, there is an open neighborhood $U$ of $x$ in $X$ such that $\phi|_{U}$ can be represented as $f/g$ for some $f, g \in \kk[u_1, \ldots, u_n]$ such that $g$ does {\em not} vanish at any point of $U$.
\end{enumerate}
\end{prop}

\begin{proof}
This follows immediately from \cref{exercise:basic-Zariski-equivalence} and the observation that if $f,g$ are homogeneous polynomials of same degree $d$ in $(x_0, \ldots, x_n)$, then $f/g = (f/x_j^d)/(g/x_j^d)$.
\end{proof}
%
%

\begin{prop} \label{prop:regular-equivalence}
Let $X$ be a subvariety of $U_j$. The set of regular functions on $X$ (when $X$ is regarded as a quasiprojective variety) are precisely the restrictions of polynomials in $(u_1, \ldots, u_n)$. In particular,
\begin{align}
\kk[X] \cong \kk[u_1, \ldots,u_n] / I(X) \label{eqn:K[X]-affine-equivalence}
\end{align}
where\footnote{Note that $I(X)$ denotes somewhat different objects depending on whether $X$ is a projective variety or an affine variety. We hope that the intended meaning will always be clear from context.} $I(X)$ is the set of polynomials in $(u_1, \ldots, u_n)$ that are identically zero on $X$.
\end{prop}

\begin{proof}
\footnote{This proof is inspired by \cite[Proof of Proposition 1.11]{mummetry}.}Note that identity \eqref{eqn:K[X]-affine-equivalence} is an immediate consequence of the first assertion and identity \eqref{eqn:K[X]-affine} from \cref{regular-section-aff}. Consequently we only prove the first assertion. \Woutlog\ we may assume that $j = 0$. We will also identify a point on $U_0$ with the corresponding point in $\kk^n$ via \eqref{map:U_j-to-k^n-1}; in other words, we identify $a := (a_1, \ldots, a_n) \in \kk^n$ with the point $\hat a := [1: a_1: \cdots : a_n] \in U_0$. Since each $u_i$ is a regular function on $X$, we only need to show that every regular function on $X$ is the restriction of a polynomial in $(u_1, \ldots, u_n)$. Given a regular function $\rho$ on $X$, let $\qqq$ be the ideal of $\kk[u_1, \ldots, u_n]$ consisting of all polynomials $g$ such that $g\rho$ is the restriction of a polynomial in $(u_1, \ldots, n)$ on some nonempty open subset of $X$. It is clear that $I(X) \subseteq \qqq$, so that $V(\qqq) \subseteq X$ (where $V(\qqq)$ is the subvariety of $U_0$ determined by $\qqq$). On the other hand, \cref{prop:rational-on-Uj} implies that for each $a = (a_1, \ldots, a_n) \in X$, there is $g \in \qqq$ such that $g(a) \neq 0$. Taken together, these observations imply that $V(\qqq) = \emptyset$. The Nullstellensatz (\cref{thm:Hilbert-nulls}) then implies that $1 \in \qqq$, i.e.\ $\rho$ agrees with a polynomial on some nonempty open subset of $X$. In the case that $X$ is {\em irreducible}, it then follows that $\rho$ is a polynomial on all of $X$ (\cref{prop:regular-zeroes}) and the proposition is true. In the general case, since $X$ has finitely many {\em irreducible components} (due to \cref{thm:irred-components-affine,prop:Zariski-equivalence}), we can proceed by induction on the number of irreducible components of $X$. So assume that proposition is true for all subvarieties of $U_0$ with at most $k$ irreducible components, where $k \geq 1$, and that $X$ has $k+1$ irreducible components. Let $\rho$ be a regular function on $X$. Let $\rrr$ be the ideal of all polynomials $h \in \kk[u_1, \ldots, u_n]$ such that $h\rho$ is the restriction of a polynomial on {\em all of} $X$.

\begin{proclaim} \label{claim:regular-equivalence-2}
For each $a = (a_1, \ldots, a_n) \in X$, there is $g \in \rrr$ such that $g(a) \neq 0$.
\end{proclaim}

\begin{proof}
Let $X_1$ be an irreducible component of $X$, and $X_2$ be the union of the other irreducible components of $X$. Due to the induction hypothesis there are polynomials $f_1, f_2 \in \kk[u_1, \ldots, u_n]$ such that $\rho$ agrees with $f_i$ on $X_i$, $i = 1, 2$. \Cref{prop:rational-on-Uj} implies that there are polynomials $h',g' \in \kk[u_1, \ldots, u_n]$ and an open neighborhood $U'$ of $a$ in $X$ such that $g'$ does not vanish at any point on $U'$ and $g'\rho = h'$ on $U'$. Pick $g'' \in \kk[u_1, \ldots, u_n]$ such that $g''(a) \neq 0$ and $U \supseteq X \setminus V(g'')$ (this is possible, e.g., due to \cref{exercise:morphism-zariski-continuous}). Then $g'f_i - h' \equiv 0$ on $X_i \setminus V(g'')$, which implies that $g''(g'f_i - h') \in I(X_i)$. Consequently $(g''g'\rho)|_{X_i} = (g''h')|_{X_i}$ for each $i$, and the claim holds with $g := g'g''$.
\end{proof}

Since $I(X) \subseteq \rrr$, it follows (as in the case of $\qqq$) that $V(\rrr) = \emptyset$, and therefore $1 \in \rrr$. But then $\rho$ is the restriction of a polynomial in $(u_1, \ldots, u_n)$, as required.
\end{proof}

\begin{example} \label{example:Pn-constant}
\Cref{prop:regular-equivalence} in particular implies that $\kk[U_j] = \kk[x_0/x_j, \ldots, x_n/x_j]$ for each $j$. We now show that $\kk[\pp^n] = \kk$, i.e.\ the only regular functions on $\pp^n$ are the constants (in fact, we will see in \cref{complete-section} that this is true for {\em all} irreducible projective varieties). Write $y_i := x_i/x_0$, $i = 1, \ldots, n$. A regular function $\rho$ on $\pp^n$ restricts to $f \in \kk[x_1/x_0, \ldots, x_n/x_0] = \kk[y_1, \ldots, y_n]$ on $U_0$, and to $g \in \kk[x_0/x_1, x_2/x_1, \ldots, x_n/x_1] = \kk[1/y_1, y_2/y_1, \ldots, y_n/y_1]$ on $U_1$. Since $U_0 \setminus V(y_1) \subseteq U_1$, it follows that $f(y_1, \ldots, y_n) = g(1/y_1, y_2/y_1, \ldots, y_n/y_1)$ in the field of fractions of $\kk[y_1, \ldots, y_n]$ (\cref{exercise:rational-equality-Kn}). Equating the degrees\footnote{Recall that the {\em degree} of $f_1/f_2$, where $f_1, f_2$ are polynomials, is defined as $\deg(f_1) - \deg(f_2)$.} of both sides shows that $\deg(f) = 0$, i.e.\ $f$ is a {\em constant}. Since $\pp^n$ is irreducible (\cref{exercise:Pn-irreducible}), it follows that $\rho$ is constant on all of $\pp^n$.
\end{example}

\subsection{Exercises}

\begin{exercise}\label{exercise:regular-zeroes}
Prove \cref{prop:regular-zeroes}.
\end{exercise}

\begin{exercise}\label{exercise:basic-Zariski-equivalence}
Let $f$ be a homogeneous polynomial of degree $d$ in $(x_0, \ldots, x_n)$ and $U_j = \pp^n \setminus V(x_j)$, $0 \leq j \leq n$, be a basic open subset of $\pp^n$. Show that
\begin{enumerate}
\item $f/x_j^d$ can be expressed as a polynomial in $(x_i/x_j)_{i \neq j}$.
\item $V(f) \cap U_j = \{x \in U_j: f/x_j^d = 0\}$.
\item Deduce that every Zariski closed subset of $U_j$ is the set of zeroes on $U_j$ of a collection of polynomials in $(x_i/x_j)_{i \neq j}$.
\end{enumerate}
\end{exercise}

\begin{exercise} \label{exercise:Zariski-continuous}
Due to \cref{prop:Zariski-equivalence} we can treat $\kk$ as a quasiprojective variety with Zariski topology induced from $\pp^1$. Prove that every regular function $f: X \to \kk$ on a quasiprojective variety $X$ is \index{Zariski!continuous}{\em Zariski continuious}, i.e.\ continuous with respect to Zariski topology.
\end{exercise}

\begin{exercise} \label{exercise:quasiaffine}
Let $X$ be a quasiprojective subset of $\kk^n$ (when $\kk^n$ is treated as a quasiprojective variety)\footnote{These sets are sometimes referred to as {\em quasiaffine} varieties.}. Show that $X$ is an open subset of a subvariety of $\kk^n$.
\end{exercise}

\begin{exercise} \label{exercise:rational-equality-Kn}
Let $f_1, g_1, f_2, g_2 \in k[x_1, \ldots, x_n]$, where $k$ is an infinite field and $x_1, \ldots, x_n$ are indeterminates over $k$. Assume $f_1(a)/g_1(a) = f_2(a)/g_2(a)$ for all $a \in k^n \setminus \{x: g_1(x)g_2(x) = 0\}$. Show that $f_1/g_1 = f_2/g_2$ as {\em rational functions} in $(x_1, \ldots, x_n)$, i.e.\ as elements of the field of fractions of $k[x_1, \ldots, x_n]$. [Hint: use \cref{exercise:pol-zero-on-Kn}.]
\end{exercise}

\section{Morphisms of quasiprojective varieties; affine varieties as quasiprojective varieties}
A map $\phi: Y \to X$ of quasiprojective varieties is called a \index{Morphism}{\em morphism} if for every Zariski open subset $U$ of $X$ and every regular function $h$ on $U$, the pullback $h \circ \phi$ is a regular function on $\phi^{-1}(U)$. An \index{Isomorphism of varieties}{\em isomorphism} is a bijective morphism whose inverse is also a morphism; we write $X \cong Y$ to denote that $X$ and $Y$ are {\em isomorphic}, i.e.\ there is an isomorphism between $X$ and $Y$. We say that a morphism $\phi:Y \to X$ is a \index{Closed!embedding of varieties}{\em closed embedding} if $\phi(Y)$ is a (closed) subvariety of $X$ and $\phi$ is an isomorphism between $Y$ and $\phi(Y)$. In \cref{quasi-regular-section} we have seen that the identification of $\kk^n$ with a basic open set of $\pp^n$ imbues every subvariety of $\kk^n$ with the structure of a quasiprojective variety. From now on we treat every affine variety as a quasiprojective variety via this identification. We saw in \cref{quasi-regular-section} that this identification is compatible with the topology or regular functions on affine varieties. In \cref{prop:morphism-equivalence} below we show that this is also compatible with morphisms, i.e.\ a map between affine varieties is a morphism (respectively, isomorphism) of affine varieties if and only if it is a morphism (respectivel, isomorphism) of quasiprojective varieties. First we go over some examples and basic properties of morphisms of quasiprojective varieties.

\begin{example}
Morphisms from a quasiprojective variety $Y$ to $\kk$ are precisely regular functions on $Y$.
\end{example}

\begin{example} \label{example:morphism-Pn-to-affine}
If $\phi: \pp^n \to \kk^m$ is a morphism for some $m, n$, then the pullback of each coordinate on $\kk^m$ is a regular function on $\pp^n$, and therefore is a constant (\cref{example:Pn-constant}). It follows that the only possible morphism from a projective space to an affine space is a trivial map which maps everything to a point. In fact we will see that this property holds for {\em all} irreducible projective varieties (\cref{exercise:morphism-proj-to-affine}).
\end{example}

The following result is immediate from the definitions; we leave its proof as an exercise.

\begin{prop} \label{prop:quasi-morphism-pullback}
A morphism $\phi: Y \to X$ of quasiprojective varieties induces via pullback a $\kk$-algebra homomorphism $\phi^*: \kk[X] \to \kk[Y]$. If $\phi$ is an isomorphism of quasiprojective varieties, then $\phi^*$ is an isomorphism of $\kk$-algebras.
\end{prop}

\begin{example}
A constant morphism from $\pp^n$ to a point for $n \geq 1$ shows that the converse of the last assertion of \cref{prop:quasi-morphism-pullback} is {\em not} true in general. Contrast this with the case of morphisms between {\em affine} varieties (\cref{exercise:affine-pullback}).
\end{example}

Our next result shows that for quasiprojective subsets of affine spaces, morphisms can be characterized in terms of the affine coordinates. It is a counterpart of \cref{prop:rational-on-Uj} which characterizes regular functions on these spaces in terms of the affine coordinates. The proof is straightforward using \cref{prop:rational-on-Uj} and the definitions of rational functions and morphisms on quasiprojective varieties - we leave it as an exercise.


\begin{prop} \label{prop:morphism-on-Uj}
Let $X$ be a quasiprojective subset of $\kk^n$ with coordinates $(x_1, \ldots, x_n)$, and $Y$ be a quasiprojective variety. Given a map $\phi: Y \to X$, the following are equivalent:
\begin{enumerate}
\item \label{morphism-on-Uj:0} $\phi$ is a morphism of quasiprojective varieties;
\item \label{morphism-on-Uj:1} for each $i$, $x_i \circ \phi$ is a regular function on $Y$.
\end{enumerate}
If $Y$ is a quasiprojective subset of an affine space $\kk^m$ with coordinates $(y_1, \ldots, y_m)$, then the preceding properties are equivalent to the following:
\begin{enumerate}[resume]
\item \label{morphism-on-Uj:2} for each $i$, $x_i \circ \phi$ can be locally represented on $Y$ by rational functions in $(y_1, \ldots, y_m)$.  \qed
\end{enumerate}
\end{prop}

\begin{cor} \label{prop:morphism-equivalence}
A map between affine varieties is a morphism (respectively, isomorphism) of quasiprojective varieties if and only if it is a morphism (respectively, isomorphism) of affine varities.
\end{cor}

\begin{proof}
Since regular functions on affine varieties are precisely polynomials in the coordinate functions (\cref{prop:regular-equivalence}), \cref{prop:morphism-on-Uj} implies that a map $Y \to X$ between affine varieties is a morphism of quasiprojective varieties if and only if the pullback of all regular functions on $X$ is a regular function on $Y$. The latter property is precisely what defines a morphism of affine varieties.
\end{proof}

\begin{example}
\Cref{prop:morphism-equivalence} in particular implies that all the examples of morphisms of affine varieties from \cref{regular-section-aff} are also examples of morphisms between quasiprojective varieties. E.g.\ the map $\kk \to \kk^2$ given by $x \mapsto (x, x^2)$ defines an isomorphism between $\kk$ and the {\em parabola} $V(y - x^2)\subseteq \kk^2$. Since the parabola is closed in $\kk^2$, this is actually a closed embedding.
\end{example}

\begin{example}[Continuation of {\cref{example:Kn-minus-V(f)}}] \label{example:affine-complement}
Let $X = V(\qqq) \subseteq \kk^n$, where $\qqq$ is an ideal in $\kk[x_1, \ldots, x_n]$. Given a polynomial $f$ in $(x_1, \ldots, x_n)$, define $Y := V(\qqq, x_{n+1}f - 1) \subseteq \kk^{n+1}$. The projection onto $(x_1, \ldots, x_n)$-coordinates maps $Y$ bijectively onto the Zariski open subset $X \setminus V(f)$ of $X$, with its inverse given by $(x_1, \ldots, x_n) \mapsto (x_1, \ldots, x_n, 1/f)$. \Cref{prop:morphism-on-Uj} implies that both these maps are morphisms, so that $X \setminus V(f) \cong Y$ as quasiprojective varieties; in particular, $X \setminus V(f)$ is isomorphic to an {\em affine variety}. It follows (due to \cref{prop:quasi-morphism-pullback}) that
\begin{align}
\kk[X \setminus V(f)]
	\cong \kk[Y]
	\cong \kk[x_1, \ldots, x_n] /\langle I(X), x_{n+1}f - 1 \rangle
	\cong \kk[X]_f
\label{eqn:affine-complement}
\end{align}
where $\kk[X]_f$ is the {\em localization} (see \cref{local-ring-section}) of $\kk[X]$ at $f$. Taking $n = 1$, $\qqq =$ the zero ideal, and $f = x_1$ yields that $\kk\setminus \{0\} \cong V(x_1x_2 -1)$ and $\kk[\kk \setminus \{0\}] \cong \kk[x_1, 1/x_1]$ (\cref{fig:hyperbola}).
\end{example}

\begin{figure}[htb]
\begin{center}
\def\colorzero{blue}
\def\colorone{red}
\begin{tikzpicture}[scale=0.45]

\def\nsamples{103}
\def\comphy{0}
\def\compvx{0}
\def\maxx{4}
\def\maxy{4}

\def\xmin{-\maxx}
\def\xmax{\maxx}
\def\ymin{-\maxy}
\def\ymax{\maxy}
\def\xminleft{\xmin}
\def\xmaxleft{-0.1}
\def\xminright{0.1}
\def\xmaxright{\xmax}

\pgfmathsetmacro{\minx}{1/\maxx};

\draw[\colorzero, thick, smooth, samples=\nsamples, domain=\xminleft:-\minx] plot(\x,{1/\x});
\draw[\colorzero, thick, smooth, samples=\nsamples, domain=\minx:\xmaxright] plot(\x,{1/\x});
\draw[samples=2, domain=\xmin:\xmax] plot ({\x} ,{\comphy});
\draw[samples=2, domain=\ymin:\ymax] plot ({\compvx}, {\x});
%
\end{tikzpicture}
\caption{The projection onto $x_1$-axis induces an isomorphism between the hyperbola $x_1x_2 - 1 = 0$ and $\kk \setminus \{0\}$}
\label{fig:hyperbola}
\end{center}
\end{figure}

Using \cref{example:affine-complement} we can show that all quasiprojective varieties are built by ``gluing'' finitely many affine varieties:

\begin{prop} \label{affine-cover}
Every quasiprojective variety has a finite open cover by open subsets isomorphic to affine varieties.
\end{prop}

\begin{proof}
Let $W$ be the open subset of a subvariety of $\pp^n$. Then $W$ is covered by $W \cap U_j$, where $U_j$, $j = 0, \ldots, n$, are the basic open subsets of $\pp^n$. The correspondence between $U_j$ and $\kk^n$ identifies $W \cap U_j$ with an open subset of an affine variety, which has a finite open covering by subsets isomorphic to $X \setminus V(f)$ for an affine variety $X$ and a regular function $f$ on $X$ (\cref{exercise:morphism-zariski-continuous}). Since $X\setminus V(f)$ is isomorphic to an affine variety (\cref{example:affine-complement}), the proof is complete.
\end{proof}

The identification of $\kk^n$ with the basic open subset $U_0$ of $\pp^n$ identifies all subvarieties of $\kk^n$ with quasiprojective subsets of $\pp^n$. We next describe the closure of these sets in $\pp^n$. For each $f \in \kk[x_1, \ldots, x_n]$, the \index{Homogenization}{\em homogenization} of $f$ in $(x_0, \ldots, x_n)$ is the homogeneous polynomial $\tilde f := x_0^{\deg(f)}f(x_1/x_0, \ldots, x_n/x_0)$. The homogenization of an ideal $\qqq$ of $\kk[x_1, \ldots, x_n]$ in $(x_0, \ldots, x_n)$ is the ideal of $\kk[x_0, \ldots, x_n]$ generated by the homogenization of all $f \in \qqq$.

\begin{prop}  \label{prop:projective-closure}
Let $\qqq$ be the any ideal of $\kk[x_1, \ldots, x_n]$. The Zariski closure of $V(\qqq) \subseteq \kk^n$ in $\pp^n$ is the subvariety $V(\tilde \qqq)$ of $\pp^n$ defined by the homomgenization of $\qqq$.
\end{prop}

\begin{proof}
Recall that the identification between $\kk^n$ and $U_0$ is given by the map $(x_1, \ldots, x_n) \mapsto [1: x_1: \cdots : x_n]$. It follows that the closure of $X := V(\qqq)$ in $\pp^n$ is
\begin{align*}
\bar X
	&= \bigcap\{ V(h): h\ \text{homogeneous in } (x_0, \ldots, x_n),\ h(1, a_1, \ldots, a_n) = 0\ \text{for all}\ (a_1, \ldots, a_n) \in X\} \\
	&= \bigcap\{ V(h): h\ \text{homogeneous in } (x_0, \ldots, x_n),\ h(1, x_1, \ldots, x_n) \in \sqrt \qqq\} \\
	&= \bigcap\{ V(h): h\ \text{homogeneous in } (x_0, \ldots, x_n),\ h^k(1, x_1, \ldots, x_n) \in \qqq\ \text{for some}\ k \geq 0\} \\
	&= \bigcap\{ V(h): h\ \text{homogeneous in } (x_0, \ldots, x_n),\ h^k = x_0^m \tilde f\ \text{for some}\ f \in \qqq,\ k, m \geq 0\} \\
	&= \bigcap\{ V(h): h\ \text{homogeneous in } (x_0, \ldots, x_n),\ h^k \in \tilde \qqq\ \text{for some}\ k \geq 0\} \\
	&= V(\tilde \qqq)
\end{align*}
where the last equality follows from the correspondence between subvarieties of projective varieties and radical homogeneous ideals of $\kk[x_0, \ldots, x_n]$ (\cref{exercise:projective-nullspondence}).
\end{proof}

\begin{example} \label{example:principal-closure}
If $\qqq$ is a principal ideal generated by $f \in \kk[x_1, \ldots, x_n]$, then $\tilde \qqq$ is also a principal ideal generated by $\tilde f$, which has the same degree as $f$. It follows that if $X$ is the hypersurface (respectively, \index{Hyperplane}{\em hyperplane}\footnote{A {\em hyperplane} is a hypersurface defined by a linear polynomial.}) of $\kk^n$ determined by $f$, then $\bar X$ is the hypersurface (respectively, hyperplane) of $\pp^n$ determined by $\tilde f$. Note that the set of ``points at infinity'' on $\bar X$ is $V(\tilde f) \setminus U_0 = V(\tilde f) \cap V(x_0) = \{[0: x_1: \cdots : x_n] \in \pp^n: \ld(f)(x_1, \ldots, x_n) = 0\}$ where $\ld(f)$ is the {\em leading form}\footnote{The leading form of a polynomial $f = \sum c_\alpha x^\alpha$ of degree $d$ is the sum of all monomial terms $c_\alpha x^\alpha$ of degree $d$.} of $f$;  this is the relation which yields the correspondence \eqref{bezout-infinity-correspondence} from \cref{brief-chapter}.  
\end{example}

\begin{example}
In general, if $X = V(f_1, \ldots, f_k)$, then $\bar X \subsetneq V(\tilde f_1, \ldots, \tilde f_k)$. E.g.\ if $f_1 = x_1$ and $f_2 = x_2 - x_1^2$, then $X := V(f_1, f_2) \subset \kk^2$ is the singleton consisting of the origin, so that $\bar X = X$, whereas $V(\tilde f_1, \tilde f_2) = V(x_1, x_0x_2 - x_1^2) = \{[1:0:0], [0:0:1]\} \subset \pp^2$ (cf.\ \cref{bezout-counterexample}).
\end{example}

\subsection{Exercises}
\begin{exercise} \label{exercise:morphism-zariski}
\begin{enumerate}
\item Show that every morphism is {\em Zariski continuous}, i.e.\ continuous with respect to the Zariski topology.
\item Show that every bijection of $\pp^1$ is Zariski continuous.
\item Deduce that there are Zariski continuous maps from $\pp^1$ to $\pp^1$ which are not morphisms.
\end{enumerate}
\end{exercise}

\begin{exercise} \label{exercise:quasi-morphism-pullback}
Prove \cref{prop:quasi-morphism-pullback}.
\end{exercise}

\begin{exercise}\label{exercise:morphism-on-Uj}
Prove \cref{prop:morphism-on-Uj}. [Hint: use \cref{prop:rational-on-Uj}.]
\end{exercise}

\begin{exercise} \label{exercise:affine-complement-1}
Verify the isomorphisms in \eqref{eqn:affine-complement} from \cref{example:affine-complement}.
\end{exercise}

\begin{exercise}[Converse to \cref{prop:projective-closure}] \label{exercise:basic-open-equation}
Let $V$ be a a subvariety of $\pp^n$ defined by a homogeneous ideal $\qqq$ of $\kk[x_0, \ldots, x_n]$. The basic open subset $U_0$ of $\pp^n$ is isomorphic to $\kk^n$ with coordinates $(x_1/x_0, \ldots, x_n/x_0)$. Let
\begin{align*}
\qqq_0 := \{f/x_0^{\deg(f)}: f\ \text{is a homogeneous polynomial in}\ \qqq\}
\end{align*}
Show that
\begin{enumerate}
\item $\qqq_0$ is an ideal of $\kk[x_1/x_0, \ldots, x_n/x_0]$.
\item $V \cap U_0$ is the subvariety of $U_0$ defined by $\qqq_0$. In particular, the restriction of a projective hypersurface to $U_0$ is an affine hypersurface. 
\end{enumerate}
\end{exercise}

\begin{exercise} \label{exercise:ruled-0-regular}
Let $X := V(x_0x_3 - x_1x_2) \setminus V(x_1, x_3) \subseteq \pp^3$. 
\begin{enumerate}
\item Identify $U_3 := \pp^3 \setminus V(x_3)$ with $\kk^3$ with coordinates $(x,y,z) := (x_0/x_3, x_1/x_3, x_2/x_3)$. Show that $X \cap U_3 = V(x - yz)$. [Hint: use \cref{exercise:basic-open-equation}.]
\item Deduce that $\kk[X \cap U_3] \cong \kk[y, z]$.
\item Similarly show that $\kk[X \cap U_1] \cong \kk[u, v]$, where $(u,v) = (x_0/x_1, x_3/x_1)$.  
\item Deduce that $\kk[X] \cong \kk[z] \cong \kk[u]$. 
\end{enumerate}
\end{exercise}

\begin{exercise}[Local nature of morphisms] \label{exercise:quasi-morphisms-local-defn}
Given a map $\phi: Y \to X$ between quasiprojective varieties, show that the following are equivalent:
\begin{enumerate}
\item $\phi$ is a morphism.
\item for every $a \in X$ and for each open neighborhood $U$ of $a$ in $X$, there is an open neighborhood $U'$ of $a$ in $X$ such that $U' \subseteq U$ and $h \circ \phi$ is a regular function on $\phi^{-1}(U')$ for all $h \in \kk[U']$.
\end{enumerate}
\end{exercise}

\begin{exercise} \label{exercise:check-morphism}
Let $[x_0: \cdots : x_n]$ be homogeneous coordinates on $\pp^n$ and $U_j := \pp^n \setminus V(x_j)$, $j = 0, \ldots, n$, be the basic open subsets of $\pp^n$. Let $X$ be a quasiprojective subset of $\pp^n$ and $Y$ be a quasiprojective variety. Given a map $\phi: Y \to X$, show that the following are equivalent:
\begin{enumerate}
\item \label{check-morphism:0} $\phi$ is a morphism.
\item \label{check-morphism:1} $(x_i/x_j) \circ \phi$ is a regular function on $\phi^{-1}(U_j \cap X)$ for each $i, j$.
\end{enumerate}
[Hint: use \cref{exercise:quasi-morphisms-local-defn,prop:morphism-on-Uj}.]
\end{exercise}

\begin{exercise} \label{exercise:morphism-dense=equal}
Let $\phi_j:X \to Y$ be morphisms of varieties, $j = 1, 2$. Assume $\phi_1(x) = \phi_2(x)$ for all $x$ on a dense subset of $X$. Then show that $\phi_1 = \phi_2$. [Hint: use \cref{exercise:quasi-morphisms-local-defn,prop:regular-zeroes}.]
\end{exercise}

\begin{exercise} \label{exercise:Pn-linear-auto}
Consider an $(n+1) \times (n+1)$-matrix $A$ over $\kk$ with entries $a_{i,j}$, $i,j = 0, \ldots, n$. If $\det(A) \neq 0$, then show that the linear isomorphism of $\kk^{n+1}$ given by $x \mapsto Ax$ (where $x$ is regarded as a column vector with $(n+1)$-rows) induces an {\em automorphism}\footnote{Recall that an {\em autmorphism} of $X$ is an isomorphism from $X$ to itself. In fact it turns out that {\em all} automorphisms of $\pp^n$ are of this form (see e.g.\ \cite[Exercise III.1.17]{shaf1}).} of $\pp^n$. [Hint: use \cref{exercise:check-morphism}.]
\end{exercise}

\begin{exercise} \label{exercise:veronese-2}
Show that the map $\phi: \pp^1 \to \pp^2$ given by $[x_0: x_1] \to [x_0^2: x_0x_1: x_2^2]$ induces an isomorphism between $\pp^1$ and the subvariety $X := V(x_0x_2 - x_1^2)$ of $\pp^2$. Since $X$ is closed in $\pp^3$, $\phi$ is in fact a {\em closed embedding} of $\pp^1$ into $\pp^2$.
\end{exercise}

\begin{exercise} \label{exercise:projective-closure-example}
Let $a = (a_1, \ldots, a_n), b = (b_1, \ldots, b_n) \in \kk^n$, $b \neq 0$, and $L$ be the line $\{a + bt: t \in \kk\} \subset \kk^n$. We consider $L$ as a subest of $\pp^n$ by identifying $\kk^n$ with the basic open subset $U_0 = \pp^n \setminus V(x_0)$ of $\pp^n$. In this exercise you will calculate the closure $\bar L$ of $L$ in $\pp^n$.
\begin{enumerate}
\item Show that there is an automorphism $\phi$ of $\pp^n$ such that $\phi(L) = \{[1: t: 0: \cdots: 0]: t \in \kk\}$. [Hint: use \cref{exercise:Pn-linear-auto}.]
\item Show that the homogeneous ideal generated by all polynomials vanishing on $\phi(L)$ is generated by $x_2, \ldots, x_n$. Compute the closure of $\phi(L)$ in $\pp^n$.
\item Show that the closure $\bar L$ of $L$ in $\pp^n$ is $L \cup \{[0: b_1: \cdots :b_n]\}$.
\item Conclude that two lines in $\kk^n$ intersect at infinity in $\pp^n$ if and only if they are parallel, i.e.\ have the same ``direction vector.''
\end{enumerate}
\end{exercise}

\section{Rational functions and rational maps on irreducible varieties} \label{rational-section}
Let $X$ be an irreducible quasiprojective variety. A \index{Rational!function!on a variety}{\em rational function} on $X$ is a regular function on a nonempty Zariski open subset of $X$. Formally, define a binary relation $\sim$ on the collection of pairs $(f, U)$, where $U$ is a nonempty Zariski open subset of $X$ and $f$ is a regular function on $U$, as follows: 
\begin{align*}
(f, U) \sim (f',U')\ \text{if and only if $f$ and $f'$ agree on}\ U \cap U'
\end{align*}
we will see in \cref{prop:k(X)-0} below that $\sim$ is an {\em equivalence relation}; a {\em rational function} is an equivalence class of $\sim$. We write $\kk(X)$ for the set of rational functions on $X$.


\begin{prop} \label{prop:k(X)-0}
Let $X$ be an irreducible variety.
\begin{enumerate}
\item \label{k(X)-0:equivalence} $\sim$ is an equivalence relation. In particular, rational functions on $X$ are well defined. 
\item \label{k(X)-0:open} If $U$ is a nonempty open subset of $X$, then $\kk(X) \cong \kk(U)$.
\item \label{k(X)-0:field} $\kk(X)$ is a field.
\item \label{k(X)-0:affine} If $X$ is also affine, then $\kk(X)$ is isomorphic to the field of fractions of $\kk[X]$.
\end{enumerate}
\end{prop}

\begin{proof}
It  is clear that $(f, U) \sim (f, U)$, and if $(f, U) \sim (f', U')$ then $(f', U') \sim (f, U)$. Now, given $(f, U) \sim (f', U')$ and $(f', U') \sim (f'', U'')$, it is clear that $f = f''$ on $U \cap U' \cap U''$. Since $X$ is irreducible, it follows that $U \cap U' \cap U''$ is nonempty, and it is dense in $U \cap U''$ (\cref{prop:irreducible-properties}). Consequently $f = f''$ on $U \cap U''$ (\cref{prop:regular-zeroes}) and $(f, U) \sim (f'', U'')$. This proves the first assertion. Similarly, given open subsets $U, U'$ of $X$, $U \cap U' \neq \emptyset$, and for every regular function $f'$ on $U'$, $(f', U') \sim (f'|_{U \cap U'}, U \cap U')$. This immediately implies that $\kk(X) \cong \kk(U)$. If $(f, U) \in \kk(X)$ and $f$ is not identically zero on $U$, then it is nonzero on a nonempty open subset $U''$ of $X$ (\cref{prop:regular-zeroes}) and $(1/f, U'')$ is the multiplicative inverse of $(f, U)$ in $\kk(X)$. This implies that $\kk(X)$ is a field. The last assertion then follows from the fact that regular functions on a nonempty open subset of an affine variety $X$ can be represented as quotients of regular functions on $X$ (\cref{prop:rational-on-Uj}).
\end{proof}

\begin{example}
Assertion \eqref{k(X)-0:affine} of \cref{prop:k(X)-0} implies that the rational functions on $\kk^n$ are the quotients of polynomials in $(x_1, \ldots, x_n)$, i.e.\ the notion of ``rational functions'' is compatible with its classical usage.
\end{example}

\begin{example}
Assertion \eqref{k(X)-0:open} of \cref{prop:k(X)-0} implies that $\kk(\pp^n) \cong \kk(\kk^n) \cong \kk(x_1, \ldots, x_n)$. Since $\kk[\pp^n] = \kk$ (\cref{example:Pn-constant}), assertion \eqref{k(X)-0:affine} of \cref{prop:k(X)-0} in general does not hold when $X$ is not affine.
\end{example}

A \index{Rational!map}{\em rational map} $\phi$ from an irreducible variety $X$ to a variety $Y$ is a morphism from a nonempty open subset $U$ of $X$ to $Y$; usually a rational map is denoted by a broken arrow
\begin{align*}
\phi : X \dashrightarrow Y
\end{align*}
It is called \index{Dominant!rational map}{\em dominant} if the image of $U$ is dense in $Y$, and \index{Birational}{\em birational} if $Y$ is irreducible and there is a rational map $\psi: Y \dashrightarrow X$ such that $\psi \circ \phi$ (respectively, $\phi \circ \psi$) restricts to an automorphism of a nonempty open subset of $X$ (respectively, $Y$). We say that {\em $X$ and $Y$ are birational} or {\em $X$ is birational to $Y$} if there is a birational map from $X$ to $Y$. A \index{Rational!variety}\index{Variety!rational}{\em rational variety} is a variety birational to $\kk^n$ for some $n \geq 0$.

\begin{example}
Every morphism is trivially a rational map. The mapping $x \mapsto (x, 1/x)$ defines a rational map $\kk \dashrightarrow \kk^2$ which is {\em not} a morphism, and {\em not} dominant. This map induces an isomorphism between $\kk \setminus \{0\}$ and the hyperbola $H := V(xy - 1) \subseteq \kk^2$ (\cref{example:affine-complement}). It follows that $\kk$ and $H$ are birational, and $H$ is a rational variety. 
\end{example}

\begin{example}
An irreducible variety is trivially birational to all of its nonempty open subsets. In particular, $\pp^n$ is a rational variety. 
\end{example}

\begin{example}
A basic source of rational maps which are not morphisms is the projection from $\pp^n$ to a coordinate subspace. E.g., given $m < n$, the map $[x_0: \cdots : x_n] \mapsto [x_0: \cdots : x_m]$ defines a rational map $\pp^n \dashrightarrow \pp^m$ which is not defined on the subvariety $Z := V(x_0, \ldots, x_m)$ of $\pp^n$. Note that $Z \cong \pp^{n - m - 1}$. 
\end{example}

We outline the proof of the following properties of birational maps in \crefrange{exercise:dominant-k(X)}{exercise:birational-hypersurface} below.

\begin{prop} \label{prop:birational-properties}
Two varieties are birational if and only if their fields of rational functions are isomorphic. Every irreducible variety is birational to a hypersurface. 
\end{prop} 

\begin{example}
For $X := V(x^3 - y^2) \subseteq \kk^2$, we have that $\kk[X] \cong \kk[x, y]/\langle x^3 - y^2 \rangle$. Since $x = (y/x)^2$ in $\kk(X)$, it follows that $\kk(X) \cong \kk(y/x)$. \Cref{prop:birational-properties} therefore implies that $X$ is birational to $\kk$; in particular $X$ is rational. Indeed, we saw in \cref{example:k-to-cusp-cubic} that the map $\phi:t \mapsto (t^2, t^3)$ is a bijective morphism from $\kk$ to $X$. It is straightforward to check that $\phi$ induces an isomorphism between $\kk \setminus \{0\}$ and $X \setminus \{(0,0)\}$ with the inverse given by $t = y/x$. 
\end{example}

Cubic curves $V(y^2 - x(x-1)(x- \lambda)) \subseteq \kk^2$ are {\em not} birational to $\kk$ for $\lambda \neq 0, 1$ and $\character(\kk) \neq 2$. Usually this is proven using ``diffrential forms,'' which we do not develop in this book; see \cite[Section 2.2]{reid-undergrad} for an elementary proof without using differential forms. 

\subsection{Exercises}

\begin{exercise} \label{exercise:k(X)-reducible}
Let $X$ be a reducible affine variety. Show that
\begin{enumerate}
\item There are nonzero regular functions $f, g$ on $X$, such that $fg = 0$ on $X$.
\item The relation $\sim$ from the definition of rational functions is {\em not} an equivalence relation when applied to $X$. [Hint: $(f, X \setminus V(f)) \sim (g, X \setminus V(g)) \sim (\alpha f, X \setminus V(f))$, where $\alpha \neq 0, 1$.]
\item The above problem can be rectified if in the definition of rational functions we only allow open subsets which have nonempty inersection with {\em every} irreducible component of $X$. 
\item The resulting set of ``rational functions'' on $X$ is the localization of $\kk[X]$ at the set of non zero-divisors; in other words, $\kk(X) \cong \{f/g: f, g \in \kk[X],\ g$ is {\em not} a zero-divisor in $\kk[X]\}$.  
\end{enumerate}
\end{exercise}
 
\begin{exercise} \label{exercise:dominantly-irreducible}
Let $\phi: X \to Y$ be a dominant rational map. Show that if $X$ is irreducible, then so is $Y$. [Hint: morphisms are continuous maps (\cref{exercise:morphism-zariski}); use \cref{prop:irreducible-properties}.]
\end{exercise}

\begin{exercise} \label{exercise:dominant-k(X)}
Let $\phi$ be a rational map from $X$ to $Y$, both irreducible varieties. Show that
\begin{enumerate}
\item If $\phi$ is dominant, then $\phi$ induces an injection $\phi^*: \kk(Y) \into \kk(X)$. [Hint: if $f$ is a nonzero regular function on $Y$, then show that $f$ can not be identically zero on the image of $X$.]
\item Conversely, if $\phi$ induces a well defined map from $\kk(Y) \to \kk(X)$, then $\phi$ must be dominant.   
\end{enumerate}
\end{exercise}

\begin{exercise} \label{exercise:birational-fields}
Let $X$ and $Y$ be irreducible varieties. In this exercise you will show that $X$ is birational to $Y$ if and only if $\kk(X) \cong \kk(Y)$. 
\begin{enumerate}
\item Show that if $X$ is birational to $Y$ then $\kk(X) \cong \kk(Y)$. [Hint: use \cref{exercise:dominant-k(X)}.]
\item Assume there is an isomorphism $F: \kk(Y) \to \kk(X)$. Choose an affine open subset of $Y$ isomorphic to a subvariety $Y'$ of $\kk^m$ with coordinates $(y_1, \ldots, y_m)$. Show that
\begin{enumerate}
\item The correspondence $x \mapsto (F(y_1)(x),  \ldots, F(y_m)(x))$ defines a well-defined rational map $\phi: X \dashrightarrow \kk^m$. 
\item The image of $\phi$ is a dense subset of $Y'$. [Hint: a polynomial $g$ in $(y_1, \ldots, y_m)$ vanishes on the image of $\phi$ if and only if $g \in I(Y')$.]
\item $\phi$ induces a birational map from $X$ to $Y$. 
\end{enumerate}
\end{enumerate} 
\end{exercise}

\begin{exercise} \label{exercise:birational-hypersurface}
Let $X$ be an irreducible variety. In this exercise you will show that there is $n \geq 1$ and a polynomial $f$ in $(x_1, \ldots, x_n)$ such that $X$ is birational to $V(f) \subset \kk^n$.
\begin{enumerate}
\item Show that $\kk(X)$ is a finitely generated field extension of $\kk$. [Hint: it suffices to consider the case that $X$ is affine. Use identity \eqref{eqn:K[X]-affine} and \cref{prop:k(X)-0}.]
\item Use \cref{thm:Schmidt-infinite} to show that there are $x_1, \ldots, x_n \in \kk(X)$ such that
\begin{enumerate}
\item $x_1, \ldots, x_{n-1}$ are algebraically independent over $\kk$,
\item $x_n$ is algebraic over $\kk(x_1, \ldots, x_{n-1})$,
\item $f(x_n) = 0$ for some irreducible $f \in \kk[x_1, \ldots, x_n]$, and
\item $\kk(X) = \kk(x_1, \ldots, x_n)$.
\end{enumerate}
\item Compute the field of rational functions of $Y := V(f) \subset \kk^n$ and use \cref{exercise:birational-fields} to conclude that $X$ is birational to $Y$.
\end{enumerate}
\end{exercise}

\begin{exercise}
Given a rational map $\phi:X \to Y$ between irreducible varieties, show that the following are equivalent: 
\begin{enumerate}
\item $\phi$ is a birational map;
\item there are open subsets $U'$ of $U$ and $V$ of $Y$ such that $\phi|_{U'}: U' \to V$ is an isomorphism.
\end{enumerate}
 \end{exercise}

\section{Product spaces, Segre map, Veronese embedding} \label{product-section}
\subsection{Product spaces, Segre map}
\label{segre-subsection}
Usually the topology considered on products of topological spaces is the ``product topology,'' whose open sets are unions of products of open subsets of each factor. The product topology on algebraic varieties however is very restrictive:

\begin{example} \label{example:product-top}
Since the only proper closed subsets of $\kk$ are finite sets of points (\cref{example:subvarieties-aff1}), the proper closed subsets of $\kk^2 \cong \kk \times \kk$ under the product topology are finite unions of sets of the form $S_1 \times S_2$, where at least one of the $S_i$ is finite. Given a polynomial $f$ in $(x,y)$, it follows that $V(f)$ on $\kk^2$ is closed in the product topology if and only if $f$ can be expressed as $g(x)h(y)$ for polynomials $g,h$ in one variable (\cref{exercise:product-top}). In particular the product topology on $\kk^2$ is {\em different} from the Zariski topology. 
\end{example}

A more natural topology on products of varieties is constructed as follows: consider projective spaces $\pp^m, \pp^n$ with homogeneous coordinates respectively $[x_0 : \cdots : x_m]$ and $[y_0: \cdots : y_n]$. A polynomial $f$ in $(x_0, \ldots, x_m, y_0, \ldots, y_n)$ is called \index{Bi-homogeneous polynomial}{\em bi-homogeneous} of \index{Bi-degree}{\em bi-degree} $(d,e)$ in the $x_i$ and $y_j$ if each monomial that appears in $f$ has degree $d$ in the $x$-variables and $e$ in the $y$-variables. The set of zeroes of a bi-homogeneous polynomial on $\pp^m \times \pp^n$ is well-defined, and the Zariski topology on $\pp^m \times \pp^n$ is by definition the (unique) topology whose basic closed subsets are intersections of zero sets of bi-homogeneous polynomials. If $V, W$ are quasiprojective varieties with their Zariski closures being subvarieties respectively of $\pp^m, \pp^n$, $j = 1, 2$, then the Zariski topology on $V \times W$ is the topology induced from the Zariski topology on $\pp^m \times \pp^n$. 

\begin{example} \label{example:Zariski-product-example-0}
If $X, Y$ are Zariski closed subsets of resepctively $\pp^m$ and $\pp^n$, then $X \times Y$ is Zariski closed in $\pp^m \times \pp^n$, since every homogeneous polynomial in $(x_0, \ldots, x_m)$ is trivially bi-homogeneous in  $(x_0, \ldots, x_m, \linebreak[1] y_0, \ldots, y_n)$. In particular, the Zariski topology on $\pp^m \times \pp^n$ is at least as fine as the product topology. 
\end{example}

\begin{example}\label{example:Zariski-affine-product}
The Zariski topology on $\kk^m \times \kk^n$ comes from identifying $\kk^m$ (respectively, $\kk^n$) with the basic open subset $\pp^m \setminus V(x_0)$ (respectively, $\pp^n \setminus V(y_0)$), so that $\kk^m \times \kk^n$ is identified with the open subset $U_{00} := (\pp^m \times \pp^n) \setminus V(x_0y_0)$ of $\pp^m \times \pp^n$. Write $x'_i := x_i/x_0$ and $y'_j := y_j/y_0$. It is straightforward to check (\cref{exercise:Zariski-affine-product}) that there is a one-to-one correspondence between the following collections of sets:
\begin{align}
\parbox{0.4\textwidth}{
intersections of $U_{00}$ and sets of zeroes on $\pp^m \times \pp^n$ of bi-homogeneous polynomials in $(x_0, \ldots, x_m, y_0, \ldots, y_n)$}
\quad
\longleftrightarrow
\quad
\parbox{0.35\textwidth}{
sets of zeroes on $\kk^{m + n}$ of polynomials in $(x'_1, \ldots, x'_m, y'_1, \ldots, y'_n)$}
\label{identification:Zariski-affine-product}
\end{align}
It follows that the Zariski topology on $\kk^m \times \kk^n$ is the same as the Zariski topology on $\kk^{m+n}$. Combined with the preceding examples, we see that the Zariski topology on $\kk^m \times \kk^n$ is finer than the product topology.
\end{example}

An equivalent formulation of the Zariski topology on product spaces can be given via the \index{Segre map}{\em Segre map} $s: \pp^m \times \pp^n \to \pp^N$, where $N := {(m+1)(n+1) -1} = mn + m + n$. It is defined as follows: let $z_{ij}$, $0 \leq i \leq m$, $0\leq j \leq n$, be a system of homogeneous coordinates on $\pp^N$; then
\begin{align}
s:([x_0 : \cdots : x_m], [y_0: \cdots : y_n]) \mapsto [x_iy_j]_{i,j}
\label{segre-defn}
\end{align}
where the right hand side denotes the point on $\pp^N$ with homogeneous coordinates $z_{ij} = x_iy_j$. The fundamental result in this section is the following result; we outline its proof in \cref{exercise:segre}. 

\begin{prop} \label{prop:segre}
The image of $\pp^m \times \pp^n$ under $s$ is the subvariety $Z$ of $\pp^N$ defined by (homogeneous) quadrics of the form $z_{ij}z_{kl} - z_{il}z_{kj}$, and $s$ induces a homeomorphism between $\pp^m \times \pp^n$ and $Z$. 
\end{prop}

\begin{example} \label{example:ruled}
When $m = n = 1$, $N = mn + m + n = 3$, and the image of $\pp^1 \times \pp^1 \into \pp^3$ is the hypersurface $V(x_0x_3 - x_1x_2)$ (we already encountered this variety in \cref{example:ruled-0}). It is called a \index{Ruled surface}{\em ruled surface}, since it comes with ``rulings'' given by lines of the form $\{a\} \times \pp^1$ and $\pp^1 \times \{b\}$ for $a, b \in \pp^1$. Under the usual identification of basic open subsets with affine spaces, the intersection of the ruled surface with $\kk^3$ is the hypersurface $V(z - xy)$ (\cref{exercise:basic-open-equation}) and the rulings are given by the lines $\{(a, t, at): t \in \kk\}$ and $\{(t,b, bt): t \in \kk\}$ for $a, b \in \kk$ (see \cref{fig:ruled}). 
\end{example}

\begin{center}

\begin{figure}[htb]
\newcommand\surfcolor{red!50!green}
\newcommand\lineonecolor{blue}
\newcommand\linetwocolor{red}
\newcommand\curvecolor{black}
\tikzstyle{sdot} = [yellow, circle, minimum size = 3pt, inner sep = 0pt, fill]
\pgfmathsetmacro\scalefactor{1}
\pgfmathsetmacro\axiscalefactor{1}
\pgfmathsetmacro\hview{60} 
\pgfmathsetmacro\vview{15} 
\pgfmathsetmacro\samplenum{50}
\pgfmathsetmacro\sampley{50}
\pgfmathsetmacro\opacity{0.3}
\pgfmathsetmacro\fillopacity{0.3}

\begin{tikzpicture}[scale = 0.75]

\pgfmathsetmacro\c{1}
\pgfmathsetmacro\xplus{1}
\pgfmathsetmacro\xminus{1}
\pgfmathsetmacro\yplus{1}
\pgfmathsetmacro\yminus{1}

\pgfmathsetmacro\cone{0}
\pgfmathsetmacro\ctwo{-0.75}
\pgfmathsetmacro\cthree{-0.25}

\pgfmathsetmacro\done{0}
\pgfmathsetmacro\dtwo{-0.75}
\pgfmathsetmacro\dthree{-0.25}

\begin{axis}[
    axis lines = middle,
    axis equal,
    view/h = \hview, view/v = \vview,
    scale = \axiscalefactor,
    xticklabels = {,,},
    yticklabels = {,,},
    zticklabels = {,,},
    x label style={anchor = north},
    y label style={anchor = north},
    z label style={anchor = south},
    xlabel = {\picfontsize $x$},
    ylabel = {\picfontsize $y$},
    zlabel = {\picfontsize $z$}
]
\addplot3 [
    surf,
    color = \surfcolor,
    opacity = \opacity,
    fill opacity = \fillopacity,
    faceted color = \surfcolor,
    samples = \samplenum,
    samples y = \sampley,
    domain = -\xminus : \xplus,
    y domain = -\yminus : \yplus
] (
	{x},
    {y},
   	{x*y*\c}
);

\addplot3 [
    samples = 2,
    samples y = 0,
    domain = -\xminus : \xplus,
    \lineonecolor, thick
] (
	{x},
	{\ctwo},
    {x*\ctwo*\c}
);
\addplot3 [
    samples = 2,
    samples y = 0,
    domain = -\xminus : \xplus,
    \lineonecolor, thick
] (
	{x},
	{\cthree},
    {x*\cthree*\c}
);

%
\addplot3 [
    samples = 2,
    samples y = 0,
    domain = -\yminus : \yplus,
    \linetwocolor, thick
] (
	{\dtwo},
	{x},
    {x*\dtwo*\c}
);
\addplot3 [
    samples = 2,
    samples y = 0,
    domain = -\yminus : \yplus,
    \linetwocolor, thick
] (
	{\dthree},
	{x},
    {x*\dthree*\c}
);

\addplot3 [
    samples = \samplenum,
    samples y = 0,
    domain = -\xminus : \xplus,
    \curvecolor, thick
] (
	{x},
	{x*x},
    {x*x*x*\c}
);


\end{axis}
\end{tikzpicture}
\caption{The rational normal curve and some rulings on the ruled surface $z = xy \subseteq \rr^3$ }
\label{fig:ruled}
\end{figure}
\end{center}

\Cref{prop:segre} implies that the Zariski topology on a product space agrees with the topology induced by the Zariski topology on the projective space upon identification of the product space with its image under the Segre map. From now on we view $\pp^m \times \pp^n$ as a projective variety by identifying it with the subvariety $s(\pp^m \times \pp^n)$ of $\pp^{(m+1)(n+1) - 1}$.

\subsection{Veronese embedding}
\label{veronese-subsection}

Let $d$ be a positive integer and $\scrV_d := \{\alpha = (\alpha_0, \ldots, \alpha_n) \in \zzeroo{n+1}: \sum_{j=0}^n \alpha_j = d\}$ be the set of exponents of monomials of degree $d$ in $(x_0, \ldots, x_n)$. The \index{Veronese map}{\em degree-$d$ Veronese map} $\nu_d: \pp^n \to \pp^{|\scrV_d|-1}$ is given by
\begin{align}
\nu_d &: [x_0: \cdots : x_n] \mapsto [ x^\alpha: \alpha \in \scrV_d] \label{veronese-defn}
\end{align}
where $x^\alpha$ is a shorthand for $x_0^{\alpha_0} \cdots x_n^{\alpha_n}$. 
\begin{prop} \label{veronese-embedding}
$\nu_d$ is a {\em closed embedding}, i.e.\ $\nu_d(\pp^n)$ is a subvariety of $\pp^{|\scrV_d|-1}$, and $\nu_d$ is an isomorphism between $\pp^n$ and $\nu_d(\pp^n)$.
\end{prop}

\begin{proof}
It is straightforward to see that $\nu_d$ can be expressed as a composition
\begin{align*}
\pp^n
	\overset{\delta}{\rightarrow}
		\pp^n \times \cdots \times \pp^n
	\overset{s_d}{\rightarrow}
		\pp^{(n+1)^d - 1}
	\overset{\pi}{\dashrightarrow}
		\pp^{|\scrV_d|-1}	
\end{align*}
where $\delta:x \mapsto (x, \ldots, x)$ is the diagonal map to the product of $d$ copies of $\pp^n$, $s_d$ is the Segre map on the $d$ factors, and $\pi$ is a projection which omits ``redundant'' coordinates of $s_d \circ \delta$ (i.e.\ for each $\alpha \in \scrV_d$, $\pi$ retains only one of the coordinates of $s_d \circ \delta$ equalling $x^\alpha$). Since $s_d$ and $\delta$ are closed embeddings (\cref{prop:segre,exercise:product-morphisms}), so is $s_d \circ \delta$. It is then straightforward to check that $\pi \circ s_d \circ \delta$ is a closed embedding as well.
\end{proof}

\begin{example} \label{example:rational-normal}
For $n = 1$, $\nu_d$ maps $[x_0: x_1] \mapsto [x_0^d: x_0^{d-1}x_1: \cdots : x_1^d] \in \pp^d$ and its image is called the \index{Rational!normal curve}{\em rational normal curve of degree $d$}. Under the identifications of basic open subsets of projective spaces with affine spaces, this becomes the map $\kk \to \kk^d$ given by $t \mapsto (t, t^2, \ldots, t^d)$ so that the rational normal curve of degree $d$ in $\kk^d$ is the subvariety defined by $x_2 = x_1^2, x_3 = x_1^3, \ldots, x_d = x_1^d$. Note that for $d = 3$, the rational curve lies on the ruled surface from \cref{example:ruled} (see \cref{fig:ruled}). 
\end{example}

We use \cref{veronese-embedding} to show that the complement of a projective hypersurface in a projective variety is affine. It is the projective analogue of the fact that the complement of a hypersurface in an affine variety is also affine (\cref{example:affine-complement}). 

\begin{prop} \label{projective-complement}
Let $X$ be a subvariety of $\pp^n$ and $f$ be a non-constant homogeneous polynomial in $(x_0, \ldots, x_n)$. Then $X \setminus V(f)$ is isomorphic to an affine variety.
\end{prop}

\begin{proof}
Let $d := \deg(f) \geq 1$. If $d = 1$, then after a suitable linear automorphism of $\pp^n$ from \cref{exercise:Pn-linear-auto}, we may assume that $f = x_0$. But then $X \setminus V(f)$ is a Zariski closed subset of the basic open set $U_0 := \pp^n \setminus V(x_0)$. Since $U_0 \cong \kk^n$, it follows that $X \setminus V(f)$ is isomorphic to an affine variety, as required. In the case that $d > 1$, let $Y := \nu_d(X) \subset \pp^{|\scrV_d|-1}$, then $\nu_d$ maps $X \setminus V(f)$ isomorphically onto $Y \setminus V(h)$ for a linear polynomial $h$ on $\nu_d$, and the result follows from the $d = 1$ case. 
\end{proof}

\subsection{Exercises}

\begin{exercise} \label{exercise:product-top}
Given a polynomial $f$ in $(x,y)$, show that the set $f(x,y) = 0$ on $\kk \times \kk$ is closed in the product topology if and only if $f$ can be expressed as $g(x)h(y)$ for polynomials $g,h$ in one variable. 
\end{exercise}

\begin{exercise} \label{exercise:Zariski-affine-product}
Show that the correspondence between the collections of sets in \eqref{identification:Zariski-affine-product} is bijective. 
\end{exercise}

\begin{exercise} \label{exercise:Zariski-affine-times-projective}
Show that the Zariski topology on $\kk^m \times \pp^n$ can be described as follows: let $(x_1, \ldots, x_m)$ be a system of coordinates  polynomial coordinates of $\kk^m$ and $[y_0: \cdots : y_n]$ be homogeneous coordinates on $\pp^n$. Then the closed sets on $\kk^m \times \pp^n$ are (finite) intersections of zero sets of polynomials $f \in \kk[x_1, \ldots, x_m, y_0, \ldots, y_n]$ which are {\em homogeneous in the $y$-variables}, i.e.\ of the form
\begin{align*}
f = \sum_{\beta_0 + \cdots + \beta_n = d} (\text{coefficient}) x_1^{\alpha_1} \cdots x_m^{\alpha_m}  y_0^{\beta_0} \cdots y_n^{\beta_n}
\end{align*}
for some $d \geq 0$ (in this case $d$ is said to be the {\em degree} of $f$ in $(y_0, \ldots, y_n)$).
\end{exercise}

\begin{exercise} \label{exercise:segre}
In this exercise you will prove \cref{prop:segre}.
\begin{enumerate}
\item Show that the image of $s$ is contained in $Z$. 
\item Show that to prove the first assertion of \cref{prop:segre} it suffices to prove that for each $i,j$, $s$ induces a homeomorphism between $s^{-1}(Z \setminus V(z_{ij}))$ and $Z \setminus V(z_{ij})$. 
\item It suffices to consider the case $i = j = 0$. Show that $s^{-1}(Z \setminus V(z_{00})) = (\pp^m \setminus V(x_0)) \times (\pp^n \setminus V(y_0))$. 
\item Write 
\begin{alignat*}{3}
& u_i &&:= z_{i0}/z_{00}, \quad && i = 1, \ldots, m, \\
& v_j &&:= z_{0j}/z_{00}, \quad && j = 1, \ldots, n, \\
& w_{ij} &&:= z_{ij}/z_{00}, \quad && i = 1, \ldots, m,\ j = 1, \ldots, n.
\end{alignat*}
Show that $Z \setminus V(z_{00})$ is isomorphic to the subvariety of $\kk^{m + n + mn}$ with coordinates $(u_1, \ldots, u_m,\linebreak[1] v_1, \ldots, v_n, \linebreak[1] w_{11}, \ldots, w_{mn})$ determined the equations $w_{ij} = u_iv_j$, $i = 1, \ldots, m$, $j = 1, \ldots, n$ [Hint: use \cref{exercise:basic-open-equation}.]. Deduce that $Z\setminus V(z_{00}) \cong \kk^{m+n}$. 
\item Write $x'_i := x_i/x_0$, $i = 1, \ldots, n$, $y'_j := y_j/y_0$, $j = 1, \ldots, n$. Show that the restriction of $s$ to $(\pp^m \setminus V(x_0)) \times (\pp^n \setminus V(y_0))$ can be described as the map from $\kk^m \times \kk^n \to \kk^{m + n + mn}$ given by 
\begin{alignat*}{3}
& u_i &&= x'_i, \quad && i = 1, \ldots, m, \\
& v_j &&= y'_j, \quad && j = 1, \ldots, n, \\
& w_{ij} &&= x'_iy'_j, \quad && i = 1, \ldots, m,\ j = 1, \ldots, n.
\end{alignat*}

\item Deduce that $s$ restricts to a homemorphism between $\kk^m \times \kk^n$ and $Z \setminus V(z_{00}) \cong \kk^{m+n}$. [Hint: use \cref{example:Zariski-affine-product}.] 
\end{enumerate}
\end{exercise}

\begin{exercise} \label{exercise:affine-product}
If $X$ is a subvariety of $\kk^m$ with coordinates $(x_1, \ldots, x_m)$ and $Y$ is a subvariety of $\kk^n$ with coordinates $(y_1, \ldots, y_n)$, show that $X \times Y$ is isomorphic to a subvariety of $\kk^{n+m}$ with coordinates $(x_1, \ldots, x_m, \linebreak[1] y_1, \ldots, y_n)$, and $I(X \times Y)$ in $\kk[x_1, \ldots, x_m, \linebreak[1] y_1, \ldots, y_n]$ is generated by $I(X) \subseteq \kk[x_1, \ldots, x_m]$ and $I(Y) \subseteq \kk[y_1, \ldots, y_n]$. [Hint: identify $\kk^m$ and $\kk^n$ with basic open subsets respectively of $\pp^m$ and $\pp^n$. Follow the steps of the isomorphism in \cref{exercise:segre}.] 
\end{exercise}

\begin{exercise} \label{exercise:product-projection}
Let $X, Y$ be quasiprojective varieties and $\pi_X: X \times Y \to X$ be the natural projection. 
\begin{enumerate}
\item In the case that $X$ is a quasiprojective subset of $\kk^m$ with coordinates $(x_1, \ldots, x_m)$, show that $x_i \circ \pi_X$ is a regular function on $X \times Y$. [Hint: reduce to the case that $Y$ is affine. Use \cref{exercise:affine-product}.]
\item Deduce that $\pi_X$ is a morphism. [Hint: morphisms can be checked locally (\cref{exercise:check-morphism}). Use the preceding assertion and \cref{prop:morphism-on-Uj}]. 
\end{enumerate}
\end{exercise}

\begin{exercise}\label{exercise:product-morphisms}
Let $\phi: X \to Y$ be a morphism of varieties. 
\begin{enumerate}
\item \label{product-morphisms:graph} Show that the {\em graph} of $\phi$, i.e.\ the set $\graph(\phi) := \{(x,\phi(x)): x \in X\}$ is a subvariety (i.e.\ a Zariski closed subset) of $X \times Y$ and it is isomorphic to $X$. [Hint: use \cref{exercise:product-projection} to express $\graph(\phi)$ as the set of zeroes of regular functions on $X \times Y$.]
\item  \label{product-morphisms:m1} If $\psi:X \to Z$ is another morphism of varieties, then show that the map $X \to Y \times Z$ given by $x \mapsto (\phi(x), \psi(x))$ is a morphism. [Hint: reduce to the case that $Y, Z$ are affine. Use \cref{exercise:affine-product}.]
\item \ \label{product-morphisms:m2} Deduce that if $\phi':X' \to Y'$ is another morphism of varieties, then the map $X \times X' \to Y \times Y'$ given by $(x,x') \mapsto (\phi(x), \phi'(x'))$ is a morphism. [Hint: $\phi \circ \pi_X$ and $\phi' \circ \pi_{X'}$ are morphisms from $X \times X'$ to respectively $Y$ and $Y'$.]
\end{enumerate}
\end{exercise}


\section{Completeness and compactification} \label{complete-section}
While studying manifolds, it is often necessary to compactify them in order that the intersection theory of submanifolds is well behaved. However, the Noetherianness of finitely generated $\kk$-algebras (Hilbert's basis theorem) implies that {\em all} varieties satisfy the usual definitions of compactness (\cref{exercise:noetherian-to-compact}). The property which plays in the case of varieties the role similar to that of compactness in the case of manifolds is {\em completeness}. A subset of a variety $X$ is called \index{Variety!complete}\index{Complete!variety}{\em complete} if for every variety $Y$, the projection map $\pi: X \times Y \to Y$ is closed, i.e.\ it maps closed sets on to closed sets. 

\begin{example} 
If $n \geq 1$, the projection onto the $x_{n+1}$-coordinate maps the hypersurface $V(x_1x_{n+1} - 1)$ of $\kk^{n+1} = \kk^n \times \kk$ to $\kk \setminus \{0\}$, which is not closed in $\kk$. It follows that $\kk^n$ is {\em not} complete for any $n \geq 1$. 
\end{example}

\begin{example} \label{example:regular-on-complete}
We can push the reasoning in the preceding example a bit further. If $f$ is a regular function on a complete variety $X$, consider the subset $Z := V(ft- 1) \subseteq X \times \kk$, where $t$ is the coordinate on $\kk$. Since $Z$ closed in $X \times \kk$, its projection $\pi(Z)$ is closed in $\kk$. Since $\pi(Z) \subsetneq \kk$, it must be finite, i.e.\ $f$ takes {\em finitely many} values. This implies, e.g., that the only complete subsets of affine varieties are finite sets of points. See \cref{exercise:regular-on-complete} for some other immediate consequences.
\end{example}

For many ``natural'' topological spaces (e.g.\ the Euclidean topology), compactness is equivalent to completeness - this is discussed in \crefrange{exercise:first-countable}{exercise:compact=complete}. The following properties of complete sets are straightforward to verify and left as exercises (\cref{exercise:closed-sub-complete}): 

\begin{prop}\label{prop:closed-sub-complete}
\begin{enumerate}
\item A complete subset of a variety is also Zariski closed.
\item Every subvariety of a complete variety is complete.
\item The image of a complete variety under a morphism is also complete. \qed
\end{enumerate}
\end{prop}

The following is the main result of this section: 

\begin{thm} \label{thm:Pn-complete-0}
The projective space $\pp^n$ is complete for all $n$. 
\end{thm}

We give a proof of \cref{thm:Pn-complete-0} following \cite[Proof of Theorem 2.23]{mummetry}. In fact we prove \cref{thm:Pn-complete-1} below; the equivalence of this result with \cref{thm:Pn-complete-0} is left as \cref{exercise:Pn-complete-0=1}. 

\begin{prop} \label{thm:Pn-complete-1}
The projection map $\pi: \pp^n \times \kk^m \to \kk^m$ is closed for each $m$. 
\end{prop}

\begin{proof} 
Choose affine coordinates $(y_1, \ldots, y_m)$ on $\kk^m$ and homogeneous coordinates $[x_0: \cdots: x_n]$ on $\pp^n$. Let $Z$ be a closed subset of  $\pp^n \times \kk^m$. Then $Z = V(f_1, \ldots, f_k)$ for polynomials $f_j (x,y) \in \kk[x_0, \ldots, x_n, \linebreak[1], y_1, \ldots, y_m]$ which are homogeneous of degree $d_j \geq 0$ in the $x$-variables (\cref{exercise:Zariski-affine-times-projective}). For each $b \in \kk^m$, it follows from the correspondence between subvarieties of $\pp^n$ and homogeneous ideals of $\kk[x_0, \ldots, x_n]$ (\cref{exercise:projective-nullspondence}) that $y \not\in \pi(Z)$ if and only if there is $d \geq 1$ such that $\mmm_+^d$ is contained in the ideal of $\kk[x_0, \ldots, x_n]$ generated by $f_1(x,b), \ldots, f_k(x,b)$, where $\mmm_+$ is the {\em irrelevant} maximal ideal of $\kk[x_0, \ldots, x_n]$ generated by $x_0, \ldots, x_n$. Therefore it suffices to show that for each $d \geq 1$, the set $\{b \in \kk^m:  \mmm_+^d$ is contained in the ideal of $\kk[x_0, \ldots, x_n]$ generated by $f_1(x,b), \ldots, f_k(x,b)\}$ is Zariski open in $\kk^m$. For each integer $e$, let $V_e$ be the vector space of homogeneous polynomials of degree $e$ in $(x_0, \ldots, x_n)$ (with $V_e = 0$ if $e < 0$) and let $m_e := \dim_\kk(V_e)$. For each $b \in \kk^m$, and $d \geq 1$, consider the linear map 
\begin{align*}
T^d(b): V_{d - d_1} \dsum \cdots \dsum V_{d-d_k} \to V_d
	\qquad 
	(g_1(x), \ldots, g_k(x)) \mapsto \sum_{j=1}^k f_j(x,b)g_j(x) 
\end{align*}
Let $[T^d(b)]$ be the matrix of $T^d(b)$ with respect to fixed bases of $\dsum_j V_{d-d_j}$ and $V_d$. It is straightforward to check that $T^d(b)$ is surjective if and only if there is an $m_d \times m_d$-minor of $[T^d(b)]$ with nonzero determinant. It follows that the set $\{b \in \kk^m:  T^d(b)$ is surjective$\}$ is open in $\kk^m$, as required. 
\end{proof}

A \index{Compactification!of a variety}{\em compactification}\footnote{We would have liked to use ``completion'' instead of ``compactification''; however, it might have suggested a (misleading) connection with the notion of completion of local rings (discussed in \cref{completion-section}).} of a variety $X$ is a complete variety $\bar X$ containing an dense open subset isomorphic to $X$. \Cref{thm:Pn-complete-0} implies that every quasiprojective variety has a compactification:

\begin{prop} \label{prop:proj-complete}
\begin{enumerate}
\item Every projective variety is complete. In particular, if $X$ is a quasiprojective subset of $\pp^n$, then the closure $\bar X$ of $X$ in $\pp^n$ is a compactification of $X$. 
\item Let  $\phi: X \to Y$ be a morphism of varieties such that $\phi(X)$ is dense in $Y$. Then $\phi$ can be extended to a {\em surjective} morphism $\phi': X' \to Y$ where $X'$ is a variety which contains (an isomorphic copy of) $X$ as a dense open subset. 
\end{enumerate}
\end{prop}

\begin{proof}
The first assertion is immediate from \cref{prop:closed-sub-complete,thm:Pn-complete-0}. For the second assertion, let $\bar X$ be a compactification of $X$, and take $X'$ to be the closure in $\bar X \times Y$ of the graph $\graph(\phi) := \{(x, \phi(x)): x \in X\}$ of $\phi$, and $\phi'$ to be the natural projection from $X'$ to $Y$. The completeness of $\bar X$ implies that $\phi'$ is surjective. It remains to prove that $\graph(\phi)$ is open in $X'$. Indeed, since $\graph(\phi)$ is already closed in $X \times Y$, it follows that $\graph(\phi) = X' \cap (X \times Y)$. Since $X \times Y$ is open in $\bar X \times Y$, it follows that $\graph(\phi)$ is open in $X'$, as required. 
\end{proof}

\subsection{Exercises}

\begin{exercise} \label{exercise:regular-on-complete}
Use \cref{example:regular-on-complete} to prove  that
\begin{enumerate}
\item Every regular function on a complete irreducible variety is constant, i.e.\ $\kk[X] = \kk$.
\item If $f$ is a morphism from an irreducible complete variety to the affine space, then the image of $f$ is a point. 
\item If $\bar X$ is a compactification of an affine variety $X$ and $C$ is a complete subset of $\bar X$ consisting of infinitely many points, then $C \cap (\bar X \setminus X) \neq \emptyset$.
\end{enumerate}
\end{exercise}

\begin{exercise}\label{exercise:closed-sub-complete}
Show that
\begin{enumerate}
\item A complete subset of a variety is also (Zariski) closed. [Hint: Given $Z \subseteq X$, the diagonal subset $Z' := \{(z,z): z \in Z\}$ is closed in $Z \times X$ (\cref{exercise:product-morphisms}).] 
\item Every subvariety of a complete variety is complete. [Hint: if $Z$ is a subvariety (i.e.\ Zariski closed subset) of $X$, then $Z \times Y$ is closed in $X \times Y$.] 
\item The image of a complete variety under a morphism is also complete. [Hint: for every variety $Y$, every morphism $\phi: Z \to Z'$ of varieties lifts to a morphism $\phi': Z \times Y \to Z' \times Y$ (\cref{exercise:product-morphisms}).]
\end{enumerate}
\end{exercise}

\begin{exercise} \label{exercise:Pn-complete-0=1}
Show that in order to prove that $\pp^n$ is complete, it suffices to prove that the projection map $\pi: \pp^n \times \kk^m \to \kk^m$ is closed for each $m$. [Hint: every variety has a finite open covering by closed subsets of the affine space (\cref{affine-cover}).]
\end{exercise}

\begin{exercise} \label{exercise:morphism-proj-to-affine} 
Show that a morphism from an irreducible projective variety to an affine variety must be a constant map. [Hint: use \cref{prop:closed-sub-complete,prop:proj-complete}.]
\end{exercise}

\begin{exercise} \label{exercise:noetherian-to-compact}
Recall that a topological space is called {\em compact} if each of its open covers has a finite subcover, and it is called {\em sequentially compact} if every infinite sequence of points has a convergent subsequence. Show that
\begin{enumerate}
\item Every variety is compact. [Hint: use \cref{exercise:Zariski-quirk:compact,affine-cover}.]
\item Every variety is sequentially compact. [Hint (the following strategy requires the notion of {\em dimension} covered in \cref{dimension-section}): start with a sequence of points on a quasiprojective variety $X$ and take its closure $Z$ in $X$. Reduce to the case that $Z$ is irreducible and affine. Then consider two cases: (i) there is a nonzero regular function on $Z$ that vanishes on an infinite subsequence, and (ii) every nonzero regular function vanishes on only finitely many points from that sequence.]
\end{enumerate}
\end{exercise}

\begin{exercise} \label{exercise:first-countable}
We say that a topological space $X$ is {\em first-countable} if for every point $x \in X$, there is a sequence of open neighborhoods $\{U_j\}_j$ of $x$ in $X$ such that every open neighborhood $U$ of $x$ in $X$ contains some $U_j$. Show that the following are equivalent:
\begin{enumerate}
\item $X$ is first countable.
\item \label{strongly-first-countable} For every point $x \in X$, there is a sequence of open neighborhoods $\{U_j\}_j$ of $x$ in $X$ such that for every open neighborhood $U$ of $x$ in $X$, there is $N$ such that $U$ contains $U_j$ for each $j \geq N$.
\end{enumerate}
\end{exercise}

\begin{exercise} \label{exercise:compact=complete}
Let $X$ be a first countable topological space. Consider the following properties:
\begin{enumerate}
\item \label{seq-compactness} $X$ is sequentially compact.
\item \label{completeness-general} For every first countable topological space $Y$, the projection map $X \times Y \to Y$ is closed with respect to the product topology on $X \times Y$.
\end{enumerate}
Show that \eqref{seq-compactness} implies \eqref{completeness-general}. If $X$ is a first countable $T_1$-space\footnote{A topological space is $T_1$ if for every pair of distinct points, each has an open neighborhood not containing the other.}, show that \eqref{completeness-general} implies \eqref{seq-compactness}.
\end{exercise}

\section{Image of a morphism: Part I}
 Since projective varieties are complete (\cref{prop:proj-complete}), the image of a morphism from a projective variety must be a subvariety of the target space, and in addition must be a {\em finite set} if the target space is an affine variety. The image of a morphism from an affine variety can be wilder - we have seen that it does not have to be closed or open (\cref{example:constructible-0}); the following example\footnote{It was motivated by a comment by R.\ Borcherds to the answer \cite{mathoverflow-algeom-examples-aff-to-proj} by A.\ Bayer on {\em MathOverflow}.} shows that it can be a non-trivial projective variety: 
 
 \begin{example} \label{example:aff-to-proj}
 Consider the morphism $\phi:\kk \to \pp^1$ given by $x \mapsto [(x-a_1)(x-a_2): (x-b_1)(x-b_2)]$, where $a_i, b_j \in \kk$ such that $a_i \neq b_j$ for any $i, j = 1, 2$. Then $\phi$ is surjective if (and only if) $a_1 + a_2 \neq b_1 + b_2$ (\cref{exercise:aff-to-proj}).
 \end{example}

A morphism $\phi: X \to Y$ of varieties is said to be \index{Dominant!morphism}\index{Morphism!dominant}{\em dominant} or {\em dominating} if $\phi(X)$ is Zariski dense in $Y$, or equivalently, if $\phi(X)$ is not contained in any proper subvariety of $Y$. Any morphism $\phi: X \to Y$ turns into a dominant map if the target space is changed from $Y$ to the Zariski closure of $\phi(X)$ in $Y$. 

\begin{prop} \label{prop:star-injective}
If a morphism $\phi: X \to Y$ of quasiprojective varieties is dominant (in particular, if it is surjective), then the pullback $\phi^*: \kk[X] \to \kk[Y]$ is injective.
\end{prop}

\begin{proof}
Given a nonzero regular function $f$ on $Y$, the set $Y \setminus V(f)$ is nonempty and, due to \cref{prop:regular-zeroes}, Zariski open in $Y$. Since $\phi$ is dominant, $\phi(X)$ intersects $Y \setminus V(f)$, so that $U := \phi^{-1}(Y \setminus V(f))$ is a nonempty Zariski open subset of $X$. Since $\phi^*(f)$ is nonzero on $U$, the proposition follows.
\end{proof}

\begin{example}
The converse of \cref{prop:star-injective} is {\em not} true. Indeed, since $\kk[\pp^n] = \kk$, if $\phi: \pp^n \to \pp^n$ is a constant morphism that maps every point to a fixed point on $\pp^n$, then $\phi^*$ is injective even though $\phi$ is far from being dominant if $n \geq 1$. 
\end{example}

Let $X$ be an affine variety and $\qqq$ be an ideal of $\kk[X]$. Given generators $f_1, \ldots, f_N$ of $\qqq$, the \index{Blow up}{\em blow up $\bl_\qqq(X)$ of $X$ at $\qqq$} is the closure in $X \times \pp^{N-1}$ of the graph of the map $X \setminus V(\qqq) \to \pp^{N-1}$ given by $x \mapsto [f_1(x) : \cdots : f_N(x)]$. The following result shows that $\bl_\qqq(X)$ depends only on $\qqq$, not on the choice of the $f_j$.

\begin{prop} \label{prop:blow-up-well-defined}
Let $g_1, \ldots, g_q$ (respectively, $h_1, \ldots, h_r$) be  generators of $\qqq$ in $\kk[X]$. Write $Y$ (respectively, $Z$) for the closure in $X \times \pp^{q-1}$ (respectively, $X \times \pp^{r-1}$) of the graph of the map $\phi_g: x \mapsto [g_1(x) : \cdots :g_q(x)]$ (respectively, $\phi_h: x \mapsto [h_1(x) : \cdots : h_r(x)]$). Then $Y \cong Z$.
\end{prop}

\begin{proof}
Here we only consider the special case that $r= q+1$ and $h_j = g_j$, $j = 1, \ldots, q$; the general case is left as \cref{exercise:blow-up-well-defined}. There are $g'_1, \ldots, g'_q \in \kk[X]$ such that $h_{q+1} = g'_1g_1 + \cdots + g'_qg_q$. Let $\psi: X  \times \pp^{q-1} \to X \times \pp^q$ be the map that sends
\begin{align*}
(x,[z_1: \cdots : z_q])
\mapsto
(x, [z_1: \cdots : z_q : g'_1(x)z_1 + \cdots + g'_q(x)z_q])
\end{align*}
It is straightforward to check that $\psi$ is a morphism and, if $\pi: X  \times \pp^q \dashrightarrow X \times \pp^{q-1}$ is the natural projection given by $(x, [z_1: \cdots : z_{q+1}]) \mapsto (x, [z_1: \cdots : z_q])$, then $\pi \circ \psi$ is identity on $\graph(\phi_g)$, and therefore it is identity everywhere on $Y$ (\cref{exercise:morphism-dense=equal}). On the other hand, $\graph(\phi_h) \subseteq V(z_{q+1} - \sum_{j=1}^q g'_j(x)z_j)$, since the latter is closed in $X \times  \pp^q$, it follows that $Z \subseteq V(z_{q+1} - \sum_{j=1}^q g'_j(x)z_j)$. This implies that $\pi$ is well-defined on $Z$. Since $\psi \circ \pi$ is identity on $\graph(\phi_h)$, it then follows that $\psi \circ \pi$ is identity on $Z$. Consequently, $\phi$ induces an isomorphism $Y \cong Z$.
\end{proof}

We next study the restriction to $\bl_{\qqq}(X)$ of the natural projection from $X \times \pp^{N-1}$ to $X$; often this map is what one means by the ``blow up of $X$ at $\qqq$.''  The next result shows that $\sigma$ induces an isomorphism on $X' := X \setminus V(\qqq)$ and its image is the closure $\bar X'$ of $X'$ in $X$. (Note that if $X$ is irreducible and $\qqq$ is not the zero ideal, then $\bar X' = X$.)

\begin{prop} \label{prop:blow-up-isomorphism}
Let $\qqq$ be an ideal of the coordinate ring of an affine variety $X$ and $\sigma: \bl_\qqq(X) \to X$ be the blow up of $X$ at $\qqq$. Let $
\bar X'$ be the closure of $X' := X \setminus V(\qqq)$ in $X$. Then 
\begin{enumerate}
\item $\sigma$ maps $\bl_\qqq(X)$ onto $\bar X'$, and
\item $\sigma$ induces an isomorphism $\bl_\qqq(X) \setminus \sigma^{-1}(V(\qqq)) \cong X \setminus V(\qqq) = \bar X' \setminus V(\qqq) = X'$. In particular, if $X$ is irreducible and $\qqq$ is a nonzero ideal, then $\bl_\qqq(X)$ is irreducible and birational to $X$.
\end{enumerate}
\end{prop}

\begin{proof}
Since $\pp^{N-1}$ is complete and $\bl_\qqq(X)$ is closed in $\bar X' \times \pp^{N-1}$, it follows that its projection to $\bar X'$ is closed. Since $\sigma(\bl_\qqq(X))$ clearly contains $X'$, it follows that $\sigma(\bl_\qqq(X)) = \bar X'$, as required for the first assertion. The second assertion follows directly from assertion \eqref{product-morphisms:graph}  of \cref{exercise:product-morphisms}. 
\end{proof}

If $V$ is an irreducible subvariety of $X$ not contained in $V(\qqq) \subset X$, then the \index{Strict transform}{\em strict transform} of $V$ on $\bl_\qqq(X)$ is the closure in $\bl_\qqq(X)$ of $\sigma^{-1}(V \setminus V(\qqq))$. In the case that $\qqq$ is the ideal $I(Y)$ of $\kk[X]$ consisting of regular functions vanishing on a subvariety $Y$ of $X$, we also write $\bl_Y(X)$ for $\bl_\qqq(X)$, and call it the ``blow up of $Y$ in $X$.'' 

\begin{example}
If $Y$ is a point in $X := \kk^n$, then $\sigma: \bl_Y(X) \to X$ is an isomophism away from $Y$, and $\sigma^{-1}(Y) \cong \pp^{n-1}$ (\cref{exercise:blow-up-origin}). In particular, for $n \geq 2$, the point is ``blown up'' to a hyperplane in $\pp^n$. This is probably the origin of the name ``blow up''  (see \cref{fig:blow-up}). In \cref{weighted-blow-up-section} we revisit the blow up at a point of $\kk^n$ from the perspective of ``toric varieties'' as a special case of ``weighted blow ups'' at a point. 
\end{example}

\begin{center}
\begin{figure}[htb]
\begin{tikzpicture}
\pgfmathsetmacro\scalefactor{1}
\pgfmathsetmacro\bottomarrowstartx{1.71*\scalefactor}
\pgfmathsetmacro\bottomarrowstarty{0.2*\scalefactor}
\pgfmathsetmacro\bottomarrowendy{-0.4*\scalefactor}
\pgfmathsetmacro\bottomfigdist{2.5*\scalefactor}

\pgfmathsetmacro\rightarrowstartx{4*\scalefactor}
\pgfmathsetmacro\rightarrowstarty{1.5*\scalefactor}
\pgfmathsetmacro\rightarrowendx{6.5*\scalefactor}
\pgfmathsetmacro\rightfigdist{7.25*\scalefactor}
\pgfmathsetmacro\rightfigtopdist{0.4*\scalefactor}

\pgfmathsetmacro\diagarrowstartx{\rightarrowendx}
\pgfmathsetmacro\diagarrowstarty{\bottomarrowstarty}
\pgfmathsetmacro\diagarrowendx{\rightarrowstartx}
\pgfmathsetmacro\diagarrowendy{\bottomarrowendy}

\pgfmathsetmacro\axiscalefactor{0.5}
	
\begin{scope}
				
\begin{axis}[
    hide axis,
    view={0}{30},
    scale=\axiscalefactor
]

\addplot3 [
    samples=40,
    samples y=0,
	domain=150:-30,
	red, ultra thick
] (
    {-cos(-x)},
    {-sin(-x)},
    {(-(x)^3)/90}
);

\addplot3 [
    surf, shader=faceted interp,
    point meta=x,
    samples=40,
    samples y=5,
    z buffer=sort,
    domain=150:-30,
    y domain=-1:1
] (
    {y*cos(-x)},
    {y*sin(-x)},
    {(-(x)^3)/90});

\addplot3 [
    samples=50,
    domain=150:135, 
    samples y=0,
    green, ultra thick
] (
    {0},
    {0},
    {-x^3/90});

\addplot3 [
    samples=50,
    domain=135:100, 
    samples y=0,
    green, ultra thick,
    dashed
] (
    {0},
    {0},
    {-x^3/90});

\addplot3 [
    samples=50,
    domain=90:-30, 
    samples y=0,
    green, ultra thick
] (
    {0},
    {0},
    {-x^3/90});

\addplot3 [
    samples=40,
    samples y=0,
	domain=-1:1,
] (
    {x*cos(30)},
    {x*sin(30)},
    {-150^3/90});

\addplot3 [
    samples=40,
    samples y=0,
	domain=-1:1,
] (
    {x*cos(30)},
    {x*sin(30)},
    {30^3/90});

\addplot3 [
    samples=40,
    samples y=0,
	domain=150:-30,
	red, ultra thick
] (
    {cos(-x)},
    {sin(-x)},
    {(-(x)^3)/90}
);
\end{axis}

\begin{scope}[shift={(0,-\bottomfigdist)}]

\begin{axis}[
    hide axis,
    view={0}{30},
    scale=\axiscalefactor
]

\addplot3 [
    surf, shader=faceted interp,
    point meta=x,
    samples=40,
    samples y=5,
    z buffer=sort,
    domain=150:-30,
    y domain=-1:1
] (
    {y*cos(-x)},
    {y*sin(-x)},
    {-190^3/90});

\addplot3 [
    samples=2,
    samples y=0,
	domain=-1:1,
] (
    {x*cos(30)},
    {x*sin(30)},
    {-190^3/90});

\addplot3 [
    samples=40,
    samples y=0,
   domain=-180:180,
   red, ultra thick
] (
    {cos(x)},
    {sin(x)},
    {-190^3/90});

\addplot3 [
    surf, shader=faceted interp,
    point meta=x,
    colormap = {greenmap}{rgb(0cm)=(0,1,0); rgb(1cm)=(0,1,0)},
    samples=40,
    samples y=2,
    z buffer=sort,
    domain=150:-30,
    y domain=-0.05:0.05
] (
    {y*cos(-x)},
    {y*sin(-x)},
    {-190^3/90});

\end{axis}
\end{scope}

\begin{scope}[shift={(\rightfigdist,\rightfigtopdist)}]
\begin{axis}[
    hide axis,
    view={70}{-10},
    scale = \axiscalefactor
]
\addplot3 [
    surf, shader=faceted interp,
    point meta=x,
    samples=40,
    samples y=5,
    z buffer=sort,
    domain=0:360,
    y domain=-0.5:0.5
] (
	{0.5*y*sin(x/2)},
	{(1+0.5*y*cos(x/2))*sin(x)},
    {(1+0.5*y*cos(x/2))*cos(x)}
    );

\addplot3 [
    samples=50,
    domain=130:190, 
    samples y=0,
    dashed,
    green, ultra thick
] (
	{0},
    {sin(x)},
    {cos(x)}
    );

\addplot3 [
    samples=50,
    samples y=0,
    domain=0:360,
    red, ultra thick
] (
	{0.25*sin(x/2)},
	{(1+0.25*cos(x/2))*sin(x)},
    {(1+0.25*cos(x/2))*cos(x)}
    );

\addplot3 [
    samples=50,
    samples y=0,
    domain=0:123,
    red, ultra thick
] (
	{-0.25*sin(x/2)},
	{(1-0.25*cos(x/2))*sin(x)},
    {(1-0.25*cos(x/2))*cos(x)}
    );


\addplot3 [
    samples=50,
    samples y=0,
    domain=145:360,
    red, ultra thick
] (
	{-0.25*sin(x/2)},
	{(1-0.25*cos(x/2))*sin(x)},
    {(1-0.25*cos(x/2))*cos(x)}
    );

\addplot3 [
    samples=50,
    domain=-190:130, 
    samples y=0,
    green, ultra thick
] (
	{0},
    {sin(x)},
    {cos(x)}
    );
    		
\end{axis}
\end{scope}

\draw [dashed, ->] (\bottomarrowstartx,\bottomarrowendy) -- (\bottomarrowstartx,\bottomarrowstarty);

\draw [->] (\rightarrowstartx,\rightarrowstarty) -- (\rightarrowendx, \rightarrowstarty) node [above, pos=0.5] {\scriptsize{identify two ends}};

\draw [->] (\diagarrowstartx,\diagarrowstarty) -- (\diagarrowendx, \diagarrowendy) node [above, pos=0.5] {\scriptsize{blow up map}};
\end{scope}

\end{tikzpicture}
\caption{Blow up of a point on a disk in $\rr^2$ is a {\em M\"obius strip}}
\label{fig:blow-up}
\end{figure}
\end{center}

\begin{example}  \label{example:linear-blow-up}
More generally, if $Y$ is a ``linear subspace of dimension $m$'' in $X := \kk^n$, then $\sigma: \bl_Y(X) \to X$ is an isomorphism away from $Y$ and $\sigma^{-1}(y) \cong \pp^{n-m - 1}$ for each $y \in Y$ (\cref{exercise:linear-blow-up}).
\end{example}

\begin{example} \label{example:projective-blow-up}
If $X$ is an open affine subset of $\bar X$ and $Z := V(\qqq) \subseteq X$ consists of finitely many points, then the blow up map $\sigma: \bl_\qqq(X) \to X$ canonically extends to a map $\bar \sigma: \bar Y \to \bar X$ such that
\begin{prooflist}
\item $\bar Y$ contains $Y := \bl_\qqq(X)$ as an open subset, 
\item $\bar \sigma|_Y = \sigma|_Y$, 
\item $\bar \sigma$ induces an isomorphism $\bar Y \setminus E \cong \bar X \setminus Z$, where $E := \sigma^{-1}(Z)$ is the ``exceptional divisor'' of $\sigma$. 
\end{prooflist}
Indeed, since $\bl_\qqq(X)$ is a closed subset of $X \times \pp^N$ for some $N$, and $\sigma$ is the projection onto $X$, it is straightforward to check that we can take $\bar Y$ to be the closure of $\bl_\qqq(X)$ in $\bar X \times \pp^N$ and $\bar \sigma$ to be the projection onto $\bar X$. In particular, this gives a construction of the blow up of $\pp^n$ at a point from the blow up of $\kk^n$ at a point.  
\end{example}

The next two results describe coordinate rings of closures of morphisms from affine varieties. These results are used in \cref{toric-intro} to compute coordinate rings of  ``toric varieties.''

\begin{prop} \label{affine-embedding}
Let $\phi: X \to \kk^N$ be the morphism from an affine variety $X$ given by $x \mapsto (f_1(x), \linebreak[1] \ldots, \linebreak[1] f_N(x))$, where $f_1, \ldots, f_N$ are regular functions $X$. Let $Z$ be the closure of $\phi(X)$ in $\kk^N$ and $R$ be the $\kk$-subalgebra of $\kk[X]$ generated by $f_1, \ldots, f_N$.
\begin{enumerate}
\item $\kk[Z] \cong R$.
\item \label{subspace-assertion} Let $\pi: \kk^N \to \kk^M$ be the natural projection onto the first $M$ coordinates, where $M \leq N$, and let $H$ be the coordinate subspace of $\kk^N$ spanned by the first $M$ coordinates. Then $Z \cap H$ is contained in the closure of $\pi \circ \phi(X)$.
\item If $\phi(X) \subseteq \nktoruss{N}$ (where $\kk^* := \kk \setminus \{0\}$), then the coordinate ring of the closure of $\phi(X)$ in $\nktoruss{N}$ is isomorphic to $R_{f_1\cdots f_N} := \kk[f_1,  f_1^{-1}, \ldots, f_N, f_N^{-1}]$.
\end{enumerate}
\end{prop}

\begin{proof}
Let $\phi^*: \kk[y_1, \ldots, y_N] \to R$ be the ring homomorphism which maps each $y_i \mapsto f_i$. Let $I(Z)$ be the ideal of $\kk[y_1, \ldots, y_N]$ consisting of polynomials that vanish on $Z$. Then $g \in I(Z)$ if and only if $g$ vanishes on $\phi(X)$ if and only if $g(f_1, \ldots, f_N)(x)$ is zero for all $x \in X$ if and only if $\phi^*(g) = 0$ in $\kk[X]$. Therefore $I(Z) = \ker \phi^*$, and $R \cong \kk[y_1, \ldots, y_N]/I(Z) \cong \kk[Z]$ (see \cref{eqn:K[X]-affine} for the last isomorphism). The second assertion follows from a general property of retractions - see \cref{exercise:retraction-intersection}.
Finally, since $\nktoruss{N} = \kk^n \setminus V(y_1 \cdots y_N)$, the closure of $\phi(X)$ in $\nktoruss{N}$ is $Z\setminus V(y_1 \cdots y_N)$ (\cref{exercise:closure-in-open-subset}), so that the coordinate ring of $Z\setminus V(y_1 \cdots y_N)$ is $\kk[Z]_{y_1 \cdots y_N}$ (\cref{example:affine-complement}), which by the first assertion is isomorphic to $R_{f_1 \cdots f_N}$, as required for the third assertion.
\end{proof}

\begin{cor} \label{projective-embedding}
Let $f_0, \ldots, f_N$ be regular functions on an {\em irreducible} affine variety $X$ such that no $f_i$ is identically zero on $X$. Let $Z$ be the closure in $\pp^N$ of the image of the map $\phi: X \setminus V(f_0, \ldots, f_N) \to \pp^N$ defined by $x \mapsto [f_0(x): \cdots : f_N(x)]$. Let $[z_0: \cdots : z_N]$ be homogeneous coordinates on $\pp^N$ and $U_i := \pp^N\setminus V(z_i)$, $0 \leq i \leq N$. Then
\begin{enumerate}
\item \label{projective-i} $\kk[Z \cap U_i] \cong \kk[f_0/f_i, \ldots, f_N/f_i]$,
\item \label{projective-ij} $\kk[Z \cap U_i \cap U_j] \cong \kk[f_0/f_i, \ldots, f_N/f_i, f_i/f_j]$.
\item \label{rational-combination} Assume there exists $b \in \zplus$ and $b_1, \ldots, b_N \in \zzero$ such that $b = \sum_{i=1}^n b_i$ and $f_0^{b} = \prod_{i \geq 1} f_i^{b_i}$. Then $Z \cap U_0 \subseteq \bigcup_{i \geq 1} (Z \cap U_i)$.
\item \label{subspace-intersection}  Let $H = V(z_{M+1}, \ldots, z_N)$, $H' := V(z_0, \ldots, z_M)$ for some $M \leq N$, and $\pi: \pp^N\setminus H' \to H$ be the natural projection onto the first $M$ coordinates. Then $Z \cap H$ is contained in the closure of $\pi \circ \phi(X \setminus V(f_0, \ldots, f_M))$.
\end{enumerate}
\end{cor}

\begin{proof}
Since each $U_i \cong \kk^N$, assertion \eqref{projective-i} follows from \cref{affine-embedding}. Write $Z_i := Z \cap U_i$ and $Z_{ij} := Z \cap U_i \cap U_j$. Since $Z_i$ is an affine variety, $f_j/f_i$ is a regular function on $Z_i$, and $Z_{ij} = Z_{i} \setminus V(f_j/f_i)$, it follows that $\kk[Z_{ij}] = \kk[Z_i]_{f_j/f_i} =  \kk[f_0/f_i, \ldots, f_N/f_i, f_i/f_j]$ (see \cref{example:affine-complement} for the first equality), which proves assertion \eqref{projective-ij}. For assertion \eqref{rational-combination}, note that $z_0^b - \prod_{i \geq 1} z_i^{b_i}$ is identically zero on $Z$. It follows that $Z \cap U_0 = Z \setminus V(z_0^b) = Z \setminus V(\prod_{i \geq 1} z_i^{b_i}) \subseteq \bigcup_{i \geq 1} (Z \cap U_i)$, as required. For the last assertion it suffices to show that $Z \cap H \cap U_i$ is contained in the closure of $\pi \circ \phi(X \setminus V(f_i))$  for each $i = 0, \ldots, M$. This follows from assertion \eqref{subspace-assertion} of \cref{affine-embedding}.
\end{proof}

\subsection{Exercises}

\begin{exercise} \label{exercise:aff-to-proj}
Prove the assertion from \cref{example:aff-to-proj}. 
\end{exercise}

\begin{exercise} \label{exercise:blow-up-well-defined} 
Show that \cref{prop:blow-up-well-defined} is true in general if it holds for the special case that  $r= q+1$ and $h_j = g_j$, $j = 1, \ldots, q$.
\end{exercise}

\begin{exercise} \label{exercise:blow-up-origin}
Let $a = (a_1, \ldots, a_n) \in \kk^n$. Recall that the ideal $I(a)$ of all polynomials vanishing at $a$ is generated by $x_1 - a_1, \ldots, x_n - a_n$, so that $\bl_a(\kk^n)$ is the closure in $\kk^n \times \pp^{n-1}$ of the graph $\gr(\phi)$ of the map $\phi: \kk^n \setminus \{a\} \to \pp^{n-1}$ given by 
\begin{align*}
(x_1, \ldots, x_n) \mapsto [x_1 - a_1: \cdots : x_n - a_n]
\end{align*}
\begin{enumerate}
\item Given $b := (b_1, \ldots, b_n) \in \kk^n \setminus \{0\}$, let $L_b := \{(a_1 + tb_1, \ldots, a_n + tb_n): t \in \kk\}$ be the line through $a$ in the ``direction'' of $b$. Show that the closure of $\graph(\phi|_{L_b \setminus \{a\}})$ is $\graph(\phi|_{L_b \setminus \{a\}}) \cup \{[b_1: \cdots : b_n]\}$. 
\item Deduce that if $\sigma: \bl_a(\kk^n) \to \kk^n$ is the blow up map, then $\sigma^{-1}(a) = \pp^{n-1}$.
\end{enumerate}
\end{exercise}


\begin{exercise}  \label{exercise:linear-blow-up}
Let $Y = V(x_1, \ldots, x_k) \subseteq X := \kk^n$. Show that $\bl_Y(X)$ is the closure in $\kk^n \times \pp^{k-1}$ of the graph $\gr(\phi)$ of the map $\phi: \kk^n \setminus Y \to \pp^{k-1}$ given by 
\begin{align*}
(x_1, \ldots, x_n) \mapsto [x_1 : \cdots : x_k]
\end{align*}
If $\sigma: \bl_Y(X) \to X$ is the blow up map, then show that $\sigma^{-1}(y) = \{y\} \times \pp^{k-1}$ for each $y \in Y$. [Hint: follow the steps of \cref{exercise:blow-up-origin}. In particular, for for each $a  = (0, \ldots, 0, a_{k+1}, \ldots, a_n) \in Y$ and each $b := (b_1, \ldots, b_k) \in \kk^k \setminus \{0\}$, consider the restriction of $\phi$ to the line $L_b := \{(tb_1, \ldots, tb_k, \linebreak[1], a_{k+1}, \ldots, a_n): t \in \kk\}$.] 
\end{exercise}

\begin{exercise} \label{exercise:retraction-intersection}
Let $Y$ be a topological space, and $r:Y \to H$ be a {\em retract} onto a subset $H \subset Y$ (which means $r$ restricts to the identity map on $H$). Given $W \subset Y$, show that $\overline{W} \cap H \subset \overline{r(W)}$.
\end{exercise}

\begin{exercise}\label{exercise:closure-in-open-subset}
Let $U$ be an open subset of a topological space $Y$ and $W \subseteq U$. Show that the closure of $W$ in $U$ is the intersection of $U$ and the closure of $W$ in $Y$.

\end{exercise}

\section{Dimension}\label{dimension-section}
The \index{Dimension!of a variety}{\em dimension} $\dim(X)$ of an irreducible quasiprojective variety $X$ is the transcendence degree over $\kk$ of the field $\kk(X)$ of rational functions on $X$. In general $\dim(X)$ is the maximum of the dimensions of all its irreducible components. We say that $X$ has \index{Dimension!pure}\index{Pure dimension!of a variety}{\em pure dimension} $d$ if each of its irreducible components has dimension $d$. If $a \in X$ and there is a Zariski open neighborhood $U$ of $a$ in $X$ which has pure dimension $d$, we say that $X$ has pure dimension $d$ {\em near} $a$. The \index{Codimension}{\em codimension} of a subvariety $Y$ of $X$ is the integer $\dim(X) - \dim(Y)$.

\begin{example}
Since $\kk(\kk^n) \cong \kk(x_1,\ldots, x_n)$, it follows that $\dim(\kk^n) = n$, as expected. Since for irreducible varieties rational functions are determined by any nonempty open subset, it follows that $\dim(\pp^n) = n$ as well. A variety is zero dimensional if and only if it is a finite set (\cref {exercise:dimension-0}). 
\end{example}

\begin{example}
Let $X := V(xyz) \subseteq \kk^3$ be the union of the three coordinate planes. Each of these planes is isomorphic to $\kk^2$, and therefore has dimension two. Cosequently $X$ has pure dimension two. Let $Y$ be the union of $X$ with the ``diagonal'' line $L = \{(t,t,t): t \in \kk\}$. It is straightforward to check that $t \mapsto (t,t,t)$ induces an isomorphism $L \cong \kk$, so that $\dim(L) = 1$. It follows that $\dim(Y) = 2$ and $Y$ is {\em not} pure dimensional. However, since $X$ and $L$ interescts only at the origin, $Y$ is pure dimensional near every point except for the origin. 
\end{example}

The following fact underlies all essential properties of dimension of varieties. It is a special case of Krull's Principal Ideal Theorem (\cref{thm:principal-ideal}). 

\begin{prop} \label{prop:pure-dimension-1-0}
If $f$ is a non-constant polynomial in $(x_1, \ldots, x_n)$, then the hypersurface $V(f)$ of $\kk^n$ has pure dimension $n-1$. 
\end{prop}

\begin{proof}
Recall that the irreducible components of $X := V(f)$ are hypersurfaces $V(f_i)$, where $f_i$ are irreducible factors of $f$ (\cref{example:affine-hypersurface}). Therefore \woutlog\ we may assume $f$ is irreducible, and moreover, that $x_n$ appears in $f$ with a nonzero coefficient. There are nonnegative integers $d_1 > d_2 > \cdots$ with $d_1 \geq 1$ and nonzero polynomials $c_j$ in $(x_1, \ldots, x_{n-1})$ such that $f = \sum_j c_j(x_1, \ldots, x_{n-1}) x_n^{d_j}$. The irreducibility of $f$ implies that $I(X) = \langle f \rangle$ (\cref{example:affine-hypersurface}), so that no polynomial in $(x_1, \ldots, x_{n-1})$ identically vanishes on $X$. This implies that 
\begin{prooflist}
\item $x_1, \ldots, x_{n-1}$ are algebraically independent over $\kk$ in $\kk(X)$, i.e.\ $\trd_\kk(\kk(X)) \geq n - 1$;
\item $c_j \neq 0 \in \kk(X)$ for any $j$, so that $x_n$ is algebraic over $\kk(x_1, \ldots, x_n)$ in $\kk(X)$ via the relation
$$\sum_j c_j(x_1, \ldots, x_{n-1})x_n^{d_j} = 0$$
\end{prooflist}
It follows that $\trd_\kk(\kk(X)) = n - 1$, as required.
\end{proof}

\begin{example}
A \index{Curve}{\em curve} is a variety of pure dimension one and a \index{Surface}{\em surface} is a variety of pure dimension two. \Cref{prop:pure-dimension-1-0} implies that for all nonconstant polynomials, $V(f) \subseteq \kk^n$ is a curve if $n = 2$ and a surface if $n = 3$. 
\end{example}

The dimension of a variety is an analogue of the dimension of smooth manifolds. In particular, when $\kk = \cc$ and $X$ is irreducible, it follows from the results in \cref{nonsingular-section} that there is an open Zariski dense subset of $X$ (consisting of {\em nonsingular} points of $X$) which is a smooth manifold of (complex) dimension $\dim(X)$. Dimensions of subvarieties of a variety, however, in some ways behave rather differently from dimensions of submanifolds of a manifold:


\begin{thm} \label{thm:dim-smaller}
Every proper subvariety of an irreducible variety $X$ has smaller dimension than that of $X$.
\end{thm}

\begin{proof}
Let $Y$ be a proper subvariety of $X$. Since the field of rational functions does not change after restricting to a (nonempty) open subset, taking an open neighborhood of a point on $Y$ if necessary we may assume \woutlog\ that $X,Y$ are affine. Since $X,Y$ are also irreducible, it follows that $X = V(\ppp)$ and $Y = V(\qqq)$ for prime ideals $\ppp \subsetneq \qqq$ of $\kk[x_1, \ldots, x_n]$, $n \geq 1$. Let $d := \dim(X)$. \Woutlog\ we may assume that $x_1, \ldots, x_d$ are algebraically independent over $\kk$ in $\kk(X)$. It suffices to show [why?] that $x_1, \ldots, x_d$ are {\em not} algebraically independent over $\kk$ in $\kk(Y)$. Pick $q \in \qqq \setminus \ppp$. Since $x_1, \ldots, x_d, q$ are algebraically dependent over $\kk$ in $\kk(X)$, there is an irreducible polynomial $F(x_1, \ldots, x_d, y) \in \kk[x_1, \ldots, x_d, y]$ such that $F(x_1, \ldots, x_d, q) = 0$ in $\kk[X]$. Since $q \neq 0$ in $\kk[X]$, it follows that $F$ is {\em not} a multiple of $y$. Therefore $F(x_1, \ldots, x_d, 0)$ is a nonzero polynomial in $\kk[x_1, \ldots, x_d]$ which represents an algebraic relation of $x_1, \ldots, x_d$ in $\kk[Y]$, as required.
\end{proof}

\begin{example} \label{example:curve-in-curve}
Let $Z$ be a curve and $Z' \subseteq Z$ be such that $Z'$ contains a dense open subset of $Z$. Then \cref{thm:dim-smaller} implies that $Z \setminus Z'$ has finitely many points. Consequently $Z'$ is also open in $Z$ and therefore $Z'$ is a quasiprojective variety. This in particular implies that the image of a morphism from a curve is a quasiprojective variety. Note that this is {\em not} true for varieties of dimension greater than one, e.g.\ the image of the morphism $\kk^2 \to \kk^2$ given by $(x,y) \mapsto (x, xy)$ is not a quasiprojective variety (\cref{example:constructible-0}). 
\end{example}

We presented \cref{thm:dim-smaller} above due to the simple proof (which follows \cite[Proof of Proposition 1.14]{mummetry}). However, if one is willing to use more elaborate machinery, namely Krull's principal ideal theorem (\cref{thm:principal-ideal}), then we can obtain the following more precise version of \cref{thm:dim-smaller} in a straightforward way - its proof is left as \cref{exercise:pure-dimension}. 

\begin{thm} \label{thm:pure-dimension}
Let $f_1, \ldots, f_k$ be regular functions on a variety $X$ of pure dimension $n$. 
\begin{enumerate}
\item If $f_1$ does not identically vanish on any irreducible component of $X$, then $V(f_1) \subseteq X$ is either empty or has pure dimension $n-1$.
\item The dimension of every irreducible component of $V := V(f_1, \ldots, f_k) \subseteq X$ is at least $n - k$. In other words, to define a subvariety of codimension $m$ one needs at least $m$ equations. 
\item If $\dim(V) = n- k$, then for all $I \subseteq [k]$, $V(f_i: i \in I)$ has pure dimension $n - |I|$ near every point of $V$. \qed
\end{enumerate}
\end{thm}

\begin{example}
\Cref{prop:pure-dimension-1-0} is a special case of \cref{thm:pure-dimension}. In general the dimension of $V(f_1, \ldots, f_k)$ may be greater than $n-k$, e.g.\ $V(x(x-y), y(x-y))$ is the diagonal line $x = y$ on $\kk^2$, and therefore has dimension $1 > 2 -2$.
\end{example}

\begin{example} \label{example:noncomplete-intersection}
Consider $\phi: \kk \to \kk^3$ given by $t \mapsto (t^3, t^4, t^5)$. In \cref{exercise:noncomplete-intersection} you will show that 
\begin{enumerate}
\item $X := \phi(\kk)$ is a codimension two subvareity of $\kk^3$ which is the set of zeroes of two polynomials (in other words, $X$ is a ``set theoretic complete intersection'' in $\kk^3$),
\item but the ideal of polynomials vanishing on $X$ can not be generated by two polynomials (in other words, $X$ is {\em not} a ``complete intersection'' in $\kk^3$).
\end{enumerate}
\end{example}

%
%
%


\begin{cor}[Curve selection lemma I] \label{closure-curve-lemma}
\index{Curve!selection lemma I}
Let $U$ be a nonempty Zariski open subset of an irreducible variety $X$ and $a \in X$. If $\dim(X) > 0$, then there is an irreducible curve $C$ on $X$ containing $a$ such that $C \cap U$ is nonempty.
\end{cor}

\begin{proof}
We prove this by induction on $n := \dim(X)$. If $n= 1$, then it is trivially true, so assume $n > 1$. \Woutlog\ we may assume that $X$ is affine. 

\begin{proclaim} \label{claim:closure-curve-lemma}
There is a regular function $f$ on $X$ such that $f(a) = 0$ and $U$ intersects every irreducible component of $V(f) \subseteq X$. 
\end{proclaim}

\begin{proof}
By \cref{thm:dim-smaller} the dimension of $Z := X \setminus U$ is less than $n$. If $\dim(Z) < n -1$, then due to \cref{thm:pure-dimension} the claim is satisfied by any non-constant $f \in \kk[X]$ with $f(a) = 0$.  So assume $\dim(Z) = n -1 > 0 = \dim(\{a\})$. Then for each irreducible $(n-1)$-dimensional component $Z_i$ of $Z$, the ideal $I(Z_i)$ of regular functions vanishing on $Z_i$ can not contain the ideal $I(a)$ of regular functions vanishing at $a$. It then follows from the ``prime avoidance'' phenomenon (assertion \eqref{mp:prime-avoidence} of \cref{thm:mp}) that there is $f \in I(a)$ such that $f$ does not identically vanish on any $(n-1)$-dimensional component of $Z$. Then $f$ satisfies the claim.
\end{proof}

Let $X'$ be an irreducible component of $V(f)$ containing $a$. Then $U' := U \cap X \neq \emptyset$. Since $\dim(X') = n - 1$ (\cref{thm:pure-dimension}), we are done by induction.
\end{proof}

Intuitively, the image of a map can not have bigger dimension than the source. The following result makes it precise:

\begin{prop} \label{dominant-dimension}
Let $\phi:X \to Y$ be a dominant morphism of varieties. Then $\dim(X) \geq \dim(Y)$.
\end{prop}

\begin{proof}
This is immediate from the observation that if $X$ and $Y$ are irreducible, then $\phi^*: \kk(Y) \to \kk(X)$ is an injection (\cref{exercise:dominant-k(X)}). 
\end{proof}

Now we study the dimension of the fibers of a morphism. Consider e.g.\ the map $\phi: \kk^2 \to \kk^2$ from \cref{example:constructible-0} given by $(x,y) \mapsto (x,xy)$. The image of $\phi$ is $(\kk^2 \setminus V(x)) \cup \{(0,0)\}$. For most points on the image, in fact for all $(u,v) \in \phi(\kk^2) \setminus \{(0,0)\}$, the fiber $\phi^{-1}(u,v)$ over $(u,v)$ is zero dimensional, and consists of a single point $(u,v/u)$. On the other hand, $\phi^{-1}(0,0)$ is all of the $y$-axis, and hence has dimension one. In particular, there is a nonempty Zariski open subset of the image over which the dimension of fibers is constant, and it is possible that some fibers have higher dimension. Our next result shows that this is true in general.

\begin{thm}
\label{fiber-dimension}
Let $X$ be an irreducible variety and $\phi: X \to Y$ be a surjective morphism. Then 
\begin{enumerate}
\item \label{fiber-dimension:geq} For every $y \in Y$ and every irreducible component $V$ of $\phi^{-1}(y)$, $\dim (V) \geq \dim(X) - \dim(Y)$. 
\item \label{fiber-dimension:eq} Moreover, there is a nonempty Zariski open subset $U$ of $Y$ such that $\phi^{-1}(y)$ is of pure dimension $\dim(X) - \dim(Y)$ for each $y \in U$.
\end{enumerate}
\end{thm}

\begin{proof}[Proof {\cite[Proof of Threorem I.7]{shaf1}}]
We may assume \woutlog\ that $Y$ is affine. Let $m := \dim(Y)$. Take $y \in Y$. The following can be proved via an induction on $m$, we leave its proof as an exercise (\cref{exercise:complementary-intersection}).

\begin{proclaim} \label{claim:complementary-intersection}
There are regular functions $h_1, \ldots, h_m$ on $Y$ such that $V(h_1, \ldots, h_m) \cap U = \{y\}$ for some open neighborhood $U$ of $y$ in $Y$. \qed
\end{proclaim}

Pick $h_1, \ldots, h_m, U$ as in \cref{claim:complementary-intersection}. Since $X$ is irreducible, it follows that $\dim(\phi^{-1}(U)) = \dim(X)$ (\cref{exercise:closure-dimension}). Since $\phi^{-1}(y) \cap \phi^{-1}(U)$ is the zero-set of $h_j \circ \phi$, $j = 1, \ldots, m$, assertion \eqref{fiber-dimension:geq} follows from \cref{thm:pure-dimension}. For the second assertion we proceed by induction on $n := \dim(X)$. Clearly it is true for $n = 0$. Let $X'$ be an arbitrary nonempty open {\em affine} subset of $\phi^{-1}(Y)$. Due to the first assertion it suffices [why?] to prove that
\begin{align}
\parbox{.54\textwidth}{
there is a nonempty Zariski open subset $Y'$ of $Y$ such that $\dim(\phi^{-1}(y) \cap X') \leq n - m$ for each $y \in Y' \cap \phi(X')$.
} \label{fiber-reduction}
\end{align}
Since $X'$ is dense in $X$ and $\phi$ is continuous (\cref{exercise:Zariski-continuous}) and surjective, it follows that $\phi(X')$ is dense in $Y$. Therefore $\phi$ induces an injection $\kk[Y] \to \kk[X']$ (\cref{prop:star-injective}); we will consider $\kk[Y]$ as a subring of $\kk[X']$. Let $g_1, \ldots, g_M$ (respectively $f_1, \ldots, f_N$) be $\kk$-algebra generators of $\kk[Y]$ (respectively $\kk[X']$). Since $\dim(X') = \dim(X) = n$, the transcendence degree of $\kk(X')$ over $\kk(Y)$ is $n-m$. Therefore we may assume \woutlog\ that $f_1, \ldots, f_{n-m}$ are algebraically independent over $\kk(Y)$ and $f_j$ is algebraically dependent over $\kk(Y)(f_1, \ldots, f_{n-m})$ for each $j > n-m$. Fix $j$, $n-m+1 \leq j \leq  N$, and pick a nonzero polynomial $F_j(u_1, \ldots, u_{n-m},v)$ in variables $u_1, \ldots, u_{n-m}, v$ with coefficients in $\kk[Y]$ such that $F_j(f_1, \ldots, f_{n-m},f_j) = 0$. Write
\begin{align*}
F_j(u_1, \ldots, u_{n-m},v)
	= \sum_{i = 0}^{d_j} F_{j,i}(u_1, \ldots, u_{n-m})v^i
\end{align*}
where $d_j$ is the degree of $F_j$ in $v$, and write
\begin{align*}
F_{j,d_j}(u_1, \ldots, u_{n-m})
	= \sum_{\alpha = (\alpha_1, \ldots, \alpha_{n-m})} g_{j,\alpha} u_1^{\alpha_1} \cdots u_{n-m}^{\alpha_{n-m}}
\end{align*}
where each $g_{j,\alpha} \in \kk[Y]$. Fix $\alpha^j = (\alpha^j_1, \ldots, \alpha^j_{n-m})$ such that $g_{j,\alpha^j}$ is a {\em nonzero} regular function on $Y$. Let $\tilde Y:= Y \setminus V(\prod_j g_{j,\alpha^j})$. \Cref{thm:dim-smaller} implies that $\tilde Y$ is a nonempty open subset of $Y$.

\begin{proclaim} \label{fiber-non-containment}
There is a nonempty open subset $Y'$ of $\tilde Y$ such that for each $y \in Y'$, no irreducible component of $\phi^{-1}(y)$ is contained in $V := \bigcup_{j = n-m+1}^N V(F_{j,d_j}) \cap X'$.
\end{proclaim}

\begin{proof}
If $\phi(V)$ is not dense in $Y$, then we can simply take $Y'$ to be the complement in $\tilde Y$ of the closure of $\phi(V)$. So assume $\phi(V)$ is dense in $Y$. Then $\phi|_V$ can be extended to a surjective morphism $\phi':V' \to Y$, where $V'$ is a variety containing $V$ as a dense subset (\cref{prop:proj-complete}). Since $\dim(V') = \dim(V)$ (\cref{exercise:closure-dimension}), \cref{thm:dim-smaller} implies that $\dim(V') < n$. The inductive hypothesis then implies that there is a nonempty open subset $Y'$ of $\tilde Y$ such that $\dim(\phi'^{-1}(y)) = \dim(V') - m$ for each $y \in Y'$. The first assertion of \cref{fiber-dimension} therefore implies that for each $y \in Y'$, no component of $\phi^{-1}(y)$ can be contained in $V'$, as required.
\end{proof}

Let $Y'$ be as in \cref{fiber-non-containment}. Fix $y \in Y' \cap \phi(X')$ and an irreducible component $W$ of $\phi^{-1}(y) \cap X'$. Let $\bar f_j := f_j|_W$, $j = 1, \ldots, N$. Now fix $j$, $n-m < j \leq N$. Let $\bar F_j(u_1, \ldots, u_{n-m},v)$ be the polynomial in $\kk[u_1, \ldots, u_{n-m}, v]$ constructed from $F_j$ by evaluating each coefficient from $\kk[Y]$ at $y$. Since $\bar F_j(\bar f_1, \ldots, \bar f_{n-m}, \bar f_j) = 0$, \cref{fiber-non-containment} implies that $\bar f_j$ is algebraically dependent on $\kk(\bar f_1, \ldots, \bar f_{n-m})$ [why?]. Since $\kk[W] = \kk[\bar f_1, \ldots, \bar f_N]$, and since $\kk(W)$ is the quotient field of $\kk[W]$ (\cref{prop:k(X)-0}), it follows that $\dim(W) \leq n-m$, as required.
\end{proof}

We next show that given a surjective morphism $\phi: X \to Y$ of irreducible varieties, the set $Y_0 :=  \{y \in Y: \dim(\phi^{-1}(y)) = \dim(X) - \dim(Y)\}$ does {\em not} have to be open in $Y$ (however, we will see in the next section that it is a {\em constructible} subset of $Y$). 

\begin{example}[{\cite{mathoverflow-algeom-examples-roy-smith}}] \label{example:fiber-dimension-nopen}
The map $\psi: \kk^3 \to \kk^3$ given by $(x,y,z) \mapsto (x, y, z^2)$ is surjective. For each $\xi \in \kk\setminus \{0\}$, let $L_{\xi}$ be the line $V(y, z - \xi)$ on $\kk^3$. Assume $\character(\kk) \neq 2$. Then $\psi^{-1}(L_{\xi^2})$ is the disjoint union of $L_{\xi}$ and $L_{-\xi}$. Consider the blow up $\sigma: \bl_{L_{\xi}}(\kk^3) \to \kk^3$ of $\kk^3$ at $L_{\xi}$. Recall that $\sigma$ is isomorphism on $\sigma^{-1}(\kk^3 \setminus L_{\xi})$ and $\sigma^{-1}(c) \cong \pp^1$ for each $c \in L_{\xi}$ (\cref{example:linear-blow-up}). Now fix $c = (\rho, 0, \xi) \in L_{\xi}$. Let $X := \bl_{L_{\xi}}(\kk^3) \setminus \sigma^{-1}(c)$ and $\phi$ be the restriction of $\psi \circ \sigma$ on $X$. Then it is straightforward to check that $\phi: X \to \kk^3$ is surjective, and the set $Y_0 := \{a \in \kk^3: \dim(\phi^{-1}(a)) = 0\}$ is $(\kk^3 \setminus L_{\xi^2}) \cup \{(\rho, 0, \xi^2)\}$, and $Y_1 :=  \{a \in \kk^3: \dim(\phi^{-1}(a)) = 0\}$ is $L_{\xi^2} \setminus \{(\rho, 0, \xi^2)\}$. In particular, neither $Y_0$ nor $Y_1$ is open or closed in $Y := \kk^3$. 
\end{example}

\subsection{Exercises}

\begin{exercise} \label{exercise:dimension-0}
Show that a zero dimensional variety is the union of finitely many points, in other words, a positive dimensional variety contains infinitely many points.
\end{exercise}

\begin{exercise} \label{exercise:closure-dimension}
Let $X \subseteq Y$ be irreducible quasiprojective varieties, and let $X'$ be the closure of $X$ in $X'$. Show that $\dim(X') = \dim(X)$.
\end{exercise}

\begin{exercise} \label{exercise:pure-dimension}
Prove \cref{thm:pure-dimension} [Hint: use Krull's Principal Ideal Theorem (\cref{thm:principal-ideal})].
\end{exercise}

\begin{exercise} \label{exercise:noncomplete-intersection}
Let $\phi: \kk \to \kk^3$ be the morphism given by $t \mapsto (t^3, t^4, t^5)$. Show that 
\begin{enumerate}
\item $\phi(\kk)$ is defined by equations $y^3 - x^4 = 0$ and $z^3 - x^5 = 0$ on $\kk^3$. [Hint: if $(a,b, c)$ is a solution to these equations, and $\alpha$ is a third root of $a$, then $b = \zeta_1 \alpha^4$ and $c = \zeta_2 \alpha^5$ for third roots $\zeta_1, \zeta_2$ of $1$. One can choose a third root $\zeta$ of $1$ such that $\phi(\zeta \alpha) = (a, b, c)$.] 
\item $X$ is irreducible and $\dim(X) = 1$. [Hint: $\kk$ is birational to $\phi(\kk)$.]
\item The ideal $I(X)$ of polynomials vanishing on $X$ can not be generated by less than three polynomials. [Hint: for all monomials of degree $\leq 3$ in $(x,y,z)$, compute its pullback by $\phi$. Show that the only degree $2$ polynomial in $I(X)$ is $g_1 := y^2 - xz$, and modulo multiples of $g_1$, the only degree $3$ polynomials in $I(X)$ are $g_2 := x^3 - yz$ and $g_3 := z^2 - x^2y$. Conclude that there can not be two polynomials $f_1, f_2 \in I(X)$ such that $g_1, g_2, g_3$ are in the ideal generated by $f_1, f_2$.]
\item $\dim(X) = 1$. 
\end{enumerate}
\end{exercise}

\begin{exercise} \label{exercise:noempty-intersection-hypersurface}
Let $X$ be a subvariety of $\pp^n$ of dimension $d$.
\begin{enumerate}
\item If $d > 0$, show that $X$ intersects every hypersurface on $\pp^n$. [Hint: the complement of a hypersurface in $\pp^n$ is affine, and every projective variety is complete.]
\item Deduce that $\dim(X \cap V(f)) \geq d -1$ for every homogeneous polynomial $f$ (in the homogeneous coordinates on $\pp^n$). [Hint: \cref{thm:pure-dimension}.]
\item Deduce that $X \cap V(f_1, \ldots, f_d) \neq \emptyset$ for all homogeneous $f_1, \ldots, f_d$. 
\item Show that there are homogeneous polynomials $f_1, \ldots, f_d$ such that $X \cap V(f_1, \ldots, f_d)$ has finitely many elements. [Hint: given a subvariety $Y$ of $\pp^n$, fix a point $a_j$ on each irreducible component of $Y$, and choose a homogeneous polynomial $f$ which does not vanish at any $a_j$. Then $\dim(Y \cap V(f)) = \dim(Y) - 1$.]
\end{enumerate}
\end{exercise}

\begin{exercise} \label{exercise:noncomplete}
Let $X$ be the ``ruled surface'' $\pp^1 \times \pp^1$ and $Y$ be a ruling on $X$ (i.e.\ $Y$ is of the form $\{a\} \times \pp^1$ or $\pp^1 \times \{a\}$ for some $a \in \pp^1$). We treat $X$ as a subvariety of some projective space $\pp^n$ (recall from \cref{example:ruled} that we may take $n = 3$) with homogeneous coordinates $[x_0: \cdots : x_n]$. Show that 
\begin{enumerate}
\item $Y$ is a subvariety of pure codimension one in $X$. 
\item There is no homogeneous polynomial $f$ in $(x_0, \ldots, x_n)$ such that $Y = X \cap V(f)$. [Hint: $X \setminus V(f)$ is affine (\cref{projective-complement}). There are rulings on $X \setminus Y$ which are isomorphic to $\pp^1$ and therefore complete.]
\end{enumerate}
\end{exercise}

\begin{exercise} \label{exercise:complementary-intersection}
Let $X$ be an affine variety in $\kk^n$ of dimension $m \geq 1$ and $a \in X$. Fix $d \geq 1$. Show that it is possible to find $m$ polynomials $h_1, \ldots, h_m$ in $(x_1, \ldots, x_n)$ of degree $d$ such that $a$ is an \index{Isolated point}{\em isolated point} of $V(h_1, \ldots, h_m) \cap X$, i.e.\ there is a Zariski open subset $U$ of $a$ in $\kk^n$ such that $U \cap X \cap V(h_1, \ldots, h_m) = \{a\}$. [Hint: use induction on dimension; at every step choose a polynomial of degree $d$ which vanishes at $a$, but does not identically vanish on any of the irreducible components of the variety.]
\end{exercise}

\begin{exercise} \label{exercise:fiber-dimension-nopen}
Verify the assertions from \cref{example:fiber-dimension-nopen}. 
\end{exercise}

\begin{exercise} \label{exercise:curve-to-curve}
Let $\phi: X \to Y$ be a morphism between curves. Assume $\phi$ does not map any irreducible component of $X$ to a point. Show that $\phi$ is a finite-to-one map. 
\end{exercise}

\begin{exercise} \label{exercise:affine-inverse}
Let $C \subset \pp^n$ be a projective curve, and $\phi: C \to \pp^1$ be a surjective morphism which is not constant on any of the irreducible components of $C$. Let $T \subset \pp^1$ be a finite set. In this exercise you will show that
\begin{align}
\parbox{0.66\textwidth}{
	there is a finite set $S \subset \pp^1 \setminus T$ such that $\phi^{-1}(\pp^1 \setminus S)$ is an affine curve.
}\label{display:affine-inverse}
\end{align}
\begin{enumerate}
\item Pick a point on $\pp^1 \setminus T$; denote it by $\infty$. Show that there is a hypersurface $X$ of $\pp^n$ containing $\phi^{-1}(\infty)$ such that $X \cap \phi^{-1}(T) = \emptyset$ and $|X \cap C| < \infty$. [Hint: use \cref{exercise:curve-to-curve,exercise:complement-hypersurface}.]
\item Show that $C' := C \setminus X$ is an affine curve. [Hint: use \cref{projective-complement}.]
\item Identifying $\pp^1 \setminus \{\infty\}$ with $\kk$, show that $\phi|_{C'}$ is induced by a regular function $f$ on $C'$. [Hint: use \cref{prop:morphism-on-Uj}.]
\item Let $T' := \phi( X \cap C) \setminus \{\infty\} = \phi(C\setminus C') \cap \kk$. Show that $C'' := C' \setminus V(\prod_{t \in T'}(f - t))$ is an affine curve. [Hint: use \cref{example:affine-complement}.]
\item Conclude that \eqref{display:affine-inverse} holds with $S := T' \cup \{\infty\}$.
\end{enumerate}
\end{exercise}

\section{Image of a morphism: Part II - Constructible sets} \label{constructible-section}
A \index{Constructible set}{\em constructible} subset of a topological space is a finite union of open subsets of its closed subsets. In particular, constructible subsets of quasiprojective varieties are simply {\em finite unions of quasiprojective subsets}. The relevance of constructible sets in algebraic geometry stems from the fundamental result of C.\ Chevalley that images of morphisms of varieties are constructible sets in Zariski topology (see \cref{example:constructible-0}). Before we prove this result, we state a ``constructible version'' of the ``curve selection lemma''; it is a straightforward consequence of curve selection lemma I (\cref{closure-curve-lemma}) and its proof is left as \cref{exercise:curve-selection}. 

\begin{prop}[Curve selection lemma II] \label{closure-curve-prop}
\index{Curve!selection lemma II}
Let $W$ be a constructible subset of a variety $X$ and $\bar W$ be the Zariski closure of $W$ in $X$. Pick an irreducible component $Z$ of $\bar W$ and a point $a \in Z$. If $\dim(Z) \geq 1$, then there is an irreducible curve $C$ on $Z$ containing $a$ such that $C \cap W$ is nonempty and (Zariski) open in $C$. \qed
\end{prop}

The main result of this section is the following result of C.\ Chevalley:

\begin{thm}[Chevalley's theorem]\label{Chevalley}
\index{Chevalley's theorem}
The image $\phi(X)$ of a morphism $\phi:X \to Y$ of varieties is a constructible subset of $Y$.
\end{thm}

We give a proof of \cref{Chevalley} following \cite[Section 2C]{mummetry}. In fact we only prove \cref{thm:constructible-projection} below, and leave it as \cref{exercise:constructible-projection} to show that this is equivalent to Chevalley's theorem. Note that the statement of \cref{thm:constructible-projection} is precisely what you get from substituting ``closed subsets'' by ``constructible subsets'' in the definition of complete varieties.

\begin{thm} \label{thm:constructible-projection}
Let $X,Y$ be varieties. Then the projection map $X \times Y \to Y$ maps constructible sets to constructible sets.
\end{thm}

\begin{proof}
It suffices [why?] to consider the case that $X = \kk^n$ and $Y = \kk^m$. Then taking compositions we can further reduce it to the case of the projection $\pi: \kk \times \kk^m \to \kk^m$. It is straightforward to check, and we leave it as an exercise (\cref{exercise:constructible-projection-2}) to show, that it suffices to prove the following statement:
\begin{align}
\parbox{0.8\textwidth}{
If $V$ is an irreducible subvariety of $\kk \times \kk^m \cong \kk^{m+1}$ and $U$ is a nonempty open subset of $V$, then $\pi(U)$ contains a nonempty open subset of the closure of $\pi(V)$ in $\kk^m$.
}  \label{constructible-projection-claim}
\end{align}
We now prove \eqref{constructible-projection-claim}. Let $W$ be the closure of $\pi(V)$ in $\kk^m$. Then $\pi$ induces an injective map $\kk[W] \into \kk[V]$ (\cref{prop:star-injective}), so that $\kk[V] \cong \kk[W][x_1]/\ppp$ for some ideal $\ppp$ of the polynomial ring $\kk[W][x_1]$ in one variable over $\kk[W]$. At first consider the case that $\ppp = 0$. Then $V = \kk \times W$. Let $(a_1, \ldots, a_{m+1})$ be any point of $U$. Then $W' :=  \{a_1\} \times W$ is a subvariety of $V$ and therefore $U \cap W'$ is a nonempty open subset of $W'$ whose projection is open in $W$. Now consider the remaining case that $\ppp \neq 0$. Then \cref{thm:dim-smaller} and \cref{dominant-dimension} imply that $\dim(V) = \dim(W)$. Let $\bar V$ be the closure of $V$ in $\pp^1 \times \kk^m$ and $V' := \bar V \setminus U$. Since $\pp^1$ is complete, it follows that $\pi(\bar V)$ and $\pi(V')$ are closed in $\kk^m$. Since $\pi(\bar V)$ is closed, it follows that $\pi(\bar V) \supset W$. On the other hand $\dim(V') < \dim(\bar V) = \dim(V) = \dim(W)$ (\cref{thm:dim-smaller,exercise:closure-dimension}), and therefore $\pi(V')$ cannot contain $W$ (\cref{dominant-dimension}). Since $\pi(U) \supset \pi(\bar V) \setminus \pi(V') \supset W \setminus \pi(V')$, the claim follows.
\end{proof}

We now extend Chevalley's theorem and show that the set of all points $y$ in the target space of a morphism $\phi$ such that $\phi^{-1}(y)$ has a given dimension is constructible (see \cref{example:fiber-dimension-nopen}). 

\begin{cor} \label{constructible-fiber-dimension}
Let $\phi:X \to Y$ be a morphism of varieties and $k$ be a nonnegative integer. Then $Y_k := \{y \in Y: \dim(\phi^{-1}(y)) < k\}$ is a constructible subset of $Y$.
\end{cor}

\begin{proof}
We proceed by double induction on $m_\phi := \dim(\phi(X))$ and $k$. Due to \cref{fiber-dimension,Chevalley} the corollary is true whenever $m_\phi = 0$ or $k \leq \dim(X) - m_\phi$. Now assume it is true for all $\phi$ such that $m_{\phi} < m$. Pick $\phi$ with $m_\phi = m$ and $k'$ such that the corollary holds for $\phi$ and $k'$. We will show that it holds for $\phi$ and $k' + 1$. By the inductive hypothesis $Y \setminus Y_{k'}$ is constructible, and therefore is a union of quasiprojective varieties. Let $Y^0_{k'}$ be an irreducible component of $Y \setminus Y_{k'}$, and $X'^0_{k'}$ be an irreducible component of $\phi^{-1}(Y^0_{k'})$. Note that both $Y^0_{k'}$ and $X'^0_{k'}$ are quasiprojective varieties, so that $Y'^0_{k'} := \phi(X'^0_{k'})$ is constructible due to Chevalley's theorem (\cref{Chevalley}). It suffices to show that $Y'^0_{k'} \cap Y_{k'+1}$ is constructible. Let $\phi'^0_{k'}$ be the restriction of $\phi$ to $X'^0_{k'}$. Then $\phi'^0_{k'}: X'^0_{k'} \to Y'^0_{k'}$ is surjective and by construction $\dim((\phi'^0_{k'})^{-1}(y)) \geq k'$ for each $y \in Y'^0_{k'}$. If $k' < \dim(X'^0_{k'}) - \dim(Y'^0_{k'})$, then \cref{fiber-dimension} implies that $Y_{k'+1} \cap Y'^0_{k'} = \emptyset$, which is trivially constructible. Otherwise \cref{fiber-dimension} implies that $k' =  \dim(X'^0_{k'}) - \dim(Y'^0_{k'})$ and there is a nonempty Zariski open subset $U'^0_{k'}$ of $Y'^0_{k'}$ such that $\dim((\phi'^0_{k'})^{-1}(y))= k'$ for each $y \in U'^0_{k'}$. Let $Y''^0_{k'} := Y'^0_{k'} \setminus U'^0_{k'}$ and $\phi''^0_{k'}$ be the restriction of $\phi$ to $\phi^{-1}(Y''^0_{k'})$. Since $\dim(Y''^0_{k'}) < \dim(Y) = m$, the corollary is true for $\phi''^0_{k'}$ (and {\em all} values of $k$), so that $Y_{k'+1} \cap Y''^0_{k'}$ is constructible. Therefore $Y_{k'+1} \cap Y'^0_{k'} = U'^0_{k'} \cup (Y_{k'+1} \cap Y''^0_{k'})$ is constructible as well.
\end{proof}

\subsection{Exercises}

\begin{exercise} \label{exercise:complementally-constructible}
If $Y$ is a constructible subset of $X$, show that $X \setminus Y$ is also a constructible subset of $X$. 
\end{exercise}

\begin{exercise} \label{exercise:transitively-constructible}
If $Y$ is a constructible subset of $X$ and $Z$ is a constructible subset of $Y$, then show that $Z$ is a constructible subset of $X$. 
\end{exercise}

\begin{exercise}\label{exercise:curve-selection}
Prove \cref{closure-curve-prop}.
\end{exercise}

\begin{exercise} \label{exercise:constructibly-open}
Let $\phi: X \to Y$ be a dominant morphism of varieties and $U$ be a constructible subset of $Y$. Show that the following are equivalent:
\begin{enumerate}
  \item $U$ contains a nonempty Zariski open subset of $Y$.
  \item $\phi^{-1}(U)$ contains a nonempty Zariski open subset of $X$.
\end{enumerate}
\end{exercise}

\begin{exercise} \label{exercise:constructible-inverse}
Let $\phi: X \to Y$ be a morphism of varieties and $U$ be a constructible subset of $X$. Show that for each constructible subset $V$ of $Y$, $\phi^{-1}(V) \cap U$ is a constructible subset of $X$.
\end{exercise}

\begin{exercise} \label{exercise:dimstructibly-open}
The \index{Dimension!of a constructible subset}{\em dimension} of a constructible subset $U$ of a variety $X$ is simply the dimension of the Zariski closure of $U$ in $X$. Let $\phi: X \to Y$ be a dominant morphism of varieties. If $X$ is irreducible, then show that the following are equivalent:
\begin{enumerate}
  \item $U$ contains a nonempty Zariski open subset of $X$.
  \item there is a nonempty Zariski open subset $Y'$ of $Y$ such that for each $y \in Y'$, $\dim(\phi^{-1}(y) \cap U) = \dim(X) - \dim(Y)$.
\end{enumerate}
\end{exercise}

\begin{exercise}\label{exercise:constructible-projection}
Show that \cref{Chevalley} is equivalent to \cref{thm:constructible-projection}. 
\end{exercise}

\begin{exercise}\label{exercise:constructible-projection-2}
Show that it suffices to prove \eqref{constructible-projection-claim} in order to prove that the projection $\kk \times \kk^m \to \kk^m$ maps constructible sets to constructible sets. [Hint: use induction on dimension of the constructible subset of $\kk \times \kk^m$.]
\end{exercise}

\section{Tangent space, singularities, local ring at a point}\label{nonsingular-section}
\subsection{The case of affine varieties}
Consider a straight line $L = \{a + tv: t \in \kk\}$ through a point $a = (a_1, \ldots, a_N) \in \kk^N$, where $v = (v_1, \ldots, v_N) \in \kk^N$ determines the ``direction'' of $L$. Assume $f(a) = 0$, where $f$ is a polynomial in $(x_1, \ldots, x_N)$. We say that $L$ is {\em tangent} to $V(f)$ at $a$ if $\ord_t(f(a+tv)) > 1$, or equivalently, if
\begin{align}
\sum_{i=1}^N v_i\partialxifrac{f}(a) = 0 \label{partial-tangency}
\end{align}
More generally, $L$ is \index{Tangent!line}{\em tangent at $a$} to an {\em affine variety $X$} containing $a$ if \eqref{partial-tangency} holds for {\em all} $f$ vanishing on $X$. The \index{Tangent!space}{\em tangent space} $\tangent{X}{a}$ to $X$ at $a$ is the union of all tangent lines to $X$ at $a$. It is straightforward to check (\cref{exercise:tangent-equation}) that
\begin{align}
\tangent{X}{a} = V \left(\sum_{i=1}^N (x_i - a_i)\partialxifrac{f}(a) : f \in I(X) \right) \label{tangent-equation}
\end{align}
where $I(X)$ is the ideal in $\kk[x_1, \ldots, x_N]$ of polynomials vanishing on $X$. It is clear from \eqref{tangent-equation} that $T_a(X)$ is of the form $V + a$ where $V$ is a linear subspace (through the origin) of $\kk^N$; the {\em dimension} of $\tangent{X}{a}$ is simply the dimension of $V$ (as a vector space over $\kk$). Let $f_1, \ldots, f_s$ be a set of generators of $I(X)$, then \cref{exercise:tangent-equation} implies that
\begin{align}
\dim \tangent{X}{a}
	&= N - \rank \left(
				\partialxiifrac{f_i}{j}(a)
		   \right)_{\substack{1 \leq i \leq s \\ 1 \leq j \leq N}}
	\label{tangent-dim}
\end{align}
where $\rank(\cdot)$ is the rank over $\kk$ of the corresponding matrix.

\begin{prop} \label{tangent-upper-semicontinuous}
For each integer $k$, the set $X_{\geq k} := \{a \in X: \dim \tangent{X}{a} \geq k\}$ is Zariski closed in $X$; in other words, the map $X \mapsto \zz$ given by $a \mapsto \dim\tangent{X}{a}$ is {\em upper semicontinuous}.
\end{prop}

\begin{proof}[Proof {\cite[Section 1A]{mummetry}}]
Let $\qqq_k$ be the ideal of $\kk[x_1, \ldots, x_N]$ generated by determinants of $(N-k+1) \times (N-k+1)$-minors of the matrix $(\partialxiifrac{f_i}{j}(a))$. Identity \eqref{tangent-dim} implies that $X_{\geq k} = V(I(X) + \qqq_k)$.
\end{proof}

Let $d := \min \{\dim \tangent{X}{a}: a \in X\}$. Assume $X$ is irreducible. Then we say that $a \in X$ is a \index{Singular!point of a variety}{\em singular} (respectively \index{Nonsingular!point of a variety}{\em nonsingular}) point if $\dim \tangent{X}{a} > d$ (respectively $\dim \tangent{X}{a} = d$). \Cref{tangent-upper-semicontinuous} implies that the set of nonsingular points of $X$ is a nonempty Zariski open subset of $X$.

\begin{example} \label{example:curve-singularities}
Every line or a quadric curve on $\kk^2$ is everywhere nonsingular (\cref{exercise:nonsingular<3}). Both the curves $C_1 = \{x^2 = y^3\}$ and $C_2 = \{x^2 = y^2 - y^3\}$ are singular at the origin (\cref{exercise:singular3}); see \cref{fig:curve-singularities}. The singularity of $C_1$ is called a \index{Cusp}{\em cusp}, and, when $\character(\kk) \neq 2$, that of $C_2$ is called a \index{Node}{\em node}.
\end{example}

\begin{center}
\begin{figure}[htb]
\begin{subfigure}[b]{0.2\textwidth}
\begin{tikzpicture}[scale=0.5]
\def\alabelx{2.6}
\def\alabely{2.7}
\def\nsamples{103}
\def\tmin{-3}
\def\tmax{3}

\def\xmin{-1.2}
\def\xmax{1}
\def\ymin{-1.3}
\def\ymax{1.3}
\def\initialshift{3}

\begin{axis}[
xmin = \xmin, xmax=\xmax, ymin = \ymin, ymax= \ymax,
axis equal=true, axis equal image=true, hide axis
]
\addplot[blue, thick, domain=\tmin:\tmax, samples=\nsamples] ({x^3} ,{x^2});
\end{axis}
\draw (\alabelx, \alabely) node [below] {\picfontsize $O$};

\end{tikzpicture}
\caption{$x^2 = y^3$}
\label{fig:cusp}
\end{subfigure}\hspace{0.1\textwidth}
\begin{subfigure}[b]{0.2\textwidth}
\begin{tikzpicture}[scale=0.5]
\def\alabelx{2.6}
\def\alabely{2.7}
\def\nsamples{103}
\def\tmin{-3}
\def\tmax{3}

\def\xmin{-1.2}
\def\xmax{1}
\def\ymin{-1.3}
\def\ymax{1.3}
\def\initialshift{3}

\begin{axis}[
xmin = \xmin, xmax=\xmax, ymin = \ymin, ymax= \ymax,
axis equal=true, axis equal image=true, hide axis
]
\addplot[blue, thick, domain=\tmin:\tmax, samples=\nsamples] ({x-x^3} ,{1 - x^2});
\end{axis}
\draw (\alabelx, \alabely) node [below] {\picfontsize $O$};

\end{tikzpicture}
\caption{$x^2 = y^2 - y^3$}
\label{fig:node}
\end{subfigure}
\caption{Some curve singularities}
\label{fig:curve-singularities}
\end{figure}
\end{center}

\subsection{Intrinsicness of the tangent space; tangent spaces and singularities on arbitrary varieties} \label{intrinsiction}
The definitions of tangent spaces and (non-)singular points given above applies only to irreducible affine varieties; moreover, they depend on the defining equations of the affine variety, and a priori it is not clear if they are preserved by isomorphisms. In this section we extend these notions to arbitrary varieties. We need two kinds of objects for this; the first one is a ``derivation'': given a point $a$ of a subvariety $X$ of $\kk^N$, a map $D: \kk[X] \to \kk$ is called a \index{Derivation}{\em derivation centered at $a$} if $D$ satisfied the following properties:
\begin{defnlist}
\item $D$ is $\kk$-linear,
\item \label{derivative-property} $D(fg) = f(a)D(g) + g(a)D(f)$ for all $f, g \in \kk[X]$,
\item $D(\alpha) = 0$ for all $\alpha \in \kk$.
\end{defnlist}
The following is straightforward to see; we leave the proof as an exercise.

\begin{prop} \label{prop:poly-derivations}
Given $a = (a_1, \ldots, a_N) \in \kk^N$, the derivations $D: \kk[x_1, \ldots, x_N] \to \kk$ centered at $a$ are in one-to-one correspondence with $\kk^N$ given by:
\begin{align}
\lambda = (\lambda_1, \ldots, \lambda_N) \mapsto D_\lambda,\ \text{where}\
D_\lambda(f) = \sum_{i=1}^N \lambda_i \partialxifrac{f}(a),\ \text{for all}\ f \in \kk[x_1, \ldots, x_N] \label{D-lambda} \qed
\end{align}
\end{prop}

The \index{Local ring!of a variety at a point}``local ring'' of a variety $X$ at $a \in X$ is the set of regular functions on arbitrarily small neighborhoods of $a$ in $X$; more precisely, consider a binary relation $\sim$ on the collection of pairs $(f, U)$, where $U$ is an open neighborhood of $a$ in $X$ and $f$ is a regular function on $U$, as follows: 
\begin{align*}
(f, U) \sim (f',U')\ \text{if and only if $f$ and $f'$ agree on}\ U \cap U'
\end{align*}
(Note the similarity to the definition of rational functions.) It is clear that $\sim$ is an {\em equivalence relation}; the {\em local ring $\local{X}{a}$ of $X$ at $a$} is the set of equivalence classes of $\sim$ with the natural $\kk$-algebra structure. The following result compiles basic properties of $\local{X}{a}$ - its proof is straightforward and left as \cref{exercise:locally-basic}. 

\begin{prop} \label{prop:locally-basic}
Let $X$ be a quasiprojective variety and $a \in X$. 
\begin{enumerate}
\item \label{locally-basic:open-restriction} If $U$ is an open neighborhood of $a$ in $X$, then $\local{X}{a} \cong \local{U}{a}$.
\item A morphism $\phi: X \to Y$ of varieties induces by pullback a $\kk$-algebra morphism $\phi^*: \local{Y}{\phi(a)} \to \local{X}{a}$. If $\phi$ is an isomorphism, then $\phi^*$ is an isomorphism of $\kk$-algebras. 
\item Let $\mmm_a$ be the ideal of $\local{X}{a}$ generated by all regular functions on neighborhoods of $a$ on $X$ which vanish at $a$. Then $\mmm_a$ is the {\em unique} maximal ideal of $\local{X}{a}$; in other words, $\local{X}{a}$ is a {\em local ring}.
\item \label{locally-basic:sublocal} Assume $a \in Z \subseteq X$, where $Z$ is a subvariety of $X$. Let $\qqq_Z$ be the ideal of $\local{X}{a}$ generated by all elements of the form $(f, U)$ such that $f|_{Z \cap U} \equiv 0$. Then $\local{Z}{a} \cong \local{X}{a}/\qqq_Z$. \qed
\end{enumerate}
\end{prop}

\begin{example} \label{example:Kn-local}
Given $a = (a_1, \ldots, a_n) \in \kk^n$, $\local{\kk^n}{a} = \{f/g: f, g \in \kk[x_1, \ldots, x_n],\ g(a) \neq 0\}$; in other words, $\local{\kk^n}{a}$ is the {\em localization} of $\kk[x_1, \ldots, x_n]$ at the ideal generated generated by polynomials vanishing at $a$. Assertion \eqref{locally-basic:open-restriction} of \cref{prop:locally-basic} implies that $\local{\pp^n}{a} \cong \local{\kk^n}{a}$. Given any other point $b = (b_1, \ldots, b_n) \in \kk^n$, the automorphism of $\kk^n$ given by the translation $x \mapsto x + (b-a)$ induces a $\kk$-algebra isomorphism between $\local{\kk^n}{a}$ and $\local{\kk^n}{b}$. 
\end{example}

\begin{example}
Consider the map $\phi: \kk^2 \to \kk^2$ from \cref{example:constructible-0} given by $(x,y) \mapsto (x, xy)$. Then $\phi$ maps the origin to itself and $\phi^*: \local{\kk^2}{\origin} \to \local{\kk^2}{\origin}$ is not surjective, in particular it is {\em not} an isomorphism. However, $\phi^*$ does induce an isomorphism between $\local{\kk^2}{\phi(a_1, a_2)}$ and $\local{\kk^2}{(a_1, a_2)}$ whenver $a_1 \neq 0$ (\cref{exercise:locally-x-to-xy}). Note that $\phi$ restricts to an automorphism on $\kk^2 \setminus V(x)$. In general the following is true: for $a \in X$, $b \in Y$, $\local{X}{a} \cong \local{Y}{b}$ as $\kk$-algebras if and only if there are open neighborhoods $U$ of $a$ in $X$ and $V$ of $b$ in $Y$ such that $U \cong V$ (\cref{exercise:local-isomorphism}). 
\end{example}

\Cref{example:Kn-local}, together with assertion \eqref{locally-basic:sublocal}  of \cref{prop:locally-basic}, immediately implies that local rings of affine varieties are {\em localizations} of the coordinate ring. 

\begin{prop} \label{prop:affine-local}
Assume $X$ is a subvariety of $\kk^N$ with coordinates $(x_1, \ldots, x_N)$. Then $\local{X}{a}$ is the ring of rational functions $\{f/g: f, g \in \kk[x_1, \ldots, x_N],\ g(a) \neq 0\}$ modulo the ideal generated by polynomials vanishing on $X$; in other words $\local{X}{a}$ is the localization of $\kk[X]$ with respect to the multiplicative set of regular functions {\em not} vanishing at $a$. \qed
\end{prop} 

\begin{cor} \label{prop:irreducible-local}
The local ring of an irreducible variety at a point is an integral domain.
\end{cor}

\begin{proof}
Since the coordinate ring of an irreducible affine variety is an integral domain (\cref{exercise:prime=integral}), the claim follows immediately from \cref{prop:affine-local}. 
\end{proof}

The following result describes the relation among derivations, local rings and tangent spaces: 

\begin{thm} \label{thm:tangent-space}
Given $a \in X$, where $X$ is an affine variety, the following spaces are isomorphic as vector spaces over $\kk$:
\begin{enumerate}
\item \label{txa} $\tangent{X}{a}$ with the vector space structure on $\tangent{X}{a}$ induced from that of $\tangent{X}{a} - a$ (in particular, $a$ is the origin of $\tangent{X}{a}$),
\item \label{dxa} the space of derivations $\kk[X] \to \kk$ centered at $a$,
\item \label{mxa} the space $(\mmm_a/\mmm_a^2)^*$ of linear functions from $\mmm_a/\mmm_a^2$ to $\kk$, where $\mmm_a$ is the maximal ideal of $\local{X}{a}$ generated by functions vanishing at $a$.
\end{enumerate}
\end{thm}

\begin{proof}[Proof {\cite[Section 1A]{mummetry}}]
A derivation $\kk[X] \to \kk$ centered at $a$ is the same as a derivation $\kk[x_1, \ldots, x_N] \to \kk$ centered at $a$ which vanishes on $I(X)$. Derivations $\kk[x_1, \ldots, x_N] \to \kk$ centered at $a$ are of the form $D_\lambda$ from the map defined in \eqref{D-lambda}. Given $\lambda \in \kk^n$, identity \eqref{tangent-equation} implies that $D_\lambda$ vanishes on all $f \in I(X)$ if and only if $\lambda + a \in \tangent{X}{a}$, which proves the isomorphisms between \eqref{txa} and \eqref{dxa}. \Cref{exercise:local-derivation} below shows that every derivation $\kk[X] \to \kk$ centered at $a$ defines an element of $(\mmm_a/\mmm_a^2)^*$. Conversely, given a linear map $\phi: \mmm_a \to \kk$ such that $\phi|_{\mmm_a^2} \equiv 0$, let $\lambda_\phi := (\phi((x_1-a_1)|_X), \ldots, \phi((x_N - a_N)|_X)) \in \kk^N$, and define $D_\phi := D_{\lambda_\phi}$ as in \eqref{D-lambda}. Let $f \in I(X)$. The Taylor series expansion of $f$ shows that $f = f(a) + \sum_{i=1}^n (x_i-a_i) \partialxifrac{f}(a) + f'$, where $f'$ is of order two or higher in the $(x_i - a_i)$. Since $f(a) = 0$ and $\phi(f|_X) = 0 = \phi(f'|_X)$ [why?], it follows that $\sum_{i=1}^N \phi((x_i-a_i)|_X) \partialxifrac{f}(a) = 0$, i.e.\ $D_\phi(f) = 0$. It follows that $D_\phi$ is a derivation $\kk[X] \to \kk$. It is straightforward to check that the maps we defined between \eqref{dxa} and \eqref{mxa} are inverse to each other, and induce an isomorphism of vector spaces.
\end{proof}

The \index{Tangent!space}{\em tangent space} to a variety $X$ at a point $a \in X$ is by definition $\tangent{X'}{a}$ for any affine open neighborhood $X'$ of $a$ in $X$. \Cref{thm:tangent-space} shows that $\tangent{X}{a} = (\mmm_a/\mmm_a^2)^*$, where $\mmm_a$ is the maximal ideal of $\local{X}{a}$ generated by functions vanishing at $a$; in particular, $\tangent{X}{a}$ is well defined for an arbitrary quasiprojective variety.

\begin{prop} \label{tangent-dimension}
Assume $X$ is an irreducible variety. Then $\min \{\dim \tangent{X}{a}: a \in X\} = \dim(X)$.
\end{prop}

\begin{proof}
Let $d := \min \{\dim \tangent{X}{a}: a \in X\}$ and $U := \{x \in X: \dim \tangent{X}{a} = d\}$. \Cref{tangent-upper-semicontinuous} implies that $U$ is open and dense in $X$. Since $X$ is birational to a hypersurface (\cref{exercise:birational-hypersurface}), and since dimension is invariant under birational maps (\cref{exercise:birational-fields}), \woutlog\ we may assume $X = V(f) \subset \kk^N$ for some polynomial $f$ in $(x_1, \ldots, x_N)$. If $f = 0$, then an easy computation shows that $U = \kk^N$ (\cref{exercise:tangent-k^n}) and $d = N$, as required. Otherwise $f$ is a nonconstant irreducible polynomial, and since $\kk$ is algebraically closed, not all the partial derivatives $\partial f/\partial x_j$ vanish identically on $X$ (\cref{exercise:nonzero-derivative}). It then follows the definition of the tangent space (or identity \eqref{tangent-dim}) that $d = N - 1 = \dim(X)$ (the last equality uses \cref{prop:pure-dimension-1-0}). 
\end{proof}

Now we can extend the notion of (non-)singular points to arbitrary (possibly reducible) varieties. Let $X$ be a variety and $a \in X$. Define $\dim_a(X)$ to be the maximum of the dimensions of the irreducible components of $X$ containing $a$. We say that $a$ is a \index{Singular!point of a variety}{\em singular} (respectively \index{Nonsingular!point of a variety}{\em nonsingular}) point of $X$ if $\dim \tangent{X}{a} > \dim_a(X)$ (respectively  $\dim \tangent{X}{a} = \dim_a(X)$). A variety is called \index{Singular!variety}{\em singular} if it has a singular point; otherwise it is called \index{Nonsingular!variety}{\em nonsingular}. The following is an immediate consequence of 
 \Cref{tangent-upper-semicontinuous,tangent-dimension}: 
 
 \begin{prop} \label{prop:openly-nonsingular}
 The set of nonsingular points of a quasiprojective variety $X$ is Zariski open and has a nonempty intersection with every irreducible component of $X$. \qed
 \end{prop}
 
\begin{example} \label{example:curve-singularities-1}
 \Cref{prop:openly-nonsingular} implies that a curve has at most finitely many singular points. Consider the affine curves $C_1 = V(x^2 - y^3)$ and $C_2 = V(x^2 - y^2 - y^3)$ from \cref{example:curve-singularities}. Embedding $\kk^2$ into $\pp^2$ with homogeneous coordinates $[x : y : z]$ via the map $(x, y) \mapsto [x:y:1]$ and homogenizing with respect to the $z$-coordinate shows that the closures of $C_j$ in $\pp^2$ are $\bar C_1 := V(x^2z - y^3)$ and $\bar C_2 := V(x^2z - y^2z - y^3)$ (\cref{example:principal-closure}). It follows that for both $j$, $\bar C_j \setminus C_j$ has only one point, namely $P := [1:0:0]$. Identifying the basic open subset $U := \pp^2 \setminus V(x)$ with $\kk^2$ with coordinates $(u, v) := (y/x, z/x)$, we see that $C_1 \cap U = V(v - u^3)$ and $C_2 \cap U = V(v - u^2 - u^3)$ (\cref{exercise:basic-open-equation}). Since $P = (0,0)$ with respect to $(u,v)$-coordinates, it follows that both $\bar C_j$ are nonsingular at $P$. 
\end{example}

\subsection{Equations near a nonsingular point}
Recall that for a subvariety $X$ of codimension $k$ in $\kk^N$, one needs at least $k$ polynomials to generate the ideal $I(X)$ of polynomials vanishing on it (\cref{thm:pure-dimension}). We say that a subvariety $X$ is a \index{Complete!intersection}{\em complete intersection} if $I(X)$ can be generated by $k$ polynomials. Not all varieties are complete intersections. Indeed, we have seen in \cref{example:noncomplete-intersection} that the image $X$ of the morphism $\kk \to \kk^3$ given by $t \mapsto (t^3, t^4, t^5)$ is a codimension two subvariety of $\kk^3$ and it takes at least three polynomials to generate $I(X)$. On the other hand, \cref{exercise:no-to-yes-comeplete-intersection} below shows that if $U = \kk^3 \setminus V(x)$ (where $(x,y,z)$ are coordinates on $\kk^3$), then the ideal of $X \cap U$ in $\kk[U] = \kk[x,y,z, 1/x]$ can be generated by two regular functions on $U$, namely $y - (y/x)^4$, $z - (y/x)^5$, and in addition, $X \cap U$ is {\em nonsingular}. In this section we show that in general every nonsingular point on a variety has an affine neighborhood which is a complete intersection. We follow the approach of \cite[Proof of Theorem 1.16]{mummetry}. Let $R := \local{\kk^N}{0}$ and $\hat R := \kk[[x_1, \ldots, x_N]]$ be the ring of formal power series in $(x_1, \ldots, x_N)$ over $\kk$. Recall that $\hat R$ has only one maximal ideal, namely the ideal $\hat \mmm$ generated by all polynomials with zero constant term, and we can view $R$ as a subring of $\hat R$ via the expansion $(1 - \sum_{i=1}^N x_ig_i)^{-1} = 1 + \sum_{j \geq 1} (x_ig_i)^j$ for any $g_1, \ldots, g_N \in \kk[x_1, \ldots, x_N]$ (see \cref{formal-power-section}). 

\begin{lemma} \label{hat-intersection}
If $\qqq$ is an ideal of $R$, then $\qqq \hat R \cap R = \qqq$.
\end{lemma}

\begin{proof}
Given $f \in \qqq \hat R \cap R$, it suffices to show that $f \in \qqq$. Indeed, write $f = \sum_j \phi_j f_j$, where $f_j$ are polynomials which generate $\qqq$ and $\phi_j$ are power series in $(x_1, \ldots, x_n)$. For each $k > \deg(f)$, if $\phi_{j,k}$ are the {\em polynomials} consisting of all monomial terms of $\phi_j$ of order at most $k$, then $g_k := \sum_j f_j \phi_{j,k} \in \qqq$, so that $f = g_k + \sum_j f_j(\phi - \phi_{j,k}) \in \qqq + \hat \mmm^k$, where $\hat \mmm$ is the (unique) maximal ideal of $\hat R$. Let $\mmm := \hat \mmm \cap R$ be the (unique) maximal ideal of $R$. Since $\hat \mmm^k \cap R = \mmm^k$ (\cref{prop:m-hat-m}), it follows that $f \in \bigcap_{k \geq 0} (\qqq + \mmm^k) = \qqq$ (\cref{thm:Krull}).
\end{proof}

\begin{cor} \label{locally-prime}
Let $f_1, \ldots, f_r$ be polynomials in $(x_1, \ldots, x_N)$ with no constant term and linearly independent (over $\kk$) linear terms. Then the ideal $\ppp$ generated by $f_1, \ldots, f_r$ in $\local{\kk^N}{0}$ is prime. 
\end{cor}

\begin{proof}
$f_1, \ldots, f_r$ generate a prime ideal $\qqq$ in $\hat R$ (\cref{cor:hat-prime}). Now apply \cref{hat-intersection}.
\end{proof}

\begin{thm} \label{local-equations}
Every nonsingular point on a variety has an affine open neighborhood which is irreducible and a complete intersection. 
\end{thm}

\begin{proof}
Let $a$ be a nonsingular point of a variety $X$, with $\dim_a(X) = n$. We may assume \woutlog\ that $X$ is a subvariety of $\kk^m$ for some $m \geq n$, and $a$ is the origin in $\kk^m$. Identity \eqref{tangent-dim} implies that there are $f_1, \ldots, f_{m-n} \in I(X)$ with no constant term and linearly independent linear terms. Let $\qqq$ be the ideal of $\kk[x_1, \ldots, x_m]$ generated by $f_1, \ldots, f_{m-n}$. \Cref{locally-prime} implies that the ideal $\qqq \local{\kk^m}{a}$ generated by $\qqq$ in $\local{\kk^m}{a}$ is prime, so that $\qqq' := (\qqq \local{\kk^m}{a}) \cap \kk[x_1, \ldots, x_m]$ is also prime. Let $Z := V(\qqq)$ and $Z' := V(\qqq') \subseteq \kk^m$. Note that $Z'$ is irreducible and $a \in Z' \subseteq Z$.

\begin{proclaim} \label{local-equations:equal'}
$\dim(Z') = n$. There is a polynomial $g$ such that $g(a) \neq 0$ and $X \setminus V(g) = Z' \setminus V(g) =  Z \setminus V(g)$.
\end{proclaim}

\begin{proof}
Let $h'_1, \ldots, h'_s$ be a set of generators of $\qqq'$. Then each $h'_j$ can be expressed as as $h_j/g_j$ for some polynomials $g_j,h_j$ such that $h_j \in \qqq$ and $g_j(a) \neq 0$. If $g := \prod_j g_j$, then it follows that $Z' \setminus V(g) = Z \setminus V(g) \supseteq X \setminus V(g)$. Since $Z$ is defined by $m-n$ equations in $\kk^m$, it follows that $\dim (Z'\setminus V(g)) = \dim (Z \setminus V(g))  \geq n$ (\cref{thm:pure-dimension}). On the other hand the assumptions on linear parts of $f_j$ and \cref{tangent-dimension} imply that $\dim (Z' \setminus V(g)) \leq n$. It follows that $\dim (Z' \setminus V(g))= n \leq \dim(X \setminus V(g))$. On the other hand, since $Z'$ is irreducible, every proper subvariety of $Z' \setminus V(g)$ has dimension smaller than $n$ (\cref{thm:dim-smaller}). It follows that $Z' \setminus V(g) = X \setminus V(g)$, which completes the proof.
\end{proof}

Since $Z \setminus V(g)$ is isomorphic to the subvariety of $\kk^{m+1}$ defiend by $f_1, \ldots, f_{m-n}, gx_{m+1} - 1$ (\cref{example:affine-complement}), it is a complete intersection. The proof is now complete due to \cref{local-equations:equal'}. 
\end{proof}

\begin{cor} \label{prop:nonsingular-local-irreducible}
The local ring of a variety at a nonsingular point is an integral domain.
\end{cor}

\begin{proof}
Since the local ring of an irreducible variety at a point is an integral domain (\cref{prop:irreducible-local}), this follows directly from \cref{local-equations}. 
\end{proof}

\subsection{Parametrizations of a curve at a nonsingular point} \label{curve-dvr-section}
The local ring at the origin of $\kk$ is $\local{\kk}{\origin} = \{f/g: f, g \in \kk[t],\ g(\origin) \neq 0\}$. Recall that the \index{Order!of a rational function in one variable}{\em order} of a polynomial $f \in \kk[t]$, denoted $\ord(f)$, is the smallest integer $d$ such that the coefficient of $t^d$ in $f$ is nonzero. One can extend $\ord$ uniquely to $\kk(t)$ by defining 
\begin{align*}
\ord(f/g) := \ord(f) - \ord(g)
\end{align*}
It is straightforward to check that $\ord$ satisfies the following properties: $\ord(fg) = \ord(f) + \ord(g)$ and $\ord(f+g) \geq \min\{\ord(f), \ord(g)\}$, and $\local{\kk}{\origin} = \{h \in \kk(t): \ord(h) \geq 0\}$. In other words, $\ord$ is a {\em discrete valuation} on $\kk(t)$ and $\local{\kk}{\origin}$ is a {\em discrete valuation ring} (see \cref{discrete-valuection} for a discussion on discrete valuations). In this section we will see that the local ring of any curve at a nonsingular point is a discrete valuation ring. Let $a$ be a point on a curve $C$ and $\mmm_a$ be the (unique) maximal ideal of $\local{C}{a}$.

\begin{prop} \label{prop:curve-local-radical}
If $f \in \local{C}{a}$ is not identically zero on any of the irreducible components of $C$ containing $a$, then the radical of the ideal generated by $f$ in $\local{C}{a}$ is either $\local{C}{a}$ itself or $\mmm_a$.
\end{prop}

\begin{proof}
\Woutlog\ we may assume $C$ is affine. Then $f  = f_1/f_2$ for some $f_1, f_2 \in \kk[C],\ f_2(a) \neq 0$ (\cref{prop:affine-local}). If $f_1(a) \neq 0$, then $f$ is invertible in $\local{C}{a}$. Otherwise \cref{thm:dim-smaller} implies that $V(f_1) \subset C$ consists of finitely many points excluding $a$. Choose a polynomial $g$ which does not vanish at $a$ but vanishes at every other point of $V(f_1)$. For any $h \in \kk[C]$ such that $h(a) = 0$, the Nullstellensatz (\cref{thm:Hilbert-nulls}) implies that $gh \in \sqrt{f_1} \subseteq \kk[C]$. Since $g$ is invertible in $\local{C}{a}$, it follows that $h \in \sqrt{f} \subseteq \local{C}{a}$. This implies that $\sqrt{f} = \mmm_a$, as required.
\end{proof}

\begin{prop} \label{prop:curve-dvr}
Assume $C$ is nonsingular at $a$. Fix $t \in \mmm_a \setminus \mmm_a^2$ (that $\mmm_a \setminus \mmm_a^2$ is nonemepty is a consequence of \cref{thm:tangent-space}). 
\begin{enumerate}[resume]
\item \label{curve-dvr:m-principal} $\mmm_a$ is the principal ideal generated by $t$. 
\item \label{curve-dvr:nu-defn} Let $\nu: \local{C}{a} \to \zz_{\geq 0} \cup \{\infty\}$ be the map given by $g \mapsto \inf\{m \geq 0: t^m \in \langle g \rangle\}$. Then for all $g \in \local{C}{a} \setminus \{0\}$, 
\begin{enumerate}
\item \label{curve-dvr:u} $g = u t^{\nu(g)}$ for some unit $u \in \local{C}{a}$, 
\item \label{curve-dvr:uu} $g = c t^{\nu(g)} + g'$ for some $c \in \kk\setminus\{0\}$ and  $g' \in \mmm_a$ such that $\nu(g') > \nu(g)$.
\end{enumerate}
\item \label{curve-dvr:nu-dvr} $\nu$ extends to a discrete valuation on the field of fractions of $\local{C}{a}$ (recall that $\local{C}{a}$ is an integral domain due to \cref{prop:nonsingular-local-irreducible}) and its valuation ring is $\local{C}{a}$. 
\end{enumerate}
\end{prop}

\begin{proof}
Since $\dim(\tangent{C}{a}) = 1$, the image of $t$ generates $\mmm_a/\mmm_a^2$ over $\kk$ (assertion \eqref{mxa} of \cref{thm:tangent-space}). Assertion \eqref{curve-dvr:m-principal} then follows directly from \cref{Nakayama-generation}, which is a corollary of Nakayama's lemma (\cref{Nakayama}). For assertion \eqref{curve-dvr:nu-defn} pick $g \in \local{C}{a}$. \Cref{prop:curve-local-radical} implies that $\nu(g)$ is well-defined. Let $m := \nu(g)$. By definiton of $\nu$, $t^m = u'g$ for some $u \in \local{C}{a}$. If $u' \in \mmm_a$, then assertion \eqref{curve-dvr:m-principal} would imply that $u' = th$ for some $h \in \local{C}{a}$, which would in turn imply that $t^{m-1} = hg$ (since $\local{C}{a}$ is an integral domain), contradicting the minimality of $\nu(g)$. Therefore $u'$ is a unit in $\local{C}{a}$, proving \eqref{curve-dvr:u} with $u := u'^{-1}$. Let $c := u(a) \neq 0$. Then  $g = ut^m = ct^m + (u - c)t^m$. Since $u - c \in \mmm_a$, assertion \eqref{curve-dvr:m-principal} implies that $u-c \in \langle t \rangle$, which in turn implies that $\nu((u - c)t^m) > m$, proving \eqref{curve-dvr:uu}. Assertion \eqref{curve-dvr:nu-dvr} then follows in a straightforward way by extending $\nu$ to the field of fractions of $\local{C}{a}$ by defining $\nu(f/g) := \nu(f) - \nu(g)$ for $f, g \in \local{C}{a}$; we leave the details as an exercise. 
\end{proof}

If $C$ is nonsingular at $a$, \cref{prop:curve-dvr} implies that any $t \in \mmm_a \setminus \mmm_a^2$ is a {\em parameter} of the discrete valuation ring $\local{X}{a}$; we say that $t$ is a \index{Parameter!of a curve at a nonsingular point}{\em parameter} of $C$ at $a$.

\begin{cor} \label{prop:curve-nonsingular-morphextension}
Let $C$ be an irreducible curve and $f: C \dashrightarrow \pp^N$ be a rational map. Assume there is $C' \subseteq C$ such that $f|_{C'}$ is a morphism and $C$ is nonsingular at every point of $C \setminus C'$. Then $f$ extends to a morphism $C \to \pp^N$. 
\end{cor}

\begin{proof}
Note that $C \setminus C'$ is finite. Fix $a \in C \setminus C'$. It suffices to show that $f$ can be extended to a morphism on a neighborhood of $a$. \Woutlog\ we may assume that 
\begin{prooflist}
\item $C$ is an (affine) subvariety of $\kk^n$ with coordinates $(x_1, \ldots, x_n)$, and
\item there is an open neighborhood $U$ of $a$ in $C$ such that $f$ is a morphism from $U \setminus \{a\} \to \pp^N$ given by $x \mapsto [h_0(x): \cdots : h_N(x)]$. 
\end{prooflist}
Pick a parameter $t$ of $C$ at $a$. Then each $h_j = u_jt^{m_j}$ for some units $u_j \in \local{C}{a}$ and $m_j := \nu(u_j)$. Choose an open neighborhood $U'$ of $a$ in $U$ such that each $u_j$ is a regular function on $U'$. If $m := \min\{m_j\}_j$, then $f$ uniquely extends to $U': \to \pp^N$ given by $[u_0t^{m_0-m}: \cdots : u_Nt^{m_N - m}]$. 
\end{proof}

\begin{example}
Let $X$ be the image of $\phi: \kk \to \kk^3$ given by $t \mapsto (t^3, t^4, t^5)$. \Cref{exercise:no-to-yes-comeplete-intersection} shows that $X \setminus \{(0,0,0)\}$ is nonsingular. We now use \cref{prop:curve-nonsingular-morphextension} to show that $X$ is singular at $O := (0,0,0)$\footnote{To prove this directly using the definitions would require computation of $I(X)$ on an open neighborhood of $O$, which is a relatively complicated task.}. Indeed, $\phi$ is a birational map, and by the usual identification of $\kk$ with $\pp^1 \setminus V(x_0)$, we see that $\phi^{-1}: X \setminus \{O\} \dashrightarrow \pp^1$ is a well-defined morphism given by $(x,y,z) \mapsto [1: y/x]$. Now assume $X$ is nonsingular at $O$. Then $\phi^{-1}$ extends to a morphism $\psi: X \to \pp^1$. Since $\psi \circ \phi: \kk \to \pp^1$ is a morphism which is identity on $\kk \setminus \{0\}$, the Zariski-continuity of morphisms (\cref{exercise:morphism-zariski}) implies that it is identity everywhere on $\kk$. This implies that $\phi$ induces an isomorphism $\kk \cong X$, and consequently, an isomorphism $\local{X}{O} \cong \local{\kk}{0}$. However, it is clear that for any polynomial $f$ in $(x, y, z)$, the order of $f \circ \phi$ in $\kk(t)$ is $\geq 3$. This shows that $\phi^*$ can {\em not} be an isomorphism between $\local{X}{O}$ and $\local{\kk}{0}$, which gives the required contradiction. 
\end{example}

\Cref{prop:curve-nonsingular-morphextension} in general fails if $C \setminus C'$ has singular points - see \cref{exercise:curve-morph-nonextension} for an example. 

\subsection{Exercises}

\begin{exercise}\label{exercise:tangent-equation}
Prove identity \eqref{tangent-equation}. If $f_1, \ldots, f_s$ generate $I(X)$, then show that $\sum_{i=1}^N (x_i - a_i)\partialxifrac{f_j}(a)$, $j = 1, \ldots, s$, generate the ideal of polynomials vanishing on $\tangent{X}{a}$.
\end{exercise}

\begin{exercise} \label{exercise:tangent-k^n}
Show that $\kk^N$ is nonsingular everywhere. Given $a \in \kk^N$, compute $\tangent{\kk^N}{a}$.
\end{exercise}

\begin{exercise} \label{exercise:nonsingular<3}
Let $C = V(f) \subset \kk^2$, where $f$ is an irreducible polynomial of degree $1$ or $2$. Show that $C$ is everywhere nonsingular. [Hint: since $\kk$ is algebraically closed, a homogeneous polynomial of degree $2$ can be written as $ax^2 + by^2$, $a, b \in \kk$, after an appropriate change of coordinates on $\kk^2$. If $\deg(f) = 2$, then use this fact to reduce to the following cases: (1) $f = ax^2 + by^2 + c$, $a,b,c \neq 0$, and (2) $f = ax^2 + by$, $a, b \neq 0$.]
\end{exercise}

\begin{exercise} \label{exercise:singular3}
Let $C = V(f) \subset \kk^2$. If either $f = x^2 - y^3$ or $f = x^2 - y^2 + y^3$, show that the origin is the only singular point of $C$.
\end{exercise}

\begin{exercise} \label{exercise:poly-derivations}
Prove \cref{prop:poly-derivations}.
\end{exercise}

\begin{exercise} \label{exercise:locally-basic}
Prove \cref{prop:locally-basic}. [Hint: for assertion \eqref{locally-basic:sublocal} it suffices to show that the map $\local{X}{a} \to \local{Z}{a}$ given by restriction to $Z$ is surjective. Any open neighborhood of $a$ in $Z$ is of the form $U \cap Z$ for some open neighborhood of $a$ in $X$. Choose an open {\em affine} neighborhood $U'$ of $a$ in $X$ such that $U' \subseteq U$. Then $U' \cap Z$ is also affine and therefore all regular functions on $U' \cap Z$ are restrictions of regular functions on $U' \cap X$.]   
\end{exercise}

\begin{exercise}\label{exercise:locally-x-to-xy}
Let $\phi: \kk^2 \to \kk^2$ from \cref{example:constructible-0} given by $(x,y) \mapsto (x, xy)$. Given $a = (a_1, a_2) \in \kk^2$, consider the induced map $\phi^*: \local{\kk^2}{\phi(a)} \to \local{\kk^2}{a}$. Show that
\begin{enumerate}
\item $\phi^*$ is not surjective when $a$ is the origin.
\item $\phi^*$ is an isomorphism if $a_1 \neq 0$. [Hint: $x$ is invertible in $\local{\kk^n}{(b_1, b_2)}$ if $b_1 \neq 0$.]
\end{enumerate}
\end{exercise}

\begin{exercise} \label{exercise:local-isomorphism}
Given varieties $X, Y$ and points $a \in X$, $b \in Y$, show that the following are equivalent: 
\begin{enumerate}
\item $\local{X}{a} \cong \local{Y}{b}$ as $\kk$-algebras,
\item there are open neighborhoods $U$ of $a$ in $X$ and $V$ of $b$ in $Y$ such that $U \cong V$. 
\end{enumerate}
[Hint: The $(\Leftarrow)$ implication follows directly from \cref{prop:locally-basic}. For the $(\im)$ implication, suffices to consider the case that $X,Y$ are subvarieties respectively of $\kk^m$ with coordinates $(x_1, \ldots, x_m)$ and $\kk^n$ with coordinates $(y_1, \ldots, y_n)$. Given a $\kk$-algebra isomorphism $\Phi: \local{Y}{b} \to \local{X}{a}$ there are $g, f_1, \ldots, f_n \in \kk[x_1, \ldots, x_m]$ such that $g(a) \neq 0$ and $\Phi$ maps $y_j \mapsto f_j/g \in \local{X}{a}$. Then $\phi:x \mapsto (f_1(x)/g(x), \ldots, f_n(x)/g(x))$ is a morphism from $X \setminus V(g)$ to $Y$. Similarly, $\Phi^{-1}$ induces a morphism $\psi: Y \setminus V(q) \to X$ for some $q \in \kk[y_1, \ldots, y_n]$ such that $q(b) \neq 0$. Then $\psi \circ \phi$ and $\phi \circ \psi$ must be identity near respectively $a$ and $b$.]
\end{exercise}

\begin{exercise}\label{exercise:local-derivation}
Let $X$ be an affine variety, $a \in X$, and $D:\kk[X] \to \kk$ be a derivation centered at $a$. If $\mmm_a$ is the maximal ideal of $\local{X}{a}$ generated by polynomials vanishing at $a$, show that
\begin{enumerate}
\item $D$ extends to a linear map $\mmm_a \to \kk$ given by
\begin{align*}
D(f/g) := D(f)/g(a)\ \text{for all}\ f,g \in \kk[X]\ \text{such that}\ g(a) \neq 0
\end{align*}
\item $D(h) = 0$ for all $h \in \mmm_a^2$.
\end{enumerate}
\end{exercise}


\begin{exercise} \label{exercise:nonzero-derivative}
Let $f$ be an irreducible (nonzero) polynomial in $(x_1, \ldots, x_N)$. 
\begin{enumerate}
\item Show that there is $i$ such that $\partialxi{f}$ does {\em not} identically vanish on $V(f)$. [Hint: $I(V(f)) = \langle f \rangle$.]
\item Give examples to show that the preceding assertion may not hold if $f$ is not irreducible or if $\kk$ is not algebraically closed.
\end{enumerate}
\end{exercise}

\begin{exercise} \label{exercise:no-to-yes-comeplete-intersection}
Let $\phi: \kk \to \kk^3$ be the morphism given by $t \mapsto (t^3, t^4, t^5)$. \Cref{exercise:noncomplete-intersection} shows that $X := \phi(\kk)$ is a one dimensional subvariety of $\kk^3$. Let $U = \kk^3 \setminus V(x)$ (where $(x,y,z)$ are coordinates on $\kk^3$), so that$\kk[U] = \kk[x,y,z, 1/x]$.
\begin{enumerate}
\item Show that the ideal of $X \setminus V(x)$ in $\kk[U]$ is generated by $g_1 := x - (y/x)^3$, $g_2 := y - (y/x)^4$, and $g_3 := z - (y/x)^5$. [Hint: modulo the ideal generated by $g_1, g_2, g_3$, every polynomial $f$ in $(x, y, z)$ restricts to a polynomial $\bar f$ in $y/x$. Show that $f|_{X \cap U} \equiv 0$ if and only if $\bar f \equiv 0$.]
\item Show that $g_1$ is in the ideal generated by $g_2$ and $g_3$ in $\kk[U]$. Deduce that the ideal of $X \setminus V(x)$ in $\kk[U]$ is generated by $g_2$ and $g_3$.  
\item Show that $X \setminus V(x) \cong \kk \setminus \{0\}$; in particular, $X \setminus V(x)$ is nonsingular. 
\end{enumerate}
\end{exercise}

\begin{exercise} \label{exercise:curve-morph-nonextension}
This exercise shows that the conclusion of  \cref{prop:curve-nonsingular-morphextension} might fail if $\bar C \setminus C$ has singular points with ``more than one branch.'' Assume $\character(\kk) \neq 2$. Let $C$ be the curve on $\kk^2$ defined by the equation $y^3 = x(y^2 - 1)$. Let $[x_0: x_1: x_2]$ be homogeneous coordinates on $\pp^2$. The map $(x,y) \mapsto [1:x: y]$ identifies $\kk^2$ with the basic open subset $U_0 := \pp^2 \setminus V(x_0)$ of $\pp^2$. Let $\bar C$ be the closure of $C$ in $\pp^2$.
\begin{enumerate}
\item Show that $\bar C = V(x_2^3 - x_1(x_2^2 - x_0^2)) \subset \pp^2$ [Hint: use \cref{example:principal-closure}] and $\bar C \setminus C = \{[0:1:1], [0:1:0]\}$.
\item Show that the projection from $C$ to $y$-axis is one-to-one, and the inverse of this map extends to a morphism $\phi: \pp^1 \to \bar C$ which is generically one-to-one, and $\phi^{-1}([0:1:0])$ consists of two points.
\item Conclude that the conclusion of \cref{prop:curve-nonsingular-morphextension} fails with $C,\bar C$ and $f := y$.
\item Show that $O := [0:1:0]$ has an affine neighborhood in $\bar C$ isomorphic to the plane curve $v^3 = (v-w)(v+w)$. A drawing of this curve makes apparent the two ``branches'' (see \cref{branchion-0}) at $O$ with ``tangents'' $v - w = 0$ and $v + w = 0$ (see \cref{fig:node}).
\end{enumerate}
\end{exercise}

\section{Completion of the local ring at a point} \label{completion-section}
To study local properties of a variety $X$ near a point, sometimes one needs to pass to the ring of formal power series in affine coordinates at the point. We have seen an example of this in the proof of \cref{locally-prime}. In \cref{degree-section,mult-chapter} we study different notions of {\em multiplicites} at a point, and power series expansions in coordinates at the point play a fundamental role in our study. The usefulness of these computations depends on the fact that they do {\em not} depend on the chosen coordinates, i.e.\ the ``rings of formal power series associated to a point on a variety'' are isomorphic under local isomorphisms. One way to see this is through the theory of ``completions of local rings," which we describe now. Given an ideal $I$ of a ring $R$ and $f \in R$, consider the following property of subsets $S$ of $R$: 
\begin{align}
S \supseteq f + I^m\ \text{for some}\ m \geq 0
\label{I-adic-condition}
\end{align}
It is straightforward to check that there is a unique topology on $R$ in which a subset $S$ of $R$ is an open neighborhood of $f \in R$ if and only if satisfies condition \eqref{I-adic-condition} (\cref{exercise:I-adic}); this is called the {\em $I$-adic topology} on $R$. A {\em Cauchy sequence} in $R$ is a sequence of elements $(f_j)_{j \geq 0}$ of $R$ such that for any open neighborhood $U$ of $0$, there is $m \geq 0$ with the property that $f_i - f_j \in U$ for all $i,j \geq m$. Two Cauchy sequences $(f_i)_i$ and $(g_j)_j$ are {\em equivalent} if the sequence $(f_i - g_i)_i$ converges\footnote{Recall that on a topological space $X$, a sequence $(x_i)_i$ {\em converges} to a point $x$ if for every open neighborhood $U$ of $x$ in $X$ there is an integer $N$ such that $x_i \in U$ for all $i \geq N$.} to $0$ in $R$. The \index{Completion of a ring}{\em $I$-adic completion} $\hat R$ of $R$ is the set of equivalence classes of Cauchy sequences.

\begin{example} \label{example:0-completion}
If $I$ is the zero ideal, then all Cauchy sequences $(f_j)_{j \geq 0}$ are eventually constant, i.e.\ equivalent to a constant sequence of the form $(f,f, \ldots)$ for some $f \in R$, and moreover, $(f,f, \ldots)$ is equivalent to $(g, g,\ldots)$ if and only if $f = g$. It follows that $\hat R \cong R$. 
\end{example}

\begin{example} \label{example:1-completion}
On the other extreme, if $I = R$, then all Cauchy sequences are equivalent to $(0, 0,\ldots)$, so that $\hat R$ is the zero ring.
\end{example}

\begin{example}
\label{example:polynomial-completion}
Let $R := \kk[x_1, \ldots, x_n]$ and $I$ be the ideal of $R$ generated by all polynomials vanishing at some point $a = (a_1, \ldots, a_n) \in \kk^n$. Let $(f_j)_j$ in $R$ be a Cauchy sequence with respect to the $I$-adic topology. We treat each $f_j$ as a polynomial in $y_i := x_i - a_i$, $i = 1, \ldots, n$. For each $k \geq 0$, there is $M_k$ such that the degree (with respect to $(y_1, \ldots, y_n)$-coordinates) of $f_i - f_j$ is greater than $k$ for each $i, j \geq M_k$. In particular, with respect to $(y_1, \ldots, y_n)$-coordinates the homogeneous components $f_{j,k}$ of degree $k$ of all $f_j$ agree with each other whenever $j \geq M_k$; write $F_k := f_{M_k, k}$. It is straightforward to check that the map
\begin{align*}
(f_j)_j \mapsto \sum_k F_k
\end{align*}
induces a $\kk$-algebra isomorphism between $\hat R$ and $\kk[[x_1 - a_1, \ldots, x_n - a_n]]$.
\end{example}

\begin{example} \label{example:polynomial-local-completion}
We have seen (in the discussion preceding \cref{hat-intersection}) that the local ring $\local{\kk^n}{a}$ of $\kk^n$ at a point $a = (a_1, \ldots, a_n)$ can be treated as a subring of $\kk[[x_1 - a_1, \ldots, x_n - a_n]]$. We leave it as an exercise (\cref{exercise:polynomial-local-completion}) to check that $\hatlocal{\kk^n}{a} \cong \kk[[x_1 - a_1, \ldots, x_n - a_n]]$; in particular, it follows from \cref{example:polynomial-completion} that the completion of the local ring at a point of $\kk^n$ with respect to its maximal ideal is isomorphic to the completion of the coordinate ring of $\kk^n$ with respect to the maximal ideal at a point.
\end{example}

\begin{prop} \label{prop:complete-inclusion}
Let $\hat R$ be the completion of a ring $R$ with respect to an ideal $I$. Given $f \in R$, write $\hat f$ for the equivalence class of the constant sequence $(f, f, \ldots)$. Let $\phi : R \to \hat R$ be the map which sends $f \mapsto \hat f$.
\begin{enumerate}
\item $\ker \phi = \bigcap_{m \geq 1} I^m$.
\item \label{complete-inclusion:local} If $R$ is a Noetherian local ring, and $I$ is a proper ideal of $R$, then $\phi$ is injective.
\end{enumerate}
\end{prop}

\begin{proof}
Since the zero element of $\hat R$ is the equivalence class of $(0, 0, \ldots)$, the first assertion is immediate from the definition of completion. The second assertion then follows from \cref{thm:Krull}.
\end{proof}

\begin{rem}
\Cref{example:1-completion} shows that assertion \eqref{complete-inclusion:local} of \cref{prop:complete-inclusion} does not hold if $I = R$.
\end{rem}

If $X$ is a variety and $a \in X$, then $\local{X}{a}$ has a unique maximal ideal $\mmm_a$, namely the ideal generated by polynomials vanishing at $a$ (\cref{prop:locally-basic} ). We write $\hatlocal{X}{a}$ for the $\mmm_a$-adic completion of $\local{X}{a}$. The ring $\hatlocal{X}{a}$ captures ``very local'' information about $X$ at $a$. \Cref{example:polynomial-local-completion} shows that if $X$ is the affine space, then $\hatlocal{X}{a}$ is the ring of formal power series expansions centered at $a$; in general it is a quotient of a power series ring:

\begin{prop} \label{complete-quotient}
Let $a = (a_1, \ldots, a_N)$ be a point of a subvariety $X$ of $\kk^N$, and let $\hat R := \kk[[x_1 - a_1, \ldots, x_N - a_N]]$. Then $\hatlocal{X}{a} \cong \hat R/(I(X)\hat R)$.
\end{prop}

\begin{proof}
Since $\local{X}{a} \cong \local{\kk^N}{a}/(I(X)\local{\kk^N}{a})$ (assertion \eqref{locally-basic:sublocal}  of \cref{prop:locally-basic}), the result follows from \cref{example:polynomial-local-completion} and the exactness of completion (\cref{thm:exactly-complete}).
\end{proof}

\begin{thm} \label{smooth-completion}
Let $a$ be a nonsingular point of an irreducible variety $X$ of dimension $n$. Then
\begin{enumerate}
\item \label{complete-isomoprhism} $\hatlocal{X}{a} \cong \kk[[x_1, \ldots, x_n]]$.
\end{enumerate}
Assume $X$ is a subvariety of $\kk^N$ with coordinates $(x_1, \ldots, x_N)$. Pick $g_1, \ldots, g_n \in \kk[x_1, \ldots, x_N]$ such that $g_i(a) = 0$ for each $i$, and $(\partialxiifrac{g_i}{1}(a), \ldots, \partialxiifrac{g_i}{N}(a))$, $i = 1, \ldots, n$, generate $\tangent{X}{a}$ as a vector space over $\kk$ (see assertion \eqref{txa} of \cref{thm:tangent-space}). Then
\begin{enumerate}[resume]
\item \label{complicit-isomorphism} it is possible to choose the isomorphism from assertion \eqref{complete-isomoprhism} such that $\hat {\bar{g_i}} \to x_i$, $i = 1, \ldots, n$ (where $\bar g_i := g_i|_X \in \kk[X]$ and $\hat {\bar{g_i}}$ are defined as in \cref{prop:complete-inclusion}).
\end{enumerate}
\end{thm}

\begin{proof}
Taking an appropriate open neighborhood of $a$ we may assume that $X$ is an irreducible subvariety of $\kk^N$, and $I(X)$ is generated by $f_1, \ldots, f_{N - n}$ (the second property can be ensured due to \cref{local-equations}). By a change of coordinates if necessary, we may assume $a$ is the origin in $\kk^N$. Then $\hatlocal{X}{a} \cong \hat R/(I(X)\hat R)$, where $\hat R := \kk[[x_1, \ldots, x_N]]$ (\cref{complete-quotient}). The nonsingularity of $X$ at $a$ implies that the linear parts of the $f_j$ are linearly independent over $\kk$, and we may choose $g_1, \ldots, g_n \in \kk[X]$ which satisfy the hypothesis of assertion \eqref{complicit-isomorphism} (namely, take polynomials $g_1, \ldots, g_n$ which vanish at $a$ and are such that the linear parts of $f_1, \ldots, f_{N-n}, g_1, \ldots, g_n$ are linearly independent over $\kk$). It then follows from \cref{cor:hat-prime} that $\hat R/(I(X)\hat R)$ is isomorphic to the ring of power series in $n$-variables over $\kk$ via an isomorphism that maps $\hat {\bar {g_i }} \to x_i$, $i = 1, \ldots, n$.
\end{proof}

\begin{cor} \label{smooth-1-completion}
Let $a$ be a nonsingular point of an irreducible curve $C$, and $t$ be a parameter of $\local{C}{a}$. Then $\hatlocal{C}{a} \cong \kk[[t]]$.
\end{cor}

\begin{proof}
Pick $g_1$ as in \cref{smooth-completion}. By definition of a parameter, $g_1 = ut^k$ for some $k \geq 0$ and a unit $u$ of $\local{C}{a}$. The properties of $g_1$ implies that $k = 1$. \Cref{prop:curve-dvr} then implies that $g_1 = ct + g'_1$ for some $g'_1 \in \mmm_a^2$ (where $\mmm_a$ is the maximal ideal of $\local{C}{a}$), it follows that $\partialxi{g'_1} = 0$ and $\partialxi{g_1} = c\partialxi{t}$ for each $i = 1, \ldots, N$ (where $(x_1, \ldots, x_N)$ are coordinates on an affine open neighborhood of $a$), i.e.\ the linear parts of $g_1$ and $t$ are proportional. The corollary then follows from the arguments of the proof of \cref{smooth-completion}.
\end{proof}

\subsection{Exercises}
\begin{exercise} \label{exercise:I-adic}
Let $R$ be a ring, $f \in R$, and $I$ be an ideal of $R$. Let $\{S_j\}_{j \in \scrJ}$ be a collection of subsets of $R$ such that each $S_j$ satisfies \eqref{I-adic-condition}. Show that 
\begin{enumerate}
\item $\bigcup_{j \in \scrJ} S_j$ satisfies \eqref{I-adic-condition}. 
\item If $\scrJ$ is finite, then $\bigcap_{j \in \scrJ} S_j$ satisfies \eqref{I-adic-condition}.  
\end{enumerate}
\end{exercise}

\begin{exercise} \label{exercise:polynomial-local-completion}
Let $a = (a_1, \ldots, a_n) \in \kk^n$, and $\mmm_a$ be the unique maximal ideal of $\local{\kk^n}{a}$ generated by polynomials vanishing at $a$. Show that the $\mmm_a$-adic completion $\hatlocal{\kk^n}{a}$ of $\local{\kk^n}{a}$ is isomorphic as a $\kk$-algebra to $\kk[[x_1 - a_1, \ldots, x_n - a_n]]$.
\end{exercise}

\begin{exercise} \label{exercise:finite-quotient}
Let $a \in X$ and $\qqq$ be an ideal of $\local{X}{a}$ such that $\local{X}{a}/\qqq\local{X}{a}$ is a finite dimensional vector space over $k$. Show that $\local{X}{a}/\qqq\local{X}{a} \cong \hatlocal{X}{a}/\qqq\hatlocal{X}{a}$. [Hint: if $\mmm$ is the maximal ideal of $\local{X}{a}$, then $\qqq \supset \mmm^q \local{X}{a}$ for some $q > 0$. Then it follows from the definition of completion that
$%
	(\local{X}{a}/\mmm^q) /(\qqq\local{X}{a}/\mmm^q)
	\cong (\hatlocal{X}{a}/\mmm^q\hatlocal{X}{a}) /(\qqq\hatlocal{X}{a}/\mmm^q\hatlocal{X}{a})
$]
\end{exercise}

%
%

\section{Degree of a dominant morphism} \label{degree-section}
Let $\phi:X \to Y$ be a dominant morphism of irreducible varieties. \Cref{fiber-dimension} implies that
\begin{prooflist}
\item if $\dim(X) > \dim(Y)$, then $|\phi^{-1}(y)|$ is infinite for each $y$ in a dense open subset of $Y$, and
\item \label{case:surjective-fiber-eq-dim} if $\dim(X) = \dim(Y)$, then $|\phi^{-1}(y)|$ is finite for each $y$ in a dense open subset of $Y$. 
\end{prooflist}
Whenever case \ref{case:surjective-fiber-eq-dim} arises for a continuous (or differentiable) map in topology, it turns out that for ``almost all'' $y \in Y$, some ``measure''\footnote{E.g.\ $|\phi^{-1}(y)|$, or $|\phi^{-1}(y)|$ modulo $2$, etc.} of the number of elements in $\phi^{-1}(y)$ is constant, and that number is called the ``degree'' of $\phi$. In this section we will see that this remains true for morphisms of algebraic varieties as well. Indeed, since $\phi$ is a dominant, it induces an inclusion $\kk(Y) \into \kk(X)$ (\cref{prop:star-injective}). In case \ref{case:surjective-fiber-eq-dim} $\kk(Y)$ is a {\em finite extension}\footnote{See \cref{algebraic-extension-section,separable-section} for a discussion of field extensions and related notions.} of $\kk(X)$; the \index{Degree!of a morphism}{\em degree} $\deg(\phi)$ of $\phi$ is the degree $[\kk(X): \kk(Y)]$ of the induced extension of fields. 

\begin{example} \label{example:separable-morphism-0}
Let $\phi: \kk \to \kk$ be the morphism $x \mapsto x^d$, $d > 0$. If $t$ is the coordinate on the target, then the induced extension $\kk(t) \into \kk(x)$ is given by $t \mapsto x^d$, so that $\deg(\phi) = [\kk(x): \kk(x^d)] = d$. If $p := \character(\kk) = 0$, then indeed $|\phi^{-1}(a)| = d = \deg(\phi)$ for all $a \in \kk \setminus \{0\}$. On the other hand, if $p > 0$ and $d = qp^k$ where $k \geq 1$ and $q$ is relatively prime to $p$, then $|\phi^{-1}(a)| = q < d$ for each $a \in \kk \setminus \{0\}$. However, in this case the extension $\kk(x)/\kk(t)$ is {\em not} separable, and the {\em separable degree} of $\kk(x)/\kk(t)$ is precisely $q$ (\cref{example:separable-1}).
\end{example}

Motivated by \cref{example:separable-morphism-0} we define the \index{Separable!degree!of a morphism} {\em separable degree} $\degsep(\phi)$ (resepectively, \index{Inseparable!degree!of a morphism} {\em inseparable degree} $\deginsep(\phi)$) for the separable (respectively, inseparable) degree of the field extesion $\kk(X)/\kk(Y)$ induced by $\phi$. Note that
\begin{itemize}
\item $\deg(\phi) = \degsep(\phi) \deginsep(\phi)$, and
\item if $\character(\kk) = 0$, then $\deginsep(\phi) = 1$ and $\deg(\phi) = \degsep(\phi)$. 
\end{itemize}
\Cref{degree-morphism} below states that if $\phi:X \to Y$ is a dominant morphism between irreducible varieties of the same dimension, then for all $y$ in a dense open subset of $Y$,
\begin{itemize}
\item $|\phi^{-1}(y)| = \degsep(\phi)$, and 
\item for each $x \in \phi^{-1}(y)$, the ``multiplicity'' of $\phi$ at $x$ is $\deginsep(\phi)$, so that 
\item the sum over all $x \in \phi^{-1}(y)$ of the multiplicity of $\phi$ is precisely $\deg(\phi)$. 
\end{itemize}
Before we state and prove \cref{degree-morphism} we give some applications. 

\begin{example} \label{example:distinct-roots}
In \cref{exercise:cayley-hamilton} we used the following fact: ``Given an algebraically closed field $\kk$ and a positive integer $n$, there is a nonempty Zariski open subset $U$ of $\kk^{n+1}$ such that for each $(c_0, \ldots, c_n) \in U$, the polynomial $f := c_0\lambda^n + c_1 \lambda^{n-1} + \cdots + c_n$ has $n$ distinct roots in $\kk$.'' We now prove this fact. Indeed, let $X$ be the hypersurface $V(f)$ on the affine space $\kk^{n+2}$ with coordinates $(c_0, \ldots, c_n, \lambda)$. It is straightforward to check that $f$ is irreducible as a polynomial in $(c_0, \ldots, c_n, \lambda)$, so that $X$ is irreducible. The arguments of \cref{prop:pure-dimension-1-0} implies that $\dim(X) = n+1$ and $1, \lambda, \ldots, \lambda^{n-1}$ is a basis of $\kk(X)$ over $\kk(c_0, \ldots, c_n)$, so that  $[\kk(X):\kk(c_0, \ldots, c_n)] = n$. Since $\kk$ is algebraically closed, the projection $\pi: X \to \kk^{n+1}$ in $(c_0, \ldots, c_n)$-coordinates is dominant. On the other hand, the derivative of $f$ with respect to $\lambda$ is $c_{n-1} + 2c_{n-2} \lambda + \cdots$, which is {\em not} identically zero in $\kk(c_0, \ldots, c_n)[\lambda]$. Therefore $\kk(X)$ is separable over $\kk(c_0, \ldots, c_n)$ (\cref{prop:separable-condition}). It follows that $\degsep(\pi) = \deg(\pi) = n$ and \cref{degree-morphism} implies that $|\pi^{-1}(c_0, \ldots, c_n)| = n$ for all $(c_0, \ldots, c_n)$ on a dense open subset of $\kk^{n+1}$, which proves the ``fact.''
\end{example}

\begin{example}[Degree of a projective variety]
Let $X$ be a subvariety of $\pp^n$. If $d := \dim(X)$, we will show that for ``almost all'' $(n-d)$-dimensional linear subspaces $L$ of $\pp^n$, the number of points in the intersection $L \cap X$ (counted with appropriate multiplicity) is constant (this number is called the \index{Degree!of a subvariety of the projective space}{\em degree} of $X$). Indeed, denote the homogeneous coordinates of $\pp^n$ by $x := [x_0, \ldots, x_n]$, and consider another $d$ copies of $\pp^n$ with homogeneous coordinates $\xi^i := [\xi^i_0: \cdots : \xi^i_n]$, $i = 1, \ldots, d$. Consider the subset $Z$ of $X \times ( \pp^n)^d$ (where $(\pp^n)^d$ is the $d$-fold Segre product $\pp^n \times \cdots \times \pp^n$) consisting of all $(x, \xi^1, \ldots, \xi^n)$ such that $\sum_j \xi^i_j x_j = 0$, $i = 1, \ldots, d$. It is straightforward to check that $Z$ is Zariski closed in $X \times ( \pp^n)^d$. Consider irreducible components $Z_k$ of $Z$ such that the projection $\pi: Z \to (\pp^n)^d$ in $(\xi^1, \ldots, \xi^d)$-coordinates maps $Z_k$ dominantly to $(\pp^n)^d$. Since $Z$ is complete, any such $Z_k$, if exists, must get mapped {\em surjectively} by $\pi$. \Cref{exercise:noempty-intersection-hypersurface} implies that such $Z_k$ exists, and in addition, there are $(\xi^1, \ldots, \xi^d)$ such that $|(\pi|_{Z_k})^{-1}(\xi^1, \ldots, \xi^d)| < \infty$. It then follows due to \cref{fiber-dimension} that $\dim(Z_k) = \dim((\pp^n)^d) = nd$, so that \cref{degree-morphism} applies, and shows that for all $(\xi^1, \ldots, \xi^d)$ in a dense open subset of $(\pp^n)^d$, the number of elements in $\pi^{-1}(\xi^1, \ldots, \xi^d)$ counted with appropriate multiplicity is precisely $\sum_k \deg(\pi|_{Z_k})$. 
\end{example}

\begin{thm} \label{degree-morphism}
Let $\phi: X \to Y$ be a dominant morphism between irreducible varieties of same dimension. Then there is a nonempty Zariski open subset $U$ of $Y$ such that for each $y \in U$,
\begin{enumerate}
\item $Y$ is nonsingular at $y$;
\item $|\phi^{-1}(y)| = \degsep(\phi)$;
\item for each $x \in \phi^{-1}(y)$
\begin{enumerate}
\item $X$ is nonsingular at $x$,
\item $\dim_\kk(\hatlocal{X}{x}/\mmm_y\hatlocal{X}{x}) = \deginsep(\phi)$, where $\mmm_y$ is the maximal ideal of $\local{Y}{y}$ and $\hatlocal{X}{x}$ is the completion of $\local{X}{x}$ with respect to its maximal ideal;
\end{enumerate}
\item in particular
\begin{align*}
\sum_{x \in \phi^{-1}(y)} \dim_\kk(\hatlocal{X}{x}/\mmm_y\hatlocal{X}{x}) = \deg(\phi)
\end{align*}
\end{enumerate}
\end{thm}

\begin{proof}
We may assume \woutlog\ that $X$ and $Y$ are affine and nonsingular. Let $L$ be the separable closure of $\kk(Y)$ in $\kk(X)$. Pick regular functions $z_1, \ldots,  z_k$ on $X$ such that $L$ is the field of fractions of $T := \kk[Y][z_1, \ldots, z_k]$. There is (up to isomorphism) a unique affine variety $Z$ with coordinate ring $T$ (\cref{exercise:algebra-to-affine-variety}). The chain of inclusions $\kk[Y] \into T \into \kk[X]$ induces a factorization of $\phi$ of the form:
\begin{align*}
X \xrightarrow{\phi_i} Z \xrightarrow{\phi_s} Y \subseteq \kk^N
\end{align*}
By \cref{thm:primitive-element-infinite} there is $g \in T$ which generates $L$ over $\kk(Y)$. We can factor $\phi_s$ as:
\begin{align*}
Z \xrightarrow{\psi} Y \times \kk \xrightarrow{\pi} Y
\end{align*}
where $\psi$ maps $z \mapsto (\phi_s(z), g(z))$ and $\pi$ is the projection onto $Y$. Let $G(y,t) = \sum_{i=0}^{d_s}a_i(y) t^{d_s-i} \in \kk[Y][t]$, where $t$ is an indeterminate and $d_s := \degsep(\phi)$, be the minimal polynomial of $g$ over $\kk(Y)$. The separability of $g$ over $\kk(Y)$ implies that $(\partial G/\partial t)|_{t= g}$ is a non-zero element of $\kk[Z]$. Let $U_0$ be a nonempty Zariski open subset of $Y$ contained in $Y\setminus (V(a_0) \cup \phi_s(V( (\partial G/\partial t)|_{t=g})))$ [why does such $U_0$ exist?] and $U'_0 := \{(y,t) : \sum_{i=0}^{d_s}a_i(y) t^{d_s-i}  = 0\}  \subset U_0 \times \kk$. Then $U'_0$ is irreducible (since $G$ is irreducible in $\kk(Y)[t]$), and $\psi$ induces a birational map from $Z$ to $U'_0$ [why?]. Let $Y_0$ be a nonempty Zariski open subset of $U_0$ such that $\psi$ induces an isomorphism $\phi_s^{-1}(Y_0) \cong Y'_0 := (\pi|_{U'_0})^{-1}(Y_0)$. Let $y_0 \in Y_0$. Then $Y'_0$ contains $(y_0,t_0)$ for all the roots $t_0$ of $G (y_0,t)$. Let $z_0 := \psi^{-1}(y_0,t_0)$. Since $(\partial G/\partial t) (y_0,t_0)$ equals $(\partial G/\partial t)|_{t=g}$ evaluated at $z_0$, it follows that $(\partial G/\partial t) (y_0,t_0) \neq 0$ for every root $t_0$ of $G(y_0,t)$. Consequently, $|\pi^{-1}(y_0)| = \deg(G(y_0,t)) = d_s = \degsep(\phi)$.

\begin{proclaim}
For each $(y_0,t_0) \in \pi^{-1}(y_0)$, $\dim_\kk(\hatlocal{Y'_0}{(y_0,t_0)}/\mmm_{y_0}\hatlocal{Y'_0}{(y_0,t_0)}) =1$.
\end{proclaim}

\begin{proof}
Pick $(y_0,t_0) \in \pi^{-1}(y_0)$. The image of $G(y,t)$ in $\hatlocal{\kk^{N+1}}{(y_0,t_0)} = \kk[[y_1-y_{0,1}, \ldots, y_N -y_{0,N},t-t_0]]$ (where $(y_{0,1}, \ldots, y_{0,N})$ are coordinates of $y_0$ in $\kk^N$) is
\begin{align*}
 G(y_0, t_0) + \sum_{j=1}^N (y_j - y_{0,j})\frac{\partial G}{\partial y_j} (y_0,t_0)
 	+ (t-t_0) \frac{\partial G}{\partial t}(y_0,t_0)  + \hot
\end{align*}
where $\hot$ denotes terms with order (in $(y-y_0,t-t_0)$) greater than one. Since $G(y_0,t_0) = 0$ and $(\partial G/\partial t) (y_0,t_0) \neq 0$, \cref{thm:linear-power} implies that $t - t_0$ is in the ideal of $\kk[[y_1 -y_{0,1}, \ldots, y_N -y_{0,N},t-t_0]]$ generated by $G(y,t)$ and $y_j-y_{0,j}$, $j = 1, \ldots, N$. Since $\hatlocal{Y'_0}{(y_0,t_0)}$ is the quotient of $\kk[[y_1 -y_{0,1}, \ldots, y_N-y_{0,N},t-t_0]]$ modulo the ideal generated by $G(y,t)$ (\cref{complete-quotient}), it follows that $t - t_0$ is in the ideal of $\hatlocal{Y'_0}{(y_0,t_0)}$ generated by the $y_j - y_{0,j}$, which implies the claim.
\end{proof}

Note that the above claim and the sentence preceding it proves \cref{degree-morphism} in the case that $\kk(X)$ is separable over $\kk(Y)$, in particular when $p := \character(\kk) = 0$. It remains to consider the case that $p > 0$ and $\kk(X)$ is not separable over $\kk(Y)$. Pick $x_1, \ldots, x_q \in \kk[X]$ such that $\kk[X] = T[x_1, \ldots, x_q]$. Set $T_0 := T$ and $T_j := T_{j-1}[x_j]$ for $j = 1, \ldots, q$. For each $j$, let $X_j$ be the unique affine variety with coordinate ring $T_j$. Note that each $X_j$ is irreducible (since $T_j$ is an integral domain). The inclusions $T_{j-1} \into T_j$ induces a factorization of $\phi_i:X \to Z$ as follows:
\begin{align*}
X = X_q \xrightarrow{\phi_{i,q}} X_{q-1} 
	\xrightarrow{\phi_{i, q-1}} 
	\cdots 
	\xrightarrow{\phi_{i, 2}} X_1 
	\xrightarrow{\phi_{i, 1}} X_0 
	= Z
\end{align*}
The minimal equation of $x_j$ over the field $L_{j-1}$ of fractions of $T_{j-1}$ is of the form $a_{j,0}x_j^{p^{e_j}} - a_{j,1} = 0$ for some $a_{j,0}, a_{j,1} \in T_{j-1}$ (\cref{prop:sep-insep-deg}). It follows that $L_j$ is generated by $1, x_j, \ldots, x_j^{p^{e_j}-1}$ as a vector space over $L_{j-1}$; in particular, $[L_j: L_{j-1}] = p^{e_j}$. Consequently,
\begin{align}
\deginsep(\phi) = [L_q: L_0] = p^{\sum_{j=1}^q e_j}
\label{eqn:degree-morphism-insep}
\end{align}
Choose a nonempty open affine subset $W_0$ of $X_0$ such that 
\begin{prooflist}
\item $W_0$ is nonsingular, 
\item $a_{1, 0}(x) \neq 0$ for each $x \in W_0$,
\item $W_j := (\phi_{i, 1} \circ \cdots \circ \phi_{i,j})^{-1}(W_0)$ is nonsingular for each $j$,
\item $a_{j + 1, 0}(x) \neq 0$ for each $x \in W_j$. 
\end{prooflist}
Then each $W_j$ is isomorphic to the hypersurface $V(x_j^{p^{e_j}} - a_{j,1}/a_{j,0})$ of $W_{j-1} \times \kk$, and $\phi_j|_{W_j}$ is simply the restriction of the projection $W_{j-1} \times \kk \to W_{j-1}$. It follows that $\phi_j|_{W_j}$ is one-to-one for each $j$, and consequently so is $\phi_i|_{W_q}: W_q \to W_0$. Fix $z_0 \in W_0$ and $z_q :=(\phi_i|_{W_q})^{-1}(z_0) \in W_q$. Due to \eqref{eqn:degree-morphism-insep} in order to complete the proof of \cref{degree-morphism} it suffices to show that
\begin{align}
\dim_\kk(\hatlocal{W_q}{z_q}/\mmm_{z_0}\hatlocal{W_q}{z_q}) = p^{\sum_{j=1}^q e_j}
\label{eqn:inseparable-claim}
\end{align}
Choose coordinates $(w_1, \ldots, w_m)$ on $W_0$ such that $z_0$ becomes the origin on $\kk^m$, $m\geq 1$. Replacing each $x_j$ by $x_j - c_j$ for some appropriate $c_j \in \kk$ if necessary, we may in addition assume that $x_j$ vanishes at $z_j  := (\phi_{i, 1} \circ \cdots \circ \phi_{i,j})^{-1}(z_0)$ for each $j$. This implies that 
\begin{align}
\begin{split}
\hatlocal{W_0}{z_0} 
	&\cong R_0/\ppp_0 R_0 \\
\hatlocal{W_j}{z_j} 
	&\cong R_0[[x_1, \ldots, x_j]]/
	\langle \ppp_0, x_1^{p^{e_1}} - a_{1,1}/a_{1,0}, \ldots, x_j^{p^{e_j}} - a_{j,1}/a_{j,0} \rangle
\end{split}
\label{eqn:inseparable-completion-steps}
\end{align}
where $R_0 := \kk[[w_1, \ldots, w_m]]$ and $\ppp_0 \subseteq \kk[w_1, \ldots, w_m]$ is the ideal of polynomials vanishing on $W$ (\cref{complete-quotient}). 

\begin{proclaim} \label{claim:degree-morphism-inseparable-basis}
Each element of $\hatlocal{W_q}{z_q}$ can be represented by a linear combination of $\scrG := \{\prod_j x_j^{i_j}: 0 \leq i_j < p^{e_j}\}$ with coefficients in $R_0$. 
\end{proclaim}

\begin{proof}
Pick $\rho \in \hatlocal{W_q}{z_q}$ and a power series $f$ in $(w_1, \ldots, w_m, x_1, \ldots, x_q)$ which represetnts $\rho$. Replacing $(x_q)^{ip^{e_q}}$ by $(a_{q,1}/a_{q,0})^i$ for each $i$ yields a power series $f_1$ such that all powers of $x_q$ in $f_1$ are smaller than $p^{e_q}$, and $f_1$ also represents $\rho \in \hatlocal{W_q}{z_q}$. Continuing this process with $x_{q-1}$ and so on yields a power series as claimed. 
\end{proof}

\Cref{claim:degree-morphism-inseparable-basis} implies that $\scrG$ spans $\hatlocal{W_q}{z_q}/\mmm_{z_0}\hatlocal{W_q}{z_q}$ over $\kk$. On the other hand, using \eqref{eqn:inseparable-completion-steps} it is straightforward to check that the elements of $\scrG$ are linearly independent over $\kk$ in $\hatlocal{W_q}{z_q}/\mmm_{z_0}\hatlocal{W_q}{z_q}$. Therefore $\scrG$ is a basis of $\hatlocal{W_q}{z_q}/\mmm_{z_0}\hatlocal{W_q}{z_q}$ over $\kk$. Since $|\scrG| = p^{\sum_j e_j}$, this completes the proof of \eqref{eqn:inseparable-claim} and consequently the theorem. 
\end{proof}

\chapter{$^*$Intersection multiplicity}  \label{mult-chapter}
\section{Introduction} \label{sec:intro-mult-intro}
\footnote{The asterisk in the chapter name is to indicate that most of the material of this section might be skipped in the first reading and/or in a first course of algebraic geometry. Only a small part of \cref{toric-intro} uses the results of this chapter. For the proof of Bernstein's theorem and its applications one would mainly need \cref{cor:Macaulay,order-curve,zero-sum,finite-mult-determinacy,int-mult-curve,mult-deformation-0,mult-deformation-global} - which might be explained without proof in a first course of algebraic geometry.}In this chapter we define the intersection multiplicity of $n$ hypersurfaces at a point on a nonsingular variety $X$ of dimension $n$, and prove some of its basic properties. As \cref{fig:parabola-tangent2} suggests, nontrivial considerations prop up even in the intersection of a parabola and a line. However, we do have a natural candidate for the intersection multiplicity, namely if $f_j$ are regular functions on $X$, then for each $a \in \bigcap_j \{f_j = 0\}$, we can consider the ``multiplicity'' at $a$ of the morphism $X \to \kk^n$ given by $x \mapsto (f_1(x), \ldots, f_n(x))$ suggested by \cref{degree-morphism}, i.e.\ the quantity
\begin{align}
\dim_\kk(\hatlocal{X}{a}/\langle f_1, \ldots, f_n \rangle)
\label{eqn:int-mult-motivation}
\end{align}
This is indeed the definition we are going to use (see \cite[Section 3.3]{fulturve} for a wonderful axiomatic motivation for this definition). However, in \crefrange{bkk-chapter}{milnor-chapter} we would in addition need to use a method of computing the intersection multiplicity via ``parametrization.'' Consider e.g.\ the case of two plane curves $f = 0$ and $g = 0$ on $\kk^2$. Their intersection multiplicity at the origin, as given by \eqref{eqn:int-mult-motivation}, is the dimension (as a vector space) over $\kk$ of the quotient of the power series ring $\kk[[x,y]]$ by the ideal generated by $f,g$. The ``parametric'' procedure on the other hand is as follows: find a ``parametrization'' $\phi(t)$ of the curve $g = 0$ such that $\phi(0) = 0$, and compute the order of $f(\phi(t))$, which measures how fast $f$ is vanishing along the curve $g = 0$. 
\begin{center}
\begin{figure}[htb]
\begin{tikzpicture}[
    thick,
    dot/.style = {
      draw,
      fill,
      circle,
      inner sep = 0pt,
      minimum size = 3pt
    },
    xscale =0.66,
    yscale=0.33]

  \def\x{2.7}
  \draw [\colorP] (-\x,\x^2) parabola [bend pos=0.5] bend +(0,-\x^2) +(2*\x,0);
   \node [\colorP, right] at (-\x+0.1,\x^2-0.1) {
   		\small
   		$y - x^2 = 0$
   	};
   	\def\ybottom{-0.9}
   	\draw (-2, 0) -- (2, 0);
    \def\ymiddle{0}
    \def\xright{2.7}
   \foreach \x in {1.2, 2.4} {
   		\draw ({\x*(\ybottom - \ymiddle)/(\x^2 - \ymiddle)}, \ybottom)
  		   		-- (\xright,{\xright*(\x^2 - \ymiddle)/\x + \ymiddle});
	}
	\def\x{2.7}
	\def\y{2.7*2.7}
	\def\txtwdt{6.3cm}
	\node [below right, text width= \txtwdt,align=justify] at (\x,\y) {
 		\small
			As secants approach the tangent at $O$ more and more closely, both of the two points of intersection move arbitrarily close to $O$.
		};
\end{tikzpicture}
\caption{A tangent line intersects a parabola at a point with multiplicity two} \label{fig:parabola-tangent2}
\end{figure}
\end{center}

For example, consider the situation of \cref{fig:parabola-tangent2}, i.e.\ $f = y - x^2$ and $g = y - mx$. Then $\phi(t) := (t, mt)$ parametrizes the line $g = 0$, and $f(\phi(t)) = mt - t^2$. Consequently $\ord_t(f(\phi(t))) = 1$ if $m \neq 0$. If $m = 0$, i.e.\ $g = 0$ is a horizontal line, then $\ord_t(f(\phi(t))) = 2$, as expected. However, to use this approach in practice one needs to define order, parametrization etc.\ even in the case that the equations are ``not reduced''; consider e.g.\ the case that $f = y - x^2$ and $g = (y - mx)^2$. The algebraic quantity \eqref{eqn:int-mult-motivation} gives the expected answer (which is $2$ if $m \neq 0$, and $4$ if $m = 0$), but what is geometrically the object $(y - mx)^2 = 0$? The underlying space is still the same line $y = mx$, but the defining equation is different, and the ``coordinate ring'' $\kk[x,y]/\langle (y - mx)^2 \rangle$ is ``non-reduced,'' (i.e.\ it has a nonzero nilpotent, namely the image of $y - mx$). In particular, in order to make the geometric approach more generally applicable one needs to
\begin{enumerate}
\item build a theory of ``non-reduced varieties,'' and
\item \label{task:order} define the notion of order at a point of a ``non-reduced curve.''
\end{enumerate}
In \cref{subscheme-section} we introduce the notion of ``closed subschemes'' which are the correct candidates for ``non-reduced subvarieties,'' and in \cref{reduced-noncurve-section} we extend the notion of order to non-reduced curves. In \cref{mult-point-section,complete-mult-section} we apply these notions to the study of intersection multiplicity.

%
\section{Closed subschemes of a variety} \label{subscheme-section}
\subsection{Closed subschemes of an affine variety} \label{subscheme-ssection-affine}
Let $X$ be an affine variety and $\qqq$ be a (not necessarily radical) ideal of $\kk[X]$. The \index{Closed!subscheme}\index{Subscheme!closed}{\em closed subscheme} of $X$ determined by $\qqq$, which by an abuse of notation\footnote{This is an abuse of notation since we also use $V(\cdot)$ to denote {\em subvarieties} of a given variety. We will try to ensure that the intended meaning is clear from the context.} we denote by $V(\qqq)$ is the pair $(Z', R)$, where $Z'$ is the {\em subvariety} of $X$ determined by $\qqq$, and $R = \kk[X]/\qqq$. We say that $Z'$ is the \index{Support!of a closed subscheme}{\em support} of $V(\qqq)$. One can picture $V(\qqq)$ as a ``thickened'' version of $Z'$. For example, if $X = \kk^2$ and $\qqq = \langle x, y^2 \rangle \subseteq \kk[x,y]$, then $V(\qqq)$ is supported at the origin. The image of $g = a + bx + cy + dx^2 + exy + fy^2 + \cdots$ in $\kk[x,y]/\qqq$ is determined by $a = g(0,0)$ and $c = \D{y}{g}(0,0)$ so that $V(\qqq)$ is the origin coupled with the vertical tangent line at the origin (\cref{fig:vertical-origin}). In general one can picture $V(\qqq)$ as a union of thickened varieties corresponding to {\em primary decompositions} (\cref{appendix:primary-decomposition}) of $\qqq$, and different primary decompositions may lead to different pictures of the same closed subscheme - see \cite[Section 3.8]{eisenview} for an illuminating exposition.   

\begin{center}
\begin{figure}[htb]
\begin{tikzpicture}[
    thick,
    dot/.style = {
      draw,
      fill,
      circle,
      inner sep = 0pt,
      minimum size = 6pt
    },
    xscale =0.66,
    yscale=0.33]

	\draw (0,0) node[dot, red] {};
	\draw [->, >= stealth, red, thick] (0,0) --(0,1.5);
\end{tikzpicture}
\caption{The subscheme of $\kk^2$ corresponding to $\qqq = \langle x, y^2 \rangle$ is the origin coupled with a vertical tangent} \label{fig:vertical-origin}
\end{figure}
\end{center}

\begin{example}
If $Z'$ is a subvariety of $X$, then in general there are infinitely many closed subschemes of $X$ supported at $Z'$. However, there is a canonical one among these, namely the subscheme $V(I(Z'))$, where $I(Z')$ is the ideal in $\kk[X]$ consisting of all regular functions that vanish on $Z'$; this is called the \index{Reduced!subscheme}{\em reduced} subscheme structure on $Z'$ (since the ``coordinate ring'' $\kk[X]/I(Z')$ is reduced). 
\end{example}

\begin{example}
Following an answer on {\em MathOverflow} \cite{mathoverflow-algeom-examples-nilpotent} we now present an example where different closed subschemes with the same support appear ``naturally.'' Recall that an $n \times n$ matrix $A$ over $\kk$ is {\em nilpotent} if $A^k = 0$ for some $k \geq 1$. The space $\scrM_n$ of $n \times n$ matrices can be naturally identified with $\kk^{n^2}$ with coordinates $x_{ij}$, $1 \leq i, j \leq n$. For $A \in \scrM_n$ it is a standard result from linear algebra that each of the following properties are equivalent to $A$ being nilpotent:
\begin{prooflist}
\item \label{nilpotent:n=0} $A^n = 0$, 
\item \label{nilpotent:char=0} $\det(A - \lambda \id_n) = \lambda^n$, where $\id_n$ is the $n \times n$ identity matrix, and $\lambda$ is an indeterminate over the $x_{ij}$. 
\end{prooflist}
Each of these conditions identifies the set $\scrN_n$ of nilpotent matrices as the set of zeroes of a system of polynomials in $(x_{ij})_{i,j}$. Let $\qqq_1, \qqq_2$ be the ideals generated by respectively these two systems of polynomials. Then both $V(\qqq_j)$ are supported at $\scrN_n$. We now show that $\qqq_1 \neq \qqq_2$ when $n > 1$. Indeed, since the trace of a nilpotent matrix is zero, it follows that $f := x_{11} + \cdots + x_{nn} \in I(\scrN_n) = \sqrt{\qqq_1} = \sqrt{\qqq_2}$. Since all polynomials that arise from condition \ref{nilpotent:n=0} are {\em homogeneous} of degree $n$, it follows that $f \not\in \qqq_1$ if $n > 1$. On the other hand, the coefficient of $\lambda^{n-1}$ in $\det(A - \lambda \id_n)$ is $\pm f$, so that $f \in \qqq_2$. 
\end{example}

\subsection{Closed subschemes of a quasiprojective variety} \label{subscheme-ssection-qproj}
Defining subschemes on an arbitrary quasiprojective variety is a bit more complicated than the case of affine varieties, since the same set of equations can look different in different coordinate charts. We want the subscheme to be a notion which would keep track of equations. Giving a set of equations on a neighborhood of a point $x$ is essentially same as giving an ideal of the local ring $\local{X}{x}$ of $X$ at $x$. However, the equations at different points need to be compatible, i.e.\ the equations at all points on a sufficiently small neighborhood must ``come from the same set of equations.'' This leads to the following definition: a \index{Sheaf of ideals}{\em sheaf $\scrI$ of ideals} on a quasiprojective variety $X$ is a product $\prod_{x \in X} \scrI_x$, where each $\scrI_x$ is an ideal of $\local{X}{x}$, such that
\begin{align}
\parbox{0.6\textwidth}{
each $x \in X$ has a nonempty open affine neighborhood $U$ in $X$ and an ideal $I$ of $\kk[U]$ such that $\scrI_{x'} = I\local{X}{x'}$ for each $x' \in U$. 
}\label{subscheme-compatibility}
\end{align}
For each $x \in X$, we say that $\scrI_x$ is the \index{Stalk of a sheaf of ideals}{\em stalk of $\scrI$ at $x$}. For us the \index{Closed!subscheme}\index{Subscheme!closed}{\em closed subscheme} $V(\scrI)$ of $X$ determined by $\scrI$ would be the product $\prod_{x \in X}\local{X}{x}/\scrI_x$ of quotient rings\footnote{If you are already familiar with schemes you will note that we are identifying a closed subscheme with its structure sheaf.}. An \index{Embedded!affine chart}{\em embedded affine chart} of $V(\scrI)$ is an affine open subset $U$ of $X$ which satisfies condition \eqref{subscheme-compatibility}. Let $\{U_j\}$ be an open covering of $X$ by embedded affine charts of $V(\scrI)$, i.e.\ for each $j$, condition \eqref{subscheme-compatibility} is satisfied with $U = U_j$ and $I = I_j$ for some ideal $I_j$ of $\kk[U_j]$. Then the union of the subvarieties $Z'_j$ of $X$ determined by $I_j$ is in fact a subvariety $Z'$ of $X$. We say that $Z'$ is the \index{Support!of a closed subscheme}{\em support} of $V(\scrI)$, and write $Z' = \supp(V(\scrI))$. 

\begin{example}
Let $X$ be an affine variety. Every ideal $\qqq$ of $\kk[X]$ canonically corresponds to the sheaf $\scrI_\qqq := \prod_{x \in X}(\qqq\local{x}{X})$ of ideals on $X$. We identify the closed subscheme $V(\qqq)$ defined in \cref{subscheme-ssection-affine} with the closed subscheme $V(\scrI_\qqq)$ of $X$. 
\end{example}

\begin{example} \label{example:zero-sheaf}
Every variety can be regarded as a closed subscheme of itself corresponding to the sheaf of ``zero ideals.''
\end{example}

\begin{example} \label{example:proj-sheaf-1}
Let $X$ be a subvariety of $\pp^n$ with homogeneous coordinates $[x_0: \cdots : x_n]$. If $f$ is any homogeneous polynomial in $(x_0, \ldots, x_n)$ of degree $d$, then $f/x_j^d$ is a regular function on the basic open subset $U_j := \pp^n \setminus V(x_j)$ for each $j = 0, \ldots, n$. Moreover, on $U_i \cap U_j$, $f/x_j^d = u_{ji}f/x_i^d$, where $u_{ji} = x_i^d/x_j^d$ is a unit in $\kk[U_i \cap U_j]$. It follows that $f$ defines a sheaf of ideals on $\pp^n$ whose stalk at $x \in X \cap U_j$ is the ideal of $\local{X}{x}$ generated by $f/x_j^d$; we write $V(f)$ for the corresponding closed subscheme of $X$. It is straightforward to check that $\supp(V(f))$ is precisely the subvariety of $X$ determined by $f$, and $X \cap U_j$ is an embedded affine chart of $V(f)$ for each basic open subset $U_j$ of $\pp^n$.
\end{example}

\begin{example} \label{example:proj-sheaf-N}
The arguments from \cref{example:proj-sheaf-1} can be generalized in a straightforward way to show that any finite collection $f_1, \ldots, f_N$ of homogeneous polynomials in $(x_0, \ldots, x_n)$ determines a sheaf of ideals $\scrI(f_1, \ldots, f_N)$ on a quasiprojective subset $X$ of $\pp^n$ such that for each $x \in X \cap U_j$, the stalk of $\scrI(f_1, \ldots, f_N)$ at $x$ is the ideal of $\local{X}{x}$ generated by $f_1/x_j^{\deg(f_1)}, \ldots, f_N/x_j^{\deg(f_N)}$; the support of the corresponding closed subscheme $V(f_1, \ldots, f_N)$ of $X$ is precisely the subvariety of $X$ determined by $f_1, \ldots, f_N$. If $I$ is a homogeneous ideal of $\kk[x_0, \ldots, x_n]$, and $f_1, \ldots, f_N$ are homogeneous generators of $I$, then it is straightforward to check that the sheaf of ideals $\scrI(f_1, \ldots, f_N)$ does {\em not} depend on $f_1, \ldots, f_N$ (i.e.\ if $g_1, \ldots, g_M$ are homogeneous generators of $I$, then $\scrI(f_1, \ldots, f_N) = \scrI(g_1, \ldots, g_M)$); we denote the corresponding closed subschme on $X$ by $V(I)$. 
\end{example}

Given a sheaf $\scrI$ of ideals on a variety $X$, let $Z := V(\scrI)$, and $Z' := \supp(Z)$. Given a point $x \in Z'$, we often abuse the notation and say $x \in Z$. The \index{Local ring!of a closed subscheme at a point}{\em local ring} $\local{Z}{x}$ of $Z$ at $x$ is the quotient $\local{X}{x}/\scrI_x$. We say that $Z$ has \index{Dimension!of a closed subscheme}\index{Pure dimension!of a closed subscheme}(pure) dimension $k$ if and only if $Z'$ has (pure) dimension $k$. If $U$ is an open subset of $X$, then $\scrI|_U := \prod_{x \in U} \scrI_x$ is a sheaf of ideals on $U$; we denote the corresponding closed subscheme of $U$ by the \index{Scheme theoretic intersection}``scheme-theoretic intersection'' $Z \cap U$, and say that it is an \index{Open subscheme}\index{Subscheme!open}{\em open subscheme} of $Z$. There is also a scheme-theoretic intersection of two closed subschemes: if $Y = V(\scrJ)$ is a closed subscheme of $X$ corresponding to a sheaf $\scrJ$ of ideals on $X$, then the scheme-theoretic intersection $Y \cap Z$ is the closed subscheme of $X$ corresponding to the sheaf of ideals $\scrI + \scrJ := \prod_{x \in X}(\scrI_x + \scrJ_x)$. We identify the variety $X$ with its closed subscheme defined by the sheaf of zero ideals. This in particular implies that the scheme-theoretic intersection $Z \cap X$ is simply $Z$, as expected.

\begin{example}
Given a quasiprojective subset $X$ of $\pp^n$ and homogeneous polynomials $f_1, \ldots, f_N$ in $(x_0, \ldots, x_n)$, if $V(f_j)$ are closed subschemes of $X$ constructed in \cref{example:proj-sheaf-1}, then the scheme-theoretic intersection $\bigcap V(f_j)$ is precisely $V(f_1, \ldots, f_N)$ constructed in \cref{example:proj-sheaf-N}. 
\end{example}

\subsection{Rational functions} \label{scheme:rational-functions}
\index{Rational!function!on a closed subscheme}
Usually the notion of rational functions is considered only for irreducible varieties. \Cref{exercise:k(X)-reducible} gives an hint that defining a rational function on a reducible variety can get tricky due to the presence of nonzero regular functions which are {\em zero-divisors}, i.e.\ which vanish identically on some irreducible component. Consider e.g.\ $X = V(xy) \subset \kk^2$, i.e.\ $X$ is the union of $x$ and $y$-axes on $\kk^2$. In this case both $x$ and $y$ are zero-divisors (in fact $xy = 0$ on $X$), and it would be difficult to give geometric interpretation of a ring containing $1/x$ and $1/y$ (what would be the ``value'' of $1/x + 1/y$ at a point on $X$?). A standard solution is therefore not to allow zero-divisors in the denominator: let $Z$ be a subscheme of $X$. For each $x \in Z$, let $S_x$ be the localization of $\local{Z}{x}$ at the set of its non zero-divisors, i.e.\ $S_x := \{f_1/f_2: f_1, f_2 \in \local{Z}{x}$, $f_2$ is {\em not} a zero-divisor in $\local{Z}{x}\}$. Let $Z' := \supp(Z)$. A {\em rational function} on $Z$ is an element $f = (f_x: x \in Z') \in \prod_{x \in Z'} S_x$ such that 
\begin{align}
\parbox{0.63\textwidth}{
each $x \in Z'$ has a nonempty open affine neighborhood $U$ in $Z'$ and $f_1, f_2 \in \kk[U]$ such that $f_{x'} = f_1/f_2 \in S_{x'}$ for each $x' \in U$. 
}\label{rational-compatibility}
\end{align}
A rational function $f = (f_x: x \in Z')$ on $Z$ is a \index{Regular function!on a closed subscheme}{\em regular function} if each $f_x \in \local{Z}{x}$, and it is an {\em invertible rational} function if $1/f$ is also a rational function, i.e.\ if \eqref{rational-compatibility} holds with the additional condition that no $f_i$ is a zero-divisor in $\local{Z}{x'}$ for any $x' \in U \cap Z'$. 

\begin{example}
Let $X$ be a variety. Regardless of whether we treat $X$ as a variety or the closed subscheme of itself determined by the zero ideal, the set of regular functions on $X$ remain the same, and in addition, if $X$ is irreducible, then the set of rational functions on $X$ remains the same (recall that in \cref{rational-section} we did {\em not} define the rational functions on a reducible variety). 
\end{example}

\begin{example}
Let $X$ be the {\em subvariety} $V(xy)$ of $\kk^2$, so that $\kk[X] = \kk[x,y]/ \langle xy \rangle$. Since both $x$ and $y$ are zero divisors in $\kk[X]$, it follows that the set of rational functions on $X$ (when we treat $X$ as its closed subscheme defined by the zero ideal) can be identified with $\{f/g: f, g \in \kk[X],\ g(\origin) \neq 0\} \cong \local{\kk^2}{\origin}/\langle xy \rangle$. 
\end{example}

\subsection{Completeness and compactification of schemes} \label{scheme:compactification}
Let $\phi : Y \to X$ be a morphism of varieties. If $\scrI = \prod_{x \in X}\scrI_x$ is a sheaf of ideals on $X$, then $\phi^*\scrI := \prod_{y \in Y} \phi^*(\scrI_{\phi(y)})\local{Y}{y}$ is a sheaf of ideals on $Y$. If $\phi$ is an {\em isomorphism of varieties}, then for each $y \in Y$, $\scrI_{\phi(y)}$ is naturally isomorphic as an $\local{X}{\phi(y)}$-module to $(\phi^*\scrI)_y = \phi^*(\scrI_{\phi(y)})\local{Y}{y}$; we say that $\phi^*: V(\scrI) \to V(\phi^*(\scrI))$ is an \index{Embedded!isomorphism of closed subschemes}{\em embedded isomorphism of closed subschemes}. A basic example of embedded isomorphism arises in the following context. 

\begin{example} \label{example:open-embedded-affine}
Given a polynomial $g \in \kk[x_1, \ldots, x_n]$, recall that $X := \kk^n \setminus V(g)$ is isomorphic to the subvariety $Y := V(x_{n+1}g - 1)$ of $\kk^{n+1}$ (\cref{example:affine-complement}). Now assume $Z$ is the closed subscheme of $\kk^n$ determined by an ideal $I$ of $\kk[x_1, \ldots, x_n]$ such that $g$ does not identically vanish on $\supp(Z)$, i.e.\ $g \not\in \sqrt{I}$. Then the open subscheme $V := Z \cap X$ of $Z$ is a closed subscheme of $X$ with nonempty support. The isomorphism $Y \to X$ induces an embedded isomorphism between $V$ and the closed subscheme of $Y$ defined by the ideal generated by $I$. 
\end{example}

 A \index{Compactification!of a closed subscheme}{\em compactification} of $Z := V(\scrI)$ is a closed subscheme $\bar Z$ of a compactification $\bar X$ of $X$ such that $Z$ is embedded isomorphic to an open subscheme of $\bar Z$. A closed subscheme of a variety is \index{Complete!subscheme}{\em complete} if its support is complete; note that $\supp(\bar Z)$, being a closed subvariety of a complete variety, is complete (\cref{prop:closed-sub-complete}). A fundamental result of Nagata states that every closed subscheme $Z$ of a given variety $X$ can be compactified to a closed subscheme of a given compactification $\bar X$ of $X$. In this chapter we will use Nagata's result for the case that $\dim(Z) = 1$, which we now prove. 
 
 \begin{rem}
Our proof of \cref{nagata1} below would have been much shorter if we had used the fact that if $C$ is an irreducible curve and $S$ is a finite nonempty subset of $C$, then $C \setminus S$ is an affine curve. But with the tools developed in \cref{var-chapter}  we could prove only an approximate version of it, namely $C \setminus S'$ is affine for some $S' \supseteq S$, where $S'$ is a finite set, and given any finite set $S'' \subseteq C \setminus S$, we can ensure that $S' \cap S'' = \emptyset$. 
 \end{rem}

\begin{thm} \label{nagata1}
Let $Z$ be a one dimensional closed subscheme of a quasiprojective variety $X$ and $\bar X$ be a projective compactification of $X$. Then there is a closed subscheme $\bar Z$ of $\bar X$ such that 
\begin{enumerate}
\item \label{nagata1:bar} $Z$ is embedded isomorphic to $\bar Z \cap X$.
\item \label{nagata1:bar-support} $\supp(\bar Z)$ is the closure in $\bar X$ of $\supp(Z)$.
\item \label{nagata1:bar-rational} every rational function on $Z$ extends to a rational function on $\bar Z$. 
\item \label{nagata1:bar-invrational} every invertible rational function on $Z$ extends to an invertible rational function on $\bar Z$. 
\end{enumerate}
\end{thm}

\begin{proof}
Since $Z' := \supp(Z)$ has dimension one, $S := \bar Z' \setminus Z'$, where $\bar Z'$ is the closure of $Z'$ in $\bar X$, is {\em finite}.

\begin{proclaim} \label{claim:nagata11}
There is an affine open subset $\bar W$ of $\bar X$ containing $S$ and there is a non zero-divisor $g \in \kk[\bar W]$ such that $V(g) \cap \bar Z' = S$, $W := \bar W \setminus V(g) \subset X$. In addition one can ensure that 
\begin{defnlist}
\item $W$ is an {\em embedded affine chart} of $Z$, and
\item \label{prop:W-unmixed} the ideal $I$ of $\kk[W]$ defining $Z \cap W$ is {\em unmixed}, i.e.\ has no zero-dimensional\footnote{The dimension of an ideal $J$ of $\kk[W]$ is the dimension of the subvariety $V(J)$ of $W$.} associated prime ideals.
\end{defnlist}
\end{proclaim}

\begin{center}
\begin{figure}[htb]
\includegraphics[width=50mm]{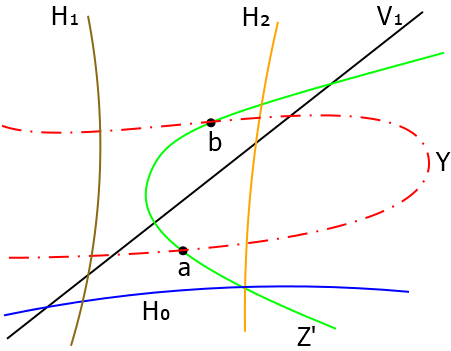}
\caption{$Y := \bar X \setminus X$,  $S = \{a, b\}$, $\bar W := \bar X \setminus (V_1 \cup H_1 \cup H_2)$, $W := \bar W \setminus (Y \cup H_0)$.}
\label{fig:nagata11}
\end{figure}
\end{center}

\begin{proof}
Choose a closed embedding $\bar X \into \pp^N$ with homogeneous coordinates $[x_0: \cdots : x_N]$. After composing $\phi$ with an appropriate Veronese embedding (see \cref{veronese-subsection}) followed by a linear change of coordinates if necessary we may ensure that  
\begin{prooflist}
\item $\bar X \setminus X \subseteq V(x_0)$, 
\item \label{prop:irr-notin-Vx0} $V(x_0)$ does {\em not} contain any irreducible component of $\bar X$ or any irreducible component of $\bar Z'$,
\item \label{prop:S-notin-Vx1} $S \cap V(x_1) = \emptyset$.
\end{prooflist}
(In \cref{fig:nagata11} $V(x_1)$ is denoted by $V_1$, and $V(x_0)$ is the union of $Y := \bar X \setminus X$ and possibly some other variety, say $H_0$, not containing any irreducible component of $\bar Z'$.) Let $U_1 := \pp^N \setminus V(x_1)$. Property \ref{prop:S-notin-Vx1} implies that $Z' \cap U_1 \neq \emptyset$, so that we can choose an affine open subset $W_1$ of $X \cap U_1$ such that 
\begin{prooflist}[resume]
\item \label{prop:W1-embedded} $W_1$ is also an embedded affine chart of $Z$, and
\item \label{prop:W1-unmixed} the ideal $I_1$ of $\kk[W_1]$ defining $Z \cap W_1$ is unmixed. 
\end{prooflist}
(In \cref{fig:nagata11} $W_1$ is the complement of $V_1 \cup H_1$.) Choose a regular function $h$ on $U_1$ such that 
\begin{prooflist}[resume]
\item \label{prop:h-vanishing} $V(h) \supseteq \left( (Z' \cap W_1 \cap V(x_0)) \cup ((X \cap U_1) \setminus W_1) \right)$
\item \label{prop:h-non-vanishing}  $h$ does {\em not} vanish at any point of $S \subset U_1$. 
\end{prooflist}
(In \cref{fig:nagata11} $V(h) \cap U_1$ is the union of $H_1$ and possibly some other variety, say $H_2$, not containing any point of $S$, but containing all points of $H_0 \cap Z'$.) Let $\bar W := (\bar X \cap U_1) \setminus V(h)$ and $W := \bar W \setminus V(x_0)$. Then $\bar W$ is an affine open subset of $\bar X$ (\cref{exercise:basic-open-equation}). Property \ref{prop:h-vanishing} implies that $W = W_1 \setminus V(hx_0/x_1)$ (note that $x_0/x_1$ is a regular function on $U_1 \supset W_1$). Since $W_1$ is an embedded affine chart of $Z$, it is straightforward to check that so is $W$. Properties \ref{prop:irr-notin-Vx0}, \ref{prop:W1-unmixed}, \ref{prop:h-vanishing} and \ref{prop:h-non-vanishing} then imply that the claim holds with $g = x_0/x_1$. 
\end{proof}

Choose $\bar W, W, g$ as in \cref{claim:nagata11}. Let $I$ be the ideal of $\kk[W]$ defining $Z \cap W$. Since $g$ is a non zero-divisor in $\kk[\bar W]$,
\begin{align}
\parbox{0.6\textwidth}{
the natural map $\kk[\bar W] \to \kk[\bar W]_g = \kk[W]$ is injective.
} \label{prop:Winjective}
\end{align}
Define $\bar I := I \cap \kk[\bar W]$. Then $I$ is generated by $\bar I$ in $\kk[W]$ (\cref{prop:localization-basic}) and $g$ is a non zero-divisor in $\kk[\bar W]/\bar I$ (\cref{prop:zero-union-minimally-prime}). If $\scrI$ is the sheaf of ideals on $X$ defining $Z$, then it follows that 
\begin{align*}
\prod_{x \in \bar W} \bar I\local{\bar X}{x} \times \prod_{x \in X \setminus W} \scrI_x  \times \prod_{x \in Y \setminus \bar W} \local{\bar X}{x}
&= \prod_{x \in X}\scrI_x \times \prod_{x \in S}\bar I\local{\bar X}{x} \times \prod_{x \in Y \setminus X} \local{\bar X}{x}
\end{align*}
(where $Y := \bar X \setminus X)$ is a sheaf of ideals on $\bar X$ and the corresponding subscheme $\bar Z$ of $\bar X$ satisfies assertions \eqref{nagata1:bar} and \eqref{nagata1:bar-support}. Now let $f = (f_x: x \in Z')$ be a rational function on $Z$. 

\begin{proclaim} \label{claim:nagata12}
There is an open subset $U'$ of $Z'$ which intersects every irreducible component of $Z'$, and in addition \eqref{rational-compatibility} is satisfied with $U = U'$.
\end{proclaim}

\begin{proof}
For each irreducible component $Z'_j$ of $Z'$, choose an open subset $U'_j$ of $Z'$ such that \eqref{rational-compatibility} is satisfied with $U = U'_j$, and $U'_j$ does {\em not} intersect any irreducible component of $Z'$ other than $Z'_j$. Then $U' := \bigcup_j U'_j$ satisfies the claim.
\end{proof}

Let $U'$ be as in \cref{claim:nagata12}. Due to assertion \ref{prop:W-unmixed} of \cref{claim:nagata11} one can find a regular function $q$ on $W$ such that $Z' \cap W \setminus V(q) \subseteq U$, and $q$ is a non zero-divisor in both $\kk[W]$ and $\kk[W]/I$. Note that $Z$ is defined in $W^* := W \setminus V(q)$ by the ideal generated by $I$. Moreover, there is $f_1, f_2 \in \kk[W^*]$ such that $f_2$ is a non zero-divisor in $\kk[W^*]/I\kk[W^*]$, and $f_x = f_1/f_2 \in \local{X}{x}/I\local{X}{x}$ for each $x \in W^* \cap Z$. Since $\kk[W^*] = \kk[W]_q = \kk[\bar W]_{gq}$, it follows that $f_1/f_2$ can be represented in the total quotient ring of $\kk[\bar W]/\bar I$ as $g^q \bar q^b \bar f_1/\bar f_2$ for some integers $a,b$ and $\bar q, \bar f_1, \bar f_2 \in \kk[\bar W]$ such that $\bar q$ and $\bar f_2$ are non zero-divisors in $\kk[\bar W]/\bar I$. This proves assertion \eqref{nagata1:bar-rational}. For assertion \eqref{nagata1:bar-invrational} note that if $f_1$ is a non zero-divisor in $\kk[W^*]/I\kk[W^*]$, then $\bar f_1$ is also a non zero-divisor in $\kk[\bar W]/\bar I$. 
\end{proof}

\subsection{Irreducible components, local rings} \label{scheme:local-section}
\index{Irreducible!component!of a closed subscheme}
Let $Z$ be a closed subscheme of a variety $X$. The {\em irreducible components} of $Z$ are simply the irreducible components of $\supp(Z)$. Note that every irreducible component of a closed subscheme $Z$ of $X$ is a {\em subvariety} of $X$. Let $Y$ be a closed irreducible {\em subvariety} of $\supp(Z)$. The \index{Local ring!of a closed subscheme at a subvariety}{\em local ring $\local{Z}{Y}$ of $Z$ at $Y$} is the set of the equivalence classes of pairs $(h, U)$, where $U$ is an open subset of $X$ such that $U \cap Y \neq \emptyset$ and $h$ is a {\em regular function} on the open subscheme $Z \cap U$ of $Z$, and the equivalence relation $\sim$ is defined as follows: $(h_1, U_1) \sim (h_2, U_2)$ if and only if $h_1 = h_2$ in $\local{Z}{x}$ for each $x \in U_1 \cap U_2$. It is straightforward to check that $\local{Z}{Y}$ is a $\kk$-algebra. A more explicit realization of $\local{Z}{Y}$ is as follows: pick an {\em embedded affine chart} $U$ of $Z$ such that $U \cap Y \neq \emptyset$. Then $U \cap Z$ is the closed subscheme of $U$ defined by an ideal $\qqq$ of $\kk[U]$. On the other hand, since $U \cap Y$ is an irreducible subvariety of $\supp(Z) \cap U$, it corresponds to a prime ideal $\ppp$ of $\kk[U]$ containing $\qqq$. Then $\local{Z}{Y}$ can be identified with the localization $(\kk[U]/\qqq)_\ppp$ of $\kk[U]/\qqq$ at the ideal generated by $\ppp$. This in particular implies that $\local{Z}{Y}$ is a {\em local ring}. If $Y = a$ is a point of $Z$, then it is straightforward to check that this definition of $\local{Z}{Y}$ agrees with the earlier definition of $\local{Z}{a}$ from \cref{subscheme-ssection-qproj}.

\begin{lemma} \label{reduced-invertible}
Let $Z$ be a closed subscheme of an affine variety $X$, and $Y$ be an irreducible subvariety of $\supp(Z)$. Given $f \in \kk[X]$, if the image of $f$ is invertible in $\local{\supp(Z)}{Y}$, then it is also invertible in $\local{Z}{Y}$. 
\end{lemma}

\begin{proof}
This immediately follows from the following observation: if $\phi: S \to T$ is a ring homomorphism such that $\ker(\phi)$ is contained in the nilradical of $S$, and if $u \in S$ is such that $\phi(u)$ is invertible in $T$, then $u$ is invertible in $S$.  
\end{proof}

\subsection{Cartier divisors} \label{scheme:cartier-section}
\index{Cartier divisor}
A {\em Cartier divisor} is a closed subscheme generated locally by single {\em non zero-divisors}; it is the natural scheme-theoretic analogue of a ``hypersurface.'' More precisely, a closed subscheme $Z = V(\scrI)$ of a variety $X$ is called a Cartier divisor if each $x \in X$ has a nonempty open affine neighborhood $U$ in $X$ and an element $g \in \kk[U]$ which is not a zero-divisor in $\kk[U]$ such that $\scrI_{x'} = g\local{X}{x'}$ for each $x' \in U$. It is straightforward to check that defining a Cartier divisor on $X$ is equivalent to prescribing a collection $\{(U_i, g_i)\}_i$ of pairs such that 
\begin{defnlist}
\item $\{U_i\}$ is an open affine covering of $X$, and 
\item $g_i \in \kk[U_i]$ are such that 
\begin{defnlist}
\item $g_i$ is a non zero-divisor in $\kk[U_i]$ for each $i$, and
\item $g_i/g_j$ is invertible in $\kk[U_i \cap U_j]$ for each $i,j$.
\end{defnlist}
\end{defnlist}
If $X$ is of pure dimension $n$, it follows from \cref{thm:pure-dimension} that the support of a Cartier divisor on $X$ is either empty or has pure dimension $n-1$. 

\begin{example}
Let $X = V(I)$ be the {\em irreducible subvariety} of $\pp^n$ determined by a {\em prime} homogeneous ideal $I$ of $\kk[x_0, \ldots, x_n]$. If $f$ is a homogeneous polynomial which is {\em not} in $I$, then the closed subscheme $V(f)$ of $X$ constructed in \cref{example:proj-sheaf-1} is a Cartier divisor on $X$.
\end{example} 

\section{Possibly non-reduced curves} \label{reduced-noncurve-section}
\index{Curve!reduced}\index{Curve!possibly non-reduced}\index{Reduced!curve}\index{Possibly non-reduced curve}
This section is devoted to task \eqref{task:order} outlined in \cref{sec:intro-mult-intro}. A {\em reduced curve} is a variety of pure dimension one. By a {\em possibly non-reduced curve} we mean a pure dimension one closed subscheme $Z$ of a variety. Unless explicitly stated otherwise, a {\em curve} will mean a reduced curve. Our convention of identifying a variety with its subscheme defined by the zero ideal sheaf implies that a ``curve'' is indeed a special case of a ``possibly non-reduced curve.'' In \cref{rcurve-section} we describe some properties of curves we are going to use without proof. In \cref{nrcurve-section} we define the notion of order at a point of a possibly non-reduced curve and describe some of its properties whose proofs are deferred to \cref{order-app}. 

\subsection{(Reduced) Curves} \label{rcurve-section}
Curves are in a sense the simplest nontrivial algebraic varieties. \Cref{curve-desingularization} below states one of their basic properties, namely that they can be {\em desingularized}. In particular, the map $\pi: \tilde C \to C$ from \cref{curve-desingularization} is called the \index{Desingularization of a curve}\index{Resolution of singularities of a curve}{\em desingularization} of $C$. It is the one-dimensional case of {\em resolution of singularities}, which is still an open problem for dimension greater than $3$ in nonzero characteristics. Proofs of \cref{curve-desingularization} can be found in many introductory algebraic geometry texts; in particular \cite{fulturve} gives an elementary (but long) proof, and \cite[Chapter 1]{kollar} contains an illuminating exposition of many different proofs. 

\begin{thm}[{\cite[Theorems II.5.6 and II.5.7]{shaf1}}] \label{curve-desingularization}
Let $C$ be an irreducible curve. Then there is a nonsingular irreducible curve $\tilde C$ and a surjective morphism $\pi: \tilde C \to C$ such that
\begin{enumerate}
\item \label{curve-desingularization:0} for each nonsingular point $a \in C$, $\pi$ restricts to an isomorphism near $\pi^{-1}(a)$;
\item \label{curve-desingularization:1} if $\phi: D \to C$ is any surjective morphism of curves with $D$ nonsingular, then there is a morphism $\tilde \phi: D \to \tilde C$ such that the following diagram commutes. 
\begin{center}
\begin{tikzcd}[column sep=tiny]
& D
\arrow[ld, swap, "\tilde \phi" ]
\arrow[rd, "\phi"]
& \\
\tilde C
\arrow[rr, "\pi"]
&
&C
\end{tikzcd}
\end{center}
\end{enumerate}
The curve $\tilde C$ is unique up to isomorphism. If $C$ is projective, then so is $\tilde C$. Moreover, condition \eqref{curve-desingularization:1} is automatically satisfied if $\tilde C$ is projective and $\pi$ satisfies condition \eqref{curve-desingularization:0}. 
\end{thm}


Let $a$ be a nonsingular point on a curve $C$. Identify an affine neighborhood of $a$ in $C$ with a curve in some $\kk^n$ with coordinates $(x_1, \ldots, x_n)$ defined by an ideal $I$ of $\kk[x_1, \ldots, x_n]$. Identity \eqref{tangent-dim} implies that there is $f \in \kk[x_1, \ldots, x_n]$ such that $f(a) = 0$ and $(\nabla f)(a) := ((\partial f/\partial x_1)(a), \ldots, (\partial f/\partial x_n)(a))$ is {\em not} in the span of $(\nabla g_1)(a), \ldots, (\nabla g_k)(a)$, where $g_1, \ldots, g_k$ are any set of generators of $I$. We say that $f$ is a \index{Parameter!of a curve at a nonsingular point}{\em parameter} of $C$ at $a$. The local ring $\local{C}{a}$ of $C$ at $a$ is a {\em discrete valuation ring} with parameter $f$ (\cref{curve-dvr-section}); we denote the discrete valuation of $\local{C}{a}$ by $\ord_a(\cdot)$, and given a rational function $g$ on $C$, we say that $\ord_a(g)$ is the \index{Order!at a nonsingular point on a curve}{\em order} of $g$ {\em at $a$}. 

\begin{example} \label{example:order-nonsing1}
Let $a$ be the origin and $C$ be the parabola $V(y-x^2) \subset \kk^2$. Since $\nabla(y-x^2)|_a = (0,1)$, it follows that $x$ is a parameter of $\local{C}{a}$. Since $y/x^2$ is invertible on $C$, it follows that $\ord_a(y|_C) = 2$.
\end{example}

\begin{thm}[{\cite[Corollary to Theorem III.2.1]{shaf1}}] \label{reduced-zero-sum}
Let $g$ be a nonzero rational function on a nonsingular curve $C$. Then there are only finitely many points $a$ on $C$ such that $\ord_a(g) \neq 0$. If $C$ is in addition projective, then 
\begin{align}
\sum_{a \in C}\ord_a(g) = 0 \label{nonsingular-order-sum}
\end{align} 
\end{thm}

Note that identity \eqref{nonsingular-order-sum} may not hold if $C$ is not projective, consider e.g.\ any non-constant polynomial on the affine line. The following result lists two basic properties of parameters. 

\begin{prop} \label{prop:nonsingular-parameter-properties}
Let $C$ be a curve on $\kk^n$ with coordinates $(x_1, \ldots, x_n)$. Then
\begin{enumerate}
\item \label{prop:nonsingular-parameter-coordinates} Given a nonsingular point $a = (a_1, \ldots, a_n)$ of $C$, there is $j$, $1 \leq j \leq n$, such that $x_j - a_j$ is a parameter of $C$ at $a$. 
\item \label{prop:nonsingular-parameter-open} Given a polynomial $f$ in $(x_1, \ldots, x_n)$, the property of $f$ being a parameter at a nonsingular point of $C$ is ``Zariski open,'' i.e.\ if $f$ is a parameter of $C$ at some nonsingular point of $C$, then it is a parameter of $C$ at each nonsingular point on a nonempty Zariski open subset of $C$. 
\end{enumerate}
\end{prop}

\begin{proof}
The first assertion follows immediately from the definition of parameters and identity \eqref{tangent-dim}. For the second assertion note that given a nonsingular point $a$ of $C$, due to \cref{local-equations} one can assume, after replacing $C$ by an appropriate open neighborhood of $a$ on $C$, that the ideal of $C$ in $\kk[x_1, \ldots, x_{n-1}]$ defined by $n-1$ polynomials $g_1, \ldots, g_{n-1}$. If $f$ is a parameter of $C$ at $a$, identity \eqref{tangent-dim} then implies that the determinant of the $n \times n$-matrix of partial derivatives of $g_1, \ldots, g_{n-1}, f$ is nonzero at $a$; therefore it is nonzero on a nonempty Zariski open subset $U$ of $C$. Then $f$ is a parameter of $C$ at each point of $U \cap C$, as required.  
\end{proof}

The following is a standard result covered in most introductory books in algebraic geometry. We use it in the proof of the main result (\cref{order-curve}) of next section. 
 
 
 \begin{thm}[{\cite[Theorem II.5.8]{shaf1}}] \label{curve-finite}
 Every non-constant morphism $\phi: C \to D$ between two {\em irreducible} projective curves is {\em finite}, i.e.\ if $U$ is any affine open subset of $D$, then $\phi^{-1}(U)$ is affine and $\kk[\phi^{-1}(U)]$ is a finite module over $\kk[U]$.
 \end{thm}
 
Consider the embedding of $\kk\setminus \{0\} \into \kk$. The coordinate ring of $\kk\setminus \{0\}$ is $\kk[x, x^{-1}]$, which is {\em not} a finite module over $\kk[x]$. This shows that the condition that $C$ and $D$ are projective is crucial in \cref{curve-finite}. Also, if $C$ is the union of the closures of $x$ and $y$ axes in $\pp^2$, then the projection map from $C$ to the $x$-axis is not finite; i.e.\ \cref{curve-finite} may fail to hold if $C$ is reducible.

\subsection{Order at a point on a possibly non-reduced curve} \label{nrcurve-section}
\index{Order!at a point on a possibly non-reduced curve}
Let $a$ be a point on a possibly non-reduced curve $C$ and $f \in \local{C}{a}$. The {\em order} $\ord_a(f)$ of $f$ at $a$ is the dimension of $\local{C}{a}/f\local{C}{a}$ as a vector space over $\kk$. Note that $\ord_a(f) = \infty$ if $f$ vanishes on any irreducible component of $\supp(C)$ containing $a$. 

\begin{example} \label{example:order-nonsing2}
Let $C$ be the parabola $V(y-x^2) \subset \kk^2$ and $f = y$. Since $y \equiv x^2$ on $C$, it follows that $\local{C}{\origin}/y\local{C}{\origin}$ is a $2$-dimensional vector space over $\kk$ generated by $1$ and $x$, so that $\ord_\origin(y|_C) = 2$.
Note that this agrees with the computation from \cref{example:order-nonsing1}. More generally, part \eqref{exercise:discrete-properties-dim} of \cref{prop:discrete-properties} shows that when $C$ is a nonsingular (reduced) curve, the two definitions of order agree. 
\end{example}

\begin{prop} \label{order-properties}
Let $a$ be a point on a possibly non-reduced curve $C$ and $f \in \local{C}{a}$. 
\begin{enumerate}
\item If $f$ is a non zero-divisor in $\local{C}{a}$, then $\ord_a(f) < \infty$.
\item $\ord_a(f) = 0$ if and only if $f$ is invertible in $\local{C}{a}$.
\item If $f$ is a non zero-divisor in $\local{C}{a}$ and $g \in \local{C}{a}$, then $\ord_a(fg) = \ord_a(f) + \ord_a(g)$.
\end{enumerate}
\end{prop}

The proof of \cref{order-properties} is given in \cref{order-app}. If $h$ is a {\em rational function} on $C$, then $h = f/g$ for some $f, g \in \local{C}{a}$ such that $g$ is a non zero-divisor in $\local{C}{a}$. We define $\ord_a(h)$ to be $\ord_a(f) - \ord_a(g)$. \Cref{order-properties} shows that $\ord_a(h)$ does not depend on the choice of $f$ or $g$. As \cref{example:order-nonsing2} suggests, it is straightforward to check using basic properties of discrete valuation rings that this definition of order agrees with the definition from \cref{rcurve-section} when both are applicable, i.e.\ $C$ is a nonsingular (reduced) curve. 

\begin{example} \label{example:order-singular}
Assume $\character(\kk) \neq 2$. Let $C' = V(x^2 - y^2 + y^3) \subset \kk^2$. We saw in \cref{example:curve-singularities} that $C'$ is singular at the origin. It is straightforward to check that $\local{C'}{\origin}/y\local{C'}{\origin}$ is a 2-dimensional vector space over $\kk$ generated by $1$ and $x$, so that $\ord_\origin(y|_{C'}) = 2$. Now we compute the order of $\pi^*(y)$ for a desingularization $\pi$ of $C'$. Consider the map $\pi: \kk \to C'$ defined as follows - given $t \in \kk\setminus \{0\}$, the straight line $\{(tu, u): u \in \kk\}$ through the origin with slope $1/t$ intersects $C'$ at a  single point other than the origin - define $\pi(t)$ to be that point (see \cref{fig:desingularize-node}). It is straightforward to compute that $\pi(t) = (t - t^3, 1 - t^2)$. It can be checked\footnote{Checking condition \eqref{curve-desingularization:0} of \cref{curve-desingularization} is straightforward. To verify condition \eqref{curve-desingularization:1}, note that $\pi$ extends to a map from $\pp^1$ to the closure $\bar C'$ of $C'$ in $\pp^2$, and use the last assertion of \cref{curve-desingularization}.} that $\pi$ is a desingularization of $C'$. Note that $\pi^{-1}(\origin)$ consists of two points $t = \pm 1$. It is straightforward to check that $\pi^*(y|_{C'}) = 1 - t^2$ is a {\em parameter} at each $\tilde a \in \pi^{-1}(\origin)$, so that $\ord_{\tilde a}(\pi^*(y|_{C'})) = 1$. It follows that 
\begin{align*}
\ord_\origin(y|_{C'}) = \sum_{\tilde a \in \pi^{-1}(\origin)} \ord_{\tilde a}(\pi^*(y|_{C'}))
\end{align*}
\end{example}

\begin{center}
\begin{figure}
\begin{tikzpicture}[scale=0.6]
\def\alabelx{2.5}
\def\alabely{2.6}
\def\nsamples{103}
\def\tmin{-3}
\def\tmax{3}
\def\rmin{-1}
\def\rmax{2.4}
\def\tzero{0.36}
\def\rtopx{3.4}
\def\rtopy{5.4}
\def\tzerox{3.15}
\def\tzeroy{4.4}

\def\xmin{-1.2}
\def\xmax{1}
\def\ymin{-1.3}
\def\ymax{1.5}

\begin{axis}[
xmin = \xmin, xmax=\xmax, ymin = \ymin, ymax= \ymax,
axis equal=true, axis equal image=true, hide axis
]
\addplot[blue, thick, domain=\tmin:\tmax, samples=\nsamples] ({x-x^3} ,{1 - x^2});
\addplot[thick, domain=\rmin:\rmax, samples=2] ({\tzero*x} ,{x});
\addplot[mark=*, red, mark size = 3pt] plot coordinates {({\tzero - \tzero^3} ,{1 - \tzero^2})};

\end{axis}
\draw (\alabelx, \alabely) node [left] {\picfontsize $O$};
\draw (\rtopx, \rtopy) node [right] {\picfontsize $x = ty$};
\draw (\tzerox, \tzeroy) node [right] {\picfontsize $\pi(t)$};
\end{tikzpicture}

\caption{Desingularization of the nodal cubic $x^2 = y^2 - y^3$}
\label{fig:desingularize-node}
\end{figure}
\end{center}

\begin{example} \label{example:order-non-reduced}
Assume $\character(\kk) \neq 2$. Let $C$ be the closed subscheme of $\kk^2$ corresponding to the ideal of $\kk[x,y]$ generated by $(x^2 - y^2 + y^3)^2$. Note that $C$ is {\em non-reduced} and $\supp(C)$ is the singular curve $C'$ from \cref{example:order-singular}. It is not hard to see that $\local{C}{\origin}/y\local{C}{\origin}$ is a 4-dimensional vector space over $\kk$ generated by $1, x, x^2, x^3$, so that $\ord_\origin(y|_C) = 4$. Combining this with the observation from \cref{example:order-singular} yields that
\begin{align*}
\ord_\origin(y|_C) = 2\sum_{\tilde a \in \pi^{-1}(\origin)} \ord_{\tilde a}(\pi^*(y|_{C'}))
\end{align*}
here the factor $2$ on the right hand side is the {\em multiplicity} of $C'$ in $C$, which we now define. 
\end{example}

Let $C$ be a closed subscheme of a variety and $D$ be an irreducible component of $\supp(C)$. Pick an embedded affine chart $U$ of $C$ such that $U \cap D \neq \emptyset$. Recall that if $\qqq$ (respectively, $\ppp$) is the ideal of $\kk[U]$ defining $C$ (respectively, $D$), then the local ring $\local{C}{D}$ of $C$ at $D$ can be identified with the localization $(\kk[U]/\qqq)_\ppp$ of $\kk[U]/\qqq$ at the ideal generated by $\ppp$. Since $D$ is an irreducible component of $C$, it follows that $\ppp$ is a {\em minimal} prime ideal containing $\qqq$. Therefore the {\em length} of $\local{C}{D}$ as a module over itself is finite (\cref{prop:minimal-localization-finite-length}); we call it the \index{Multiplicity of an irreducible component of a closed subscheme}{\em multiplicity} of $D$ in $C$ and denote it by $\mu_D(C)$. In \cref{example:order-non-reduced}, $U = \kk^2$, $\qqq$ (respectively, $\ppp$) is the ideal of $\kk[x,y]$ generated by $(x^2 - y^2 + y^3)^2$ (respectively, $x^2 - y^2 + y^3$). It is straightforward to check directly that $\local{C}{D} \supsetneq \ppp \local{C}{D} \supsetneq 0$ is a maximal chain of ideals of $\local{C}{D}$, which implies that $\mu_D(C) = 2$. The following result, whose proof is given in \cref{order-app}, is the key relation between orders at a point on a possibly non-reduced curve and the desingularizations of its irreducible components; it shows that our observation from \cref{example:order-non-reduced} holds in general. 

\begin{thm} \label{order-curve}
Let $a$ be a point on a possibly non-reduced curve $C$. Let $C_1, \ldots, C_s$ be the irreducible components of $\supp(C)$ containing $a$ and $\pi_i:\tilde C_i \to C_i$ be the desingularizations of $C_i$. If $f$ is a non zero-divisor in $\local{C}{a}$, then
\begin{align}
\ord_a(f) 
	&= \sum_i \mu_{C_i}(C) \ord_{ a}(f|_{C_i})
	= \sum_i \mu_{C_i}(C) \sum_{\tilde a \in \pi_i^{-1}(a)} \ord_{\tilde a}(\pi_i^*(f|_{C_i})) 
	\label{order-curve-identity}
\end{align}
\end{thm}

\begin{cor} \label{zero-sum}
Let $C$ be a possibly non-reduced curve and $h$ be an {\em invertible rational function} (see \cref{scheme:rational-functions}) on $C$. If $\supp(C)$ is a {\em projective} curve, then $\sum_{a \in C}\ord_a(h) = 0$. 
\end{cor}

\begin{proof}
Let $C_1, \ldots, C_s$ be the irreducible components of $C$ containing $a$ and $\pi_i:\tilde C_i \to C_i$ be the desingularizations of $C_i$. \Cref{order-curve} implies that 
\begin{align*}
\sum_{a \in C}\ord_a(h)
	&= \sum_i \mu_{C_i}(C) \sum_{\tilde a \in C_i} \ord_{\tilde a}(\pi_i^*(h|_{C_i})) 
\end{align*}
\Cref{prop:non-zero-divisor-restriction} implies that each $i$, $h|_{C_i}$ is a well-defined rational function on $C_i$. It then follows from \cref{reduced-zero-sum} that $\sum_{\tilde a \in C_i} \ord_{\tilde a}(\pi_i^*(h|_{C_i})) = 0$ for each $i$, as required.
\end{proof}

\begin{cor} \label{linear-intersection}
Let $X$ be a projective variety and $C$ be a projective (reduced) curve on $X \times \pp^1$. Assume no component of $C$ is contained in $X \times \{a\}$ for any $a \in \pp^1$. Then each component of $C$ intersects $X \times \{a\}$ for every $a \in \pp^1$. 
\end{cor}

\begin{proof}
We may assume \woutlog\ that $C$ is irreducible. Fix a point $(x_0,a_0)$, where $x_0 \in X$ and $a_0 \in \pp^1$, on $C$. Choose an arbitrary point $a \neq a_0 \in \pp^1$. We will show that $C$ intersects $X \times \{a\}$. Pick a point $\infty \in \pp^1 \setminus \{a,a_0\}$. Identify $\pp^1 \setminus \{\infty\}$ with $\kk$, so that we can treat $a,a_0$ as elements of $\kk$. Let $t$ be a coordinate on $\kk$. \Cref{zero-sum} implies that 
\begin{align*}
\sum_{(x,b) \in C} \ord_{(x,b)}((t-a)|_C)
	&= \sum_{(x,b) \in C} \ord_{(x,b)}((t-a_0)|_C) 
	= 0
\end{align*}
Since $f := (t-a)/(t-a_0)$ is regular and nonzero at all points of $X \times \{\infty\}$ [check that $f(x, \infty) = 1$ for each $x \in X$], \cref{order-properties} implies that 
\begin{align*}
\sum_{(x,b) \in C \cap (X \times \{\infty\})} \ord_{(x,b)}((t-a)|_C)
	&= \sum_{(x,b) \in C \cap (X \times \{\infty\})} \ord_{(x,b)}((t-a_0)|_C) 
\end{align*}
Note that $C \setminus (X \times \{\infty\}) = C \cap (X \times \kk)$, and for each $(x,b) \in  C \cap (X \times \kk)$, 
\begin{align*}
\ord_{(x,b)}((t-a)|_C)
	&= 
		\begin{cases}
		\text{positive} & \text{if}\ b = a,\\
		0 & \text{otherwise.}
		\end{cases}
\end{align*}
It follows that 
\begin{align*}
\sum_{(x,a) \in C \cap (X \times \{a\})} \ord_{(x,a)}((t-a)|_C)
	&= \sum_{(x,a_0) \in C \cap (X \times \{a_0\})} \ord_{(x,a_0)}((t-a_0)|_C) 
	\geq \ord_{(x_0,a_0)}((t-a_0)|_C)
	> 0
\end{align*}
This implies that $C \cap (X \times \{a\}) \neq \emptyset$, as required. 
\end{proof}

\section{Intersection multiplicity at a nonsingular point of a variety} \label{mult-point-section}
\index{Intersection multiplicity!at a nonsingular point}
\subsection{Intersection multiplicity of power series}
The {\em intersection multiplicity at the origin} of $f_1, \ldots, f_n \in \kk[[x_1, \ldots, x_n]]$ is 
\begin{align}
\multzero{f_1}{f_n} := \dim_\kk(\kk[[x_1, \ldots, x_n]]/\langle f_1, \ldots, f_n \rangle) \label{mult-0-defn}
\end{align}
Every power series can be approximated up to arbitrarily high order by polynomials. The following result shows that the intersection multiplicity of power series can also be approximated up to arbitrarily high order by (sufficiently close) polynomial approximations. Recall that the \index{Order!of a formal power series}{\em order} $\ord(f)$ of a power series $f$ is the smallest $m$ for which there is there is a monomial $x_1^{\alpha_1} \cdots x_n^{\alpha_n}$ with nonzero coefficient in $f$ such that $\sum_j \alpha_j = m$.

\begin{prop} \label{finite-mult-determinacy}
Let $f_1, \ldots, f_n \in \kk[[x_1, \ldots, x_n]]$.
\begin{enumerate}
\item If $\multfzero < \infty$, then there exists $m \geq 0$ such that $\multgzero = \multfzero$ for all $g_1, \ldots, g_n \in \kk[x_1, \ldots, x_n]$ such that $\ord(f_j - g_j) \geq m$. 
\item If $\multfzero = \infty$, then for each $N \geq 0$, there exists $m \geq 0$ such that $\multgzero \geq N$ for all $g_1, \ldots, g_n \in \kk[x_1, \ldots, x_n]$ such that $\ord(f_j - g_j) \geq m$.  
\end{enumerate}
\end{prop}

\begin{proof}
It follows immediately from \cref{finite-determinacy}, by taking $\preceq$ e.g.\ to be the graded lexicographic order (see \cref{grlex}). 
\end{proof}

\subsection{Intersection multiplicity of regular functions}
Let $a$ be a nonsingular point of a variety $X$ of dimension $n$. Recall that the completion $\hatlocal{X}{a}$ of the local ring $\local{X}{a}$ of $X$ at $a$ is isomorphic to $\kk[[x_1, \ldots, x_n]]$ (\cref{smooth-completion}). If $f_1, \ldots, f_n$ are regular functions on a neighborhood of $a$ in $X$, the {\em intersection multiplicity at $a$} of $f_1, \ldots, f_n$ is
\begin{align}
\multf{a} 
	&:= \dim_\kk(\hatlocal{X}{a}/\langle f_1, \ldots, f_n \rangle) \label{mult-a-defn-1}
\end{align}
where we identify each $f_j$ with its natural image in $\hatlocal{X}{a}$ (see \cref{prop:complete-inclusion}). \Cref{exercise:finite-quotient} implies that
\begin{align}
\multf{a} 
	&= \dim_\kk(\local{X}{a}/\langle f_1, \ldots, f_n \rangle) \label{mult-a-defn-2}
\end{align}
Now we deduce some basic properties of intersection multiplicity using results of the preceding sections. In addition, the ``unmixedness theorem'' of F.\ S.\ Macaulay (\cref{thm:Macaulay}) is used in a fundamental way in all the upcoming results of this and the following section. Given a subset $S$ of $X$ and $b \in S$, we say that $b$ is an \index{Isolated point}{\em isolated} point of $S$ if there is an open neighborhood $U$ of $b$ in $X$ such that $b$ is the only point of $S \cap U$. 

\begin{prop} \label{int-mult-curve}
Let  $V := V(f_1, \ldots, f_n) \subseteq X$. 
\begin{enumerate}
\item \label{outside-zero} $\multf{a} =0 $ if and only if $a \not\in V$. 
\item \label{isolated-finite}$0 <\multf{a} < \infty$ if and only if $a$ is an isolated point of $V$. 
\item \label{finite-to-curve} If $0 < \multf{a} < \infty$, then there is a Zariski open neighborhood $U$ of $a$ in $X$ such that $V(f_2, \ldots, f_n) \cap U$ has pure dimension one.
\item \label{mult-to-order} If there is a Zariski open neighborhood $U$ of $a$ in $X$ such that $C := V(f_2, \ldots, f_n) \cap U$ is a pure dimension one closed subscheme of $U$, then $\multf{a} = \ord_a(f_1|_C)$. 
\item \label{mult-additive} $\multp{f_1f'_1, f_2}{f_n}{a} = \multf{a} + \multp{f'_1,f_2}{f_n}{a}$. 
\end{enumerate}
\end{prop}

\begin{proof}
Assertion \eqref{outside-zero} is clear. Assertion \eqref{isolated-finite} follows from identity \eqref{mult-a-defn-2} and the fact that $a$ is an isolated point of $V$ if and only if the radical of the ideal generated by $f_1, \ldots, f_n$ in $\local{X}{a}$ is the maximal ideal of $\local{X}{a}$. The third assertion follows from assertion \eqref{isolated-finite} and \cref{thm:pure-dimension}. Since $\local{X}{a}/\langle f_1, \ldots, f_n \rangle \cong \local{C}{a}/f_1\local{C}{a}$, the fourth assertion follows from identity \eqref{mult-a-defn-2} and the definition of order. The fifth assertion is obvious in the case that either $\multf{a}$ or $\multp{f'_1,f_2}{f_n}{a}$ is zero or infinite. Otherwise assertion \eqref{finite-to-curve} implies that $C := V(f_2, \ldots, f_n)$ is a possibly non-reduced curve near $a$, and then assertion \eqref{isolated-finite} and \cref{cor:Macaulay} below imply that $f_1f'_1$ is a non zero-divisor in $\local{C}{a}$. Therefore assertion \eqref{mult-to-order} and \cref{order-properties} imply that $\multp{f_1f'_1, f_2}{f_n}{a} = \ord_a((f_1f'_1)|_C) = \ord_a(f_1|_C) + \ord_a(f'_1|_C) = \multf{a} + \multp{f'_1,f_2}{f_n}{a}$, as required.
\end{proof}

\begin{lemma} \label{cor:Macaulay}
Let $a$ be a nonsingular point of a variety $X$ of dimension $n$. Let $f_1, \ldots, f_m$, $m \leq n$, be regular functions on a neighborhood $U$ of $a$ on $X$ such that $a \in V(f_1, \ldots, f_m)$ and $V(f_1, \ldots, f_m) \cap U$ has pure dimension $n-m$. Then for every $j = 1, \ldots, m$, $f_j$ is a non zero-divisor in $\local{X}{a}/\langle f_1, \ldots, f_{j-1} \rangle$. 
\end{lemma}

\begin{proof}
If some $f_j$ is a zero-divisor in $\local{X}{a}/\langle f_1, \ldots, f_{j-1} \rangle$, then there is an open neighborhood $U'$ of $a$ in $X$ such that $f_j$ is a zero-divisor in $\kk[U']/\langle f_1, \ldots, f_{j-1} \rangle$. Due to \cref{affine-cover,local-equations} we can choose a neighborhood $U''$ of $a$ in $U'$ such that $U''$ is an affine variety in $\kk^{n+r}$ and the ideal $I(U'')$ of $U''$ in $\kk[x_1, \ldots, x_{n+r}]$ is generated by $m$ polynomials $g_1, \ldots, g_r$. Applying Macaulay's unmixedness theorem (\cref{thm:Macaulay}) to $g_1, \ldots, g_r, f_1, \ldots, f_m$, we see that $f_j$ is a non zero-divisor in $\kk[U'']/\langle f_1, \ldots, f_{j-1} \rangle$, which is a contradiction.
\end{proof}

\subsection{Intersection multiplicity in a family}
Let $X$ be a nonsingular affine variety of dimension $n$. Let $h_i := \sum_j \phi_{i,j}(t)g_{i,j}$, $i = 1, \ldots, n$, where $\phi_{i,j}$ are rational functions in an indeterminate $t$ and $g_{i,j}$ are regular functions on $X$. Let $T$ be the set of all $\epsilon \in \kk$ such that each $\phi_{i,j}$ is defined at $\epsilon$ (in other words, $T$ is the complement in $\kk$ of the poles of $\prod_{i,j} \phi_{i,j}$). For each $\epsilon \in T$, $h_{\epsilon,i} := h_i|_{t = \epsilon}$ are regular functions on $X$ which can be thought of ``deformations'' of $h_{\epsilon_0, i}$ for some fixed $\epsilon_0 \in T$. In this section we describe how, as $\epsilon$ varies in $T$, the multiplicities of $h_{\epsilon, i}$, $i = 1, \ldots, n$, change {\em locally} at a point of $X$, and {\em globally} on all of $X$. 

\begin{example}\label{example:mult-deformation-0}
Let $h(x,t) := (x+2)(t-2)(xt-2)$ (see \cref{fig:mult-deformation}) over a field $\kk$ of characteristic different from $2$ or $3$, and let $b_0 := -2 \in \kk$.  We compute the multiplicity at $x = b_0$ of $h(x, t)$ for different values of $t$. If $\epsilon \not\in \{-1, 2\}$, the multiplicity  at $x = -2$ of $h|_{t = \epsilon}$ is the same as the multiplicity of $x+2$, which is $1$. For $\epsilon = 2$, $h|_{t = \epsilon}$ is identically zero, and therefore has infinite multiplicity everywhere. On the other hand, for $\epsilon = -1$, the multiplicity at $x = -2$ is 
\begin{align*}
[h|_{t = -1}]_{x = -2} 
	&= [3(x+2)^2]_{x = -2} = 2[x+2]_{x = -2} = 2
\end{align*}
Note that the point $(x, \epsilon) = (-2, -1)$ is on two distinct irreducible components of the curve $h(x,t) = 0$, namely $x = -2$ and $xt = 2$. 
\end{example}

\begin{figure}
\begin{center}
\def\colorzero{blue}
\def\colorone{red}
\begin{tikzpicture}[scale=0.45]

\def\nsamples{103}
\def\comphy{2}
\def\compvx{-1}
\def\maxx{4}
\def\maxy{4}

\def\xmin{-\maxx}
\def\xmax{\maxx}
\def\ymin{-\maxy}
\def\ymax{\maxy}
\def\xminleft{\xmin}
\def\xmaxleft{-0.1}
\def\xminright{0.1}
\def\xmaxright{\xmax}

\pgfmathsetmacro{\minx}{1/\maxx};

\draw[\colorzero, thick, smooth, samples=\nsamples, domain=\xminleft:-\minx] plot(\x,{1/\x});
\draw[\colorzero, thick, smooth, samples=\nsamples, domain=\minx:\xmaxright] plot(\x,{1/\x});   
\draw[\colorzero, thick, samples=2, domain=\xmin:\xmax] plot ({\x} ,{\comphy});
\draw[\colorzero, thick, samples=2, domain=\ymin:\ymax] plot ({\compvx}, {\x});
\draw [\colorone, thick, dashed] (\xmin, {1/\compvx}) -- (\xmax, {1/\compvx});
\draw [\colorone, thick, dashed] (\xmin, 0) -- (\xmax, 0);

\draw [\colorone] (\xmin, 0) node [left] {\picfontsize $t = 0$};
\draw [\colorone] (\xmin, {1/\compvx}) node [left] {\picfontsize $t = -1$};
\draw [\colorzero] (\xmin, \comphy) node [left] {\picfontsize $t = \comphy$};
\draw [\colorzero] (\compvx, \ymin) node [below] {\picfontsize $x = -2$};
\end{tikzpicture}
\caption{The set of zeroes of $h(x,t) = (x+2)(t-2)(xt - 2)$}
\label{fig:mult-deformation}
\end{center}
\end{figure}

In \cref{example:mult-deformation-0} the multiplicity of $h(x,t)$ at $x = b_0$ is the same for almost all values of $t$, and it jumps to higher values for $t = \epsilon$ only when either the point $(b_0, \epsilon)$ is non-isolated in $V(h(x,t) = 0)$, or $(b_0, \epsilon)$ is also on an irreducible component of $V(h(x,t) = 0)$ different from the ``vertical line'' $\{b_0\} \times \kk$. In \cref{mult-deformation-0} below we show that this is true in general, and then in \cref{mult-deformation-global} we describe the ``global'' analogue of \cref{mult-deformation-0}. In the statements and proof of these results, $t$ will be a coordinate on $\kk$, and we will treat the $h_j(x,t)$ as regular functions on $X \times T$. Let $V := V(h_1, \ldots, h_n)$ be the closed subscheme of $X \times T$ determined by $h_1, \ldots, h_n$, and for each $\epsilon \in T$, let $Z_\epsilon$ be the subvariety $\supp(V) \cap (X \times \{\epsilon\})$; in other words, $Z_{\epsilon}$ is the set of zeroes of $h_1, \ldots, h_n, t - \epsilon$. 

\begin{thm} \label{mult-deformation-0}
Let $b_0 \in X$, and $m: T \to \zz \cup \{\infty\}$ be the function given by $\epsilon \mapsto \multp{h_{\epsilon,1}}{h_{\epsilon,n}}{b_0}$. Let $m^* := \min\{m(\epsilon): m \in T\}$ and $\tilde T := \{\epsilon \in T:$ either  $(b_0,\epsilon) \not\in Z_\epsilon$ or $(b_0, \epsilon)$ is isolated in $Z_\epsilon\}$. Then 
\begin{enumerate}
\item \label{tilde-T-open} $\tilde T$ is a Zariski open subset of $T$.
\item \label{m^*-finite-condition} $m^* < \infty$ if and only if $\tilde T$ is nonempty. 
\item \label{m^*-zero-condition} $m^* = 0$ if and only if there is $\epsilon \in T$ such that $(b_0, \epsilon) \not \in Z_\epsilon$. 
\item \label{m^*-open} Assume $m^* < \infty$. Then $\{\epsilon \in T: m(\epsilon) = m^*\}$ is a nonempty Zariski open subset of $\tilde T$.
\item \label{mult-not-min} Assume $m^* < \infty$. Then for each $\epsilon \in T$, $m(\epsilon) > m^*$ if and only if one of the following is true:
\begin{enumerate}
\item $\epsilon \not \in \tilde T$, i.e.\ $(b_0,\epsilon)$ is a non-isolated zero of $h_{\epsilon,1}, \ldots, h_{\epsilon, n}$, or
\item there is an irreducible component of $V(h_1, \ldots, h_n) \subset X \times \kk$ containing $(b_0, \epsilon)$ other than the ``vertical line'' $\{b_0\} \times \kk$. 
\end{enumerate}
\end{enumerate}
\end{thm}

\begin{proof}
The first assertion follows from the observation (due to \cref{thm:pure-dimension}) that $(b_0, \epsilon)$ is isolated in $Z_\epsilon$ if and only if $V$ has pure dimension one near $(b_0, \epsilon)$. Assertions \eqref{m^*-finite-condition} and \eqref{m^*-zero-condition} are immediate consequences of the definition of intersection multiplicity. Now we prove the last two assertions. If there is $\epsilon \in T$ such that $b_0$ is not a common zero of $h_{\epsilon, 1}, \ldots, h_{\epsilon,n}$, then $m^* = 0$, and both assertions \eqref{m^*-open} and \eqref{mult-not-min} are true. So assume $0 < m^* < \infty$. Then \cref{int-mult-curve,thm:pure-dimension} imply that $C_0 := \{b_0\} \times T$ is an irreducible component of the subscheme $V$ of $X \times T$. Let 
\begin{align*}
T_0 
	:= \{\epsilon \in T: (b_0,\epsilon)\ \text{is isolated in}\ 		  Z_\epsilon\ \text{and $C_0$ is the only component of $V$ containing}\ (b_0,\epsilon)\}
\end{align*}
Then $T_0 \subset \tilde T$ is a nonempty Zariski open subset of $T$, and there is a Zariski open neighborhood $U_0$ of $C_0 \cap (X \times T_0) = \{b_0\} \times T_0$ in $X \times T_0$ such that the scheme theoretic intersection $D_0 := V \cap  U_0$ has pure dimension one, and $\supp(D_0) = \{b_0\} \times T_0$. Pick distinct points $\epsilon_1, \epsilon_2 \in T_0$. For each $j$, \cref{cor:Macaulay} implies that $t - \epsilon_j$ is not a zero-divisor in $\local{D_0}{(b_0, \epsilon_j)} = \local{V}{(b_0, \epsilon_j)}$, which implies that $f := (t-\epsilon_1)/(t-\epsilon_2)$ is an {\em invertible rational function} on $D_0$. Embed $\kk \into \pp^1$, and let $\bar C_0 := \{b_0\} \times \pp^1$ be the closure of $C_0$ in $X \times \pp^1$. \Cref{nagata1} implies that each of $f, t - \epsilon_1, t - \epsilon_2$ can be extended to an invertible rational function on a closed subscheme $\bar D_0$ of $X \times \pp^1$ such that $\supp(\bar D_0) = \bar C_0$ and $D_0$ is (isomorphic to) an open subscheme of $\bar D_0$. The following is an immediate consequence of \cref{reduced-invertible}.

\begin{proclaim} \label{claim:invertible-p^1-infty}
$f$ is invertible in $\local{\bar D_0}{(b_0,\infty)}$, where $\infty$ is the only point of $\pp^1 \setminus \kk$. For each $j$, $t - \epsilon_j$ is invertible in $\local{\bar D_0}{(b_0, \epsilon)}$ for each $\epsilon \in \kk \setminus \{\epsilon_j\}$. \qed
\end{proclaim}

\Cref{order-properties,claim:invertible-p^1-infty} imply that 
\begin{align*}
\ord_{(b_0, \infty)}((t-\epsilon_1)|_{\bar D_0})
	= \ord_{(b_0, \infty)}((t-\epsilon_2)|_{\bar D_0})
\end{align*}
\Cref{order-curve,zero-sum,int-mult-curve,claim:invertible-p^1-infty} then imply that for each $\epsilon \in T_0$, 
\begin{align*}
\multp{h_{\epsilon,1}}{h_{\epsilon,n}}{b_0}
 	&=  \multp{t-\epsilon, h_1}{h_n}{(b_0,\epsilon)}
	=  \ord_{(b_0,\epsilon)}((t-\epsilon)|_{\bar D_0})  \\
	&=  \mu_{\bar C_0}(\bar D_0) \ord_{(b_0,\epsilon)}((t-\epsilon)|_{\bar C_0}) 
	= \mu_{\bar C_0}(\bar D_0)
\end{align*}
To complete the proof it suffices to show that $\multp{h_{\epsilon,1}}{h_{\epsilon,n}}{b_0} > \mu_{\bar C_0}(\bar D_0)$ whenever $\epsilon \in  T \setminus T_0$. This inequality is clear if $(b_0,\epsilon)$ is non-isolated in $Z_\epsilon$, so assume $(b_0,\epsilon)$ is an isolated point of $Z_\epsilon$ and there are irreducible components of $V$ containing $(b_0, \epsilon)$ other than $C_0$; denote them by $C_1, \ldots, C_k$. Choose a Zariski open neighborhood $W$ of $(b_0,\epsilon)$ in $X \times T$ such that $C_0 \cap W, \ldots, C_k \cap W$ are the only irreducible components of the open subscheme $V \cap W$ of $V$. Then 
\begin{align*}
\multp{h_{\epsilon,1}}{h_{\epsilon,n}}{b_0}
	&= \multp{t-\epsilon, h_1}{h_n}{(b_0,\epsilon)}
	=  \ord_{(b_0,\epsilon)}((t-\epsilon)|_{V \cap W})  \\
	&= \sum_{i=0}^k \mu_{C_i \cap W}(V \cap W)  \ord_{(b_0,\epsilon)}((t-\epsilon)|_{C_i \cap W}) 
	> \mu_{C_0 \cap W}(V \cap W) 
\end{align*}
Since $\local{V \cap W}{C_0 \cap W} = \local{\bar D_0}{C_0}$ (see \cref{scheme:local-section}), it follows that $\mu_{C_0 \cap W}(V \cap W) = \mu_{\bar C_0}(\bar D_0)$, which completes the proof.
\end{proof}

Now we prove a global counterpart of \cref{mult-deformation-0}. To motivate the statement of this result, we compute for different values of $\epsilon$ the number (counted with multiplicities) of isolated solutions of $h|_{t = \epsilon}$, for $h(x,t) := (x+2)(t-2)(xt-2)$ from \cref{example:mult-deformation-0}. It is straightforward to check (see \cref{fig:mult-deformation}) that if $\epsilon \not\in\{0,2\}$, then this number is $2$; indeed, if $\epsilon \not\in \{-1, 0, 2\}$, then $h|_{t = \epsilon}$ has two solutions of multiplicity one: $(-2, \epsilon), (2/\epsilon, \epsilon)$, and if $\epsilon = -1$, then $h|_{t = \epsilon}$ has one solution of multiplicity two: $(-2, -1)$. On the other hand, if $\epsilon = 2$, the polynomial $h|_{t = \epsilon}$ is identically zero on $\kk$, and therefore has {\em zero} isolated solutions. Finally, for $\epsilon = 0$ there is only one point of multiplicity one on $h|_{t = \epsilon} = 0$, namely the point $(-2, \epsilon)$; the other solution $(2/\epsilon, \epsilon)$ ``goes to infinity'' at $t = 0$. In particular, the total number of isolated solutions of $h|_{t = \epsilon}$ is equal for almost all values of $\epsilon$, and can only drop in exceptional cases when some of the solutions become non-isolated or run to infinity. \Cref{mult-deformation-global} below states that this is also the case in general; in particular, ``minimum'' in the local case (i.e.\ \cref{mult-deformation-0}) becomes ``maximum'' in the global case (i.e.\ \cref{mult-deformation-global}). The notation (in particular the meaning of $V$ and $Z_\epsilon$) below remains unchanged from \cref{mult-deformation-0}.

\begin{thm} \label{mult-deformation-global}
Let $C \subset X \times T$ be the union of all irreducible components of $V$ containing at least one isolated point of $Z_\epsilon$ for some $\epsilon \in T$. 
\begin{enumerate}
\item \label{onely-isolated} Either $C$ is empty, or it has pure dimension one.
\item \label{horizontally-finite}  $C_\epsilon := C \cap  (X \times \{\epsilon\})$ is finite for every $\epsilon \in T$. 
\end{enumerate}
Now assume $C$ is not empty. 
\begin{enumerate}[resume]
\item \label{openly-isolated} Let $T^*$ be the set of all $\epsilon \in T$ such that all points on $C_\epsilon$ are isolated in $Z_\epsilon$. Then $T^*$ is a nonempty Zariski open subset of $\kk$.
\item  \label{openly-max} For each $\epsilon \in T$, let $\tilde C_\epsilon := \{(b,\epsilon) \in X \times T: (b,\epsilon)$ is isolated in $Z_\epsilon\}$. The function $M: T \to \zz$ given by $\epsilon \mapsto \sum_{(b,\epsilon) \in \tilde C_\epsilon} \multp{h_{\epsilon,1}}{h_{\epsilon,n}}{b}$ achieves the maximum on a nonempty Zariski open subset of $T^*$.
\item \label{non-max-condition} If $\epsilon \in T$, then $M(\epsilon)$ fails to attain the maximum if and only if at least one of the following is true:
\begin{enumerate}
\item $\epsilon \not\in T^*$, i.e.\ there is a point on $C_\epsilon$ which is a non-isolated zero of $h_{\epsilon,1}, \ldots, h_{\epsilon, n}$, or
\item $C$ ``has a point at infinity at $t = \epsilon$'', i.e.\ if $\bar X$ is a projective compactification of $X$ and $\bar C$ is the closure of $C$ in $\bar X \times \pp^1$, then $\bar C \cap ((\bar X \setminus X) \times \{\epsilon\}) \neq \emptyset$. 
\end{enumerate}
\end{enumerate}
\end{thm}

\begin{proof}
If $(b, \epsilon)$ is an isolated point of $Z_\epsilon$, which is defined by $n+1$ regular functions $h_1, \ldots ,h_n, t - \epsilon$ on a variety of dimension $n+1$, \cref{thm:pure-dimension} implies that $V = V(h_1, \ldots, h_n)$ has pure dimension one near $(b, \epsilon)$, which proves assertion \eqref{onely-isolated}. If assertion \eqref{horizontally-finite} does not hold, then assertion \eqref{onely-isolated} implies that there is $\epsilon \in T$ such that $X \times \{\epsilon\}$ contains an irreducible component of $C$. But then no point on this component is isolated in $Z_\epsilon$, contradicting the definition of $C$. This proves assertion \eqref{horizontally-finite}. For assertion \eqref{openly-isolated}, let $Y$ be the union of the irreducible components of $V$ not contained in $C$. Since $\dim(C) = 1$, it follows that $C \cap Y$ is a finite set. If $\{(b'_j, \epsilon'_j)\}_j$ are the points in this intersection, then note that $T^* = T\setminus \{\epsilon'_j\}_j$. It remains to prove the last two assertions. Let 
\begin{align*}
\tilde C 
	&:= \bigcup_{\epsilon} \tilde C_\epsilon 
	= \{(b,\epsilon) \in C : (b, \epsilon)\ \text{is isolated in}\ Z_\epsilon \}
\end{align*}
Then $\tilde C$ is a Zariski open subset of $C$, i.e.\ there is a Zariski open subset $U$ of $X \times T$ such that $\tilde C = U \cap C$. Let $\tilde D := V \cap U$ be the corresponding open subscheme of $V$. Now choose a projective compactification $\bar X$ of $X$. Let $\bar C$ be the closure of $C$ in $\bar X \times \pp^1$. For every $\epsilon \in \kk$ and every $(b, \epsilon) \in \tilde C$, \cref{cor:Macaulay} implies that $t - \epsilon$ is a non zero-divisor in $\local{\tilde D}{(b, \epsilon)} = \local{V}{(b, \epsilon)}$. If $\epsilon_1, \epsilon_2$ are distinct elements of $\kk$, then \cref{nagata1} implies that there is a closed subscheme $\bar D$ of $\bar X \times \pp^1$ containing $\tilde D$ as an open subscheme such that $\supp(\bar D) = \bar C$, and each of $t - \epsilon_1, t-\epsilon_2, f := (t- \epsilon_1)/(t-\epsilon_2)$ extends to an {\em invertible} rational function on $\bar D$. The following is an immediate consequence of \cref{reduced-invertible}.

\begin{proclaim} \label{claim:invertible-over-p^1-infty}
Denote the only point of $\pp^1 \setminus \kk$ by $\infty$. For each $b \in \bar X$ such that $(b, \infty) \in \bar C$, $f$ is invertible in $\local{\bar D}{(b,\infty)}$. For each $j$, $t - \epsilon_j$ is invertible in $\local{\bar D}{(b, \epsilon)}$ for each $(b, \epsilon) \in \bar C$ such that $\epsilon \in \kk \setminus \{\epsilon_j\}$. \qed
\end{proclaim}

\Cref{order-properties,claim:invertible-over-p^1-infty} imply that for each $\epsilon_1, \epsilon_2 \in \kk$,
\begin{align*}
\sum_{b:(b,\infty) \in \bar C} \ord_{(b, \infty)}((t-\epsilon_1)|_{\bar D})
	= \sum_{b:(b,\infty) \in \bar C} \ord_{(b, \infty)}((t-\epsilon_2)|_{\bar D})
\end{align*}
\Cref{order-curve,zero-sum,int-mult-curve,claim:invertible-over-p^1-infty} then imply that for each $\epsilon \in \kk$, 
\begin{align*}
M^*	
	&:= \sum_{b:(b,\epsilon) \in \bar C} \ord_{(b, \epsilon)}((t-\epsilon)|_{\bar D})
	= \sum_{b:(b,\epsilon) \in \bar C} \sum_i \mu_{C_i}(\bar D) \ord_{(b, \epsilon)}((t-\epsilon)|_{C_i})
\end{align*}
is constant, where the $C_i$ are irreducible components of $C$. On the other hand, \cref{order-curve,int-mult-curve} imply that for all $\epsilon \in T$, 
\begin{align*}
M(\epsilon)
	&= \sum_{b:(b,\epsilon) \in \tilde C_\epsilon}  \multp{h_{\epsilon,1}}{h_{\epsilon,n}}{b}
	= \sum_{b:(b,\epsilon) \in \tilde C_\epsilon}  \multp{t-\epsilon, h_1}{h_n}{(b,\epsilon)} \\
	&= \sum_{b:(b,\epsilon) \in \tilde C_\epsilon} \ord_{(b,\epsilon)}((t-\epsilon)|_D) 
	= \sum_{b:(b,\epsilon) \in \tilde C_\epsilon}  \sum_i \mu_{C_i}(\bar D) \ord_{(b, \epsilon)}((t-\epsilon)|_{C_i})
\end{align*}
Since $t - \epsilon$ is regular and has a zero at each point of $\bar C \cap (X \times \{\epsilon\})$, it follows that $M^* \geq M(\epsilon)$; moreover, $M^* > M(\epsilon)$ if and only if $\tilde C_\epsilon \subsetneqq \bar C \cap (\bar X \times \{\epsilon\})$. By construction of $\bar C$, the latter condition is true if and only if at least one of the conditions of assertion \eqref{non-max-condition} holds. Since these conditions hold at at most finitely many points of $C$, assertion \eqref{openly-max} also holds.
\end{proof}

\section{Intersection multiplicity of complete intersections} \label{complete-mult-section}
\index{Intersection multiplicity!of complete intersections}
Let $f_1, \ldots, f_k$ be regular functions on a nonsingular variety $X$ of dimension $n \geq k$ such that $V(f_1, \ldots, f_k)$ has an irreducible component $Z$ of dimension $n-k$. Then there is a Zariski open subset $U$ of $X$ such that $Z \cap U = V(f_1, \ldots, f_k) \cap U$, and $Z \cap U$ has pure dimension $n-k$ (in particular, $Z \cap U$ is nonempty). Let $V$ be the {\em closed subscheme} of $U$ defined by (the ideal generated by) $f_1, \ldots, f_k$. Then $Z \cap U = \supp(V)$, so that we can define the {\em multiplicity} $\mu_{Z \cap U}(V)$ of $Z \cap U$ in $V$ as in \cref{nrcurve-section} (in the paragraph preceding \cref{order-curve}). It is straightforward to check that $\mu_{Z \cap U}(V)$ does not depend on $U$; we say that $\mu_{Z \cap U}(V)$ is the {\em intersection multiplicity $\multp{f_1}{f_k}{Z}$ of $f_1, \ldots, f_k$ along $Z$}. If $k = n$, then $Z$ is a singleton $\{a\}$, and $\local{V}{Z} \cong \local{X}{a}/\langle f_1, \ldots, f_n \rangle$, and \cref{prop:local-field-finite-length} implies that $\mu_Z(V) = \dim_\kk( \local{X}{a}/\langle f_1, \ldots, f_n \rangle)$, so that the definition of $\multp{f_1}{f_n}{Z}$ from this section agrees with the definition from the preceding section. 

\begin{prop} \label{mult-chain-0}
Let $f_1, \ldots, f_{n-1} \in \kk[x_1, \ldots, x_n]$. Let $Y$ be the coordinate subspace $x_1 = \cdots = x_k = 0$ of $\kk^n$. Assume 
\begin{enumerate}
\item $Y$ is an irreducible component of $V(f_1, \ldots, f_k) \subset \kk^{n}$. 
\item $V(f_{k+1}|_Y, \ldots, f_{n-1}|_Y)$ has a one dimensional irreducible component $Z$.
\item $Z$ is not contained in any irreducible component of $V(f_1, \ldots, f_k)$ other than $Y$.
\end{enumerate}
For each $\epsilon = (0, \ldots, 0, \epsilon_{k+1}, \ldots, \epsilon_n) \in Y$, and each $j = 1, \ldots, k$, we write $f_{j,\epsilon}$ for the polynomial in $(x_1, \ldots, x_k)$ obtained by substituting $\epsilon_i$ for $x_i$ for $i = k+1, \ldots, n$. Then %
$$\multp{f_1}{f_{n-1}}{Z}
 	= \multzero{f_{1, \epsilon}}
 	{f_{k, \epsilon}}
 		\multp{f_{k+1}|_Y}{f_{n-1}|_Y}{Z}$$ %
 for generic $\epsilon \in Y$.  
\end{prop}

\begin{proof}
We prove this by induction on $n-k$. At first consider the case that $n-k = 1$. Then $Z$ is the $x_n$-axis, and for generic $\epsilon \in \kk$, $Z$ is the only irreducible component of $V(f_1, \ldots, f_{n-1})$ containing $a_\epsilon := (0, \ldots, 0, \epsilon)$. Since $\ord_{a_\epsilon}((x_n - \epsilon)|_{Z}) = 1$, \cref{order-curve} implies that 
\begin{align*}
\multp{f_1}{f_{n-1}}{Z} 
	&= \multp{x_n-\epsilon, f_1}{f_{n-1}}{a_\epsilon} \\
	&=  \dim_\kk (\kk[[x_1, \ldots, x_{n-1}, x_n - \epsilon]]/\langle x_n - \epsilon, f_1, \ldots, f_{n-1} \rangle ) \\
	&= \dim_\kk (\kk[[x_1, \ldots, x_{n-1}]]/\langle f_1|_{x_n = \epsilon}, \ldots, f_{n-1}|_{x_n = \epsilon} \rangle ) \\
	&= \multzero{ f_{1,\epsilon}}{ f_{n-1,\epsilon}}
\end{align*}
as required. In the general case, pick a nonsingular point $z = (0, \ldots, 0, z_{k+1}, \ldots, z_n)$ of $Z$. Then there is $j$, $k + 1 \leq j \leq n$, such that $(x_j - z_j)|_Z$ has order one at $z$ (\cref{prop:nonsingular-parameter-properties}, assertion \eqref{prop:nonsingular-parameter-coordinates}). Pick a generic $\epsilon_j \in \kk$. Assertion \eqref{prop:nonsingular-parameter-open} of \cref{prop:nonsingular-parameter-properties} implies that the set $V(x_j - \epsilon_j) \cap Z$ is nonempty and contains a nonsingular point $a$ of $Z$ such that $a$ is {\em not} in any other irreducible component of $V(f_1, \ldots, f_{n-1})$, and $\ord_a((x_j - \epsilon_j)|_Z) = 1$. Since $a$ is an isolated zero of $x_j - \epsilon_j, f_1, \ldots, f_{n-1}$, \cref{thm:pure-dimension} implies that $x_j-\epsilon_j, f_1, \ldots, f_{n-2}$ defines a possibly non-reduced curve $W$ near $a$. Let $W_1, \ldots, W_s$ be the irreducible components of $W$ and $\pi_i: \tilde W_i \to W_i$ be the desingularization. \Cref{order-curve,int-mult-curve} imply that
\begin{align*}
\multp{f_1}{f_{n-1}}{Z} 
	&= \multp{x_j-\epsilon_j, f_1}{f_{n-1}}{a}  
	= \ord_{a}(f_{n-1}|_{W}) \\
	&= \sum_i \mu_{W_i}(W) 
	\sum_{\tilde a \in \pi_i^{-1}(a)} \ord_{\tilde a}(\pi_i^*(f_{n-1}|_{W_i})) \\
	&= \sum_i \multp{f_1|_{x_j = \epsilon_j}}{f_{n-2}|_{x_j = \epsilon_j}}{W_i}
	\sum_{\tilde a \in \pi_i^{-1}(a)} \ord_{\tilde a}(\pi_i^*(f_{n-1}|_{W_i})) 
\end{align*}
Let $Y_{\epsilon_j} := Y \cap V(x_j - \epsilon_j)$. Then the inductive hypothesis implies that
\begin{align*}
\multp{f_1}{f_{n-1}}{Z} 
	&= \multzero{f_{1,\epsilon}}{f_{k,\epsilon}}
	\sum_i 	\multp{f_{k+1}|_{Y_{\epsilon_j}}}{f_{n-2}|_{Y_{\epsilon_j}}}{W_i} 
	\sum_{\tilde a \in \pi_i^{-1}(a)} \ord_{\tilde a}(\pi_i^*(f_{n-1}|_{W_i})) \\
	 &=  \multzero{f_{1,\epsilon}}{f_{k,\epsilon}}
	 	 	\ord_{a}(f_{n-1}|_{V}) 
\end{align*}
where $V$ is the closed subscheme $W \cap Y$ of $Y$. It follows that 
\begin{align*}
\multp{f_1}{f_{n-1}}{Z} 
	 &=  \multzero{f_{1,\epsilon}}{f_{k,\epsilon}}
	 	 	\multp{x_j - \epsilon_j, f_{k+1}|_Y}{f_{n-1}|_Y}{a} \\
	 &=  \multzero{f_{1,\epsilon}}{f_{k,\epsilon}}
	 	 	 	\multp{f_{k+1}|_Y}{f_{n-1}|_Y}{Z}
\end{align*}
as required.
\end{proof}

\begin{cor}\label{mult-chain-1}
Let the assumptions be as in \cref{mult-chain-0}. Let $a \in Y$ be such that $V := V(f_1, \ldots, f_{n-1})$ is purely one dimensional near $a$ and no irreducible component of $V$ containing $a$ is contained in any irreducible component of $V(f_1, \ldots, f_k)$ other than $Y$. Then for all $f_n \in \kk[x_1, \ldots, x_n]$
\begin{align*}
\multp{f_1}{f_n}{a} = \multzero{f_{1, \epsilon}}
 	{f_{k, \epsilon}}
 		\multp{f_{k+1}|_Y}{f_n|_Y}{a}
\end{align*}
for generic $\epsilon \in Y$.  
\end{cor}

\begin{proof}
Follows from \cref{int-mult-curve,mult-chain-0}. 
\end{proof}

\chapter{Convex polyhedra} \label{appolytopes}
\def\scalefactor{.3}
\def\colorAB{blue}
\def\colorBC{magenta}
\def\colorCA{black}
\def\colorA{orange}
\def\colorB{gray}
\def\colorC{yellow}
\def\opazero{0.5}
\def\colordot{red}
\def\colorzero{green}
\def\colorone{blue}
\def\colortwo{black}
\def\colornu{purple}
\def\picfontsize{\small}

A ``polytope'' has two equivalent definitions: a convex hull of finitely many points, or a bounded intersection of finitely many ``half spaces.'' In \cref{basic-convex-section,poly-characterizection} we prove the equivalence of these definitions after introducing the basic terminology. The rest of the chapter is devoted to different properties of polytopes which are implicitly or explicitly used in the forthcoming chapters. 

\section{Basic notions} \label{basic-convex-section}
In this chapter we treat the spaces $\rr^n$, $n \geq 0$, as vector spaces over $\rr$ equipped with the Euclidean topology, and deal only with ``affine maps'' between these spaces. Recall that a map $\phi: \rr^n \to \rr^m$ is \index{Affine!map}{\em affine} if there is $\beta \in \rr^m$ such that $\phi(\cdot) = \beta + \phi_0(\cdot)$ for some {\em linear} map $\phi_0: \rr^n \to \rr^m$. Given $S \subset \rr^n$ and $T \subset \rr^m$, and a map $\phi: S \to T$, we say that $\phi$ is {\em affine} if it is the restriction of an affine map from $\rr^n \to \rr^m$; we say that $\phi$ is an {\em affine isomorphism} if $\phi$ is affine and bijective. An \index{Affine!subspace}{\em affine subspace} $A$ of $\rr^n$ is a subset of $\rr^n$ which is the image of an affine map; in other words it is simply a translation of a linear subspace $L$ of $\rr^n$. The \index{Dimension!of an affine subspace}{\em dimension} $\dim(A)$ of $A$ is the dimension of $L$ as a vector space over $\rr$. A \index{Hyperplane}{\em hyperplane} in $\rr^n$ is an affine subspace of dimension $n-1$. The \index{Affine!hull}{\em affine hull} $\aff(S)$ of a set $S$ of $\rr^n$ is the smallest affine subspace of $\rr^n$ containing $S$; alternatively, if $L$ is the linear subspace of $\rr^n$ spanned by all elements of the form $\alpha - \beta$ such that $\alpha, \beta \in S$, then $\aff(S) = L + \alpha$ for any $\alpha \in S$. In \cref{fig:min-max}, it is straightforward to check that $\aff(S) = \rr^2$, $\aff(\{C\}) = \{C\}$, and since $A, B, E$ are collinear, $\aff(\{A,B, E\})$ is the (unique) line $L$ through these points.

\begin{figure}[h]
\def\xmin{-0.5}
\def\xmax{16.5}
\def\ymin{-4.5}
\def\ymax{5.5}
\def\tx{0}
\def\ty{5}
\def\tw{5cm}

\tikzstyle{dot} = [\colordot, circle, minimum size=4pt, inner sep = 0pt, fill]

\begin{center}

\begin{tikzpicture}[scale=\scalefactor]

\node[dot] (A) at (1,-4) {};
\node[dot] (B) at (16,1) {};
\node[dot] (C) at (10,5) {};
\draw[thick, \colorone, fill=\colorzero, opacity=\opazero ] (A.center) --  (B.center) -- (C.center) -- cycle;
\node[dot] (D) at (11,1) {};
\node[dot] (E) at (7,-2) {};
\coordinate (nu) at (1,-3);
\coordinate (nuperp) at (3,1);

\node[anchor = north] at (A) {\picfontsize $A$};
\node[anchor = north] at (B) {\picfontsize $B$};
\node[anchor = south] at (C) {\picfontsize $C$};
\node[anchor = east] at (D) {\picfontsize $D$};
\node[anchor = north] at (E) {\picfontsize $E$};

\draw[thick, \colortwo] ($(A)-(nuperp)$) --  ($(B)+(nuperp)$);
\coordinate (O) at ($(A)!0.5!(B)$);
\draw [thick, \colornu, ->] (O) -- ($(O) + 0.4*(nu)$);
\node[anchor = north] at ($(O) + 0.4*(nu)$) {\picfontsize $\nu$};
\node[anchor = north] at ($(B)+(nuperp)$) {\picfontsize $L$};

\draw[thick, \colortwo, dashed] ($(A)-(nuperp)-1.8*(nu)$) --  ($(B)+(nuperp)-1.8*(nu)$);
\draw [thick, \colornu, ->] ($(O)-1.8*(nu)$) -- ($(O) - 1.4*(nu)$);
\end{tikzpicture}

\caption{
	$S = \{A, B, C, D, E\}$, $\In_\nu(S) = \{C\}$, $\ld_\nu(S) = \{A, E, B\}$
} \label{fig:min-max}
\end{center}
\end{figure}
It is straightforward to check that an affine map preserves affine subspaces, affine hulls, and if the map is injective, then also the dimension of affine hulls (\cref{exercise:affine-image}). Given $\alpha \in \rr^n$ and $\nu \in \rnstar$, we write $\langle \nu, \alpha \rangle$ for the ``value of $\nu$ at $\alpha$,'' and write $\nu^\perp := \{\alpha \in \rr^n: \langle \nu, \alpha \rangle = 0\}$. Define
\begin{align*}
\min_S(\nu)
	&:= \min\{\langle \nu, \alpha \rangle: \alpha \in S\},\ \text{provided the minimum exists.} \\
\max_S(\nu)
	&:= \max\{\langle \nu, \alpha \rangle: \alpha \in S\},\ \text{provided the maximum exists.}  \\
\In_\nu(S)
	&:= \{\alpha \in S: \langle \nu, \alpha \rangle =  \min_S(\nu)\},\ \text{provided $\min_S(\nu)$ exists.}  \\
\ld_\nu(S)
	&:= \{\alpha \in S: \langle \nu, \alpha \rangle =  \max_S(\nu)\},\ \text{provided $\max_S(\nu)$ exists.}
\end{align*}
See \cref{fig:min-max} for an illustration of these notions for a planar set. Note that in \cref{fig:min-max} we depicted $\nu \in \rnnstar{2}$ on $\rr^2$ by identifying it with an element (modulo ``parallel translations'') in $\rr^2$ via the ``dot product.'' A set is \index{Convex!set}{\em convex} if it contains the line segment joining any two points in it. The \index{Convex!hull}{\em convex hull} $\conv(S)$ of $S$ is the smallest convex set containing $S$. In \cref{fig:min-max} the convex hull of the $5$ points is the green triangle. Given $\alpha_1, \ldots, \alpha_k \in \rr^n$, an expression of the form $\sum_{j=1}^k \epsilon_j \alpha_j$, where the $\epsilon_j$ are nonnegative numbers whose sum is $1$, is called a \index{Convex!combination}{\em convex combination} of the $\alpha_j$. The set of convex combinations of two elements in $\rr^n$ is precisely the line segment joining them, and from this observation it can be shown that the convex hull of a set consists of the convex combinations of its points (\cref{exercise:conv-rep}), i.e.\
\begin{align}
\conv(S) = \{\sum_{j=1}^k \epsilon_j \alpha_j:
			k \geq 0,\
			\alpha_j \in S,\ \epsilon_j \geq 0\
			\text{for each}\ j,\ 1 \leq j \leq k,\
			\text{and}\ \sum_{j=1}^k \epsilon_j  = 1			
		\}
\label{conv-rep}
\end{align}
Given a nonnegative real number $r$ and $S \subset \rr^n$, we define $rS := \{r\alpha : \alpha \in S\}$, i.e.\ $rS$ is the ``dilation of $S$ by a factor of $r$.'' We say that $S$ is a \index{Cone}{\em cone} if it is ``dilation invariant,'' i.e.\ if $rS \subseteq S$ for each $r \geq 0$. The \index{Cone!convex}\index{Convex!cone}{\em convex cone generated by $S$} is defined to be
\begin{align}
\cone(S) := \{\sum_{j=1}^k \epsilon_j \alpha_j:
			k \geq 0,\
			\alpha_j \in S\
			\text{and}\ \epsilon_j  \geq 0\
			\text{for each}\ j,\ 1 \leq j \leq k
		\}
\label{cone-rep}
\end{align}
We say that $\cone(S)$ is \index{Finitely generated cone}{\em finitely generated} if $S$ is finite. A \index{Polyhedron}\index{Convex!polyhedron}{\em convex polyhedron} $\scrP$ is a subset of $\rr^n$ defined by finitely many linear inequalities, i.e.\ inequalities of the form $a_0 +  a_1x_1 + \cdots + a_nx_n \geq 0$. Geometrically, a polyhedron is a finite intersection of ``half-spaces,'' where a ``half-space'' is the set of all points on one side of a hyperplane - see \cref{fig:two-hedra}. If $\scrP$ is bounded, we call it a convex \index{Polytope}\index{Convex!polytope}{\em polytope}, and if it is a cone, we call it a convex \index{Polyhedral cone}\index{Convex!polyhedral cone}{\em polyhedral cone}. A convex polyhedral cone is equivalently a set defined by finitely many linear inequalities with {\em zero} constant term (\cref{exercise:cone-inequalities}). In this book we only consider convex polyhedra, and therefore we will simply write ``polyhedra,'' ``polytopes,'' ``cones,'' ``polyhedral cones'' to mean ``convex polyhedra,'' ``convex polytopes,'' ``convex cones,'' ``convex polyhedral cones'' respectively. The \index{Dimension!of a convex polyhedron}{\em dimension} $\dim(\scrP)$ of a polyhedron $\scrP$ is the dimension of its affine hull. \Cref{fig:two-hedra} depicts a few two dimensional convex polyhedra.

\begin{center}
\begin{figure}[h]

\def\scalefactor{.3}

\def\tinyshift{0.05} 
\def\xmin{-2}
\def\xmaxwithoutshift{8}
\pgfmathsetmacro{\xmax}{\xmaxwithoutshift+\tinyshift}
\def\ymin{-5}
\def\ymax{8}

\def\tx{0}
\def\ty{\ymin}

\def\xshift{5}
\xdefinecolor{c1}{HTML}{D4ACEE}
\xdefinecolor{c2}{HTML}{EEc6AC}
\xdefinecolor{c3}{HTML}{EDEEAC}

\tikzstyle{dot} = [\colordot, circle, minimum size=4pt, inner sep = 0pt, fill]

\begin{tikzpicture}[scale=\scalefactor]
\pgfmathsetmacro\shiftone{\xshift + \xmax -\xmin};

\coordinate (A) at (-2, -1);
\coordinate (B) at (3,-3);
\coordinate (C) at (7,1);
\coordinate (D) at (5,6);
\coordinate (E) at (1,4);

\draw[fill=\colorzero, opacity=\opazero ] (A) --  (B) -- (C) -- (D) -- (E) -- cycle;

\draw[thick, \colorone ] (A) --  (B) -- (C) -- (D) -- (E) -- cycle;

\draw (\xmax, 0) -- (\xmin, 0);
\draw (0, \ymax) -- (0, \ymin);

\node[dot] (O) at (0,0){};
\node[anchor = south west] at (O) {\picfontsize $(0,0)$};
\node[anchor = west] at (\xmax, 0) {\picfontsize $x$};
\node[anchor = south] at (0, \ymax) {\picfontsize $y$};

\node [below right, text width= 0.75cm, align=left] at (\tx,\ty) {
	\picfontsize
	Polytope
};	

\begin{scope}[shift = {(\shiftone,0)}]
	\coordinate (A) at (-2, -1);
	\pgfmathsetmacro\Cy{-1 -0.4*(\xmax+2)};
	\coordinate (C') at (\xmax, \Cy);
	\coordinate (C') at (8,-5);
	\pgfmathsetmacro\Dy{4 + 0.5*(\xmax - 1)};
	\coordinate (D') at (\xmax, \Dy);
	\coordinate (E) at (1,4);
	
	\fill[\colorzero, opacity=\opazero ] (A) -- (C') -- (D') -- (E) -- cycle;
	
	\draw[thick, \colorone ] (A) -- (C');
	\draw[thick, \colorone ] (A) -- (E) -- (D');
	
	\draw (\xmax, 0) -- (\xmin, 0);
	\draw (0, \ymax) -- (0, \ymin);
	
	\node[dot] (O) at (0,0){};
	\node[anchor = south west] at (O) {\picfontsize $(0,0)$};
	\node[anchor = west] at (\xmax, 0) {\picfontsize $x$};
	\node[anchor = south] at (0, \ymax) {\picfontsize $y$};
	
	\node [below right, text width= 1.8cm, align=left] at (\tx,\ty) {
		\picfontsize
		Unbounded polyhedron
	};	
\end{scope}

\begin{scope}[shift = {(2*\shiftone,0)}]
	\pgfmathsetmacro\Cy{-0.4*\xmax};
	\coordinate (C') at (\xmax,\Cy);
	\pgfmathsetmacro\Dy{0.5*\xmax};
	\coordinate (D') at (\xmax, \Dy);
	\coordinate (O) at (0,0);
	

	
	\pgfmathsetmacro\xmindiffed{\xmin + 2*\tinyshift};

	\foreach \x in {\xmin,\xmindiffed,...,\xmaxwithoutshift} {
		\draw[color=c2, ultra thin] (\x, \ymax) -- (\x,{-0.4*\x});
		\draw[color=c1, ultra thin] ({\x + \tinyshift}, \ymin) -- ({\x + \tinyshift}, {0.5*(\x + \tinyshift)});
	}
	\draw[color=c2] (\xmin, {-0.4*\xmin}) -- (\xmax, {-0.4*\xmax});
	\draw[color=c1] (\xmin, {0.5*\xmin}) -- (\xmax, {0.5*\xmax});
	
	\draw[thick, \colorone ] (O) -- (C');
	\draw[thick, \colorone ] (O) -- (D');
	
	\draw (\xmax, 0) -- (\xmin, 0);
	\draw (0, \ymax) -- (0, \ymin);
	
	\node[dot] (extra) at (O) {};
	\node[anchor = north] at (O) {\picfontsize $(0,0)$};
	\node[anchor = west] at (\xmax, 0) {\picfontsize $x$};
	\node[anchor = south] at (0, \ymax) {\picfontsize $y$};
	
	\node [below right, text width= 2.7cm, align=left] at (\tx,\ty) {
		\picfontsize
		Polyhedral cone
	};	
\end{scope}

\end{tikzpicture}

\caption{
	Some planar convex polyhedra; ``half-planes'' defining the cone are also depicted.
} \label{fig:two-hedra}

\end{figure}
\end{center}

A \index{Strongly convex cone}\index{Cone!strongly convex}{\em strongly convex cone} $\scrC$ is a convex cone which does not contain any line through the origin, equivalently, for all $\alpha \in \scrC\setminus \{0\}$, $-\alpha \not\in \scrC$. The \index{Minkowski sum}{\em Minkowski sum} of two subsets $\scrP, \scrQ$ of $\rr^n$ is $\scrP + \scrQ := \{\alpha + \beta : \alpha \in \scrP,\ \beta \in \scrQ\}$. It is straightforward to check that the Minkowski sum of convex sets is also convex (\cref{exercise:minkonvex}); see \cref{fig:minkowski-sum} for some examples in $\rr^2$.  We now show that every convex polyhedral cone or finitely generated cone\footnote{We will see in \cref{poly-characterizection} that convex polyhedral cones and finitely generated cones are the same. However, we will use \cref{strongly-convex-prop} in the proof of this equivalence.} can be represented as the Minkowski sum of a linear subspace and a strongly convex cone.

\begin{figure}[h]
\begin{center}
\begin{tikzpicture}[scale=0.4]
\def\shifttwo{5}
\def\minx{-1.5}
\def\maxx{1.5}
\def\maxxx{1.5}
\def\miny{-1.5}
\def\maxy{1.5}
\def\tx{\minx}
\def\ty{\miny}
\def\bigshiftdhor{3}
\def\bigshiftdver{3}

\draw [gray,  line width=0pt] (\minx,\miny) grid (\maxx,\maxy);
\draw[ultra thick, \colorone] (0,0) -- (-1,1);
\node[anchor = west] at (\maxx, 0) {\picfontsize $x$};
\node[anchor = south] at (0, \maxy) {\picfontsize $y$};
\draw (\tx,\ty) node [below right] {\picfontsize $\scrP$};  	

\begin{scope}[shift={(\shifttwo,0)}]
	\def\tx{0}
	\draw [gray,  line width=0pt] (\minx,\miny) grid (\maxx,\maxy);
	\draw[ultra thick, \colorone] (0,0) -- (1,-1);
	\node[anchor = west] at (\maxx, 0) {\picfontsize $x$};
	\node[anchor = south] at (0, \maxy) {\picfontsize $y$};
	\draw (\tx,\ty) node [below right] {\picfontsize $\scrQ$};    	
\end{scope}

\begin{scope}[shift={(2*\shifttwo,0)}]
	\draw [gray,  line width=0pt] (\minx,\miny) grid (\maxx,\maxy);
	\draw[ultra thick, \colorone] (-1,1) -- (1,-1);
	\node[anchor = west] at (\maxx, 0) {\picfontsize $x$};
	\node[anchor = south] at (0, \maxy) {\picfontsize $y$};
	\draw (\tx,\ty) node [below right] {\picfontsize $\scrP+ \scrQ$};
\end{scope}

\begin{scope}[shift={(2*\shifttwo+\maxxx - \minx +\bigshiftdhor,0)}]
	\def\minx{-0.5}
	\def\maxx{2.5}
	\def\maxxx{4.5}
	
	\draw [gray,  line width=0pt] (\minx,\miny) grid (\maxx,\maxy);
	\draw[ultra thick, \colorone] (0,1) -- (1,0);
	\node[anchor = west] at (\maxx, 0) {\picfontsize $x$};
	\node[anchor = south] at (0, \maxy) {\picfontsize $y$};
	\draw (\tx,\ty) node [below right] {\picfontsize $\scrP$};  	
	
	\begin{scope}[shift={(\shifttwo,0)}]
		\def\tx{1}
		\draw [gray,  line width=0pt] (\minx,\miny) grid (\maxx,\maxy);
		
		\draw[ultra thick, \colorone] (1,0) -- (2,0);
		\node[anchor = west] at (\maxx, 0) {\picfontsize $x$};
		\node[anchor = south] at (0, \maxy) {\picfontsize $y$};
		\draw (\tx,\ty) node [below right] {\picfontsize $\scrQ$};    	
	\end{scope}
	
	\begin{scope}[shift={(2*\shifttwo,0)}]
		\def\tx{1}
		\draw [gray,  line width=0pt] (\minx,\miny) grid (\maxxx,\maxy);		
		\node[anchor = west] at (\maxxx, 0) {\picfontsize $x$};
		\node[anchor = south] at (0, \maxy) {\picfontsize $y$};
		\draw[fill=\colorzero, opacity=\opazero] (1,1) --  (2,1) -- (3,0) -- (2,0) -- cycle;
		\draw[thick, \colorone ] (1,1) --  (2,1) -- (3,0) -- (2,0) -- cycle;
		\draw (\tx,\ty) node [below right] {\picfontsize $\scrP+ \scrQ$};
	\end{scope}
\end{scope}		

\begin{scope}[shift={(0,-\maxy +\miny -\bigshiftdver)}]

	\draw [gray,  line width=0pt] (\minx,\miny) grid (\maxx,\maxy);
	\node[anchor = west] at (\maxx, 0) {\picfontsize $x$};
	\node[anchor = south] at (0, \maxy) {\picfontsize $y$};
	\draw (\tx,\ty) node [below right] {\picfontsize $\scrP$};  	
	\draw[fill=\colorzero, opacity=\opazero] (-1,0) --  (-1,1) -- (0,0) -- cycle;
	\draw[thick, \colorone ] (-1,0) --  (-1,1) -- (0,0) -- cycle;

	\begin{scope}[shift={(\shifttwo,0)}]
		\def\tx{0}
		\draw [gray,  line width=0pt] (\minx,\miny) grid (\maxx,\maxy);
		\node[anchor = west] at (\maxx, 0) {\picfontsize $x$};
		\node[anchor = south] at (0, \maxy) {\picfontsize $y$};
		\draw[fill=\colorzero, opacity=\opazero] (1,0) --  (1,-1) -- (0,0) -- cycle;
		\draw[thick, \colorone ] (1,0) --  (1,-1) -- (0,0) -- cycle;
		\draw (\tx,\ty) node [below right] {\picfontsize $\scrQ$};    	
	\end{scope}
	
	\begin{scope}[shift={(2*\shifttwo,0)}]
		\draw [gray,  line width=0pt] (\minx,\miny) grid (\maxx,\maxy);
		\node[anchor = west] at (\maxx, 0) {\picfontsize $x$};
		\node[anchor = south] at (0, \maxy) {\picfontsize $y$};
		\draw[fill=\colorzero, opacity=\opazero] (-1,0) --  (-1,1) -- (0,1) -- (1,0) -- (1,-1) -- (0, -1) -- cycle;
		\draw[thick, \colorone ] (-1,0) --  (-1,1) -- (0,1) -- (1,0) -- (1,-1) -- (0, -1) -- cycle;
		\draw (\tx,\ty) node [below right] {\picfontsize $\scrP+ \scrQ$};
	\end{scope}

	\begin{scope}[shift={(2*\shifttwo+\maxxx - \minx +\bigshiftdhor,0)}]
	
		\def\minx{-0.5}
		\def\maxx{2.5}
		\def\maxxx{4.5}
		
		\draw [gray,  line width=0pt] (\minx,\miny) grid (\maxx,\maxy);
		\node[anchor = west] at (\maxx, 0) {\picfontsize $x$};
		\node[anchor = south] at (0, \maxy) {\picfontsize $y$};
		\draw[fill=\colorzero, opacity=\opazero] (0,-1) --  (1,2) -- (2,-1) -- (1,0) -- cycle;
		\draw[thick, \colorone ] (0,-1) --  (1,2) -- (2,-1) -- (1,0) -- cycle;
		\draw (\tx,\ty) node [below right] {\picfontsize $\scrP$};  	
		
		\begin{scope}[shift={(\shifttwo,0)}]
			\def\tx{1}
			\draw [gray,  line width=0pt] (\minx,\miny) grid (\maxx,\maxy);
			
			\draw[ultra thick, \colorone] (1,0) -- (2,0);
			\node[anchor = west] at (\maxx, 0) {\picfontsize $x$};
			\node[anchor = south] at (0, \maxy) {\picfontsize $y$};
			\draw (\tx,\ty) node [below right] {\picfontsize $\scrQ$};    	
		\end{scope}
		
		\begin{scope}[shift={(2*\shifttwo,0)}]
			\def\tx{1}
			\draw [gray,  line width=0pt] (\minx,\miny) grid (\maxxx,\maxy);		
			\node[anchor = west] at (\maxxx, 0) {\picfontsize $x$};
			\node[anchor = south] at (0, \maxy) {\picfontsize $y$};
			\draw[fill=\colorzero, opacity=\opazero] (1,-1) --  (2,2) -- (3,2) -- (4,-1) -- (3,-1) -- (2.5, -0.5) -- (2,-1) -- cycle;
			\draw[thick, \colorone ] (1,-1) --  (2,2) -- (3,2) -- (4,-1) -- (3,-1) -- (2.5, -0.5) -- (2,-1) -- cycle;
			\draw (\tx,\ty) node [below right] {\picfontsize $\scrP+ \scrQ$};
		\end{scope}
	\end{scope}		
\end{scope}
\end{tikzpicture}
\caption{Minkowski sums of planar sets}  \label{fig:minkowski-sum}
\end{center}
\end{figure}

\begin{prop} \label{strongly-convex-prop}
Let $\scrC$ be a convex cone in $\rr^n$.
\begin{enumerate}
\item \label{strongly-convex:polyhedral} If $\scrC$ is a polyhedral cone, then there is a strongly convex polyhedral cone $\scrC'$ and a linear subspace $L$ of $\rr^n$ such that $\scrC = \scrC' + L$ and $\aff(\scrC') \cap L = \{0\}$.
\item  \label{strongly-convex:finite} If $\scrC$ is a finitely generated cone, then there is a strongly convex finitely generated cone $\scrC'$ and a linear subspace $L$ of $\rr^n$ such that $\scrC = \scrC' + L$ and $\aff(\scrC') \cap L = \{0\}$.
\item \label{dsum-linear} Let $L'$ be a linear subspace of $\rr^n$ such that $L' \cap \aff(C) = \{0\}$. Then $\scrC$ is a polyhedral cone (respectively, finitely generated cone) if and only if $\scrC + L'$ is a polyhedral cone (respectively, finitely generated cone).
\end{enumerate}
\end{prop}

\begin{proof}
Let $L$ be the (unique) maximal linear subspace of $\rr^n$ contained in $\scrC$ (see \cref{exercise:maxlinconvex}). After a linear change of coordinates if necessary, we may assume that $L$ is the coordinate subspace spanned by the first $k$ coordinates. At first assume $\scrC$ is a polyhedral cone. \Cref{exercise:cone-inequalities} implies that it is defined by finitely many inequalities of the form $a_{i,1}x_1 + \cdots + a_{i,n}x_n \geq 0$, $1 \leq i \leq m$. For each $j = 1, \ldots, n$, write $e_j$ for the $j$-th standard unit vector in $\rr^n$. Since $re_j \in L \subseteq \scrC$ for each $j = 1, \ldots, k$, and each $r \in \rr$, it follows that $a_{i,1} = \cdots = a_{i,k} = 0$ for each $i$. Let $\scrC'$ be the polyhedral cone on the $(n-k)$-dimensional coordinate subspace of $\rr^n$ spanned by $e_{k+1}, \ldots, e_n$ defined by the inequalities $a_{i,k+1}x_{k+1} + \cdots + a_{i,n}x_n = 0$. Then assertion \eqref{strongly-convex:polyhedral} holds with $\scrC'$ and $L$. Now assume $\scrC$ is the cone generated by finitely many elements $\alpha_1, \ldots, \alpha_m \in \rr^n$. Let $\alpha_i := (\alpha_{i,1}, \ldots, \alpha_{i,n})$. Since $\tilde \alpha_i := (-\alpha_{i,1}, \ldots, -\alpha_{i,k}, 0, \ldots, 0) \in L \subseteq \scrC$, it follows that $\alpha'_i := \alpha_i + \tilde \alpha_i = (0, \ldots, 0, \alpha_{i, k+1}, \ldots, \alpha_{i,n}) \in \scrC$ for each $i$. Then assertion \eqref{strongly-convex:finite} holds with $L$ and $\scrC' := \cone(\alpha'_1, \ldots, \alpha'_m)$. The proof of assertion \eqref{dsum-linear} follows via the same arguments as in the proof of the preceding assertions by choosing a system of coordinates on $\rr^n$ such that $L'$ is the coordinate subspace spanned by the first $k'$ coordinates, where $k' := \dim(L')$.
\end{proof}

\begin{rem} \label{equivalent-multiplication}
If $\scrA \subset \rr^n$ and $d$ is a positive integer, then $d\scrA$ has two natural interpretations: one is the dilation (which is how we defined it), and the other is the Minkowski sum of $d$ copies of $\scrA$. If $\scrA$ is convex, then these are equivalent (\cref{exercise:integer-multiplication}).
\end{rem}

\subsection{Exercises}

\begin{exercise}\label{exercise:affine-image}
Let $\phi: \rr^n \to \rr^m$ be an affine map and let $S \subseteq \rr^n$. Show that
\begin{enumerate}
\item \label{affine-image-affine} If $A$ is an affine subspace of $\rr^n$, then $\phi(A)$ is an affine subspace of $\rr^m$, and $\dim(\phi(A)) \leq \dim(A)$.
\item \label{affine-image-hulls} $\phi(\aff(S)) = \aff(\phi(S))$ and $\phi(\conv(S)) = \conv(\phi(S))$.
\item \label{affine-image-cone} Assume $\phi$ is in addition a {\em linear map}, i.e.\ $\phi(0) = 0$. Then $\phi(\cone(S)) = \cone(\phi(S))$.
\end{enumerate}
Now assume $\phi$ is in addition injective. Then show that
\begin{enumerate} [resume]
\item The inequality in assertion \eqref{affine-image-affine} holds with equality.
\item The converses of assertions \eqref{affine-image-affine}, \eqref{affine-image-hulls} and \eqref{affine-image-cone} are also true, i.e.\
\begin{enumerate}
\item if $A \subseteq \rr^n$ is such that $\phi(A)$ is an affine subspace of $\rr^m$, then $A$ is also an affine subspace of $\rr^n$,
\item $\phi^{-1}(\aff(S)) = \aff(S)$ and $\phi^{-1}(\conv(\phi(S)) = \conv(S)$,
\item if $\phi$ is a linear map, then $\phi^{-1}(\cone(\phi(S)) = \cone(S)$.
\end{enumerate}
\end{enumerate}
\end{exercise}

\begin{exercise}\label{exercise:conv-rep}
Show that the line segment joining $\alpha, \beta \in \rr^n$ is precisely the set of their convex combinations. Use it to prove identity \eqref{conv-rep}. If $\nu \in \rnstar$ and $S \subset \rr^n$ are such that $\min_S (\nu)$ exists, then use identity \eqref{conv-rep} to show that $\min_{\conv(S)}(\nu) = \min_S(\nu)$.
\end{exercise}

\begin{exercise}\label{exercise:pol-is-convex}
Check that a convex polyhedron as defined in \cref{basic-convex-section} is actually convex.
\end{exercise}

\begin{exercise}\label{exercise:zerone-dim-pol}
Let $S \subseteq \rr^n$.
\begin{enumerate}
\item Show that $\aff(S)$ is zero dimensional if and only if $S$ is a ``singleton'' (i.e.\ it consists of a single point). Deduce that zero dimensional polyhedra are precisely singletons.
\item Assume $S$ is convex and closed in $\rr^n$, and $\aff(S)$ is one dimensional. Show that there is an affine isomorphism $\phi: \rr \to \aff(S)$ such that $\phi^{-1}(S)$ is a (possibly unbounded) closed interval of $\rr$.
\end{enumerate}
\end{exercise}

\begin{exercise}\label{exercise:2-cone}
Show that every finitely generated two dimensional cone $S$ in $\rr^2$ can be generated as a cone by two nonzero elements of $\rr^2$ (these correspond to ``edges'' of $S$). If $S$ is in addition strongly convex, show that the generators of $S$ are linearly independent.
\end{exercise}

\begin{exercise} \label{exercise:cone-inequalities}
Let $\scrP \subseteq \rr^n$.
\begin{enumerate}
\item Assume $\scrP$ is ``dilation invariant,'' i.e.\ for every $\alpha \in \scrP$ and every $r \geq 0$, $r\alpha \in \scrP$. Show that for every $\nu \in \rnstar$, $\inf\{ \langle \nu, \alpha \rangle: \alpha \in \scrP\}$ is either $0$ or $-\infty$.
\item Deduce that $\scrP$ is a convex polyhedral cone if and only if it can be defined by finitely many inequalities with zero constant term, i.e.\ inequalities of the form $a_1x_1 + \cdots + a_nx_n \geq 0$.
\end{enumerate}
\end{exercise}

\begin{exercise} \label{exercise:linear-sum-invariance}
If $S, T \subseteq \rr^m$ and $\phi: \rr^m \to \rr^n$ is a linear map, show that $\phi(S+T) = \phi(S) + \phi(T)$.
\end{exercise}

\begin{exercise}\label{exercise:minkonvex}
Let $\scrP, \scrQ$ be subsets of $\rr^n$.
\begin{enumerate}
\item If $\scrP = \conv(S)$ and $\scrQ = \conv(T)$, then show that $\scrP+ \scrQ= \conv(S + T)$. [Hint: if $\sum_i \delta_i = \sum_j \epsilon_j = 1$, then $\sum_i \delta_i a_i + \sum_j \epsilon_j b_j = \sum_{i,j}\delta_i\epsilon_j (a_i + b_j)$.]
\item Deduce that if $\scrP$ and $\scrQ$ are convex (respectively, convex hulls of finitely many points, finitely generated cones) then so is $\scrP + \scrQ$.
\end{enumerate}
\end{exercise}

\begin{exercise} \label{exercise:maxlinconvex}
Let $\scrC$ be a convex subset of $\rr^n$ containing the origin. Show that $\scrC$ contains a unique maximal linear subspace of $\rr^n$, i.e.\ there is a linear subspace $L$ of $\rr^n$ such that $L \subseteq \scrC$, and if $L'$ is any linear subspace of $\rr^n$ contained in $\scrC$, then $L' \subseteq L$. [Hint: if $L_1, L_2$ are two linear subspaces of $\rr^n$ contained in $\scrC$, then $L_1 + L_2 \subseteq \scrC$.]
\end{exercise}

\begin{exercise} \label{exercise:integer-multiplication}
Let $\scrA \subset \rr^n$ and $d$ be a positive integer.
\begin{enumerate}
\item If $\scrA$ is convex, then show that $\{d\alpha: \alpha \in \scrA\} = \{\alpha_1 + \cdots + \alpha_d: \alpha_j \in \scrA$ for each $j\}$.
\item Give an example to show that the above equality may not hold if $\scrA$ is not convex.
\end{enumerate}
\end{exercise}

\begin{exercise} \label{exercise:eq-dim-embedding}
Let $S \subset \rr^m$ and $n: = \dim(\aff(S))$.
\begin{enumerate}
\item \label{eq-dim-embedding:affine} Show that there is $T \subseteq \rr^n$ and an injective affine map $\phi:\rr^n \to \rr^m$ such that $\aff(T)$ is $\rr^n$, and $\phi(T) = S$.  [Hint: $\aff(S)$ is a translation of an $n$-dimensional linear subspace of $\rr^m$.]
\item If $S$ contains the origin (this is the case when e.g.\ $S$ is a cone) then show that it is possible to ensure in assertion \eqref{eq-dim-embedding:affine} that $\phi$ is in addition a {\em linear} map.
\end{enumerate}
\end{exercise}

\begin{exercise} \label{exercise:affine-invariance}
Show that each of the following properties is invariant under injective affine maps, i.e.\ if $\phi: \rr^{n_1} \to \rr^{n_2}$ is an injective affine map then each of the following properties holds with $S = S_1$ and $n = n_1$ if and only if it holds with $S = \phi(S_1)$ and $n = n_2$.
\begin{enumerate}
\item $S$ is an affine subspace of $\rr^n$,
\item $S$ is a convex subset of $\rr^n$,
\item $S$ is the convex hull of finitely many points in $\rr^n$,
\item $S$ is convex polyhedron in $\rr^n$,
\item $S$ is convex polytope in $\rr^n$.
\end{enumerate}
\end{exercise}

\begin{exercise}\label{exercise:linear-invariance}
Show that each of the following properties is invariant under injective linear maps, i.e.\ if $\phi: \rr^{n_1} \to \rr^{n_2}$ is an injective linear map then each of the following properties holds with $S = S_1$ and $n = n_1$ if and only if it holds with $S = \phi(S_1)$ and $n = n_2$.
\begin{enumerate}
\item $S$ is a linear subspace of $\rr^n$,
\item $S$ is a convex cone in $\rr^n$,
\item $S$ is a strongly convex cone in $\rr^n$,
\item $S$ is a finitely generated convex cone in $\rr^n$,
\item $S$ is a convex polyhedral cone in $\rr^n$.
\end{enumerate}
\end{exercise}

\begin{exercise} \label{exercise:dual-extension}
Let $H$ be a linear subspace of $\rr^n$.
\begin{enumerate}
\item \label{dual-extension-0} Show that every linear map $\nu:H \to \rr$ can be extended to a map $\rr^n \to \rr$.

\item \label{dual-extension-positive} Assume $H \subseteq \eta^\perp$ for some $\eta \in \rnstar$. Let $S \subset \rr^n$ be a finite set such that $\langle \eta, \alpha \rangle >0$ for each $\alpha \in S$, and let $c$ be an arbitrary real number. Show that in assertion \eqref{dual-extension-0} it can be ensured that $\langle \nu', \alpha \rangle > c$ for each $\alpha \in S$. [Hint: after a linear change of coordinates on $\rr^n$ we may assume that $H$ is spanned by the first $k$-coordinates and $\eta$ is the projection onto the $(k+1)$-th coordinate.]
\end{enumerate}
\end{exercise}

\begin{exercise}  \label{exercise:injectively-positive-dual}
Given $S \subset \rr^n$, consider the following property:
\begin{align}
\parbox{0.6\textwidth}{
there is $\nu \in \rnstar$ such that $\langle \nu, \alpha \rangle > 0$ for each $\alpha \in S\setminus \{0\}$.
}
\label{positive-dual-property}
\end{align}
Let $\phi: \rr^{n_1} \to \rr^{n_2}$ be an injective linear map and $S_1 \subseteq \rr^{n_1}$. Show that property \eqref{positive-dual-property} holds with $S = S_1$ and $n = n_1$ if and only if it holds with $S = \phi(S_1)$ and $n = n_2$. [Hint: use assertion \eqref{dual-extension-0} of \cref{exercise:dual-extension}.]
\end{exercise}


\begin{exercise} \label{exercise:0-line}
Let $S \subseteq \rr^n\setminus\{0\}$. Show that the origin is in the convex hull of $S$ if and only if $\cone(S)$ contains a line through the origin. [Hint: $0$ can be written as a convex combination of elements from $S$ if and only if there is $\alpha \in \cone(S)$, $\alpha \neq 0$, such that $-\alpha$ is also in $\cone(S)$.]
\end{exercise}


\begin{exercise} \label{exercise:strongly-convex-projection}
Let $\scrC$ be a convex cone in $\rr^n$ and $\pi:\rr^n \to \rr^m$ be a linear map such that $\ker(\pi) \cap \scrC = \{0\}$. Show that $\scrC$ is a strongly convex cone in $\rr^n$ if and only if $\pi(\scrC)$ is a strongly convex cone in $\rr^m$. [Hint: use \cref{exercise:0-line}.]
\end{exercise}

\section{Characterization of convex polyhedra} \label{poly-characterizection}
In this section we describe a characterization of convex polytopes and polyhedral cones, and use it to show that every polyhedron is the Minkowski sum of a polytope and a polyhedral cone. The proofs we give are elementary and geometric, but not the most ``efficient''; see e.g.\ \cite[Section 7.2]{schrijver} for a quicker proof (which is somewhat less intuitive in the beginning), and \cite[Section 1.3]{ziegtopes} for a more algorithmic proof.

\begin{thm}[Farkas (1898, 1902), Minkowski (1896), Weyl (1935)]  \label{poly-basic-thm}
Let $\scrP$ be a convex subset of $\rr^n$.
\begin{enumerate}
\item \label{characterize-polytope} $\scrP$ is a polytope if and only if it is the convex hull of finitely many points.
\item \label{characterize-cone} $\scrP$ is a polyhedral cone if and only if it is a finitely generated cone.
\item \label{characterize-strongly-convex} Assume $\scrP$ is a polyhedral cone. Then it is strongly convex if and only if there is $\nu \in \rnstar$ such that $\langle \nu, \alpha \rangle > 0$ for each $\alpha \in \scrP \setminus \{0\}$.
\end{enumerate}
\end{thm}


\begin{proof}
We are going to prove assertions \eqref{characterize-cone}, \eqref{characterize-strongly-convex} and $(\Leftarrow)$ implication of assertion \eqref{characterize-polytope}. The $(\im)$ direction of assertion \eqref{characterize-polytope} follows from these assertions; it is left as an exercise (\cref{exercise:pol-to-hull}). We start with the proof of assertion \eqref{characterize-strongly-convex}. For the $(\Leftarrow)$ direction note that if $\scrP$ contains both $\alpha$ and $-\alpha$ for some nonzero $\alpha \in \rr^n$, then for all $\nu \in \rnstar$, either $\langle \nu, \alpha \rangle \leq 0$ or $\langle \nu, -\alpha \rangle < 0$. For the $(\im)$ direction we proceed by induction on $\dim(\scrP)$. Due to \cref{exercise:eq-dim-embedding,exercise:linear-invariance,exercise:injectively-positive-dual} we may assume \woutlog\ that $\dim(\scrP) = n$. Since the only strongly positive cones in $\rr$ are $\rr_{\geq 0}$ and $\rr_{\leq 0}$, it holds for $n = 1$. In the general case $\scrP$ is defined by finitely many inequalities of the form $\langle \nu, x \rangle \geq 0$ for nonzero $\nu \in \rnstar$ (\cref{exercise:cone-inequalities}). Take one such $\nu$. Since $\nu^\perp \cap \scrP$ is a smaller dimensional strongly convex polyhedral cone, it follows by induction that there is $\eta \in \rnstar$ which is positive on $(\nu^\perp \cap \scrP) \setminus \{0\}$.

\begin{proclaim} \label{claim:positive-dual-by-induction}
For sufficiently small $\epsilon > 0$, $\nu + \epsilon \eta$ is positive on $\scrP\setminus \{0\}$.
\end{proclaim}

\begin{proof}
If the claim is false, then we can find a sequence of positive numbers $\epsilon_k \to 0$ and $\alpha_k \in S^{n-1} \cap \scrP$, where $S^{n-1}$ is the unit sphere centered at the origin in $\rr^n$, such that $\langle \nu + \epsilon_k \eta, \alpha_k \rangle < 0$. Since $\scrP$ is closed and $S^{n-1}$ is compact, we may assume that the $\alpha_k$ converge to $\alpha \in \scrP \cap S^{n-1}$. Then $\langle \nu, \alpha \rangle = \lim_{k \to \infty} \langle \nu + \epsilon_k \eta, \alpha_k \rangle \leq  0$, so that $\langle \nu, \alpha \rangle = 0$. It follows that $\alpha \in  S^{n-1} \cap \nu^\perp \cap \scrP$, so that $\langle \eta, \alpha \rangle > 0$. By continuity $\eta$ is positive on an open neighborhood $U$ of $\alpha$, which means that $\langle \nu + \epsilon_k \eta, \alpha_k \rangle > 0$ for sufficiently large $k$. This contradiction proves the claim.
\end{proof}

\Cref{claim:positive-dual-by-induction} finishes the proof of assertion \eqref{characterize-strongly-convex}. We now prove the $(\im)$ direction of assertion \eqref{characterize-cone} by induction on $\dim(\scrP)$. Due to  \cref{strongly-convex-prop,exercise:eq-dim-embedding,exercise:linear-invariance} we may assume \woutlog\ that $\scrP$ is an $n$-dimensional strongly convex polyhedral cone in $\rr^n$. As in the proof of assertion \eqref{characterize-strongly-convex}, the $n = 1$ case follows directly from the observation that the only strongly convex polyhedral cones in $\rr$ are $\rr_{\geq 0}$ and $\rr_{\leq 0}$. Now consider the case $n \geq 2$. Assertion \eqref{characterize-strongly-convex} implies that there is $\nu \in \rnstar$ such that $\langle \nu, \alpha \rangle > 0$ for each $\alpha \in \scrP \setminus \{0\}$.

\begin{proclaim} \label{bounded-nu=1}
$\scrP' := \{\alpha \in \scrP : \langle \nu, \alpha \rangle = 1\}$ is bounded.
\end{proclaim}

\begin{proof}
Indeed, otherwise take a sequence $\alpha_k \in \scrP'$ such that $\norm{\alpha_k}$ is unbounded, where $\norm{\cdot}$ is the Euclidean norm on $\rr^n$. Note that each $\alpha_k/\norm{\alpha_k}$ is in the intersection of $\scrP$ and the $(n-1)$-dimensional unit sphere $S^{n-1}$. Since $S^{n-1}$ is compact and $\scrP$ is closed, we may assume \woutlog\ that $\alpha_k/\norm{\alpha_k}$ converge to $\alpha \in \scrP \cap S^{n-1}$. But then $\langle \nu, \alpha \rangle = \lim_{k \to \infty} \langle \nu, \alpha_k/\norm{\alpha_k} \rangle = 0$, which contradicts the choice of $\nu$.
\end{proof}

\Cref{exercise:cone-inequalities} implies that $\scrP$ is defined by inequalities of the form $\langle \nu_j, x \rangle \geq 0$, $j = 1, \ldots, N$, for some $\nu_1, \ldots, \nu_N \in \rnstar$. For each $j$, let $\scrP_j := \scrP \cap \nu_j^\perp$. By the inductive hypothesis, each $\scrP_j$ is generated by finite sets $S_j \subset \rr^n$. We claim that $\scrP$ is generated by $\bigcup_j S_j$. Indeed, take $\alpha \in \scrP \setminus \{0\}$. There is $r > 0$ such that $r\alpha \in \scrP'$. Take any straight line through $r\alpha$ on the hyperplane $\{\alpha: \langle \nu, \alpha \rangle = 1\}$. \Cref{bounded-nu=1} implies that either $r\alpha \in \scrP_j$ for some $j$, or each end of the line intersects one of the $\scrP_j$. In any event, the inductive hypothesis implies that $\alpha$ is a nonnegative linear combination of elements from $\bigcup_j S_j$, as required to complete the proof of $(\im)$ implication of assertion \eqref{characterize-cone}. Now we prove $(\Leftarrow)$ implications of assertions \eqref{characterize-polytope} and \eqref{characterize-cone} by induction on $\dim(\aff(\scrP))$. We will first show that it suffices to prove only the implication from assertion \eqref{characterize-cone}.

\begin{figure}[h]
\begin{center}

\xdefinecolor{c1}{HTML}{D4ACEE}
\xdefinecolor{c2}{HTML}{EEc6AC}
\xdefinecolor{c3}{HTML}{EDEEAC}
\def\tinyshifth{0.09}
\def\tinyshiftv{0.09}

\def\xmin{-2}
\def\xmax{8}
\def\ymin{-3}

\tikzstyle{dot} = [red, circle, minimum size=4pt, inner sep = 0pt, fill]

\begin{tikzpicture}[scale=\scalefactor]

\coordinate (A) at (0,0);
\coordinate (B) at (4,2);
\coordinate (C) at (3,3);

\foreach \x in {0,\tinyshifth,...,\xmax} {
	\draw[color=c1, ultra thin] (\x, {0.5*\x}) -- (\x,\x);
}
\pgfmathsetmacro\xlim{4-\xmin}
\foreach \x in {0,\tinyshifth,...,\xlim} {
	\draw[color=c2, ultra thin] ({4-\x}, {2-0.5*\x}) -- ({4-\x},{2+\x});
}
\pgfmathsetmacro\ylim{3-\ymin}
\foreach \y in {0,\tinyshiftv,...,\ylim} {
	\draw[color=c3, ultra thin] ({3-\y}, {3-\y}) -- ({3+\y}, {3-\y});
}

\node[dot] at (A) {};
\node[dot] at (B) {};
\node[dot] at (C) {};

\draw[thick, \colorAB] (A) -- (B);
\draw[thick, \colorBC] (B) -- (C);
\draw[thick, \colorCA] (C) -- (A);
\node[anchor = north] at (A) {$\alpha_1$};
\node[anchor = west] at (B) {$\alpha_2$};
\node[anchor = south] at (C) {$\alpha_3$};


\def\xmin{-5.5}
\def\xshift{9}
\def\yshift{2.5}
\pgfmathsetmacro\shiftone{\xshift + \xmax -\xmin};
\def\xmax{5.5}
\def\ymin{-5.5}
\def\ymax{5.5}

\begin{scope}[shift={(\shiftone,\yshift)}]
\coordinate (O) at (0,0);
\coordinate (BA) at (\xmin, \xmin/2);
\coordinate (AB) at (\xmax, \xmax/2);
\coordinate (CA) at (\ymin, \ymin);
\coordinate (AC) at (\ymax, \ymax);
\coordinate (BC) at (-\ymax, \ymax);
\coordinate (CB) at (-\ymin, \ymin);

\foreach \x in {0,\tinyshifth,...,\xmax} {
	\draw[color=c1, ultra thin] (\x, {0.5*\x}) -- (\x,\x);
	\draw[color=c2, ultra thin] (-\x, -{0.5*\x}) -- (-\x,\x);
}

\foreach \y in {0,-\tinyshiftv,...,\ymin} {
	\draw[color=c3, ultra thin] (\y, \y) -- (-\y,\y);
}


\draw [thick, \colorAB] (AB) -- (BA);
\draw [thick, \colorBC] (CB) -- (BC);
\draw [thick, \colorCA] (AC) -- (CA);

\node at (-3.5,0.5) {\picfontsize $\scrC_2$};
\node at (4.5,3.5) {\picfontsize $\scrC_1$};
===\node at (0.5,-3.5) {\picfontsize $\scrC_3$};
\end{scope}
\end{tikzpicture}
\caption{Convex hull of finitely many points is a finite intersection of translations of cones}  \label{fig:polytope=polyhedral-intersection}
\end{center}
\end{figure}

\begin{proclaim} \label{fin2cone-to-hull2pol}
Let $k \geq 1$. Assume the $(\Leftarrow)$ implication of assertion \eqref{characterize-cone} holds whenever $\dim(\aff(\scrP)) \leq k$. Then the $(\Leftarrow)$ implication of assertion \eqref{characterize-polytope} holds whenever $\dim(\aff(\scrP)) \leq k$.
\end{proclaim}

\begin{proof}
Let $\scrP$ be the convex hull of a finite set $S$, and $m := \dim(\aff(\scrP)) \leq k$. \Woutlog\ we may assume that $\scrP\subset \rr^m$ (\cref{exercise:eq-dim-embedding,exercise:affine-invariance}). Let $\alpha_1, \ldots, \alpha_s$ be the elements of $S$. For each $j$, let $\scrC_j$ be the cone generated by $S - \alpha_j = \{\alpha_i - \alpha_j: 1 \leq i \leq s\}$. The hypothesis of the claim implies that each $\scrC_j$ is a polyhedral cone, which implies in turn that $\scrC_j + \alpha_j$ is a polyhedron (\cref{exercise:affine-invariance}). It follows that $\scrP' := \bigcap_j (\scrC_j + \alpha_j)$ is also a polyhedron. Therefore it suffices to show that $\scrP' = \scrP$ (\cref{fig:polytope=polyhedral-intersection}). By \eqref{conv-rep} every element $\alpha \in \scrP$ can be expressed as $\sum_i \epsilon_i\alpha_i$ where $\epsilon_i$ are nonnegative real numbers with $\sum_i \epsilon_i = 1$. But then $\alpha = \alpha_j + \sum_i \epsilon_i(\alpha_i - \alpha_j) \in \scrC_j + \alpha_j$ for each $j$. Therefore $\scrP\subseteq \scrP'$. Now take $\alpha \not\in \scrP$. To complete the proof of the claim it suffices to show that $\alpha \not \in \scrP'$. Since $\alpha \not\in \scrP$, it follows that $\cone(\scrP - \alpha)$ is strongly convex (\cref{exercise:0-line}). By the hypothesis of the claim $\cone(\scrP - \alpha)$ is polyhedral, so that due to assertion \eqref{characterize-strongly-convex} there is $\nu \in \rnstar$ which is positive on $ \scrP - \alpha$. Pick $j$ such that $\langle \nu, \alpha_j\rangle = \min \{\langle \nu, \alpha_i\rangle: 1 \leq i \leq s\}$. Then $\nu$ is nonnegative on $\scrC_j$, whereas $\langle \nu, \alpha - \alpha_j \rangle $ is negative. Therefore $\alpha \not\in \scrC_j + \alpha_j$, and consequently $\alpha \not\in \scrP'$, as required.
\end{proof}

Now we start the proof of $(\Leftarrow)$ direction of assertion \eqref{characterize-cone} by induction on $\dim(\aff(\scrP))$. Due to \cref{strongly-convex-prop,exercise:eq-dim-embedding,exercise:linear-invariance} we may assume \woutlog\ that $\scrP$ is a strongly convex cone generated by a finite set $S \subset \rr^n$ and $n = \dim(\aff(\scrP))$. For $n = 1$ the possibilities for $\scrP$ are $\rr_{\geq 0}$ and $\rr_{\leq 0}$, both of which are polyhedral. For general $n$ we proceed by induction on $|S|$. The case $|S| = 1$ is also covered by the case $n = 1$. So assume $|S| \geq 2$.

\begin{proclaim}
There is $\nu \in \rnstar$ which is positive on $S \setminus \{0\}$.
\end{proclaim}

\begin{proof}
 Pick $\alpha \in S \setminus \{0\}$. By the inductive hypothesis $\scrP_1 := \cone(S \setminus \{\alpha\})$ is polyhedral, so that $\scrP_1$ is defined by finitely many inequalities of the form $\langle \nu_j, x \rangle \geq 0$, $j = 1, \ldots, N$ (\cref{exercise:cone-inequalities}). We claim that $\langle \nu_j, \alpha \rangle > 0$ for some $j$. Indeed, otherwise $\langle \nu_j, -\alpha \rangle \geq 0$ for each $j$, so that $-\alpha \in \scrP_1$. But then $\scrP$ contains the line through the origin and $\alpha$, contradicting the strong convexity of $\scrP$. So we can pick $j$ such that $\langle \nu_j, \alpha \rangle > 0$. Now let $\scrP_2 := \scrP \cap \nu_j^\perp$. Then $\dim(\aff(\scrP_2)) < n$, so that $\scrP_2$ is polyhedral by the inductive hypothesis. Assertion \eqref{characterize-strongly-convex} then implies that there is $\nu \in \rnstar$ which is positive on $\scrP_2 \setminus \{0\}$. We claim that if $\epsilon$ is a sufficiently small positive number, then $\nu_j + \epsilon \nu$ is positive on $\scrP \setminus \{0\}$. Indeed, simply take $\epsilon$ such that $\epsilon | \langle  \nu, \beta \rangle | < \langle \nu_j, \beta \rangle$ for all (the finitely many) elements of $S \setminus \scrP_2$.
\end{proof}

Let $\scrP' := \{\alpha \in \scrP: \langle \nu, \alpha \rangle =1\}$. If $\alpha_1, \ldots, \alpha_s$ are elements of $S$, then $\alpha'_j := \alpha_j/\langle \nu, \alpha_j\rangle \in \scrP'$ for each $j$, and it is straightforward to check that $\scrP'$ is the convex hull of $\alpha'_1, \ldots, \alpha'_s$. Since $\dim(\aff(\scrP')) < n$, the inductive hypothesis and \cref{fin2cone-to-hull2pol} imply that $\scrP'$ is a polytope. It is then straightforward to check that $\scrP$ is a polyhedral cone (\cref{exercise:pol-to-cone}), as required.
\end{proof}

\begin{cor}  \label{cor:polyhedral-sum}
The Minkowski sum of two convex polytopes (respectively, polyhedral cones) is a convex polytope (respectively, polyhedral cone).
\end{cor}

\begin{proof}
Follows immediately from \cref{poly-basic-thm,exercise:minkonvex}.
\end{proof}

The next corollary of \cref{poly-basic-thm} shows that each convex polyhedron $\scrQ$ has a representation of the form $\scrQ= \scrP+ \scrC$, where $\scrC$ is a polyhedral cone and $\scrP$ is a polytope. In general the decomposition is {\em not} unique, see \cref{fig:non-unique-sum}. However, $\scrC$ is uniquely determined from $\scrQ$ (\cref{exercise:minkomposition-uniqueness}). It is also possible to find $\scrP$ which is ``minimal'' (modulo translations in the case that $\scrC$ is not strongly convex), but we will not get into that.

\begin{center}
\begin{figure}[h]

\def\scalefactor{.3}

\def\xmin{-2}
\def\xmax{8}
\def\ymin{-5}
\def\ymax{8}

\def\tx{0}
\def\ty{\ymin}

\def\xshiftone{6}
\def\xshifttwo{5}

\def\yminsmall{-2.5}
\def\ymaxsmall{4.5}

\def\ygap{5}
\def\yshiftup{6}
\def\yshiftdown{-5}
\tikzstyle{dot} = [\colordot, circle, minimum size=4pt, inner sep = 0pt, fill]

\begin{tikzpicture}[scale=\scalefactor]
\pgfmathsetmacro\shiftone{\xshiftone + \xmax -\xmin};
\pgfmathsetmacro\shifttwo{\xshifttwo + \xmax -\xmin};

\coordinate (A) at (-2, -1);
\coordinate (C') at (8,-5);
\coordinate (D') at (8,7.5);
\coordinate (E) at (1,4);

\fill[\colorzero, opacity=\opazero ] (A) -- (C') -- (D') -- (E) -- cycle;

\draw[thick, \colorone ] (A) -- (C');
\draw[thick, \colorone ] (A) -- (E) -- (D');

\draw (\xmax, 0) -- (\xmin, 0);
\draw (0, \ymax) -- (0, \ymin);

\node[dot] (O) at (0,0){};
\node[anchor = south west] at (O) {\picfontsize $(0,0)$};
\node[anchor = west] at (\xmax, 0) {\picfontsize $x$};
\node[anchor = south] at (0, \ymax) {\picfontsize $y$};

\node [below right, text width= 1.8cm, align=left] at (\tx,\ty) {
	\picfontsize
	$\scrQ$
};
	
\def\ymax{\ymaxsmall}
\def\ymin{\yminsmall}

\begin{scope}[shift = {(\shiftone,\yshiftup)}]

\coordinate (A) at (-2, -1);
\coordinate (B) at (3,-3);
\coordinate (C) at (7,1);
\coordinate (D) at (5,6);
\coordinate (E) at (1,4);

\draw[fill=\colorzero, opacity=\opazero ] (A) --  (B) -- (C) -- (D) -- (E) -- cycle;

\draw[thick, \colorone ] (A) --  (B) -- (C) -- (D) -- (E) -- cycle;

\draw (\xmax, 0) -- (\xmin, 0);
\draw (0, \ymax) -- (0, \ymin);

\node[dot] (O) at (0,0){};
\node[anchor = south west] at (O) {\picfontsize $(0,0)$};
\node[anchor = west] at (\xmax, 0) {\picfontsize $x$};
\node[anchor = south] at (0, \ymax) {\picfontsize $y$};

\node [below right, text width= 0.75cm, align=left] at (\tx,\ty) {
	\picfontsize
	$\scrP_1$
};	

\end{scope}

\begin{scope}[shift = {(\shiftone + \shifttwo,\yshiftup)}]
	\coordinate (C') at (8,-3.2);
	\coordinate (D') at (8, 4);
	\coordinate (O) at (0,0);
	
	\fill[\colorzero, opacity=\opazero ] (O) -- (C') -- (D') -- cycle;
	
	\draw[thick, \colorone ] (O) -- (C');
	\draw[thick, \colorone ] (O) -- (D');
	
	\draw (\xmax, 0) -- (\xmin, 0);
	\draw (0, \ymax) -- (0, \ymin);
	
	\node[dot] (extra) at (O) {};
	\node[anchor = north] at (O) {\picfontsize $(0,0)$};
	\node[anchor = west] at (\xmax, 0) {\picfontsize $x$};
	\node[anchor = south] at (0, \ymax) {\picfontsize $y$};
	
	\node [below right, text width= 2.7cm, align=left] at (\tx,\ty) {
		\picfontsize
		$\scrC$
	};	
\end{scope}

\begin{scope}[shift = {(\shiftone,\yshiftdown)}]

\coordinate (A) at (-2, -1);
\coordinate (E) at (1,4);

\draw[thick, \colorone ] (A) -- (E);

\draw (\xmax, 0) -- (\xmin, 0);
\draw (0, \ymax) -- (0, \ymin);

\node[dot] (O) at (0,0){};
\node[anchor = south west] at (O) {\picfontsize $(0,0)$};
\node[anchor = west] at (\xmax, 0) {\picfontsize $x$};
\node[anchor = south] at (0, \ymax) {\picfontsize $y$};

\node [below right, text width= 0.75cm, align=left] at (\tx,\ty) {
	\picfontsize
	$\scrP_2$
};	

\end{scope}

\begin{scope}[shift = {(\shiftone + \shifttwo,\yshiftdown)}]
	\coordinate (C') at (8,-3.2);
	\coordinate (D') at (8, 4);
	\coordinate (O) at (0,0);
	
	\fill[\colorzero, opacity=\opazero ] (O) -- (C') -- (D') -- cycle;
	
	\draw[thick, \colorone ] (O) -- (C');
	\draw[thick, \colorone ] (O) -- (D');
	
	\draw (\xmax, 0) -- (\xmin, 0);
	\draw (0, \ymax) -- (0, \ymin);
	
	\node[dot] (extra) at (O) {};
	\node[anchor = north] at (O) {\picfontsize $(0,0)$};
	\node[anchor = west] at (\xmax, 0) {\picfontsize $x$};
	\node[anchor = south] at (0, \ymax) {\picfontsize $y$};
	
	\node [below right, text width= 2.7cm, align=left] at (\tx,\ty) {
		\picfontsize
		$\scrC$
	};	
\end{scope}

\end{tikzpicture}

\caption{
$\scrQ = \scrP + \scrC$ for $\scrP = \scrP_1$ and $\scrP= \scrP_2$. The ``minimal'' possible choice for $\scrP$ is the line segment $\scrP_2$.
} \label{fig:non-unique-sum}

\end{figure}
\end{center}

\begin{cor}[Motzkin (1936)] \label{poly-basic-cor}
A subset of $\rr^n$ is a convex polyhedron if and only if it is the Minkowski sum of a convex polyhedral cone and a convex polytope.
\end{cor}

\begin{proof}
At first we prove the $(\im)$ implication. Let $\scrQ$ be a polyhedron in $\rr^n$ defined by inequalities $a_{q,0} + a_{q,1}x_1 + \cdots + a_{q,n}x_n \geq 0$, $q = 1, \ldots, N$. We will show that it is the sum of a polyhedral cone and a polytope.
Consider the polyhedral cone $\scrQ'$ in $\rr^{n+1}$ defined by the ``homogenizations'' $a_{q,0}x_0 +  a_{q,1}x_1 + \cdots + a_{q,n}x_n \geq 0$, $q = 1, \ldots, N$, and $x_0 \geq 0$. Note that $\scrQ = \scrQ' \cap (\{1\} \times \rr^n)$. \Cref{poly-basic-thm} implies that $\scrQ'$ is a finitely generated cone. We may choose a set of generators of $\scrQ'$ of the form $(1, \alpha_1), \ldots, (1, \alpha_k), (0, \beta_1), \ldots, (0, \beta_l)$, where $\alpha_1, \ldots, \alpha_k, \beta_1, \ldots, \beta_l \in \rr^n$. Let $\scrP$ be the convex hull of $\alpha_1, \ldots, \alpha_k$, and $\scrC$ be the cone in $\rr^n$ generated by $\beta_1, \ldots, \beta_l$. \Cref{poly-basic-thm} implies that $\scrP$ is a polytope and $\scrC$ is a polyhedral cone.

\begin{proclaim} \label{poly-basic-cor:cla-im}
$\scrQ = \scrP + \scrC$.
\end{proclaim}

\begin{proof}
Let $\gamma \in \scrQ$. Then $(1, \gamma) \in \scrQ'$, and therefore $(1, \gamma) = \sum_i a_i (1, \alpha_i) + \sum_j b_j (0, \beta_j)$ where the $a_i, b_j$ are nonnegative numbers. Then it follows that $\sum_i a_i = 1$ and $\gamma = \alpha + \beta$, where $\alpha = \sum_i a_i \alpha_i$ and $\beta = \sum_j b_j \beta_j$. It is clear that $\beta \in \scrC$. Since $\sum_i a_i = 1$, it follows that $\alpha$ is a convex combination of the $\alpha_i$, so that $\alpha \in \scrP$. Therefore $\scrQ \subset \scrP+ \scrC$. Now we check the opposite inclusion. Let $\alpha \in \scrP$ and $\beta \in \scrC$. Then $\alpha = \sum_i a_i \alpha_i$ where the $a_i$ are nonnegative numbers such that $\sum_i a_i = 1$. Since $(1, \alpha_i) \in \scrQ'$ for each $i$, it follows that $(1, \alpha) = \sum_i a_i (1,\alpha_i) \in \scrQ'$. On the other hand, it follows from the construction of $\scrC$ that $\{0\} \times \scrC = \scrQ' \cap (\{0\} \times \rr^n)$, so that $(0,\beta) \in \scrQ'$. Since the sum of two elements of a convex cone is also in that cone, it follows that $(1, \alpha + \beta) = (1, \alpha) + (0, \beta) \in \scrQ'$. Therefore, $\alpha + \beta \in \scrQ$, as required.
\end{proof}
\Cref{poly-basic-cor:cla-im} finishes the proof of $(\im)$ implication of the corollary. Now we prove the $(\Leftarrow)$ implication. Let $\scrP$ be a convex polytope and $\scrC$ be a convex polyhedral cone in $\rr^n$. We will show that $\scrP+ \scrC$ is a convex polyhedron. By \cref{exercise:affine-invariance} we may assume \woutlog\ that $\scrP$ contains the origin. Identify $\rr^{n+1}$ with $\rr \times \rr^n$. Let $\scrP'$ be the cone in $\rr^{n+1}$ generated by $\{1\} \times \scrP$, and $\scrC'$ be the closure of the cone in $\rr^{n+1}$ generated by $\{1\} \times \scrC$. Either \cref{poly-basic-thm} or \cref{exercise:pol-to-cone} implies that $\scrP'$ is a polyhedral cone, and \cref{exercise:cone-to-cone} implies that $\scrC'$ is a polyhedral cone. \Cref{cor:polyhedral-sum} then implies that $\scrP' + \scrC'$ is a polyhedral cone, so that $\scrQ := (\scrP' + \scrC') \cap (\{1\} \times \rr^n)$ is a polyhedron.

\begin{proclaim} \label{claim:cone-hedral-crossectoin}
$\scrQ = \{1\} \times (\scrP + \scrC)$.
\end{proclaim}

\begin{proof}
Let the inequalities defining $\scrP$ be $a_{i,0} + a_{i,1}x_1 + \cdots + a_{i,n}x_n \geq 0$, $i = 1, \ldots, k$, and the inequalities defining $\scrC$ be $b_{j,1}x_1 + \cdots + b_{j,n}x_n \geq 0$, $j = 1, \ldots, l$. \Cref{exercise:pol-to-cone} shows that $\scrP'$ is defined by the inequalities $a_{i,0}x_0 + a_{i,1}x_1 + \cdots + a_{i,n}x_n \geq 0$, $i = 1, \ldots, k$, and $x_0 \geq 0$. On the other hand \cref{exercise:cone-to-cone} implies that $\scrC'$ is defined by the same inequalities as $\scrC$ together with the inequality $x_0 \geq 0$. If $\alpha \in \scrP$ and $\beta \in \scrC$, then it follows that $(1, \alpha) \in \scrP'$ and $(0, \beta) \in \scrC'$, so that $(1, \alpha + \beta) \in \scrP' + \scrC'$, which proves that $\scrQ \supset \{1\} \times (\scrP + \scrC)$. For the opposite inclusion, pick $(1, \gamma) \in \scrQ$. We will show that $\gamma \in \scrP+ \scrC$. Write $(1, \gamma) = (a, \alpha) + (b, \beta)$, where $(a, \alpha) \in \scrP'$ and $(b, \beta) \in \scrC'$, and $a, b$ are nonnegative numbers such that $a + b = 1$. An examination of the inequalities defining $C'$ shows that $\beta \in \scrC$. If $a = 0$, then assertion \eqref{pol-to-cone-0} of \cref{exercise:pol-to-cone} implies that $\alpha = 0$. Since the origin is in $\scrP$, it follows that $\gamma = 0 + \beta \in \scrP + \scrC$. On the other hand, if $a \neq 0$, then $(1/a)\alpha \in \scrP$. Since $0 < a \leq 1$, and since $0 \in \scrP$, it follows by convexity of $\scrP$ that $\alpha \in \scrP$. Therefore $\gamma = \alpha + \beta \in \scrP+ \scrQ$, as required.
\end{proof}

\Cref{claim:cone-hedral-crossectoin} finishes the proof of $(\Leftarrow)$ implication, and therefore the proof of \cref{poly-basic-cor}.
\end{proof}

\begin{cor} \label{polyhedron-minkomposition}
Let $\scrP$ be a polyhedron in $\rr^n$.
\begin{enumerate}
\item $\scrP$ has a decomposition of the form $\scrP = \scrP_0 + \scrC + L$, where $\scrP_0$ is a polytope, $\scrC$ is strongly convex polyhedral cone, and $L$ is a linear subspace of $\rr^n$ such that $\scrC \cap L = \{0\}$. Moreover, $L$ and $\scrC + L$ are uniquely determined by $\scrP$.
\item If $\scrP$ is defined by inequalities $a_{i,0} + a_{i,1}x_1 + \cdots + a_{i,n}x_n \geq 0$, $i = 1, \ldots, N$, then
\begin{enumerate}
\item $L = \{x:a_{i,1}x_1 + \cdots + a_{i,n}x_n = 0,\ i = 1, \ldots, N\}$.
\item $\scrC + L = \{x:a_{i,1}x_1 + \cdots + a_{i,n}x_n \geq 0,\ i = 1, \ldots, N\}$.
\item For each $r \gg 1$, the set $\scrP_r := \{x:a_{i,1}x_1 + \cdots + a_{i,n}x_n = 0,\ i = 1, \ldots, N,$ and $-r \leq x_j \leq r$, $j = 1, \ldots, n\}$ is a polytope, and we can take $\scrP_0 = \scrP_r$ for any such $r$.
\end{enumerate}
\end{enumerate}
\end{cor}

\begin{proof}
Combine \cref{strongly-convex-prop,poly-basic-cor}.
\end{proof}

\begin{cor} \label{cor:polyhedral-sum-general}
The Minkowski sum of finitely many convex polyhedra is a convex polyhedron.
\end{cor}

\begin{proof}
Follows directly from \cref{cor:polyhedral-sum} and the arguments of the proof of \cref{poly-basic-cor}.
\end{proof}

\begin{cor}
The image of a polyhedron under an affine map is also a polyhedron.
\end{cor}

\begin{proof}
By \cref{exercise:affine-invariance} it suffices to consider the case of linear maps. Let $\scrP$ be a polyhedron in $\rr^n$, and $\phi:\rr^n \to \rr^m$ be a linear map. Write $\scrP = \scrP_0 + \scrC + L$ as in \cref{polyhedron-minkomposition}. Then $\phi(\scrP) = \phi(\scrP_0) + \phi(\scrC) + \phi(L)$ (\cref{exercise:linear-sum-invariance}). \Cref{poly-basic-thm,exercise:affine-image} imply that $\phi(\scrP)$ is a polytope, and $\phi(\scrC)$ is a polyhedral cone. Since $\phi(L)$ is a linear subspace of $\rr^m$ (due to linearity of $\phi$), \cref{cor:polyhedral-sum-general} implies that $\phi(\scrP)$ is a polyhedron.
\end{proof}

\subsection{Exercises}

\begin{exercise} \label{exercise:pol-to-cone}
Let $\scrP$ be a nonempty polytope in $\rr^n$ defined by inequalities $a_{i,0} + a_{i,1}x_1 + \cdots + a_{i,n}x_n \geq 0$, $i = 1, \ldots, N$. Let $\scrC_0 := \{x \in \rr^n: a_{i,1}x_1 + \cdots + a_{i,n}x_n \geq 0,\ 1 \leq i \leq N\}$, and $\scrC$ be the cone generated by $\{1\} \times \scrP$ in $\rr^{n+1} = \rr \times \rr^n$.
\begin{enumerate}
\item \label{pol-to-cone-0} Show that $\scrC_0 = \{0\}$. [Hint: $\scrC_0$ is a cone, $\scrC_0 + \scrP\subset\scrP$, and $\scrP$ is bounded.]
\item Deduce that $\scrC$ is a polyhedral cone in $\rr^{n+1}$ defined by the inequalities $x_0 \geq 0$ and $a_{i,0}x_0 + a_{i,1}x_1 + \cdots + a_{i,n}x_n \geq 0$, $i = 1, \ldots, N$.
\item Show that the assumption ``$\scrP$ is nonempty'' is necessary in the preceding statements. In particular, show by an example that $\scrC_0$ may be a nontrivial cone if $\scrP = \emptyset$.
\end{enumerate}
\end{exercise}

\begin{exercise} \label{exercise:pol-to-hull}
Show that if the $(\im)$ implication of assertion \eqref{characterize-cone} of \cref{poly-basic-thm} holds for all $n$, then the $(\im)$ direction of assertion \eqref{characterize-polytope} of \cref{poly-basic-thm} also holds for all $n$. [Hint: use \cref{exercise:pol-to-cone}.]
\end{exercise}

\begin{exercise} \label{exercise:cone-to-cone}
Let $\scrC$ be a polyhedral cone in $\rr^n$. \Cref{exercise:cone-inequalities} implies that $\scrC$ is defined by finitely many inequalities of the form $a_1x_1 + \cdots + a_nx_n \geq 0$. Show that the closure in $\rr^{n+1}$ of the cone generated by $\{1\} \times \scrC$ is a polyhedral cone defined by the same inequalities as the ones defining $\scrC$ coupled with $x_0 \geq 0$.
\end{exercise}

\begin{exercise} \label{exercise:minkomposition-uniqueness}
Let $\scrP \subset \rr^n$.
\begin{enumerate}
\item Assume $\scrP = \scrP_1 + \scrC_1 = \scrP_2 + \scrC_2$, where $\scrP_j$ are bounded and $\scrC_j$ are closed cones. Prove that $\scrC_1 = \scrC_2$. Show by examples that it is possible to have $\scrP_1 \neq \scrP_2$.
\item Assume $\scrP = \scrC_1 + L_1 =  \scrC_2 + L_2$ where $\scrC_j$ are strongly convex polyhedral cones and $L_j$ are linear subspaces of $\rr^n$ such that $\scrC_j \cap L_j = \{0\}$. Prove that $L_1 = L_2$. Show by examples that it is possible to have $\scrC_1 \neq \scrC_2$.
\item Assume $\scrP = \scrP_1 + \scrC_1 + L_1 = \scrP_2 + \scrC_2 + L_2$ where $\scrP_j$ are polytopes, $\scrC_j$ are strongly convex polyhedral cones and $L_j$ are linear subspaces of $\rr^n$ such that $\scrC_j \cap L_j = \{0\}$. Prove that $L_1 = L_2$ and $\scrC_1 + L_1 = \scrC_2 + L_2$.
\end{enumerate}
[Hint: we may assume \woutlog\ that $\scrP$ contains origin. Then for the first assertion pick $\alpha \in \scrP$ such that $r\alpha \in \scrP$ for each $r \geq 0$. Express $\alpha$ as a sum of elements of $\frac{1}{r}\scrP_j$ and $\scrC_j$, and then take the limit as $r \to \infty$ to show that each $\scrC_j$ must be the largest cone contained in $\scrP$. For the second assertion prove that each $L_j$ is the largest linear subspace of $\rr^n$ contained in $\scrP$ as follows: pick $\alpha \in \rr^n$ such that both $\alpha$ and $-\alpha$ are in $\scrP$. Express both $\alpha$ and $-\alpha$ as elements of $\scrC_j + L_j$, add them up, and use the condition that $\scrC_j \cap L_j = \{0\}$ to show that $\alpha \in L_j$.]
\end{exercise}

\section{Basic properties of convex polyhedra} \label{poly-basic-section}
In this section we establish a few basic properties of convex polyhedra. Throughout this section $\scrP$ will denote a {\em nonempty} convex polyhedron in $\rr^n$. The first property we state follows directly from the definition of a polyhedron:

\begin{prop} \label{poly-separation}
For each $\alpha \in \rr^n \setminus \scrP$, there is $\nu \in \rnstar$ such that $\min_\scrP(\nu)$ exists and $\langle \nu, \alpha \rangle < \min_\scrP(\nu)$. \qed
\end{prop}

Geometrically, \cref{poly-separation} states that every point in the complement of $\scrP$ is separated from $\scrP$ by a hyperplane - see \cref{fig:pols-0}. A \index{Face!of a convex polyhedron}{\em face} of $\scrP$ is a subset of the form $\In_\nu(\scrP)$ for some $\nu \in \rnstar$.

\begin{figure}[h]
\def\xmin{-0.5}
\def\xmax{16.5}
\def\ymin{-4.5}
\def\ymax{5.5}
\def\tx{0}
\def\ty{5}
\def\tw{5cm}

\tikzstyle{dot} = [\colordot, circle, minimum size=4pt, inner sep = 0pt, fill]

\begin{center}

\begin{tikzpicture}[scale=\scalefactor]

\coordinate (A) at (1,-4);
\coordinate (B) at (16,1);
\coordinate (C) at (10,5);
\fill[\colorzero, opacity=\opazero ] (A) --  (B) -- (C) -- cycle;
\draw[thick, \colorone] (A) --  (B) -- (C) -- cycle;

\node[dot] at (A) {};
\node[dot] at (B) {};
\node[dot] at (C) {};

\coordinate (nu) at (0,-1);

\coordinate (eta) at (-2,-3);
\coordinate (etaperp) at (-3,2);

\node[anchor = north] at (A) {\picfontsize $A$};
\node[anchor = north] at (B) {\picfontsize $B$};
\node[anchor = south] at (C) {\picfontsize $C$};

\coordinate (A1) at (1,5);
\coordinate (B1) at (16,5);
\coordinate (O) at (6.5,5);
\draw[thick, \colortwo, dashed] (A1) -- (B1);
\draw [thick, \colornu, ->] (O) -- ($(O) + (nu)$);
\node[anchor = north] at ($(O) + (nu)$) {\picfontsize $\nu$};
\node[anchor = south] at ($(B1)$) {\picfontsize $L$};

\coordinate (out) at (3, 6);
\node[dot] at (out) {};
\node[anchor = west] at (out) {\picfontsize $\alpha$};

\coordinate (O) at ($(B)!0.5!(C)$);
\draw[thick, \colortwo, dashed] ($(B)-(etaperp)$) --  ($(C)+(etaperp)$);
\draw [thick, \colornu, ->] (O) -- ($(O) + 0.3*(eta)$);
\node[anchor = east] at ($(O) + 0.3*(eta)$) {\picfontsize $\eta$};
\node at (9,{2/3}) {\picfontsize $\scrP$};
\end{tikzpicture}

\caption{
	$\In_\nu(\scrP) = \{C\}$, $\In_\eta(\scrP) = BC$. The line $L$ separates $\scrP$ from $\alpha$.
} \label{fig:pols-0}
\end{center}
\end{figure}

\begin{prop} \label{facedron}
Every face of $\scrP$ is a convex polyhedron. If $\scrP$ is a polytope (respectively polyhedral cone), then every face of $\scrP$ is also a polytope (respectively, polyhedral cone).
\end{prop}

\begin{proof}
The case of a convex polyhedron and that of a polytope are direct consequences of the definitions of polyhedra, polytopes and faces. The case of a polyhedral cone follows from combining \cref{exercise:cone-inequalities,exercise:minconelinear}.
\end{proof}

We also note the following property whose proof is left as an exercise:

\begin{prop} \label{prop:face-invariance}
Let $\phi: \scrP\to \scrQ$ be an affine isomorphism of convex polyhedra. Then $\phi$ induces a bijection of faces, i.e.\ $\scrR \subseteq \scrP$ is a face of $\scrP$ if and only if $\phi(\scrR)$ is a face of $\scrQ$.
\end{prop}

\begin{proof}
This is \cref{exercise:face-invariance}.
\end{proof}

\Cref{facedron} implies that every face of $\scrP$ is a polyhedron, and therefore has a well-defined dimension. A \index{Vertex!of a convex polyhedron}\index{Edge!of a convex polyhedron}\index{Facet!of a convex polyhedron}{\em vertex} (respectively, {\em edge, facet}) is a face of dimension zero (respectively one, $\dim(\scrP) - 1$). In \cref{fig:pols-0} the vertices of $\scrP$ are $A,B,C$, and the edges are the three sides. \Cref{prop:finite-vertedges} below gives a more precise description of vertices and edges. The part of the proposition about finiteness of numbers of vertices and edges is a special case of the more general result (\cref{prop:finite-faces}) that a polyhedron can have only finitely many faces. We will also see later that every polytope is the convex hull of its vertices (as it is evident for $\scrP$ from \cref{fig:pols-0}).

\begin{prop} \label{prop:finite-vertedges}
Every vertex of a polyhedron is a singleton, and up to an affine isomorphism, every edge is a (possibly unbounded) closed interval of $\rr$. Every polytope has finitely many vertices. Every strongly convex polyhedral cone has only one vertex (namely, the origin) and finitely many edges.
\end{prop}

\begin{proof}
The first statement follows directly from \cref{exercise:zerone-dim-pol}. By \cref{poly-basic-thm} every polytope $\scrP$ is the convex hull of a finite set. The first statement together with \cref{exercise:conv-rep} then implies that $\scrP$ has only finitely many vertices. Now let $\scrC$ be a strongly convex polyhedral cone. \Cref{poly-basic-thm} implies that the origin is a vertex of $\scrC$. Since every face of $\scrC$ is also a cone (\cref{facedron}), and since the only zero-dimensional cone is the origin, it follows that the origin is the only vertex of $\scrC$. The finiteness of the numbers of edges of $\scrC$ follows from \cref{exercise:cone-to-pol} and finiteness of the numbers of vertices of polytopes.
\end{proof}

Note in \cref{fig:pols-0} every pair of edges of $\scrP$ intersects in a vertex. The following result shows that this corresponds to a general phenomenon:

\begin{prop} \label{prop:facection}
The intersection of finitely many faces of $\scrP$ is also a face of $\scrP$, provided the intersection is nonempty.
\end{prop}

\begin{proof}
Let $\scrQ := \In_\nu(\scrP)$ and $\scrR := \In_\eta(\scrP)$ be faces of $\scrP$ with $\scrQ \cap \scrR \neq \emptyset$. It is then straightforward to check that $\scrQ \cap \scrR = \In_{\nu + \eta}(\scrP)$. [Where is it used that $\scrQ \cap \scrR \neq \emptyset$? See \cref{exercise:facection-counter-example}.]
\end{proof}

Every polyhedron $\scrP$ is a face of itself, since $\scrP = \In_0 (\scrP)$. A \index{Proper!face of a convex polyhedron}{\em proper face} of $\scrP$ is a face which is properly contained in $\scrP$, and the \index{Relative interior of a convex polyhedron}{\em relative interior} $\relint(\scrP)$ of $\scrP$ is the complement in $\scrP$ of the union of its proper faces. In \cref{fig:pols-0} the relative interior is the complement in $\scrP$ of the union of the three edges. 

\begin{prop} \label{prop:relative-containment}
Let $\scrQ, \scrR$ be faces of $\scrP$.
\begin{enumerate}
\item \label{relatively-smaller-dim} If $\scrQ$ is a proper face of $\scrP$, then $\dim(\scrQ) < \dim(\scrP)$.
\item If $\scrQ \cap \relint(\scrR) \neq \emptyset$, then $\scrQ \supset \scrR$.
\end{enumerate}
\end{prop}

\begin{proof}
The first assertion is an immediate consequence of the observation that if $\nu \in \rnstar$ is non-constant on an affine subspace $H$ of $\rr^n$, then for every $r \in \rr$, $\{ \alpha \in H: \langle \nu, \alpha \rangle = r\}$ is an affine subspace of $\rr^n$ whose dimension is one less than that of $H$. For the second assertion, note that if $\scrQ \cap \scrR \neq \emptyset$, then by \cref{prop:facection} $\scrQ \cap \scrR$ is a face of $\scrR$, and it can not be a proper face of $\scrR$ if it contains a point from $\relint(\scrR)$.
\end{proof}

Note that in \cref{fig:pols-0} the relative interior of $\scrP$ is also its \index{Topological!interior}{\em topological interior}\footnote{A point $x \in S \subseteq \rr^n$ is in the topological interior of $S$ if and only if $S$ contains an open neighborhood of $x$ in $\rr^n$.} in $\rr^2$. We will now see that this is a manifestation of a general property of polyhedra.


\begin{prop} \label{prop:relint-inequalities}
$\scrP$ is the closure of its relative interior. In particular, $\relint(\scrP)$ is nonempty, and it is precisely the topological interior of $\scrP$ in $\aff(\scrP)$ (where $\aff(\scrP)$ is equipped with the relative topology from $\rr^n$). In particular, $\relint(\scrP)$ is a nonempty relatively open subset of $\aff(\scrP)$. In the case that $\scrP \subset \rr^n$ is ``full dimensional'' (i.e.\ $\dim(\scrP) = n$), and $a_{i,0} +  a_{i,1}x_1 + \cdots + a_{i,n}x_n \geq 0$, $i = 1, \ldots, N$, are a set of ``nontrivial inequalities''\footnote{``Nontrivial inequalities'' means that we do not allows inequalities with $a_{i,0} = \cdots = a_{i,n} = 0$.} defining $\scrP$, then $\relint(\scrP)$ is the nonempty open set of $\rr^n$ defined by the strict inequalities $a_{i,0} +  a_{i,1}x_1 + \cdots + a_{i,n}x_n > 0$, $i = 1, \ldots, N$.
\end{prop}

\begin{proof}
Due to \cref{exercise:eq-dim-embedding,exercise:affine-invariance,exercise:face-invariance} it suffices to consider the case that $\scrP$ is full dimensional. For each $i$, let $H_i$ be the hyperplane $a_{i,0} +  a_{i,1}x_1 + \cdots + a_{i,n}x_n = 0$. Let $\scrP^0 := \{x \in \rr^n:a_{i,0} +  a_{i,1}x_1 + \cdots + a_{i,n}x_n > 0,\ i = 1, \ldots, N\} = \scrP \setminus \bigcup_i H_i$. \Cref{exercise:top-bdry} implies that $\scrP^0$ is the topological interior of $\scrP$. We claim that $\scrP^0 \neq \emptyset$. Indeed, for each $i = 1, \ldots, N$, there is $\alpha_i \in \scrP$ such that $a_{i,0} +  a_{i,1}x_1 + \cdots + a_{i,n}x_n > 0$, since otherwise $\scrP$ would be contained in $H_i$, and $\dim(\scrP)$ would be less than $n$. Since $\scrP$ is convex, $\alpha := (1/N)\sum_{i=1}^N \alpha_i \in \scrP$. It is straightforward then to check that $\alpha \in \scrP^0$, as required. We now show that $\scrP^0 = \relint(\scrP)$. The full-dimensionality of $\scrP$ also implies that $H_i \cap \scrP$ is a proper face of $\scrP$ for each $i$, so that $\scrP^0 \supset \relint(\scrP)$. On the other hand, if $\alpha \in \scrP^0$, then by openness of $\scrP^0$, for every $\nu \in \rnstar \setminus \{0\}$, there is $\alpha' \in \scrP^0$ such that $\langle \nu, \alpha' \rangle < \langle \nu, \alpha \rangle$, so that $\alpha \not\in \In_\nu(\scrP)$. This shows that $\scrP^0 \subset \relint(\scrP)$, and therefore $\scrP^0 = \relint(\scrP)$. Finally, to see that $\scrP$ is the closure of its relative interior, let $\beta := (\beta_1, \ldots, \beta_n) \in \relint(\scrP)$. Given $\gamma = (\gamma_1, \ldots, \gamma_n) \in \scrP$, let $L := \{(1-\epsilon)\beta + \epsilon\gamma: 0 \leq \epsilon \leq \beta \}$ be the line segment from $\beta$ to $\gamma$. We claim that $L \setminus \{\gamma\} \subseteq \relint(\scrP)$. Indeed, for each $i = 1, \ldots,N$, let
\begin{align*}
\epsilon_i
	&:= \sup\{\epsilon: 0 \leq \epsilon \leq 1,\
		a_{i,0} +  a_{i,1}((1-\epsilon) \beta_1 + \epsilon\gamma_1) +  \cdots +  a_{i,n}((1-\epsilon) \beta_n + \epsilon\gamma_n)
		> 0 \}
\end{align*}
Note that each of the functions $a_{i,0} +  a_{i,1}((1-\epsilon) \beta_1 + \epsilon\gamma_1) +  \cdots +  a_{i,n}((1-\epsilon) \beta_n + \epsilon\gamma_n)$ is linear in $\epsilon$, and it is strictly positive at $\epsilon = 0$. Therefore if $\epsilon_i < 1$, then it must be zero at $\epsilon_i$ and negative on the interval $(\epsilon_i, 1]$, contradicting the fact that $L \subset \scrP$. Therefore $\epsilon_i = 1$ for each $i$, and consequently, $L\setminus \{\gamma\} \subset \scrP^0 = \relint(\scrP)$. This implies that $\scrP$ is the closure of $\relint(\scrP)$, as required.
\end{proof}

\begin{prop} \label{prop:sum-faces}
Let $\scrP_1, \ldots, \scrP_s$ be polyhedra in $\rr^n$, and $\nu \in \rnstar$. Then
\begin{enumerate}
\item \label{sum-faces:existence} $\min_{\sum_j \scrP_j}(\nu)$ exists if and only if $\min_{\scrP_j}(\nu)$ exists for each $j$.
\end{enumerate}
Now assume $\min_{\sum_j \scrP_j}(\nu)$ exists. Then
\begin{enumerate}[resume]
\item \label{sum-faces:min} $\min_{\sum_j \scrP_j}(\nu) = \sum_j \min_{\scrP_j}(\nu)$.
\item \label{sum-faces:sum} $\In_\nu(\sum_j \scrP_j) = \sum_j \In_\nu(\scrP_j)$.
\item \label{sum-faces:uniqueness} $\In_\nu(\scrP_j)$ are the ``unique maximal'' subsets of $\scrP_j$ whose sum is $\In_\nu(\sum_j \scrP_j)$, i.e.\ if $\alpha_j \in \scrP_j$ are such that $\sum_j \alpha_j \in \In_\nu(\sum_j \scrP_j)$, then $\alpha_j \in \In_\nu(\scrP_j)$ for each $j$.
\item \label{sum-faces:polytope-uniqueness} If each $\scrP_j$ is a polytope, then $\In_\nu(\scrP_j)$ are in fact unique faces of $\scrP_j$ whose sum is $\In_\nu(\sum_j \scrP_j)$, i.e.\ if $\scrQ_j$ are faces of $\scrP_j$ such that $\sum_j \scrQ_j = \In_\nu(\sum_j \scrP_j)$, then $\scrQ_j = \In_\nu(\scrP_j)$ for each $j$.
\item \label{sum-faces:conelinear} Assume there is $j$ such that
\begin{enumerate}
\item \label{conelinear:cone} either $\scrP_j$ is a cone, or
\item \label{conelinear:linear} $\scrP_j$ is a linear subspace of $\rr^n$.
\end{enumerate}
Then $\In_\nu(\scrP_j)$ contains the origin and $\min_{\scrP_j}(\nu) = 0$. In addition in case \eqref{conelinear:linear}, $\scrP_j \subset \nu^\perp$, and $\In_\nu(\scrP_j) = \scrP_j$.
\end{enumerate}
\end{prop}

\begin{proof}
It suffices to treat the case $s = 2$. Assume at first $\min_{\scrP_1}(\nu)$ does not exist. Then there are $\alpha_k \in \scrP_1$ such that $\langle \nu, \alpha_k \rangle \to -\infty$. If $\beta$ is an arbitrary element in $\scrP_2$, then $\langle \nu, \alpha_k + \beta \rangle \to -\infty$ as well, so that $\min_{\sum_j \scrP_j}(\nu)$ does not exist. On the other hand, if $\min_{\scrP_j}(\nu)$ exists for each $j$, then pick $\alpha_j \in \In_\nu(\scrP_j)$ for each $j$, and note that for all $\beta_j \in \scrP_j$, $\langle \nu, \sum_j \beta_j \rangle \geq \sum_j \langle \nu, \alpha_j \rangle =  \langle \nu, \sum_j \alpha_j\rangle$. This simultaneously proves assertions \eqref{sum-faces:existence} to \eqref{sum-faces:sum}. For assertion \eqref{sum-faces:uniqueness} note that if $\alpha_1 \not\in \In_\nu(\scrP_1)$, then $\langle \nu, \alpha_1 \rangle > \min_{\scrP_1}(\nu)$. It follows that for each $\alpha_2 \in \scrP_2$, $\langle \nu, \alpha_1 + \alpha_2 \rangle > \min_{\scrP_1}(\nu) + \min_{\scrP_2}(\nu)$, so that $\alpha_1 + \alpha_2 \not\in \In_\nu(\scrP_1 + \scrP_2)$. For assertion \eqref{sum-faces:polytope-uniqueness} assume $\scrQ_1 \subsetneq \In_\nu(\scrP_1)$. Since $\scrQ_1$ is a polytope, there is $\eta \in \rnstar$ such that $\min_{\scrQ_1}(\eta) > \min_{\In_\nu(\scrP_1)} (\eta)$. On the other hand, $\In_\nu(\scrP_j)$ is a polytope for each $j$, and therefore by assertion \eqref{sum-faces:existence},
\begin{align*}
\min_{\scrQ_1 + \In_\nu(\scrP_2)}(\eta)
	&= \min_{\scrQ_1}(\eta)  + \min_{\In_\nu(\scrP_2)}(\eta)
	> \min_{\In_\nu(\scrP_1)}(\eta)  + \min_{\In_\nu(\scrP_2)}(\eta)
	= \min_{\In_\nu(\scrP_1) + \In_\nu(\scrP_2)}(\eta)
\end{align*}
It follows that $\scrQ_1 + \In_\nu(\scrP_2) \neq \In_\nu(\scrP_1)(\eta)  + \In_\nu(\scrP_2)$, as required. Finally, assertion \eqref{sum-faces:conelinear} follows directly from \cref{exercise:minconelinear}.
\end{proof}

The example from \cref{fig:non-unique-sum} shows that assertion \eqref{sum-faces:polytope-uniqueness} of \cref{prop:sum-faces} may not be true if some of the $\scrP_j$ are not bounded [what goes wrong with the proof?]. On the other hand, if $\scrP_j$ are polytopes, then it is not too hard to show that it remains true even if $\scrQ_j$ are allowed to be arbitrary convex subsets of $\scrP_j$, i.e.\ $\In_\nu(\scrP_j)$ are unique {\em convex} subsets of $\scrP_j$ whose sum is $\In_\nu(\sum_j \scrP_j)$. \Cref{exercise:non-unique-sum-faces} asks you to show that convexity is necessary for uniqueness.

\begin{cor} \label{prop:finite-faces}
Every polyhedron has finitely many distinct faces.
\end{cor}

\begin{proof}
Let $\nu \in \rnstar$. If $\scrP_0 = \conv(S)$ for some finite set $S \subset \rr^n$, \cref{exercise:positively-convex} implies that $\In_\nu(\scrP_0) = \conv(\In_\nu(S))$; in particular, every face of the polytope $\scrP_0$ is the convex hull of a subset of $S$, and therefore $\scrP_0$ has finitely many distinct faces. Now consider the case that $\scrC = \cone(T)$ for a finite subset $T$ of $\rr^n$. If $\min_\scrC(\nu)$ exists, then \cref{prop:sum-faces} implies that $\min_\scrC(\nu) = 0$ and $\In_\nu(\scrC) = \nu^\perp \cap \scrC = \cone(\nu^\perp \cap T)$. Therefore, the number of distinct faces of $\scrC$ is bounded by the number of subsets of $T$. By \cref{polyhedron-minkomposition} every polyhedron $\scrP$ be can be expressed as $\scrP_0 + \scrC + L$, where $\scrP_0$ is a polytope, $\scrC$ is a strongly convex polyhedral cone, and $L$ is a linear subspace of $\rr^n$. \Cref{prop:sum-faces} then implies that every face of $\scrP$ is of the form $\scrP'_0 + \scrC' + L$, where $\scrP'_0$ (respectively, $\scrC'$) is a face of $\scrP_0$ (respectively, $\scrC$). The result then follows from the cases of polytopes and polyhedral cones.
\end{proof}

\begin{center}
\begin{figure}[h]
\def\viewx{60}
\def\viewy{30}
\begin{tikzpicture}[scale=0.75]
\pgfplotsset{view={\viewx}{\viewy}}

\begin{axis}[axis lines = none]
	\addplot3 [draw, thick, fill=\colorzero,opacity=\opazero] coordinates{(0,0,1) (0,1,0) (1,0,0) (0,0,1)};
	\addplot3 [draw, thick, fill=\colorone,opacity=\opazero] coordinates{(0,0,0) (1,0,0) (0,0,1) (0,0,0)};
	\addplot3 [dashed, \colortwo, thick] coordinates{(0,0,0) (0,1,0)};
	\addplot3 [\colortwo, thick] coordinates{(0.3,0,0.7) (0,0.8,0.2)};
\end{axis}
\draw (0,1.5) node {\picfontsize $O$};
\draw (0,4.75) node {\picfontsize $C$};
\draw (2.5,-0.1) node {\picfontsize $A$};
\draw (4.6,2.75) node {\picfontsize $B$};

\end{tikzpicture}

\caption{$ABC$ is the minimal face containing edges $AC$ and $BC$. The relative interior of the line segment from a point in $\relint(AC)$ to a point in $\relint(BC)$ is contained in $\relint(ABC)$.} \label{fig:positively-relative}
\end{figure}
\end{center}

As a corollary of the finiteness of the number of faces, we now prove a more technical result that we use in \cref{normal-section}. It shows that given two faces $\scrQ_1, \scrQ_2$ of a polyhedron $\scrP$, there is a ``smallest'' face $\scrQ$ of $\scrP$ that contains both $\scrQ_j$, and given two points in the relative interiors of the $\scrQ_j$, the relative interior of the line segment joining these points is contained in the relative interior of $\scrQ$, see \cref{fig:positively-relative}.

\begin{cor} \label{prop:positively-relative}
Let $\scrQ_1, \scrQ_2$ be faces of a convex polyhedron $\scrP$.
\begin{enumerate}
\item There is a (unique) face $\scrQ$ of $\scrP$ such that $\scrQ_1 \cup \scrQ_2 \subset \scrQ$, and $\scrQ \subset \scrQ'$ for every face $\scrQ'$ of $\scrP$ such that $\scrQ_1 \cup \scrQ_2 \subset \scrQ'$.
\item Pick $\alpha_i \in \relint(\scrQ_i)$, $i = 1,2$. Then for every $\epsilon \in (0,1)$, the positive convex combination $\epsilon \alpha_1 + (1- \epsilon)\alpha_2$ of $\alpha_1$ and $\alpha_2$ is in the relative interior of $\scrQ$.
\end{enumerate}
\end{cor}

\begin{proof}
Let $\scrQ$ be the intersection of all faces of $\scrP$ containing both $\scrQ_j$. \Cref{prop:facection,prop:finite-faces} imply that $\scrQ$ is a face of $\scrP$, and it is clearly as in the first assertion. The second assertion follows from \cref{exercise:positively-convex} and the definition of relative interior.
\end{proof}

The following proposition shows that the property of being a face is transitive. Note that this is clear in \cref{fig:pols-0,fig:positively-relative}: every vertex of every edge of $\scrP$ is also a vertex of $\scrP$.

\begin{prop} \label{prop:face-of-face}
Every face of a face of $\scrP$ is also a face of $\scrP$.
\end{prop}

\begin{proof}
Let $\scrQ := \In_\nu(\scrP)$ and $\scrR := \In_\eta(\scrQ)$, $\nu, \eta \in \rnstar$. Let $\scrP = \scrP_0 + \scrC + L$ be a decomposition as in \cref{polyhedron-minkomposition}. \Cref{prop:sum-faces} then implies that
\begin{prooflist}
\item \label{face-of-face:Q} $\scrQ = \In_\nu(\scrP_0) + \In_\nu(\scrC) + L$, where $\nu|_L \equiv 0$ and $\min_\scrC(\nu) = 0$;
\item \label{face-of-face:R} $\scrR = \In_\eta(\In_\nu(\scrP_0)) + \In_\eta(\In_\nu(\scrC)) + L$, where $\eta|_L \equiv 0$, and $\min_{\In_\nu(\scrC)}(\eta) = 0$.
\end{prooflist}
Due to \cref{prop:face-invariance} after a translation of $\scrP$ if necessary we may further assume that $\In_\nu(\scrP_0)$ contains the origin, which implies that
\begin{prooflist}[resume]
\item \label{face-of-face:min-P0-nu} $\min_{\scrP_0}(\nu) = 0$, and
\item \label{face-of-face:Q-containment} $\scrQ \supset \In_\nu(\scrP_0) \cup \In_\nu(\scrC) \cup L$.
\end{prooflist}
By \cref{poly-basic-thm} there are finite sets $S,T \subset \rr^n$ such that $\scrP_0 = \conv(S)$ and $\scrC = \cone(T)$. By observations \ref{face-of-face:Q} and \ref{face-of-face:min-P0-nu}, $\nu$ is nonnegative on $S \cup T$, and therefore we may choose $r > 0$ such that $\langle r\nu + \eta, \alpha \rangle > 0$ for each $\alpha \in (S \cup T) \setminus \nu^\perp$.

\begin{proclaim} \label{face-of-face:claim:r}
$\scrR = \In_{r\nu+ \eta}(\scrP)$.
\end{proclaim}

\begin{proof}
This is left as \cref{exercise:claim:r}.
\end{proof}
The proposition follows immediately from \cref{face-of-face:claim:r}.
\end{proof}

\begin{prop} \label{prop:vertex-hull}
Let $\scrP$ be a polytope and $\scrV$ be the set of vertices of $\scrP$.
\begin{enumerate}
\item $\scrP = \conv(\scrV)$.
\item $\scrV$ is the unique minimal set whose convex hull is $\scrP$.
\item \label{vertex-hull:relatively-positive} $\relint(\scrP)$ is the set of {\em positive} convex combinations of its vertices, i.e.\ $\alpha \in \relint(\scrP)$ if and only if $\alpha = \sum_{\beta \in \scrV} \epsilon_\beta \beta$, where $\epsilon_\beta$ are positive real numbers such that $\sum_{\beta \in \scrV} \epsilon_\beta = 1$.
\end{enumerate}
\end{prop}

\begin{proof}
For the first assertion proceed by induction on $\dim(\scrP)$. It is evident when $\dim(\scrP) = 0$. Now consider the case that $\dim(\scrP) \geq 1$. Take $\alpha \in \scrP$. If $\alpha \not\in \relint(\scrP)$, then it is on a proper face $\scrQ$ of $\scrP$. Since $\dim(\scrQ) < \dim(\scrP)$, by induction $\alpha$ is in the convex hull of vertices of $\scrQ$. \Cref{prop:face-of-face} implies that every vertex of $\scrQ$ is also a vertex of $\scrP$, so that $\alpha \in \conv(\scrV)$, as required. If $\alpha \in \relint(\scrP)$, then take a line $L$ through $\alpha$ on $\aff(\scrP)$. Since $\scrP$ is bounded, \cref{prop:relint-inequalities} implies that each end of $L$ intersects a proper face of $\scrP$. Since we already showed that every proper face of $\scrP$ is in $\conv(\scrP)$, it follows that $\alpha \in \conv(\scrV)$ and completes the proof of the first assertion. The second assertion and the $(\Leftarrow)$ implication of the third assertion follow from the first assertion and \cref{exercise:positively-convex}. Now pick $\alpha \in \relint(\scrP)$. It remains to show that $\alpha$ is a positive convex combination of the vertices of $\scrP$. By \cref{prop:finite-vertedges} we may assume $\dim(\scrP) \geq 1$. Since $\relint(\scrP)$ is relatively open in $\aff(\scrP)$ (\cref{prop:relint-inequalities}), $\alpha - \epsilon \sum_{\beta \in \scrV} \beta \in \relint(\scrP)$ for sufficiently small positive $\epsilon$. By the first assertion then $\alpha - \epsilon \sum_{\beta \in \scrV} \beta \in \conv(\scrV)$. It then follows that $\alpha$ is a positive convex combination of elements from $\scrV$, as required.
\end{proof}

\begin{cor} \label{prop:edge-hull}
If $\scrC$ is a strongly convex polyhedral cone, then
\begin{enumerate}[resume]
\item $\scrC$ is the cone generated by its edges.
\item $\relint(\scrC)$ is the set of {\em positive} linear combinations of nonzero elements of its edges, i.e.\ if $\scrE$ is a set consisting of one nonzero element from each of the edges of $\scrC$, then $\alpha \in \relint(\scrC)$ if and only if $\alpha = \sum_{\beta \in \scrE} r_\beta \beta$, where $r_\beta$ are positive real numbers.
\end{enumerate}
\end{cor}

\begin{proof}
Combine \cref{prop:vertex-hull,exercise:cone-to-pol}.
\end{proof}

%

\begin{center}
\begin{figure}[h]

\def\scalefactor{.3}
\def\xmin{-2}
\def\xmax{9}
\def\ymin{-6}
\def\ymax{3}
\def\nux{6}

\tikzstyle{dot} = [\colordot, circle, minimum size=4pt, inner sep = 0pt, fill]

\begin{tikzpicture}[scale=\scalefactor]

\coordinate (O) at (0,0);
\coordinate (A) at (3,-2);
\coordinate (B) at (4,-1);
\coordinate (C) at (4,4/3);
\coordinate (A') at (9,-6);
\coordinate (C') at (9,3);

\fill [\colorzero, opacity=\opazero ] (O) --  (A') -- (C') -- cycle;

\draw[thick, \colorone ] (O) --  (A');
\draw[thick, \colorone ] (O) --  (C');

\draw (\xmax, 0) -- (\xmin, 0);
\draw (0, \ymax) -- (0, \ymin);

\node[dot] at (O){};
\node[dot] at (A){};
\node[dot] at (B){};
\node[dot] at (C){};
\node[anchor = west] at (\xmax, 0) {\picfontsize $x_{k+1}$};
\node[anchor = south] at (0, \ymax) {\picfontsize $x_{k+2}$};

\coordinate (nu1) at (2,3);
\coordinate (nu2) at (1,-3);
\pgfmathsetmacro\nuratio{sqrt(1.3)};
\def\nuonescale{0.4}

\pgfmathsetmacro\nutwoscale{\nuonescale*\nuratio};
\coordinate (nu1base) at (\nux, {-2*\nux/3});
\coordinate (nu1head) at ($(nu1base) + \nuonescale*(nu1)$);
\draw [thick, \colornu, ->] (nu1base) -- (nu1head);
\node[anchor = west] at (nu1head) {\picfontsize $\eta_1$};

\coordinate (nu2base) at (\nux, {\nux/3});
\coordinate (nu2head) at ($(nu2base) + \nuonescale*\nuratio*(nu2)$);
\draw [thick, \colornu, ->] (nu2base) -- (nu2head);
\node[anchor = west] at (nu2head) {\picfontsize $\eta_2$};

\end{tikzpicture}

\caption{
$\cone(\pi(S \cup T))$ is a two dimensional cone contained in the half-plane $x_{k+1} \geq 0$, and intersects $x_{k+2}$-axis only at the origin.
} \label{fig:pos-2jection}

\end{figure}
\end{center}

\begin{prop} \label{prop:facet-containment}
Every proper face $\scrQ$ of a polyhedron $\scrP$ is the intersection of facets of $\scrP$ containing $\scrQ$. In particular, every proper face is contained in a facet.
\end{prop}

\begin{proof}
Let $\scrP = \scrP_0 + \scrC+L$ be a decomposition of $\scrP$ as in \cref{polyhedron-minkomposition}. Let $\scrQ = \In_\nu(\scrP)$, $\nu \in \rnstar$, be a proper face of $\scrP$. \Cref{prop:sum-faces} implies that $\scrQ = \In_\nu(\scrP_0) + \In_\nu(\scrC) + L$, where $\nu|_L \equiv 0$ and $\min_\scrC(\nu) = 0$. Due to \cref{prop:face-invariance} after an affine isomorphism of $\scrP$ if necessary we may further assume that $\dim(\scrP) = n$, and $\aff(\scrQ)$ is the coordinate subspace of $\rr^n$ spanned by the first $k$-coordinates (in particular, $\min_{\scrP}(\nu) = \min_{\scrP_0}(\nu) = 0$), where $k := \dim(\scrQ) < n$, and $\nu$ is the projection on to the $(k+1)$-th coordinate. If $n=k+1$, then $\scrQ$ is a facet, and we are done. So assume $n - k \geq 2$. Then after a linear isomorphism of $\scrP$ we can also ensure that each point of $\scrQ$ has zero $(k+2)$-th coordinate. Choose finite sets $S, T$ such that $\scrP_0 = \conv(S)$ and $\scrC = \cone(T)$. Let $\pi$ be the projection map from $\rr^n$ onto the two-dimensional coordinate subspace spanned by the $(k+1)$-th and $(k+2)$-th coordinate. By construction every nonzero element of $\pi(S \cup T)$ has {\em positive} $x_{k+1}$-coordinate. We claim that the cone generated by $\pi(S \cup T)$ is two dimensional, i.e.\ it is as in \cref{fig:pos-2jection}. Indeed, by construction $\pi(L) = 0$, so that \cref{exercise:affine-image,exercise:linear-sum-invariance} imply that $\pi(\scrP) = \pi(\scrP_0) + \pi(\scrC) = \conv(\pi(S)) + \cone(\pi(T)) \subseteq \cone(\pi(S \cup T))$. Therefore, if $\cone(\pi(S \cup T))$ is contained in a line $L$, and $\scrP$ would be contained in $\pi^{-1}(L)$, contradicting the full-dimensionality of $\scrP$. It follows that $\dim(\cone(\pi(S\cup T))) = 2$, and $\cone(\pi(S \cup T))$ has two edges (\cref{exercise:2-cone}). As in \cref{fig:pos-2jection}, let $\eta_1, \eta_2$ be the elements in $\rnnstar{2}$ which attains their minima on the two edges of $\cone(\pi(S \cup T))$. It is straightforward to check that
\begin{itemize}
\item for each $j$, $\scrR_j := \In_{\pi^*(\eta_j)}(\scrP)$ is a face of $\scrP$ such that $\scrR_j \supset \scrQ$ and $\dim(\scrR_j) > \dim(\scrQ)$,
\item If $n = k+2$, then each $\scrR_j$ is a facet of $\scrP$, and $\scrQ = \scrR_1 \cap \scrR_2$.
\end{itemize}
The proposition follows from these observations by a straightforward induction on $\dim(\scrP) - \dim(\scrQ)$.
\end{proof}

\begin{cor} \label{prop:minimal-polynequalities}
Let $\scrP$ be an $n$-dimensional polyhedron in $\rr^n$. Assume $\scrP \neq \rr^n$. Then $\scrP$ is determined by the affine hyperplanes corresponding to facets of $\scrP$. More precisely,
\begin{enumerate}
\item up to multiplications by nonzero real numbers, there is a unique minimal set of inequalities determining $\scrP$.
\end{enumerate}
Let $a_{i,0} +  a_{i,1}x_1 + \cdots + a_{i,n}x_n \geq 0$, $i = 1, \ldots, N$, be the minimal set of inequalities defining $\scrP$. Then
\begin{enumerate}[resume]
\item $\scrP \cap \{x : a_{i,0} +  a_{i,1}x_1 + \cdots + a_{i,n}x_n = 0\}$ is a facet of $\scrP$ for each $i$;
\item for every face $\scrQ$ of $\scrP$, there is $I \subseteq \{1, \ldots, N\}$ such that $\scrQ = \scrP \cap \{x: a_{i,0} +  a_{i,1}x_1 + \cdots + a_{i,n}x_n = 0$ for each $i \in I\}$.
\end{enumerate}
\end{cor}

\begin{proof}
Let $a_{i,0} +  a_{i,1}x_1 + \cdots + a_{i,n}x_n \geq 0$, $i = 1, \ldots, N$, be a minimal set of inequalities defining $\scrP$. It suffices to show that $\scrQ_i := \scrP \cap \{x : a_{i,0} +  a_{i,1}x_1 + \cdots + a_{i,n}x_n = 0\}$, $i = 1, \ldots, N$, are the facets of $\scrP$. Indeed, if $\scrQ$ is a facet of $\scrP$, then \cref{prop:relint-inequalities} implies that $\scrQ \subset \bigcup_i \scrQ_i$. \Cref{exercise:face-union} and assertion \eqref{relatively-smaller-dim} of \cref{prop:relative-containment} then imply that $\scrQ = \scrQ_i$ for some $i$. It remains to show that every $\scrQ_i$ is a facet. Indeed, reorder the inequalities so that $\scrQ_i$ are facets for $i = 1, \ldots, M$. Let $\scrP'$ be the polytope determined by  $a_{i,0} +  a_{i,1}x_1 + \cdots + a_{i,n}x_n \geq 0$, $i = 1, \ldots, M$. If $M < N$, then $\scrP \subsetneq \scrP'$ and \cref{prop:relint-inequalities} implies that there is $\alpha \in \relint(\scrP') \setminus \scrP$. Pick $\beta \in \relint(\scrP)$ and let $L$ be the line segment from $\alpha$ to $\beta$. Then \cref{prop:relative-containment,prop:relint-inequalities} imply that
\begin{itemize}
\item for each $i = 1, \ldots, M$, and each $x \in L$, $a_{i,0} +  a_{i,1}x_1 + \cdots + a_{i,n}x_n > 0$.
\item $L$ contains a point $\gamma$ on the \index{Topological!boundary}{\em topological boundary}\footnote{Topological boundary of $S \subset \rr^n$ is the complement in $S$ of the topological interior of $S$.} of $\scrP$.
\end{itemize}
But then \cref{prop:relative-containment} implies that $\gamma$ is a point of a face of
 $\scrP$ which is not on any facet of $\scrP$, contradicting \cref{prop:facet-containment}.
\end{proof}

The theory of toric varieties (to be discussed in \cref{toric-intro}) exploits many similarities between polytopes and algebraic varieties. Here we note some of the more obvious analogues between a polyhedron $\scrP$ and an irreducible variety $X$.

\begin{align*}
\begin{array}{lcl}
\parbox{0.36\textwidth}{
	The dimension of every proper irreducible subvariety of $X$ is less than $\dim(X)$ (\cref{thm:dim-smaller}).
}
&
\leftrightarrow
&
\parbox{0.3\textwidth}{
	The dimension of every proper face of $\scrP$ is less than $\dim(\scrP)$ (\cref{prop:relative-containment}).
}\\
&
&\\
\parbox{0.36\textwidth}{
	If an irreducible subvariety of $X$ is contained in the union of finitely many subvarieties of $X$, then it is contained in one of them (\cref{prop:irreducible-properties}).
}
&
\leftrightarrow
&
\parbox{0.3\textwidth}{
	If a face of $\scrP$ is contained in the union of finitely many faces of $\scrP$, then it is contained in one of them (\cref{exercise:face-union}).
}
\end{array}
\end{align*}

\subsection{Exercises}
\begin{exercise} \label{exercise:afface-hull}
Let $\scrQ$ be a face of $\scrP$. Show that $\scrQ = \scrP \cap \aff(\scrQ)$. [Hint: if $\scrQ = \In_\nu(\scrP)$, then $\langle \nu, \cdot \rangle$ is constant on $\aff(\scrQ)$.]
\end{exercise}

\begin{exercise} \label{exercise:minconelinear}
Let $\scrP$ be a convex polyhedron and $\nu \in \rnstar$. Assume $\inf\{\langle \nu, \alpha \rangle: \alpha \in \scrP \} > -\infty$, and
\begin{enumerate}
\item \label{minconelinear:cone} either $\scrP$ is a cone, or
\item \label{minconelinear:linear} $\scrP$ is a linear subspace of $\rr^n$.
\end{enumerate}
Then $\In_\nu(\scrP)$ contains the origin and $\min_{\scrP}(\nu) = 0$. In addition, in case \eqref{minconelinear:linear} $\scrP \subset \nu^\perp$, and $\In_\nu(\scrP) = \scrP$. [Hint: in either case $\scrP$ contains the origin, so that $\min_{\scrP}(\nu) \leq 0$. If there is $\alpha \in \scrP$ such that $\langle \nu, \alpha \rangle < 0$, then $\langle \nu, r\alpha \rangle$ approaches $-\infty$ as $r \to \infty$.]
\end{exercise}

\begin{exercise} \label{exercise:face-invariance}
Prove \cref{prop:face-invariance}. [Hint: use assertion \eqref{dual-extension-0} of \cref{exercise:dual-extension}.]
\end{exercise}

\begin{exercise} \label{exercise:top-bdry}
Let $S \subset \rr^n$ and $\nu$ be a nonzero element in $\rnstar$ such that $m := \inf\{\langle \nu, \alpha \rangle : \alpha \in S\} > -\infty$. Show that no point of the hyperplane $\{x \in \rr^n: \langle \nu, \alpha \rangle = m\}$ is in the topological interior of $S$ in $\rr^n$.
\end{exercise}


\begin{exercise}\label{exercise:cone-to-pol}
Let $\scrC$ be a positive dimensional strongly convex polyhedral cone. By \cref{poly-basic-thm} there is $\nu \in \rnstar$ such that $\nu$ is positive on $\scrC \setminus \{0\}$. Let $\scrP :=  \{\alpha \in \scrP: \langle \nu, \alpha \rangle = 1\}$. \Cref{bounded-nu=1} implies that $\scrP$ is a polytope.
\begin{enumerate}
\item If $\scrQ$ is a face of $\scrP$, show that $\cone(\scrQ)$ is a face of $\scrC$ of dimension $\dim(\scrQ) + 1$.
\item Show that the correspondence $\scrQ \mapsto \cone(\scrQ)$ induces a bijection between faces of $\scrP$ and positive dimensional faces of $\scrC$,
\item \label{cone-to-pol:relint} Show that the above correspondence also preserves the relative interiors of the faces, i.e.\ If $\scrQ$ is a face of $\scrP$, then $\relint(\cone(\scrQ)) = \cone(\relint(\scrQ))$.
\end{enumerate}
[Hint: change coordinates on such that $\nu$ is the projection onto the last coordinate.]
\end{exercise}


\begin{exercise} \label{exercise:facection-counter-example}
Let $\scrP$ be a polyhedron, and $\scrQ = \In_\nu(\scrP)$ and $\scrR = \In_\eta(\scrP)$ be faces of $\scrP$. If $\scrQ \cap \scrR = \emptyset$, show by an example that it may not be true that $\In_{\nu + \eta}(\scrP) = \scrQ \cap \scrR$.
\end{exercise}

\begin{exercise} \label{exercise:non-unique-sum-faces}
Find examples of polytopes $\scrP_1, \scrP_2$ in $\rr^n$ and faces  $\scrQ_j$ of $\scrP_j$ such that $\scrQ_1 + \scrQ_2$ is a face of $\scrP_1 + \scrP_2$, and there are $\scrQ'_j \subsetneq \scrQ_j$ such that $\scrQ'_1 + \scrQ'_2 = \scrQ_1 + \scrQ_2$. [Hint: there are examples with $n = 1$ with $\scrQ_j = \scrP_j$ for each $j$.]
\end{exercise}

\begin{exercise} \label{exercise:positively-convex}
Let $\alpha = \sum_{j=1}^k \epsilon_j \alpha_j$ be a convex combination of $\alpha_1, \ldots, \alpha_k \in \rr^n$. Assume each $\epsilon_j$ is positive.
\begin{enumerate}
\item If $\nu \in \rnstar$, then show that $\langle \nu, \alpha \rangle \geq \min_j \langle \nu, \alpha_j \rangle$, with equality if and only if $\langle \nu, \alpha_1 \rangle = \cdots = \langle \nu, \alpha_k \rangle$.
\item Conclude that if the $\alpha_j$ are points of a polyhedron $\scrP$, and $\scrQ$ is a face of $\scrP$, then $\scrQ$ contains $\alpha$ if and only if $\scrQ$ contains each $\alpha_j$.
\end{enumerate}
\end{exercise}

\begin{exercise} \label{exercise:relatively-positive-interior}
Let $S$ be a finite subset of $\rr^n$. Let $\scrP := \conv(S)$ and $\scrC$ be the cone generated by $S$. Show that
\begin{enumerate}
\item $\alpha \in \relint(\scrP)$ if and only if $\alpha = \sum_{\beta \in S} \epsilon_\beta \beta$, where $\epsilon_\beta$ are positive real numbers such that $\sum_{\beta \in S} \epsilon_\beta = 1$.
\item $\alpha \in \relint(\scrC)$ if and only if $\alpha = \sum_{\beta \in S} r_\beta \beta$, where $r_\beta$ are positive real numbers.
\end{enumerate}
[Hint: follow the arguments from the proof of assertion \eqref{vertex-hull:relatively-positive} of \cref{prop:vertex-hull}.]
\end{exercise}

\begin{exercise} \label{exercise:face-union}
Let $\scrQ_1, \ldots, \scrQ_k$ be faces of a polyhedron $\scrP$. If $\scrQ$ is a convex subset (e.g.\ a face) of $\scrP$ such that $\scrQ \subseteq \bigcup_j \scrQ_j$, then show that $\scrQ \subseteq \scrQ_j$ for some $j$. [Hint: otherwise for each $j$, there is $\alpha_j \in \scrQ \setminus \scrQ_j$. Apply \cref{exercise:positively-convex} to a positive convex combination of the $\alpha_j$.]
\end{exercise}

\begin{exercise} \label{exercise:edge-graph}
Given distinct vertices $\alpha, \alpha'$ of a polytope $\scrP$, show that there are vertices $\beta_0 = \alpha, \beta_2, \ldots, \beta_k = \alpha'$ of $\scrP$ such that there is an edge of $\scrP$ connecting $\beta_{j-1}$ to $\beta_j$ for each $j = 1, \ldots, k$; in other words, the ``edge-graph'' of $\scrP$ is connected. [Hint: reduce to the case that $\scrP$ is full dimensional. Pick $\nu, \nu' \in \rnstar$ such that $\In_\nu(\scrP) = \{\alpha\}$ and $\In_{\nu'}(\scrP) = \{\alpha'\}$. Consider $\scrP_\epsilon := \In_{\nu + \epsilon \nu'}(\scrP)$ for $\epsilon \geq 0$. As $\epsilon$ goes from $0$ to $\infty$, $\scrP_\epsilon$ transitions from $\{\alpha\}$ to $\{\alpha'\}$ in finitely many steps. Apply induction on dimension to each $\scrP_\epsilon$.]
\end{exercise}

\begin{exercise} \label{exercise:claim:r}
Prove \cref{face-of-face:claim:r}. [Hint: every $\alpha \in \scrP$ can be expressed as $\sum_{\alpha \in S} \epsilon_\alpha \alpha + \sum_{\beta \in T} r_\beta \beta + \gamma$, where $\gamma \in L$, $r_\beta$ are nonnegative real numbers, and $\epsilon_\alpha$ are nonnegative real numbers such that $\sum_{\alpha \in S} \epsilon_\alpha = 1$. Define $\alpha' := \sum_{\alpha \in S \cap \nu^\perp} \epsilon_\alpha \alpha$, $\alpha'' := \sum_{\alpha \in S \setminus \nu^\perp} \epsilon_\alpha \alpha$, $\beta' := \sum_{\beta \in T \cap \nu^\perp} r_\beta \beta$ and $\beta'' := \sum_{\beta \in T \setminus \nu^\perp} r_\beta \beta$. Examine the value of $r\nu + \eta$ on each of $\alpha', \alpha'', \beta', \beta''$. Show that $\langle r\nu + \eta, \alpha \rangle \geq \min_{\In_\nu(\scrP_0)}(\eta)$, with equality if and only if $\alpha'' = \beta'' = 0$, $\alpha' \in \In_\eta(\In_\nu(\scrP_0))$ and $\beta' \in \In_\eta(\In_\nu(\scrC))$.]
\end{exercise}

\section{Normal fan of a convex polytope} \label{normal-section}
A \index{Fan}{\em fan}\footnote{Our definition of a fan differs from the definition in standard texts on toric varieties (e.g.\cite{fultoric,littlehalcox}) in that we do not require the cones in a fan to be strongly convex.} in $\rr^n$ is a collection $\Sigma$ of convex polyhedral cones in $\rr^n$ such that
\begin{enumerate}
\item \label{cone-face} Each face of a cone in $\Sigma$ is also a cone in $\Sigma$.
\item \label{cone-intersection} The intersection of any two cones in $\Sigma$ is a face of each of them.
\end{enumerate}
Any (finite dimensional) vector space $V$ over real numbers can be identified with $\rr^n$ after choosing a basis $\scrB$, and thus the notions of convex polyhedra, cones, polytopes, fans, etc., can be extended to subsets of $V$. In this section we take $V = \rnstar$ and $\scrB$ to be the basis dual to the standard basis of $\rr^n$. \Cref{exercise:affine-image} implies that convex polyhedra (and therefore cones, polytopes, fans, etc) in $\rnstar$ which are defined in this way remain so after linear changes of coordinates on $\rr^n$. Let $\scrP$ be a convex polytope in $\rr^n$. For each face $\scrQ$ of $\scrP$, define
\begin{align}
\sigma_\scrQ := \{\nu \in \rnstar: \In_\nu(\scrP) \supseteq \scrQ\} \subset \rnstar
\label{sigma-Q}
\end{align}
It is straightforward to see that $\sigma_\scrQ$ is a convex polyhedral cone in $\rnstar$ (\cref{exercise:dual-cone}); it is called the \index{Normal!cone}{\em normal cone} of $\scrQ$. Let $\Sigma_\scrP := \{\sigma_\scrQ: \scrQ$ is a face of $\scrP\}$. We will show in \cref{normal-fan} below that $\Sigma_\scrP$ is a fan in $\rnstar$; this is called the \index{Fan!normal}\index{Normal!fan}{\em normal fan} of $\scrP$. Given a face $\scrQ$ of $\scrP$, let $\sigma^0_\scrQ := \{\nu \in \rnstar: \In_\nu(\scrP) = \scrQ\} \subset \sigma_\scrQ$. We show in \cref{normal-faces} below that $\sigma^0_\scrQ$ is the relative interior of $\sigma_\scrQ$.

\begin{figure}[h]
\begin{center}
\def\xmin{-0.5}
\def\xmax{16.5}
\def\ymin{-5.5}
\def\ymax{5.5}

\tikzstyle{dotbig} = [circle, minimum size=5pt, inner sep = 0pt, fill]

\begin{tikzpicture}[scale=\scalefactor]
\coordinate (A) at (1, -4);
\coordinate (B) at (16, 1);
\coordinate (C) at (10, 5);
\fill[\colorzero, opacity=\opazero ] (A) --  (B) -- (C);
\draw[ultra thick, \colorAB] (A) -- (B);
\draw[ultra thick, \colorBC] (B) -- (C);
\draw[ultra thick, \colorCA] (C) -- (A);

\node[dotbig, \colorA] at (A) {};
\node[dotbig, \colorB] at (B) {};
\node[dotbig, \colorC] at (C) {};
\node[anchor = north] at (A) {\picfontsize $A$};
\node[anchor = west] at (B) {\picfontsize $B$};
\node[anchor = south] at (C) {\picfontsize $C$};

\node at (9.5,0.5) {\picfontsize $\scrP$};

\def\xmin{-8.5}
\def\xshift{5}
\pgfmathsetmacro\shiftone{\xshift + \xmax -\xmin};
\def\xmax{8.5}

\begin{scope}[shift={(\shiftone,0)}]
\coordinate (O) at (0,0);

\coordinate (ABperp) at (-\ymax/3, \ymax);
\coordinate (BCperp) at (\ymin*2/3, \ymin);
\coordinate (CAperp) at (-\ymin, \ymin);

\fill[\colorB, opacity=\opazero ] (O) --  (ABperp) -- (\xmin,\ymax) -- (\xmin,\ymin) -- (BCperp);
\fill[\colorC, opacity=\opazero ] (O) --  (BCperp) -- (CAperp);
\fill[\colorA, opacity=\opazero ] (O) --  (CAperp) -- (\xmax, \ymin) -- (\xmax,\ymax) -- (ABperp);

\draw [ultra thick, \colorAB] (O) -- (ABperp);
\draw [ultra thick, \colorBC] (O) -- (BCperp);
\draw [ultra thick, \colorCA] (O) -- (CAperp);


\node[dotbig, \colorzero] at (O) {};

\node at (-4.5,0.5) {\picfontsize $\sigma^0_B$};
\node at (4.5,1.5) {\picfontsize $\sigma^0_A$};
\node at (.5,-3.5) {\picfontsize $\sigma^0_C$};

\node [anchor= north] at (BCperp) {\picfontsize $\sigma^0_{BC}$};
\node [anchor=north] at (CAperp) {\picfontsize $\sigma^0_{CA}$};
\node [anchor=south] at (ABperp) {\picfontsize $\sigma^0_{AB}$};

\end{scope}
\end{tikzpicture}
\caption{Normal fan of a triangle $\scrP$ in $\rr^2$. $\sigma_\scrP = \sigma^0_\scrP$ is the origin.}  \label{fig:normal-fan}
\end{center}
\end{figure}

\begin{prop} \label{prop:dual-interior}
Let $\scrR_1, \ldots,\scrR_k$ be the faces of $\scrP$ containing $\scrQ$. Then $\sigma^0_\scrQ = \sigma_\scrQ \setminus \bigcup_j \sigma_{\scrR_j}$.
\end{prop}

\begin{proof}
This follows directly from the definition of a face and the definition of $\sigma_\scrQ$.
\end{proof}

Let $L_\scrQ$ be the linear subspace of $\rr^n$ spanned by all vectors of the form $\alpha - \beta$ with $\alpha, \beta \in \scrQ$ and $L_\scrQ^\perp$ be the linear subspace of $\rnstar$ consisting of all $\nu \in \rnstar$ such that $\nu|_{L_\scrQ} \equiv 0$.

\begin{prop} \label{affine-dimension}
$L_\scrQ^\perp = \aff(\sigma_\scrQ) = \aff(\sigma^0_\scrQ)$. In particular, $\dim(\sigma_\scrQ) = n - \dim(\scrQ)$.
\end{prop}

\begin{proof}
Since $\scrQ$ is a face of $\scrP$, there is $\nu' \in \rnstar$ such that $\scrQ= \min_\scrP(\nu')$. It follows that $\nu' \in \sigma^0_\scrQ$; in particular $\sigma^0_\scrQ$ is nonempty.

\begin{proclaim} \label{interior-claim}
Let $\nu \in \sigma^0_\scrQ$. Then for each $\mu \in L_\scrQ^\perp$, $\nu + \epsilon \mu \in \sigma^0_\scrQ$ for all sufficiently small positive $\epsilon$.
\end{proclaim}

\begin{proof}
Indeed, let $m := \min_\scrP(\nu)$ and pick $\epsilon$ such that for all vertices $\alpha$ of $\scrP$ not on $\scrQ$, $\epsilon \langle -\mu, \alpha \rangle < \langle \nu, \alpha \rangle  - m$.
\end{proof}

 \Cref{interior-claim} implies that $\aff(\sigma^0_\scrQ) \supset L_\scrQ^\perp$. On the other hand, if $\nu \in \sigma^0_\scrQ$, then for each $\alpha, \beta \in \scrQ$, $\langle \nu, \alpha - \beta \rangle = 0$, so that $\nu \in L_\scrQ^\perp$. It follows that $\aff(\sigma^0_\scrQ) = L_\scrQ^\perp$. Since $\sigma_\scrQ$ is the union of the $\sigma^0_\scrR$ over all faces $\scrR$ of $\scrP$ containing $\scrQ$, and since $L_\scrR^\perp \subset L_\scrQ^\perp$ whenever $\scrR$ contains $\scrQ$, it follows that $\aff(\sigma_\scrQ) = L_\scrQ^\perp$ as well, as required.
\end{proof}

\begin{cor} \label{normal-faces}
Let $\scrQ$ be a face of $\scrP$.
\begin{enumerate}
\item \label{face-pol-to-normal} If $\scrR$ is a face of $\scrP$ containing $\scrQ$, then $\sigma_\scrR = L_\scrR^\perp \cap \sigma_\scrQ$. In particular, $\sigma_\scrR$ is a face of $\sigma_\scrQ$.
\item \label{normal-interior} $\sigma^0_\scrQ$ is the relative interior of $\sigma_\scrQ$.
\item \label{face-normal-to-pol}  Every face of $\sigma_\scrQ$ is of the form $\sigma_\scrR$ for some face $\scrR$ of $\scrP$ containing $\scrQ$.
\end{enumerate}
\end{cor}

\begin{proof}
In this proof we will identify $\rr^n$ with $(\rnstar)^*$ in the usual way by treating $\alpha \in \rr^n$ as the linear function on $\rnstar$ given by $\nu \mapsto \langle \nu, \alpha \rangle$. \Cref{affine-dimension} implies that $\sigma_\scrR \subset L_\scrR^\perp \cap \sigma_\scrQ$.
Now pick $\nu \in L_\scrR^\perp \cap \sigma_\scrQ$. Since $\nu \in L_\scrR^\perp$, it follows that $\nu$ is constant on $\scrR$. On the other hand, since $\nu \in \sigma_\scrQ$, it follows that $\min_\scrP(\nu)$ is achieved on $\scrQ$. Since $\scrQ \subset \scrR$, it follows that $\min_\scrP(\nu)$ is achieved on all of $\scrR$, and $\nu \in \sigma_\scrR$. Therefore $\sigma_\scrR = L_\scrR^\perp \cap \sigma_\scrQ$. For assertion \eqref{face-pol-to-normal} it remains to show that $\sigma_\scrR$ is a face of $\sigma_\scrQ$.

\begin{proclaim}\label{min-beta-alpha}
Let $\alpha \in \scrQ$ and $\beta \in \scrR$. Then $\beta - \alpha$ induces a nonnegative linear function on $\sigma_\scrQ$. In particular, $\min_{\sigma_\scrQ}(\beta - \alpha) = 0$.
\end{proclaim}

\begin{proof}
If $\langle \nu, \beta - \alpha \rangle < 0$ for some $\nu \in \sigma_\scrQ$, then $\alpha \not\in \In_\nu(\scrP)$, which contradicts the definition of $\sigma_\scrQ$.
\end{proof}

Choose a basis of $L_\scrR$ of the form $\beta_j - \alpha$, $j = 1, \ldots, \dim(L_\scrR)$, where $\alpha \in \scrQ$, and each $\beta_j \in \scrR$. Then \cref{min-beta-alpha} implies that $\sigma_\scrR = L_\scrR^\perp \cap \sigma_\scrQ = \bigcap_j( \sigma_\scrQ \cap (\beta_j - \alpha)^\perp) = \bigcap_j(\In_{\beta_j - \alpha}(\sigma_\scrQ ))$ is an intersection of finitely many faces of $\sigma_\scrQ$, and therefore is also a face of $\sigma_\scrQ$. This completes the proof of assertion \eqref{face-pol-to-normal}, and due to \cref{prop:dual-interior} also implies that $\sigma^0_\scrQ$, being the complement of a union of faces of $\sigma_\scrQ$, contains $\relint(\sigma_\scrQ)$. On the other hand, \cref{interior-claim} implies that for each $\nu \in \sigma_\scrQ$, one can fit inside $\sigma_\scrQ$ a $\dim (\sigma_\scrQ)$-dimensional ball $B_\nu$ centered at $\nu$. \Cref{prop:relint-inequalities} then implies that $\sigma^0_\scrQ \subset \relint(\sigma_\scrQ)$. Consequently, $\sigma^0_\scrQ = \relint(\sigma_\scrQ)$, which proves assertion \eqref{normal-interior}. By assertions \eqref{face-pol-to-normal} and \eqref{normal-interior}, and \cref{prop:finite-faces,exercise:face-union}, every proper face of $\sigma_\scrQ$ is a face of $\sigma_\scrR$ for some face $\scrR$ of $\scrP$ properly containing $\scrQ$. Assertion \eqref{face-normal-to-pol} therefore follows from a straightforward induction on $\dim(\scrP) - \dim(\scrQ)$.
\end{proof}

\begin{cor} \label{normal-fan}
$\Sigma_\scrP$ is a fan in $\rnstar$.
\end{cor}

\begin{proof}
\Cref{normal-faces} shows that $\Sigma_\scrP$ satisfies property \eqref{cone-face} of a fan. Let $\scrQ_1$ and $\scrQ_2$ be faces of $\scrP$, and let $\scrQ$ from \cref{prop:positively-relative} be the smallest face of $\scrP$ containing both $\scrQ_j$. It is clear that $\sigma_\scrQ \subset \sigma_{\scrQ_1} \cap \sigma_{\scrQ_2}$. On the other hand, since $\scrQ$ is contained in every face of $\scrP$ containing both $\scrQ_j$, it follows that for every $\nu \in \sigma_{\scrQ_1}  \cap \sigma_{\scrQ_2}$, $\In_\nu(\scrP) \supset \scrQ$, so that $\nu \in \sigma_\scrQ$. Therefore $\sigma_{\scrQ_1}  \cap \sigma_{\scrQ_2}= \sigma_\scrQ$, which is a face of each $\sigma_{\scrQ_j}$ (\cref{normal-faces}). This shows $\Sigma_\scrP$ satisfies property \eqref{cone-intersection} of a fan as well.
\end{proof}

Recall that in our definition the cones in a fan do not have to be strongly convex. They do however turn out to be so in an important case.

\begin{prop} \label{stritctly-convex-fan}
If $\scrP$ is an $n$-dimensional convex polytope in $\rr^n$, then each cone of $\Sigma_\scrP$ is strongly convex.
\end{prop}

\begin{proof}
Indeed, if $\dim(\scrP) = n$, then for every $\nu \in \rnstar \setminus \{0\}$, $\scrP \setminus \In_\nu(\scrP)$ is nonempty. If $\nu \in  \sigma_\scrQ$ for some face $\scrQ$ of $\scrP$, it then follows from \eqref{sigma-Q} that $\In_{-\nu}(\scrP) \cap \scrQ = \emptyset$, i.e.\ $-\nu \not\in \sigma_\scrQ$, as required.
\end{proof}

\subsection{Exercises}
\begin{exercise} \label{exercise:deformation-relint}
Let $\scrP$ be an $n$-dimensional polyhedron in $\rr^n$, $\scrQ$ be a facet of $\scrP$, and $\alpha \in \relint(\scrQ)$. Let $\scrQ = \In_\nu(\scrP)$, $\nu \in \rnstar$. If $\beta$ is an arbitrary element of $\rr^n$ such that $\langle \nu, \beta \rangle \geq 0$, then show that $r\alpha + \beta \in r\scrP$ for all $r \gg 1$. [Hint: use \cref{prop:minimal-polynequalities} to show that $\alpha + \epsilon \beta \in \scrP$ if $\epsilon$ is a sufficiently small positive number.]
\end{exercise}

\begin{exercise} \label{exercise:dual-cone}
Let $\scrQ$ be a face of a polytope $\scrP$. Show that $\sigma_\scrQ$ defined in \cref{normal-section} is a convex polyhedral cone in $\rnstar$. [Hint: express $\scrP$ and $\scrQ$ as convex hulls of finite sets. Show that the condition defining $\sigma_\scrQ$ can be expressed in terms of finitely many linear inequalities as in \cref{exercise:cone-inequalities}.]
\end{exercise}

\section{Rational polyhedra} \label{rational-polysection}
We say that a convex polyhedron $\scrP$ is \index{Rational!polyhedron}\index{Polyhedron!rational}{\em rational} if it can be defined in $\rr^n$ by finitely many inequalities of the form $a_0 +  a_1x_1 + \cdots + a_nx_n \geq 0$ where each $a_i$ is a rational number (or equivalently, an integer). In particular, a {\em rational affine hyperplane} is the set of zeroes in $\rr^n$ of an equation $a_0  + a_1x_1 + \cdots + a_nx_n = 0$, where each $a_i \in \qq$, and $(a_1, \ldots, a_n) \neq (0, \ldots, 0)$. \Cref{poly-rational-thm} below gives a characterization of rational polyhedra. In particular, it implies that a polytope is rational if and only if its vertices have rational coordinates, and a polyhedral cone is rational if and only if each of its edges have a nonzero element with rational coordinates. We say that a polytope in $\rr^n$ is \index{Integral!polytope}\index{Polytope!integral}{\em integral} if all its vertices are in $\zz^n$.

%

\begin{prop} \label{poly-rational-thm}
Let $\scrP \subset \rr^n$.
\begin{enumerate}
\item \label{characterize-ratiotope} $\scrP$ is a rational polytope if and only if it is the convex hull of finitely many points with rational coordinates.
\item \label{characterize-ratione} $\scrP$ is a rational polyhedral cone if and only if it is the cone generated by finitely many points with rational coordinates.
\item \label{characterize-ratiohedron} $\scrP$ is a rational polyhedron if and only if it is the Minkowski sum of a rational convex polyhedral cone and a rational convex polytope.
\end{enumerate}
\end{prop}

\begin{proof}
At first assume $\scrP = \conv(S) + \cone(T)$ for $S, T \subset \qq^n$. We claim that $\scrP$ is a rational polyhedron. Due to \cref{exercise:rational-reduction} the claim reduces to the case that $\dim(\aff(\scrP)) = n$. By the results from \cref{poly-basic-section}, $\scrP$ is an $n$-dimensional polyhedron in $\rr^n$, and if $\scrQ$ is a facet of $\scrP$, then $\scrQ = \conv(S') + \cone(T')$ for $S' \subset S$ and $T' \subset T$. It is then straightforward to see that $\aff(\scrQ)$ is a rational affine hyperplane. Since an $n$-dimensional polyhedron in $\rr^n$ is determined by the hyperplanes containing its facets (\cref{prop:minimal-polynequalities}), it follows that $\scrP$ is a rational polyhedron, as required, and proves the $(\Leftarrow)$ implications of all three assertions of \cref{poly-rational-thm}. We now prove the $(\im)$ implications of all three assertions simultaneously by induction on $\dim(\scrP)$. Due to \cref{exercise:rational-reduction} we may assume \woutlog\ that $\dim(\scrP) = n$. The cases of $\dim(\scrP) = 0$ and $\dim(\scrP) = 1$ are then obvious. Now assume $\dim(\scrP) \geq 2$. Since the vertices and edges of $\scrP$ are also rational polyhedra (\cref{prop:minimal-polynequalities}), the inductive hypothesis and \cref{prop:vertex-hull,prop:edge-hull} imply that if $\scrP$ is a polytope, then $\scrP = \conv(S)$ for a finite $S \subset \qq^n$, and if $\scrP$ is a strongly convex polyhedral cone, then $\scrP = \cone(T)$ for a finite $T \subset \qq^n$. In general, \cref{polyhedron-minkomposition} implies that $\scrP = \scrP_0 + \scrC$, where $\scrP_0$ is a rational convex polytope and $\scrC$ is a rational convex polyhedron. Due to what we already proved, it suffices to show that $\scrC = \cone(T)$ with finite $T \subset \qq^n$, even if $\scrC$ is {\em not} strongly convex. Indeed, if $\scrC$ is not strongly convex, then \cref{strongly-convex-prop,polyhedron-minkomposition} imply that $\scrC = \scrC' + L$ where $\scrC'$ is a cone and $L$ is a {\em positive} dimensional rational linear subspace of $\rr^n$. Let $k := \dim(L)$. Then we can choose a surjective linear map $\pi: \rr^n \to \rr^{n-k}$ such that $\ker(\pi) = L$, and $\pi$ maps $\qq^n$ onto $\qq^{n-k}$. \Cref{exercise:cone-rational-projection} implies that $\pi(\scrC)$ is a rational polyhedral cone, and therefore, by the inductive hypothesis $\pi(\scrC) = \cone(T')$ for some finite $T' \subseteq \qq^{n-k}$. Another application of \cref{exercise:cone-rational-projection} then implies that $\scrC = \cone(T)$ for a finite $T \subseteq \qq^n$, as required.
\end{proof}

\begin{cor}
If $\scrP$ is a convex rational polytope, then each cone in the normal fan of $\scrP$ is also rational (with respect to the basis on $\rnstar$ which is dual to the standard basis of $\rr^n$).
\end{cor}

\begin{proof}
Let $\scrQ$ be a face of $\scrP$. We will show that the cone $\sigma_\scrQ$ from \cref{normal-section} is rational. \Cref{prop:minimal-polynequalities} implies that $\scrQ$ is also a convex rational polytope. Therefore by \cref{poly-rational-thm} there are finite subsets $S, T$ of $\qq^n$ such that each $\scrP$ is $\scrP = \conv(S)$ and $\scrQ = \conv(T)$. Since $\sigma_\scrQ = \{\nu \in \rnstar: \langle \nu, \beta - \beta' \rangle = 0$ for all $\beta, \beta' \in T$, and $\langle \nu, \alpha - \beta \rangle \geq 0$ for all $\alpha \in S$, $\beta \in T\}$, it follows that $\sigma_\scrQ$ is a rational polyhedral cone, as required.
\end{proof}

Recall that the \index{Hausdorff distance}{\em Hausdorff distance} of $\scrP, \scrQ \subset \rr^n$ is $\sup_{\alpha \in \scrP}\{\inf_{\beta \in \scrQ} \norm{\alpha - \beta}\}$, where $\norm{\cdot}$ is the Euclidean norm on $\rr^n$.

\begin{cor} \label{poly-rational-appxn}
Every polytope can be approximated arbitrarily closely (with respect to the Hausdorff distance) by rational polytopes.
\end{cor}

\begin{proof}
This follows from \cref{poly-basic-thm,poly-rational-thm}, since to approximate the convex hull of $\alpha_1, \ldots, \alpha_N \in \rr^n$, it suffices to take the convex hull of $\beta_1, \ldots, \beta_N$, where each $\beta_j$ has rational coordinates, and is sufficiently close to $\alpha_j$.
\end{proof}

Many results from the theory of rational polyhedra (including \cref{poly-rational-thm}) are ultimately based on the following basic fact from linear algebra.

\begin{lemma} \label{rational-solutions}
Consider a linear system of equations of the form $a_{i,0} + a_{i,1}x_1 + \cdots + a_{i,n}x_n  = 0$, $i = 1, \ldots, k$. Assume each $a_{i,j}\in \qq$. Let $H_\qq$ and $H_\rr$ be the set of common solutions of this system respectively in $\qq^n$ and $\rr^n$. Then $\dim_\qq(H_\qq) = \dim_\rr(H_\rr)$. In particular, the system has a common solution over $\rr$ if and only if it has a common solution over $\qq$.
\end{lemma}

\begin{proof}
If the system has a solution, then it can be found by Gaussian elimination, which produces a solution in $\qq$ (since each $a_{i,j} \in \qq$). Moreover, in that case the dimension of the space of solutions is precisely $n$ minus the rank of the $k \times n$ matrix $(a_{i,j})$, $1 \leq i \leq k$, $1 \leq j \leq n$. Since the rank over $\qq$ is the same as the rank over $\rr$, this completes the proof.
\end{proof}

\begin{cor}[Cf.\ \cref{prop:vertex-hull}] \label{positive-linear-comb-lemma}
Let $\scrP$ be a convex rational polyhedron and $\alpha$ be a point with rational coordinates in the relative interior of $\scrP$.
\begin{enumerate}
\item  \label{positive-linear-comb:pol} If $\scrP$ is a polytope with vertices $\alpha_1, \ldots, \alpha_k$, then there are positive rational numbers $q_1, \ldots, q_k$ such that $\sum_j q_j = 1$ and $\alpha = \sum_j  q_j\alpha_j$.
\item \label{positive-linear-comb:cone} If $\scrP$ is a cone generated by $\alpha_1, \ldots, \alpha_k \in \zz^n$, then there are positive rational numbers $q_1, \ldots, q_k$ such that $\alpha = \sum_j  q_j\alpha_j$.
\end{enumerate}
\end{cor}

\begin{proof}
Let $\phi: \rr^k \to \rr^n$ be the map which sends $(r_1, \ldots, r_k) \mapsto \sum_{j=1}^k r_j\alpha_j - \alpha$. If $\scrP$ is a cone, then \cref{prop:edge-hull}  implies that $\phi^{-1}(0) \cap \rpluss{k} \neq \emptyset$. \Cref{rational-solutions} then implies that $\phi^{-1}(0) \cap \qpluss{k} \neq \emptyset$, which together with \cref{prop:edge-hull} proves assertion \eqref{positive-linear-comb:cone}. Assertion  \eqref{positive-linear-comb:pol} follows  by the same arguments from \cref{prop:vertex-hull,rational-solutions} by considering the map $\phi': \rr^k \to \rr^{n+1}$ given by $(r_1, \ldots, r_k) \mapsto (\phi(r_1, \ldots, r_k), r_1 + \cdots + r_k - 1)$.
\end{proof}

The following result is one of the foundations of the theory of toric varieties we will encounter in \cref{toric-intro}.

\begin{lemma}[Gordan's lemma] \label{gordan}
\index{Gordan's lemma}
If $\sigma$ is a rational convex polyhedral cone in $\rr^n$, then $\sigma \cap \zz^n$ is a finitely generated semigroup.
\end{lemma}

\begin{proof}
It is straightforward to check that $\sigma \cap \zz^n$ is a semigroup. Pick $\alpha_1, \ldots, \alpha_s \in S_\sigma$ which generate $\sigma$ as a cone. Let $K := \{\sum_{i=1}^s t_i\alpha_i: t_i \in \rr, 0 \leq t_i \leq 1\} \subset \rr^n$. Since $K$ is compact, $K \cap \zz^n$ is a finite set. Now pick $\alpha \in \sigma \cap \zz^n$. Then $\alpha = \sum_{i=1}^s r_i \alpha_i$, where each $r_i \geq 0$. Write $\alpha = \sum_{i=1}^s \lfloor r_i \rfloor \alpha_i + \beta$, where $\lfloor r_i \rfloor$ is the greatest integer less than or equal to $r_i$ for each $i$, and $\beta := \sum_i (r_i - \lfloor r_i \rfloor) \alpha_i \in K \cap \zz^n$. It follows that $\sigma \cap \zz^n$ is generated as a semigroup by all the $\alpha_j$ together with $K \cap \zz^n$.
\end{proof}

\subsection{Exercises}
\begin{exercise} \label{exercise:rational-reduction}
Assume either $\scrP$ is a rational polyhedron or $\scrP = \conv(S) + \cone(T)$, where $S, T \subset \qq^n$. If $\dim(\aff(\scrP)) < n$, then show that there is a rational affine hyperplane in $\rr^n$ containing $\scrP$. [Hint: if $\scrP$ is a polyhedron in $\rr^n$ with $\dim(\scrP) < n$, then the set of points for which each of the inequalities defining $\scrP$ is strict must be empty, for otherwise $\dim(\scrP)$ would be $n$. In the other case, $\scrP = \conv(S) + \cone(T)$, with $S, T \subset \qq^n$. Let $L_\scrP$ be the linear subspace of $\rr^n$ generated by all elements of the form $\alpha - \beta$, where $\alpha - \beta \in \scrP$. It is possible to choose a generator of $L_\scrP$ consisting of elements of the form $\alpha_1 - \alpha_2 + \beta_1 - \beta_2$, where $\alpha_1, \alpha_2 \in S$ and $\beta_1, \beta_2$ are scalar multiples of elements from $T$. Since each such element is in $\qq^n$, it follows that $L_\scrP$ is contained in a rational hyperplane. $\scrP$ is contained in a translation of $L_\scrP$ by an element of $\qq^n$.]
\end{exercise}

\begin{exercise} \label{exercise:cone-rational-projection}
Let $\pi:\rr^n \to \rr^m$ be a linear map such that $\pi$ maps $\qq^n$ into $\qq^m$. Given $\scrC \subseteq \rr^m$, show that
\begin{enumerate}
\item $\scrC = \cone(T)$ for some (finite) $T \subset \qq^m$ if and only if $\pi^{-1}(\scrC) = \cone(T')$ for some (finite) $T' \subset \qq^n$.
\item $\scrC$ is a rational polyhedral cone in $\rr^m$ if and only if $\pi^{-1}(\scrC)$ is a rational polyhedral cone in $\rr^n$. [Hint: use \cref{exercise:cone-inequalities}.]
\end{enumerate}
\end{exercise}

\section{$^*$Volume of convex polytopes}
\label{rational-vol-section}
\footnote{The asterisk in the section name is to indicate that the material of this section is not going to be used in \cref{toric-intro}. It is only used in the proof of Bernstein's theorem in \cref{bkk-chapter}.}Given linearly independent elements $\alpha_1, \ldots, \alpha_d \in \rr^n$, the \index{Parallelotope}{\em parallelotope generated by the $\alpha_j$} is the set
\begin{align*}
\scrP := \{\lambda_1 \alpha_1 + \cdots + \lambda_d \alpha_d: 0 \leq \lambda_j \leq 1,\ 1 \leq j \leq d\}
\end{align*}
It is straightforward to check that $\scrP$ is a $d$-dimensional polytope in $\rr^n$ (\cref{exercise:parallelotope}). We denote by $\vol_n$ the usual $n$-dimensional volume (i.e.\ the Lebesgue measure) on $\rr^n$. Given an affine subspace $H$ of $\rr^n$, we write $\vol_H$ for the measure induced on $H$ by $\vol_n$. Here is a precise definition of $\vol_H$: let $H_0$ be the (unique) linear subspace of $\rr^n$ which is a translate of $H$ (i.e.\ $H_0 = H - \alpha$ for any $\alpha \in H$), $d := \dim(H)$, and $\beta_1, \ldots, \beta_{n-d}$ be an orthonormal set (with respect to the dot product) of elements in $\rr^n$ such that $H_0$ is precisely the set of elements in $\rr^n$ whose dot product with each $\beta_j$ is zero. Let $\scrQ$ be the parallelotope generated by $\beta_1, \ldots, \beta_d$. If $\scrR$ is any subset of $H$, then $\scrR$ is measurable with respect to $\vol_H$ if and only if $\scrQ +\scrR$ is measurable with respect to $\vol_n$, and $\vol_H(\scrR) = \vol_n(\scrQ + \scrR)$.

\begin{center}
\begin{figure}[h]

\def\scalefactor{.6}
\def\xmin{-2.5}
\def\xmax{7.5}
\def\ax{6}
\def\ay{-1}
\def\nux{-\ay}
\def\nuy{\ax}
\pgfmathsetmacro\nulength{sqrt(\nux*\nux + \nuy*\nuy)}
\def\nubasex{-1.5}

\begin{tikzpicture}[scale=\scalefactor]

\coordinate (leftL) at (\xmin, {\xmin*\ay/\ax});
\coordinate (rightL) at (\xmax, {\xmax*\ay/\ax});

\draw [\colortwo] (leftL) -- (rightL);

\coordinate (O) at (0,0);
\coordinate (A) at (\ax,\ay);
\coordinate (nu) at ({\nux/\nulength}, {\nuy/\nulength});
\coordinate (B) at ($(A) + (nu)$);

\fill [\colorzero, opacity=\opazero ] (O) --  (A) -- (B) -- (nu) -- cycle;
\draw[thick, \colorone] (O) -- (A);

\coordinate (nubase) at (\nubasex, {\nubasex*\ay/\ax});
\coordinate (nutip) at ($(nubase) + (nu)$);
\draw [thick, \colornu, ->] (nubase) -- (nutip);

\coordinate (nu2) at ({\ax/\nulength}, {\ay/\nulength});
\coordinate (nu2tip) at ($(nubase) + (nu2)$);
\draw [thick, \colornu, ->] (nubase) -- (nu2tip);

\coordinate (beta1) at ($(nubase)!0.5!(nutip)$);
\coordinate (beta2) at ($(nubase)!0.5!(nu2tip)$);
\node [anchor = east] at (beta1) {\picfontsize $\beta_1$};
\node [anchor = north] at (beta2) {\picfontsize $\beta_2$};

\coordinate (Q) at ($(O)!0.5!(A)$);
\node [anchor = north] at (Q) {\picfontsize $\scrQ$};

\node [anchor = west] at (rightL) {\picfontsize $H$};

\end{tikzpicture}

\caption{
If $\scrQ$ is a segment of a line $H \subset \rr^2$, then $\vol_H(\scrQ)$ is the length of $\scrQ$, and equals the area of the rectangle with base $\scrQ$ and height one.
} \label{fig:length=area}

\end{figure}
\end{center}
If $\scrP \subset \rr^n$ is such that the $\dim(\aff(\scrP)) \leq d$, then we write $\vol_d(\scrP)$ for $\vol_H(\scrP)$, where  $H$ is any $d$-dimensional affine subspace of $\rr^n$ containing $\scrP$; this is well defined since $\vol_H(\scrP)$ does {\em not} depend on $H$ (\cref{exercise:vol-d}). The following properties of $\vol$ follow from basic analysis; we use these without proof.

\begin{thm} \label{vol-basics}
Let $\scrP$ be a polytope in $\rr^n$.
\begin{enumerate}
\item If $\dim(\scrP) \leq d \leq n$, then the map $\lambda \in \rr_{\geq 0} \mapsto \vol_d(\lambda \scrP)$ is {\em homogeneous of order $d$}, i.e.\ $\vol_d(\lambda\scrP) = \lambda^d \vol_d(\scrP)$ for all $\lambda \geq 0$.
\item As a real valued function from the set of the polytopes in $\rr^n$, $\vol_n$ is continuous with respect to the Hausdorff distance.
\end{enumerate}
\end{thm}

Let $\scrP$ be an $n$-dimensional polytope, and $O$ be an arbitrary point of $\rr^n$. For each facet $\scrQ$ of $\scrP$, let $\scrS_{\scrQ,O}$ be the convex hull of $\scrQ \cup \{O\}$, i.e.\ $\scrS_{\scrQ,O}$ is the ``cone with base $\scrQ$ and apex $O$.'' Then it is intuitively clear that the volume of $\scrP$ can be computed in terms of the volumes of the $\scrS_{\scrQ,O}$, see \cref{fig:volume-facet}. This leads to our next result. To state it we introduce a notation: if $\scrQ$ is a facet of an $n$-dimensional polytope $\scrP \subset \rr^n$, then (up to a positive multiple) there is a unique $\nu \in \rnstar \setminus \{0\}$ such that $\scrQ = \In_\nu(\scrP)$. Define
\begin{align*}
\sign_\scrQ(\scrP, O)
	&=
	\begin{cases}
		1 & \text{if}\ \langle \nu, O \rangle > \min_\scrP(\nu),\\
		-1 & \text{otherwise.}
	\end{cases}
\end{align*}
Geometrically, $\sign_\scrQ(\scrP, O)$ is $1$ if and only if $O$ and $\scrP$ are on the same side of $\aff(\scrQ)$ - see \cref{fig:volume-facet}.

\begin{center}
\begin{figure}[h]
\def\xmin{-0.5}
\def\xmax{16.5}
\def\ymin{-5.5}
\def\ymax{5.5}

\tikzstyle{dot} = [\colordot, circle, minimum size=4pt, inner sep = 0pt, fill]
\begin{subfigure}[b]{0.36\textwidth}
\begin{tikzpicture}[scale=\scalefactor]
\coordinate (A) at (1, -4);
\coordinate (B) at (16, 1);
\coordinate (C) at (10, 5);
\coordinate (O) at (9.5, 0.5);
\coordinate (nu) at (-1,1);
\pgfmathsetmacro\nuscale{1}
\coordinate (acperpbase) at ($(A)!(O)!(C)$);
\coordinate (nutop) at ($(O) + \nuscale*(nu)$);

\draw[thick, \colorAB] (A) -- (B);
\draw[thick, \colorBC] (B) -- (C);
\draw[thick, \colorCA] (C) -- (A);
\fill[\colorAB, opacity=\opazero ] (O) --  (A) -- (B) -- cycle;
\fill[\colorBC, opacity=\opazero ] (O) --  (B) -- (C) -- cycle;
\fill[\colorCA, opacity=\opazero ] (O) --  (C) -- (A) -- cycle;

\draw[dashed, thick, \colorA] (O) -- (acperpbase);
\draw[thick, \colornu, ->] (O) -- (nutop);

\node[dot] at (O) {};
\node[dot] at (acperpbase) {};
\node[anchor = north] at (A) {\picfontsize $A$};
\node[anchor = west] at (B) {\picfontsize $B$};
\node[anchor = south] at (C) {\picfontsize $C$};
\node[anchor = north] at (O) {\picfontsize $O$};
\node[anchor = west] at (nutop) {\picfontsize $\nu_{CA}$};
\node[anchor = east] at (acperpbase) {\picfontsize $\alpha_{CA}$};
\end{tikzpicture}
\caption{$\vol_2(ABC) = \vol_2(OAB) + \vol_2(OBC) + \vol_2(OCA)$}
\end{subfigure}\hspace{0.18\textwidth}
\begin{subfigure}[b]{0.39\textwidth}
\begin{tikzpicture}[scale=\scalefactor]
\coordinate (A) at (1, -4);
\coordinate (B) at (16, 1);
\coordinate (C) at (10, 5);
\coordinate (O) at (3, 4);

\coordinate (nu) at (-1,1);
\pgfmathsetmacro\nuscale{1}
\coordinate (acperpbase) at ($(A)!(O)!(C)$);
\coordinate (nutop) at ($(O) + \nuscale*(nu)$);

\draw[thick, \colorAB] (A) -- (B);
\draw[thick, \colorBC] (B) -- (C);
\draw[thick, \colorCA] (C) -- (A);
\fill[\colorAB, opacity=\opazero ] (O) --  (A) -- (B) -- cycle;
\fill[\colorBC, opacity=\opazero ] (O) --  (B) -- (C) -- cycle;
\fill[\colorCA, opacity=\opazero ] (O) --  (C) -- (A) -- cycle;

\draw[dashed, thick, \colorA] (O) -- (acperpbase);
\draw[thick, \colornu, ->] (O) -- (nutop);

\node[dot] at (O) {};
\node[dot] at (acperpbase) {};
\node[anchor = north] at (A) {\picfontsize $A$};
\node[anchor = west] at (B) {\picfontsize $B$};
\node[anchor = south] at (C) {\picfontsize $C$};
\node[anchor = east] at (O) {\picfontsize $O$};
\node[anchor = east] at (nutop) {\picfontsize $\nu_{CA}$};
\node[anchor = west] at (acperpbase) {\picfontsize $\alpha_{CA}$};
\end{tikzpicture}
\caption{$\vol_2(ABC) = \vol_2(OAB) + \vol_2(OBC) - \vol_2(OCA)$}
\label{fig:volume-facet-2}
\end{subfigure}
\caption{If $\scrP$ is a triangle in $\rr^2$ and $O \in \rr^2$, then for each side $\scrQ$ of $\scrP$, $\scrS_{\scrQ,O}$ is the triangle formed by $\scrQ$ and $O$. In \cref{fig:volume-facet-2} $\sign_{CA}(ABC, O) = -1$ since $O$ and $ABC$ are on different sides of the line containing $CA$.}  \label{fig:volume-facet}
\end{figure}
\end{center}

\begin{thm} \label{volume-facet}
Let $\scrP$ be an $n$-dimensional polytope in $\rr^n$ and $O \in \rr^n$. For each facet $\scrQ$ of $\scrP$, let $d(O,\aff(\scrQ))$ denote the distance between $O$ and the affine hull of $\scrQ$. Then
\begin{align}
\vol_n(\scrP) = \frac{1}{n}\sum_\scrQ \sign_\scrQ(\scrP, O) d(O, \aff(\scrQ)) \vol_{n-1}(\scrQ)
\label{volume-facet-formula}
\end{align}
where the sum is over all facets $\scrQ$ of $\scrP$. In particular,
\begin{align}
\vol_n(\scrP)
	= \frac{1}{n}\sum_{\substack{\nu \in \rnstar\\ \norm{\nu} = 1}}
	 			\max_\scrP(\nu) \vol_{n-1}(\ld_\nu(\scrP))
	= -\frac{1}{n}\sum_{\substack{\nu \in \rnstar\\ \norm{\nu} = 1}}
 			\min_\scrP(\nu) \vol_{n-1}(\In_\nu(\scrP))
 	\label{volume-facet-unit-formula}
\end{align}

\end{thm}

\begin{rem}
The norm $\norm{\cdot}$ on $\rnstar$ in \eqref{volume-facet-unit-formula} is the Euclidean norm induced from $\rr^n$ upon identification of $\rnstar$ and $\rr^n$ via the basis dual to the standard basis of $\rr^n$.
\end{rem}

\begin{proof}
Let $\scrI_+$ (respectively, $\scrI_-$) be the collection of facets $\scrQ$ of $\scrP$ such that $\sign_\scrQ(\scrP,O) = 1$ (respectively, $-1$). The results of \cref{poly-basic-section} imply that
\begin{prooflist}
\item \label{volume-facet:P-in-cones} $\scrP \subset \bigcup_{\scrQ \in \scrI_+} \scrS_\scrQ(\scrP,O)$ (\cref{exercise:volume-facet:P-in-cones});
\item \label{volume-facet:minus-interior} if $\scrQ \in \scrI_-$ then $\scrS_{\scrQ,O} \cap \relint (\scrP) = \emptyset$ and $\scrS_{\scrQ,O} \subset \bigcup_{\scrQ' \in \scrI_+} \scrS_{\scrQ'}(\scrP,O)$ (\cref{exercise:volume-facet:minus-interior});
\item \label{volume-facet:plus-exterior} if $\scrQ \in \scrI_+$ then $\scrS_{\scrQ,O} \setminus \scrP \subset \bigcup_{\scrQ' \in \scrI_-} \scrS_{\scrQ'}(\scrP,O)$ (\cref{exercise:volume-facet:plus-exterior});
\item \label{volume-facet:intersect-cones} if $\scrQ$ and $\scrQ'$ are distinct facets of $\scrP$ such that $\sign_\scrQ(\scrP, O) =  \sign_{\scrQ'}(\scrP, O)$, then $\dim(\scrS_{\scrQ,O} \cap \scrS_{\scrQ',O}) \leq n-1$ (\cref{exercise:volume-facet:intersect-cones}).
\end{prooflist}
These observations immediately imply that
\begin{align*}
\vol_n(\scrP)
	&= \sum_\scrQ  \sign_\scrQ(\scrP, O) \vol_n(\scrS_{\scrQ,O})
\end{align*}
Since every cross section of $\scrS_{\scrQ,O}$ parallel to the hyperplane $\aff(\scrQ)$ is a dilation of $\scrQ$, it follows that
\begin{align*}
\vol_n(\scrS_{\scrQ,O})
	&= \int_{r=0}^{d(O,\aff(\scrQ))} \vol_{n-1}(r\scrQ)dr
	= \vol_{n-1}(\scrQ) \int_{r=0}^{d(O,\aff(\scrQ))}  r^{n-1} dr
	=\frac{1}{n} {d(O,\aff(\scrQ))} \vol_{n-1}(\scrQ),
\end{align*}
where the second equality follows from \cref{vol-basics}. This completes the proof of identity \eqref{volume-facet-formula}. Now for each facet $\scrQ$ of $\scrP$, let $\nu_\scrQ \in \rnstar$ be the {\em outward facing} unit normal to $\scrQ$, i.e.\ $\scrQ = \ld_\nu(\scrP)$ and $\norm{\nu} = 1$. We now apply identity \eqref{volume-facet-formula} with $O$ being the origin of $\rr^n$. If we identify $\rnstar$ with $\rr^n$ via the basis dual to the standard basis of $\rr^n$, then $\alpha_\scrQ := \sign_\scrQ(\scrP, O) d(O, \aff(\scrQ)) \nu_\scrQ $ is a point on $\aff(\scrQ)$ (see e.g.\ \cref{fig:volume-facet}). It follows that $\max_\scrP(\nu_\scrQ)
 = \langle \nu_\scrQ, \alpha_\scrQ  \rangle  = \sign_\scrQ(\scrP, O) d(O, \aff(\scrQ))$. The first equality of identity \eqref{volume-facet-unit-formula} now follows from identity \eqref{volume-facet-formula}. The second equality follows from the first by replacing $\nu$ by $-\nu$.
\end{proof}

\subsection{Exercises}
\begin{exercise} \label{exercise:parallelotope}
Let $\scrP := \{\lambda_1 \alpha_1 + \cdots + \lambda_d \alpha_d: 0 \leq \lambda_j \leq 1$ for each $j = 1, \ldots, d\}$, where $\alpha_1,  \ldots, \alpha_d \in \rr^n$. Show that
\begin{enumerate}
\item $\scrP$ is the convex hull of the set consisting of the origin and all elements of the form $\alpha_{i_1} + \cdots + \alpha_{i_k}$, where $1 \leq i_1 < \cdots < i_k \leq d$.
\item $\aff(\scrP)$ is the linear subspace of $\rr^n$ spanned by $\alpha_1, \ldots, \alpha_d$.
\end{enumerate}
\end{exercise}

\begin{exercise}\label{exercise:vol-d}
Let $\scrP \subset \rr^n$ and $k := \dim(\scrP)$. Fix $d$, $k \leq d \leq n$. If $H_1, H_2$ are $d$-dimensional affine subspaces of $\rr^n$ containing $\scrP$, show that
\begin{enumerate}
\item $\scrP$ is measurable with respect to $\vol_{H_1}$ if and only if it is measurable with respect to $\vol_{H_2}$.
\item If $\scrP$ is measurable with respect to either of them, then $\vol_{H_1}(\scrP) = \vol_{H_2}(\scrP)$.
\end{enumerate}
\end{exercise}

\begin{exercise} \label{exercise:volume-facet:P-in-cones}
Prove observation \ref{volume-facet:P-in-cones} from the proof of \cref{volume-facet}. [Hint: it suffices to prove that $\relint(\scrP) \subset \bigcup_{\scrQ \in \scrI_+} \scrS_\scrQ(\scrP,O)$. If $\alpha \in \relint(\scrP)\setminus \{O\}$, the line through $\alpha$ and $O$ intersects the facets at two points. One of these is contained in facets from $\scrI_+$.]
\end{exercise}

\begin{exercise} \label{exercise:volume-facet:minus-interior}
Prove observation \ref{volume-facet:minus-interior} from the proof of \cref{volume-facet}. [Hint: if $\scrI_- \neq \emptyset$, then either $O \not\in \scrP$ or $O$ is on a facet of $\scrP$. Pick $\alpha \in \scrQ$, where $\scrQ \in \scrI_-$. If $\scrQ = \In_\nu(\scrP)$, for the first part show that $\langle \nu, \beta \rangle \leq \min_\scrP(\nu)$ for each $\beta$ on the line segment from $O$ to $\alpha$. For the second part it suffices to show $\relint(\scrS_{\scrQ,O}) \subset \bigcup_{\scrQ' \in \scrI_+}$. Every point of $\relint(\scrS_{\scrQ,O})$ is on the line segment between $O$ and a point on $\relint(\scrQ)$. Extending this line segment hits another point of the topological boundary of $\scrP$ which belongs to a facet from $\scrI_+$.]
\end{exercise}

\begin{exercise} \label{exercise:volume-facet:plus-exterior}
Prove observation \ref{volume-facet:plus-exterior} from the proof of \cref{volume-facet}. [Hint: it suffices to consider the case that $O \not\in \scrP$. If $\alpha \in \relint(\scrQ)$ and $\scrQ \in \scrI_+$, then the line segment $L$ from $O$ to $\alpha$ intersects the boundary of $\scrP$ at a point $\beta$ ``in between'' $O$ and $\alpha$. Any facet $\scrR$ of $\scrP$ containing $\beta$ is in $\scrI_-$ and $L\setminus \scrP \subset \scrS_\scrR(\scrP,O)$.]
\end{exercise}

\begin{exercise} \label{exercise:volume-facet:intersect-cones}
Prove observation \ref{volume-facet:intersect-cones} from the proof of \cref{volume-facet}. [Hint: pick distinct facets $\scrQ_1, \scrQ_2$ of $\scrP$ such that $\scrS_{\scrQ_j,O}$ are full dimensional and $\relint(\scrS_{\scrQ_1,O}) \cap \relint(\scrS_{\scrQ_2,O}) \neq \emptyset$. It suffices to show that $\sign_{\scrQ_j}(\scrP, O)$ have different signs for $j=1$ and $j = 2$. Indeed, if $\alpha$ is a point in the intersection, then the line through $O$ and $\alpha$ intersects $\scrP$ at $\beta_j \in \relint(\scrQ_j)$, $j = 1, 2$. Show that one of the $\beta_j$ is ``in between'' $O$ and the other $\beta_j$.]
\end{exercise}

\begin{exercise} \label{exercise:omegadm}
Given positive $\omega_1, \ldots, \omega_n$ and nonnegative $d, m_1, \ldots, m_p$, where $0 \leq p \leq n$, let $\scrQ(\vec \omega, d, \vec m)$ be the polytope in $\rr^n$ determined by the following inequalities:
\begin{align*}
&x_i \geq 0\ \text{for each}\ i = 1, \ldots, n,\\
&\sum_{i=1}^n \omega_i x_i \leq d,\\
&x_i \leq m_i,\ i = 1, \ldots, p.
\end{align*}
Let $\scrI(\vec \omega, d, \vec m)$ be the collection of subsets $I$ of $\{1, \ldots, p\}$ such that $\sum_{i \in I} \omega_im_i \leq d$.
\begin{enumerate}
\item Show that the vertices of $\scrQ(\vec \omega, d, \vec m)$ are precisely the elements $\alpha_{I} = (\alpha_{I,1}, \ldots, \alpha_{I,n})$ and $\beta_{I,j} = (\beta_{I,j,1}, \ldots, \beta_{I,j,k}) \in \rr^n$, indexed by $I \in \scrI(\vec \omega, d, \vec m)$ and $j \in \{1, \ldots, n\} \setminus I$, and defined as follows:
\begin{align*}
\alpha_{I,k}
	&:=	\begin{cases}
			m_k & \text{if}\ k \in I,\\
			0	& \text{if}\ k \not\in I,
		\end{cases}\\
\beta_{I,j,k}
	&:=	\begin{cases}
			m_k
				& \text{if}\ k \in I,\\
			\frac{d - \sum_{i \in I}\omega_im_i}{\omega_j}
				& \text{if}\ k = j,\\
			0	& \text{if}\ k \not\in I \cup \{j\}.
		\end{cases}
\end{align*}
\item \label{sumQ} Let $d_i, m_{i,1}, \ldots, m_{i,p}$, $i = 1, 2$, be nonnegative real numbers such that $\scrI(\vec \omega, d_1, \vec m_1) = \scrI(\vec \omega, d_2, \vec m_2) = \scrI(\vec \omega, d_1 + d_2, \vec m_1 + \vec m_2)$. Then show that $\scrQ(\vec \omega, d_1+d_2, \vec m_1 + \vec m_2) = \scrQ(\vec \omega, d_1, \vec m_1) + \scrQ(\vec \omega, d_2, \vec m_2)$.
\item With $\vec \omega = (1,1,1)$, and $p = 2$, show that $\scrQ(\vec \omega, 3, (1, 1)) + \scrQ(\vec \omega, 3, (3, 3)) \subsetneq \scrQ(\vec \omega, 6, (4, 4))$, i.e.\ the conclusion of assertion \eqref{sumQ} may not be true in the absence of its assumption.
\end{enumerate}
\end{exercise}

\begin{exercise} \label{exercise:vol-omegadm}
\begin{enumerate}
\item Given polytopes $\Delta_1, \ldots, \Delta_p \subset \rr^n$, show that
\begin{align*}
\vol\left( \bigcup_{q=1}^p \Delta_q \right)
	&= \sum_{q=1}^p (-1)^{q-1} \sum_{I \subseteq [p],\ |I| = q}
		\vol\left(\bigcap_{i \in I} \Delta_i \right)
\end{align*}
where $[p]$ denotes the set $\{1, \ldots, p\}$ (this is the so called ``inclusion-exclusion principle,'' it holds for all ``measurable'' sets with finite volume, i.e.\ as long as the volume of each intersection on the right hand side is well defined and finite).

\item Given $\vec \omega, d, \vec m$ as in \cref{exercise:omegadm}, define
\begin{align*}
\Delta_0
	&:= \{(x_1, \ldots, x_n): x_i \geq 0,\ i = 1, \ldots, n,\ \sum_{i=1}^n \omega_i x_i \leq d\} \\
\Delta_q
	&:= \{(x_1, \ldots, x_n): x_i \geq 0,\ i = 1, \ldots, n,\ x_q \geq m_q, \sum_{i=1}^n \omega_i x_i \leq d\},\ q = 1, \ldots, p.
\end{align*}
Show that the polytope $\scrQ(\vec \omega, d, \vec m)$ from \cref{exercise:omegadm} equals $\Delta_0 \setminus \bigcup_{q=1}^p \Delta_p$. Conclude that
\begin{align*}
\vol(\scrQ(\vec \omega, d, \vec m))
	&=  \frac{1}{n!\omega_1 \cdots \omega_n} 	
	   \sum_{q = 0}^p (-1)^q
	   \sum_{\substack{I \subseteq [p],\ |I| = q\\ \sum_{i \in I} \omega_i m_i < d}}
	   (d - \sum_{i \in I}\omega_i m_i)^n
\end{align*}
\end{enumerate}
\end{exercise}

\section{$^*$Volume of special classes of polytopes}
\footnote{The asterisk in the section name is to indicate that the material of this section is not going to be used in \cref{toric-intro}. It is only used in the proof of Bernstein's theorem in \cref{bkk-chapter}.}In \cref{minkowski-section} we study the dependence of the volume of Minkowski sums of polytopes on its summands, and in \cref{rational-volume-section} we give a formula of the volume of rational polytopes in terms of ``lattice volumes'' of its facets. 

\subsection{Minkowski sums} \label{minkowski-section}
\Cref{vol-basics} implies that volume interacts well with Minkowski addition, in the sense that given compact convex subsets $\scrP, \scrQ$ of $\rr^n$, the function from $\rzero$ to $\rzero$ given by $\lambda \mapsto \vol_n(\scrP + \lambda \scrQ)$ is continuous. However, it turns out that this function is much more than a continuous function, it is a {\em polynomial}. In this section we are going to prove this result for the case of polytopes. At first we need the following result.

\begin{lemma} \label{sum-lemma-2}
Let $\scrP_1, \ldots, \scrP_s$ be subsets of $\rr^n$ and $\lambda := (\lambda_1, \ldots, \lambda_s) \in \rpluss{s}$. Then for different $\lambda$, the affine hull $A_\lambda$ of $\lambda_1\scrP_1 + \cdots + \lambda_s \scrP_s$ are translations of each other. In particular, $\dim(A_\lambda)$ is independent of $\lambda$.
\end{lemma}

\begin{proof}
Fix an arbitrary element $\alpha_i$ of $\scrP_i$, $i = 1, \ldots, n$. \Woutlog\ we may replace $\scrP_i$ by $\scrP_i - \alpha_i$ and assume that each $\scrP_i$ contains the origin. For each $\lambda \in \rpluss{s}$, it then follows that $A_\lambda$ contains each $\scrP_i$, and therefore it is simply the linear subspace of $\rr^n$ spanned by elements in $\bigcup_i \scrP_i$.
\end{proof}

\begin{thm} \label{volume-thm}
Let $\scrP_1, \ldots, \scrP_s$ be convex polytopes in $\rr^n$. Then there are nonnegative real numbers $v_\alpha(\scrP_1, \ldots, \scrP_s)$ for all $\alpha \in \scrE_s:= \{(\alpha_1, \ldots, \alpha_s) \in \zzeroo{s}: \alpha_1 + \cdots  + \alpha_s = n\}$ such that for all $\lambda = (\lambda_1, \ldots, \lambda_s) \in \rzeroo{s}$,
\begin{align*}
\vol_n(\lambda_1\scrP_1 + \cdots + \lambda_s \scrP_s) = \sum_{\alpha \in \scrE_s} v_\alpha(\scrP_1, \ldots, \scrP_s) \lambda_1^{\alpha_1} \cdots \lambda_s^{\alpha_s}
\end{align*}
where $\vol_n$ is the $n$-dimensional volume.
\end{thm}

\begin{proof}
We proceed by induction on $n$. If $n = 1$ each $\scrP_i$ is of the form $[a_i, b_i]$, so that
\begin{align*}
\vol_1(\lambda_1\scrP_1 + \cdots + \lambda_s \scrP_s)
	&= \vol_1([\lambda_1a_1 + \cdots + \lambda_s a_s, \lambda_1b_1 + \cdots + \lambda_s b_s]) \\
	&= \lambda_1(b_1 - a_1) + \cdots + \lambda_s(b_s - a_s) \\
	&= \sum_i \lambda_i\vol_1(\scrP_i)
\end{align*}
Now assume it is true for convex polytopes in $\rr^{n-1}$. Pick convex polytopes $\scrP_1, \ldots, \scrP_n$ in $\rr^n$. Since the volume is translation invariant, we may assume that
\begin{align}
\parbox{.6\textwidth}{the origin is in the relative interior of each $\scrP_j$.}
\tag{$*$} \label{origin-interior}
\end{align}
Let $P_\lambda := \lambda_1\scrP_1 + \cdots + \lambda_s \scrP_s$. Due to \cref{vol-basics} it suffices to consider the case that each $\lambda_i$is {\em positive}. If $\dim(P_\lambda) \leq n - 1$, then due to \cref{sum-lemma-2} the result is true with all $\nu_\alpha$ being zero. So assume $\dim(P_\lambda) = n$. Then \cref{prop:sum-faces,sum-lemma-2} imply that the number of facets of $\scrP_\lambda$ does not depend on $\lambda$, and moreover, if $\scrP_{\lambda, 1}, \ldots, \scrP_{\lambda, N}$ are the facets of $\scrP_\lambda$, then for each $j$, there are faces $\scrP_{i,j}$ of $\scrP_i$, $i = 1, \ldots, s$, such that
\begin{align*}
\scrP_{\lambda, j} = \lambda_1 \scrP_{1, j} + \cdots + \lambda_s \scrP_{s,j}
\end{align*}
For each $i,j$, pick an arbitrary $\alpha_{i, j} \in \scrP_{i,j}$. Let $\nu_j$ be the outward pointing unit normal to $\scrQ_{\lambda, j}$. Identity \eqref{volume-facet-unit-formula} implies that
\begin{align*}
\vol_n(\scrP_\lambda)
	&= \frac{1}{n}\sum_j
 			\max_{\scrP_\lambda}(\nu_j) \vol_{n-1}(\scrP_{\lambda, j})\\
 	&= \frac{1}{n}\sum_j
 	 			\langle \nu_j, \lambda_1 \alpha_{1, j} + \cdots + \lambda_s \alpha_{s,j}\rangle \vol_{n-1}(\lambda_1 \scrP_{1, j} + \cdots + \lambda_s \scrP_{s,j})
\end{align*}
Condition \eqref{origin-interior} implies that  $\langle \nu_j, \alpha_{i,j} \rangle$ is nonnegative for each $i,j$. Since for each $j$, all the $\scrP_{i,j}$ can be identified (via a volume preserving affine map from $\aff(\scrP_{\lambda,j})$ to $\rr^{n-1}$) with polytopes in $\rr^{n-1}$, the result then follows from the inductive hypothesis.
\end{proof}

\subsection{Rational polytopes} \label{rational-volume-section}

%

Let $H$ be a $d$-dimensional rational affine subspace of $\rr^n$. If $\beta \in H \cap \zz^n$, \cref{rational-solutions} implies that $G_H := (H - \beta) \cap  \zz^n$ is isomorphic (as an abelian group) to $\zz^d$. A \index{Fundamental lattice parallelotope}\index{Parallelotope!fundamental lattice}{\em fundamental lattice parallelotope} on $H$ is a polytope of the form $\scrP + \beta$, where $\scrP$ is a ($d$-dimensional) parallelotope generated by $d$ elements from $G_H$ which generate $G_H$ as an abelian group. We write $\fund(H) := \vol_H(\scrP)$. \Cref{prop:fund-vol} below shows that $\fund(H)$ is well defined. In this section we identify $\rr^n$ with $\rnstar$ via the dot product, and given $\alpha, \beta \in \rr^n$, write $\langle \alpha, \beta \rangle$ for the dot product of $\alpha$ and $\beta$. Similarly we write $\beta^\perp := \{\gamma \in \rr^n: \langle \beta, \gamma \rangle = 0\}$.

\begin{center}
\begin{figure}[h]
\def\xmin{-3.5}
\def\xmax{5.5}
\def\ymin{-3.5}
\def\ymax{3.5}
\def\scalefactor{0.45}
\tikzstyle{dotsmall} = [\colordot, circle, minimum size=3pt, inner sep = 0pt, fill]

\begin{tikzpicture}[scale=\scalefactor]
\draw [gray,  line width=0pt] (\xmin, \ymin) grid (\xmax,\ymax);
\draw [<->] (0, \ymax) |- (\xmax, 0);

\coordinate (A') at (-3.5, -2.5);
\coordinate (B') at (2.5, 3.5);
\coordinate (A) at (-2, -1);
\coordinate (B) at (1, 2);

\coordinate (C') at (-3.5, -3.25);
\coordinate (D') at (5.5, 1.25);
\coordinate (C) at (-1, -2);
\coordinate (D) at (3, 0);

\draw (A') -- (B');
\draw (C') -- (D');
\draw[thick, \colorAB] (A) -- (B);
\draw[thick, \colorBC] (C) -- (D);

\node[dotsmall] at (A) {};
\node[dotsmall] at (B) {};
\node[dotsmall] at (C) {};
\node[dotsmall] at (D) {};

\node[anchor = east] at (A) {\picfontsize $A$};
\node[anchor = east] at (B) {\picfontsize $B$};
\node[anchor = north] at (C) {\picfontsize $C$};
\node[anchor = north] at (D) {\picfontsize $D$};
\node[anchor = south] at (B') {\picfontsize $L$};
\node[anchor = west] at (D') {\picfontsize $M$};
\end{tikzpicture}

\caption{$\fund(L)= \sqrt{2}$, $\fund(M) = \sqrt{5}$. $\vol'_L(AB) = 3$, $\vol'_M(CD) = 2$.}
\label{fig:fundamentally-normal-vol}
\end{figure}
\end{center}

\begin{prop} \label{prop:fund-vol}
Let $H$ be a rational affine subspace of $\rr^n$. If $\scrP_1, \scrP_2$ are two fundamental lattice parallelotopes of $H$, then $\vol_H(\scrP_1) = \vol_H(\scrP_2)$.
\end{prop}

\begin{proof}
By translating $H$ and the $\scrP_j$ if necessary we may assume $H$ is a linear subspace of $\rr^n$, and each $\scrP_j$ is the parallelotope generated by $\alpha_{j,1}, \ldots, \alpha_{j, d} \in H \cap \zz^n$, where $d := \dim(H)$. Pick an orthonormal set $\beta_1, \ldots, \beta_{n-d} \in \rr^n$ such that $H = \bigcap_i \beta_i^\perp$. For each $j = 1, 2$, let $\scrB_j$ be the basis of $\rr^n$ consisting of $\beta_1, \ldots, \beta_{n-d}, \alpha_{j,1}, \ldots, \alpha_{j,d}$. By definition $\vol_H(\scrP_j) = \vol_n(\scrQ_j)$, where $\scrQ_j$ is the parallelotope generated by the elements of $\scrB_j$. Let $\phi: \rr^n \to \rr^n$ be the linear map which changes coordinates with respect to $\scrB_1$ to that of $\scrB_2$. Since $\alpha_{j,1}, \ldots, \alpha_{j,d}$ generate $H \cap \zz^n$ (as an abelian group), it follows that the matrices of both $\phi$ and $\phi^{-1}$ have only integer entries. This means that the determinant of the matrix of $\phi$ is $\pm 1$, and therefore $\phi$ preserves $\vol_n$. Since $\phi$ maps one of the $\scrQ_j$ to the other, this completes the proof.
\end{proof}

Let $H'$ be a rational affine subspace of $\rr^n$ such that $H' \supset H$ and $\dim(H') = d+1$. We now describe the relation between $\fund(H)$ and $\fund(H')$. Pick $\beta \in H$. Since $(H- \beta) \cap \qq^n \subset (H'- \beta) \cap \qq^n$ is an inclusion of vector spaces over $\qq$, it follows from the elementary theory of vector spaces that there is $\eta' \in (H'- \beta) \cap \qq^n$ such that $H - \beta = \eta'^\perp \cap (H' - \beta)$. Pick $r \in \qq \setminus \{0\}$ such that $\eta := r\eta'$ is ``primitive,'' i.e.\ there is no integer $k > 1$ such that $\eta = k\eta''$ for some $\eta'' \in \zz^n$. Let $\norm{\cdot}$ denote the Euclidean norm on $\rr^n$.

\begin{prop} \label{fundamental-cor'}
\begin{align*}
\fund(H) =  \frac{\norm{\eta}}{\min\{|\langle \eta, \alpha\rangle|: \alpha \in ((H' - \beta) \cap \zz^n) \setminus (H - \beta)\}}\fund(H')
\end{align*}
In particular, if $d = n-1$ and $H' = \rr^n$, then $\fund(H) = \norm{\eta}$.
\end{prop}

\begin{proof}
Replacing $H$ by $H - \beta$ if necessary we may assume that $\beta= 0$. \Cref{basis-lemma} implies that we may pick $\alpha_1, \ldots, \alpha_{d+1}, u_{d+1}, \ldots, u_n \in \rr^n$ such that
\begin{itemize}
\item $\alpha_1, \ldots, \alpha_d$ generate $H \cap \zz^n$,
\item $\alpha_1, \ldots, \alpha_{d+1}$ generate $H' \cap \zz^n$,
\item $u_{d+1} \in H'$, $H = u_{d+1}^\perp \cap H'$, $\norm{u_{d+1}} = 1$,
\item $u_{d+2}, \ldots, u_n$ are orthonormal (with respect to dot product), and $H' = \bigcap_{j = d+2}^n u_j^\perp$.
\end{itemize}
Let $M$ be the matrix with column vectors $\alpha_1, \ldots, \alpha_d, u_{d+1}, \ldots, u_n$ and $M'$ be the matrix with column vectors $\alpha_1, \ldots, \alpha_{d+1}, u_{d+2}, \ldots, u_n$. Write $u_{d+1} = c_{d+1}\alpha_{d+1} + u'$, where $u' \in H$. Then
\begin{align*}
\fund(H) = |\det(M)| = |c_{d+1} \det(M')| =| c_{d+1}|\fund(H')
\end{align*}
Note that $u_{d+1} = \pm \eta/\norm{\eta}$. It follows that $\langle \eta, u_{d+1} \rangle = \pm \norm{\eta}$. On the other hand, $\langle \eta, u_{d+1} \rangle = \langle \eta, c_{d+1}\alpha_{d+1} + u' \rangle = c_{d+1} \langle \eta, \alpha_{d+1} \rangle$. Since for each $\alpha \in H' \cap \zz^n$, $\langle  \eta, \alpha \rangle$ is an integer multiple of $ \langle \eta, \alpha_{d+1} \rangle$, the result follows.
\end{proof}

\begin{defn} \label{vol'}
Let $H$ be a rational affine subspace of $\rr^n$ of dimension $d$. The {\em $H$-normalized volume} is
\begin{align*}
\vol'_{H}(\cdot) := \vol_d(\cdot)/\fund(H)
\end{align*}
see \cref{fig:fundamentally-normal-vol} for some examples with $n =2$. An \index{Integral!element of $\rr^n$}{\em integral} element of $\rr^n$ is an element with integral coordinates; an integral element of $\rnstar$ is one which has integral coordinates with respect to the basis which is dual to the standard basis of $\rr^n$. An integral element $\eta$ of $\rr^n$ or $\rnstar$ is \index{Primitive!integral element of $\rr^n$}{\em primitive}, if it is not of the form $k\eta'$, where $k > 1$ and $\eta'$ is also integral. If $\nu$ is an integral element of $\rnstar$, then we write $\vol'_\nu$ for $\vol'_{\nu^\perp}$.
\end{defn}

\begin{cor} \label{volume-facet-rational}
Let $\scrP$ be a convex rational polytope in $\rr^n$. Then
\begin{align}
\vol_n(\scrP)
	&= \frac{1}{n}\sum_\nu \max_\scrP(\nu) \vol'_\nu(\ld_\nu(\scrP))
	= -\frac{1}{n}\sum_\nu \min_\scrP(\nu) \vol'_\nu(\In_\nu(\scrP)) \label{rational-volume}
\end{align}
where the sum is over all primitive integral $\nu  \in \rnstar$.
\end{cor}

\begin{proof}
Since $\scrP$ is rational, every facet of $\scrP$ is determined by integral elements of $\rnstar$. Therefore the result follows from combining \cref{volume-facet-unit-formula,fundamental-cor'}.
\end{proof}

\begin{rem} \label{primitive-remark}
If $\scrQ$ is a facet of an $n$-dimensional rational polytope $\scrP$ in $\rr^n$, then the \index{Primitive!inner normal}\index{Primitive!outer normal}{\em primitive inner} (respectively {\em outer}) {\em normal} to $\scrF$ is the unique primitive integral $\nu \in \rnstar$ such that $\scrQ = \In_\nu(\scrP)$ (respectively $\scrQ = \ld_\nu(\scrP)$). Note that the sum in \eqref{rational-volume} is practically finite: the only non-zero contributions come from those $\nu$ which are primitive outer normal to $(n-1)$-dimensional faces of $\scrP$.
\end{rem}

\part{Number of zeroes on the torus} \label{toric-part}
\chapter{Toric varieties over algebraically closed fields} \label{toric-intro}
\def\scalefactor{.3}

This chapter introduces {\em toric varieties}, which are the setting of all the subsequent results of this book. Our treatment will be mostly based on the results from \cref{var-chapter,appolytopes}; only in \cref{xone-section} we use the notion of {\em closed subschemes} discussed in \cref{subscheme-section}. Unless explicitly stated otherwise, from this chapter onward $\kk$ denotes an algebraically closed field (of arbitrary characteristic), and $\kk^*$ denotes $\kk \setminus \{0\}$.

\section{Algebraic torus}
%

If $(x_1, \ldots, x_n)$ are coordinates on $\kk^n$, $n \geq 1$, then $\nktorus = \kk^n \setminus V(x_1 \cdots x_n)$. This implies that $\nktorus$ is an affine variety, and its coordinate ring is the ring $\kk[x_1, x_1^{-1}, \ldots, x_n, x_n^{-1}]$ of \index{Laurent polynomial}{\em Laurent polynomials} in $x_1, \ldots, x_n$ over $\kk$ (\cref{example:affine-complement}). An \index{Algebraic!torus}\index{Torus}{\em (algebraic) torus} is a variety $X$ isomorphic to $\nktorus$ for some $n \geq 1$. A \index{System of coordinates on a torus}{\em system of coordinates} on $X$ is an ordered collection $(x_1, \ldots, x_n)$ of regular functions on $X$ such that the coordinate ring $\kk[X]$ of $X$ is $\kk[x_1, x_1^{-1}, \ldots, x_n, x_n^{-1}]$. A basic property of a torus is that {\em every} morphism between two tori is a group homomorphism and also a {\em monomial} map with respect to {\em every} set of coordinates. Indeed, let $\phi: X \to Y \cong \nktoruss{N}$ be a morphism. Choose coordinates $(y_1, \ldots, y_N)$ on $Y$. If $\phi(x) = (\phi_1(x), \ldots, \phi_N(x))$, then each $\phi_j$ must be a monomial in $(x_1, \ldots, x_n)$, for otherwise it will be zero at some point of $X$ (\cref{exercise:vanishing-monomial}). This shows that $\phi$ is a monomial map. Write $\phi_j(x) = x^{\alpha_j}$, $\alpha_1, \ldots, \alpha_N \in \zz^n$. If $x = (x_1, \ldots, x_n), x' = (x'_1, \ldots, x'_n) \in X$, then it follows that $\phi(x \cdot x') = \phi(x_1x'_1, \ldots, x_nx'_n) = ((x \cdot x')^{\alpha_1}, \ldots, (x \cdot x')^{\alpha_N}) = \phi(x) \cdot \phi(x')$, so that $\phi$ is indeed a homomorphism. This implies in particular that the multiplication with respect to every set of (algebraic) coordinates on a torus induces the same group structure on it, and the image of a morphism between two tori is a subgroup of the target. \Cref{torus-image} below shows that it is in addition a Zariski closed subset of the target\footnote{Contrast this to the case of $\kk^n$: the additive group structures on $\kk^n$ with respect to different systems of algebraic coordinates are in general different, and the image of a morphism from $\kk^n$ to $\kk^n$ is in general neither a subgroup nor a closed subvariety of the target (see \cref{example:constructible-0,constructible-section}).}. We use the following notation in \cref{torus-image}: given a monomial map $\phi:x \mapsto (x^{\alpha_1}, \ldots, x^{\alpha_N})$ between tori with fixed systems of coordinates, we write $\mat\phi$ for the $N \times n$ matrix whose rows are the $\alpha_i$. Some basic properties of $\mat{\phi}$ are established in \cref{exercise:mat-properties}.

\begin{prop} \label{torus-image}
 Let $G$ be the subgroup of $\zz^n$ generated by the $\alpha_i$, and $\bar G := \{\alpha \in \zz^n: k\alpha \in G$ for some $k \geq 1\}$ be the ``saturation'' of $G$ in $\zz^n$. Let $q$ be the index of $G$ in $\bar G$. Let $r$ be the rank of $\mat{\phi}$ as a matrix over $\qq$. Then
\begin{enumerate}
\item $\phi(X)$ is a torus and a closed subvariety of $Y$ of dimension $r$.
\item $\bar G/G \cong \prod_{j=1}^r \zz/q_j\zz$ for positive integers $q_1, \ldots, q_r$ such that $q = \prod_j q_j$.
\item $\ker(\phi)$ is an $(n-r)$-dimensional subgroup of $X$ isomorphic to $(\bar G/G) \times \nktoruss{n-r}$. In particular, if $r=n$, then the degree of $\phi$ (as a map from $X$ to $\phi(X)$) is $q$.
\item Pick a basis $(\beta_1, \ldots, \beta_{n-r})$ of $\ker \mat{\phi} \subseteq \zz^n$, and let $\eta: \nktoruss{n-r} \to \nktorus$ be the morphism such that the column vectors of $[\eta]$ are the $\beta_j$. Then the irreducible component of $\ker(\phi)$ containing $(1, \ldots, 1)$ is the image of $\eta$.
\end{enumerate}
\end{prop}

\begin{proof}
Let $\phi': \zz^n \to \zz^N$ be the map corresponding to multiplication by $\mat{\phi}$. With respect to appropriate coordinates on $\zz^n$ and $\zz^N$, the matrix of $\phi'$ is of the form
\begin{align*}
\left[
	\begin{array}{c|c}
	D & 0\\
	\hline
	0 & 0
	\end{array}
\right]
\end{align*}
where $D$ is a diagonal matrix with positive integers as diagonal entries (\cref{change-of-basis-lemma}). This means that we can choose coordinates $(x_1, \ldots, x_n)$ on $X$ and $(y_1, \ldots, y_N)$ on $Y$ with respect to which $\phi$ takes the form $(x_1, \ldots, x_n) \mapsto (x_1^{q_1}, \ldots, x_r^{q_r}, 1, \ldots, 1)$ for positive integers $q_1, \ldots, q_r$. All the assertions are now straightforward; their proofs are left as \cref{exercise:torus-image}.
\end{proof}

\subsection{Exercises}
\begin{exercise} \label{exercise:vanishing-monomial}
Show that every polynomial in $(x_1, \ldots, x_n)$ which is not a monomial vanishes at some points on $\nktorus$. [Hint: use \cref{exercise:nonempty-hypersurface-aff}.]
\end{exercise}

\begin{exercise} \label{exercise:mat-properties}
Let $\phi: X \to Y$ be a morphism between two tori. Let $(x_1, \ldots, x_n)$ (respectively, $(y_1, \ldots, y_N)$) be coordinates on $X$ (respectively, $Y$) and $\mat{\phi}$ be the corresponding matrix of $\phi$. Show that
\begin{enumerate}
\item  for each $\beta = (\beta_1, \ldots, \beta_N) \in \zz^N$, $(\phi(x))^\beta = x^{\beta \mat{\phi}}$, where $\beta \mat{\phi}$ is the product of $\beta$ (regarded as a $1 \times N$ matrix) and $[\phi]$;
\item $\phi$ is an isomorphism if and only if $N = n$ and $\mat\phi$ is invertible over $\zz$;
\item \label{mat:composition=prod} if $Z$ is a torus and $\psi: Y \to Z$ is a morphism, then $\mat{\psi \circ \phi} = \mat{\psi}\mat{\phi}$.
\end{enumerate}
\end{exercise}

\begin{exercise} \label{exercise:torus-image}
Complete the proof of \cref{torus-image}.
\end{exercise}

\section{Toric varieties from finite subsets of $\zz^n$} \label{xa-section}
A \index{Toric variety!from a finite subset of $\zz^n$}\index{Variety!toric}{\em toric variety} is a variety $X$ which contains an algebraic torus as a dense open subset such that the (multiplicative) action of the torus on itself extends to an action on all of $X$. Given a finite subset $\scrA = \{\alpha_0, \ldots, \alpha_N\}$ of $\zz^n$, we write $\phi_\scrA: \nktorus \to \pp^N$ for the map given by
\begin{align}
x \mapsto [x^{\alpha_0}: \cdots: x^{\alpha_N}] \label{phi_A}
\end{align}
We write $\xzeroa$ for the image of $\phi_\scrA$ and $\xa$ for the closure of $\xzeroa$ in $\pp^N$. We will now show that $\xa$ is a toric variety with torus $\xzeroa$. Denote the homogeneous coordinates of $\pp^N$ by $[z_{\alpha_0}: \cdots : z_{\alpha_N}]$. Let $U_\alpha := \pp^N \setminus V(z_\alpha)$, $\alpha \in \scrA$, be the basic open subsets of $\pp^N$.

\begin{prop} \label{xa-toric}
$\xa$ is a toric variety. More precisely,
\begin{enumerate}
\item \label{xa:torus} $\xzeroa$ is a torus and $\xzeroa = \xa \cap \bigcap_{\alpha \in \scrA} U_\alpha$.
\item \label{xa:U-alpha} For each $\alpha \in A$, $\xa \cap U_\alpha$ is an affine variety with coordinate ring $\kk[x^\beta: \beta \in S_\alpha]$, where $S_\alpha$ is the subsemigroup of $\zz^n$ generated by $\scrA - \alpha := \{\beta - \alpha: \beta \in \scrA\}$.
\item \label{xa:dimension} The dimension of $\xa$ (and equivalently, of $\xzeroa$) equals the dimension (as a polytope) of the convex hull of $\scrA$ in $\rr^n$.
\item \label{xa:action} $\xzeroa$ acts on $\xa$ via the multiplicative action of $\xzeroa$ on $\pp^N$ given by:
\begin{align}
[y_{\alpha_0}: \cdots: y_{\alpha_N}] \cdot [z_{\alpha_0}: \cdots: z_{\alpha_N}]:= [y_{\alpha_0}z_{\alpha_0}: \cdots: y_{\alpha_N}z_{\alpha_N}] \label{xa-action}
\end{align}
for all $[y_{\alpha_0}: \cdots: y_{\alpha_N}] \in \xzeroa$ and $[z_{\alpha_0}: \cdots: z_{\alpha_N}]\in \pp^N$.
\end{enumerate}
\end{prop}

\begin{proof}
Since $\xzeroa \subset U := \bigcap_{\alpha \in \scrA} U_\alpha$, and since $U \cong \nktoruss{N}$ via the map $[z_{\alpha_0}: \cdots : z_{\alpha_N}] \mapsto (z_{\alpha_1}/z_{\alpha_0}, \ldots, z_{\alpha_N}/z_{\alpha_0})$, it follows that $\xzeroa$ is the image in $\nktoruss{N}$ of the map $x \mapsto (x^{\alpha_1-\alpha_0}, \ldots, x^{\alpha_n - \alpha_0})$. \Cref{torus-image} then implies that $\xzeroa$ is a torus, and also implies assertion \eqref{xa:dimension}. It also says that $\xzeroa$ is a closed subset of $U$, which implies that $\xzeroa = \xa \cap U$, and proves assertion \eqref{xa:torus}. Assertion \eqref{xa:U-alpha} follows directly from \cref{projective-embedding}. Finally, since for a fixed $y \in \xzeroa$, the action of $y$ on $\pp^N$ given by \eqref{xa-action} is an isomorphism (\cref{exercise:continuous-action}) and since $\xzeroa$ is closed under this action, it follows that $\xa$ is also closed under it (\cref{exercise:action-closure}), as required to prove assertion \eqref{xa:action}.
\end{proof}

\Cref{xa-toric} states in particular that $\xa$ is ``equivariantly embedded'' in $\pp^N$, i.e.\ the action of the torus on $\xa$ extends to all of $\pp^N$. Conversely every equivariantly embedded projective toric variety is essentially of the form $\xa$ for some appropriate $\scrA$ (see e.g.\ \cite[Proposition 5.1.5]{gkz}); we will not use this result. We now show that $\xzeroa$ and $\xa$ depend only on the affine geometry of the set $\scrA$.

\begin{prop}[{\cite[Proposition 5.1.2]{gkz}}] \label{xa-invariant}
Let $\scrA \subset \zz^n$, $\scrB \subset \zz^m$, and $T: \zz^n \to \zz^m$ be an injective integer affine transformation such that $T(\scrA) = \scrB$. Then $\xzeroa = \xzerob$ and $\xa = \xb$ as subsets of $\pp^N$, where $N = |\scrA| - 1$.
\end{prop}

\begin{proof}
Let $A = \{\alpha_0, \ldots, \alpha_N\}$, $\scrB = \{\beta_0, \ldots, \beta_N\}$, where $\beta_j = T(\alpha_j)$, $j = 0, \ldots, N$. By definition there is $\lambda = (\lambda_1, \ldots, \lambda_m) \in \zz^m$ and an $n \times  m$ matrix $M$ such that for each $\gamma = (\gamma_1, \ldots, \gamma_n) \in \zz^n$, $T(\gamma) = \lambda + \gamma M$. Let $\mu_j \in \zz^m$ be the $j$-th row vector of $M$, $j = 1, \ldots, n$, and $T^*: \nktoruss{m} \to \nktorus$ be the map defined by $T^*(y) = (y^{\mu_1}, \ldots, y^{\mu_n})$. \Cref{exercise:mat-properties} implies that $T^*(y)^{\alpha_j} = y^{\alpha_j M} = y^{\beta_j - \lambda}$ for each $j = 0, \ldots, N$. Since the rank of $M = \mat{T^*}$ is $n$, \cref{torus-image} implies that $T^*$ is surjective, and therefore
\begin{align*}
\xzerob
	&= \{[y^{\beta_0}: \cdots : y^{\beta_N}]: y \in \nktoruss{m}\}
	= \{[y^{-\lambda}y^{\beta_0}: \cdots : y^{-\lambda}y^{\beta_N}]: y \in \nktoruss{m}\}\\
	&= \{[T^*(y)^{\alpha_0}: \cdots : T^*(y)^{\alpha_N}]: y \in \nktoruss{m}\}
	= \{[x^{\alpha_0}: \cdots : x^{\alpha_N}]: x \in \nktorus\}
	= \xzeroa
\end{align*}
which completes the proof.
\end{proof}

\subsection{Exercises}

\begin{exercise} \label{exercise:continuous-action}
Show that for a fixed $y \in \xzeroa$, the action of $y$ on $\pp^N$ given by \eqref{xa-action} is an isomorphism from $\pp^N$ to itself.
\end{exercise}

\begin{exercise} \label{exercise:action-closure}
Let $W$ be a subset of a topological space $X$ and $\phi: X \to X$ be a continuous map. Show that $\phi(\overline{W}) \subseteq \overline{\phi(W)}$ (where the ``bar'' indicates closure in $X$). Deduce that if $\phi(W) \subseteq W$, then $\phi(\overline{W}) \subseteq \overline{W}$. 
\end{exercise}

\begin{exercise} \label{exercise:xa-binomial}
In this exercise you will show that $\xa$ is a {\em binomial variety}, i.e.\ $\xa$ is defined in $\pp^N$ by binomial equations.  Given $q = (q_0, \ldots, q_N) \in \zz^N$, write $z^q := \prod_{j=0}^N z_{\alpha_j}^{q_j}$. Let $J$ be the ideal of $R := \kk[z_{\alpha_0}, \ldots, z_{\alpha_N}]$ generated by binomials of the form $z^{q_1} - z^{q_2}$, where $q_i := (q_{i,0}, \ldots, q_{i,N}) \in \zz^n$, $i = 1, 2$, are such that $\sum_{j=0}^N q_{1,j} = \sum_{j=0}^N q_{2,j}$ and $\sum_{j=0}^N q_{1,j} \alpha_j = \sum_{j=0}^N q_{2,j} \alpha_j$.  Let $I$ be the ideal of $R$ consisting of all homogeneous polynomials that vanish on $\xa$.
\begin{enumerate}
\item Show that $J \subset I$.
\item Let $f = \sum_{q \in \zzeroo{N}} c_q z^q \in I$. Use the fact that $f(x^{\alpha_0}, \ldots, x^{\alpha_N})$ is identically zero on $\nktorus$ to show that for each $\alpha \in \zz^n$, $\sum_{q \in \scrS_\alpha}c_q = 0$, where $\scrS_\alpha := \{(q_0, \ldots, q_N) \in \zzeroo{N}:  \sum_{j=0}^N q_j\alpha_j = \alpha\}$.
\item Given $\alpha \in \zz^n$, write $f_\alpha := \sum_{q \in \scrS_\alpha}c_qz^q$. Use the preceding step to show that $f_\alpha \in J$. Deduce that $I = J$.
\end{enumerate}
\end{exercise}

\section{Examples of toric varieties}
\label{toric-exaction}
Given a subset $Y$ of a topological space $X$, in this section we write $\cl_X(Y)$ for the closure of $Y$ in $X$. We also write $e_1, \ldots, e_n$ for the standard unit vectors of $\rr^n$, i.e.\ for each $i,j$, the $j$-th coordinate of $e_i$ is $0$ if $j \neq i$, and $1$ if $j = i$.

\begin{example} \label{example:p^n}
If $\scrA = \{0, 1\} \subset \zz$, then $\xa = \cl_{\pp^2}(\{[1: x]: x \in \kk^*\} = \pp^1$. More generally, if $\scrA = \{0, e_1, \ldots, e_n\}$, then $\xa = \cl_{\pp^n}(\{[1:x_1: \cdots : x_n]: (x_1, \ldots, x_n) \in \nktorus\}) = \pp^n$.
\end{example}

\begin{example} \label{example:cube}
If $\scrA = \{0,1\}^2 = \{(0,0), (1,0), (0,1), (1,1)\}$, then $\xa = \cl_{\pp^3}(\{[1:x:y:xy]: x,y \in \kk^*\}) =  V(z_1z_2 - z_0z_3) \subset \pp^3$ (the last equality follows by a dimension count). Recall (from \cref{example:ruled}) that $\xa$ is the ruled surface isomorphic to $\pp^1 \times \pp^1$. More generally, if $\scrA = \{0,1\}^n$, assertion \eqref{toric-product} of \cref{segre-veronese} below implies that $\xa \cong (\pp^1)^n$.
\end{example}

Given a set $\scrA \subset \rr^n$ and a positive integer $d$, in this chapter we write $d\scrA := \{\alpha_1 + \cdots + \alpha_d: \alpha_j \in \scrA$ for each $j\} \subset \rr^n$. If $\scrA,\scrB$ are finite subsets of $\zz^n$, we write $\delta_{\scrA, \scrB}$ for the {\em diagonal} map from $\kk^n$ to $\pp^{|\scrA| -1} \times \pp^{|\scrB| - 1}$ given by
\begin{align}
x \mapsto (\phi_{\scrA}(x), \phi_{\scrB}(x))
\label{delta_AB}
\end{align}
where $\phi_{\scrA}$ and $\phi_{\scrB}$ are defined as in \eqref{phi_A}.

\begin{prop} \label{segre-veronese}
Let $\scrA_1, \scrA_2$ be finite subsets of $\zz^{n_i}$ and $d_1, d_2$ be positive integers.
\begin{enumerate}
\item \label{toric-veronese} $\xaa{\scrA_1} \cong \xaa{d_1\scrA_1}$.
\item \label{toric-product} $\xaa{\scrA_1} \times \xaa{\scrA_2} \cong \xaa{d_1\scrA_1 \times d_2\scrA_2}$.
\item \label{toric-diagonal} Assume $n_1 = n_2 = n$. Let $X_{d_1, d_2}$ be the closure in $\xaa{d_1\scrA_1} \times \xaa{d_2\scrA_2}$ of $\delta_{d_1\scrA_1,d_2\scrA_2}(\nktorus)$. Then $\xaa{\scrA_1 + \scrA_2} \cong \xaa{d_1\scrA_1+d_2\scrA_2} \cong X_{d_1, d_2} \cong X_{1,1}$.
\qed
\end{enumerate}
\end{prop}

\begin{proof}
We may assume \woutlog\ that each $\scrA_i$ contains the origin. If $\nu_{d_i}$ is the degree-$d_i$ {\em Veronese map} (see \cref{veronese-subsection}), then it follows that $\phi_{d_i\scrA_i} = \pi \circ \nu_{d_i} \circ \phi_{\scrA_i}$, where $\pi$ is a projection which omits ``redundant'' coordinates of $\nu_{d_i} \circ \phi_{\scrA_i}$ (i.e.\ for each $\alpha \in d_i\scrA_i$, $\pi$ retains only one of the coordinates of $\nu_{d_i} \circ \phi_{\scrA_i}$ equalling $x^\alpha$). Assertion \eqref{toric-veronese} then follows from \cref{veronese-embedding}. For assertion \eqref{toric-product}, let $\phi_{\scrA_1} \times \phi_{\scrA_2}: \nktoruss{n_1+n_2} \to \xaa{\scrA_1} \times \xaa{\scrA_2}$ be the morphism which maps $(x_1,x_2) \mapsto (\phi_{\scrA_1}(x_1), \phi_{\scrA_2}(x_2))$, where $x_i \in \nktoruss{n_i}$, $i = 1, 2$. If $s: \pp^{|\scrA_1|-1} \times \pp^{|\scrA_2|-1} \to \pp^{|\scrA_1||\scrA_2|-1}$ is the Segre map, then assertion \eqref{toric-product} follows from assertion \eqref{toric-veronese} and the observation that $s \circ (\phi_{\scrA_1} \times \phi_{\scrA_2}) = \phi_{\scrA_1 \times \scrA_2}$. The proof of assertion \eqref{toric-diagonal} is left as \cref{exercise:toric-diagonal}.
\end{proof}

\begin{example} \label{example:rational-normal-curve}
Let $d$ be a positive integer. If $\scrA = \{0,1,2, \ldots, d\} \subset \zz$, then assertion \eqref{toric-veronese} of \cref{segre-veronese} and \cref{example:p^n} imply that $\xa \cong \pp^1$. Note that $\xa$ is the {\em rational normal curve of degree $d$} from \cref{example:rational-normal}. 
\end{example}

\begin{example} \label{example:blow-up-toric}
Let $\scrA_1 = \{0, e_1, \ldots, e_n\}$ and $\scrA_2 = \{e_1, \ldots, e_n\}$, so that $\scrA_1 + \scrA_2 = \{e_1, \ldots, e_n\} \cup \{e_i + e_j\}_{i,j}$. Since $\xaa{\scrA_1} \cong \pp^n$ and $\xaa{\scrA_2} \cong \pp^{n-1}$ (\cref{example:p^n}), it follows from \cref{segre-veronese} that $\xaa{\scrA_1 + \scrA_2}$ is the closure in $\pp^n \times \pp^{n-1}$ of the image of the map from $\nktorus \to \pp^n \times \pp^{n-1}$ given by 
\begin{align*}
(x_1, \ldots, x_n) \mapsto ([1: x_1: \cdots : x_n], [x_1: \cdots : x_n])
\end{align*}
It follows that $\xaa{\scrA_1 + \scrA_2}$ is precisely the {\em blow up} of $\pp^n$ at the point $[1: 0: \cdots: 0]$ (see \cref{example:projective-blow-up}). \Cref{fig:proj-blow-up} shows the convex hull of $\scrB := \scrA_1 + \scrA_2$ for $n = 2, 3$. 
\end{example}

\begin{center}
\begin{figure}[htb]
\def\scalefactor{1}
\def\xmin{-0.5}
\def\xmax{2.5}
\def\ymin{-0.5}
\def\ymax{2.5}
\def\opazero{0.5}
\def\colorzero{green}
\def\colorg{orange}

\begin{subfigure}[b]{0.2\textwidth}
\begin{tikzpicture}[scale=\scalefactor]
\draw [gray,  line width=0pt] (\xmin, \ymin) grid (\xmax,\ymax);
\draw [<->] (0, \ymax) |- (\xmax, 0);

\coordinate (A) at (1,0);
\coordinate (B) at (0,1);
\coordinate (C) at (0,2);
\coordinate (D) at (2,0);

\fill[green, opacity=\opazero ] (A) --  (B) -- (C) -- (D) -- cycle;
\draw[thick] (A) -- (B) -- (C) -- (D) -- cycle;
\end{tikzpicture}
\caption{$n = 2$}  
\end{subfigure}\hspace{0.1\textwidth}
\def\viewx{60}%
\def\viewy{30}%
\begin{subfigure}[b]{0.2\textwidth}
\def\viewx{60}
\def\viewy{30}
\begin{tikzpicture}[scale=0.66]
\pgfplotsset{every axis title/.append style={at={(0,-0.2)}}, view={\viewx}{\viewy}, axis lines=middle, enlargelimits={upper}}

\begin{axis}
\addplot3[ultra thick, draw, fill=\colorg,opacity=\opazero] coordinates{(0,0,1) (0,1,0) (1,0,0)};
\addplot3[ultra thick, draw, fill=\colorg,opacity=\opazero] coordinates{(0,0,1) (0,1,0) (0,2,0) (0,0,2)};
\addplot3[ultra thick, draw, fill=green,opacity=\opazero] coordinates{(0,0,2) (0,2,0) (2,0,0)};
\addplot3[ultra thick, draw, fill=red,opacity=\opazero] coordinates{(0,0,2) (2,0,0) (1,0,0) (0,0,1)};
\end{axis}
\end{tikzpicture}
\caption{$n = 3$} 
\end{subfigure}
\caption{Convex hull of $\scrB$ such that $\xb$ is the blow up of $\pp^n$ at a point}  \label{fig:proj-blow-up}
\end{figure}
\end{center}

\begin{example}
If $\scrA = \{0,2, 3\} \subset \zz$, then $\xa = \cl_{\pp^2}(\{[1:x^2:x^3]: x \in \kk^*\}) = V(z_0z_2^2 - z_1^3) \subset \pp^2$. If $(x,y) = (z_1/z_0, z_2/z_0)$ are coordinates on the basic open subset $U_0 = \pp^2 \setminus V(z_0)$ of $\pp^2$, then $\xa \cap U_0$ is the curve defined by $y^2 =x^3$. Note that the curve is singular at the origin (\cref{example:curve-singularities}).
\end{example}

A toric variety of dimension $n$ by definition contains an open subset isomorphic to $\nktorus$. The following result desribes a class of sets $\scrA$ such that $\xa$ contains open subsets isomorphic to $\kk^n$, i.e.\ $\xa$ is a {\em compactification} of $\kk^n$ - its proof is left as \cref{exercise:toric-k^n}.

\begin{prop} \label{toric-k^n}
Let $\scrA$ be a finite subset of $\zz^n$.
\begin{enumerate}
\item \label{toric-k^n:origin} Assume $\scrA \subset \znzero$ and $\scrA \supset \{\origin, e_1, \ldots, e_n\}$. Show that $U_\origin \cap \xa \cong \kk^n$.
\item More generally, assume there is a vertex $\alpha$ of $\scrA$ such that there are precisely $n$ edges $\scrE_1, \ldots, \scrE_n$ of $\conv(\scrA)$ containing $\alpha$ and on each $\scrE_i$, there is an element of $\scrA$ of the form $\alpha + \beta_i$, where $\beta_1, \ldots, \beta_n$ is a basis of $\zz^n$. Then show that $U_\alpha \cap \xa \cong \kk^n$. \qed
\end{enumerate}
\end{prop}

\subsection{Exercises}
\begin{exercise} \label{exercise:toric-diagonal}
Prove assertion \eqref{toric-diagonal} of \cref{segre-veronese}.
\end{exercise}

\begin{exercise} \label{exercise:toric-k^n}
Prove \cref{toric-k^n}. [Hint: the first assertion follows from \cref{xa-toric}. For the second assertion and the definition of ``vertex'' read the next section and use \cref{xa-thm}.]
\end{exercise}

\section{Structure of $\xa$}
We continue with the notation and the set up of \cref{xa-section}. In this section we study the complement of $\xzeroa$ in $\xa$ (i.e.\ the subvariety of $\xa$ ``at infinity'') and catch a glimpse of its beautiful combinatorial structure. Let $\scrP$ be the convex hull of $\scrA$. Then $\scrP$ is a rational polytope (\cref{poly-rational-thm}). A \index{Face!of a finite set}{\em face} of $\scrA$ is by definition a set $\scrB$ of the form $\scrQ \cap \scrA$ where $\scrQ$ is a face of $\scrP$. We say that $\scrB$ is a \index{Facet!of a finite set}\index{Vertex!of a finite set}\index{Edge!of a finite set}{\em facet} (respectively {\em vertex, edge}) of $\scrA$ if $\scrQ$ is a facet (respectively vertex, edge) of $\scrP$, see \cref{fig:face-set}.

\begin{center}
\begin{figure}[h]
\def\xmin{-0.5}
\def\xmax{13.5}
\def\ymin{-4.5}
\def\ymax{5.5}
\def\tx{0}
\def\ty{5}
\def\tw{5cm}
\def\opazero{0.5}
\def\colorzero{green}
\def\colorone{blue}

\tikzstyle{dot} = [red, circle, minimum size=4pt, inner sep = 0pt, fill]

\begin{subfigure}[b]{0.45\textwidth}
\begin{tikzpicture}[scale=\scalefactor]
\draw [gray,  line width=0pt] (\xmin, \ymin) grid (\xmax,\ymax);
\draw [<->] (0, \ymax) |- (\xmax, 0);

\node[dot] (A) at (1,-4) {};
\node[dot] (B) at (13,2) {};
\node[dot] (C) at (10,5) {};
\draw[thick, \colorone, fill=\colorzero, opacity=\opazero ] (A.center) --  (B.center) -- (C.center) -- cycle;
\node[dot] (D) at (10,2) {};
\node[dot] (E) at (7,-1) {};

\node[anchor = north] at (A) {\picfontsize $A$};
\node[anchor = west] at (B) {\picfontsize $B$};
\node[anchor = south] at (C) {\picfontsize $C$};
\node[anchor = east] at (D) {\picfontsize $D$};
\node[anchor = north] at (E) {\picfontsize $E$};

\node [below right, text width= \tw, align=left] at (\tx,\ty) {
	\picfontsize
};	
\end{tikzpicture}
\caption{
$\scrA = \{A, B, C, D, E\}$\\
Facets: $\{A, E, B\}$, $\{B, C\}$, $\{C, A\}$\\
Vertices: $A$, $B$, $C$
} \label{fig:A}
\end{subfigure}	
\begin{subfigure}[b]{0.45\textwidth}
\def\xmin{-3.5}
\def\xmax{7.5}
\def\tx{-8}
\begin{tikzpicture}[scale=\scalefactor]
\draw [gray,  line width=0pt] (\xmin, \ymin) grid (\xmax,\ymax);
\draw [<->] (0, \ymax) |- (\xmax, 0);

\def\singularPtwo#1#2#3{
\begin{scope}[shift={#1}]
\node[dot] (A) at (0,0) {};
\node[dot] (B) at (#2,0) {};
\node[dot] (C) at (0,#3) {};
\draw[thick, \colorone, fill=\colorzero, opacity=\opazero ] (A.center) --  (B.center) -- (C.center) -- cycle;
\node[dot] (D) at (0,1) {};
\node[dot] (E) at (#2-1,0) {};

\node[anchor = north] at (A) {\picfontsize $A'$};
\node[anchor = north] at (B) {\picfontsize $B'$};
\node[anchor = south] at (C) {\picfontsize $C'$};
\node[anchor = east] at (D) {\picfontsize $D'$};
\node[anchor = north] at (E) {\picfontsize $E'$};
\end{scope}
}
\singularPtwo{(-3,-4)}{9}{9}

\end{tikzpicture}
\caption{
	$\scrA' = \{A', B', C', D', E'\}$ \\
	Facets: $\{A', E', B'\}$, $\{B', C'\}$, $\{C', D', A'\}$\\
	Vertices: $A'$, $B'$, $C'$
} \label{fig:A'}
\end{subfigure}	

\caption{Faces of some planar sets}  \label{fig:face-set}

\end{figure}
\end{center}

\begin{prop} \label{point-face}
Let $z \in \xa$. Define $\scrA_z := \{\alpha \in \scrA: z \in U_\alpha\}$. Then $\scrA_z$ is a face of $\scrA$.
 \end{prop}

\begin{proof}
Pick $\alpha_1, \alpha_2 \in \scrA_z$. For each $i$, let $\scrQ_i$ be the (unique) face of $\scrP$ which contains $\alpha_i$ in its interior. For every $\epsilon \in (0,1)$, the convex combination $\alpha_\epsilon := \epsilon \alpha_1 + (1- \epsilon)\alpha_2$ of $\alpha_1$ and $\alpha_2$ belongs to the relative interior of a face $\scrQ$ of $\scrP$ containing both $\scrQ_i$ (\cref{prop:positively-relative}). It suffices to show that $\scrQ \cap \scrA \subset \scrA_z$. Indeed, let $\epsilon$ be a rational number in $(0,1)$. Let $\beta_1, \ldots, \beta_k$ be the vertices of $\scrQ$. There are positive integers $N_1, \ldots, N_k$ such that $(\sum_{j=1}^k N_j)\alpha_\epsilon = \sum_{j=1}^k N_j\beta_j$ (\cref{positive-linear-comb-lemma}). Write $N := \sum_{j=1}^k N_j$. Multiplying the $N_j$ by some appropriate integer we may ensure that $N\epsilon$ is a positive integer. Then $\prod_{j = 1}^k z_{\beta_j}^{N_j} = z_{\alpha_1}^{N\epsilon}z_{\alpha_2}^{N-N\epsilon}$ on $\xa$. Since both $\alpha_i$ are in $\scrA_z$, it follows that $z_{\beta_j}|_z \neq 0$ for each $j$. Since each $\beta \in \scrQ \cap \scrA$ is a convex rational linear combination of the $\beta_j$, it follows by the same reasoning that $z_\beta|_z \neq 0$ for each $\beta \in \scrQ \cap \scrA$, as required.
\end{proof}

Given $\scrB' \subseteq \scrB \subseteq \scrA$, we write $\scrB' \preceq \scrB$ (respectively, $\scrB' \precneqq \scrB$) to denote that $\scrB'$ is a face (respectively, a proper face) of $\scrB$.

\begin{thm} \label{xa-thm}
For each face $\scrB$ of $\scrA$, define $\corbit{\scrB} := \xa \setminus \bigcup_{\alpha \not\in \scrB} U_\alpha$ and $\orbit{\scrB} := \corbit{\scrB} \cap \bigcap_{\beta \in \scrB} U_\beta$.
\begin{enumerate}
\item \label{union-of-orbits} $\corbit{\scrB} = \bigcup_{\scrB' \preceq \scrB} \orbit{\scrB'}$. In particular,
\begin{align}
\xa\setminus \xzeroa
	= \bigcup_{\scrB \precneqq \scrA} \orbit{\scrB}
	= \bigcup_{\scrB \precneqq \scrA} \corbit{\scrB}
	= \bigcup_{j=1}^k \corbit{\scrA_j}
\label{eq:union-of-orbits}
\end{align}
where $\scrA_1, \ldots, \scrA_k$ are the facets of $\scrA$.
\item \label{orbits} There is a one-to-one correspondence between the collection of $\orbit{\scrB}$ for $\scrB \preceq \scrA$ and the set of orbits of $\xzeroa$ on $\xa$. In particular, each $\corbit{\scrB}$ is invariant under the action of $\xzeroa$.
\item \label{toric-faces} Each $\corbit{\scrB}$ is a toric variety with torus $\orbit{\scrB}$. More precisely, the pair $(\corbit{\scrB}, \orbit{\scrB})$ is isomorphic to $(\xaa{\scrB}, \xzeroaa{\scrB})$. The isomorphism is given by the projection map $\pi_\scrB: V_\scrB\to \pp^{|\scrB|-1}$ which ``drops'' all the coordinates $z_\alpha$ such that $\alpha \not\in \scrB$; in other words, $\pi_\scrB([z_\alpha: \alpha \in \scrA]) = [z_\beta: \beta \in \scrB]$.
\item \label{compatible-action} The action of $\orbit{\scrB}$ on $\corbit{\scrB}$ is compatible with the action of $\xzeroa$. More precisely, assume $y_\scrA \in \xzeroa$ and $y_\scrB \in \orbit{\scrB} \cong \xzerob$ correspond to the same $x \in \nktorus$, i.e.\ $y_\scrA = [x^{\alpha}: \alpha \in \scrA]$ and $y_\scrB = [x^{\beta}: \beta \in \scrB]$. Then for all $z \in \corbit{\scrB}$,
\begin{align*}
y_\scrA \cdot_\scrA z = y_\scrB \cdot_\scrB z
\end{align*}
%
where we write $\cdot_\scrA$ (respectively, $\cdot_\scrB$) to denote the action of $\xzeroa$ (respectively, $\orbit{\scrB}$) on $\corbit{\scrB}$.
\end{enumerate}
\end{thm}

\begin{proof}
The first statement of assertion \eqref{union-of-orbits} follows from \cref{point-face}, and this in turn implies the first two equalities of \eqref{eq:union-of-orbits}. Since $\corbit{\scrB'} \subseteq \corbit{\scrB}$ whenever $\scrB' \preceq \scrB$, and since every proper face is contained in a facet (\cref{prop:facet-containment}), the last equality of \eqref{eq:union-of-orbits} follows. We now prove the remaining assertions. Given $\scrB \preceq \scrA$, let $H_\scrB := V(z_\alpha: \alpha \not\in \scrB) \subset \pp^N$ be the coordinate subspace containing $\corbit{\scrB}$. Let $\pi_\scrB: \pp^N \setminus V(z_\beta: \beta \in \scrB) \to H_\scrB$ be the natural projection and $\zb$ be the closure in $\pp^N$ of $\pi_\scrB(\xzeroa)$.

\begin{proclaim} \label{corbit-claim}
$\corbit{\scrB} = \zb$.
\end{proclaim}

\begin{proof}
The inclusion $\corbit{\scrB} \subseteq \zb$ follows from a general property of morphisms (assertion \eqref{subspace-intersection} of \cref{projective-embedding}). For the opposite inclusion it suffices (due to the definition of $\corbit{\scrB}$) to show that $\xa \supset \zb$. Pick $\beta \in \scrB$. We will show that $\xa \cap U_\beta \supseteq \zb \cap U_\beta$. We may assume $\beta = \alpha_0$. Write $z'_i := z_{\alpha_i}/z_{\alpha_0}$, so that $U_\beta \cong \kk^N$ with coordinates $(z'_1, \ldots, z'_N)$. Let $f(z'_1, \ldots, z'_N) = \sum_{\gamma \in \zz^N} c'_\gamma z'^\gamma$ be a polynomial in $(z'_1, \ldots, z'_N)$ which vanishes on $\xzeroa \cap U_\beta$. It suffices to show that $f$ vanishes on $\zb \cap U_\beta$ as well.  Note that $\xzeroa \cap U_\beta = \{(x^{\alpha'_1}, \ldots, x^{\alpha'_N}): x \in \nktorus\}$, where $\alpha'_i := \alpha_i - \alpha_0$, $i = 1, \ldots, N$. Write $f = f' + f''$, where the monomials in $f'$ consist solely of the $z'_i$ such that $\alpha_i \in \scrB$ and each monomial in $f''$ contains at least one $z'_i$ such that $\alpha_i \not\in \scrB$. If $B$ is the {\em affine hull} of $\{\alpha_i - \alpha_0: \alpha_i \in \scrB\}$, \cref{exercise:face-semi-affine-hull} implies that the exponent of each monomial in $f'(x^{\alpha'_1}, \ldots, x^{\alpha'_N})$ is on $B$, whereas the exponent of no monomial in $f''(x^{\alpha'_1}, \ldots, x^{\alpha'_N})$ is on $B$; in particular, the monomials that appear in $f'(x^{\alpha'_1}, \ldots, x^{\alpha'_N})$ are distinct from those appearing in $f''(x^{\alpha'_1}, \ldots, x^{\alpha'_N})$. Since $f(x^{\alpha'_1}, \ldots, x^{\alpha'_N})$ is identically zero on $\nktorus$, it follows that  $f'(x^{\alpha'_1}, \ldots, x^{\alpha'_N})$ is also identically zero on $\nktorus$. This implies that $f \circ \pi_\scrB(x^{\alpha'_1}, \ldots, x^{\alpha'_N}) = f'(x^{\alpha'_1}, \ldots, x^{\alpha'_N}) = 0$ for all $x \in \nktorus$, so that $f$ vanishes on $\pi_\scrB(\xzeroa)$, as required.
\end{proof}

It is evident that $\zb$ can be identified to $\xb$ by ``forgetting'' the coordinates $z_\alpha$ for all $\alpha \not\in \scrB$. Since $\corbit{\scrB} = \zb$, \cref{xa-toric} implies that this induces an identification of $\xzerob$ with $\orbit{\scrB}$, which proves assertion \eqref{toric-faces}. Assertion \eqref{compatible-action} then follows from identity \eqref{xa-action}. Since assertion \eqref{compatible-action} in particular implies that $\orbit{\scrB}$ is an orbit of $\xzeroa$, this proves assertion \eqref{orbits} as well.
\end{proof}

\Cref{cartier-infty} below is an immediate corollary of \cref{xa-thm} - we leave its proof as \cref{exercise:cartier-infty}. The second statement of \cref{cartier-infty} in particular implies that the complement of the torus is locally defined by a single equation on $\xa$; in the terminology of \cref{scheme:cartier-section}, $\xa\setminus \xzeroa$ is the ``support of a Cartier divisor'' on $\xa$.

\begin{cor} \label{cartier-infty}
If $\scrV$ is the set of vertices of $\scrA$, then $\xa \subset \bigcup_{\alpha \in \scrV} U_\alpha$. In particular, $\xa \setminus \xzeroa = V(\prod_{\alpha \in \scrV} z_\alpha) \cap \xa$. \qed
\end{cor}

\subsection{Exercises}

\begin{exercise} \label{exercise:face-semi-affine-hull}
Let $\scrQ$ be a face of a polyhedron $\scrP \subset \rr^n$. Assume that the origin is in $\aff(\scrQ)$. Let $\alpha = \sum_{j= 1}^k r_j \alpha_j$, where $r_j > 0$ and $\alpha_j \in \scrP$ for each $j$. Show that $\alpha \in \aff(\scrQ)$ if and only if $\alpha_j \in \scrQ$ for each $j$. [Hint: if $\scrQ =\In_\nu(\scrP)$, consider the value of $\langle \nu, \cdot \rangle$ on $\scrQ$, $\scrP$ and $\alpha$.]
\end{exercise}

\begin{exercise} \label{exercise:cartier-infty}
Prove \cref{cartier-infty}. [Hint: use \cref{prop:face-of-face}.]
\end{exercise}

\section{Toric varieties from polytopes} \label{toric-polysection}

%

If $\scrA \subset \scrA'$ are finite subsets of $\zz^n$, then the natural projection $\pp^{|\scrA'|-1} \dashrightarrow \pp^{|\scrA|-1}$ restricts to a rational map $\pi_{\scrA, \scrA'}:\xaa{\scrA'} \dashrightarrow \xa$. This map is in general not defined everywhere on $\xaa{\scrA'}$. For example, if $\scrA = \{(0,0), (1,0)\}$ and $\scrA' = \{(0,0), (1,0), (0,1)\}$, then \cref{example:p^n} implies that $\xaa{\scrA'} \cong \pp^2$, $\xa \cong \pp^1$, and $\pi_{\scrA,\scrA'}$ maps $[x_0:x_1:x_2] \mapsto [x_0:x_1]$, which is not defined at the point $[0:0:1]$. However, if in addition $\scrA$ and $\scrA'$ have the same convex hull in $\rr^n$, then \cref{cartier-infty} implies that $\pi_{\scrA,\scrA'}$ is well-defined everywhere on $\xaa{\scrA'}$ (\cref{exercise:equi-hull-projection}). If $\scrP$ is the convex hull of $\scrA$ in $\rr^n$, this observation shows that there is a natural morphism $X_{\scrP \cap \zz^n} \to \xa$, and for every positive integer $k$ there is a natural morphism
\begin{align*}
\xaa{(k+1)\scrP \cap \zz^n} \to \xaa{(k\scrP \cap \zz^n) + (\scrP \cap \zz^n)}
\end{align*}
(For subsets $\scrS$ of $\rr^n$ and a positive integer $d$, in \cref{appolytopes} we defined $d\scrS$ as a ``dilation,'' whereas in \cref{toric-exaction} we defined it as the sum of $d$-copies of $\scrS$. This does not lead to any conflict for the case of convex polytopes - see \cref{equivalent-multiplication}.) \Cref{segre-veronese} implies that $\xaa{(k\scrP \cap \zz^n) + (\scrP \cap \zz^n)}$ is isomorphic to a subset of $\xaa{k\scrP \cap \zz^n} \times \xaa{\scrP \cap \zz^n}$, so that the projection onto the first factor induces a morphism $\xaa{(k\scrP \cap \zz^n) + (\scrP \cap \zz^n)} \to \xaa{k\scrP \cap \zz^n}$. Consequently, there is a sequence of morphisms
\begin{align}
 \cdots
 	\xrightarrow{\pi_{3,4}} \xaa{3\scrP \cap \zz^n}
 	\xrightarrow{\pi_{2,3}} \xaa{2\scrP \cap \zz^n} \xrightarrow{\pi_{1,2}} \xaa{\scrP \cap \zz^n}
 	\xrightarrow{\pi_1} \xa
 \label{xp-to-xa}
\end{align}

\begin{prop} \label{xprop}
For $k$ sufficiently large, $\pi_{k,k+1}:\xaa{(k+1)\scrP \cap \zz^n} \to \xaa{k\scrP \cap \zz^n}$ is an isomorphism.
\end{prop}

\begin{proof}
Let $\alpha_0, \ldots, \alpha_s$ be the vertices of $\scrP$, so that $k\alpha_j$, $j = 0, \ldots, s$, are the vertices of $k\scrP$. \Cref{cartier-infty} implies that $\xaa{k\scrP \cap \zz^n}$ is the union of affine open sets $\xaa{k\scrP \cap \zz^n} \cap U_{k\alpha_j}$, $j = 0, \ldots, s$. For each $j$, \cref{xa-toric} implies that the coordinate ring of $\xaa{k\scrP \cap \zz^n} \cap U_{k\alpha_j}$ is the semigroup algebra $\kk[S_{j,k}]$, where $S_{j,k}$ is the subsemigroup of $\zz^n$ generated by $\{\alpha - k\alpha_j: \alpha \in k\scrP \cap \zz^n \}$. Since $\alpha - k\alpha_j = (\alpha + \alpha_j) - (k+1)\alpha_j$, it follows that $S_{j,k} \subset S_{j, k+1}$, and $\pi_{k,k+1}$ is simply the morphism induced by this inclusion. Note that each $S_{j,k}$ is a subsemigroup of $\scrC_{\alpha_j} \cap \zz^n$, where $\scrC_{\alpha_j}$ is the {\em rational} convex polyhedral cone in $\rr^n$ generated by $\{\alpha_i - \alpha_j: i = 0, \ldots, s\}$. Gordan's lemma (\cref{gordan}) implies that $\scrC_{\alpha_j}  \cap \zz^n$ is finitely generated. Due to \cref{exercise:C-alpha-j} below we may choose an integer $K$ such that for each $j$, $S_{j,k}$ contains each of these generators for each $k > K$. For each $k > K$ and each $j$, it follows that $S_{j,k} = \scrC_{\alpha_j}  \cap \zz^n$ is independent of $k$; consequently $\xaa{(k+1)\scrP \cap \zz^n} \cong \xaa{k\scrP \cap \zz^n}$.
\end{proof}

\begin{figure}[h]
\begin{center}

\def\xmin{-0.5}
\def\xmax{13.5}
\def\ymin{-5.5}
\def\ymax{5.5}
\def\opazero{0.5}
\def\colorzero{green}

\tikzstyle{dot} = [red, circle, minimum size=4pt, inner sep = 0pt, fill]
\tikzstyle{smalldot} = [black, circle, minimum size=3pt, inner sep = 0pt, fill]

\begin{tikzpicture}[scale=\scalefactor]
\draw [gray,  line width=0pt] (\xmin, \ymin) grid (\xmax,\ymax);
\draw [<->] (0, \ymax) |- (\xmax, 0);
\node[dot] (A) at (1,-4) {};
\node[dot] (B) at (13,2) {};
\node[dot] (C) at (10,5) {};
\fill[\colorzero, opacity=\opazero ] (A.center) --  (B.center) -- (C.center);
\draw[thick, \colorAB] (A) -- (B);
\draw[thick, \colorBC] (B) -- (C);
\draw[thick, \colorCA] (C) -- (A);
\node[anchor = north] at (A) {$A$};
\node[anchor = west] at (B) {$B$};
\node[anchor = south] at (C) {$C$};

\node at (8,1) {\picfontsize $\scrP$};

\def\xmin{-5.5}
\def\xshift{5}
\pgfmathsetmacro\shiftone{\xshift + \xmax -\xmin};
\def\xmax{5.5}

\begin{scope}[shift={(\shiftone,0)}]
\draw [gray,  line width=0pt] (\xmin, \ymin) grid (\xmax,\ymax);
\coordinate (O) at (0,0);
\coordinate (BA) at (\xmin, \xmin/2);
\coordinate (AB) at (\xmax, \xmax/2);
\coordinate (CA) at (\ymin, \ymin);
\coordinate (AC) at (\ymax, \ymax);
\coordinate (BC) at (-\ymax, \ymax);
\coordinate (CB) at (-\ymin, \ymin);

\fill[\colorzero, opacity=\opazero ] (O) --  (BA) -- (BC);
\fill[\colorzero, opacity=\opazero ] (O) --  (AC) -- (AB);
\fill[\colorzero, opacity=\opazero ] (O) --  (CA) -- (CB);

\draw [thick, \colorAB] (AB) -- (BA);
\draw [thick, \colorBC] (CB) -- (BC);
\draw [thick, \colorCA] (AC) -- (CA);

\node at (-3.5,0.5) {\picfontsize $\scrC_B$};
\node at (4.5,3.5) {\picfontsize $\scrC_A$};
\node at (0.5,-3.5) {\picfontsize $\scrC_C$};

\node[smalldot] at (-1, -1) {};
\node[smalldot] at (0, -1) {};
\node[smalldot] at (1, -1) {};

\node[smalldot] at (2,1) {};
\node[smalldot] at (1,1) {};

\node[smalldot] at (-1,1) {};
\node[smalldot] at (-1,0) {};
\node[smalldot] at (-2,-1) {};

\end{scope}
\end{tikzpicture}
\caption{Cones of regular functions on basic open subsets of $\xp$. The dots in each cone mark generators of the semigroup of integral points in the cone.}  \label{fig:xp-support-algebra}
\end{center}
\end{figure}

Given a convex integral polytope $\scrP \subseteq \rr^n$, we write \index{Toric variety!from a convex integral polytope}$\xp$ for the toric variety isomorphic to $\xaa{k \scrP \cap \zz^n}$ for sufficiently large $k$ and $\xzerop$ for the torus $\xzeroaa{k\scrP \cap \zz^n}$ of $\xp$; we also write $\phi_\scrP$ for the morphism $\phi_{k\scrP \cap \zz^n}: \nktorus \to \xzerop$ defined as in \eqref{phi_A}. The arguments from the proof of \cref{xprop} show that $\xp$ is the union of open affine subsets $U'_\alpha := \xp \cap U_\alpha$ corresponding to vertices $\alpha$ of $\scrP$, and the coordinate ring of each $U'_\alpha$ is generated by the monomials whose exponents belong to the cone $\scrC_\alpha$ generated by $\scrP - \alpha$. See \cref{fig:xp-support-algebra} for the cones corresponding to the vertices of a triangle. The following proposition summarizes some basic properties of $\xp$.

\begin{prop} \label{xp-thm}
Let $\scrA$ be a finite subset of $\zz^n$ and $\scrP$ be the convex hull of $\scrA$ in $\rr^n$.  Let $G_\scrP$ be the subgroup of $\zz^n$ generated by pairwise differences of integral elements in the affine hull of $\scrP$ and $G_\scrA$ be the subgroup of $G_\scrP$ generated by the pairwise differences of elements from $\scrA$.
\begin{enumerate}
\item \label{phi-P} If $k$ is such that $G_\scrP$ is generated by $\{\alpha - \beta: \alpha, \beta \in k\scrP \cap \zz^n\}$, then $\xp \cong \xaa{k \scrP \cap \zz^n}$.
\item \label{torus-identification} The dimension of $\xp$ is the same as the dimension of $\scrP$. If $\dim(\scrP) = n$, then $\phi_\scrP$ is an isomorphism between $\nktorus$ and $\xzerop$; i.e.\ $\xp$ is a compactification of $\nktorus$.
\item There is a natural finite-to-one morphism $\phi_{\scrP, \scrA}: \xp \to \xa$ defined as the composition of the following maps:
\begin{align}
\xp
	\cong \xaa{k\scrP \cap \zz^n}
 	\xrightarrow{\pi_{k\scrP \cap \zz^n, k\scrA}} \xaa{k\scrA}
 	\xrightarrow{\phi_{k\scrA, \scrA}} \xa
 \label{xp-to-xa-short}
\end{align}
where $k$ is as in assertion \eqref{phi-P}, $\pi_{k\scrP \cap \zz^n, k\scrA}$ is the projection which drops coordinates corresponding to elements in $(k\scrP \cap \zz^n) \setminus k\scrA$, and $\phi_{k\scrA, \scrA}$ is the inverse of the Veronese isomorphism of degree $k$ (see \cref{segre-veronese}). The degree of $\phi_{\scrP, \scrA}$ is the index of $G_\scrA$ in $G_\scrP$. In particular, if $G_\scrA= G_\scrP$, then $\phi_{\scrP, \scrA}$ restricts to an isomorphism between $\xzerop$ and $\xzeroa$.
\item \label{C-S-alpha} There is an open cover $\{U'_\alpha\}_{\alpha \in  \scrP \cap \zz^n}$ of $\xp$ such that the coordinate ring of each $U'_\alpha$ is the semigroup algebra $\kk[S_\alpha]$, where $S_\alpha$ is the semigroup of integral points in the convex polyhedral cone $\scrC_\alpha$ generated by $\{\alpha' - \alpha: \alpha' \in \scrP\}$.
\item There is a natural one-to-one correspondence between the faces of $\scrP$ and the orbits of $\xzerop$ on $\xp$. For each face $\scrQ$ of $\scrP$, let $\orbit{\scrQ}$ be the corresponding orbit and $\corbit{\scrQ}$ be the closure of $\orbit{\scrQ}$ in $\xp$. Then $\corbit{\scrQ}$ is naturally isomorphic to $\xq$, and the isomorphism identifies $\orbit{\scrQ}$ with $\xzeroq$.
\item In particular, $\xp \setminus \xzerop$ is the union of the $\corbit{\scrQ}$ for the facets $\scrQ$ of $\scrP$.
\item \label{xp:diagonal} Let $\scrP'$ be a convex integral polytope in $\rr^n$. Then $\xaa{\scrP+\scrP'}$ is isomorphic to the closure in $\xp \times \xaa{\scrP'}$ of the image of the diagonal map $\nktorus \to \xp \times \xaa{\scrP'}$ which sends $x \mapsto (\phi_{\scrP}(x), \phi_{\scrP'}(x))$. Moreover, $\xaa{k\scrP+l\scrP'} \cong \xaa{\scrP+ \scrP'}$ for every pair of positive integers $k,l$.
\end{enumerate}
\end{prop}

\begin{proof}
Due to \cref{xa-invariant} we may assume $\dim(\scrP) = n$ and $\scrP$ contains the origin. If $k$ is as in assertion \eqref{phi-P}, then after another application of \cref{xa-invariant} we may assume that each standard unit vector is in $k\scrP \cap \zz^n$, which immediately implies that $\phi_\scrP$ is an isomorphism and proves the first two assertions. The last assertion is a straightforward corollary of \cref{segre-veronese,xprop}.  The remaining statements follow from \cref{torus-image,xa-toric,xa-thm}. We leave it as an exercise to complete the proof.
\end{proof}

\begin{example} \label{example:deg2}
Let $\scrA := \{0, 2\} \subset \zz$. Since $x \mapsto 2x$ maps $\{0,1\}$ onto $\scrA$, \cref{xa-invariant,example:rational-normal-curve} imply that $\xa \cong \pp^1$. If $\scrP = \conv(\scrA) \subset \rr$, then $G_\scrP = \zz$ and assertion \eqref{phi-P} of \cref{xp-thm} is satisfied with $k = 1$, i.e.\ $\xp$ is the closure in $\pp^2$ of the image of $\kk^*$ under the map $x \mapsto [1:x:x^2]$. \Cref{example:rational-normal-curve} implies that $\xp$ is also isomorphic to $\pp^1$. The map $\phi_{\scrP, \scrA}: \xp \to \xa$ is the restriction to $\xp$ of the projection $[z_0:z_1:z_2] \to [z_0: z_2]$, and therefore on the level of the tori it is simply the map $x \mapsto x^2$. Note that $\deg(\phi_{\scrP, \scrA}) = 2$ is also the index of $G_\scrA = 2\zz$ in $\zz = G_\scrP$.
\end{example}

\begin{example} \label{example:xp^n}
Let $\scrS^n$ be the $n$-dimensional simplex in $\rr^n$ with vertices at the origin and at the elements of the standard unit basis of $\rr^n$, and let $\scrK^n := [0,1]^n \subset \rr^n$. Then $G_{\scrS ^n}= G_{\scrK^n} = \zz^n$, and both $\scrS^n$ and $\scrK^n$ satisfy assertion \eqref{phi-P} of \cref{xp-thm} with $k = 1$. Therefore it follows from \cref{example:p^n,example:cube} that $\xpp{\scrS^n} \cong \pp^n$ and $\xpp{\scrK^n} \cong (\pp^1)^n$.
\end{example}

\begin{example} \label{A'-example}
Consider $\scrA'$ from \cref{fig:A'}. The convex hull $\scrP'$ of $\scrA$ is a translation of $9\scrS^2$, so that \cref{segre-veronese,xp-thm,example:xp^n} imply that $\xpp{\scrP'} \cong \xpp{\scrS^2} \cong \pp^2$. Note that $G_{\scrP'} = G_{\scrA'} = \zz^2$, i.e.\ the map $\phi_{\scrP', \scrA'}: \xpp{\scrP'} \to \xpp{\scrA'}$ has degree one. However, it is {\em not} an isomorphism. Indeed, \cref{xa-toric} implies that the coordinate ring of $\xpp{\scrA'} \cap U_{A'}$ is $\kk[x^8, x^9, y] \cong \kk[u,v,w]/\langle u^9 - v^8 \rangle$ (\cref{exercise:desingxample2}). Therefore $\xpp{\scrA'} \cap U_{A'}$ is isomorphic to the hypersurface in $\kk^3$ defined by $u^9 - v^8 = 0$, which is {\em singular} at all points of the $w$-axis, and $\phi_{\scrP', \scrA'}: \xpp{\scrP'} \to \xpp{\scrA'}$ is a {\em desingularization} of $\xpp{\scrA'}$.
\end{example}

Even though $\xp$ is by definition isomorphic to $\xaa{k\scrP \cap \zz^n}$ for $k \gg 1$, in each of the preceding examples it suffices to take $k = 1$. In general it suffices to take $k \geq \dim(\scrP) - 1$, see e.g.\ \cite[Theorem 2.2.12]{littlehalcox}. In particular, to find examples for which one needs to take $k > 1$ requires polytopes with dimensions at least $3$.

\subsection{Exercises}
\begin{exercise}\label{exercise:equi-hull-projection}
Let $\scrA \subset \scrA'$ be finite subsets of $\zz^n$ such that $\conv(\scrA) = \conv(\scrA')$. Show that $\pi_{\scrA,\scrA'}$ is well-defined everywhere on $\xaa{\scrA'}$. [Hint: $\scrA$ and $\scrA'$ have the same set of vertices. Use \cref{cartier-infty}.]
\end{exercise}

\begin{exercise} \label{exercise:C-alpha-j}
In the notation of the proof of \cref{xprop}, show that for each $\beta \in \scrC_{\alpha_j}$, there is $K \geq 0$ such that $\beta \in k\scrP - k\alpha_j$ for each $k \geq K$.
\end{exercise}

\begin{exercise} \label{exercise:xp}
Complete the proof of \cref{xp-thm}.
\end{exercise}

\begin{exercise} \label{exercise:desingxample2}
Show that $\kk[x^8, x^9, y] \cong \kk[u,v,w]/\langle u^9 - v^8 \rangle$. Prove that the hypersurface $V(u^9 - v^8)$ in $\kk^3$ is singular at all points of the $w$-axis.
\end{exercise}

\begin{exercise} \label{exercise:singxp}
In the set up of \cref{fig:xp-support-algebra}, let $S_A$ (respectively, $S_B,S_C$) be the subsemigroup of $\zz^n$ consisting of integral elements in $\scrC_A$ (respectively, $\scrC_B,\scrC_C$).
\begin{enumerate}
\item Show that $S_A$ is generated as a semigroup by $(1,1), (2,1)$; $S_B$ is generated as a semigroup by $(-1,1), (-1,0), (-2,-1)$; $S_C$ is generated as a semigroup by $(-1,-1), (0,-1), (1,-1)$.
\item Deduce that $\xp \cap U_A \cong \kk^2$, $\xp \cap U_B \cong V(uv - w^3) \subset \kk^3$, and $\xp \cap U_C \cong V(uv - w^2) \subset \kk^3$. [Hint: $\kk[\xp \cap U_A] = \kk[xy, x^2y]$, $\kk[\xp \cap U_B] = \kk[x^{-1}y, x^{-1},x^{-2}y^{-1}]$, $\kk[\xp \cap U_C] = \kk[x^{-1}y^{-1}, y^{-1},xy^{-1}]$.]
\item Conclude that $\xp$ has precisely two singular points.
\end{enumerate}
\end{exercise}

\begin{exercise} \label{exercise:invertible-on-intersection}
Let $\scrP$ be a convex integral polytope in $\rr^n$ and $\alpha, \alpha'$ be vertices of $\scrP$.
\begin{enumerate}
\item \label{invertible:0} Show that $x^{\alpha' - \alpha}$ is an invertible regular function on $U_\alpha \cap U_{\alpha'} \cap \xp$.
\item  Let $\scrQ$ be the smallest face of $\scrP$ containing both $\alpha$ and $\alpha'$, and $\scrQ_\alpha$ be the cone generated by $\{\beta - \alpha: \beta \in \scrQ\}$. Show that $\alpha' - \alpha$ is in the relative interior of $\scrQ_\alpha$.
\item Let $H$ be the linear subspace of $\rr^n$ spanned by $\{\beta - \alpha: \beta \in \scrQ\}$. For each $\gamma \in H \cap \zz^n$, show that $x^\gamma$ is an invertible regular function on $U_\alpha \cap U_{\alpha'} \cap \xp$. [Hint: Choose $\beta_1, \ldots, \beta_k \in \scrQ \cap \zz^n$ such that $\beta_j - \alpha$, $j = 1, \ldots, k$, generate $\scrQ_\alpha$ as a cone. Use \cref{exercise:relatively-positive-interior} and the preceding assertions to show that for each $j$, $x^{\beta_j - \alpha}$ is a regular function on $U_\alpha \cap \xp$ which does not vanish at any point of $U_\alpha \cap U_{\alpha'} \cap \xp$.]
\item Deduce that if $\scrQ = \scrP$, then $U_\alpha \cap U_{\alpha'} \cap \xp = \xzerop$.
\end{enumerate}
\end{exercise}

\section{Nonsingularity in codimension one on $\xp$} \label{U_Q-section}
\Cref{A'-example,exercise:singxp} above show that toric varieties $\xp$ from polytopes might be singular. However, we will see in this section that the $\xp$ is nonsingular outside a subvariety of dimension at most $\dim(\xp) - 2$, i.e.\ $\xp$ is ``nonsingular in codimension one.'' Note that this is in general {\em not} true for varieties $\xa$ (see \cref{A'-example} above). We continue to use the notation of \cref{toric-polysection}. Let $\scrQ$ be a facet of $\scrP$. We will show that $\xp$ is nonsingular at all points of $\orbit{\scrQ}$. Due to \cref{xa-invariant} we may assume \woutlog\ that $\scrP$ is {\em full dimensional}, i.e.\ $\scrP$ is an $n$-dimensional polytope in $\rr^n$. Let $\nu$ be the {\em primitive inner normal} (see \cref{primitive-remark}) to $\scrQ$. Let $\znnuzero := \{\beta \in \zz^n: \langle \nu, \beta \rangle \geq 0\}$ and $\znnuperp := \{\beta \in \zz^n: \langle \nu, \beta \rangle = 0\}$. Choose an arbitrary element $\alpha_\nu \in \zz^n$ such that $\langle \nu, \alpha_\nu \rangle = 1$. We write $\zzero \langle \alpha_\nu \rangle := \{k\alpha_\nu: k \in \zz,\ k \geq 0\}$. The following is a straightforward implication of \cref{basis-lemma}.

\begin{lemma} \label{lemma:nuzero-decomposition}
There is a basis of $\zz^n$ of the form $\beta_1, \ldots, \beta_n$, where $\beta_n = \alpha_\nu$ and $\beta_1, \ldots, \beta_{n-1}$ constitute a basis of $\znnuperp$. In particular, $\znnuperp \cong \zz^{n-1}$, and as a semigroup $\znnuzero$ is isomorphic to $\znnuperp + \zzero \langle \alpha_\nu \rangle \cong \zz^{n-1} \times \zzero$. \qed
\end{lemma}

\begin{prop} \label{U_Q}
Let $U_\scrQ := \xzerop \cup \orbit{\scrQ}$.
\begin{enumerate}
\item \label{U_Q:open} $U_\scrQ$ is an open affine neighborhood of $\orbit{\scrQ}$ in $\xp$.
\item \label{U_Q:coordinates} $\kk[U_\scrQ] = \kk[x^\beta : \beta \in \znnuzero] \cong \kk[\zz^{n-1} \times \zzero]$. In particular, $U_\scrQ \cong \nktoruss{n-1} \times \kk$.
\item \label{O_Q-nonsingular} $\orbit{\scrQ} = V(x^{\alpha_\nu}) \subset U_\scrQ$ and $\kk[\orbit{\scrQ}] \cong \kk[U_\scrQ]/\langle x^{\alpha_\nu} \rangle \cong \kk[\znnuperp]$. In particular, the embedding $\orbit{\scrQ} \into U_\scrQ$ is isomorphic to the embedding $\nktoruss{n-1} \times \{0\} \into  \nktoruss{n-1} \times \kk$.
\end{enumerate}
\end{prop}

\begin{proof}
Let $\scrA := k\scrP \cap \zz^n$, where $k$ is large enough so that $\xp \cong \xa$ and there is an integral element $\alpha_0$ in the relative interior of $k\scrQ$. Since $\scrQ$ is the smallest face of $\scrP$ containing $\alpha_0$ (\cref{prop:relative-containment}),  \cref{point-face} implies that $U_\scrQ = \xp \cap U_{\alpha_0}$, which proves assertion \eqref{U_Q:open}. Let $\beta_1, \ldots, \beta_n = \alpha_\nu$ be a basis of $\zz^n$ as in \cref{lemma:nuzero-decomposition}. Choosing a large enough $k$ we can also ensure that $\alpha_0 + \beta_j \in k\scrP \cap \zz^n$ for each $j = 1, \ldots, n$ (\cref{exercise:deformation-relint}). It is then straightforward to check that the semigroup generated by $\scrA - \alpha_0$ is precisely $\znnuzero$. Assertions \eqref{U_Q:coordinates} and \eqref{O_Q-nonsingular} then follow from \cref{xp-thm,lemma:nuzero-decomposition} in a straightforward manner - we leave the proof as an exercise.
\end{proof}

\begin{cor}[Nonsingularity of $\xp$ in codimension one]
\index{Nonsingularity in codimension one}
The set $\sing(\xp)$ of singular points of $\xp$ is contained in $\bigcup_i  \orbit{\scrQ'_{i}}$, where $\scrQ'_i$ are faces of $\scrP$ of dimension $\leq n - 2$. In particular, $\dim(\sing(\xp)) \leq n - 2$.
\end{cor}

\begin{proof}
$\xp \setminus \bigcup_i  \orbit{\scrQ'_{i}}$ is the union of $U_{\scrQ}$ over all facets $\scrQ$ of $\scrP$, which is nonsingular due to \cref{U_Q}. Since $\dim(\bigcup_i  \orbit{\scrQ'_{i}}) = n - 2$, the result follows.
\end{proof}

\def\colorg{orange}
\def\gp#1#2#3{
\begin{scope}[shift={(#1,#2)}]
	\def\gh{1}
	\draw[fill=\colorg, opacity = \opazero] (0,0) -- (-2*\gh,0) -- (-4*\gh,-2*\gh) -- cycle;
	\node[anchor=north] at (-2*\gh,-2*\gh) {\picfontsize #3};
\end{scope}
}

\begin{figure}[h]
\begin{center}

\def\xmin{-0.5}
\def\xmax{13.5}
\def\ymin{-5.5}
\def\ymax{5.5}
\def\opazero{0.5}
\def\colorzero{green}

\tikzstyle{dot} = [red, circle, minimum size=4pt, inner sep = 0pt, fill]

\begin{tikzpicture}[scale=\scalefactor]
\draw [gray,  line width=0pt] (\xmin, \ymin) grid (\xmax,\ymax);
\draw [<->] (0, \ymax) |- (\xmax, 0);
\node[dot] (A) at (1,-4) {};
\node[dot] (B) at (13,2) {};
\node[dot] (C) at (10,5) {};
\fill[\colorzero, opacity=\opazero ] (A.center) --  (B.center) -- (C.center);
\draw[thick, \colorAB] (A) -- (B);
\draw[thick, \colorBC] (B) -- (C);
\draw[thick, \colorCA] (C) -- (A);
\node[anchor = north] at (A) {\picfontsize $A$};
\node[anchor = west] at (B) {\picfontsize $B$};
\node[anchor = south] at (C) {\picfontsize $C$};
\node[anchor = west] at ($(B)!0.6!(C)$) {\picfontsize $\scrQ$};

\node at (8,1) {\picfontsize $\scrP$};

\gp{4}{4}{$\np(g)$}

\def\xmin{-5.5}
\def\xshift{5}
\pgfmathsetmacro\shiftone{\xshift + \xmax -\xmin};
\def\xmax{5.5}
\begin{scope}[shift={(\shiftone,0)}]
\draw [gray,  line width=0pt] (\xmin, \ymin) grid (\xmax,\ymax);

\coordinate (BC) at (\xmin, -\xmin);
\coordinate (CB) at (\xmax, -\xmax);
\draw [thick, \colorBC] (BC) -- (CB);

\coordinate (tshift) at (0,0.2);

\coordinate (perp) at (-1,-1);
\draw [thick, \colorBC, ->] (0,0) -- (perp);
\node [anchor = north] at ($(perp) + (tshift)$) {\picfontsize $\nu = (-1,-1)$};

\def\alphax{-4}
\def\alphay{4}
\gp{\alphax}{\alphay}{$\np(\talphanu(g))$}
\node[dot] (alpha) at (\alphax, \alphay) {};
\coordinate (tshift) at (-0.15,0.3);
\node [anchor=west] at ($(alpha) + (tshift)$) {\picfontsize $\alpha_\nu = (-1,0)$};

\def\alphax{4}
\def\alphay{-4}
\gp{\alphax}{\alphay}{$\np(\talphanu(g))$}
\node[dot] (alpha) at (\alphax, \alphay) {};
\node [anchor= west] at ($(alpha) + (tshift)$) {\picfontsize $\alpha_\nu = (0,-1)$};

\end{scope}
\end{tikzpicture}
\caption{Different choices for $\talphanu(g)$ when $g = 2y^2 - 3x^2y^3 + 7x^4y^4 - 6x^2y^4$}  \label{fig:g-nu}
\end{center}
\end{figure}

Let $\nu$ and $\alpha_\nu$ be as in \cref{U_Q}. Given $g  = \sum_\beta c_\beta x^\beta \in \kk[x_1, x_1^{-1}, \ldots, x_n, x_n^{-1}]$, the \index{Support!of a Laurent polynomial}{\em support} $\supp(g)$ of $g$ is the set of all $\beta \in \zz^n$ such that $c_\beta \neq 0$, and the \index{Newton!polytope}{\em Newton polytope} $\np(g)$ of $g$ is the convex hull of $\supp(g)$. Let $m := \min_{\supp(g)}(\nu) = \min_{\np(g)}(\nu)$. Choose an arbitrary isomorphism $\psi_\nu:  \znnuperp \cong \zz^{n-1}$. Define $\talphanu(g) := x^{-m\alpha_\nu} g$ and $\Inalphapsinu(g) :=  \sum_{ \langle \nu, \beta \rangle = m} c_{\beta}x^{\psi_\nu(\beta - m\alpha_\nu)}$ (the ``$T$'' in $\talphanu(\cdot)$ is supposed to imply ``translation,'' and ``$\In$'' in $\Inalphapsinu(\cdot)$ is to suggest ``initial form''). See \cref{fig:g-nu} for an example with $\scrP$ from \cref{fig:xp-support-algebra}. The following result is an immediate corollary of \cref{U_Q}. Its proof is left as an exercise.

\begin{cor} \label{O_Q}
\begin{enumerate}
\item $\talphanu(g)$ is a regular function on $U_\scrQ$ for each Laurent polynomial $g$.
\item The correspondence $\talphanu(g)|_{\orbit{\scrQ}} \mapsto \Inalphapsinu(g)$ induces an isomorphism $\psi_\nu^*: \nktoruss{n-1} \cong \orbit{\scrQ}$.
\item If $\alpha'_\nu$ is another element in $\zz^n$ such that $\langle \nu, \alpha'_\nu \rangle = 1$, then $\talphanu(g)/\talphanuu{\alpha'_\nu}(g) = x^{m(\alpha'_\nu - \alpha_\nu)}$ is invertible on $U_\scrQ$.
\end{enumerate}
\end{cor}

A basic property of the varieties $\xp$ is that they are ``normal,'' and nonsingularity in codimension one follows from normality. In this book we do not treat this notion - see \cite[Section 2.1]{fultoric} or \cite[Section 2.4]{littlehalcox} for an exposition of this and other fundamental properties of $\xp$ including the following result (which we do not use): assume $\scrP \subset \rr^n$ is full dimensional. Then $\xp$ is nonsingular if and only if both of the following are true for every vertex $\alpha$ of $\scrP$:
\begin{itemize}
\item $\scrC_\alpha$ has precisely $n$ edges, and
\item the primitive integral elements of the edges of $\scrC_\alpha$ form a basis of $\zz^n$.
\end{itemize}
(It then follows due to \cref{toric-k^n} that if $\xp$ is nonsingular, then it is a compactification of $\kk^n$.) 
\subsection{Exercises}

\begin{exercise}
Complete the proof of \cref{U_Q}.
\end{exercise}

\begin{exercise}
Prove \cref{O_Q}.
\end{exercise}

\section{Extending closed subschemes of the torus to $\xp$} \label{xone-section}
Recall (from \cref{subscheme-section}) that a {\em closed subscheme} $V$ of a variety $X$ is essentially a Zariski closed subset $V'$ of a variety together with a {\em sheaf of ideals} $\scrI$ such that $V'$ is precisely the set of zeroes of elements in $\scrI$. If $\bar X$ is a variety containing $X$ as a Zariski open subset, and $\bar V$ is a closed subscheme of $\bar X$, we say that $\bar V$ \index{Extension of a closed subscheme}{\em extends} $V$ if the scheme-theoretic intersection $\bar V \cap X$ is precisely $V$. In this section we will study the case that $X = \nktorus$ and $\bar X$ is the toric variety $\xp$ corresponding to an $n$-dimensional convex integral polytope $\scrP$. We are specially interested in the case that
\begin{defnlist}
\item \label{Cartier-extension} both $V$ and $\bar V$ are {\em Cartier divisors} (see \cref{scheme:cartier-section}), and
\item \label{tight-extension} $\supp(\bar V)$ is the closure of $\supp(V)$.
\end{defnlist}
Each Laurent polynomial $g \in \kk[x_1, x_1^{-1}, \ldots, x_n, x_n^{-1}]$ defines a Cartier divisor $V(g)$ on $\nktorus$ (in fact it is not hard to see, using the fact that the ring of Laurent polynomials is a UFD, that {\em every} Cartier divisor on $\nktorus$ is of the form $V(g)$ for some Laurent polynomial $g$ - but we will not use it). There are many ways to extend $V(g)$ to $\xp$, e.g.\ if $\np(g) \subset \scrP$, then \cref{exercise:in-extension} below prescribes a way to extend $V(g)$ to a Cartier divisor on $\xp$ which satisfies property \ref{tight-extension} if and only if $\np(g)$ intersects each facet of $\scrP$. However, we will shortly see that for some $\scrP$ there are Cartier divisors on $\nktorus$ which can {\em not} be extended to $\xp$ in a way to satisfy property \ref{tight-extension}. Then we will define an open subset $\xonep$ of $\xp$ such that $\xp \setminus \xonep$ has dimension $\leq n -2$, and every Cartier divisor on $\nktorus$ does admit extensions to $\xonep$ satisfying both properties \ref{Cartier-extension} and \ref{tight-extension}.

\def\colorh{purple}
\def\hp#1#2#3{
\begin{scope}[shift={(#1,#2)}]
	\draw[fill=\colorh, opacity = \opazero] (0,0) -- (2,0) -- (2,2) -- cycle;
	\node[anchor=north] at (1,0) {\picfontsize #3};
\end{scope}
}

\begin{center}
\begin{figure}[h]
\def\xmin{-0.5}
\def\xmax{13.5}
\def\ymin{-5.5}
\def\ymax{5.5}
\tikzstyle{dot} = [red, circle, minimum size=4pt, inner sep = 0pt, fill]

\begin{tikzpicture}[scale=\scalefactor]
\draw [gray,  line width=0pt] (\xmin, \ymin) grid (\xmax,\ymax);
\draw [<->] (0, \ymax) |- (\xmax, 0);
\node[dot] (A) at (1,-4) {};
\node[dot] (B) at (13,2) {};
\node[dot] (C) at (10,5) {};
\fill[\colorzero, opacity=\opazero ] (A.center) --  (B.center) -- (C.center);
\draw[thick, \colorAB] (A) -- (B);
\draw[thick, \colorBC] (B) -- (C);
\draw[thick, \colorCA] (C) -- (A);
\node[anchor = north] at (A) {\picfontsize $A$};
\node[anchor = west] at (B) {\picfontsize $B$};
\node[anchor = south] at (C) {\picfontsize $C$};

\node at (8,1) {\picfontsize $\scrP$};
\gp{4}{4}{$\np(g)$}
\hp{9}{-3}{$\np(h)$}

\def\xmin{-5.5}
\def\xshift{5}
\pgfmathsetmacro\shiftone{\xshift + \xmax -\xmin};
\def\xmax{5.5}
\begin{scope}[shift={(\shiftone,0)}]

\draw [gray,  line width=0pt] (\xmin, \ymin) grid (\xmax,\ymax);
\coordinate (O) at (0,0);
\coordinate (BA) at (\xmin, \xmin/2);
\coordinate (AB) at (\xmax, \xmax/2);
\coordinate (CA) at (\ymin, \ymin);
\coordinate (AC) at (\ymax, \ymax);
\coordinate (BC) at (-\ymax, \ymax);
\coordinate (CB) at (-\ymin, \ymin);

\draw [thick, \colorAB] (AB) -- (BA);
\draw [thick, \colorBC] (CB) -- (BC);
\draw [thick, \colorCA] (AC) -- (CA);

\gp{0}{0}{}
\gp{1}{-1}{}
\gp{4}{2}{}

\fill[\colorzero, opacity=\opazero ] (0,0) -- (-2,0) -- (-4,-2) -- (BA) -- (BC);
\fill[\colorzero, opacity=\opazero ] (0,0) --  (-1,-1) -- (1,-1);
\fill[\colorzero, opacity=\opazero ] (1,-1) -- (-3,-3) -- (CA) -- (CB);
\fill[\colorzero, opacity=\opazero ] (2,2) -- (4,2) -- (AB) -- (AC);

\node at (-4.5,0.5) {\picfontsize $\scrC_B$};
\node at (5.5,3.5) {\picfontsize $\scrC_A$};
\node at (-0.5,-4.5) {\picfontsize $\scrC_C$};
\end{scope}

\pgfmathsetmacro\shiftwo{\shiftone + \xshift + \xmax -\xmin};
\begin{scope}[shift={(\shiftwo,0)}]

\draw [gray,  line width=0pt] (\xmin, \ymin) grid (\xmax,\ymax);
\coordinate (O) at (0,0);
\coordinate (BA) at (\xmin, \xmin/2);
\coordinate (AB) at (\xmax, \xmax/2);
\coordinate (CA) at (\ymin, \ymin);
\coordinate (AC) at (\ymax, \ymax);
\coordinate (BC) at (-\ymax, \ymax);
\coordinate (CB) at (-\ymin, \ymin);

\draw [thick, \colorAB] (AB) -- (BA);
\draw [thick, \colorBC] (CB) -- (BC);
\draw [thick, \colorCA] (AC) -- (CA);

\hp{-2}{-2}{}
\hp{2}{2}{}
\hp{-10/3}{-2/3}{}

\fill[\colorzero, opacity=\opazero ] (0,0) -- (-4/3,-2/3) -- (-4/3,4/3);
\fill[\colorzero, opacity=\opazero ] (-4/3,4/3) -- (-10/3,-2/3) -- (-4/3,-2/3) -- (BA) -- (BC);
\fill[\colorzero, opacity=\opazero ] (0,0) --  (0,-2) -- (-2,-2) -- (CA) -- (CB);
\fill[\colorzero, opacity=\opazero ] (0,0) -- (4,2) -- (2,2);
\fill[\colorzero, opacity=\opazero ] (4,4) -- (4,2) -- (AB) -- (AC);

\node at (-4.5,0.5) {\picfontsize $\scrC_B$};
\node at (5.5,3.5) {\picfontsize $\scrC_A$};
\node at (-0.5,-4.5) {\picfontsize $\scrC_C$};
\end{scope}
\end{tikzpicture}
\caption{Extending Cartier divisors in dimension two}  \label{fig:sub2extension}
\end{figure}
\end{center}

\begin{example} \label{example:2-extension}
Consider $\scrP$ and $g$ from \cref{fig:g-nu}. \Cref{xa-thm,xp-thm} imply that $\xp \subset U_A \cup U_B \cup U_C$, and for each $P \in \{A,B,C\}$, the coordinate ring of $U_P \cap \xp$ consists of Laurent polynomials supported at the cone $\scrC_P$; consider the (unique) translation of $\np(g)$ which is contained in $\scrC_P$ and touches both sides of $\scrC_P$ (see the middle panel of \cref{fig:sub2extension}). More precisely, let
\begin{align*}
g_P :=
	\begin{cases}
		2 - 3x^2y + 7x^4y^2 - 6x^2y^2               & \text{if}\ P = A, \\
		x^{-4}y^{-2}(2 - 3x^2y + 7x^4y^2 - 6x^2y^2) & \text{if}\ P = B, \\
		x^{-3}y^{-3}(2 - 3x^2y + 7x^4y^2 - 6x^2y^2) & \text{if}\ P = C.
	\end{cases}
\end{align*}
Then it is straightforward to check that the pairs $(U_P \cap \xp, g_P)$, $P \in \{A, B,C\}$, defines a Cartier divisor $D$ on $\xp$ such that $D$ extends the Cartier divisor $V(g)$ and $\supp(D)$ is the closure in $\xp$ of $V(g) \subset \nktoruss{2}$ (see \cref{exercise:simplicial-extension} for a more general result). On the other hand, if $h$ is any Laurent polynomial with $\np(h)$ as in the left panel of \cref{fig:sub2extension}, then the translation of $\np(h)$ that is contained in $\scrC_B$ and touches both sides of $\scrC_B$ is {\em not} integral - see the right panel of \cref{fig:sub2extension}. It follows that if $E$ is any Cartier divisor on $\xp$ which extends $V(h) \subset (\kk^*)^2$, then $\supp(E)$ must contain either $V_{AB}$ or $V_{BC}$, so that $\supp(E)$ is larger than the closure in $\xp$ of $V(h) \subset (\kk^*)^2$ (\cref{exercise:non-2-extension}); in particular, $E$ does not satisfy property \ref{tight-extension} of an extension.
\end{example}

\begin{center}
\begin{figure}[h]
\def\scalefactor{0.4}
\def\xmin{-3.5}
\def\xmax{3.5}
\def\ymin{-3.5}
\def\ymax{3.5}
\def\opazero{0.5}
\def\colorzero{green}

\tikzstyle{dot} = [red, circle, minimum size=4pt, inner sep = 0pt, fill]

\begin{subfigure}[b]{0.3\textwidth}
\begin{tikzpicture}[scale=\scalefactor]
\draw [gray,  line width=0pt] (\xmin, \ymin) grid (\xmax,\ymax);
\draw [<->] (0, \ymax) |- (\xmax, 0);

\coordinate (A) at (-2,2);
\coordinate (B) at (-2,-2);
\coordinate (C) at (2,-2);
\coordinate (D) at (2,2);

\coordinate (E) at (-2,1);
\coordinate (F) at (-2,-1);
\coordinate (G) at (2,-1);
\coordinate (H) at (2,1);

\fill[\colorg, opacity=\opazero ] (E) --  (F) -- (G) -- (H);
\fill[\colorzero, opacity=\opazero ] (A) --  (E) -- (H) -- (D);
\fill[\colorzero, opacity=\opazero ] (C) --  (G) -- (F) -- (B);

\draw[thick] (A) -- (B) -- (C) -- (D);


\node [anchor = north] at (0,-2) {\picfontsize $\scrP$};
\node at (0,0) {\picfontsize $\np(g)$};

\end{tikzpicture}
\caption{Crossections of $\scrP$ and $\np(g)$}  \label{fig:subnonextension-cross}
\end{subfigure}
\def\viewx{60}%
\def\viewy{30}%
\begin{subfigure}[b]{0.3\textwidth}
\begin{tikzpicture}[scale=0.66]
\pgfplotsset{every axis title/.append style={at={(0,-0.2)}}, view={\viewx}{\viewy}, axis lines=middle, enlargelimits={upper}}

\begin{axis}
\addplot3[ thick, draw, fill=\colorzero,opacity=\opazero] coordinates{(-2,2,0) (-2,-2,0) (2,-2,0) (2,2,0)};
\addplot3[ thick, draw, fill=red,opacity=\opazero] coordinates{(0,0,2) (2,-2,0) (2,2,0)};
\addplot3[ thick, draw, fill=blue,opacity=\opazero] coordinates{(0,0,2) (2,-2,0) (-2,-2,0)};
\addplot3[draw, thick, dashed] coordinates {(0,0,2) (-2,2,0)};
\end{axis}
\end{tikzpicture}
\caption{$\scrP$}  \label{fig:subnonextension-P}
\end{subfigure}
\begin{subfigure}[b]{0.3\textwidth}
\def\viewx{60}
\def\viewy{30}
\begin{tikzpicture}[scale=0.66]
\pgfplotsset{every axis title/.append style={at={(0,-0.2)}}, view={\viewx}{\viewy}, axis lines=middle, enlargelimits={upper}}

\begin{axis}
\addplot3[ thick, draw, fill=\colorg,opacity=\opazero] coordinates{(-2,1,0) (-2,-1,0) (2,-1,0) (2,1,0)};
\addplot3[ thick, draw, fill=\colorg,opacity=\opazero] coordinates{(-2,1,2) (-2,-1,2) (2,-1,2) (2,1,2)};
\addplot3[ thick, draw, fill=red,opacity=\opazero] coordinates{(2,-1,2) (2,1,2) (2,1,0) (2,-1,0)};
\addplot3[ thick, draw, fill=blue,opacity=\opazero] coordinates{(2,-1,2) (-2,-1,2) (-2,-1,0) (2,-1,0)};
\addplot3[draw, thick, dashed] coordinates {(-2,1,2) (-2,1,0)};

\addplot3[draw] coordinates {(0,0,0) (0,-2,0)};
\addplot3[draw] coordinates {(0,0,0) (0,2,0)};
\end{axis}
\end{tikzpicture}
\caption{$\np(g)$}  \label{fig:subnonextension-g}
\end{subfigure}
\caption{Obstruction to extension of Cartier divisors in dimension $\geq 3$}  \label{fig:subnonextension}
\end{figure}
\end{center}

\begin{example} \label{example:3-non-extension}
Let $\scrP$ be the polytope in $\rr^3$ with vertices $(2,2,0)$, $(-2,2,0)$, $(-2, -2,0)$, $(2, -2,0)$, $(0, 0,2)$, and let $g \in \kk[x,y,z]$ be a polynomial such that $\np(g)$ is the ``cuboid'' with vertices $(2,1,0)$, $(-2,1,0)$, $(-2, -1,0)$, $(2, -1,0)$, $(2,1,2)$, $(-2,1,2)$, $(-2, -1,2)$, $(2, -1,2)$ (see \cref{fig:subnonextension}). It is straightforward to check that for each $m \geq 1$, if a translation $\scrQ$ of $m\np(g)$ is contained in $\scrC_P$, then $\scrQ$ can not touch all the facets of $\scrC_P$, and therefore every extension to $\xp$ of $V(g^m) \subset \nktoruss{3}$ by a Cartier divisor fails property \ref{tight-extension} (\cref{exercise:non-3-extension}).
\end{example}

In general, to find a Cartier divisor on $\xp$ which extends $V(g) \subset \nktorus$ and satisfies property \ref{tight-extension} requires a ``modification'' of $\xp$ by finding a ``common refinement'' of the normal fans of $\scrP$ and $\np(g)$. This leads to beautiful combinatorial geometry, but we will not get into this. We will rather find an open subset $\xonep$ of $\xp$ so that every Cartier divisor of $\nktorus$ can be extended to a Cartier divisor on $\xonep$ and the extension also satisfies property \ref{tight-extension}.

\subsection{The subset $\xonep$} \label{xonep-section}
Let $\scrP$ be an $n$-dimensional convex integral polytope in $\rr^n$. \Cref{xp-thm} implies that we can identify $\nktorus$ with $\xzerop$. Denote the facets of $\scrP$ by $\scrQ_1, \ldots, \scrQ_s$, and the faces of $\scrP$ of dimension $\leq n - 2$ by $\scrQ'_1, \ldots, \scrQ'_{s'}$. Define
\begin{align*}
\xonep
	:= \xzerop \cup \bigcup_{i=1}^s \orbit{\scrQ_i}
	= \bigcup_{i=1}^s U_{\scrQ_i}
	= \xp \setminus \bigcup_{i'=1}^{s'} \orbit{\scrQ'_{i'}}
	= \xp \setminus \bigcup_{i'=1}^{s'} \corbit{\scrQ'_{i'}}
\end{align*}
where the equalities are implications of \cref{xa-thm,xp-thm,U_Q}. In particular, it follows that $\xonep$ is Zariski open in $\xp$, and it is the union of torus orbits in $\xp$ of ``codimension smaller than one'' (this is the motivation for the ``$1$'' in the notation $\xonep$). Let $\nu_i$ be the primitive inner normal to $\scrQ_i$, $i = 1, \ldots, s$. Pick $\alpha_{\nu_i}$ such that $\langle \nu_i, \alpha_{\nu_i} \rangle = 1$. Given $g \in \kk[x_1, x_1^{-1}, \ldots, x_n, x_n^{-1}]$, \cref{O_Q} implies that $\talphanuu{\alpha_{\nu_i}}(g)$ is a regular function on $U_{\scrQ_i}$ for each $i$. Since $U_{\scrQ_i} \cap U_{\scrQ_j} = \xzerop \cong \nktorus$ for $i \neq j$ (\cref{U_Q}) and since $\talphanuu{\alpha_{\nu_i}}(g)/\talphanuu{\alpha_{\nu_j}}(g)$ is a monomial in $x_1, \ldots, x_n$, it follows that $\{(U_{\scrQ_i}, \talphanuu{\alpha_{\nu_i}}(g))\}_i$ defines a Cartier divisor $\vonep(g)$ on $\xonep$. In \cref{exercise:vonep-Cartier} you will check that $\vonep(g)$ is an extension to $\xonep$ of the Cartier divisor $V(g)$ of $\nktorus$, and it also satisfies property \ref{tight-extension}. This construction can be carried out for arbitrary closed subschemes of $\nktorus$. Indeed, given $g_1, \ldots, g_k \in \kk[x_1, x_1^{-1}, \ldots, x_n, x_n^{-1}$, for each $i$, define $I_{\scrQ_i}$ to be the ideal of $\kk[U_{\scrQ_i}]$ generated by $\talphanuu{\alpha_{\nu_i}}(g_j)$, $j = 1, \ldots, k$. \Cref{exercise:vonep-subscheme} implies that the ideals $I_{\scrQ_i}$ can be glued over their intersection to form a {\em sheaf of ideals} $\scrI$ of $\sheaf_{\xp}$. We write $\vonep(g_1, \ldots, g_k)$ for the corresponding closed subscheme $V(\scrI)$ of $\xp$. Assertion \eqref{alpha-non-dependence} of \cref{exercise:vonep-subscheme} implies that $\vonep(g_1, \ldots, g_k)$ does not depend on the choice of the $\alpha_{\nu_i}$.

\begin{example} \label{example:different-vonep}
\Cref{exercise:vonep-Cartier} implies that if $k = 1$, then $\vonep(g_1, \ldots, g_k)$ depends only on the ideal $I$ of $\kk[x_1, x_1^{-1}, \ldots, x_n, x_n^{-1}]$ generated by $g_1, \ldots, g_k$. We now show that this is in general false if $k > 1$. Let $\scrP$ be the triangle in $\rr^2$ with vertices $(0,0), (1,0), (0,1)$, so that $\xp \cong \pp^2$ (\cref{example:xp^n}), and with respect to corresponding homogeneous coordinates $[z_{0,0}:z_{1,0}:z_{0,1}]$ on $\pp^2$, $\xzerop = \pp^2 \setminus V(z_{0,0}z_{1,0}z_{0,1})$ and $\xonep = \xp \setminus \{[0:0:1], [0:1:0], [1:0:0]\}$. Note that $(x,y) := (z_{1,0}/z_{0,0}, z_{0,1}/z_{0,0})$ is a system of coordinates on $\xzerop$. Let $g_1 := x - 1$ and $g_2 := y - 1$. It is straightforward to check that $\vonep(g_1, g_2) \cap (\xonep \setminus \xzerop) = \emptyset$, and as a subscheme of $\xzerop$, $\vonep(g_1, g_2) = V(g_1, g_2)$. Now let $h_1 := x - y$ and $h_2 := x + y - 2 + (x-y)^2$, so that the ideal generated by $h_1, h_2$ in $\kk[x,x^{-1}, y, y^{-1}]$ is the same as the ideal generated by $g_1, g_2$. It can be checked that $\vonep(h_1, h_2)$ contains the point $[0:1:1] \in  \xonep \setminus \xzerop$, so that $\vonep(h_1, h_2) \neq \vonep(g_1, g_2)$.
\end{example}

Let $\scrQ$ be a facet of $\scrP$ with primitive inner normal $\nu$, and $\alpha_\nu$ be an arbitrary element of $\zz^n$ such that $\langle \nu, \alpha_\nu \rangle = 1$. \Cref{U_Q} implies that the ideal of $\kk[U_\scrQ]$ consisting of elements vanishing on $\orbit{\scrQ}$ is generated by $x^{\alpha_\nu}$. Therefore $\orbit{\scrQ}$ is precisely the support of the closed {\em subscheme} $Z_{\alpha_\nu} := V(x^{\alpha_\nu})$ of $U_\scrQ$. Since $\orbit{\scrQ}$ is Zariski closed in $\xonep$ (\cref{exercise:orbit-closed-in-xonep}), it follows that $Z_{\alpha_\nu}$ is in fact a closed subscheme of $\xonep$. We now determine the {\em embedded isomorphism} (see \cref{scheme:compactification}) type of the ``scheme-theoretic intersection'' $\vonep(g_1, \ldots, g_k) \cap Z_{\alpha_\nu}$. Let $\psi_\nu^*: \nktoruss{n-1} \cong \orbit{\scrQ}$ be as in \cref{O_Q}.

\begin{prop} \label{vonepinfinity}
As a closed subscheme of $\orbit{\scrQ}$, the scheme-theoretic intersection $\vonep(g_1, \ldots, g_k) \cap Z_{\alpha_\nu}$ is embedded isomorphic via $\psi_\nu^*$ to the closed subscheme $V(\Inbetapsinu(g_1), \ldots, \Inbetapsinu(g_k))$ of $\nktoruss{n-1}$.
\end{prop}

\begin{proof}
The ideal $I_\scrQ$ of $\vonep(g_1, \ldots, g_k) \cap Z_{\alpha_\nu}$ is generated in $\kk[U_\scrQ]$ by $x^{\alpha_\nu}, \talphanuu{\alpha_\nu}(g_1), \ldots, \talphanuu{\alpha_\nu}(g_k)$. \Cref{U_Q,O_Q} imply that $\kk[U_\scrQ]/\langle x^{\alpha_\nu},\talphanuu{\alpha_\nu}(g_1), \ldots, \talphanuu{\alpha_\nu}(g_k) \rangle$ is isomorphic via $\psi_\nu^*$ to $\kk[\nktoruss{n-1}]/\langle \Inbetapsinu(g_1), \ldots, \Inbetapsinu(g_k)\rangle$, which directly implies the result.
\end{proof}

\subsection{Exercises}

\begin{exercise} \label{exercise:in-extension}
Let $\scrP$ be an $n$-dimensional convex integral polytope in $\rr^n$ and $g$ be a Laurent polynomial with $\np(g) \subseteq \scrP$.
\begin{enumerate}
\item Show that the collection $(U_\alpha \cap \xp, x^{-\alpha}g)$, where $\alpha$ varies over the vertices of $\scrP$, defines a Cartier divisor $D$ on $\xp$ such that $D$ extends the Cartier divisor $V(g)$ of $\nktorus$. [Hint: use \cref{exercise:invertible-on-intersection}.]
\item Show that the following are equivalent:
\begin{enumerate}
\item $\supp(D)$ is the closure in $\xp$ of $V(g) \subset \nktorus$.
\item $\np(g) \cap \scrQ \neq \emptyset$ for each facet $\scrQ$ of $\scrP$.
\end{enumerate}
\item \label{in-extension:k^n} Assume in addition that $\scrA := \scrP \cap \zz^n$ satisfies the hypothesis of assertion \eqref{toric-k^n:origin} of \cref{toric-k^n} so that $U_\origin \cap \xp \cong \kk^n$. Show that the following are equivalent:
\begin{enumerate}
\item $\supp(D)$ is the closure in $\xp$ of $V(g) \subset \kk^n$.
\item $\np(g) \cap \scrQ \neq \emptyset$ for each facet $\scrQ$ of $\scrP$ which is not contained in any coordinate hyperplane of $\rr^n$.
\end{enumerate}
\end{enumerate}
\end{exercise}

\begin{exercise} \label{exercise:simplicial-extension}
Let $\scrP$ be an $n$-dimensional convex integral polytope in $\rr^n$. For each vertex $\alpha$ of $\scrP$, let $\scrC_\alpha$ be the cone in $\rr^n$ generated by $\{\alpha' - \alpha: \alpha' \in \scrP\}$. Assume there is $g \in \kk[x_1, x_1^{-1}, \ldots, x_n, x_n^{-1}]$ such that for each vertex $\alpha$ of $\scrP$, there is $\beta_\alpha \in \zz^n$ such that $\beta_\alpha + \np(g) \subset \scrC_\alpha$ and $\beta_\alpha + \np(g)$ touches every edge of $\scrC_\alpha$.
\begin{enumerate}
\item Let $\alpha, \alpha'$ be vertices of $\scrP$, and $\scrQ$ be the smallest face of $\scrP$ containing both $\alpha$ and $\alpha'$. Let $H$ be the linear subspace of $\rr^n$ spanned by $\{\beta - \alpha: \beta \in \scrQ\}$. Show that $\beta_\alpha - \beta_{\alpha'} \in H$. [Hint: use \cref{exercise:edge-graph} to reduce to the case that there is an edge of $\scrP$ connecting $\alpha$ and $\alpha'$. In that case show that $\beta_\alpha - \beta_{\alpha'}$ is on the line through the origin and $\alpha - \alpha'$.]
\item Deduce that the collection of pairs $(U_\alpha \cap \xp, x^{\beta_\alpha}g)$, where $\alpha$ varies over the vertices of $\scrP$, defines a Cartier divisor $D$ on $\xp$ such that $D$ extends the Cartier divisor $V(g)$ of $\nktorus$, and $\supp(D)$ is the closure in $\xp$ of $V(g) \subset \nktorus$. [Hint: use \cref{exercise:invertible-on-intersection}.]
\item If $\scrP$ is {\em simplicial}, i.e.\ each vertex of $\scrP$ is connected to precisely $n$ distinct edges, then show that for every $g \in \kk[x_1, x_1^{-1}, \ldots, x_n, x_n^{-1}]$, there is $m \geq 1$ such that $V(g^m) \subset \nktorus$ extends to a Cartier divisor $D$ on $\xp$ such that $\supp(D)$ is the closure in $\xp$ of $V(g^m) \subset \nktorus$.
\end{enumerate}
\end{exercise}

\begin{exercise} \label{exercise:non-2-extension}
Let $\scrP$ and $h$ be as in \cref{example:2-extension}. If $D$ is any Cartier divisor on $\xp$ which extends $V(h^m) \subset \nktoruss{2}$, where $m$ is a positive integer, then show that $D$ satisfies property \ref{tight-extension} if and only if $m$ is a multiple of $3$.
\end{exercise}

\begin{exercise} \label{exercise:non-3-extension}
Prove the claims made in \cref{example:3-non-extension}.
\end{exercise}

\begin{exercise} \label{exercise:orbit-closed-in-xonep}
Let $\scrQ$ be a facet of an $n$-dimensional convex integral polytope $\scrP$ in $\rr^n$. Show that $\orbit{\scrQ}$ is a closed subvariety of $\xonep$.
\end{exercise}

\begin{exercise} \label{exercise:vonep-Cartier}
Let $g \in \kk[x_1, x_1^{-1}, \ldots, x_n, x_n^{-1}]$ and $\scrP$ be a convex integral polytope of dimension $n$. Show that
\begin{enumerate}
\item The Cartier divisor $\vonep(g)$ of $\xonep$ does not depend on the choice of the $\alpha_{\nu_i}$.
\item If $g/h$ is a monomial in $x_1, \ldots, x_n$, then $\vonep(g) = \vonep(h)$.
\item $\supp(\vonep(g))$ is the closure in $\xonep$ of $V(g) \subset \nktorus$.
\end{enumerate}
\end{exercise}

\begin{exercise} \label{exercise:vonep-subscheme}
Let $g_1, \ldots, g_k \in \kk[x_1, x_1^{-1}, \ldots, x_n, x_n^{-1}]$ and $\scrP$ be a convex integral polytope of dimension $n$.
\begin{enumerate}
\item Let $\scrQ_1, \scrQ_2$ be facets of $\scrP$, $\nu_i$ be the primitive inner normal to $\scrQ_i$, and $\alpha_{\nu_i} \in \zz^n$ be such that $\langle \nu_i, \alpha_{\nu_i} \rangle = 1$. Show that for each $x \in U_{\scrQ_1} \cap U_{\scrQ_2}$ the ideal of $\local{\xp}{x}$ generated by $\talphanuu{\alpha_{\nu_1}}(g_1), \ldots, \talphanuu{\alpha_{\nu_1}}(g_k)$ is the same as the ideal generated by $\talphanuu{\alpha_{\nu_2}}(g_1), \ldots, \talphanuu{\alpha_{\nu_2}}(g_k)$.
\item \label{alpha-non-dependence} If $g_j/h_j$ is a monomial in $x_1, \ldots, x_n$ for each $j$, then show that for each $i$ and each $x \in U_{\scrQ_i}$, the ideal of $\local{\xp}{x}$ generated by $\talphanuu{\alpha_{\nu_i}}(g_1), \ldots, \talphanuu{\alpha_{\nu_i}}(g_k)$ is the same as the ideal generated by $\talphanuu{\alpha_{\nu_i}}(h_1), \ldots, \talphanuu{\alpha_{\nu_i}}(h_k)$.
\end{enumerate}
\end{exercise}

\begin{exercise} \label{exercise:different-vonep}
Verify the claims made in \cref{example:different-vonep}.
\end{exercise}

\section{Branches of curves on the torus} \label{branchion}
\subsection{Branch of a curve on a variety} \label{branchion-0}
Let $C$ be a {\em curve}, i.e.\ a variety of pure dimension one. Fix a desingularization $\pi: C' \to C$ of $C$ and a nonsingular compactification $\bar C'$ of $C'$. A \index{Branch}{\em branch} of $C$ is the germ of a point in $\bar{C'}$. Equivalently, consider the equivalence relation $\sim$ on the collection of pairs $\{(Z,z): z \in \bar{C'}$ and $Z$ is an open neighborhood of $z$ in $\bar{C'}\}$ defined as follows: $(Z, z) \sim (Z',z')$ if and only if $z = z'$. Then a branch of $C$ is a equivalence class of $\sim$. Let $X$ be a variety containing $C$ as a subvariety. Let $\bar X$ be an arbitrary projective compactification of $X$ and $\bar C$ be the closure of $C$ in $\bar X$. Then $\pi$ extends to a map $\bar C' \to \bar C$ (\cref{prop:curve-nonsingular-morphextension}), which we also denote by $\pi$. If $B := (Z,z)$ is a branch of $C$ and $y := \pi(z)$, we say that $y$ is the \index{Center of a branch}{\em center} of $B$ on $\bar X$, or equivalently, $B$ is a {\em branch of $C$ at $y$}. If $y \not\in X$, we say that (with respect to $X$) $B$ is a \index{Branch!at infinity}\index{Branch!centered at infinity}{\em branch at infinity}, or that it is {\em centered at infinity}. Since $\bar C'$ is nonsingular, $\local{Z}{z}$ is a discrete valuation ring (\cref{prop:curve-dvr}), and corresponds to a {\em unique} discrete valuation $\ord_z(\cdot)$ on the field of fractions of $\local{Z}{z}$ (assertion \ref{uniquely-dvr} of \cref{prop:discrete-properties}). If $f$ is a regular function on a neighborhood of $z$ on $X$, then we write $\ord_z(f)$ for $\ord_z(\pi^*(f))$.

\begin{example} \label{example:branch-order} 
Assume $\character{\kk} \neq 2$. Let $C$ be the curve on $\kk^2$ given by $x^2 = y^2 - y^3$. In \cref{example:order-singular} we have seen that $\pi: \kk \to C$ given by $t \mapsto (t - t^3, 1 - t^2)$ is a desingularization of $C$. Let $P = (0,0) \in C$. Then $\pi^{-1}(P) = \{1, -1\}$, and $(1, \kk)$ and $(-1, \kk)$ represent the two branches of $C$ at $P$ (see \Vref{fig:desingularize-node}). Note that both these branches are centered at infinity with respect to $\nktoruss{2}$ (since $P \not\in \nktoruss{2}$). In \cref{exercise:branch-order} you are asked to compute $\ord_1(f|_C)$ and $\ord_{-1}(f|_C)$ for different $f \in \kk[x,y]$. 
\end{example} 

\subsection{Weights of a branch on the torus} \label{branch-weight-section}
Fix a system of coordinates $(x_1, \ldots, x_n)$ on $\nktorus$. Let $B = (Z,z)$ be a branch of a curve on $\nktorus$. The \index{Weight corresponding to a branch}{\em weight} of $x_j$ corresponding to $B$ is $\ord_z(x_j|_Z)$. By $\nu_B$ we denote the element in $\rnstar$ with coordinates $(\ord_z(x_1), \ldots, \ord_z(x_n))$ with respect to the basis dual to the standard basis of $\rr^n$. The proof of the following result is left as an exercise. 

\begin{lemma} \label{branch-lemma-0}
\begin{enumerate}
\item  \label{dot-property} For each $\alpha \in \zz^n$, $x^\alpha$ restricts to a well-defined rational function on $Z$ and $\ord_z(x^{\alpha}) = \langle \nu_B, \alpha \rangle$.
\item $\nu_B \neq 0$ if and only if $B$ is centered at infinity with respect to $\nktorus$. \qed
\end{enumerate}
\end{lemma}


\begin{example}
Recall that the curve $C$ from \cref{example:branch-order} has two branches $B_1$ and $B_2$ at the origin. \Cref{exercise:branch-order} implies that both ``weight vectors'' $\nu_{B_1}$ and $\nu_{B_2}$ are $(1, 1)$. Note that
\begin{prooflist}
\item $B_j$ are centered at infinity with respect to $\nktoruss{2}$, and $\nu_{B_j}$ are nonzero. 
\item $\nu_B$ is {\em not} an invariant of $B$, but the embedding $C \cap \nktoruss{2} \into \nktoruss{2}$. \Cref{exercise:branch-nu-B'} presents a curve $C' \cong C$ via a map that sends the origin to itself, such that the branches at the origin correspond to {\em different} weight vectors. 
\end{prooflist}
\end{example} 

\subsection{Weighted order on Laurent polynomials} \label{wt-order-section}
Let $\nu$ be an integral element $\nu$ of $\rnstar$. The \index{Weighted!order}{\em weighted order} corresponding to $\nu$ is an integer valued map, which by an abuse of notation we also denote by $\nu$, on the ring of Laurent polynomials defined as follows: given a Laurent polynomial $f = \sum_\alpha c_\alpha x^\alpha$,
\begin{align*}
\nu(f) := \min\{ \langle \nu, \alpha \rangle: c_\alpha \neq 0\}
\end{align*}
In particular, if $f$ is the zero polynomial, then $\nu(f)$ is defined to be $\infty$. The \index{Initial!form!with respect to a weighted order}{\em initial form} $\In_\nu(f)$ of $f$ with respect to $\nu$ is the sum of all $c_\alpha x^{\alpha}$ such that $\langle \nu, \alpha \rangle= \nu(f)$. Given a subset $\scrS$ of $\rr^n$, we say that $f$ is \index{Supported at a set}{\em supported at} $\scrS$ if $\supp(f) \subseteq \scrS$, and we write $\scrL(\scrS)$ for the set of all Laurent polynomials supported at $\scrS$. In the case that $\scrS \cap \zz^n$ is a finite set, we equip $\scrL(\scrS)$ with the Zariski topology by identifying it with $\kk^{|\scrS \cap \zz^n|}$ via the map
\begin{align*}
\sum_{\alpha \in \scrS \cap \zz^n}  c_\alpha x^\alpha \mapsto (c_\alpha)_{\alpha \in \scrS \cap \zz^n}
\end{align*}
The result we will now prove ties these notions with those from \cref{branch-weight-section}; we will encounter many of its variants in the forthcoming chapters. Let $B = (Z,z)$ be a branch of a curve on $\nktorus$, and $\nu_B$ be as in \cref{branch-weight-section}. By an abuse of notation, we denote by $\nu_B$ also the weighted order corresponding to $\nu_B$.

\begin{prop}
Assume $\scrS \cap \zz^n$ has finitely many elements. There is a nonempty Zariski open subset $U$ of $\scrL(\scrS)$ such that $\ord_z(f) = \nu_B(f)$ for each $f \in U$. More precisely, define
\begin{align*}
\scrL^*(\scrS) := \{f \in \scrL(\scrS): \ord_z(f) = \nu_B(f)\}
\end{align*}
Then $\scrL^*(\scrS)$ is a constructible subset of $\scrL(\scrS)$ of dimension $|\scrS|$.
\end{prop}

\begin{proof}
Pick a parameter $\rho$ of $\local{Z}{z}$. Then $\ord_z(\rho) = 1$. For each $i = 1, \ldots, n$, if $m_i := \nu_B(x_i)$, then there is a representation of the form $x_i = c_i\rho^{m_i} + h_i$ where $c_i \in \kk^*$, and $h_i \in \local{Z}{z}$ such that $\ord_z(h_i) > m_i$ (assertion \eqref{curve-dvr:uu} of \cref{prop:curve-dvr}). For each $\scrA \subseteq \scrS \cap \zz^n$, let $\scrL^*_B(\scrA) :=  \{f \in \scrL^*(\scrA): \supp(\In_{\nu_B}(f)) \subset \In_{\nu_B}(\scrA) \}$. It is straightforward to check that
\begin{align*}
\scrL^*_B(\scrA)
	&=  \{( c_\alpha)_{\alpha \in \scrS \cap \zz^n}:
		c_\alpha = 0\ \text{if}\ \alpha \not\in \scrA,\
		\sum_{\alpha \in \In_{\nu_B}(\scrA)} \prod_{i=1}^n (c_i)^{\alpha_i} c_\alpha \neq 0 \}
\end{align*}
which implies that $\scrL^*_B(\scrA)$ is a {\em nonempty} open subset of a closed subset of $\scrL(\scrS)$. Since $\scrL^*(\scrS)$ is the union of $\scrL^*_B(\scrA)$ over all subsets $\scrA$ of $\scrS \cap \zz^n$, it follows that it is a constructible subset of $\scrL(\scrS)$. The remaining parts of the proposition follows from taking $\scrA = \scrS \cap \zz^n$.
\end{proof}

\subsection{Exercises}
\begin{exercise} \label{exercise:branch-order} 
Assume $\character{\kk} \neq 2$. Consider the desingularization $\pi: \kk \to C$ from \cref{example:branch-order} given by $t \mapsto (t - t^3, 1 - t^2)$. Show that 
\begin{enumerate}
\item $1 - t^2$ is a parameter at $\local{\kk}{z}$ for both $z = 1$ and $z = -1$.
\item $\ord_z(x|_C) = \ord_z(y|_C) = 1$ for both $z = 1$ and $z = -1$. 
\item $\ord_1((y - x)|_C) = 2$, but $\ord_{-1}((y - x)|_C) = 1$.  
\end{enumerate}
\end{exercise}

\begin{exercise} \label{exercise:branch-nu-B'} 
Assume $\character{\kk} \neq 2$. Let $C$ be the curve from \cref{example:branch-order}, $C' := V((x+y)^2 - y^2 + y^3) \subseteq \kk^2$, and $P' := (0, 0) \in C'$. 
\begin{enumerate}
\item the map $\phi: (x, y) \mapsto (x-y, y)$ induces an isomorphism $C \cong C'$ and maps $P \to P'$, where $P := (0,0) \in C$.
\item there are two branches $B'_1, B'_2$ of $C'$ at the origin and the corresponding ``weight vectors'' $\nu_{B'_1}, \nu_{B'_2} \in \rnstar$ are different. [Hint: use \cref{example:branch-order,exercise:branch-order}.]
\end{enumerate}
\end{exercise}

\begin{exercise} \label{exercise:branch-order-1}
For each of the following curves $C  \subseteq \kk^n$, show that $C$ has a single branch $B$ at the origin, and compute $\nu_B \in \rnstar$: 
\begin{enumerate}
\item $C = V(x^2 - y^3)$. [Hint: the desingularization of $C$ is given by $t \mapsto (t^3, t^2)$.]
\item $C = V(y^3 - x^4, z^3 - x^5) $. [Hint: use \cref{exercise:noncomplete-intersection}.]
\end{enumerate}
\end{exercise} 

\begin{exercise}
Prove \cref{branch-lemma-0}.
\end{exercise}

\section{Points at infinity on toric varieties} \label{toric-infinity-section}
\subsection{Centers of branches at infinity on the torus} \label{toric-center-section}
Let $\scrA:= \{\alpha_0, \ldots, \alpha_N\}$ be a finite subset of $\zz^n$ and $B = (Z,z)$ be a branch centered at infinity on $\nktorus$. If $\phi_\scrA: \nktorus \to \xzeroa$ is the map from \eqref{phi_A}, then $\phi_\scrA(B)$ is a branch centered at infinity on $\xzeroa$. We now determine the torus orbit of $\xa$ that contains the center $o_B$ of $\phi_\scrA(B)$ on $\xa$; we will see that this orbit is completely determined by the element $\nu_B \in \rnstar$ defined in \cref{branch-weight-section}. Let $\scrP$ be the convex hull of $\scrA$ and $\scrB$ be a face of $\scrA$. Then the convex hull $\scrQ$ of $\scrB$ is a face of $\scrP$. As in \cref{normal-section} we write $\Sigma_\scrP$ for the {\em normal fan} of $\scrP$, and denote the {\em normal cone} of $\scrQ$ by $\sigma_\scrQ$ and the relative interior of $\sigma_\scrQ$ by $\sigma^0_\scrQ$.

\begin{figure}[h]
\begin{center}
\def\xmin{-0.5}
\def\xmax{13.5}
\def\ymin{-5.5}
\def\ymax{5.5}
\def\opazero{0.5}
\def\colorzero{green}

\tikzstyle{dot} = [circle, minimum size=5pt, inner sep = 0pt, fill]

\begin{tikzpicture}[scale=\scalefactor]
\draw [gray,  line width=0pt] (\xmin, \ymin) grid (\xmax,\ymax);
\draw [<->] (0, \ymax) |- (\xmax, 0);
\node[dot, \colorA] (A) at (1,-4) {};
\node[dot, \colorB] (B) at (13,2) {};
\node[dot, \colorC] (C) at (10,5) {};
\fill[\colorzero, opacity=\opazero ] (A.center) --  (B.center) -- (C.center);
\draw[thick, \colorAB] (A) -- (B);
\draw[thick, \colorBC] (B) -- (C);
\draw[thick, \colorCA] (C) -- (A);
\node[anchor = north] at (A) {\picfontsize $A$};
\node[anchor = west] at (B) {\picfontsize $B$};
\node[anchor = south] at (C) {\picfontsize $C$};

\node at (8,1) {\picfontsize $\scrP$};

\def\xmin{-7.5}
\def\xshift{5}
\pgfmathsetmacro\shiftone{\xshift + \xmax -\xmin};
\def\xmax{7.5}

\begin{scope}[shift={(\shiftone,0)}]
\draw [gray,  line width=0pt] (\xmin, \ymin) grid (\xmax,\ymax);
\coordinate (O) at (0,0);

\coordinate (ABperp) at (-\ymax/2, \ymax);
\coordinate (BCperp) at (\ymin, \ymin);
\coordinate (CAperp) at (-\ymin, \ymin);

\fill[\colorB, opacity=\opazero ] (O) --  (ABperp) -- (\xmin,\ymax) -- (\xmin,\ymin) -- (BCperp);
\fill[\colorC, opacity=\opazero ] (O) --  (BCperp) -- (CAperp);
\fill[\colorA, opacity=\opazero ] (O) --  (CAperp) -- (\xmax, \ymin) -- (\xmax,\ymax) -- (ABperp);

\draw [thick, \colorAB] (O) -- (ABperp);
\draw [thick, \colorBC] (O) -- (BCperp);
\draw [thick, \colorCA] (O) -- (CAperp);


\node[dot, \colorzero] at (O) {};

\node at (-4.5,0.5) {\picfontsize $\sigma^0_B$};
\node at (4.5,1.5) {\picfontsize $\sigma^0_A$};
\node at (.5,-3.5) {\picfontsize $\sigma^0_C$};

\node [anchor= north] at (BCperp) {\picfontsize $\sigma^0_{BC}$};
\node [anchor=north] at (CAperp) {\picfontsize $\sigma^0_{CA}$};
\node [anchor=south] at (ABperp) {\picfontsize $\sigma^0_{AB}$};

\end{scope}
\end{tikzpicture}
\caption{Normal fan of $\scrP$}  \label{fig:xp-normal-fan}
\end{center}
\end{figure}

\begin{prop} \label{curve-to-orbit}
$o_B \in \orbit{\scrB}$ if and only if $\nu_B \in \sigma^0_\scrQ$. In particular, if $\scrB$ is a facet of $\scrA$, then $o_B \in \orbit{\scrB}$ if and only if $\nu_B$ is a positive multiple of the {\em primitive inner normal} to $\scrQ$.
\end{prop}

\begin{proof}
Pick $\beta \in \scrB$. \Cref{xa-thm} implies that $x^{\alpha - \beta} = z_\alpha/z_\beta$ is a regular function on $\orbit{\scrB}$ for each $\alpha \in \scrA$. It also implies that $o_B \in \orbit{\scrB}$ if and only if the following holds: ``$x^{\alpha-\beta}|_{o_B} \neq 0$ if and only if $\alpha \in \scrB$.'' Due to assertion \eqref{dot-property} of \cref{branch-lemma-0} the latter condition is  equivalent to the condition that $\In_{\nu_B}(\scrA) =  \scrB$, which is in turn equivalent to the condition that $\In_{\nu_B}(\scrP) = \scrQ$, as required.
\end{proof}

We now describe the coordinates of $o_B$. As in \cref{wt-order-section}, we write $\nu_B$ also for the {\em weighted order} on $\kk[x_1, x_1^{-1}, \ldots, x_n, x_n^{-1}]$ induced by $\nu_B$. Fix an arbitrary element $\rho_B \in \local{Z}{z}$ such that $\ord_z(\rho_B|_Z) = 1$. We say that $\rho_B$ is a \index{Parameter!of a branch of a curve}{\em parameter} of $B$. Define
\begin{align}
\begin{split}
\In_B(x_j)
	&:=
		\left. \frac{x_j}{(\rho_B)^{\nu_B(x_j)}} \right|_z \in \kk^*, \\
\In(B)
	&:= (\In_B(x_1), \ldots, \In_B(x_n)) \in \nktorus
\end{split}
\label{in-B-torus}
\end{align}
That $\In_B(x_j)$ and $\In(B)$ are well defined follows from \cref{prop:curve-dvr}. Note that $\In(B)$ depends on the choice of $\rho_B$. In all cases considered in this book, whenever a branch $B$ of a curve is considered, a corresponding parameter $\rho_B$ is assumed to be fixed from the beginning. Let $\scrB$ be the face of $\scrA$ such that $o_B \in \orbit{\scrB}$ (\cref{curve-to-orbit} shows that $\scrB$ is uniquely determined by $\nu_B$). Let $\phi_\scrB : \nktorus \to \xzerob$ be the map from \eqref{phi_A}. \Cref{xa-thm} implies that we may think of $\phi_\scrB$ as a map from $\nktorus$ to $\orbit{\scrB} \subset \xa$, simply by adjoining a zero in place of each coordinate $z_\alpha$ on $\pp^{N}$ such that $\alpha \not\in \scrB$.

\begin{prop} \label{branch-center}
$o_B = \phi_\scrB(\In(B))$. In particular, $\phi_\scrB(\In(B))$ does not depend on the choice of $\rho_B$ (even though $\In(B)$ does).
\end{prop}

\begin{proof}
Pick $\alpha, \beta \in \scrB$. \Cref{curve-to-orbit} implies that $(z_\alpha/z_\beta)|_{o_B} = x^{\alpha - \beta}|_{o_B} = (\In(B))^{\alpha - \beta}$. The result then follows immediately from \cref{xa-thm}.
\end{proof}

\subsection{Closure of subvarieties of the torus}
Let $W$ be a closed subvariety of $\nktorus$ defined by Laurent polynomials $f_1, \ldots, f_m$ in $\nktorus$. Let $\scrA$ be a finite subset of $\zz^n$ and $\phi_\scrA: \nktorus \to \xzeroa$ be the map from \eqref{phi_A}. Write $\bar W'$ for the closure in $\xa$ of $W' := \phi_\scrA(W) \subset \xzeroa$. In this section we give a partial description of the points in $\bar W'$. \\


\begin{lemma} \label{branch-lemma-1}
Let $B$ be a branch of a curve contained in $W$. Then $\In (B)$ is a common zero of $\In_{\nu_B}(f_i)$, $i = 1, \ldots, m$.
\end{lemma}

\begin{proof}
Let $B = (Z,z)$. Pick a parameter $\rho_B$ of $B$. \Cref{prop:curve-dvr} implies that for each $j$, $x_j/\rho_B^{\nu_B(x_j)}$ is a regular function on a neighborhood of $z$ on $Z$, and it can be expressed as $\In_B(x_j) + g_j$, where $g_j \in \local{Z}{z}$, $\ord_z(g_j) > 0$. Consequently, for each $i$, $f_i/\rho_B^{\nu_B(f_i)}$ can be expressed in $\local{Z}{z}$ as $\In_{\nu_B}(f_i)(\In(B)) + h_i$, where $h_i \in \local{Z}{z}$, $\ord_z(h_i) > 0$. Since $f_i/\rho_B^{\nu_B(f_i)}$ maps to the zero element in $\local{Z}{z}$, it follows that $\In_{\nu_B}(f_i)(\In(B))= 0$, as required.
\end{proof}

Let $\nu$ be an integral element of $\rnstar$; we write $\nu$ also for the corresponding weighted order on $\kk[x_1,  x_1^{-1},\ldots, x_n, x_n^{-1}]$ and define $W_\nu(f_1, \ldots, f_m) := V(\In_\nu(f_1), \ldots, \In_\nu(f_m)) \subset \nktorus$.
Let $\scrB$ be a face of $\scrA$. As in \cref{branch-center}, we regard the map $\phi_\scrB : \nktorus \to \xzerob$ from \eqref{phi_A} as a map from $\nktorus$ to $\orbit{\scrB} \subset \xa$. As in \cref{curve-to-orbit} we write $\scrQ$ for the convex hull of $\scrB$ and $\sigma^0_\scrQ$ for the relative interior of the corresponding cone of the normal fan of the convex hull of $\scrA$.

\begin{cor} \label{closure-containment}
$\bar W' \cap  \orbit{\scrB} \subset \bigcup_{\nu \in \sigma^0_\scrQ} \phi_\scrB(W_\nu(f_1, \ldots, f_m))$.
\end{cor}

\begin{proof}
Let $w \in \bar W' \cap \orbit{\scrB}$. If $w \in W'$, then we must have that $\scrB = \scrA$. In that case $0 \in \sigma^0_\scrQ$. Since $W_0(f_1, \ldots, f_m) =  W$, it follows that $w \in \phi_\scrB(W_0(f_1, \ldots, f_m)) = W'$, as required. So assume $w \in \bar W' \setminus W'$. Then there is an irreducible curve $C' \subset W'$ such that $w$ is in the closure of $C'$ (\cref{closure-curve-prop}). Pick a branch $B = (Z,z)$ of $\phi_\scrA^{-1}(C)$ such that $z \mapsto w$ under the morphism induced by $\phi_\scrA$. \Cref{curve-to-orbit} implies that $\nu_B \in \sigma^0_\scrQ$ and \cref{branch-center} implies that $w = \phi_\scrB(\In(B))$. Since $\In(B) \in W_{\nu_B}(f_1, \ldots, f_m)$ (\cref{branch-lemma-1}), the result follows.
\end{proof}

For each $i = 1, \ldots, m$, there are only finitely many possibilities for $\In_\nu(f_i)$ as $\nu$ varies over $\rnstar$. It follows that the union in the statement of \cref{closure-containment} can be regarded as being over a {\em finite} collection of $\nu \in \rnstar$. \Cref{exercise:closure-proper-containment} shows that the containment of \cref{closure-containment} is in general proper. However, if $W$ is a hypersurface (i.e.\ $m = 1$), then \cref{exercise:closure-hypersurface} shows that \cref{closure-containment} holds with ``$=$'' in place of $\subset$.

\subsection{Exercises}

\begin{exercise} \label{exercise:closure-proper-containment}
Let $\scrA := \{(0,0), (1,0), (0,1)\} \subset \rr^2$ and $h_1, h_2$ be as in \cref{example:different-vonep}. Let $W := V(h_1, h_2) = \{(1,1)\} \in \nktoruss{2}$ and $\scrB := \{(1,0), (0,1)\}$. Note that $\scrQ := \conv(\scrB)$ is an edge of $\scrP := \conv(\scrA)$.
\begin{enumerate}
\item Show that $\bar W' \cap  \orbit{\scrB} = \emptyset$.
\item Let $\nu$ be the primitive inner normal to $\scrQ$. Show that $W_\nu(h_1, h_2) \neq \emptyset$.
\end{enumerate}
\end{exercise}

\begin{exercise} \label{exercise:closure-hypersurface}
Let $f \in \kk[x_1, x_1^{-1}, \ldots, x_n, x_n^{-1}]$ and $W := V(f) \subset \nktorus$. Show that in the notation of \cref{closure-containment}, $\bar W' \cap  \orbit{\scrB} = \bigcup_{\nu \in \sigma^0_\scrQ} \phi_\scrB(W_\nu(f))$ for each face $\scrB$ of $\scrA$. [Hint: use \cref{exercise:vonep-Cartier} to prove it in the case that $\scrQ$ is a facet of $\scrP$. Then use induction on $\dim(\scrP)$.]
\end{exercise}

\section{$^*$Weighted projective spaces} \label{weighted-toric-section}
\footnote{The asterisk in the section name is to indicate that the material of this section is not going to be used in the proof of Bernstein's theorem. It is used for the first time in \cref{bezout-chapter}.} Recall that the $n$-dimensional projective space is the space of straight lines through the origin in $\kk^{n+1}$. A \index{Weighted!projective space}{\em weighted projective space} is constructed in the same way, using {\em weighted rational curves} in place of straight lines. Let $\omega$ be an integral element of $\rnnstar{n+1}$ with coordinates $(\omega_0, \ldots, \omega_n)$ with respect to the basis dual to the standard basis of $\rr^{n+1}$. Assume each $\omega_j$ is {\em positive}.
The corresponding weighted projective space, which we denote by $\wpn$ or $\wwpn{\omega_0, \ldots, \omega_n}$, is the set of curves in $\kk^{n+1}$ of the form $C_a := \{(a_0t^{\omega_0}, \ldots, a_nt^{\omega_n}): t \in \kk\}$, where $a  := (a_0, \ldots, a_n) \in \kk^{n+1} \setminus \{0\}$. The \index{Weighted!homogeneous!coordinates}{\em weighted homogeneous coordinates} of $C_a$ are $[a_0: \cdots : a_n]$. Note that the projective space $\pp^n$ is the special case of $\wpn$ for $\omega = (1, \ldots, 1)$. Like $\pp^n$, each $\wpn$ can be given the structure of a complete algebraic variety. In this section we give two
(equivalent) realizations of $\wpn$ as a toric variety.

\begin{figure}[h]
\begin{center}
\def\xmin{-0.5}
\def\xmax{13.5}
\def\ymin{-5.5}
\def\ymax{5.5}
\def\opazero{0.5}
\def\colorzero{green}

\tikzstyle{dot} = [circle, minimum size=5pt, inner sep = 0pt, fill]

\begin{tikzpicture}[scale=\scalefactor]
\draw [gray,  line width=0pt] (\xmin, \ymin) grid (\xmax,\ymax);
\draw [<->] (0, \ymax) |- (\xmax, 0);
\node[dot, \colorA] (A) at (1,-4) {};
\node[dot, \colorB] (B) at (13,2) {};
\node[dot, \colorC] (C) at (10,5) {};
\fill[\colorzero, opacity=\opazero ] (A.center) --  (B.center) -- (C.center);
\draw[thick, \colorAB] (A) -- (B);
\draw[thick, \colorBC] (B) -- (C);
\draw[thick, \colorCA] (C) -- (A);
\node[anchor = north] at (A) {\picfontsize $A$};
\node[anchor = west] at (B) {\picfontsize $B$};
\node[anchor = south] at (C) {\picfontsize $C$};

\node at (8,1) {\picfontsize $\scrP$};

\def\xmin{-7.5}
\def\xshift{5}
\pgfmathsetmacro\shiftone{\xshift + \xmax -\xmin};
\def\xmax{7.5}

\begin{scope}[shift={(\shiftone,0)}]

\coordinate (O) at (0,0);
\coordinate (ABperp) at (-1, 2);
\coordinate (BCperp) at (-1, -1);
\coordinate (CAperp) at (1, -1);

\draw [thick, \colorAB, ->] (O) -- (ABperp);
\draw [thick, \colorBC, ->] (O) -- (BCperp);
\draw [thick, \colorCA, ->] (O) -- (CAperp);

\node [anchor = north east] at (BCperp) {\picfontsize $\nu_{BC} = (-1,-1)$};
\node [anchor = north west] at (CAperp) {\picfontsize $\nu_{CA} = (1,-1)$};
\node [anchor = south] at (ABperp) {\picfontsize $\nu_{AB} = (-1,2)$};

\node [anchor = north] at (0,-4) {\picfontsize $\nu_{BC} + 2\nu_{AB} + 3\nu_{BC} = 0$};
\end{scope}
\end{tikzpicture}
\caption{$\xp \cong \pp^2(1,2,3)$}  \label{fig:wt-proj-n-construction}
\end{center}
\end{figure}

\subsection{$\wpn$ via polytopes in $\rr^n$} \label{implicit-wpn}
Pick integral elements $\nu_0, \ldots, \nu_n$ of $\rnstar$ such that $\nu_0, \ldots, \nu_n$ span $\znstar$, and $\sum_{j=0}^n \omega_j \nu_j = 0$. Let $\scrP$ be an $n$-dimensional integral simplex in $\rr^n$ such that its inner facet normals are $\nu_0, \ldots, \nu_n$; note that $\scrP$ is uniquely determined by the $\nu_j$ up to translation and scaling - see \cref{fig:wt-proj-n-construction}. We will show that $\wpn$ can be identified with $\xp$.

\begin{center}
\begin{figure}[h]

\def\shiftone{7.5}
\def\colorzero{green}
\def\colorone{red}
\def\opazero{0.5}
\def\viewx{75}
\def\viewy{30}
\def\titlex{1}
\def\titley{-1}
\begin{tikzpicture}[scale=0.6]
\pgfplotsset{every axis/.append style = {view={\viewx}{\viewy}, axis lines=middle, enlargelimits={upper}}}

\begin{axis}
\addplot3[ultra thick, draw=\colorone, fill=\colorzero,opacity=\opazero] coordinates{(6,0,0) (0,3,0) (0,0,2) (6,0,0)};
\draw (axis cs:3,0.5,1) node {$\scrP'$};
\end{axis}
\end{tikzpicture}

\caption{A polytope $\scrP'$ in $\rr^3$ such that $\xpp{\scrP'} \cong \pp^2(1,2,3)$} \label{fig:wt-proj-n+1-construction}
\end{figure}
\end{center}

\subsection{$\wpn$ via polytopes in $\rr^{n+1}$} \label{explicit-wpn}
Let $p := \lcm(\omega_0, \ldots, \omega_n)$ and $\scrP'$ be the $n$-dimensional simplex in $\rr^{n+1}$ with vertices $\beta_j := (p/\omega_j)e_j$, $j = 0, \ldots, n$, where $e_0, \ldots, e_n$ are the standard unit vectors in $\rr^{n+1}$ (\cref{fig:wt-proj-n+1-construction}). We will show that $\xpp{\scrP'} \cong \xp$, where $\scrP$ is from \cref{implicit-wpn}.

\subsection{Equivalence of the constructions} \label{wt-equivalence-section}
Let $\scrP$ be as in \cref{implicit-wpn}. For each $j = 0, \ldots, n$, let $\scrQ_j$ be the facet of $\scrP$ with inner normal $\nu_j$, and $\alpha_j$ be the (unique) vertex of $\scrP$ which is not on the facet $\scrQ_j$. Pick $j$, $0 \leq j \leq n$. \Cref{xp-thm} implies that the coordinate ring of $U_{\alpha_j} \cap \xp$ is $\kk[S_j]$, where $S_j$ is the semigroup of integral points in the polyhedral cone $\scrC_{\alpha_j} \subset \rr^n$ generated by $\alpha_i - \alpha_j$, $i = 0, \ldots, n$. \Cref{exercise:equivalent-cones} below implies that $\alpha \in S_j$ if and only if $\alpha \in \zz^n$ and $\langle \nu_i, \alpha \rangle \geq 0$ for each $i = 0, \ldots, \hat j, \ldots, n$. Consider the map $\phi: \zz^n \mapsto \zz^{n+1}$ given by $\alpha \mapsto (\langle \nu_0, \alpha \rangle, \ldots, \langle \nu_n, \alpha \rangle)$.

\begin{claim}
$\phi$ induces an isomorphism between $S_j$ and $S'_j := \{\beta = (\beta_0, \ldots, \beta_n) \in \zz^{n+1}: \langle \omega, \beta \rangle = 0,\ \beta_i \geq 0$ for each $i = 0, \ldots, \hat j,\ldots,  n\}$.
\end{claim}

\begin{proof}
\Cref{exercise:equivalent-cones} implies that $S_j = \phi^{-1}(S'_j)$. We now show that $S'_j = \phi(S_j)$. Indeed, let $H_j$ and $H'_j$ be the subgroups of $\zz^{n+1}$ generated respectively by $\phi(S_j)$ and $S'_j$. It suffices to show that $H'_j = H_j$. Since $H_j \subset H'_j$ are subgroups of $\zz^n$ of the same rank $n$, we have to show that if $k\beta \in H_j$ for some positive integer $k$ and $\beta = (\beta_0, \ldots, \beta_n) \in \zz^{n+1}$, then $\beta \in H_j$. Indeed, if $k\beta_j = \langle \nu_j, \alpha \rangle$ for each $j = 0, \ldots, n$, then since the $\nu_j$ span $\znstar$, it follows that $\alpha/k \in \zz^n$, and $\beta = \phi(\alpha/k) \in H_j$. Therefore $H'_j = H_j$, which proves the claim.
\end{proof}

As in \cref{explicit-wpn} let $\beta_j$ be the vertex of $\scrP'$ on the $j$-th axis. \Cref{exercise:equivalent-cones} implies that $S'_j$ is the semigroup of integral points in the polyhedral cone $\scrC'_{\beta_j} \subset \rr^{n+1}$ generated by $\beta_i - \beta_j$, $i = 0, \ldots, n$. Therefore \cref{xp-thm} implies that $\phi$ induces an isomorphism $\xp \cong \xpp{\scrP'}$.

\subsection{Identification with $\wpn$} \label{wt-identification-section}
Let $f \in \kk[x_0, \ldots, x_n]$. We say that $f$ is \index{Weighted!homogeneous!polynomial}{\em weighted homogeneous} with respect to $\omega$ (or in short, {\em $\omega$-homogeneous}) if $\omega$ is constant on $\supp(f)$. If $f$ is $\omega$-homogeneous, then the set $V(f) := \{[a_0 : \cdots : a_n]: f(a_0, \ldots, a_n) = 0\} \subset \wpn$ of zeroes of $f$ is a well defined subset of $\wpn$. As in the case of $\pp^n$, the basic open subsets of $\wpn$ are $U_j := \wpn \setminus V(x_j)$, $j = 0, \ldots, n$. If $f$ is $\omega$-homogeneous with $\omega(f)$ a multiple of $\omega_j$, say $\omega(f) = k\omega_j$, $k \geq 0$, then $(f/x_j^k)|_{C_a}$ is constant for all $C_a \in U_j$ and therefore $f/x_j^k$ is a well defined function on $U_j$. \Cref{exercise:R_j-S_j} shows that the $\kk$-algebra $R_j$ generated by all these $f/x_j^k$, $k \geq 0$, is finitely generated, and if $h_1, \ldots, h_s$ generate $R_j$ as a $\kk$-algebra, then they induce a bijection from $U_j$ to an open affine subvariety of $\xpp{\scrP'}$ (where $\scrP'$ is as in \cref{explicit-wpn}) which extends to a bijection from $\wpn \to \xpp{\scrP'}$. We use this bijection to identify $\wpn$ with $\xpp{\scrP'} \cong \xp$ (where $\scrP$ is as in \cref{implicit-wpn}). In particular, it follows that each $U_j$ is an open affine subvariety of $\wpn$ with coordinate ring $R_j$. To completely describe the identification of $\wpn$ and $\xpp{\scrP'}$ it remains to explicitly identify points in $\wpn$ with those of $\xpp{\scrP'}$. Given $a := [a_0: \cdots : a_n] \in \wpn$, we now compute the ideal $I_a$ of all $\omega$-homogeneous polynomials in $(x_0, \ldots, x_n)$ that vanish at $a$. Let $J_a$ be the ideal of $\kk[x_0, \ldots, x_n]$ generated by all the $x_i$ such that $a_i = 0$, and all binomials of the form $a^{\alpha_2}x^{\alpha_1} - a^{\alpha_1}x^{\alpha_2}$ where $\alpha_1, \alpha_2 \in \zzeroo{n+1}$ such that $\langle \omega, \alpha_1 \rangle = \langle \omega, \alpha_2 \rangle$. It is clear that $J_a \subset I_a$. The following proposition, which you will prove in \cref{exercise:ideal-prop}, shows that the converse is also true.

\begin{prop} \label{ideal-prop}
$I_a = J_a$. \qed
\end{prop}

\subsection{Exercises}

\begin{exercise}\label{exercise:equivalent-cones}
In the notation of \cref{wt-equivalence-section} show that
\begin{enumerate}
\item $\scrC_{\alpha_j} = \{\alpha \in \rr^n : \langle \nu_i, \alpha \rangle \geq 0,\ i = 0, \ldots, \hat j, \ldots, n\}$.
\item $\scrC'_{\beta_j} = \{\beta \in \rr^{n+1} : \langle \omega, \beta \rangle = 0,\ \beta_i \geq 0,\ i = 0, \ldots, \hat j, \ldots, n\}$.
\end{enumerate}
[Hint: use \cref{prop:minimal-polynequalities}.]
\end{exercise}

\begin{exercise}\label{exercise:R_j-S_j}
For each $j = 0, \ldots, n$, let $R_j$ be as in \cref{wt-identification-section} and $S'_j$ be as in \cref{wt-equivalence-section}.
\begin{enumerate}
\item Show that $R_j = \kk[S'_j]$ for each $j$.
\item Let $U'_j := U_{\beta_j} \cap \xpp{\scrP'}$ be the affine open subset of $\xpp{\scrP'}$ (where $\scrP'$ is as in \cref{explicit-wpn}) with coordinate ring $\kk[S'_j]$. If $h_1, \ldots, h_s$ generate $R_j$ as a $\kk$-algebra, then show that the map $x \mapsto (h_1(x), \ldots, h_s(x))$ induces a bijection $\phi_j$ from $U_j := \wpn \setminus V(x_j)$ to $U'_j$.
\item Show that for distinct $j$, the maps $\phi_j$ are ``compatible,'' i.e.\ for each $j,j'$, $\phi_j$ maps $U_j \cap U_{j'}$ bijectively onto $U'_j \cap U'_{j'}$.
\end{enumerate}
\end{exercise}

\begin{exercise}\label{exercise:ideal-prop}
Prove \cref{ideal-prop}. [Hint: use arguments analogous to those outlined in \cref{exercise:xa-binomial}.]
\end{exercise}

\section{$^*$Weighted blow up} \label{weighted-blow-up-section}
\footnote{The asterisk in the section name is to indicate that the material of this section is not going to be used in the proof of Bernstein's theorem. It is used for the first time in \cref{multiplicity-chapter}.}Let $\nu$ be a weighted order on $\kk[x_1, x_1^{-1}, \ldots, x_n, x_n^{-1}]$ with {\em positive} weights $\nu_j := \nu(x_j)$, $j = 1, \ldots, n$. Identify $\nu$ with the integral element of $\rnstar$ with coordinates $(\nu_1, \ldots, \nu_n)$ with respect to the basis dual to the standard basis of $\rr^n$. Fix a positive integer $k$. Let $\scrB_k := \{\alpha \in \zzeroo{n}: \langle \nu, \alpha \rangle = kp\}$, where $p := \lcm(\nu_1, \ldots, \nu_n)$, and $\qqq_k$ be the ideal of $\kk[x_1, \ldots, x_n]$ generated by $\{x^\alpha: \alpha \in \scrB_k\}$. Recall that the {\em blow up} $\bl_{\qqq_k}(\kk^n)$ of $\kk^n$ at $\qqq_k$ is the closure in $\kk^n \times \xaa{\scrB_k}$ of the graph of the map $\phi_{\scrB_k}: \nktorus \to \xaa{\scrB_k}$ defined as 
\begin{align*}
x \mapsto [x^{\alpha_0}: \cdots: x^{\alpha_N}]
\end{align*}
where $\alpha_0, \ldots, \alpha_N$ are the elements of $\scrB_k$. Let $\scrP := \{\alpha \in \rzeroo{n}: \langle \nu, \alpha \rangle = p\}$. Then $\scrP$ is the $(n-1)$-dimensional simplex in $\rr^n$ with vertices $(p/\nu_j)e_j$, $j = 1, \ldots, n$, where $e_1, \ldots, e_n$ form the standard basis of $\rr^n$. Let $\scrQ$ be the $n$-dimensional simplex in $\rr^n$ whose vertices are $0,e_1, \ldots, e_n$, so that $\kk^n$ can be naturally identified with the basic open set $U_0$ of $\xq \cong \pp^n$. If $k$ is sufficiently large, then $\xaa{\scrB_k} \cong \xp$, so that assertion \eqref{xp:diagonal} of \cref{xp-thm} implies that $\bl_{\qqq_k}(\kk^n)$ is isomorphic to the open subset of $\xaa{\scrP + \scrQ}$ which is the union of basic open subsets $U_{(p/\nu_j)e_j}$, $j = 1, \ldots, n$ (note that each $(p/\nu_j)e_j$ is a vertex of $\scrP+ \scrQ$). In particular, $\bl_{\qqq_k}(\kk^n)$ are isomorphic for all sufficiently large $k$; we call the corresponding algebraic variety the \index{Weighted!blow up}\index{Blow up!weighted}{\em $\nu$-weighted blow up of $\kk^n$} and denote it by $\blnukn$. Let $\sigma: \blnukn \to \kk^n$ be the blow up map. The \index{Exceptional divisor}{\em exceptional divisor} on $\blnukn$ is $\Enu := \sigma^{-1}(0)$. Note that $\sigma$ is an isomorphism on $\blnukn \setminus \Enu$. If $V$ is a subvariety of $\kk^n$, the \index{Strict transform}{\em strict transform} of $V$ on $\blnukn$ is the closure in $\blnukn$ of $\sigma^{-1}(V \setminus \{0\})$. Note that $\scrP$ is the ``lower'' facet of $\scrP+ \scrQ$ (see \cref{fig:wt-blow-up}) and $\Enu$ is precisely the subvariety $\corbit{\scrP}$ of $\xaa{\scrP+ \scrQ}$ corresponding to the facet $\scrP$ of $\scrP+ \scrQ$. The construction in \cref{explicit-wpn} therefore shows that $\Enu \cong \wwpnn{\nu}{n-1}$. The following proposition makes this isomorphism more explicit.

\begin{center}
\begin{figure}[h]

\def\shiftx{7.5}
\def\colorzero{green}
\def\colorone{red}
\def\colortwo{blue}
\def\colorthree{yellow}
\def\opazero{0.5}
\def\viewx{60}
\def\viewy{30}
\def\titlex{2}
\def\titley{-1}
\begin{tikzpicture}[scale=0.66]
\pgfplotsset{every axis title/.append style={at={(0,-0.2)}}, view={\viewx}{\viewy}, axis lines=middle, enlargelimits={upper}}

\begin{axis}
\addplot3[ultra thick, draw, fill=\colorzero,opacity=\opazero] coordinates {(6,0,0) (0,3,0) (0,0,2) (6,0,0)};
\end{axis}
\draw (\titlex,\titley) node {$\scrP$};


\begin{scope}[shift={(\shiftx,0)}]
\begin{axis}
	\addplot3 [draw, ultra thick, fill=\colortwo,opacity=\opazero] coordinates{(0,0,1) (0,1,0) (1,0,0) (0,0,1)};
	\addplot3 [draw, ultra thick, fill=\colorone,opacity=\opazero] coordinates{(0,0,0) (1,0,0) (0,0,1) (0,0,0)};
\end{axis}
\draw (\titlex,\titley) node {$\scrQ$};
\end{scope}


\begin{scope}[shift={(2*\shiftx,0)}]
\begin{axis}
	\addplot3[ultra thick, draw, fill=\colorzero, opacity=\opazero] coordinates {(6,0,0) (0,3,0) (0,0,2) (6,0,0)};
	\addplot3[ultra thick, draw, fill=\colorthree, opacity=\opazero] coordinates {(0,0,3) (6,0,1) (0,3,1) (0,0,3)};
	\addplot3 [draw, ultra thick, fill=\colornu, opacity=\opazero] coordinates {(6,0,1) (7,0,0) (6,1,0) (6,0,1)};
	\addplot3 [draw, ultra thick, fill=gray, opacity=\opazero] coordinates {(0,3,1) (0,4,0) (6,1,0) (6,0,1)};
\end{axis}
\draw (\titlex,\titley) node {Lower and upper faces of $\scrP + \scrQ$};
\end{scope}

\end{tikzpicture}

\caption{Construction of $\bl_{(1,2,3)}(\kk^3)$} \label{fig:wt-blow-up}
\end{figure}
\end{center}

\begin{prop} \label{Enu-Pnu}
For each $a = (a_1, \ldots, a_n) \in \kk^n \setminus \{0\}$, let $C_a := \{(a_1t^{\nu_1}, \ldots, a_nt^{\nu_n}): t \in \kk\} \subset \kk^n$. The strict transform of each $C_a$ intersects $\Enu$ at precisely one point which we denote by $[a]$. The map $[a_1: \cdots : a_n] \to [a]$ yields an isomorphism between $ \wwpnn{\nu}{n-1}$ and $\Enu$.
\end{prop}

\begin{proof}
Fix $a = (a_1, \ldots, a_n) \in \kk^n \setminus \{0\}$ and $j$ such that $a_j \neq 0$. Let $\scrC'_j$ be the convex polyhedral cone  generated by $(p/\nu_i)e_i - (p/\nu_j)e_j$, $i = 1, \ldots, \hat j, \ldots, n$, and $\scrC_j$ be the convex polyhedral cone generated by $\scrC'_j$ and $e_j$. Let $S_j$ and $S'_j$ be the semigroups of integral elements respectively in $\scrC_j$ and $\scrC'_j$. \Cref{xp-thm} implies that $\kk[U_{\alpha_j} \cap \xaa{\scrP+\scrQ}] = \kk[S_j]$ and $\kk[U_{\alpha_j} \cap \corbit{\scrP} ] = \kk[S'_j]$. Let $C'_a$ be the strict transform of $C_a$ in $\blnukn$. It turns out that (\cref{exercise:Enu-Pnu})
\begin{prooflist}
\item \label{Enu-Pnu:S_j} For every $\alpha \in S_j$, either $x^\alpha$ identically vanishes on $C_a$ or $\ord_t(x^\alpha|_{(a_1t^{\nu_1}, \ldots, a_nt^{\nu_n})})$ is nonnegative. This implies that $C'_a \cap \Enu \subset U_{\alpha_j}$.
\item \label{Enu-Pnu:S'_j} $S'_j = \{\alpha = (\alpha_1, \ldots, \alpha_n) \in \zz^n: \langle \nu, \alpha \rangle = 0,\ \alpha_i \geq 0,\ i = 1, \ldots, \hat j, \ldots, n\}$. If $\alpha \in S'_j$, it follows that
\begin{prooflist}
\item either $\alpha_i > 0$ for some $i$ such that $a_i = 0$, in which case $x^\alpha$ identically vanishes on $C'_a$,
\item or $\alpha_i \neq 0$ if and only if $a_i \neq 0$, in which case $x^\alpha$ takes the constant (nonzero) value $a^\alpha$ on $C'_a$.
\end{prooflist}
\end{prooflist}
These observations together with \cref{ideal-prop} implies that $C'_a \cap E_\nu$ corresponds precisely to the point $[a_1: \cdots : a_n] \in \pp^{n-1}(\nu)$, as required.
\end{proof}

We write $\Onu$ for the torus $\orbit{\scrP}$ of $\Enu = \corbit{\scrP}$. The isomorphism between $\Enu$ and $\pp^{n-1}(\nu)$ from \cref{Enu-Pnu} induces an isomorphism $\Onu \cong \pp^{n-1}(\nu) \setminus V(x_1 \cdots x_n) \cong \nktoruss{n-1}$. \Cref{U_Q} implies that $\Onu$ is a nonsingular hypersurface of $\blnukn$. Let $I \subset [n] := \{1, \ldots, n\}$ and $\Ki := V(x_j: j \not\in I)$ be the corresponding coordinate subspace of $\kk^n$. We identify $\Ki$ with $\kk^{|I|} $. Let $\nu'$ be a weighted order with positive weights on $\kk[x_i : i \in I]$ such that $(\nu'(x_i): i \in I)$ is proportional to $(\nu_i: i \in I)$. It follows from the definition of a weighted blow up that $\blnuprimeknn{I}$ can be identified with the strict transform of $\Ki$ on $\blnukn$. The following proposition compiles some properties of the embedding $\blnuprimeknn{I} \into \blnukn$, its proof is left as an exercise.

\begin{prop} \label{blow-up-coordinates}
Assume $\gcd(\nu_i: i \in I) =  1$. Let $k := |I|$. Then
there is a Zariski open neighborhood $U$ of $\Onuprime$ in $\blnukn$ and regular functions $(z_1, \ldots, z_n)$ on $U$ such that
\begin{enumerate}
\item $U \cong  \kk \times \nktoruss{k-1} \times \kk^{n-k}$ with coordinates $(z_1, \ldots, z_n)$,
\item $z_1, \ldots, z_k$ are monomials in $(x_i: i \in I)$,
\item $\nu(z_1) = 1$, $\nu(z_i) = 0$, $2 \leq i \leq n$,
\item for all $i' \not\in I$, there is $i'$ such that $z_{i'} = x_{i'}/z_1^{\nu_{i'}}$,
\item $\Enu \cap U = V(z_1) \cong \nktoruss{k-1} \times \kk^{n-k}$,
\item $\Onu  = (\Enu \cap U) \setminus V(z_{k+1} \cdots z_n)$,
\item $\blnuprimeknn{I} \cap U  = V(z_{k+1}, \ldots, z_n) \cong \kk \times \nktoruss{k-1}$,
\item $\Enuprime \cap U  = \Onuprime = V(z_1, z_{k+1}, \ldots, z_n) \cong \nktoruss{k-1}$. \qed
\end{enumerate}
\end{prop}

\subsection{Exercises}

\begin{exercise} \label{exercise:Enu-Pnu}
Verify observations \ref{Enu-Pnu:S_j} and \ref{Enu-Pnu:S'_j} from the proof of \cref{Enu-Pnu}.
\end{exercise}

\begin{exercise}
Prove \cref{blow-up-coordinates}. [Hint: Choose $\beta_1, \ldots, \beta_k \in \zz^I$ such that $\langle \nu, \beta_1 \rangle = 1$ and $(\beta_2, \ldots, \beta_k)$ is a basis of $\nu^\perp \cap \zz^I$. Set $\beta_j := e_j - \nu_j\beta_1$, $j = k+1, \ldots, n$. Show that $(\beta_1, \ldots, \beta_n)$ is a basis of $\zz^n$. Choose $d \geq 1$ such that the relative interior of $d\scrP \cap \zz^I$ contains a point $\beta_0$ with integral coordinates and consider the cone $\scrC_{\beta_0}$ generated by $d\scrP + \scrQ - \beta_0$. Show that $\scrC_{\beta_0} \cap \zz^I$ is generated as a semigroup by $\beta_1, \pm\beta_2, \ldots, \pm\beta_k$, and $\scrC_{\beta_0} \cap \zz^n$ is generated as a semigroup by $\beta_1, \pm\beta_2, \ldots, \pm\beta_k, \beta_{k+1}, \ldots, \beta_n$. Set $z_j := x^{\beta_j}$, $j = 1, \ldots, n$.] 
\end{exercise}

\chapter{Number of zeroes on the torus: BKK bound} \label{bkk-chapter}
\section{Introduction}
In this chapter we derive Bernstein's theorem for the number of isolated solutions of generic systems of $n$ Laurent polynomials on the algebraic torus $\nktorus$ over an algebraically closed field $\kk$, apply it to derive some properties of mixed volume, and discuss some open problems related to Bernstein's theorem. The term \index{BKK!bound}``BKK bound'' in the title of this section refers to the bound from Bernstein's theorem which is also known in the literature as ``BKK theorem'' after D.\ Bernstein, A.\ Kushnirenko and A.\ Khovanski.

\section{Mixed volume} \label{mixed-section}

The set of convex polytopes in $\rr^n$, $n \geq 1$, is a commutative semigroup under Minkowski addition (see \cref{appolytopes}). The interaction between Minkowski addition and volume gives rise to the theory of {\em mixed volumes}. The starting point of this theory is the following result proven in \cref{minkowski-section}:

\begin{repthm}{volume-thm}
Let $\scrP_1, \ldots, \scrP_s$ be convex polytopes in $\rr^n$. Then there are nonnegative real numbers $v_\alpha(\scrP_1, \ldots, \scrP_s)$ for all $\alpha \in \scrE_s:= \{(\alpha_1, \ldots, \alpha_s) \in \zzeroo{s}: \alpha_1 + \cdots  + \alpha_s = n\}$ such that for all $\lambda = (\lambda_1, \ldots, \lambda_s) \in \rzeroo{s}$,
\begin{align*}
\vol_n(\lambda_1\scrP_1 + \cdots + \lambda_s \scrP_s) = \sum_{\alpha \in \scrE_s} v_\alpha(\scrP_1, \ldots, \scrP_s) \lambda_1^{\alpha_1} \cdots \lambda_s^{\alpha_s}
\end{align*}
where $\vol_n$ is the $n$-dimensional Euclidean volume.
\end{repthm}

\begin{defn}  \label{mixed-defn}
\index{Mixed volume}The {\em mixed volume} $\mv(\scrP_1, \ldots, \scrP_n)$ of convex polytopes $\scrP_1, \ldots, \scrP_n$ in $\rr^n$ is $v_{(1, \ldots, 1)}(\scrP_1, \ldots, \scrP_n)$.
\end{defn}


\begin{thm} \label{mixed-unique}
Let $\scrK$ be any collection of convex polytopes in $\rr^n$ which is invariant under Minkowski addition\footnote{E.g.\ $\scrK$ may be the set of all convex polytopes in $\rr^n$, or the set of convex integral polytopes in $\rr^n$.}. Then $\mv: \scrK^n \to \rr$ is the unique function such that
\begin{enumerate}
\item $\mv(\scrP, \ldots, \scrP) = n!\vol_n(\scrP)$ for all $\scrP \in \scrK$,
\item $\mv$ is symmetric in its arguments, and
\item $\mv$ is {\em multiadditive}, i.e.\
\begin{align*}
\mv(k_1\scrP_1 + k'_1\scrP'_1, \scrP_2, \ldots, \scrP_n)
	= k_1\mv(\scrP_1, \scrP_2, \ldots, \scrP_n) + k'_1\mv(\scrP'_1, \scrP_2, \ldots, \scrP_n)
\end{align*}
for all $k_1, k'_1 \in \zzero$ and $\scrP_1, \ldots, \scrP_n, \scrP'_1 \in \scrK$.
\end{enumerate}
Moreover, $\mv$ can be expressed in terms of the volume (we write $[n]$ to denote $\{1, \ldots, n\}$):
\begin{align}
\mv(\scrP_1, \ldots, \scrP_n)
	&= \sum_{ \substack{I \subseteq [n]\\ I \neq \emptyset}} (-1)^{n - |I|}\vol_n\left(\sum_{i \in I} \scrP_i\right)
	\label{vol-to-mixed}
\end{align}
\end{thm}

\begin{proof}
This follows from combining \cref{volume-thm,tom-lemma,mixed-lemma}.
\end{proof}

\begin{example} \label{mixed-example}
For $n = 1$, a convex polytope is simply an interval and its mixed volume is its length. For $n = 2$, if $\scrP, \scrQ$ are convex polygons in $\rr^2$, then identity \eqref{vol-to-mixed} implies (see \cref{mixed-figure}) that
\begin{align}
\mv(\scrP, \scrQ) = \area(\scrP + \scrQ) - \area(\scrP) - \area(\scrQ) \label{vol-to-mixed-2}
\end{align}
\end{example}

\begin{figure}[h]
\begin{center}
\begin{tikzpicture}[
	dot/.style = {
      draw,
      fill,
      circle,
      inner sep = 0pt,
      minimum size = 3pt
    }, scale=0.6
    ]
    \def\shiftone{4}
    \def\shifttwo{5}
    \draw (\shiftone,1.5) node {$+$};    	
    \draw (\shiftone+\shifttwo,1.5) node {$=$};    	
	\begin{scope}[shift={(0,0)}]
		\draw [gray,  line width=0pt] (-0.5,-0.5) grid (3.5,3.5);
       	
       	\draw[ultra thick, fill=green, opacity=0.4 ] (0,0) -- (2,0) -- (2,2) -- (0,2) -- cycle;
		\draw (-0.5,-0.5) node [below right] {
     			\small
     			$\scrP$
     		};    	
	\end{scope}
	
	\begin{scope}[shift={(\shifttwo,0)}]
		\draw [gray,  line width=0pt] (-0.5,-0.5) grid (3.5,3.5);
       	
       \draw[ultra thick, fill=green, opacity=0.4 ] (0,0) -- (1,0) -- (2,1) -- (1,2) -- (0,1) -- cycle;
       \draw (-0.5,-0.5) node [below right] {
     			\small
     			$\scrQ$
     		};    	
	\end{scope}
	
	\begin{scope}[shift={(2*\shifttwo,0)}]
			\draw [gray,  line width=0pt] (-0.5,-0.5) grid (4.5,4.5);
	       	
	       	\draw[ultra thick] (0,0) -- (3,0) -- (4,1) --  (4,3) -- (3,4) -- (1,4) -- (0,3) -- cycle;
	       	\draw[fill=blue, opacity=0.4] (0,2) -- (2,2) -- (2,3) -- (3,4) -- (1,4) -- (0,3) -- cycle;
	       	\draw[fill=blue, opacity=0.4] (2,0) -- (3,0) -- (4,1) -- (4,3) -- (3,2) -- (2,2) -- cycle;
	       	\draw[dashed, ultra thick, fill=green, opacity=0.4 ] (0,0) -- (2,0) -- (2,2) -- (0,2) -- cycle;
	       	\begin{scope}[shift={(2,2)}]
	       		\draw[ultra thick, dashed, fill=green, opacity=0.4 ] (0,0) -- (1,0) -- (2,1) -- (1,2) -- (0,1) -- cycle;
	       	\end{scope}
	        \draw (-0.5,-0.5) node [below right] {
       			\small
       			$\scrP + \scrQ$
       		};    	
		\end{scope}
		
\end{tikzpicture}
\caption{$\mv(\scrP, \scrQ)$ is $8$, which is the area of the blue part of $\scrP+\scrQ$}  \label{mixed-figure}
\end{center}
\end{figure}

\begin{rem} \label{mixed-remark}
\Cref{volume-thm} implies that the mixed volume is nonnegative, and identity \eqref{vol-to-mixed} implies that $\mv$ is invariant under volume preserving transformations of $\rr^n$. Moreover, \cref{vol-basics} coupled with identity \eqref{vol-to-mixed} implies that mixed volume is continuous with respect to the Hausdorff distance on polytopes. In \cref{bkk-convex-applisection} we use Bernstein's theorem to deduce some other basic properties of mixed volume.
\end{rem}
%
%

\section{Theorems of Kushnirenko and Bernstein} \label{bkk-section}
\subsection{The bound}
Given Laurent polynomials  $f_1, \ldots, f_n$ in indeterminates $x_1, \ldots, x_n$, the number of their zeroes counted with multiplicity is
\begin{align*}
\multfntorus := \sum_{a \in \nktorus} \multf{a} 
\end{align*}
where $\multf{a}$ is the {\em intersection multiplicity} of $f_1, \ldots, f_n$ at $a$. If $V := V(f_1, \ldots, f_n)$ has non-isolated points, then $\multfntorus = \infty$ (\cref{int-mult-curve}). We write $\multfntorusiso$ for the sum of intersection multiplicities of $f_1, \ldots, f_n$ at all the isolated points of $V$. Recall that the {\em support} $\supp(f)$ of a Laurent polynomial $f = \sum_{\alpha \in \scrA}c_\alpha x^\alpha$ is the set of all $\alpha \in \zz^n$ such that $c_\alpha \neq 0$. Given a subset $\scrS$ of $\rr^n$, we say that $f$ is {\em supported at} $\scrS$ if $\supp(f) \subseteq \scrS$, and we write $\scrL(\scrS)$ for the set of all Laurent polynomials supported at $\scrS$. Given an ordered collection $\scrA := (\scrA_1, \ldots \scrA_m)$ of {\em finite} subsets of $\zz^n$, we write $\scrL(\scrA)$ for the set of all $m$-tuples $(f_1, \ldots, f_m)$ of Laurent polynomials such that $f_j$ is supported at $\scrA_j$ for each $j$. We say that some property holds for {\em generic} $f_j$ supported at $\scrA_j$, $j = 1, \ldots, m$, if it holds for all $(f_1, \ldots, f_m)$ in a nonempty Zariski open subset of $\scrL(\scrA) \cong \prod_{j= 1}^m \scrL(\scrA_j) \cong \kk^{\sum_j |\scrA_j|}$.

\begin{thm}[Bernstein's theorem: the bound] \label{bkk-bound-thm}
\index{Bernstein's!theorem!bound}
Let $\scrP_j$ be the convex hull of $\scrA_j$, $j = 1, \ldots, n$. If $\supp(f_j) \subset \scrA_j$, $j = 1, \ldots, n$, then
\begin{align}
\multfntorusiso \leq \mv(\scrP_1, \ldots, \scrP_n) \label{bkk-ineq}
\end{align}
Moreover, for generic $f_j$ supported at $\scrA_j$, $j = 1, \ldots, n$, we have
\begin{align}
\multfntorus = \multfntorusiso = \mv(\scrP_1, \ldots, \scrP_n) \label{bkk-eq}
\end{align}
\end{thm}
The {\em Newton polytope} $\np(f)$ of a Laurent polynomial $f$ is the convex hull (in $\rr^n$) of the support of $f$. Bernstein's theorem in particular states that the number of isolated solutions of a system of polynomials is bounded by the mixed volume of their Newton polytopes.
A.\ Kushnirenko proved Bernstein's theorem for the case that all the Newton polytopes are identical, in which case the mixed volume equals $n!$ times the volume of any of these polytopes. D.\ Bernstein found \cref{bkk-bound-thm} while trying to understand and generalize Kushnirenko's result. Kushnirenko however, not only gave the bound, but also gave a precise characterization of the collections of $f_1, \ldots, f_n$ for which the bound is achieved. There is a natural way to understand this characterization in the case that $\kk = \cc$; we describe it now.

\subsection{The non-degeneracy condition} \label{sec:bkk-non-degeneracy-statement}
Let $\scrA_1, \ldots, \scrA_n$ be finite subsets of $\zz^n$, and $f_1, \ldots, f_n$ be Laurent polynomials in $x_1, \ldots, x_n$ over $\cc$ such that $\supp(f_j) = \scrA_j$, $j = 1, \ldots, n$. Assume there are Laurent polynomials $g_1, \ldots, g_n$ such that $\supp(g_j) \subseteq \scrA_j$ for each $j$, and $\multfnctorusiso < \multgnctorusiso$. Write $h_j := (1-t)f_j + tg_j$, $j = 1, \ldots, n$. Since the $g_j = h_j|_{t=1}$ have ``more'' common zeroes on $\ntorus$ than the $f_j = h_j|_{t=0}$, intuitively we may expect that there is a curve $C(t)$ on $\ntorus$ such that $h_j(C(t)) =0$ and $\lim_{t \to 0}C(t)$ is {\em not} on $\ntorus$, i.e.\ as $t$ approaches $0$, either $C(t)$ approaches one of the coordinate hyperplanes of $\cc^n$, or $|C(t)|$ approaches infinity (see \cref{fig:h-t}).
\begin{figure}[h]
\begin{center}

\begin{tikzpicture}[scale=0.33]

   \def\xmin{0}
   \def\xmax{18}
   \def\ymin{0}
   \def\ymax{18}
   \def\s{6}
   \def\a{6}
   \def\b{1}
   \def\c{18}
   \def\xshift{-6}
   \def\tlabelx{0}
   \def\tlabely{-0.3}
   \def\nsamples{103}
   \def\tmin{-5}
   \def\tmax{4}

   \begin{scope}[shift={(0,0)}]
		  \begin{axis}[
		  	xmin = \xmin, xmax=\xmax, ymin = \ymin, ymax= \ymax,
		  	axis equal=true, axis equal image=true, hide axis
		  	]
		    \addplot[blue, thick, domain=\tmin:\tmax, samples=\nsamples]({x^2 + 1 + x} ,{x^2 + 1});
		    \addplot[red, thick, domain=\xmin:\xmax,samples=2](x,x+1/\s);
		  \end{axis}

		  \draw (\tlabelx,\tlabely) node [below right] {\small $t = 0$};    	
		
	\end{scope}
	
	\begin{scope}[shift={(3*\xshift,0)}]
	   		  \begin{axis}[
			  	xmin = \xmin, xmax=\xmax, ymin = \ymin, ymax= \ymax,
			  	axis equal=true, axis equal image=true, hide axis
			  	]
			    \addplot[blue, thick, domain=\tmin:\tmax, samples=\nsamples]({x^2 + 1 + x} ,{x^2 + 1});
			    \addplot[red, thick, domain=\xmin:\xmax,samples=2](x,{(\b*x + \c)/\a});
			  \end{axis}
			  \draw (\tlabelx,\tlabely) node [below right] {\small $t = 1$};    	
	 		
	\end{scope}
	
	\begin{scope}[shift={( 2*\xshift,0)}]
			  \def\t{0.6}
	   		  \begin{axis}[
			  	xmin = \xmin, xmax=\xmax, ymin = \ymin, ymax= \ymax,
			  	axis equal=true, axis equal image=true, hide axis
			  	]
			    \addplot[blue, thick, domain=\tmin:\tmax, samples=\nsamples]({x^2 + 1 + x} ,{x^2 + 1});
			    \addplot[red, thick, domain=\xmin:\xmax,samples=2](x,{(x*((1-\t)*\s + \t*\b)+ 1-\t + \t*\c)/((1-\t)*\s + \t*\a)});
			  \end{axis}
			
	 		  \draw (\tlabelx,\tlabely) node [below right] {\small $t = \t$};    	
		\end{scope}
	
	\begin{scope}[shift={(\xshift,0)}]
			  \def\t{0.3}
			  \begin{axis}[
			  	xmin = \xmin, xmax=\xmax, ymin = \ymin, ymax= \ymax,
			  	axis equal=true, axis equal image=true, hide axis
			  	]
			    \addplot[blue, thick, domain=\tmin:\tmax, samples=\nsamples]({x^2 + 1 + x} ,{x^2 + 1});
			    \addplot[red, thick, domain=\xmin:\xmax,samples=2](x,{(x*((1-\t)*\s + \t*\b)+ 1-\t + \t*\c)/((1-\t)*\s + \t*\a)});
			  \end{axis}
			
	 		  \draw (\tlabelx,\tlabely) node [below right] {\small $t = \t$};    	
		\end{scope}
		
		\def\textshiftx{1}
		\def\textshifty{5}
		\def\polyshiftone{1}
		\def\polyshifttwo{9}
		
		\begin{scope}[shift={(-\xshift + \textshiftx,0)}]	
			\node [
				right,
			 	text width= 5.4cm,align=justify] at (0,\textshifty) {
			 		\small
						$(f_1, f_2) = (6y-6x-1, y- (y-x)^2 - 1)$ \\
						$(g_1, g_2) = (6y - x - 18, y - (y-x)^2 - 1)$
					};
					
			\begin{scope}[shift={(\polyshiftone,0)}]					
					\draw [gray,  line width=0pt] (-0.5,-0.5) grid (2.5,2.5);
						\draw [<->] (0,2.5)  |- (2.5,0);
					      	
					       	\draw[thick, fill=green, opacity=0.4 ] (0,0) -- (1,0) -- (0,1)-- cycle;
					     	\draw (-0.5, \tlabely) node [below right] {
								\small
								$\np(f_1) = \np(g_1)$ \newline
							};
			\end{scope}
			
			\begin{scope}[shift={(\polyshiftone+\polyshifttwo,0)}]					
						\draw [gray,  line width=0pt] (-0.5,-0.5) grid (2.5,2.5);
						\draw [<->] (0,2.5)  |- (2.5,0);
      	
				       	\draw[thick, fill=green, opacity=0.4 ] (0,0) -- (2,0) -- (0,2)-- cycle;

				     	\draw (-0.5, \tlabely) node [below right] {
							\small
							$\np(f_2) = \np(g_2)$ \newline
						};
			\end{scope}
				
		\end{scope}
		
\end{tikzpicture}
\caption{\mbox{One of  the common roots of $(1-t)f_j + tg_j = 0$, $j = 1,2$, approaches infinity as $t \to 0$}}
\label{fig:h-t}
\end{center}

\end{figure}
In any event, assuming our intuition is correct, there is a punctured neighborhood $U$ of the origin on $\cc$ and a parametrization $U \to C(t)$ of the form $\gamma: t \mapsto (a_1 t^{\nu_1} + \cdots, \ldots, a_nt^{\nu_n} + \cdots)$, where $a = (a_1, \ldots, a_n) \in \ntorus$ and for each $j$, $\nu_j$ is the order (in $t$) of the $j$-th coordinate. Since $\lim_{t \to 0}C(t) \not\in \ntorus$, it follows that not all the $\nu_j$ are zero. Let $\nu$ be the element in $\rnstar$ with coordinates $(\nu_1, \ldots, \nu_n)$ with respect to the basis dual to the standard basis of $\rr^n$. Given a subset $\scrS$ of $\rr^n$ and a Laurent polynomial $g = \sum_\alpha c_\alpha x^\alpha$ supported at $\scrS$, we write
\begin{align}
\In_{\scrS, \nu}(g)
	&:= \sum_{\alpha \in \In_\nu(\scrS)} c_\alpha x^\alpha
    =
	\begin{cases}
	\In_\nu(g) & \text{if}\ \supp(g) \cap \In_\nu(\scrS) \neq \emptyset,\\
	0 & \text{otherwise.}
	\end{cases}
\label{insnu}
\end{align}
where $\In_\nu(\scrS)$ is defined as in \cref{basic-convex-section}. Let $q_j := \min_\nu(\scrA_j)$. Since $\supp(g_j) \subseteq \scrA_j = \supp(f_j)$, it follows that
\begin{align*}
h_j(\gamma(t)) 
	&= t^{q_j} \In_{\scrA_j, \nu}(f_j)(a) + t^{q_j + 1} (- \In_{\scrA_j, \nu}(f_j)(a) + t^\epsilon \In_{\scrA_j, \nu}(g_j)(a)) + \cdots
\end{align*}
where $\epsilon$ is nonnegative, and the orders in $t$ of the omitted terms are higher than $q_j$. Since $h_j(\gamma(t)) \equiv 0$, it follows that $\In_{\scrA_j, \nu}(f_j)(a) = 0$ for each $j = 1, \ldots, n$. This leads to the following definition. 

\begin{defn} \label{defn:b-non-deg}
Let $\scrA_1, \ldots, \scrA_m$ be finite subsets of $\zz^n$ and $(f_1, \ldots, f_m) \in \scrL(\scrA)$. We say that $f_1, \ldots, f_m$ are \index{Non-degeneracy!with respect to tuples of finite subsets of $\zz^n$}{\em $(\scrA_1, \ldots, \scrA_m)$-non-degenerate} if they satisfy the following condition:
\begin{align}
\parbox{0.45\textwidth}{for each $\nu \in \rnstar \setminus \{0\}$, there is no common root of $\In_{\scrA_j, \nu}(f_j)$, $j = 1, \ldots, m$, on $\nktorus$.}
\tag{$\text{ND}^*$}
\label{b-non-degeneracy}
\end{align}
We say that $f_1, \ldots, f_m$ are \index{BKK!non-degeneracy}\index{Non-degeneracy!BKK}{\em BKK non-degenerate} if they are $(\supp(f_1), \ldots, \supp(f_m))$-non-degenerate.
\end{defn}

The preceding argument suggests that  for $\kk = \cc$, $\scrA$-non-degeneracy is sufficient for the maximality of $\multfntorusiso$. It turns out that it is also necessary; Kushnirenko proved it in the case that the convex hulls of the $\scrA_j$ are identical and Bernstein treated the general case. Both necessity and sufficiency remain valid even if $\cc$ is replaced by an arbitrary algebraically closed field:

\begin{thm}[Bernstein's theorem: non-degeneracy condition] \label{bkk-non-degenerate-thm}
\index{Bernstein's!theorem!non-degeneracy}
Let $\scrP_j := \conv(\scrA_j)$, $j = 1, \ldots, n$. If the mixed volume of $\scrP_1, \ldots, \scrP_n$ is nonzero, then the bound \eqref{bkk-ineq} is satisfied with an equality if and only if $f_1, \ldots, f_n$ are $(\scrA_1, \ldots, \scrA_n)$-non-degenerate.
\end{thm}

\begin{rem} \label{wt-non-degenerate-remark}
Recall from \cref{wt-order-section} that integral elements in $\rnstar$ can be identified with {\em weighted orders} on the ring of Laurent polynomials. In the case that $\np(f_j) = \conv(\scrA_j)$ for each $j = 1, \ldots, m$, condition \eqref{b-non-degeneracy} is equivalent to the following condition:
\begin{align}
\parbox{0.48\textwidth}{for each nontrivial weighted order $\nu$, there is no common root of $\In_\nu(f_j)$, $j = 1, \ldots, m$, on $\nktorus$.}
\tag{$\text{ND}'^*$}
\label{b-non-degeneracy-wt}
\end{align}
\end{rem}

\begin{rem} \label{finite-remark}
On the face of it condition \eqref{A-non-degeneracy} consists of uncountably many conditions, one for each element in $\rnstar$. However, it is equivalent to finitely many conditions. Indeed, let $\scrP_j := \conv(\scrA_j)$, $j = 1, \ldots, m$, and $\scrP := \scrP_1 + \cdots + \scrP_m$. For each face $\scrQ$ of $\scrP$, there are unique faces $\scrQ_j$ of $\scrP_j$ such that $\scrQ = \scrQ_1 + \cdots + \scrQ_m$ (\cref{prop:sum-faces}). Given $f_j := \sum_{\alpha} c_{j, \alpha}x^\alpha$, let $f_{j,\scrQ} := \sum_{\alpha \in \scrQ_j} c_{j,\alpha} x^\alpha$ be the ``component'' of $f_j$ supported at $\scrQ_j$. Then \eqref{b-non-degeneracy} is equivalent to the following condition:
\begin{align}
\parbox{0.5\textwidth}{for each face $\scrQ$ of dimension less than $n$ of $\scrP$, there is no common root of $f_{j,\scrQ}$, $j = 1, \ldots, m$, on $\nktorus$.}
\tag{$\text{ND}''^*$}
\label{b-non-degeneracy-pol}
\end{align}
\end{rem}

Kushnirenko's theorem follows immediately by applying \cref{mixed-unique,bkk-bound-thm,bkk-non-degenerate-thm} to the case that each $f_j$ is supported at the same (finite) subset of $\zz^n$:

\begin{cor}[{Kushnirenko \cite{kush-poly-milnor}}] \label{kushnirenko}
\index{Kushnirenko's theorem!on number of solutions}
Let $f_1, \ldots, f_n$ be Laurent polynomials supported at a finite subset $\scrA$ of $\zz^n$ and $\scrP:= \conv(\scrA)$. Then $\multfntorusiso \leq n! \vol_n(\scrP)$. If $\vol_n(\scrP) > 0$, then the bound is satisfied with equality if and only if $f_1, \ldots, f_n$ are $(\scrA, \ldots, \scrA)$-non-degenerate. \qed
\end{cor}

\begin{example} \label{example:kush-inequal-np}
To attain the bound in Bernstein's and Kushnirenko's theorems it is {\em not} necessary that $\np(f_j)$ be equal to $\scrP_j$ for each (or, any!) $j$. Indeed, consider $f_1 = 1 + x^4 + x^2y^4$ and $f_2 = xy^2 + x^3y^3 + x^6$. Newton polygons of the $f_j$ are proper subsets of the triangle $\scrP$ with vertices $A = (0,0)$, $B = (6,0)$ and $C = (2, 4)$ (see \cref{fig:kush-inequal-np}). However, it is straightforward to check that if $\scrQ$ is any edge or vertex of $\scrP$,then either $f_{1,\scrQ}$ or $f_{2, \scrQ}$ is a monomial, so that $(f_1, f_2)$ satisfy the non-degeneracy condition \eqref{b-non-degeneracy-pol} from \cref{finite-remark}. Consequently, the number of solutions of $(f_1, f_2)$ on $\nktoruss{2}$ is $2\area(\scrP) = 24$. 
\end{example}

\begin{figure}[h]
\begin{center}
\def\xmin{-0.5}
\def\xmax{7.5}
\def\ymin{-0.5}
\def\ymax{5.5}
\def\opazero{0.5}
\def\colorzero{green}
\def\colorone{red}
\def\colortwo{blue}
\tikzstyle{dot} = [circle, minimum size=5pt, inner sep = 0pt, fill]
\def\scalefactor{0.6}

\begin{tikzpicture}[scale=\scalefactor]
\draw [gray,  line width=0pt] (\xmin, \ymin) grid (\xmax,\ymax);
\draw [<->] (0, \ymax) |- (\xmax, 0);
\coordinate (A) at (0, 0);
\coordinate (B) at (6, 0);
\coordinate (C) at (2, 4);
\coordinate (D) at (4, 0);
\fill[\colorzero, opacity=\opazero ] (A) --  (B) -- (C);
\draw[thick] (A) -- (C) -- (D) -- cycle;

\node[dot, \colorone] at (A) {};
\node[dot, \colorone] at (B) {};
\node[dot, \colorone] at (C) {};
\node[dot, \colortwo] at (D) {};

\node[anchor = north east] at (A) {\picfontsize $A$};
\node[anchor = north west] at (B) {\picfontsize $B$};
\node[anchor = south] at (C) {\picfontsize $C$};
\node[anchor = north west] at (D) {\picfontsize $D$};

\node at (2, 1.5) {\picfontsize $\np(f_1)$};

\def\xshift{5}
\pgfmathsetmacro\shiftone{\xshift + \xmax -\xmin};

\begin{scope}[shift={(\shiftone,0)}]
\draw [gray,  line width=0pt] (\xmin, \ymin) grid (\xmax,\ymax);
\draw [<->] (0, \ymax) |- (\xmax, 0);
\coordinate (A) at (0, 0);
\coordinate (B) at (6, 0);
\coordinate (C) at (2, 4);
\coordinate (E) at (3, 3);
\coordinate (F) at (1, 2);

\fill[\colorzero, opacity=\opazero ] (A) --  (B) -- (C);
\draw[thick] (B) -- (E) -- (F) -- cycle;

\node[dot, \colorone] at (A) {};
\node[dot, \colorone] at (B) {};
\node[dot, \colorone] at (C) {};
\node[dot, \colortwo] at (E) {};
\node[dot, \colortwo] at (F) {};

\node[anchor = north east] at (A) {\picfontsize $A$};
\node[anchor = north west] at (B) {\picfontsize $B$};
\node[anchor = south] at (C) {\picfontsize $C$};
\node[anchor = south west] at (E) {\picfontsize $E$};
\node[anchor = south east] at (F) {\picfontsize $F$};

\node at (3, 1.75) {\picfontsize $\np(f_2)$};
\end{scope}

\end{tikzpicture}
\caption{The bound in Kushnirenko's theorem is attained when $\scrA$ is the triangle $ABC$, $\np(f_1)$ is the triangle $ADC$ and $\np(f_2)$ is the triangle $BEF$}  \label{fig:kush-inequal-np}
\end{center}
\end{figure}

\Cref{example:kush-inequal-np} shows that it is possible that $\mv(\scrP_1, \ldots, \scrP_n) = \mv(\scrQ_1, \ldots, \scrQ_n)$ even if each $\scrP_j$ {\em properly} contains $\scrQ_j$. In fact using Bernstein's theorem it is possible to precisely characterize the situations in which this happens. First we need to make a definition. Convex polytopes $\scrQ_1, \ldots, \scrQ_m$ in $\rr^n$ are said to be \index{Dependence!of polytopes}\index{Polytope!dependent}{\em dependent} if there is a nonempty subset $I$ of $[m] := \{1, \ldots, m\}$ such that $\dim(\sum_{i \in I} \scrQ_i) < |I|$; otherwise they are said to be {\em independent}.

\begin{thm}
If $\scrP'_j$ are convex polytopes in $\rr^n$ such that  $\scrP_j \subseteq \scrP'_j$ for each $j$, then $\mv(\scrP_1, \ldots, \scrP_n) \leq \mv(\scrP'_1, \ldots, \scrP'_n)$. The inequailty is strict if and only if both of the following are true:
\begin{enumerate}
\item $\scrP'_1, \ldots, \scrP'_n$ are independent, and
\item there is $\nu \in \rnstar \setminus \{0\}$ such that the collection $\{\In_\nu(\scrP_j): \scrP_j \cap \In_\nu(\scrP'_j) \neq \emptyset\}$ of polytopes is independent.
\end{enumerate}
\end{thm}

We prove this result in \cref{strictly-mixed-monotone}. Using this we can give the following alternate formulation of Bernstein's theorem in which $(\scrA_1, \ldots, \scrA_n)$-non-degeneracy is replaced by a combinatorial condition plus BKK non-degeneracy; we will prove it in \cref{alternate-corollary}.

\begin{thm}[Bernstein's theorem - alternate version] \label{alternate-bernstein}
\index{Bernstein's!theorem!alternate version}
Let $\scrA_j$ be finite subsets of $\zz^n$ and $f_j$ be Laurent polynomials supported at $\scrA_j$, $j = 1, \ldots, n$. Then
$$\multfntorusiso \leq \mv(\conv(\scrA_1), \ldots, \conv(\scrA_n))$$
If $\mv(\conv(\scrA_1), \ldots, \conv(\scrA_n)) > 0$, then the bound is satisfied with an equality if and only if both of the following conditions hold:
\begin{enumerate}
\item \label{bernstein-dependence} for each nontrivial weighted order $\nu$, the collection $\{\In_\nu(\np(f_j)): \In_\nu(\scrA_j) \cap \supp(f_j) \neq \emptyset\}$ of polytopes is dependent, and
\item  $f_1, \ldots, f_n$ are BKK non-degenerate, i.e.\ they satisfy \eqref{b-non-degeneracy-wt} with $m = n$.
\end{enumerate}
\end{thm}

\begin{rem}
The Newton polytopes of polynomials in \cref{example:kush-inequal-np} satisfy the property that for each nonzero $\nu \in \rnstar$, for at least one of these polytopes the face where $\nu$ attains the minimum is a vertex. These were extensively used by A.\ Khovanskii (see e.g. \cite{khovanskii-good-compactification}) who called these {\em developed} systems of polytopes. A basic reason for the usefulness of these systems is the following property (which is straightforward to check): if the collection of Newton polytopes of a system of Laurent polynomials is developed, then this system is automatically BKK non-degenerate. In particular, if $f_1, \ldots, f_n$ are Laurent polynomials in $(x_1, \ldots, x_n)$ such that their Newton polytopes form a developed system, then $\multfntorusiso = \mv(\np(f_1), \ldots, \np(f_n))$. 
\end{rem}

\section{Proof of Bernstein-Kushnirenko non-degeneracy condition} \label{bkk-non-deg-proof-section}
Let $\scrA := (\scrA_1, \ldots, \scrA_m)$, $m \geq 1$, be a collection of finite subsets of $\zz^n$. Let $\scrL(\scrA) := \prod_{j=1}^m \scrL(\scrA_j) \cong \kk^{\sum_j |\scrA_j|}$ and $\scrN(\scrA)$ be the set of all $\scrA$-non-degenerate $(f_1, \ldots, f_m) \in \scrL(\scrA)$. In \cref{bkk-existential-section} we prove the following result:

\begin{thm} \label{bkk-0}
$\scrN(\scrA)$ is a Zariski open subset of $\scrL(\scrA)$. If $m \geq \min\{n, \dim(\sum_{j=1}^m \conv(\scrA_j)) + 1\}$, then $\scrN(\scrA)$ is nonempty.
\end{thm}


In the case that $m = n$, define
\begin{align}
\multAntorusiso := \max\{\multfntorusiso: \supp(f_j) \subseteq \scrA_j,\ j =1 , \ldots, n\},
\label{multAntorusiso}
\end{align}
and let $\scrM(\scrA)$ be the set of all $(f_1, \ldots, f_n) \in \scrL(\scrA)$ such that $\multfntorusiso = \multAntorusiso$. In \cref{bkk-sufficient-section,bkk-necessary-section} we prove the following result:

\begin{thm} \label{bkk-1}
Assume $m = n$. Then
\begin{align*}
\scrM(\scrA) &=
						\begin{cases}
							\scrL(\scrA) & \text{if}\ \multAntorusiso = 0,\\
							\scrN(\scrA) & \text{if}\ \multAntorusiso > 0.
						\end{cases}
\end{align*}
In particular, $\scrM(\scrA)$ is a nonempty Zariski open subset of $\scrL(\scrA)$.
\end{thm}
\Cref{bkk-0} is required (in \cref{bkk-necessary-section}) for the proof of \cref{bkk-1}. The reason to defer the proof of \cref{bkk-0} to \cref{bkk-existential-section} is that it is somewhat more technical.

\subsection{Sufficiency of Bernstein-Kushnirenko non-degeneracy condition} \label{bkk-sufficient-section}
In this section we prove that $\scrN(\scrA) \subseteq \scrM(\scrA)$ in the case that $m = n$. In his proof of this part of his theorem, Bernstein used the theory of Puiseux series of curves, which only works in characteristic zero. For arbitrary characteristics, the analogous role is played by the notion of {\em branches} described in \cref{branchion-0}. Let $B = (Z, z)$ be a branch on $\nktorus$. Let $\nu_B \in \rnstar$ be the corresponding integral vector of weights, $\rho_B$ be a parameter of $B$ and $\In(B) \in \nktorus$ be the corresponding $n$-tuple of ``initial coefficients'' (see \cref{toric-center-section}).

\begin{lemma} \label{branch-lemma-2}
Let $B$ be a branch of a curve contained in the common zero set of Laurent polynomials $f_1, \ldots, f_m$ on $\nktorus$. Then $\In (B)$ is a common zero of $\In_{\nu_B}(f_j)$, $j = 1, \ldots, m$. In particular, if $B$ is a branch at infinity, and $\supp(f_j) \subseteq \scrA_j$ for each $j$, then both \eqref{b-non-degeneracy} and \eqref{b-non-degeneracy-wt} are violated with $\nu = \nu_B$.
\end{lemma}

\begin{proof}
The first assertion is a direct corollary of \cref{branch-lemma-1}. If $B$ is a branch at infinity, then $\nu_B$ is a nontrivial weighted order (\cref{branch-lemma-0}), so that the second assertion follows from the first assertion and the second identity from \eqref{insnu}.
\end{proof}

\begin{cor} \label{finite-cor}
Let $f_j$ be a Laurent polynomial supported at $\scrA_j$, $j = 1, \ldots, m$. If $f_1, \ldots, f_m$ are $(\scrA_1, \ldots, \scrA_m)$-non-degenerate, then $V(f_1, \ldots, f_m) \subset \nktorus$ is finite. \qed
\end{cor}

\begin{prop}
\label{bkk-sufficiency}
Let $f_j$ be a Laurent polynomial supported at $\scrA_j$, $j = 1, \ldots, n$. If $f_1, \ldots, f_n$ are $(\scrA_1, \ldots, \scrA_n)$-non-degenerate, then $\multfntorusiso = \multAntorusiso$. In particular, $\scrN(\scrA_1, \ldots, \scrA_n) \subseteq \scrM(\scrA_1, \ldots, \scrA_n)$.
\end{prop}

\begin{proof}
Assume to the contrary that there are Laurent polynomials $g_j$ supported at $\scrA_j$ such that $\multgntorusiso > \multfntorusiso$. It suffices to show that this leads to a contradiction. Define $h_j := (1-t)f_j + tg_j$, where $t$ is a new indeterminate. Since the set of zeroes of $f_1, \ldots, f_n$ on $\nktorus$ is finite (\cref{finite-cor}), assertion \eqref{non-max-condition} of \cref{mult-deformation-global} implies that the set of zeroes of $h_1, \ldots, h_n$ in $\nktorus \times \kk$ contains a curve which has a branch $B = (Z,z)$ at infinity (with respect to $\nktorus \times \kk$) at $t = 0$ and $B \not\subseteq \nktorus \times \{0\}$. Then $t|_B \not\equiv 0$ and $B$ determines a well-defined weighted order $\nu_B$ on $\kk[x_1, x_1^{-1}, \ldots, x_n, x_n^{-1}, t]$ such that $\nu_B(t) > 0$. Let $\nu$ be the restriction of $\nu_B$ to $\kk[x_1, x_1^{-1}, \ldots, x_n, x_n^{-1}]$. Since $\nu_B(t) > 0$ and $\supp(g_j) \subset \scrA_j$, it follows from \eqref{insnu} that
\begin{align*}
\In_{\scrA_j, \nu}(f_j)
	&=
	\begin{cases}
	\In_\nu(f_j) = \In_{\nu_B}(h_j)& \text{if}\ \supp(f_j) \cap \In_\nu(\scrA_j) \neq \emptyset,\\
	0 & \text{otherwise.}
	\end{cases}
\end{align*}
Therefore \cref{branch-lemma-2} implies that $\In_{\scrA_j, \nu}(f_j)$, $j = 1, \ldots, n$, have a common zero on $\nktorus$. Since $B$ is centered at infinity with respect to $\nktorus \times \kk$ and $\nu_B(t) > 0$, it follows that $\nu$ is a nontrivial weighted order on $\kk[x_1, x_1^{-1}, \ldots, x_n, x_n^{-1}]$. This contradicts the $\scrA$-non-degeneracy of $f_1, \ldots, f_n$, as desired.
\end{proof}


\subsection{Necessity of Bernstein-Kushnirenko non-degeneracy condition} \label{bkk-necessary-section}
In this section we finish (using \cref{bkk-0}) the proof of \cref{bkk-1}. Assume $\multAntorusiso > 0$ and $(f_1, \ldots, f_n) \in \scrL(\scrA)$ is a {\em $\scrA$-degenerate system} in $\scrL$, i.e.\ there is a nontrivial weighted order $\nu$ on $\kk[x_1, x_1^{-1}, \ldots, x_n, x_n^{-1}]$ such that $\In_{\scrA_1, \nu}(f_1), \ldots, \In_{\scrA_n,\nu}(f_n)$ have a common zero on $\nktorus$. It suffices to show that
$$\multfntorusiso < \multAntorusiso$$
We show it in {\em two} ways.
%
%

\begin{claim} \label{claim:dominant}
Given any point $a \in \nktorus$, there is $(g_1, \ldots, g_n) \in \scrL(\scrA)$ such that
\begin{enumerate}
  \item $\In_{\scrA_j, \nu}(g_j)(a) \neq 0$ for each $j = 1, \ldots, n$, and
  \item there is an isolated zero $b$ of $(g_1, \ldots, g_n)$ on $\nktorus$ such that $f_j(b) \neq 0$ for each $j = 1, \ldots, n$.
\end{enumerate}
\end{claim}

\begin{proof}
Due to \cref{bkk-0} we can pick $\scrA$-non-degenerate $(h_1, \ldots, h_n) \in \scrL(\scrA)$ which automatically satisfies the first property. For each $j$, let $\alpha_j$ be an arbitrary element of $\scrA_j$. For each $\epsilon := (\epsilon_1, \ldots, \epsilon_n) \in \nktorus$, let $h_{\epsilon_j, j} := h_j - \epsilon_j x^{\alpha_j}$. Since $\scrN(\scrA)$ is Zariski open (\cref{bkk-0}), it follows that $(h_{\epsilon_1,1}, \ldots, h_{\epsilon_n, n})$ is in $\scrN(\scrA)$ for generic $\epsilon \in \nktorus$, and therefore $\multntorus{h_{\epsilon_j, 1} }{h_{\epsilon_j, n} } = \multAntorusiso > 0$ (\cref{bkk-sufficiency}); in particular the map $\Phi: \nktorus \to \nktorus$ defined by $x \mapsto (x^{-\alpha_1}h_1, \ldots, x^{-\alpha_n}h_n)$ is dominant. Therefore it suffices to take $g_j := h_{\epsilon_j, j}$ for generic $(\epsilon_1, \ldots, \epsilon_n) \in \nktorus \setminus \Phi(V(f_1 \cdots f_n))$.
%
\end{proof}

\subsubsection{Bernstein's proof of necessity of non-degeneracy (Bernstein's trick)} \label{original-proof}
\index{Bernstein's!trick}
After a linear change of coordinates of $\zz^n$ (which corresponds to a monomial change of coordinates of $\nktorus$) and translating each $\scrA_j$ by an element in $\zz^n$ (which corresponds to multiplying each $f_j$ by a monomial) if necessary, we may arrange that $\nu = (0, \ldots, 0, 1)$ and $\In_\nu(\scrA_j) = \scrA_j \cap (\zz^{n-1} \times \{0\})$, $j = 1, \ldots, n$. For each $(g_1, \ldots, g_n) \in \scrL(\scrA)$, it follows that
\begin{align}
\In_{\scrA_j, \nu}(g_j) = g_j|_{x_n = 0} \label{initial-trick}
\end{align}
If $(a_1, \ldots, a_n)$ is a common zero of the $\In_{\scrA_j, \nu}(f_j)$ on $\nktorus$, then it follows that $a' := (a_1, \ldots, a_{n-1}, 0)$ is also a common zero of the $\In_{\scrA_j, \nu}(f_j)$ on $\kk^n$. Let $g_1, \ldots, g_n$ and $b$ be as in \cref{claim:dominant}. Identity \eqref{initial-trick} implies that $g_j(a') \neq 0$ for each $j = 1, \ldots, n$. Take any map $c: \kk \to \kk^n$ such that $c(0) = a'$ and $c(1) = b$ (e.g.\ we may take $c(t) = (1-t)a' + tb$) and set $h_j(x,t) := g_j(c(t))f_j(x) - f_j(c(t))g_j(x)$, $j = 1, \ldots, n$. Then each $h_j$ vanishes on the curve $C' := \{(c(t),t): t \in \kk\} \subset \kk^{n+1}$. Since $(b, 1)$ is an isolated zero of $h_1(x,1), \ldots, h_n(x, 1)$, and since $(a',0) \in C'$ is a ``point at infinity'' with respect to $\nktorus \times \kk$, assertion \eqref{non-max-condition} of \cref{mult-deformation-global} implies that $\multntorusiso{h_1(x,0)}{h_n(x,0)} < \multntorusiso{h_1(x,\epsilon)}{h_n(x,\epsilon)}$ for generic $\epsilon \in \kk$. Since each $h_j(x,\epsilon)$ is supported at $\scrA_j$, it follows that $\multfntorusiso < \multAntorusiso$, as required. \qed
\vspace{\topsep}

\subsubsection{A modified version of Bernstein's trick} \label{modified-proof}
Here we give an alternative proof of the necessity of Bernstein-Kushnirenko non-degeneracy by adapting Bernstein's trick to produce a curve $C'$ as in his original proof {\em without} changing the coordinates on $\nktorus$; this will be useful later, e.g.\ in the proofs of the weighted B\'ezout theorem (\cref{wt-bezout-proof-section}) and the extension of the BKK bound to the affine space (\cref{non-deg-1-proof-section}). Let $a = (a_1, \ldots, a_n)$ be a common zero of the $\In_\nu(f_j)$ on $\nktorus$. As in the original proof, let $(g_1, \ldots, g_n) \in \scrL(\scrA)$ and $b = (b_1, \ldots, b_n) \in \nktorus$ be as in \cref{claim:dominant}. Fix integers $\nu'_j > \nu_j := \nu(x_j)$. Let $C$ be the rational curve on $\kk^n$ parameterized by $c(t):= (c_1(t), \ldots, c_n(t)) :\kk \to \kk^n$ given by
\begin{align}
c_j(t) := a_jt^{\nu_j} + (b_j - a_j)t^{\nu'_j},\ j = 1, \ldots, n. \label{trick-curve}
\end{align}
Note that $c(1) = b \in \nktorus$, so that $c(t) \in \nktorus$ for all but finitely many values of $\kk$. Let $m_j := \min_{\scrA_j}(\nu)$, $j = 1, \ldots, n$. Then $t^{-m_j}f_j(c(t))$ and $t^{-m_j}g_j(c(t))$ are well defined rational functions in $t$. The following claim, which follows from a straightforward computation, implies that $t = 0$ is {\em not} a pole of $t^{-m_j}f_j(c(t))$ or $t^{-m_j}g_j(c(t))$.

\begin{claim} \label{claim:trick}
Fix $j$, $1 \leq j \leq n$. If $p_j$ is a Laurent polynomial supported at $\scrA_j$ then $t^{-m_j}p_j(c(t)) \in \kk[[t]]$. The following identity holds in $\kk[[t]]$:
\begin{align*}
t^{-m_j}p_j(c(t)) = \In_{\scrA_j, \nu}(p_j)(a) + tq_j(t)
\end{align*}
for some $q_j(t) \in \kk[[t]]$. \qed
\end{claim}
Define
\begin{align}
h_j &:=
	t^{-m_j}f_j(c(t))g_j - t^{-m_j}g_j(c(t))f_j \label{trick-deformation}
\end{align}
Let $T$ be the complement in $\kk$ of all the poles of $\prod t^{-2m_j}f_j(c(t))g_j(c(t))$. Then both $0$ and $1$ are in $T$. Moreover, for each $j$,
\begin{itemize}
\item $h_j(x,1) = f_j(b)g_j(x)$ (since $f_j(b) \neq 0 = g_j(b)$), and
\item $h_j(x,0) = -\In_{\scrA_j, \nu}(g_j)(a)f_j(x)$ (this follows from \cref{claim:trick} since $\In_{\scrA_j, \nu}(f_j)(a) = 0 \neq \In_{\scrA_j, \nu}(g_j)(a)$).
\end{itemize}
In particular, $h_j(x,0)$ and $h_j(x,1)$ are (nonzero) constant multiples of respectively $f_j$ and $g_j$. Each $h_j$ vanishes on the curve $C' := \{(c(t),t): t \in T\} \subset \kk^{n+1}$. Note that $C' \cap (\nktorus \times \{1\})$ contains $(b, 1)$ which is an isolated zero of $h_1(x,1), \ldots, h_n(x, 1)$. On the other hand, since $\nu$ is nontrivial, it follows that $C'$ has a ``point at infinity'' at $t = 0$ with respect to $\nktorus \times T$. Therefore assertion \eqref{non-max-condition} of \cref{mult-deformation-global} implies that $\multntorusiso{h_1(x,0)}{h_n(x,0)} < \multntorusiso{h_1(x,\epsilon)}{h_n(x,\epsilon)}$ for generic $\epsilon \in T$. Since each such $h_j(x,\epsilon)$ is supported at $\scrA_j$, it follows that $\multfntorusiso < \multAntorusiso$, as required. \qed

\subsection{The set of non-degenerate systems} \label{bkk-existential-section}
\newcommand{\DJB}{\scrD_{J, \scrB}}
\newcommand{\DJBbar}{\bar \scrD_{J, \scrB}}
\newcommand{\DJBprime}{\scrD'_{J, \scrB}}
\newcommand{\DJBprimebar}{\bar \scrD'_{J, \scrB}}
\newcommand{\DJBzero}{\scrD^0_{J, \scrB}}
\newcommand{\DJBzerobar}{\bar \scrD^0_{J, \scrB}}
\newcommand{\DJBzeroprime}{\scrD'^0_{J, \scrB}}
\newcommand{\DJBzeroprimebar}{\bar \scrD'^0_{J, \scrB}}
\newcommand{\LJ}{\scrL_J}
\newcommand{\LJprime}{\scrL'_J}
\newcommand{\LJprimebar}{\bar \scrL'_J}
\newcommand{\PJ}{\scrP_J}
\newcommand{\PJB}{\scrP_{J, \scrB}}
\newcommand{\UJ}{\scrL^0_J}
\newcommand{\UJprime}{\scrL'^0_J}
\newcommand{\UJprimebar}{\bar \scrL'^0_J}

In this section we prove \cref{bkk-0}. We start with some notation. Let $J \subseteq [m] := \{1, \ldots, m\}$ and $\scrL_J(\scrA) := \prod_{j \in J} \scrL(\scrA_j)$. Write $c_{j,\alpha}$ for the coefficient of $x^\alpha$ for each $j \in J$, so that $(c_{j,\alpha}: j \in J,\  \alpha \in \scrA_j)$ are the coordinates on $\scrL_J(\scrA)$. For $f = (c_{j,\alpha})_{j, \alpha} \in \scrL_J(\scrA)$, we write $f_j$ for the corresponding element in $\scrL(\scrA_j)$, i.e.\ $f_j = \sum_{\alpha \in \scrA_j} c_{j, \alpha} x^\alpha$.  Let $\scrA_J := (\scrA_j : j \in J)$ and $\scrN_J(\scrA)$ be the set of all $f \in \scrL_J(\scrA)$ which are $\scrA_J$-non-degenerate. We will at first show that $\scrN_J (\scrA)$ is Zariski open in $\scrL_J(\scrA)$. If $\scrB = (\scrB_j: j \in J)$ is an ordered tuple of (finite) sets, we say that $\scrB$ is a {\em face} of $\scrA_J$, and write $\scrB  \preceq \scrA_J$, if it satisfies the following property:
\begin{align*}
\parbox{.6\textwidth}{
    ``there is $\nu \in \rnstar$ such that $\scrB_j = \In_\nu(\scrA_j)$ for each $j \in J$.''
}
\end{align*}
For each $\scrB = (\scrB_j: j \in J) \preceq \scrA_J$ and $f = (c_{j,\alpha})_{j, \alpha} \in \scrL_J(\scrA)$, we define $f_{j, \scrB_j}  :=\sum_{\alpha \in \scrB_j} c_{j, \alpha} x^\alpha$. Let $\DJB(\scrA)$ be the set of all $f \in \scrL_J$ such that there is a common root of $f_{j, \scrB_j}$, $j \in J$, on $\nktorus$, so that
\begin{align}
\scrL_J(\scrA) \setminus \scrN_J(\scrA)
	= \bigcup_{\substack{\scrB \preceq \scrA_J \\ \dim(\PJB) < n}}  \DJB(\scrA) \label{DJB}
\end{align}
where $\PJB = \sum_{j \in J} \conv(\scrB_j)$. Let $\LJprime(\scrA) := \scrL_J(\scrA) \times \nktorus$ and $\LJprimebar(\scrA) := \scrL_J(\scrA) \times \pp^n$. Let $\DJBprime(\scrA) \subset \LJprime(\scrA)$ be the collection of all $(f, a)$, where $f \in \DJB(\scrA)$ and $a \in \nktorus$ are such that $f_{j,\scrB_j}(a) = 0$ for each $j \in J$; let $\DJBprimebar$(\scrA) be the closure of $\DJBprime(\scrA) $ in $\LJprimebar$(\scrA). Let $\pi_J: \LJprimebar(\scrA) \to \LJ(\scrA)$ be the natural projection and $\DJBbar(\scrA) := \pi_J(\DJBprimebar(\scrA))$.

\begin{claim} \label{DJBbar-closed}
Let $\scrB \preceq \scrA_J$. Then $\DJBbar(\scrA) \subset \bigcup_{\scrB' \preceq \scrB} \scrD_{J, \scrB'}(\scrA)$.
\end{claim}

\begin{proof}
Let $f^0 = (c^0_{j,\alpha})_{j,\alpha} \in \DJBbar(\scrA) \setminus \DJB(\scrA)$. Pick $a^0 \in \pp^n$ such that $(f^0,a^0) \in \DJBprimebar(\scrA)$. We can find an irreducible curve $C$ on $\DJBprimebar(\scrA)$ such that $(f^0,a^0) \in C$ and $C \cap \DJBprime(\scrA)$ is nonempty and open in $C$ (\cref{closure-curve-prop}). Let $B = (Z,z)$ be a branch of $C$ at $(f^0,a^0)$. For each $j = 1, \ldots, n$, and each $\alpha \in \scrA_j$, we write $\bar x_j, \bar c_{j, \alpha}$ respectively for the restrictions of $x_j, c_{j, \alpha}$ to $C$. Then for each $j \in J$, $\bar F_j := \sum_{\alpha \in \scrB_j} \bar c_{j, \alpha}\bar x^\alpha$ is identically zero on $C$. Since $C' \cap \DJBprime(\scrA) \neq \emptyset$, it follows that no $\bar x_j$ is identically zero on $C$. Let $\nu'_B$ be the element in $\rnstar$ with coordinates $(\nu_B(\bar x_1), \ldots, \nu_B(\bar x_n))$ with respect to the basis dual to the standard basis. Fix $j \in J$. Let $\scrB'_j := \In_{\nu'_B}(\scrB_j)$ and $m_j := \min_{\scrB_J}(\nu'_B) = \langle \nu'_B, \beta \rangle$ for any $\beta \in \scrB'_j$. Then for each $\alpha \in \scrB_j$, since $\ord_z(\bar c_{j,\alpha}) \geq 0$, it follows that
\begin{align*}
\ord_z(\bar c_{j,\alpha}\bar x^\alpha) = \ord_z(\bar c_{j,\alpha}) + \ord_z(\bar x^\alpha) \geq \ord_z(\bar x^\alpha) = \langle \nu'_B, \alpha \rangle \geq m_j
\end{align*}
Moreover, $\ord_z(\bar c_{j,\alpha}\bar x^\alpha) = m_j$ if and only if $ \ord_z(\bar c_{j,\alpha}) = 0$ and $\alpha \in \scrB'_j$, i.e.\ if and only if $c^0_{j, \alpha} \neq 0$ and $\alpha \in \scrB'_j$, i.e.\ if and only if $\alpha \in \supp(f^0_{j,\scrB'_j})$. If $\rho_B$ is a parameter at $B$, it follows that $\bar F_j$ can be expanded in $\kk((\rho_B))$ as
\begin{align*}
\sum_{ \alpha \in \scrB'_j} c^0_{j,\alpha} \prod_{i=1}^n (a^0_i)^{\alpha_i} (\rho_B)^{m_j} + \cdots
\end{align*}
where $a^0_i = \In_B(\bar x_i)$, $i = 1, \ldots, n$, and the omitted terms have higher order in $\rho_B$. Since $\bar F_j$ is identically zero on $C$, it follows that $(a^0_1, \ldots, a^0_n)$ is a zero of $\sum_{\alpha \in\scrB'_j} c^0_{j,\alpha}x^\alpha = f^0_{j,\scrB'_j}$ for each $j \in J$. Therefore $f^0 = (f^0_j: j \in J)$ is $\scrA_J$-degenerate and is an element of $\scrD_{J, \scrB'}(\scrA)$, as required.
\end{proof}

\begin{cor}\label{NJ-open}
For each $\scrB \preceq \scrA_J$, the set $\bigcup_{\scrB' \preceq \scrB} \scrD_{J, \scrB'}(\scrA)$ is Zariski closed in $\scrL_J(\scrA)$. Consequently $\scrN_J(\scrA)$ is Zariski open in $\scrL_J(\scrA)$.
\end{cor}

\begin{proof}
Since $\pp^n$ is {\em complete} (see \cref{complete-section}), it follows that $\DJBbar(\scrA) := \pi_J(\DJBprimebar(\scrA))$ is Zariski closed in $\LJ(\scrA)$. Since the relation $\preceq$ is transitive, \cref{DJBbar-closed}.therefore implies the first assertion. The second assertion then follows from \cref{DJB} (and the transitivity of $\preceq$).
\end{proof}

Let $\scrP_j := \conv(\scrA_j)$, $j = 1, \ldots, n$, and $\PJ := \sum_{j \in J} \scrP_j$. For each $\scrB = (\scrB_j: j \in J) \preceq \scrA_J$, define $\PJB := \sum_{j \in J} \conv(\scrB_j)$. Note that $\PJB$ is a face of $\PJ$.

\begin{claim} \label{DJBbar-proper}
Let $\scrB  \preceq \scrA_J$ be such that $\dim(\PJB) < |J|$. Then $\DJBbar(\scrA) \subsetneqq \LJ(\scrA)$.
\end{claim}

\begin{proof}

We proceed by induction on $|J|$. If $|J| = 1$, then the assumption $\dim(\PJB) < |J|$ is valid only if $\PJB$ is a vertex of $\PJ$, and therefore $\scrB_J$ is a vertex of $\scrA_j$ for each $j$. In that case $f_{j, \scrB_j}$ is a monomial for each $f \in \LJ(\scrA)$, so that $\DJBbar(\scrA) = \emptyset$. Now assume $|J| \geq 2$.  Let $d := \dim(\PJB)$. Since $d < |J|$, there is $J' \subset J$ such that $|J'| = |J|-1$ and $\dim(\scrP_{J', \scrB_{J'}}) =  d$, where $\scrB_{J'} := (\scrB_j: j \in J')$. Define
\begin{align*}
\scrN_{J', \scrB}(\scrA) := \scrL_{J'} (\scrA) \setminus
	\bigcup_{\substack{\scrB' \preceq \scrB_{J'} \\ \dim(\scrP_{J',\scrB'}) < d}} \bar \scrD_{J', \scrB'}(\scrA)
\end{align*}
For each $\bar \scrD_{J', \scrB'}(\scrA)$ appearing in the above union,  we have $\dim(\scrP_{J', \scrB'}) < d \leq |J'|$, so that the inductive hypothesis implies that $\bar \scrD_{J', \scrB'}(\scrA)$ is a proper Zariski closed subset of $\scrL_{J'}(\scrA')$.
Since $\scrL_{J'}(\scrA)$ is irreducible\footnote{It is straightforward to check that $\scrL_{J'}(\scrA) \cong \kk^p$ for some $p\geq 0$.}, it follows that $\scrN_{J', \scrB}(\scrA)$ is a nonempty Zariski open subset of $\scrL_{J'}(\scrA)$. \\

After an integral change of coordinates of $\zz^n$ (i.e.\ a monomial change of coordinates of $\nktorus$) we may arrange that the affine hull of $\PJB$ is of the form $H +  \alpha$ where $\alpha \in \zz^n$ and $H$ is the coordinate hyperplane of $\rr^n$ spanned by the first $d$ unit vectors. Pick $(f_j: j \in J') \in \scrN_{J',\scrB}(\scrA)$. Let $\scrB' \preceq \scrB$. For each $j \in J'$, $f_{j, \scrB'_j}$ is of the form $x^{\alpha_j}g_{j, \scrB'_j}(x_1, \ldots, x_d)$ for some $\alpha_j \in \zz^n$ and Laurent polynomial $g_{j, \scrB'_j}$ in $(x_1, \ldots, x_d)$. Since for each $j \in J'$, the support of $f_{j, \scrB'_j}$ is a translation of the support of $g_{j, \scrB'_j}$, it follows from the definition of $\scrN_{J',\scrB}(\scrA)$ that $(g_{j, \scrB'_j}: j \in J')$ is a BKK non-degenerate system of Laurent polynomials in $(x_1, \ldots, x_d)$.  \Cref{finite-cor} implies that $W_{\scrB'} := V(g_{j, \scrB'_j}: j \in J')$ is a finite set of points in $\nktoruss{d}$. If $j_0$ is the unique element of $J \setminus J'$, then a generic $f_{j_0} \in \scrL(\scrA_{j_0})$ will satisfy the following: for each $\scrB' \preceq \scrB$, $f_{j_0,\scrB'_{j_0}}$ takes a nonzero value at each point of $W_{\scrB'}$. It follows that $f := (f_j: j \in J) \not\in \bigcup_{\scrB' \preceq \scrB} \scrD_{J,\scrB'}(\scrA)$. \Cref{DJBbar-closed} now implies that $f \not\in \DJBbar(\scrA)$, as required.
\end{proof}

\begin{proof}[Proof of \cref{bkk-0}]
Apply \cref{NJ-open} with $J = [m] := \{1, \ldots, m\}$ to see that $\scrN(\scrA)$ is Zariski open in $\scrL(\scrA)$. If $\dim(\sum_{j=1}^m \scrP_j) < m$, then an application of \cref{DJBbar-proper} with $J = [m]$ shows that $\DJBbar(\scrA)$ is a proper Zariski closed subset of $\scrL(\scrA)$ for each face $\scrB$ of $\scrA = (\scrA_1, \ldots, \scrA_m)$. Since $\scrL(\scrA)$ is irreducible, it follows that $\scrL(\scrA) \setminus \bigcup_{\scrB \preceq \scrA} \DJBbar(\scrA)$ is nonempty. Identity \eqref{DJB} then implies that $\scrN(\scrA)$ is nonempty. If on the other hand $m \geq n$, an application of \cref{DJBbar-proper} with $J = [m]$ shows that $\DJBbar(\scrA)$ is a proper Zariski closed subset of $\scrL(\scrA)$ for each face $\scrB$ of $\scrA$ such that $\dim(\scrP_{[m], \scrB}) \leq n - 1$, and it follows similarly that $\scrN(\scrA)$ is nonempty.
\end{proof}

\subsection{Properly non-degenerate systems} \label{properly-non-degenerate-section}
We now introduce a notion of non-degeneracy to be used in \cref{bkk-bound-section}. Let $\scrA := (\scrA_1, \ldots, \scrA_m)$ and $(f_1, \ldots, f_m) \in \scrL(\scrA)$. We say that $f_1, \ldots, f_m$ are \index{Proper!non-degeneracy}\index{Non-degeneracy!proper}{\em properly $\scrA$-non-degenerate} if for every $J \subseteq [m]$ and {\em every}\footnote{Unlike $\scrA$-non-degeneracy, $\nu$ is allowed to be the trivial weighted order $(0, \ldots, 0)$.} weighted order $\nu$ on $\kk[x_1, x_1^{-1}, \ldots, x_n, x_n^{-1}]$ such that $\dim(\In_\nu(\sum_{j \in J} \conv(\scrA_j)) < |J|$, the Laurent polynomials $\In_{\scrA_j, \nu}(f_j)$, $j \in J$, have no common zero in $\nktorus$. If $m \geq \min\{n, \dim(\sum_j \conv(\scrA_j)) + 1\}$, then a properly $\scrA$-non-degenerate system is also $\scrA$-non-degenerate; otherwise there may be properly $\scrA$-non-degenerate systems which are not $\scrA$-non-degenerate. 

\begin{prop} \label{properly-non-degenerate-claim}
The collection $\tilde \scrN(\scrA)$ of properly $\scrA$-non-degenerate systems is a nonempty Zariski open subset of $\scrL(\scrA)$. If $m \geq \min\{n, \dim(\sum_{j=1}^m \conv(\scrA_j)) + 1\}$, then $\tilde \scrN(\scrA) \subset \scrN(\scrA)$.
\end{prop}

\begin{proof}
Indeed, for each $J \subseteq [m]$, let $\tilde \pi_J$ be the natural projection from $\scrL(\scrA)$ onto $\scrL_J(\scrA)$. It is straightforward to check that in the notation of \cref{bkk-existential-section},
\begin{align*}
\tilde \scrN(\scrA) = \scrL(\scrA) \setminus \bigcup_{\substack{J \subseteq [m] \\ \dim(\scrP_J) < |J|} } \bigcup_{\scrB \preceq \scrA_J}\tilde \pi_J^{-1}(\DJB(\scrA))
\end{align*}
\Cref{DJBbar-closed} implies that $\bigcup_{\scrB \preceq \scrA_J}\tilde \pi_J^{-1}(\DJB(\scrA)) = \bigcup_{\scrB \preceq \scrA_J}\tilde \pi_J^{-1}(\DJBbar(\scrA))$, so that $\tilde \scrN(\scrA)$ is Zariski open in $\scrL(\scrA)$. \Cref{DJBbar-proper} implies that $\tilde \scrN(\scrA)$ is nonempty.
\end{proof}

\section{Proof of the BKK bound} \label{bkk-bound-section}
For each $j = 1, \ldots, n$, let $\scrA_j$ be a finite subset of $\zz^n$ and $\scrP_j$ be the convex hull of $\scrA_j$ in $\rr^n$. In this section we show that
\begin{align}
\multAntorusiso = \mv(\scrP_1, \ldots, \scrP_n) \label{mult-A=mv}
\end{align}
In \cref{toric-order-section} we prove two relevant results from the theory of toric varieties, which we use in \cref{mult-A=mv-section} to prove \eqref{mult-A=mv}. Throughout this section we follow the convention of \cref{wt-order-section} to identify weighted orders on the ring of Laurent polynomials in $(x_1, \ldots, x_n)$ with integral elements of $\rnstar$.

\subsection{Toric propositions} \label{toric-order-section}
Let $\nu$ be a primitive integral element in $\rnstar$. For each Laurent polynomial $g$ in $(x_1, \ldots, x_n)$, define $\talphanu(g)$ and $\Inalphapsinu(g)$ as in \cref{O_Q}, where $\alpha_\nu \in \zz^n$ is such that $\langle \nu, \alpha_\nu \rangle = 1$ and $\psi_\nu: \znnuperp \cong \zz^{n-1}$ is an isomorphism (of abelian groups). Let $\scrP$ be an $n$-dimensional convex integral polytope in $\rr^n$ which has a facet $\scrQ$ with primitive inner normal $\nu$. Let $\xp$ be the toric variety corresponding to $\scrP$ from \cref{toric-polysection}, and given Laurent polynomials $g_1, \ldots, g_k$, let $\vonep(g_1, \ldots, g_k)$ be the extension from \cref{xonep-section} of the closed subscheme $V(g_1, \ldots, g_k)$ of $\nktorus$ to a closed subscheme of the Zariski open subset $\xonep$ of $\xp$. Recall that a {\em possibly non-reduced curve} is a pure dimension one closed subscheme of a variety.

\begin{prop} \label{order-nu-prop}
Let $f_1, \ldots, f_n$ be Laurent polynomials in $(x_1, \ldots, x_n)$.
\begin{enumerate}
\item \label{empty-implication} If $\Inalphapsinu(f_2), \ldots,\Inalphapsinu(f_n)$ have no common zero on $\nktoruss{n-1}$, then there is a Zariski open subset $U'$ of $\xp$ containing $\orbit{\scrQ}$ such that the support of $\vonep(f_2, \ldots, f_n) \cap U'$ is empty.

\item  \label{nonempty-implication} If the number of common zeroes of $\Inalphapsinu(f_2), \ldots,\Inalphapsinu(f_n)$ on $\nktoruss{n-1}$ is nonzero and finite, then there is a Zariski open subset $U'$ of $\xp$ containing $\orbit{\scrQ}$ such that $C' := \vonep(f_2, \ldots, f_n) \cap U'$ is a possibly non-reduced curve and every irreducible component of $C'$ intersects $\nktorus$.
\item \label{nonempty-implication-1} If in addition to the assumptions of assertion \eqref{nonempty-implication} $\Inalphapsinu(f_1)$ does {\em not} vanish at any of the common zeroes of $\Inalphapsinu(f_2), \ldots,\Inalphapsinu(f_n)$ on $\nktoruss{n-1}$, then $f_1$ restricts to a nonzero rational function on $C'$ which can be represented in $\local{C'}{a}$ as a quotient of non zero-divisors for every $a \in C'$, and
\begin{align*}
\sum_{ a \in C' \cap \orbit{\scrQ}} \ord_a(f_1|_{C'}) = \nu(f_1) \multnntorus{\Inalphapsinu(f_2)}{\Inalphapsinu(f_n)}{n-1} 
\end{align*}
\end{enumerate}
\end{prop}

\begin{proof}
If $\Inalphapsinu(f_2), \ldots,\Inalphapsinu(f_n)$ have no common zero on $\nktoruss{n-1}$, then \cref{vonepinfinity} implies that $\vonep(f_2, \ldots, f_n) \cap \orbit{\scrQ} = \emptyset$, so that part \eqref{empty-implication} holds with $U' := \xonep \setminus \vonep(f_2, \ldots, f_n)$. Now assume the number of common zeroes of $\Inalphapsinu(f_2), \ldots,\Inalphapsinu(f_n)$ on $\nktoruss{n-1}$ is nonzero and finite. Since by definition $\talphanu(f_j) = x^{-\nu(f_j)\alpha_\nu}f_j $, it follows (e.g.\ due to \cref{vonepinfinity}) that $\talphanu(f_2), \ldots, \talphanu(f_n), x^{\alpha_\nu}$ are $n$ regular functions on the $n$ dimensional variety $U_\scrQ := \xzerop \cup \orbit{\scrQ}$ such that their common zero set $Z$ is finite and nonempty. \Cref{thm:pure-dimension} then implies that $\vonep(f_2, \ldots, f_n)$ has pure dimension one near each point of $Z$. Since $Z$ is precisely the set of points in $\vonep(f_2, \ldots, f_n) \cap \orbit{\scrQ}$, this implies assertion \eqref{nonempty-implication} is satisfied with some $U' \subset U_\scrQ$. Since $U_\scrQ$ is nonsingular, and the set of zeroes of $x^{\alpha_\nu}$ on $U_\scrQ$ is precisely $\orbit{\scrQ}$ (\cref{U_Q}), and since $\orbit{\scrQ}$ does not contain any irreducible component of $C'$, \cref{cor:Macaulay} implies that $x^{\alpha_\nu}|_{C'}$ is a non zero-divisor in $\local{C'}{a}$ for each $a \in C'$. Now fix $a \in C' \cap \orbit{\scrQ}$. Since $a$ corresponds to a common zero of $\Inalphapsinu(f_2), \ldots,\Inalphapsinu(f_n)$ on $\nktoruss{n-1}$ (\cref{vonepinfinity}), under the assumption of assertion \eqref{nonempty-implication-1}, $\talphanu(f_1)$ is a regular function near $a$ which does {\em not} vanish at $a$. Therefore $\talphanu(f_1)|_{C'}$ is a unit in $\local{C'}{a}$, and $f_1 = \talphanu(f_1)(x^{\alpha_\nu})^{\nu(f_1)}$ is the quotient of two non zero-divisors in $\local{C'}{a}$. \Cref{order-properties} then implies that $\ord_a(f_1|_{C'}) = \ord_a(\talphanu(f_1)|_{C'}) + \nu(f_1)\ord_a(x^{\alpha_\nu}|_{C'}) = \nu(f_1)\ord_a(x^{\alpha_\nu}|_{C'})$. On the other hand, assertion \eqref{mult-to-order} of \cref{int-mult-curve} implies that
\begin{align*}
\ord_a(x^{\alpha_\nu}|_{C'})
	&= \multp{x^{\alpha_\nu}, \talphanu(f_2)}{\talphanu(f_n)}{a}
	= \dim_\kk(\local{\xp}{a}/\langle x^{\alpha_\nu}, \talphanu(f_2), \ldots, \talphanu(f_n) \rangle)
\end{align*}
\Cref{vonepinfinity} then implies that
\begin{align*}
\ord_a(x^{\alpha_\nu}|_{C'})
	&= \dim_\kk(\local{\nktoruss{n-1}}{(\psi_\nu)_*(a)}/\langle \Inalphapsinu(f_2), \ldots, \Inalphapsinu(f_n) \rangle) \\
	&= \multp{\Inalphapsinu(f_2)}{\Inalphapsinu(f_n)}{(\psi_\nu)_*(a)}
\end{align*}
where $(\psi_\nu)_*: \orbit{\scrQ} \cong \nktoruss{n-1}$ is the inverse of the isomorphism $\psi_\nu^*: \nktoruss{n-1} \cong \orbit{\scrQ}$ from \cref{O_Q}. Assertion \eqref{nonempty-implication-1} now follows immediately.
\end{proof}

\begin{cor} \label{order-sum-nktorus}
Let $f_1, \ldots, f_n$ be Laurent polynomials in $(x_1, \ldots, x_n)$ and $\scrR$ be the sum of Newton polytopes of $f_2, \ldots, f_n$.
\begin{enumerate}
\item  \label{2-n-dim<n-1} Assume that $\Inalphapsinu(f_2), \ldots, \Inalphapsinu(f_n)$ have no common zero on $\nktoruss{n-1}$ for every primitive integral $\nu \in \rnstar$ such that $\dim(\In_\nu(\scrR)) < n-1$. Then $V(f_2, \ldots, f_n) \cap \nktorus$ is either empty or a curve.
\item \label{2-n-dim=n-1} Assume in addition that $\Inalphapsinu(f_1), \ldots, \Inalphapsinu(f_n)$ have no common zero on $\nktoruss{n-1}$ for every primitive integral $\nu \in \rnstar$ such that $\dim(\In_\nu(\scrR)) = n-1$. Then
\begin{align}
\multfntorus = \multfntorusiso= -\sum_\nu \nu(f_1) \multnntorus{\Inalphapsinu(f_2)} {\Inalphapsinu(f_n)}{n-1} \label{order-sum-nktorus-eqn}
\end{align}
where the sum is over all primitive integral $\nu \in \rnstar$ such that $\dim(\In_\nu(\scrR)) = n-1$.
\end{enumerate}
\end{cor}

\begin{proof}
Let $C'$ be the set of common zeroes of $f_2, \ldots, f_n$ on $\nktorus$. If $\dim(\scrR) < n-1$, then the assumption of assertion \eqref{2-n-dim<n-1} implies that $C' = \emptyset$ and all three sides of \eqref{order-sum-nktorus-eqn} are zero. Therefore assume that $\dim \scrR$ is $n-1$ or $n$. Let $\scrP$ be an $n$-dimensional convex integral polytope in $\rr^n$ which satisfy the following property:
\begin{align}
\parbox{.6\textwidth}{for every $\nu \in \rnstar$, if $\nu$ is an inner normal to a face of $\scrR$ of dimension $n-1$, then $\nu$ is also an inner normal to a facet of $\scrP$.} \label{all-facets-included}
\end{align}
For example, if $\dim(\scrR) = n$, then we can simply take $\scrP = \scrR$, and if $\dim(\scrR) = n-1$, we can take $\scrP$ to be any (convex integral) polytope which has two facets parallel to $\scrR$. Let $\xp$ be the toric variety corresponding to $\scrP$. Recall that we can identify $\nktorus$ with the torus $\xzerop$ of $\xp$ (assertion \eqref{torus-identification} of \cref{xp-thm}). Let $\bar C'$ be the closure of $C'$ in $\xp$. Fix a primitive integral element $\nu \in \rnstar$. Let $\scrQ := \In_\nu(\scrP)$ and $\orbit{\scrQ}$ be the corresponding torus orbit on $\xp$.

\begin{proclaim} \label{C-Q}
Let $\scrQ' := \In_\nu(\scrR)$.
\begin{enumerate}
\item \label{Q'<n-1} If $\dim(\scrQ') < n-1$, then $V(\Inalphapsinu(f_2), \ldots, \Inalphapsinu(f_n)) \cap \nktoruss{n-1} = \emptyset$, and $\bar C' \cap \orbit{\scrQ} = \emptyset$.
\item \label{Q'=n-1} If $\dim(\scrQ') = n-1$, then $\scrQ$ is a facet of $\scrP$, and $V(\Inalphapsinu(f_2), \ldots, \Inalphapsinu(f_n)) \cap \nktoruss{n-1}$ is finite.
\end{enumerate}

\end{proclaim}

\begin{proof}
If $\dim(\scrQ') < n-1$, then the hypothesis of the corollary implies that the set of common zeroes of $\Inalphapsinu(f_2), \ldots, \Inalphapsinu(f_n)$ on $\nktorus$ is empty. Assertion \eqref{Q'<n-1} of the claim then follows from \cref{closure-containment}. On the other hand, if $\dim(\scrQ') = n-1$, then property \eqref{all-facets-included} implies that $\scrQ$ is a facet of $\scrP$, and the hypothesis of the corollary implies that $\Inalphapsinu(f_2), \ldots, \Inalphapsinu(f_n)$ are BKK non-degenerate on $\nktoruss{n-1}$. \Cref{finite-cor} then implies that $V(\Inalphapsinu(f_2), \ldots, \Inalphapsinu(f_n)) \cap \nktoruss{n-1}$ is finite, as required.
\end{proof}

Let $C := \vonep(f_2, \ldots, f_n)$. Recall that $C$ is an extension of the closed subscheme of $\kk^n$ determined by $f_2, \ldots, f_n$ to the open subset of $\xp$ obtained by removing $\orbit{\scrQ}$ for all face $\scrQ$ of $\scrP$ with dimension $\leq n-2$. Assertion \eqref{Q'<n-1} of \cref{C-Q} implies that $\bar C'  \subset \supp(C)$. If $\Inalphapsinu(f_2), \ldots, \Inalphapsinu(f_n)$ have no common zero on $\nktoruss{n-1}$ for each primitive integral $\nu \in \rnstar$, then assertion \eqref{empty-implication} of \cref{order-nu-prop} implies that $\bar C' \cap \xpinfty = \emptyset$, where $\xpinfty := \xp \setminus \nktorus$. On the other hand $\bar C'$ is a {\em complete} variety if it is nonempty, so that it can {\em not} be completely contained in the affine variety $\nktorus$ (\cref{example:regular-on-complete}). It  follows that $C' = \emptyset$ and all sides of \eqref{order-sum-nktorus-eqn} are zero. So assume there is primitive integral $\nu \in \rnstar$ such that $V(\Inalphapsinu(f_2), \ldots, \Inalphapsinu(f_n)) \cap \nktoruss{n-1} \neq \emptyset$. Then assertion \eqref{Q'=n-1} of \cref{C-Q} and \cref{order-nu-prop} imply that $C$ has dimension one near $\xpinfty$, and $\supp(C) \subset \bar C'$. It follows that $\supp(C) = \bar C'$, so that $C$ is a {\em possibly non-reduced curve} such that $\supp(C)$ is {\em projective}. Since $C \cap \nktorus$ is the closed subscheme of $\nktorus$ defined by $f_2, \ldots, f_n$, \cref{int-mult-curve} implies that
\begin{align*}
\multfntorus
	&=  \sum_{a \in C \cap \nktorus} \ord_a(f_1|_C)
\end{align*}
Since $\supp(C)$ is projective, \cref{zero-sum,order-nu-prop,C-Q} then imply that
\begin{align*}
\multfntorus
	= -\sum_{a \in C\cap \xpinfty} \ord_a(f_1|_C)
	= -\sum_\nu \nu(f_1) \multnntorus{\Inalphapsinu(f_2)} {\Inalphapsinu(f_n)}{n-1}
\end{align*}
as required.
\end{proof}

\subsection{Proof of identity \eqref{mult-A=mv}} \label{mult-A=mv-section}

Define
\begin{align}
\multPntorusiso := \max\{\multfntorusiso: \supp(f_j) \subseteq \scrP_j,\ j =1 , \ldots, n\} \label{multPntorusiso}
\end{align}
Let $\scrA'_j := \scrP_j \cap \zz^n \supseteq \scrA_j$, $j = 1, \ldots, n$. It follows from the definition that $\multPntorusiso = \multntorusiso{\scrA'_1}{\scrA'_n}\geq \multAntorusiso$.  However, due to \cref{bkk-0} we may pick BKK non-degenerate Laurent polynomials $f_1, \ldots, f_n$ such that $\supp(f_j) \subseteq \scrA_j \subseteq \scrA'_j$ and $\np(f_j) = \scrP_j$ for each $j$, and then \cref{bkk-sufficiency} implies that
$$\multAntorusiso = \multfntorusiso = \multntorusiso{\scrA'_1}{\scrA'_n} = \multPntorusiso$$
In order to prove \eqref{mult-A=mv} it therefore suffices to show that
\begin{align}
\multPntorusiso = \mv(\scrP_1, \ldots, \scrP_n) \label{mult-P=mv}
\end{align}

\begin{claim} \label{additive-claim}
Let $\scrK$ be the set of convex integral polytopes in $\rr^n$ regarded as a semigroup under Minkowski addition. The function $\scrK^n \to \rr$ given by $(\scrQ_1, \ldots, \scrQ_n) \mapsto \multntorusiso{\scrQ_1}{\scrQ_n}$ is symmetric and multiadditive.
\end{claim}

\begin{proof}
The symmetry is evident, so we prove the multiadditivity. Pick $\scrQ_1, \ldots, \scrQ_n, \scrQ'_1 \in \scrK$. \Cref{bkk-0} implies that we may choose Laurent polynomials $g_1, \ldots, g_n, g'_1$ such that $\np(g_j) = \scrQ_j$, $j =1 , \ldots, n$, $\np(g'_1) = \scrQ'_1$, and both $g_1, \ldots, g_n$ and $g'_1, g_2, \ldots, g_n$ are BKK non-degenerate. But then $g_1g'_1, g_2, \ldots, g_n$ are also BKK non-degenerate; in particular, $g_1g'_1, g_2, \ldots, g_n$ only have isolated zeroes on $\nktorus$ (\cref{finite-cor}). It follows that
\begin{align*}
\multntorusiso{\scrQ_1 + \scrQ'_1, \scrQ_2}{\scrQ_n}
	&= \multntorusiso{g_1g'_1,g_2}{g_n}
	\ \text{(\cref{bkk-sufficiency})}\\
	&= \multntorusiso{g_1}{g_n} + \multntorusiso{g'_1,g_2}{g_n}
	\ \text{(\cref{int-mult-curve}, assertion \eqref{mult-additive})} \\
	&= \multntorusiso{\scrQ_1}{\scrQ_n} + \multntorusiso{ \scrQ'_1, \scrQ_2}{\scrQ_n}
	\ \text{(\cref{bkk-sufficiency})}
\end{align*}
as required.
\end{proof}

Due to \cref{mixed-unique,additive-claim}, in order to prove \eqref{mult-P=mv} it suffices to show that $\multntorusiso{\scrP}{\scrP} = n!\vol_n(\scrP)$ for each convex integral polytope $\scrP$ in $\rr^n$. We proceed by induction on $n$. It is clearly true for $n = 1$, so assume it is true for $n - 1$. Pick a convex integral polytope $\scrP$ in $\rr^n$. Let $f_1, \ldots, f_n$ be {\em properly} $\scrA$-non-degenerate Laurent polynomials (see \cref{properly-non-degenerate-section}) such that the Newton polytope of each $f_j$ is $\scrP$. In particular $f_1, \ldots, f_n$ are BKK non-degenerate (\cref{properly-non-degenerate-claim}). Therefore \cref{finite-cor,bkk-sufficiency} imply that
\begin{align*}
\multntorusiso{\scrP}{\scrP}  = \multfntorus
\end{align*}
Since the $f_j$ are properly non-degenerate, they satisfy the hypothesis of \cref{order-sum-nktorus}, and identity \eqref{order-sum-nktorus-eqn} implies that
\begin{align*}
\multfntorus = -\sum_\nu \nu(f_1)\multnntorus{\Inalphapsinu(f_2)}{\Inalphapsinu(f_n)}{n-1}
\end{align*}
where the sum is over all primitive integral $\nu \in \rnstar$. Fix one such $\nu$. The proper $\scrA$-non-degeneracy of the $f_j$ implies that $\Inalphapsinu(f_2), \ldots, \Inalphapsinu(f_n)$ are BKK non-degenerate on $\nktoruss{n-1}$. Moreover, the Newton polytope of each $\Inalphapsinu(f_j)$ is $\scrP_\nu$, which is the convex hull in $\rr^{n-1}$ of $\{\psi_\nu(\beta-m_\nu\alpha_\nu): \beta \in \In_\nu(\scrP) \cap \zz^n\}$, where $m_\nu := \min_\scrP(\nu)$. The inductive hypothesis implies that
\begin{align*}
\multnntorus{\Inalphapsinu(f_2)} {\Inalphapsinu(f_n)}{n-1} = (n-1)! \vol_{n-1}(\scrP_\nu)
\end{align*}
It follows from the definition of $\vol'_\nu(\cdot)$ (see \cref{vol'}) that $\vol_{n-1}(\scrP_\nu) = \vol'_{\nu}(\In_\nu(\scrP))$. It follows that
\begin{align*}
\multfntorus
	= -(n-1)!\sum_\nu \min_\scrP(\nu) \vol'_{\nu}(\In_\nu(\scrP))
	=  (n-1)!\sum_\nu \max_\scrP(\nu) \vol'_{\nu}(\ld_\nu(\scrP))
	= n! \vol_n(\scrP)
\end{align*}
where the last equality uses \cref{volume-facet-rational}. Therefore, $\multntorusiso{\scrP}{\scrP} = n!\vol_n(\scrP)$, as required. \qed


\section{Applications of Bernstein's theorem to convex geometry} \label{bkk-convex-applisection}
In this section we use Bernstein's theorem to deduce some properties of mixed volume. In particular we characterize the conditions for mixed volume (of $n$ convex polytopes in $\rr^n$) being zero (\cref{positively-mixed}) and the conditions under which it is ``strictly monotonic'' (\cref{strictly-mixed-monotone}). As an application back to algebraic geometry, we prove the alternate version of Bernstein's theorem (\cref{alternate-corollary}). Throughout this section $\scrP_1, \ldots, \scrP_n$ denote convex polytopes in $\rr^n$, $n \geq 1$.

\begin{prop}[Monotonicity of mixed volume] \label{mixed-monotone}
If $\scrP'_j$ are convex polytopes in $\rr^n$ such that  $\scrP_j \subseteq \scrP'_j$ for each $j$, then $\mv(\scrP_1, \ldots, \scrP_n) \leq \mv(\scrP'_1, \ldots, \scrP'_n)$.
\end{prop}

\begin{proof}
\Cref{bkk-bound-thm} implies that it holds for rational polytopes. The general case then follows from the observation that every polytope can be approximated arbitrarily closely (with respect to the Hausdorff distance) by rational polytopes (\cref{poly-rational-appxn}), and the mixed volume is continuous with respect to the Hausdorf distance (\cref{mixed-remark}).
\end{proof}


Let $\nu$ be a nonzero integral element of $\rnstar$ and $\rnnuperp := \{\alpha \in \rr^n: \langle \nu, \alpha \rangle = 0\}$. Choose an affine transformation $\psi_\nu: \rr^n \to \rr^n$ such that $\psi_\nu$ restricts to an automorphism of $\zz^n$ and maps $\rnnuperp \cap \zz^n$ onto $\zz^{n-1} \times \{0\}$. If $\scrQ_1, \ldots, \scrQ_{n-1}$ are rational polytopes in $\rr^n$ such that each $\scrQ_j$ is a translate of some polytope $\scrQ'_j \subset \rnnuperp$, then we define
\begin{align}
\mv'_\nu(\scrQ_1, \ldots, \scrQ_{n-1})
	:=
	\mv(\psi_\nu(\scrQ'_1), \ldots, \psi_\nu(\scrQ'_{n-1})) \label{mv'}
\end{align}
where the mixed volume on the right hand side is the $(n-1)$-dimensional mixed volume on $\rr^{n-1}$. \Cref{prop:fund-vol,mixed-unique} imply that $\mv'_\nu$ does not depend on the choice of $\psi_\nu$ or the translations involved.

\begin{prop}  \label{rational-mv'-prop}
Assume $\scrP_1, \ldots,  \scrP_n$ are rational polytopes. Then
\begin{enumerate}
\item \label{rational-mv'} $\mv(\scrP_1, \ldots, \scrP_n) = \sum_\nu \max_{\scrP_1}(\nu) \mv'_\nu(\ld_\nu(\scrP_2), \ldots, \ld_{\nu}(\scrP_n))$, where the sum is over all primitive integral $\nu \in \rnstar$.
\item \label{line-mv'} Assume $\scrP_1$ is a line segment in the direction of $\nu$, where $\nu$ is a primitive integral element in $\rr^n$. Let $l(\scrP_1)$ be the ``integer length'' of $\scrP_1$ (i.e.\ the Euclidean length of $\scrP_1$ is $l(\scrP_1)$ times the length of $\nu$). Identify $\nu$ with an element of $\rnstar$ via the basis dual to the standard basis of $\rr^n$. Let $\pi_\nu: \rr^n \to \rnnuperp$ be a ``lattice projection in the direction normal to $\nu$'' (i.e.\ $\pi_\nu = \psi_\nu^{-1} \circ \pi \circ \psi_\nu$, where $\pi: \rr^n \to \rr^{n-1} \times \{0\}$ is the projection in the first $(n-1)$-coordinates). Then $\mv(\scrP_1, \ldots, \scrP_n) = l(\scrP_1) \mv'_{\nu}(\pi_\nu(\scrP_2), \ldots, \pi_\nu(\scrP_n))$.
\end{enumerate}
\end{prop}

\begin{proof}
Due to the multiadditivity of the mixed volume we may assume that each $\scrP_j$ is integral. Pick Laurent polynomials $f_j$ with Newton polytope $\scrP_j$ such that $f_1, \ldots, f_n$ are properly non-degenerate. Then they are BKK non-degenerate on $\nktorus$, they satisfy the hypothesis of \cref{order-sum-nktorus}, and for each primitive integral $\nu \in \rnstar$, $\Inalphapsinu(f_2), \ldots, \Inalphapsinu(f_n)$ are BKK non-degenerate on $\nktoruss{n-1}$. Then \cref{order-sum-nktorus,bkk-bound-thm,bkk-non-degenerate-thm} implies assertion \eqref{rational-mv'}. For assertion \eqref{line-mv'}, change coordinates on $\zz^n$ (using \cref{basis-lemma}) so that $\nu = (0, \ldots, 0, 1)$. \Woutlog\ we may assume that $\scrP_1$ is the line segment bounded by the origin and $(0, \ldots, 0, l)$, where $l := l(\scrP_1)$, and in addition, $f_1 = x_n^l + 1$. Then for generic $a = (a_1, \ldots, a_n) \in \nktorus$,
\begin{align*}
\multf{a} 
	= \multp{x_n^l + 1 - a_1,f_2}{f_n}{a}
	= \sum_{\epsilon^l = a_1 - 1} \multp{f_2|_{x_n = \epsilon}}{f_n|_{x_n = \epsilon}}{(a_2, \ldots, a_n)}
\end{align*}
Since $\pi_\nu$ is simply the projection in the first $(n-1)$-coordinates, for generic $a_1$, the Newton polytope of $f_j|_{x_n = \epsilon}$ is $\pi_\nu(\scrP_j)$ for each $j$. Assertion \eqref{line-mv'} therefore follows from \cref{bkk-bound-thm}.
\end{proof}

\def\shiftone{7.5}
\def\colorzero{blue}
\def\opazero{0.5}
\def\viewx{75}
\def\viewy{30}
\def\titlex{1}
\def\titley{-1}

\begin{figure}[htp]
\begin{center}
\begin{tikzpicture}[scale=0.6]
\pgfplotsset{every axis/.append style = {view={\viewx}{\viewy}, axis lines=middle, enlargelimits={upper}}}
\begin{scope}
\begin{axis}
	\addplot3 [dashed, black, thick] coordinates{(0,1,0) (0,0,1)};
	\addplot3[fill=red,opacity=\opazero] coordinates{(1,0,0) (0,0,1) (0,0,2)};
	\addplot3[fill=\colorzero,opacity=\opazero] coordinates{(1,0,0) (0,1,0) (0,0,2)};
\end{axis}
\draw (\titlex,\titley) node {$\scrP_1$};
\end{scope}

\begin{scope}[shift={(\shiftone,0)}]
\begin{axis}
	\addplot3[fill=\colorzero,opacity=\opazero] coordinates{(3,0,0) (0,3,0) (0,1,2) (1,0,2)};
\end{axis}
\draw (\titlex,\titley) node {$\scrP_2$};
\end{scope}

\begin{scope}[shift={(2*\shiftone,0)}]
\begin{axis}
	\addplot3[fill=\colorzero,opacity=\opazero] coordinates{(2,0,0) (0,2,0) (0,1,2) (1,0,2)};
\end{axis}
\draw (\titlex,\titley) node {$\scrP_3$};
\end{scope}
\end{tikzpicture}
\end{center}
\caption{Newton polytopes of polynomials from \cref{ex-mv'}} \label{fig-ex-mv'}
\end{figure}

\begin{example}\label{ex-mv'}
Let $f_1 = a_1x + b_1y + c_1z + d_1z^2$, $f_2 = a_2x^3 + b_2xz^2 + c_2y^3 + d_2yz^2$, $f_3 = a_3x^2 + b_3xz^2 + c_3y^2 + d_3yz^2$, where $a_i, b_i, c_i, d_i$'s are generic elements of $\kk$, and let $\scrP_j := \np(f_j)$, $j = 1, 2,3$ (see \cref{fig-ex-mv'}). We compute $\multpnodots{f_1, f_2, f_3}{\nktoruss{3}} = \mv(\scrP_1, \scrP_2, \scrP_3)$ using \cref{rational-mv'-prop}. Assertion \eqref{rational-mv'} of \cref{rational-mv'-prop} implies that $ \mv(\scrP_1, \scrP_2, \scrP_3) = \sum_\nu \max_{\scrP_1}(\nu)\mv'_\nu(\ld_\nu(\scrP_2), \ld_\nu(\scrP_3))$. \Cref{positively-mixed} implies that $\mv'_\nu(\ld_\nu(\scrP_2), \ld_\nu(\scrP_3))$ is nonzero only if $\nu$ is one of the six outer normals of facets of $\scrP_2 + \scrP_3$ (\cref{fig-ex-mv'-2}). When $\nu = (-1,0,0)$ or $(0,-1,0)$, then $\max_{\scrP_1}(\nu) = 0$, so it suffices to consider the remaining four cases. The image of the leading faces of $\scrP_2, \scrP_3$ and $\scrP_2 + \scrP_3$ under (certain choices of) $\psi_\nu$ are given in \cref{fig:ex-mv'-3}, and \cref{mixed-example} implies that $\mv'(\ld_\nu(\scrP_2), \ld_\nu(\scrP_3))$ is the area of the region shaded green inside $\psi_\nu(\ld_{\nu}(\scrP_2 + \scrP_3))$. It then follows from \cref{fig:ex-mv'-3} that

\def\scalefactor{0.3}
\def\colorzero{blue}
\def\colorone{green}

\begin{align*}
\mv(\scrP_1, \scrP_2, \scrP_3)
	&= \max_{\scrP_1}(1,1,1) \cdot \area(
				\begin{tikzpicture}[scale=\scalefactor]
						\draw[fill=\colorone, opacity=\opazero] (3,0) -- (4,0) -- (2,2) -- (1,2) -- cycle;
				\end{tikzpicture}
				)
			+  \max_{\scrP_1}(-1,-1,-1) \cdot \area(
				\begin{tikzpicture}[scale=\scalefactor]
						\draw[fill=\colorone, opacity=\opazero] (3,0) -- (5,0) -- (3,2) -- (1,2) -- cycle;
				\end{tikzpicture}
				) \\
	& \qquad \qquad
			+  \max_{\scrP_1}(2,2,1) \cdot \area(
				\begin{tikzpicture}[scale=\scalefactor]
						\draw[fill=\colorone, opacity=\opazero] (2,0) -- (5,0) -- (4,1) -- (1,1) -- cycle;
				\end{tikzpicture}
				)
			+  \max_{\scrP_1}(-2,-2,-1) \cdot \area(
				\begin{tikzpicture}[scale=\scalefactor]
						\draw[fill=\colorone, opacity=\opazero] (2,0) -- (3,0) -- (2,1) -- (1,1) -- cycle;
				\end{tikzpicture}
				) \\
	&= 2 \cdot 2 - 1\cdot4 + 2\cdot 3 -1\cdot1 \\
	&= 5
\end{align*}
\end{example}

\def\shiftone{7.5}
\def\colorzero{blue}
\def\opazero{0.5}
\def\viewx{70}
\def\viewy{5}
\def\titlex{1}
\def\titley{-1}

\begin{figure}[thp]
\begin{center}
\begin{tikzpicture}[scale=0.6]
\pgfplotsset{every axis/.append style = {view={\viewx}{\viewy}}}
\begin{scope}
\begin{axis}[axis lines=middle, enlargelimits={upper}]
	\addplot3 [fill=blue,opacity=\opazero, thick] coordinates{(5,0,0) (0,5,0) (0,4,2) (4,0,2) (5,0,0)};
	\addplot3 [fill=yellow,opacity=\opazero, thick] coordinates{(4,0,2) (0,4,2) (0,2,4) (2,0,4) (4,0,2)};
	\addplot3 [fill=red,opacity=\opazero, thick] coordinates{(2,0,4) (3,0,2) (5,0,0) (4,0,2)};
	\addplot3 [dashed, black, thick] coordinates{(2,0,4) (0,2,4) (0,3,2) (3,0,2) (2,0,4)};
	\addplot3 [dashed, black, thick] coordinates{(3,0,2) (0,3,2) (0,5,0) (5,0,0) (3,0,2)};
\end{axis}
\draw (\titlex,\titley) node {$\scrP_2 + \scrP_3$};
\end{scope}

\def\factorone{3}
\def\factortwo{3}
\def\factorthree{1}
\begin{scope}[shift={(\shiftone,0)}]
\begin{axis}[hide axis, enlargelimits={0.15}]
	\addplot3 [ultra thick, blue, ->] coordinates{(0,0,0) (2,2,1)};
	\addplot3 [ultra thick, ->] coordinates{(0,0,0) (-2*\factorthree,-2*\factorthree,-\factorthree)};
	\addplot3 [ultra thick, yellow, ->] coordinates{(0,0,0) (1,1,1)};
	\addplot3 [ultra thick, ->] coordinates{(0,0,0) (-\factorthree,-\factorthree,-\factorthree)};	
	\addplot3 [ultra thick, red, ->] coordinates{(0,0,0) (0,-\factorone,0)};
	\addplot3 [ultra thick, ->] coordinates{(0,0,0) (-\factortwo,0,0)};
	\draw (axis cs:2,2,1) node [above] {(2,2,1)};
	\draw (axis cs:1,1,1) node [left] {(1,1,1)};
	\draw (axis cs:0,-\factorone,0) node [below] {(0,-1,0)};
	\draw (axis cs: -\factortwo,0,0) node [above] {(-1,0,0)};
	\draw (axis cs:-\factorthree,-\factorthree,-\factorthree) node [right] {(-1,-1,-1)};
	\draw (axis cs:-2*\factorthree,-2*\factorthree,-\factorthree) node [below] {(-2,-2,-1)};
\end{axis}
\draw (\titlex,\titley) node [right, text width=3.5cm] {outer normals to facets of $\scrP_2 + \scrP_3$};
\end{scope}

\end{tikzpicture}
\end{center}
\caption{$\scrP_2+ \scrP_3$ and the outer normals to its facets} \label{fig-ex-mv'-2}
\end{figure}

\def\colorzero{blue}
\def\colorone{green}

\begin{figure}[h]
\begin{center}
\begin{tikzpicture}[scale=0.4]
\def\shiftone{4}
\def\shifttwo{5}
\def\nux{-0.5}
\def\nuy{3.5}
\def\plusy{1.5}
\def\opazero{0.5}
\def\tx{-0.5}
\def\ty{-0.5}
\def\gridx{3.5}
\def\gridxx{5.5}
\def\gridy{2.5}
\def\bigshiftdhor{3}
\def\bigshiftdver{3}

\draw (\nux,\nuy) node [right] {\small $\nu = (1,1,1)$};
\draw (\shiftone,\plusy) node {$+$};    	
\draw (\shiftone+\shifttwo,\plusy) node {$=$};    	

\draw [gray,  line width=0pt] (-0.5,-0.5) grid (\gridx,\gridy);

\draw[ultra thick, fill=\colorzero, opacity=\opazero ] (0,0) -- (3,0) -- (1,2) -- (0,2) -- cycle;
\draw (\tx,\ty) node [below right] {\small $\ld_\nu(\scrP_2)$};    	

\begin{scope}[shift={(\shifttwo,0)}]
	\draw [gray,  line width=0pt] (-0.5,-0.5) grid (\gridx,\gridy);
	
	\draw[ultra thick, \colorzero] (0,0) -- (1,0);
	\draw (\tx,\ty) node [below right] {\small $\ld_\nu(\scrP_3)$};    	
\end{scope}

\begin{scope}[shift={(2*\shifttwo,0)}]
	\draw [gray,  line width=0pt] (-0.5,-0.5) grid (\gridxx,\gridy);
	
	\draw[ultra thick] (0,0) -- (4,0) -- (2,2) --  (0,2) --  cycle;
	\draw[fill=\colorone, opacity=\opazero] (3,0) -- (4,0) -- (2,2) -- (1,2) -- cycle;
	\draw[fill=\colorzero, opacity=\opazero, dashed, ultra thick ] (0,0) -- (3,0) -- (1,2) -- (0,2) -- cycle;
	
	\draw (\tx,\ty) node [below right] {\small $\ld_\nu(\scrP_2 + \scrP_3)$};    	
\end{scope}

\begin{scope}[shift={(2*\shifttwo+\gridxx+\bigshiftdhor,0)}]
	\draw (\nux,\nuy) node [right] {\small $\nu = (-1,-1,-1)$};
	\draw (\shiftone,\plusy) node {$+$};    	
	\draw (\shiftone+\shifttwo,\plusy) node {$=$};    	
	
	\draw [gray,  line width=0pt] (-0.5,-0.5) grid (\gridx,\gridy);
	
	\draw[ultra thick, fill=\colorzero, opacity=\opazero ] (0,0) -- (3,0) -- (1,2) -- (0,2) -- cycle;
	\draw (\tx,\ty) node [below right] {\small $\ld_\nu(\scrP_2)$};    	
	
	\begin{scope}[shift={(\shifttwo,0)}]
		\draw [gray,  line width=0pt] (-0.5,-0.5) grid (\gridx,\gridy);
		
		\draw[ultra thick, \colorzero] (0,0) -- (2,0);
		\draw (\tx,\ty) node [below right] {\small $\ld_\nu(\scrP_3)$};    	
	\end{scope}
	
	\begin{scope}[shift={(2*\shifttwo,0)}]
		\draw [gray,  line width=0pt] (-0.5,-0.5) grid (\gridxx,\gridy);
		
		\draw[ultra thick] (0,0) -- (5,0) -- (3,2) --  (0,2) --  cycle;
		\draw[fill=\colorone, opacity=\opazero] (3,0) -- (5,0) -- (3,2) -- (1,2) -- cycle;
		\draw[fill=\colorzero, opacity=\opazero, dashed, ultra thick ] (0,0) -- (3,0) -- (1,2) -- (0,2) -- cycle;
		
		\draw (\tx,\ty) node [below right] {\small $\ld_\nu(\scrP_2 + \scrP_3)$};    	
	\end{scope}
\end{scope}		

\begin{scope}[shift={(0,-\nuy-\bigshiftdver)}]
\draw (\nux,\nuy) node [right] {\small $\nu = (2,2,1)$};
\draw (\shiftone,\plusy) node {$+$};    	
\draw (\shiftone+\shifttwo,\plusy) node {$=$};    	

\draw [gray,  line width=0pt] (-0.5,-0.5) grid (\gridx,\gridy);

\draw[ultra thick, \colorzero] (0,0) -- (3,0);

\draw (\tx,\ty) node [below right] {\small $\ld_\nu(\scrP_2)$};    	

\begin{scope}[shift={(\shifttwo,0)}]
	\draw [gray,  line width=0pt] (-0.5,-0.5) grid (\gridx,\gridy);
	\draw[ultra thick, fill=\colorzero, opacity=\opazero] (0,0) -- (2,0) -- (1,1) -- (0,1) --cycle;
	\draw (\tx,\ty) node [below right] {\small $\ld_\nu(\scrP_3)$};    	
\end{scope}

\begin{scope}[shift={(2*\shifttwo,0)}]
	\draw [gray,  line width=0pt] (-0.5,-0.5) grid (\gridxx,\gridy);
	
	\draw[ultra thick] (0,0) -- (5,0) -- (4,1) --  (0,1) --  cycle;
	\draw[fill=\colorone, opacity=\opazero] (2,0) -- (1,1) -- (4,1) -- (5,0) -- cycle;
	\draw[fill=\colorzero, opacity=\opazero, dashed, ultra thick ] (0,0) -- (2,0) -- (1,1) -- (0,1) --cycle;
	
	\draw (\tx,\ty) node [below right] {\small $\ld_\nu(\scrP_2 + \scrP_3)$};    	
\end{scope}

\begin{scope}[shift={(2*\shifttwo+\gridxx+\bigshiftdhor,0)}]
	\draw (\nux,\nuy) node [right] {\small $\nu = (-2,-2,-1)$};
	\draw (\shiftone,\plusy) node {$+$};    	
	\draw (\shiftone+\shifttwo,\plusy) node {$=$};    	
	
	\draw [gray,  line width=0pt] (-0.5,-0.5) grid (\gridx,\gridy);
	
	\draw[ultra thick, \colorzero] (0,0) -- (1,0);
	\draw (\tx,\ty) node [below right] {\small $\ld_\nu(\scrP_2)$};    	
	
	\begin{scope}[shift={(\shifttwo,0)}]
		\draw [gray,  line width=0pt] (-0.5,-0.5) grid (\gridx,\gridy);
		
		\draw[ultra thick, fill=\colorzero, opacity=\opazero] (0,0) -- (2,0) -- (1,1) -- (0,1) --cycle;
		\draw (\tx,\ty) node [below right] {\small $\ld_\nu(\scrP_3)$};    	
	\end{scope}
	
	\begin{scope}[shift={(2*\shifttwo,0)}]
		\draw [gray,  line width=0pt] (-0.5,-0.5) grid (\gridxx,\gridy);
		
		\draw[ultra thick] (0,0) -- (3,0) -- (2,1) --  (0,1) --  cycle;
		\draw[fill=\colorone, opacity=\opazero] (2,0) -- (1,1) -- (2,1) -- (3,0) -- cycle;
		\draw[fill=\colorzero, opacity=\opazero, dashed, ultra thick ] (0,0) -- (2,0) -- (1,1) -- (0,1) --cycle;
		
		\draw (\tx,\ty) node [below right] {\small $\ld_\nu(\scrP_2 + \scrP_3)$};    	
	\end{scope}
\end{scope}		
\end{scope}
\end{tikzpicture}
\caption{The image under $\psi_\nu$ of leading faces of $\scrP_2, \scrP_3$ and $\scrP_2 + \scrP_3$}  \label{fig:ex-mv'-3}
\end{center}
\end{figure}

\begin{defn} \label{dependent-defn}
\index{Dependence!of polytopes}\index{Independence!of polytopes}\index{Polytope!dependent}\index{Polytope!independent}
We say that convex polytopes $\scrQ_1, \ldots, \scrQ_m$ in $\rr^n$ are {\em dependent} if there is a nonempty subset $I$ of $[m] := \{1, \ldots, m\}$ such that $\dim(\sum_{i \in I} \scrQ_i) < |I|$; otherwise we say that they are {\em independent}. In particular if $m \geq 1$ and $\scrQ_j = \emptyset$ for some $j$, then $\scrQ_1, \ldots, \scrQ_m$ are dependent. 
\end{defn}

\begin{thm}[Minkowski] \label{positively-mixed}
$\mv(\scrP_1, \ldots, \scrP_n) = 0$ if and only if they are dependent.
\end{thm}

\begin{proof}
Due to \cref{poly-rational-appxn,mixed-remark} and the multiadditivity of mixed volumes, it suffices to consider the case that each $\scrP_j$ is integral. At first assume there is $I \subseteq [n]$ such that $\dim(\sum_{i \in I} \scrP_i) < |I|$. A recursive application of assertion \eqref{rational-mv'} of \cref{rational-mv'-prop} shows that $\mv(\scrP_1, \ldots, \scrP_n)$ can be expressed as a sum such that each summand has a multiplicative factor of the form $\mv'_\nu(\ld_\nu(\scrP'_i): i \in I)$, where $\mv'$ denotes an $|I|$-dimensional mixed volume and $\scrP'_i$ is a face of $\scrP_i$ for each $i \in I$. Now pick BKK non-degenerate Laurent polynomials $f_i$, $i \in I$, such that $\np(f_i) = \scrP_i$ for each $i$ (such $f_i$ exist due to \cref{bkk-0}). Since $\dim(\sum_{i \in I} \scrP_i) < |I|$, it follows from the definition of BKK non-degeneracy that there is no common zero of $\ld_\nu(f_i)$, $i \in I$, on $\nktorus$. \Cref{bkk-bound-thm} then implies that $\mv'_\nu(\ld_\nu(\scrP'_i): i \in I) = 0$. This in turn implies that $\mv(\scrP_1, \ldots, \scrP_n) = 0$. Now assume that $\scrP_1, \ldots, \scrP_n$ are independent. We will show that $\mv(\scrP_1, \ldots, \scrP_n) > 0$. We proceed by induction on $n$. The case of $n = 1$ is obvious. In the general case, since $\dim(\scrP_n) \geq 1$, after a change of coordinates on $\zz^n$ if necessary, we may assume that it has positive length along $x_n$-axis. Let $\pi: \zz^n \to \zz^{n-1}$ be the projection in the first $(n-1)$-coordinates. We consider two cases:\\

Case 1: $\mv(\pi(\scrP_1), \ldots, \pi(\scrP_{n-1})) = 0$. Due to the inductive hypothesis, we may assume after a reordering of the $\scrP_j$ if necessary, that there is $k$, $1 \leq k \leq n-1$, such that $\dim(\sum_{j=1}^k \pi(\scrP_j)) < k$. Since $\scrP_1, \ldots, \scrP_n$ are independent, it follows that $\dim(\sum_{j=1}^k \pi(\scrP_j)) = k-1$ and $\dim(\sum_{j=1}^k \scrP_j) = k$. After a translation of one of the $\scrP_j$ if necessary, we may assume that the affine hull $\aff(\scrP_1 + \ldots + \scrP_k)$ of $\scrP_1 + \cdots + \scrP_k$ passes through the origin. Due to \cref{basis-lemma} we can change the basis of $\zz^n$ to ensure that the subgroup of $\zz^n$ generated by $\zz^n \cap \aff(\scrP_1 + \ldots + \scrP_k)$ is $\zz^k \times \{(0, \ldots, 0)\}$. Let $\pi' : \rr^n \to \rr^{n-k}$ be the projection in the last $(n-k)$-coordinates. We claim that $\pi'(\scrP_{k+1}), \ldots, \pi'(\scrP_n)$ are independent. Indeed, if $\dim( \sum_{j \in J} \pi'(\scrP_j)) < |J|$ for some $J \subseteq \{k+1, \ldots, n\}$, then setting $J' := \{1, \ldots, k\} \cup J$ will yield that $\dim (\sum_{j \in J'} \scrP_j) < |J'|$, contradicting the independence of the $\scrP_j$. Now pick generic $f_1, \ldots, f_n$ such that $\np(f_j) = \scrP_j$ for each $j$. \Cref{bkk-bound-thm,bkk-non-degenerate-thm} and the inductive hypothesis implies that the number of solutions $Z_k$ of $f_1, \ldots, f_k$ on $\nktoruss{k}$ is nonzero, and for each $a = (a_1, \ldots, a_k) \in Z_k$, the number of solutions of $f_{k+1}, \ldots,f_n$ on $\{a \} \times \nktoruss{n-k}$ is nonzero. \Cref{bkk-bound-thm} then implies that $\mv(\scrP_1, \ldots, \scrP_n) > 0$, as required.\\

Case 2: $\mv(\pi(\scrP_1), \ldots, \pi(\scrP_{n-1})) > 0$. In this case, for each generic $\epsilon \in \nktorus$, and for generic $f_j$ with $\np(f_j) = \scrP_j$, $j = 2, \ldots, n-1$,
\begin{align*}
\multntorus{x_n-\epsilon, f_1}{ f_{n-1}}
	&= \multp{f_1|_{x_n = \epsilon}}{f_{n-1}|_{x_n = \epsilon}}{\nktoruss{n-1}}
	= \mv(\pi(\scrP_1), \ldots, \pi(\scrP_{n-1}))
	> 0
\end{align*}
Computing $\multntorus{x_n-\epsilon, f_1}{ f_{n-1}}$ using \cref{order-sum-nktorus-eqn} then implies due to \cref{bkk-bound-thm} that there is a primitive integral $\nu  \in \rnstar$ such that $\langle \nu, e_n \rangle \neq 0$ (where $e_n = (0, \ldots, 0, 1) \in \zz^n$) and $\mv'_\nu(\ld_\nu(\scrP_1), \ldots, \ld_{\nu}(\scrP_{n-1})) > 0$. After a translation if necessary, we may assume the origin is in the relative interior of $\scrP_n$. Then $\max_{\scrP_n}(\nu) > 0$, and therefore assertion \eqref{rational-mv'} of \cref{rational-mv'-prop} implies that $\mv(\scrP_1, \ldots, \scrP_n) > 0$, as required.
\end{proof}

\begin{cor}[{Strict monotonicity of mixed volume, Rojas \cite[Corollary 9]{rojas-convex}}]\label{strictly-mixed-monotone}
If $\scrP'_j$ are convex polytopes in $\rr^n$ such that  $\scrP_j \subseteq \scrP'_j$ for each $j$, then $\mv(\scrP_1, \ldots, \scrP_n) \leq \mv(\scrP'_1, \ldots, \scrP'_n)$. The inequality is strict if and only if both of the following are true:
\begin{enumerate}
\item $\scrP'_1, \ldots, \scrP'_n$ are independent, and
\item \label{strictly-independent} there is $\nu \in \rnstar \setminus \{0\}$ such that the collection $\{\In_\nu(\scrP_j): \scrP_j \cap \In_\nu(\scrP'_j) \neq \emptyset\}$ of polytopes is independent.
\end{enumerate}
\end{cor}

\begin{rem}
Recall that an empty collection of polytopes is independent. Therefore condition \eqref{strictly-independent} of \cref{strictly-mixed-monotone} holds if there is $\nu \in \rnstar \setminus \{0\}$ such that $\scrP_j \cap \In_\nu(\scrP'_j) = \emptyset$ for each $j$.
\end{rem}

\begin{proof}[Proof of \cref{strictly-mixed-monotone}]
Due to \cref{mixed-monotone,positively-mixed} it suffices to prove the following statement:
\begin{align*}
\parbox{0.75\textwidth}{
if $\scrP_j \subset \scrP'_j$, $j = 1, \ldots, n$, and $\mv(\scrP'_1, \ldots, \scrP'_n) > 0$, then $\mv(\scrP'_1, \ldots, \scrP'_n) > \mv(\scrP_1, \ldots, \scrP_n)$ if and only if condition \eqref{strictly-independent} of \cref{strictly-mixed-monotone} is true.
}
\end{align*}
So assume $\mv(\scrP'_1, \ldots, \scrP'_n) > 0$. Due to \cref{poly-rational-appxn,mixed-remark} and the multiadditivity of mixed volumes, we may assume in addition that all the $\scrP_j,\scrP'_j$ are integral. Then choose BKK non-degenerate Laurent polynomials $f_j$ such that $\np(f_j) = \scrP_j$, $j = 1, \ldots, n$. \Cref{bkk-bound-thm,bkk-non-degenerate-thm} imply that
\begin{align*}
 \mv(\scrP_1, \ldots, \scrP_n) = \multfntorusiso \leq  \mv(\scrP'_1, \ldots, \scrP'_n)
\end{align*}
and the inequality is strict if and only if there is $\nu \in \rnstar\setminus \{0\}$ such that
\begin{align} \label{strictly-zero}
\parbox{0.6\textwidth}{
$\In_{\scrP'_j, \nu}(f_1), \ldots, \In_{\scrP'_n, \nu}(f_n)$ have a common zero on $\nktorus$.
}
\end{align}
Since $f_j$ are generic, then it follows from \cref{bkk-bound-thm,positively-mixed} that condition \eqref{strictly-zero} holds if and only if  $\{\In_\nu(\scrP_j): \scrP_j \cap \In_\nu(\scrP'_j) \neq \emptyset\}$ is an independent collection of polytopes.
\end{proof}


\begin{cor}[Bernstein's theorem - alternate version] \label{alternate-corollary}
Let $\scrA_j$ be finite subsets of $\zz^n$ and $f_j$ be Laurent polynomials supported at $\scrA_j$, $j = 1, \ldots, n$. Then
$$\multfntorusiso \leq \mv(\conv(\scrA_1), \ldots, \conv(\scrA_n))$$
If $\mv(\conv(\scrA_1), \ldots, \conv(\scrA_n)) > 0$, then the bound is satisfied with an equality if and only if both of the following conditions hold:
\begin{enumerate}
\item for each nontrivial weighted order $\nu$, the collection $\{\In_\nu(\np(f_j)): \In_\nu(\scrA_j) \cap \supp(f_j) \neq \emptyset\}$ of polytopes is dependent, and
\item  $f_1, \ldots, f_n$ are BKK non-degenerate, i.e.\ they satisfy \eqref{b-non-degeneracy-wt} with $m = n$.
\end{enumerate}
\end{cor}

\begin{proof}
The bound of \cref{alternate-corollary} follows from \Cref{bkk-bound-thm}. \Cref{bkk-non-degenerate-thm} implies that if $\mv(\conv(\scrA_1), \ldots, \conv(\scrA_n)) > 0$, then the bound holds with an equality if and only if $f_1, \ldots, f_n$ are BKK non-degenerate and $\mv(\np(f_1), \ldots, \mv(\np(f_n) = \mv(\conv(\scrA_1), \ldots, \conv(\scrA_n))$. Now the result follows from \cref{strictly-mixed-monotone}.
\end{proof}

\section{Some technical results} \label{bkk-algebraic-applisection}

In this section we compile a few (technical) corollaries of Bernstein's theorem that we use in later chapters.

\begin{prop} \label{positive-dimensional-prop}
Let the set up be as in \cref{bkk-bound-thm}. Assume $\mv(\scrP_1, \ldots, \scrP_n)$ is nonzero. If the set of common zeroes of $f_1, \ldots, f_n$ on $\nktorus$ has a positive dimensional component, then $\multfntorusiso < \mv(\scrP_1, \ldots, \scrP_n)$.
\end{prop}

\begin{proof}
Since it is possible to choose a curve on such a positive dimensional component (\cref{closure-curve-lemma}), the proposition follows from \cref{bkk-non-degenerate-thm,branch-lemma-2}.
\end{proof}

\begin{cor} \label{order-nu-cor}
Let the notation be as in \cref{order-nu-prop}. Assume $V(\Inalphapsinu(f_2), \ldots,\Inalphapsinu(f_n)) \cap \nktorus \neq \emptyset$, so that $C'$ is a curve. Let $\{C'_j\}_j$ be the irreducible components of $C'$, and $\scrB'_{j,\nu}$ be the collection of all branches $B$ of $C'_j$ at infinity (with respect to $\nktorus$) such that $\nu_B$ is proportional to $\nu$. If $V(\Inalphapsinu(f_1), \ldots,\Inalphapsinu(f_n)) \cap \nktoruss{n-1} = \emptyset$, then
\begin{equation}
\begin{aligned}
\sum_{ a \in C' \cap \orbit{\scrQ}} \ord_a(f_1|_{C'})
	& = \sum_j \sum_{(Z,z) \in \scrB'_{j,\nu}} \ord_z(f_1|_{C'_j}) \multp{f_2}{f_n}{C'_j} \\
	& =  \nu(f_1) \multnntorus{\Inalphapsinu(f_2)} {\Inalphapsinu(f_n)}{n-1}
\end{aligned}
\label{order-nu-sum-1}
\end{equation}
where $\multp{f_2}{f_n}{C'_j}$ are defined as in \cref{complete-mult-section}. If in addition $f'_{2, \nu}, \ldots, f'_{n,\nu}$ are BKK non-degenerate, then
\begin{align}
\sum_{ a \in C' \cap \orbit{\scrQ}} \ord_a(f_1|_{C'})
	= \min_{\np(f_1)}(\nu) \mv'_\nu(\In_\nu(\np(f_2)), \ldots, \In_\nu(\np(f_n)))
	\label{order-nu-mv'}
\end{align}
\end{cor}

\begin{proof}
\Cref{order-curve} implies that
\begin{align*}
\sum_{ a \in C' \cap \orbit{\scrQ}} \ord_a(f_1|_{C'})
	&= \sum_j  \sum_{ a \in C'_j \cap \orbit{\scrQ}}   \multp{f_2}{f_n}{C'_j} \ord_a(f_1|_{C'_j}) \\
	&=  \sum_j  \sum_{ a \in C'_j \cap \orbit{\scrQ}}   \multp{f_2}{f_n}{C'_j} \sum_{z \in \pi_j^{-1}(a)} \ord_{z}(\pi_j^*(f_1|_{C'_j}))
\end{align*}
where $\pi_j: \tilde C_j \to C'_j$ are desingularizations of $C'_j$. It follows from the definition of branches in \cref{branchion} that each $z \in  \pi_j^{-1}(a)$, where $a \in C'_j \cap \orbit{\scrQ}$, corresponds to a branch $B = (Z,z)$ at infinity of $C'_j$. Moreover, since $\scrQ$ is a facet of $\scrP$, \cref{curve-to-orbit} implies that a branch $B$ at infinity of $C'_j$ intersects $\orbit{\scrQ}$ if and only if $\nu_B$ is proportional to $\nu$. It follows that
\begin{align*}
\sum_{ a \in C' \cap \orbit{\scrQ}} \ord_a(f_1|_{C'})
	&= \sum_j \sum_{(Z,z) \in \scrB'_{j,\nu}} \ord_z(f_1|_{C'_j}) \multp{f_2}{f_n}{C'_j}
\end{align*}
The result now follows from \cref{order-nu-prop,bkk-bound-thm,bkk-non-degenerate-thm}.
\end{proof}

\section{The problem of characterizing coefficients which guarantee non-degeneracy}
Let $\scrA = (\scrA_1, \ldots, \scrA_n)$ be an $n$-tuple of finite subsets of $\zz^n$, and as in \cref{bkk-non-deg-proof-section}, let $\scrL(\scrA)$ be the collection of $n$-tuples $(f_1, \ldots, f_n)$ of Laurent polynomials such that each $f_j$ is supported at $\scrA_j$. Given $(f_1, \ldots, f_n) \in \scrL(\scrA)$, \cref{bkk-bound-thm} implies that if the coefficients of the $f_j$ are generic, then
\begin{align}
\multfntorusiso = \mv(\scrP_1, \ldots, \scrP_n) \label{multfj=mvPj}
\end{align}
where $\scrP_j = \conv(\scrA_j)$, $j = 1, \ldots, n$. On the other hand, it is straightforward to see that not {\em all} the coefficients of the $f_j$ have to be generic for the equality in \eqref{multfj=mvPj}. E.g.\ if $\scrB_j$ is the set of vertices of $\scrP_j$ and if the coefficient of $x^\alpha$ in each $f_j$ is fixed for each $\alpha \in \scrA_j \setminus \scrB_j$, \eqref{multfj=mvPj} still holds provided the coefficients of $x^\alpha$ in the $f_j$ are generic for all $\alpha \in \scrB_j$. J.\ M.\ Rojas \cite{rojas-toric} posed the problem of identifying all $(\scrB_1, \ldots, \scrB_n)$ which have this property. The precise version of \index{Rojas's problem}Rojas's problem for $\nktorus$ is as follows: let $\scrB_j \subseteq \scrA_j$, $j = 1, \ldots, n$ and $\scrB := (\scrB_1, \ldots, \scrB_n)$. We say that $\scrB$ {\em guarantees $\nktorus$-maximality on $\scrL(\scrA)$} if for all choices from $\kk$ of coefficients of $x^\alpha$ in $f_j$ for all $j$ and all $\alpha \in \scrA_j \setminus \scrB_j$, \eqref{multfj=mvPj} holds provided the coefficients of $x^\beta$ in $f_k$ are generic for all $k$ and all $\beta \in \scrB_k$. Then Rojas's problem\footnote{Rojas \cite{rojas-toric} posed the problem in a more general context (instead of $\nktorus$ he allowed for a broader class of subsets of $\kk^n$) and presented a solution.} is to classify all $\scrB$ which guarantee $\nktorus$-maximality on $\scrL(\scrA)$.

\begin{prop}[Solution of Rojas's problem for $\nktorus$] \label{torus-genericness}
Let $\scrQ_j := \conv(\scrB_j)$, $j = 1, \ldots, n$. Then the following are equivalent:
\begin{enumerate}
\item \label{maximal-guarantee} $\scrB$ guarantees $\nktorus$-maximality on $\scrL(\scrA)$;
\item \label{mixed-equality} $\mv(\scrQ_1, \ldots, \scrQ_n) = \mv(\scrP_1, \ldots, \scrP_n)$;
\item \label{mixed-equality-II} one of the following holds:
\begin{enumerate}
\item $\scrP_1, \ldots, \scrP_n$ are dependent, or
\item for each $\nu \in \rnstar \setminus \{0\}$, the collection $\{\In_\nu(\scrQ_j): \scrQ_j \cap \In_\nu(\scrP_j) \neq \emptyset\}$ of polytopes is dependent.
\end{enumerate}
\end{enumerate}
\end{prop}

\begin{proof}
It is straightforward to check that $\scrB$ guarantees $\nktorus$-maximality on $\scrL(\scrA)$ if and only if for each $(f_1, \ldots, f_n) \in \scrL(\scrA)$, there is $(g_1, \ldots, g_n) \in \scrL(\scrB)$ such that $$\multntorusiso{f_1 + g_1}{f_n + g_n} = \mv(\scrP_1, \ldots, \scrP_n)$$
The implication \eqref{maximal-guarantee} $\im$ \eqref{mixed-equality} follows from taking $(f_1, \ldots, f_n) = (0, \ldots, 0)$ and applying \cref{bkk-bound-thm,mixed-monotone}. For the opposite implication \eqref{mixed-equality} $\im$ \eqref{maximal-guarantee}, assume $\mv(\scrQ_1, \ldots, \scrQ_n) = \mv(\scrP_1, \ldots, \scrP_n)$. Pick $\scrB$-non-degenerate $(g_1, \ldots, g_n) \in \scrL(\scrB)$ such that set $h_j := (1-t)f_j + tg_j$. \Cref{bkk-non-degenerate-thm} implies that $\multntorusiso{g_1}{g_n} \geq \multntorusiso{h_1|_{t= \epsilon}}{h_n|_{t= \epsilon}}$ for each $\epsilon \in \kk$, and therefore \cref{mult-deformation-global} implies that for generic $\epsilon \in \kk$, $\multntorusiso{h_1|_{t= \epsilon}}{h_n|_{t= \epsilon}} = \multntorusiso{g_1}{g_n} = \mv(\scrP_1, \ldots, \scrP_n)$, which implies condition \eqref{maximal-guarantee}. Finally, the equivalence of conditions \eqref{mixed-equality} and \eqref{mixed-equality-II} follows from \cref{strictly-mixed-monotone}.
\end{proof}

We now describe a natural variant of the notion of $\nktorus$-maximality on $\scrL(\scrA)$. Consider the case that $n = 2$ and $\scrA_1 = \scrA_2 = \scrP \cap \zz^2$, where $\scrP$ is the bigger triangle from \cref{fig:generic}. Let $\scrB_1 = \scrB_2 = \scrQ \cap \zz^2$, where $\scrQ$ is the smaller triangle in \cref{fig:generic}. If $f_1, f_2$ are polynomials in two variables with Newton polytope $\scrP$, it is straightforward to check that \eqref{b-non-degeneracy} is satisfied if the coefficients of monomials in $f_j$ whose exponents are in $\scrB$ are generic, and therefore the number of solutions of $f_1, f_2$ on $\nktoruss{2}$ is $\mv(\scrP_1, \scrP_2) = 2 \area(\scrP) = 48$. Note however that $\mv(\scrQ_1, \scrQ_2) = 2 \area(\scrQ) = 12 < \mv(\scrP_1, \scrP_2)$. This motivates \cref{genericness-problem} below. In its statement we use the following notation: we write $\scrL^0(\scrA)$ for the collection of all $(f_1, \ldots, f_n) \in \scrL(\scrA)$ such that $\np(f_j) = \scrA_j$ for each $j$. Given $\scrB = (\scrB_1, \ldots, \scrB_n)$ with $\scrB_j \subseteq \scrA_j$, $j = 1, \ldots, n$, we denote by $\pi_\scrB: \scrL(\scrA) \to \scrL(\scrB)$ be the natural projection which ``forgets'' the coefficients corresponding to $\alpha \not\in \scrB_j$, $j = 1, \ldots, n$, and we write $\scrA \setminus \scrB := (\scrA_1 \setminus \scrB_1, \ldots, \scrA_n \setminus \scrB_n)$. We say that $\scrB$ {\em guarantees $\nktorus$-maximality on $\scrL^0(\scrA)$} if for each $(h_1, \ldots, h_n) \in \scrL(\scrA \setminus \scrB)$, there is a nonempty Zariski open subset $\scrU$ of $\scrL(\scrB)$ such that $\multfntorusiso = \mv(\scrP_1, \ldots, \scrP_n)$ for each $(f_1, \ldots, f_n) \in \scrL^0(\scrA) \cap \pi_{\scrA \setminus \scrB}^{-1}(h_1, \ldots, h_n) \cap  \pi_\scrB^{-1}(\scrU)$.

\begin{figure}[h]
\begin{center}
\def\xmin{-0.5}
\def\xmax{10.5}
\def\ymin{-2.5}
\def\ymax{4.5}
\def\opazero{0.5}
\def\colorzero{green}
\def\colorone{orange}
\def\scalefactor{0.5}

\tikzstyle{dot} = [red, circle, minimum size=5pt, inner sep = 0pt, fill]

\begin{tikzpicture}[scale=\scalefactor]
\draw [gray,  line width=0pt] (\xmin, \ymin) grid (\xmax,\ymax);
\draw [<->] (0, \ymax) |- (\xmax, 0);
\node[dot] (A) at (1,-2) {};
\node[dot] (B) at (9,-2) {};
\node[dot] (C) at (5,4) {};
\fill[\colorzero, opacity=\opazero ] (A.center) --  (B.center) -- (C.center);
\draw[thick] (A) -- (B) -- (C) -- (A);
\draw[thick, fill=\colorone, opacity=\opazero] (5,-2) -- (7,1) -- (3,1) -- cycle;

\node at (5,-2.5) [below] {\picfontsize $\scrP$};
\node at (5,0) {\picfontsize $\scrQ$};
\end{tikzpicture}
\caption{The number of solutions on $\nktoruss{2}$ of polynomials with Newton polytope $\scrP$ is maximal if the coefficients of monomials from $\scrQ$ are generic}  \label{fig:generic}
\end{center}
\end{figure}

\begin{problem} \label{genericness-problem}
Classify all $\scrB = (\scrB_1, \ldots, \scrB_n)$ with $\scrB_j \subseteq \scrA_j$, $j = 1, \ldots, n$, such that $\scrB$ guarantees $\nktorus$-maximality on $\scrL^0(\scrA)$.
\end{problem}

In the remainder of this section we give some partial answers to \cref{genericness-problem}. We start with a condition that guarantees $\nktorus$-maximality on $\scrL^0(\scrA)$. Given a weighted order $\nu$ on the ring of Laurent polynomials in $(x_1, \ldots, x_n)$ and a subset $J$ of $[n]$, let $d_{J, \nu} := \dim(\In_\nu( \sum_{j \in J} \scrP_j))$ and $e_{J, \nu} :=  |\{j \in J:  \In_\nu(\scrP_j) \cap \scrB_j \neq \emptyset\}|$.

\begin{prop} \label{sufficiently-generic}
Assume for each nontrivial weighted order $\nu$, one of the following holds:
\begin{enumerate}
\item \label{vertex-condition} either $\In_\nu(\scrP_j)$ is a vertex of $\scrP_j$ for some $j$,
\item \label{dimension-inequality} or there exists a nonempty subset $J$ of $[n]$ such that $e_{J,\nu}  \geq d_{J,\nu} < |J|$.
\end{enumerate}
Then $\scrB$ guarantees $\nktorus$-maximality on $\scrL^0(\scrA)$.
\end{prop}

\begin{proof}
Fix $(h_1, \ldots, h_n) \in \scrL(\scrA \setminus \scrB)$. It suffices to show that for generic $(g_1, \ldots, g_n) \in \scrL(\scrB)$, if $(f_1, \ldots, f_n) \in \pi_{\scrA \setminus \scrB}^{-1}(h_1, \ldots, h_n) \cap \pi_\scrB^{-1}(g_1, \ldots, g_n)$, then \eqref{b-non-degeneracy} holds provided $(f_1, \ldots, f_n) \in \scrL^0(\scrA)$. Let $\nu$ be a nontrivial weighted order on the ring of Laurent polynomials. If \eqref{vertex-condition} holds for some $j$, then $\In_\nu(f_j)$ is a monomial whenever $\np(f_j) = \scrP_j$ and therefore is nowhere zero on $\nktorus$. So \woutlog\ we may assume \eqref{dimension-inequality} holds for some nonempty subset $J$ of $[n]$. Let $J_\nu := \{j \in J: \In_\nu(\scrP_j) \cap \scrB_j \neq \emptyset\}$, so that $e_{J,\nu} = |J_\nu|$. If $e_{J,\nu} > d_{J,\nu}$, it is straightforward to check that $V(\In_\nu(f_j): j \in J_\nu) \cap \nktorus = \emptyset$ for generic $(g_j: j \in J_\nu)$. On the other hand, if $e_{J,\nu} = d_{J,\nu}$, then there is $j' \in J \setminus J_\nu$ (since $d_{J,\nu} < |J|$), and it is straightforward to check that $V(\In_\nu(f_j): j \in J_\nu) \cap V(f_{j'}) \cap \nktorus = \emptyset$ for generic $(g_j: j \in J_\nu)$. This proves the proposition.
\end{proof}

\begin{prop}[Solution to \cref{genericness-problem} for $n = 2$]
Assume $n = 2$. Then $\scrB = (\scrB_1, \scrB_2)$ guarantees $\nktoruss{2}$-maximality on $\scrL^0(\scrA)$ if and only if for each nontrivial weighted order $\nu$, one of the following holds:
\begin{enumerate}
\item either $\In_\nu(\scrP_j)$ is a vertex of $\scrP_j$ for some $j$,
\item or there is $j$ such that $\In_\nu(\scrP_j) \cap \scrB_j \neq \emptyset$.
\end{enumerate}
\end{prop}

\begin{proof}
The $(\Leftarrow)$ direction follows from \cref{sufficiently-generic}. For the opposite inclusion assume there is a nontrivial weighted order $\nu$ such that for each $j$, $\In_\nu(\scrP_j)$ is an edge of $\scrP_j$ which does not intersect $\scrB_j$. Then it is clear we can pick $f_j$ with $\np(f_j) = \scrP_j$ such that $\In_\nu(f_1) \cap \In_\nu(f_j)$ have a common zero on $\nktoruss{2}$ and for such $(f_1, f_2)$, no choice of coefficients of monomials from $\scrB_j$ would make that zero disappear.
\end{proof}

\begin{prop}[Solution to \cref{genericness-problem} for $n = 3$]
Assume $n = 3$. Then $\scrB$ guarantees $\nktoruss{3}$-maximality on $\scrL^0(\scrA)$ if and only if for each nontrivial weighted order $\nu$, one of the following holds:
\begin{enumerate}
\item \label{vertex-condition'} $\In_\nu(\scrP_j)$ is a vertex of $\scrP_j$ for some $j$,
\item \label{dimension-inequality'-1} or there are $j_1 \neq j_2$ such that $\In_\nu(\scrP_{j_1}+\scrP_{j_2})$ has dimension one and $\In_\nu(\scrP_{j_1}) \cap \scrB_{j_1} \neq \emptyset$,
\item \label{dimension-inequality'-2} or there are $j_1 \neq j_2$ such that $\In_\nu(\scrP_{j_1}+\scrP_{j_2})$ has dimension two and $\In_\nu(\scrP_{j_k}) \cap \scrB_{j_k} \neq \emptyset$ for each $k = 1, 2$,
\item \label{subdivide-condition} or there are $j_1 \neq j_2$ such that
\begin{enumerate}
\item \label{subdivide-condition-0} there is no positive dimensional polytope which is a {\em Minkowski summand}\footnote{We say that a convex polytope $\scrP$ is a {\em Minkowski summand} of a convex polytope $\scrQ$ if there is a convex polytope $\scrR$ such that $\scrQ = \scrP + \scrR$.} of $\In_\nu(\scrP_{j_k})$ for each $k = 1, 2$, 
\item $\In_\nu(\scrP_{j_k}) \cap \scrB_{j_k} = \emptyset$ for each $k = 1, 2$, and
\item $\In_\nu(\scrP_{j_3}) \cap \scrB_{j_3} \neq \emptyset$, where $j_3$ is the single element of $\{1, 2, 3\} \setminus \{j_1, j_2\}$.
\end{enumerate}
\end{enumerate}
\end{prop}

\begin{proof}
At first we prove the $(\Leftarrow)$ implication. Fix $(h_1, \ldots, h_n) \in \scrL(\scrA \setminus \scrB)$, $(g_1, \ldots, g_n) \in \scrL(\scrB)$ and $(f_1, \ldots, f_n) \in \pi_{\scrA \setminus \scrB}^{-1}(h_1, \ldots, h_n) \cap \pi_\scrB^{-1}(g_1, \ldots, g_n) \cap \scrL^0(\scrA)$. Pick a nontrivial weighted order $\nu$. If one of \eqref{vertex-condition'}, \eqref{dimension-inequality'-1} or \eqref{dimension-inequality'-2} holds, then \cref{sufficiently-generic} implies that $V(\In_\nu(f_1), \ldots, \In_\nu(f_n)) \cap \nktorus$ is empty if $g_1, \ldots, g_n$ are generic. So assume \eqref{subdivide-condition} holds. Condition \eqref{subdivide-condition-0} implies that $\In_\nu(f_{j_1})$ and $\In_\nu(f_{j_2})$ have no common non-invertible factor in $\kk[x_1, x_1^{-1}, \ldots, x_n, x_n^{-1}]$. This implies that $V(\In_\nu(f_{j_1}), \In_\nu(f_{j_2}))$ is either empty or have codimension $2$ in $\nktoruss{3}$. Therefore, by choosing generic $g_{j_3}$, it is possible to ensure that $V(\In_\nu(f_{j_1}), \In_\nu(f_{j_2})) \cap V(\In_\nu(g_{j_3})) \cap \nktoruss{3} = \emptyset$, which completes the proof of $(\Leftarrow)$ implication. For the opposite implication, assume there is a nontrivial weighted order $\nu$ such that $\In_\nu(\scrP_j)$ is positive dimensional for each $j$, and one of the following holds:
\begin{prooflist}
\item \label{empty-Bintersection} either $\In_\nu(\scrP_j) \cap \scrB_j = \emptyset$ for each $j$,
\item \label{summand-exists} or there are $j_1 \neq j_2$ such that
\begin{prooflist}
\item there is a positive dimensional polytope $\scrQ$ which is a Minkowski summand of each $\In_\nu(\scrP_{j_k})$, $k = 1, 2$,
\item $\In_\nu(\scrP_{j_k}) \cap \scrB_{j_k} = \emptyset$ for each $k = 1, 2$, and
\item $\In_\nu(\scrP_{j_3}) \cap \scrB_{j_3} \neq \emptyset$, where $j_3$ is the single element of $\{1, 2, 3\} \setminus \{j_1, j_2\}$.
\end{prooflist}
\end{prooflist}
If \ref{empty-Bintersection} holds then one can choose $(f_1, f_2, f_3) \in \scrL^0(\scrA)$ such that $\In_\nu(f_j)$ have a common zero on $\nktorus$, and the zero would be unaffected by the coefficients of monomials from the $\scrB_j$, so that $\scrB$ will not be able to ensure that \eqref{b-non-degeneracy} holds. On the other hand, if \ref{summand-exists} holds, then we can choose $f_{j_1}, f_{j_2}$ such that $\In_\nu(f_{j_1})$ and $\In_\nu(f_{j_2})$ have a common factor $g$ with Newton polytope $\scrQ$. Then for generic $g_{j_3}$ and generic $f_{j_3} \in \pi_{\scrB_3}^{-1}(g_{j_3})$, Bernstein's theorem would imply that $V(\In_\nu(f_{j_1}), \In_\nu(f_{j_2}), \In_\nu(f_{j_3})) \cap \nktoruss{3} \supseteq V(g, \In_\nu(f_{j_3})) \cap \nktoruss{3} \neq \emptyset$, so that $(f_1, f_2, f_3)$ violates \eqref{b-non-degeneracy}, as required.
\end{proof}

\section{Notes}
All major results of this chapter are well known. A.\ Khovanskii in \cite[Section 27]{burago-zalgaller} gave a simple proof of Bernstein's formula in zero characteristic. The main distinction between his proof and ours is in the handling of intersection multiplicity: he counts the number of roots of systems which are {\em nonsingular} at the points of intersection, and then argues that every system can be deformed into such systems. We avoid this approach since in positive characteristics it would involve having as a technical overhead some versions of Bertini-type theorem every time we deal with intersection multiplicity. B.\ Huber and B.\ Sturmfels \cite{hurmfels-polyhedra} gave a constructive proof of Bernstein's theorem (in zero characteristic) which has had a deep impact on the ``homotopy continuation'' method to numerically compute the solutions of polynomial systems; the techniques of their proof also give an efficient way to compute mixed volumes. Our proof that the set of BKK non-degenerate polynomials is Zariski open follows the arguments from \cite{oka-newton-topology}, and our proof that it is nonempty comes from \cite{toricstein}. 


\part{Beyond the torus}
\label{postoric-part}
\chapter{Number of zeroes on the affine space I: (Weighted) B\'ezout theorems} \label{bezout-chapter}
\chaptermark{(Weighted) B\'ezout theorems}
\newcommand{\overmat}[2]{%
	\makebox[0pt][l]{
		$\smash{
			\color{white}
			\overbrace{
				\phantom{%
					\begin{matrix}#2\end{matrix}
				}
			}^{\text{\color{black}#1}}
		}$
	}#2
}
\newcommand\bovermat[2]{%
	\makebox[0pt][l]{
		$\smash{
			\overbrace{
				\phantom{%
					\begin{matrix}#2\end{matrix}
				}
			}^{\text{#1}}
		}$
	}#2
}

In this chapter we use the results from \cref{bkk-chapter} to prove B\'ezout's theorem (\cref{bezout-1}) and two of its classical generalizations: the weighted homogeneous version (\cref{wt-bezout-1}) and the weighted multi-homogeneous version (\cref{multi-bezout}). The weighted degrees considered in these results have the property that the weight of each variable is {\em positive}. In \cref{weighted-application,weighted-multapplication} we establish more general versions of these results involving arbitrary weighted degrees as special cases of the extension of Bernstein's theorem to $\kk^n$. We continue to assume that $\kk$ is an algebraically closed field.

\section{Weighted degree}
Let $\omega$ be an integral element of $\rnstar$. The \index{Weighted!degree}{\em weighted degree} corresponding to $\omega$, which by an abuse of notation we also denote by $\omega$, is the map $\kk[x_1, \ldots, x_n] \to \zz \cup \{-\infty\}$ given by
\begin{align*}
f \mapsto \max_{\supp(f)}(\omega)
\end{align*}
Given $f = \sum_\alpha c_\alpha x^\alpha \in \kk[x_1, \ldots, x_n]$, the \index{Leading form}{\em leading form} of $f$ with respect to $\omega$ is $\ld_\omega(f) := \sum_{ \langle \omega, \alpha \rangle = \omega(f) }c_\alpha x^{\alpha}$. We say that $f$ is \index{Weighted!homogeneous!polynomial}{\em weighted homogeneous} with respect to $\omega$ (or in short, {\em $\omega$-homogeneous}) if $\omega (x^\alpha)$ are equal for all $\alpha$ such that $c_\alpha \neq 0$, or equivalently, if $\supp(\ld_\omega(f)) = \supp(f)$.

\section{$\wpn$ as a compactification of $\kk^n$ when the $\omega_j$ are positive and $\omega_0 = 1$}
Assume $\omega$ is an integral element of $\rnnstar{n+1}$ with coordinates $(\omega_0, \ldots, \omega_n)$ with respect to the basis dual to the standard basis of $\rr^{n+1}$ such that each $\omega_j$ is {\em positive}. The {\em weighted projective space} $\wpn$ corresponding to $\omega$ was constructed in \cref{weighted-toric-section}. In this section we treat the case that $\omega_0 = 1$. For each $j = 0, \ldots, n$, let $U_j := \wpn \setminus V(x_j)$, be the ``coordinate chart'' of $\wpn$ considered in \cref{wt-identification-section}. Since $\omega_0 = 1$, \cref{toric-k^n} implies that the map $(a_1, \ldots, a_n) \mapsto [1:a_1: \cdots : a_n]$ induces an isomorphism between $\kk^n$ and $U_0$, and therefore $\wpn$ is a compactification of $\kk^n$. The set of points at infinity (with respect to $U_0 \cong \kk^n$) on $\wpn$ is $V(x_0) = \{[0: a_1: \cdots :a_n] : (a_1, \ldots, a_n) \in \kk^n \setminus \{0\}\}
\cong \wwpnn{\omega_1, \ldots, \omega_n}{n-1}$. Since $ \omega_0 = 1$, for each polynomial $f \in \kk[x_1, \ldots, x_n]$, we can define its \index{Homogenization!weighted}\index{Weighted!homogenization}{\em weighted homogenization} with respect to $\omega$ as $\tilde f := x_0^{\omega(f)}f(x_1/x_0^{\omega_1}, \ldots, x_n/x_0^{\omega_n})$. It is straightforward to check that $\tilde f$ is $\omega$-homogeneous and $\ld_\omega(\tilde f) = \ld_\omega(f)$.

\begin{prop} \label{wpn-closure}
Identify $\kk^n$ with $U_0$.
\begin{enumerate}
\item \label{wpn-closure-Z} Let $f \in \kk[x_1, \ldots, x_n]$. Then $V(\tilde f)$ is the Zariski closure in $\wpn$ of $V(f) \subset \kk^n$.
\item \label{wpn-closure-infty} If $f_1, \ldots, f_k \in \kk[x_1, \ldots, x_n]$, then the following are equivalent:
\begin{enumerate}
\item $\bigcap_{j=1}^k V(\tilde f_j) \setminus \kk^n = \emptyset$.
\item there is no common zero of $\ld_\omega(f_1), \ldots, \ld_\omega(f_k)$ on $\kk^n \setminus \{\origin\}$.
\end{enumerate}
\end{enumerate}
\end{prop}

\begin{proof}
At first we prove assertion \eqref{wpn-closure-Z}. If $f$ is a nonzero constant, then both $V(f)$ and $V(\tilde f)$ are empty. So assume $f$ is a non-constant polynomial. If $n = 1$, then $V(f)$ consists of finitely many points and it is straightforward to check that $V(f) = V(\tilde f) \subset \wwpnn{\omega_0, \omega_1}{1}$. So assume $n \geq 2$. For each $j = 0, \ldots, n$, the set of points in $V(\tilde f) \cap U_j$ is precisely the set of zeroes on $U_j$ of the regular function $(\tilde f)^{\omega_j}/x_j^{\omega(f)}$. Therefore \cref{thm:pure-dimension} implies that each irreducible component of $V(\tilde f)$ has dimension $n-1$. On the other hand,
\begin{align}
V(\tilde f) \cap (\pp^n \setminus \kk^n) = V(\tilde f) \cap V(x_0) = V(\ld_\omega(f)) \cap V(x_0) \label{wpn-infty-eqn}
\end{align}
Since $\ld_\omega(f)$ is a nonzero polynomial and $V(x_0) \cong \wwpnn{\omega_1, \ldots, \omega_n}{n-1}$, \cref{thm:pure-dimension} implies that each irreducible component of $V(\tilde f) \setminus \kk^n$ has dimension $n-2$, and therefore can not be an irreducible component of $V(\tilde f)$. Consequently every irreducible component of $V(\tilde f)$ intersects $\kk^n$ and therefore contains an irreducible component of $V(f)$. Since $V(\tilde f)$ is Zariski closed in $\wpn$, this completes the proof of assertion \eqref{wpn-closure-Z}. Assertion \eqref{wpn-closure-infty} follows from identity \eqref{wpn-infty-eqn}.
\end{proof}


\section{Weighted B\'ezout theorem} \label{wt-bezout-proof-section}
We use the following notation throughout the rest of the book: for a nonnegative integer $n$ we denote by $[n]$ the set $\{1, \ldots, n\}$; if $I \subseteq [n]$ and $k$ is a field (in most cases $k$ will be either $\kk$ or $\rr$), we write
\begin{align}
\ki & := \{(x_1, \ldots, x_n) \in k^n: x_i =  0\ \text{if}\ i \not\in I\}  \cong \kii{|I|}, \label{ki}\\
\kstari &:=  \{\prod_{i\in I} x_i \neq 0\} \cap \ki \cong \kstarii{|I|} \label{kstari}
\end{align}
where $k^* := k \setminus \{0\}$. Note that $\kii{\emptyset} = \kstarii{\emptyset} = \{0\}$.

\begin{thm}[Weighted B\'ezout theorem (\cref{wt-bezout})] \label{wt-bezout-1}
\index{Weighted!B\'ezout theorem!version I}
\index{B\'ezout's theorem!weighted!version I}
Let $\omega$ be a weighted degree on $\kk[x_1, \ldots, x_n]$ with positive weights $\omega_j$ for $x_j$, $j = 1, \ldots, n$. Then the number of isolated solutions of polynomials $f_1, \ldots, f_n$ on $\kk^n$ is bounded above by $(\prod_j \omega(f_j))/(\prod_j \omega_j)$. This bound is exact if and only if the leading weighted homogeneous forms of $f_1, \ldots, f_n$ have no common solution other than $(0, \ldots, 0)$.
\end{thm}

\begin{proof}
Due to multiadditivity of intersection multiplicities (assertion \eqref{mult-additive} of \cref{int-mult-curve}), we may replace each $f_i$ by $f_i^m$ for an appropriate positive integer $m$ and assume that $m_{i,j} := \omega(f_i)/\omega_j$ is an integer for each $i,j$. Let $\scrP_i$ be the simplex in $\rr^n$ with vertices at the origin and at $m_{i,j}e_j$, $j = 1, \ldots, n$, where $e_1, \ldots, e_n$ are the standard unit vectors in $\rr^n$. Analogous to the definition of $\multPntorusiso$ in \eqref{multPntorusiso}, define
\begin{align*}
\multPkniso := \max\{\multkniso{p_1}{p_n}: \supp(p_j) \subseteq \scrP_j,\ j =1 , \ldots, n\}
\end{align*}
It suffices to show that
\begin{prooflist}
\item \label{weighted-multpknsio} $\multPkniso = (\prod_j \omega(f_j))/(\prod_j \omega_j)$, and
\item \label{weighted-non-degeneracy} $\multfkniso = \multPkniso$ if and only if there is no common zero of $\ld_\omega(f_1), \ldots, \ld_\omega(f_n)$ on $\kk^n \setminus \{0\}$.
\end{prooflist}
Assertion \ref{weighted-multpknsio} follows from Bernstein's theorem (\cref{bkk-bound-thm}) and the following claim.

\begin{proclaim} \label{wt-bezout-claim}
\begin{enumerate}
\item $\mv(\scrP_1, \ldots, \scrP_n) = (\prod_j \omega(f_j))/(\prod_j \omega_j)$.
\item \label{wt-equal-mult} $\multPkniso = \multPntorusiso$.
\end{enumerate}
\end{proclaim}

\begin{proof}
Since $\scrP_j = \frac{\omega(f_j)}{\omega(f_1)}\scrP_1$ for each $j = 1, \ldots, n$, the properties of mixed volume from \cref{mixed-unique} imply that
\begin{align*}
\mv(\scrP_1, \ldots, \scrP_n)
	&= \mv(\scrP_1, \frac{\omega(f_2)}{\omega(f_1)}\scrP_1, \ldots, \frac{\omega(f_n)}{\omega(f_1)}\scrP_1)
	= \prod_{j=2}^n \frac{\omega(f_j)}{\omega(f_1)} \mv(\scrP_1, \ldots, \scrP_1)\\
	&= \prod_{j=2}^n \frac{\omega(f_j)}{\omega(f_1)} n!\vol_n(\scrP_1)
	=  \prod_{j=2}^n \frac{\omega(f_j)}{\omega(f_1)}   \prod_{j=1}^n m_{1,j}
	=  \prod_{j=1}^n \frac{\omega(f_j)}{\omega_j}
\end{align*}
This proves the first assertion of the claim. %
For the second assertion, let $p_j$ be an arbitrary polynomial supported at $\scrP_j$, $j = 1, \ldots, n$. Since it is clear that $\multPkniso \geq \multPntorusiso$, it suffices to show that there are $p'_j$ supported at $\scrP_j$ such that $\multkniso{p_1}{p_n} \leq \multntorusiso{p'_1}{p'_n}$. Write $\scrA_j := \scrP_j \cap \zz^n$, $j = 1, \ldots, n$. Pick $(\scrA_1, \ldots, \scrA_n)$-non-degenerate $(q_1, \ldots, q_n) \in \scrL(\scrA_1, \ldots, \scrA_n)$, and set $r_i(x,t) := (1-t)p_i + tq_i$. For each $\epsilon \in \kk$, write $r_{\epsilon,i} := r_i|_{t = \epsilon}$. Note that $r_{1,i} = q_i$ and $r_{0,i} = p_i$. Pick a generic $\epsilon \in \kk$. Due to \cref{bkk-0} we may assume that $(r_{\epsilon, 1}, \ldots, r_{\epsilon,n})$ is also $(\scrA_1, \ldots, \scrA_n)$-non-degenerate. Since the intersection of an $\scrA_j$ with a coordinate subspace is a (nonempty) face of $\scrA_j$, it is then straightforward to check (e.g.\ using \cref{finite-remark}) from the definition of $(\scrA_1, \ldots, \scrA_n)$-non-degeneracy that $r_{\epsilon, 1}, \ldots, r_{\epsilon,n}$ do not have any common zero on $\kk^n \setminus \nktorus$. Assertion \eqref{openly-max} of \cref{mult-deformation-global} then implies that $\multkniso{p_1}{p_n} \leq \multkniso{r_{\epsilon, 1}}{r_{\epsilon, n}} = \multntorusiso{r_{\epsilon, 1}}{r_{\epsilon, n}}$, as required.
\end{proof}

Now we prove assertion \ref{weighted-non-degeneracy}. At first assume $\multfkniso < \multPkniso$. We will show that the leading weighted homogeneous forms $\ld_\omega(f_i)$ of $f_i$ have a common zero on $\kk^n \setminus \{0\}$. Let $\omega' := (1, \omega_1, \ldots, \omega_n)$. Embed $\kk^n$ into $\wpnprime$ via the map $(x_1, \ldots, x_n) \mapsto [1: x_1: \cdots : x_n]$. Set $p_i := f_i$, $i = 1, \ldots, n$. Let $q_i$ and $r_i := (1-t)p_i + tq_i$ be as in the proof of \cref{wt-bezout-claim}. Let $V$ be the (finite) set of common zeroes of $q_1, \ldots, q_n$ on $\kk^n$ and $C \subseteq \kk^{n+1}$ be the union of irreducible components of $V(r_1, \ldots, r_n)$ which intersect $V \times \{1\}$. Assertion \eqref{onely-isolated} of \cref{mult-deformation-global} implies that $C$ is a curve. Since $\multfkniso < \multkniso{q_1}{q_n}$, assertion \eqref{non-max-condition} of \cref{mult-deformation-global} implies that one of the following holds:
\begin{prooflist}[resume]
\item \label{possibility-non-isolated} there is a positive dimensional component of $V(f_1, \ldots, f_n) \subset \kk^n$, or
\item \label{possibility-at-infinity} $C$ ``has a point at infinity at $t = 0$'', i.e.\ if $\bar C$ is the closure of $C$ in $\wpnprime \times \pp^1$, then $\bar C \cap ((\wpnprime \setminus X) \times \{0\}) \neq \emptyset$.
\end{prooflist}
Denote the weighted homogeneous coordinates on $\wpnprime$ by $[x_0:  \cdots : x_n]$. Let $\tilde f_i$ and $\tilde q_i$ be the {\em weighted homogenization} with respect to $\omega'$ respectively of $f_i$ and $q_i$. \Cref{wpn-closure} implies that the closures of $V(f_i)$ and $V(q_i)$ in $\wpnprime$ are respectively $V(\tilde f_i)$ and $V(\tilde q_i)$. Since $\omega(f_i) = \omega(q_i)$ for each $i$, it follows that the closure of $V(r_i)$ in $\wpnprime \times \kk$ is $V(\tilde r_i)$, where $\tilde r_i := (1-t)\tilde f_i + t \tilde q_i$. If \ref{possibility-non-isolated} holds, then the closure of $V(f_1, \ldots, f_n)$ in $\wpnprime$ contains a point $a' \in \wpnprime \setminus \kk^n$. The weighted homogeneous coordinates of $a'$ are of the form $[0: a_1: \cdots : a_n]$ with $a = (a_1, \ldots, a_n) \in \kk^n \setminus \{0\}$. Then $\tilde f_i(a')  = 0$, and therefore $\ld_\omega(f_i)(a) = 0$ for each $i$. On the other hand, if \ref{possibility-at-infinity} holds, then let $(a',0) \in \bar C \cap ((\wpnprime \setminus X) \times \{0\})$. Since $\tilde r_i(a',0) = 0$ for each $i$, this again yields a common zero of the $\ld_\omega(f_i)$ on $\kk^n\setminus \{0\}$, as required. \\

It remains to show the necessity of the non-degeneracy condition. Assume the leading weighted homogeneous forms $\ld_\omega(f_i)$ of $f_i$ have a common solution $a \in \kk^n \setminus \{\origin\}$. As in \cref{bkk-necessary-section}, pick BKK non-degenerate $g_1, \ldots, g_n$ with a common zero $b \in \kk^n$ such that for each $i$, $\np(g_i) = \scrP_i$ and $\ld_\omega(g_i)(a) \neq 0 \neq f_i(b)$. Define a rational curve $C$ on $\kk^n$ via the parametrization $c(t):= (c_1(t), \ldots, c_n(t)) :\kk \to \kk^n$ from \eqref{trick-curve} with $\nu := -\omega$, i.e.
\begin{align*}
c_j(t) := a_jt^{-\omega_j} + (b_j - a_j)t^{-\omega_j + k_j},\ j = 1, \ldots, n.
\end{align*}
where each $k_j$ is a positive integer. Define $h_1, \ldots, h_n$ as in \eqref{trick-deformation} with $m_j = -\omega(f_j) = -\omega(g_j)$, i.e.\
\begin{align*}
h_j &:=
	t^{\omega(f_j)}f_j(c(t))g_j - t^{\omega(g_j)}g_j(c(t))f_j
\end{align*}
Note that $t^{\omega(f_j)}f_j(c(t))$ and $t^{\omega(g_j)}g_j(c(t))$ are {\em polynomials} in $t$. The same arguments as in the alternate proof of necessity of BKK non-degeneracy in \cref{modified-proof} show that $\multfkniso < \multgkniso$, so that the weighted homogeneous bound is {\em not} exact.
\end{proof}

\begin{cor}[B\'ezout's theorem (\cref{bezout})] \label{bezout-1}
\index{B\'ezout's theorem}
The number of isolated solutions of polynomials $f_1, \ldots, f_n$ on $\kk^n$ is bounded above by $\prod_j \deg(f_j)$. This bound is exact if and only if the leading homogeneous forms of $f_1, \ldots, f_n$ have no common solution other than $(0, \ldots, 0)$.
\end{cor}

\begin{proof}
In \cref{wt-bezout-1} take $\omega$ to be the usual degree of polynomials.
\end{proof}

Combining the arguments of the proof of assertion \eqref{wt-equal-mult} of \cref{wt-bezout-claim} with \cref{thm:pure-dimension,bkk-bound-thm,positively-mixed} gives the following characterization of polytopes $\scrP_1, \ldots, \scrP_n$ such that $\multPkniso = \multPntorusiso$. Given a coordinate subspace $H$ of $\rr^n$, we write $\thp :=  \{j: \scrP_j \cap H \neq \emptyset\} \subset [n]$.

\begin{lemma} \label{equal-multpn}
Let $\scrP_1, \ldots, \scrP_n$ be convex integral polytopes in $\rzeroo{n}$. Then the following are equivalent:
\begin{enumerate}
\item $\multPkniso = \multPntorusiso$.
\item For each proper coordinate subspace $H$ of $\rr^n$, one of the following is true:
\begin{enumerate}
\item either there is a coordinate subspace $H'$ of $\rr^n$ such that $H' \supset H$ and $|\thprimep| < \dim(H')$, or
\item $\thp$ is nonempty and $\{\scrP_j \cap H: j \in \thp\}$ is a {\em dependent} collection of polytopes. \qed
\end{enumerate}
\end{enumerate}
\end{lemma}

\section{Products of weighted projective spaces}
In this section we examine the closures of the hypersurfaces in a product of weighted projective spaces. This would be useful in the proof of the weighted multi-homogeneous B\'ezout bound (\cref{multi-bezout}). Let $\bar X := \prod_{j=1}^s \wwpnn{\omega_j}{n_j}$, where each $\omega_j$ is a weighted degree on $A_j := \kk[x_{j,0}, \ldots, x_{j,n_j}]$ such that the weight $\omega_{j,k}$ of $x_{j,k}$ is positive for each $j,k$. Let $A := \kk[x_{j,k}: 1 \leq j \leq s,\ 0 \leq k \leq n_j]$, and for each $j$, let $\tilde \omega_j$ be the trivial extension of $\omega_j$ to $A$, i.e.\
\begin{align*}
\tilde \omega_j(x_{k,l})
	&=
	\begin{cases}
	\omega_{k,l} & \text{if}\ k = j,\\
	0 & \text{otherwise.}
	\end{cases}
\end{align*}
We say that a polynomial $h\in A$ is \index{Weighted!multi-homogeneous!polynomial}{\em (weighted multi-) homogeneous with respect to $\Omega := \{\omega_1, \ldots, \omega_s\}$}, or in short, $f$ is {\em $\Omega$-homogeneous}, if it is $\tilde \omega_j$-homogeneous for each $j = 1, \ldots, s$. If $h$ is $\Omega$-homogeneous, then $V(h) := \{a \in \bar X: h(a) = 0\}$ is a well defined Zariski closed subset of $\bar X$.

\subsection{The case that \bm{$\omega_{j,0} = 1$} for each \bm{$j = 1, \ldots, s$}} \label{multi-compactisection}
In this section we consider the case that $\omega_{j,0} = 1$ for each $j$. In this case each $\wwpnn{\omega_j}{n_j}$ is a compactification of $\kk^{n_j}$, and therefore $\bar X$ is a compactification of $\kk^n$, where $n := \sum_{j=1}^s n_j$. More precisely, $\kk^n$ can be identified with $\bar X \setminus V( \prod_{j=1}^s x_{j,0})$, and we may treat $B := \kk[x_{j,k}: 1 \leq j \leq s,\ 1 \leq k \leq n_j] \subset A$ as the coordinate ring of $\kk^n$. Let $f \in B$. Then {\em $\Omega$-homogenization} $\tilde f$ of $f$ is formed by substituting $x_{j,k}/x_{j,0}$ for $x_{j,k}$ in $f$ for each $j,k$, and then clearing out the denominator. The following result is the analogue of assertion \eqref{wpn-closure-Z} of \cref{wpn-closure}, and follows from the same arguments.

\begin{prop} \label{multi-closure-Z}
Let $f \in B$. If $\kk$ is algebraically closed, then $V(\tilde f)$ is the Zariski closure in $\bar X$ of $V(f)\subset \kk^n$. \qed
\end{prop}

The set of points at infinity on $\bar X$ is $\bar X \setminus \kk^n = V(\prod_{j=1}^s x_{j,0}) = \bigcup_{J} Y_J$, where the union is over all {\em nonempty} subsets $J$ of $[s] := \{1, \ldots, s\}$, and $Y_J$ are defined as follows:
\begin{align*}
Y_J
	&=  V(x_{j,0} : j \in J) \setminus V(\prod_{j \not\in J} x_{j,0})
	\cong \prod_{j \in J} \wwpnn{\omega'_j}{n_j-1} \times \prod_{j \not\in J} \kk^{n_j}
\end{align*}
where $\omega'_j$ are the restriction of $\omega_j$ to $\kk[x_{j,1}, \ldots, x_{j,n_j}]$ for each $j$. Fix a nonempty subset $J$ of $[s]$.  Given $f \in B$, we would like to compute the points at infinity on $V(\tilde f)$ which belong to $Y_J$. If $f = \sum_\alpha c_\alpha x^\alpha$, we write $\ld_{\Omega, J}(f)$ for be sum of $c_\alpha x^\alpha$ over all $\alpha$ such that $\tilde \omega_j(x^\alpha) = \tilde \omega_j(f)$ for each $j \in J$. In other words, $\ld_{\Omega, J}(f)$ is obtained from $\tilde f$ by substituting $x_{j,0} = 0$ for each $j \in J$ and $x_{j,0} = 1$ for each $j \not\in J$.

\begin{example}
Let $s = 3$, $n_1 = n_2 = 1$ and $n_3 = 2$ (so that $n = 4$). Let $\omega_1, \omega_2$ be the usual degree in respectively $x_{1,1}$ and $x_{2,1}$ coordinates, and $\omega_3$ be the weighted degree in $(x_{3,1}, x_{3,2})$ coordinates corresponding to weights $2$ for $x_{3,1}$ and $3$ for $x_{3,2}$. Let $f = x_{1,1}^5 +  x_{2,1}^7 + x_{1,1}^5x_{2,1}^7 + x_{3,1}^3 +  x_{1,1}^5x_{3,2}^2 $. Then $\omega_1(f) = 5$, $\omega_2(f) = 7$, $\omega_3(f) = 6$, and
\begin{align*}
\ld_{\Omega, J}(f)
	&=
	\begin{cases}
		x_{1,1}^5 + x_{1,1}^5x_{2,1}^7 +  x_{1,1}^5x_{3,2}^2 & \text{if}\ J = \{1\},\\
		x_{2,1}^7  + x_{1,1}^5x_{2,1}^7 & \text{if}\ J = \{2\},\\
		x_{3,1}^3 +  x_{1,1}^5x_{3,2}^2 & \text{if}\ J = \{3\},\\
		x_{1,1}^5x_{2,1}^7 	 & \text{if}\ J = \{1,2\},\\
		x_{1,1}^5x_{3,2}^2 	& \text{if}\ J = \{1,3\},\\
		0 & \text{if $J = \{2, 3\}$ or $J = \{1,2,3\}$.}\\
	\end{cases}
\end{align*}
Note the following difference from the weighted homogeneous case: if $|J| \geq 2$, then it might happen that $\ld_{\Omega, J}(f) = 0$ even if $f$ is a nonzero polynomial.
\end{example}

Since $\ld_{\Omega, J}(f)$ is $\tilde \omega_j$-homogeneous for each $j \in J$, it defines a well defined Zariski closed subset of $Y_J$. It is straightforward to check that this set is precisely the intersection of $V(\tilde f)$ and $Y_J$, which is the content of the next result.

\begin{prop} \label{multi-closure-infty}
Assume $\kk$ is algebraically closed. If $f \in B$, then $V(\tilde f) \cap Y_J = V(\ld_{\Omega, J}(f)) \cap Y_J$. If $f_1, \ldots, f_k \in B$, then the following are equivalent
\begin{enumerate}
\item $\bigcap_{j=1}^k V(\tilde f_j) \setminus \kk^n = \emptyset$.
\item For each nonempty subset $J$ of $[s]$, there is no common zero of $\ld_{\Omega, J}(f_1), \ldots, \ld_{\Omega, J}(f_k)$ on $\prod_{j \in J} (\kk^{n_j} \setminus \{\origin\}) \times \prod_{j \not\in J} \kk^{n_j}$. \qed
\end{enumerate}
\end{prop}

\section{Weighted multi-homogeneous B\'ezout theorem}

We now generalize the weighted B\'ezout bound to the multi-projective setting. Let $\mscrI := (I_1, \ldots, I_s)$ be an ordered partition of $[n] := \{1, \ldots, n\}$, i.e.\ $[n] = \bigcup_j I_j$ and $\sum_j |I_j| = n$. For each $j = 1, \ldots, s$, let $\omega_j$ be a weighted degree on $\kk[x_k: k \in I_j]$ with positive weights $\omega_{j,k}$ for $x_k$, $k  \in I_j$. Let $f_1, \ldots, f_n$ be polynomials on $\kk^n$. Given $d_{i,j} := \omega_j(f_i)$, we would like to compute a (sharp) upper bound of the number of isolated solutions of $f_1, \ldots, f_n$. Let $n_j := |I_j|$ and $l_j$ be the least common multiple of $\omega_{j,1}, \ldots, \omega_{j,n_j}$ and $\scrS_j$ be the simplex in $\rzeroo{n_j}$ defined by $\{ \alpha: \langle \omega_j , \alpha \rangle \leq l_j\}$. Note that $\supp(f_i) \subseteq \scrP_i := \prod_{j=1}^s (d_{i,j}/l_j) \scrS_j \subset \rr^n$. By definition $\mv(\scrP_1, \ldots, \scrP_n)$ is the coefficient of $\lambda_1 \cdots \lambda_n$ in the polynomial
\begin{align*}
\vol_n(\sum_{i=1}^n \lambda_i \scrP_i)
	&= \vol_n(\sum_{i=1}^n \lambda_i \prod_{j=1}^s \frac{d_{i,j}}{l_j}  \scrS_j)
	= \vol_n(\prod_{j=1}^s (\sum_{i=1}^n \lambda_i \frac{d_{i,j}}{l_j}  \scrS_j))
	= \prod_{j=1}^s (\sum_{i=1}^n \lambda_i \frac{d_{i,j}}{l_j} )^{n_j} \vol_n(\prod_{j=1}^s \scrS_j) \\
	&= \frac{\prod_{j=1}^s\prod_{k=1}^{n_j}(l_j/\omega_{j,k})} {\prod_{j=1}^s (n_j! l_j^{n_j})}
	\prod_{j=1}^s (\sum_{i=1}^n \lambda_i d_{i,j})^{n_j}
	= \frac{ \prod_{j=1}^s (\sum_{i=1}^n \lambda_i d_{i,j})^{n_j}}
	 {\prod_{j=1}^s (n_j! \prod_{k=1}^{n_j}\omega_{j,k})}
\end{align*}
The coefficient of $\lambda_1 \cdots \lambda_n$ in $\prod_{j=1}^s (\sum_{i=1}^n \lambda_i d_{i,j})^{n_j}$ is the {\em permanent} of the following $n \times n$ matrix:
\vspace{2mm}

\begin{align*}
\begin{matrix}
D(\mscrI, \vec d)
	&:=
	  \begin{pmatrix}
	  \bovermat{$n_1$ times}{ d_{1,1}  & \cdots & d_{1,1}}
	  	&
	  \overmat{$\mkern-3.5mu\cdots$}{\cdots & \cdots}
	  	&
	  \bovermat{$n_s$ times}{ d_{1,s}  & \cdots & d_{1,s}}
	   	\\[0.5em]
	   	\vdots & & \vdots & & & \vdots & & \vdots \\
	   	d_{n, 1}  & \cdots & d_{n,1} & \cdots & \cdots & d_{n,s}  & \cdots & d_{n,s}
	   \end{pmatrix}
\end{matrix}
\end{align*}
The preceding observations together with \cref{equal-multpn} imply that
\begin{align}
\multfkniso
    \leq \multPkniso
    = \mv(\scrP_1, \ldots, \scrP_n)
	= \frac{ \perm(D(\mscrI, \vec d))}
		 {(\prod_j n_j!) (\prod_{j,k}\omega_{j,k})}
		 \label{wt-multi-homogeneous-bound}
\end{align}
Note that when $s = 1$, then $\Omega$ consists of only one weighted degree $\omega$ and $\perm(D(\mscrI, \vec d)) = n! \omega(f_1) \cdots \omega(f_n)$, so that the bound from \eqref{wt-multi-homogeneous-bound} is precisely the weighted B\'ezout bound $\prod_{j=1}^n(\omega(f_j)/\omega_j)$. Now we determine the condition for the attainment of this bound. For each $j = 1, \ldots, s$, fix a new indeterminate $u_j$. Let $B'_j$ be the ring of polynomials (over $\kk$) in $u_j$ and $x_k$, $k \in I_j$. Note that each $\xpp{\scrP_i}$ is isomorphic to $\prod_{j=1}^s \wwpnn{\omega'_j}{n_j}$, where $\omega'_j$ is the weighted degree on $B'_j$ such that $\omega'_j(u_j) = 1$ and $\omega'_j(x_k) = \omega_{j,k}$, $k \in I_j$. For each $g \in \kk[x_1, \ldots, x_n]$ and each nonempty subset $J$ of $[s]$, define $\ld_{\Omega, J}(g)$ as in \cref{multi-compactisection}. Finally, note that since $\mscrI$ is a partition of $[n]$, we may identify $\kk^n$ with $\prod_{j=1}^s \Kii{I_j}$.


\begin{thm}[Weighted multi-homogeneous B\'ezout theorem] \label{multi-bezout}
\index{Weighted!multi-homogeneous B\'ezout theorem!version I}
\index{B\'ezout's theorem!weighted multi-homogeneous!version I}
The number of isolated solutions of polynomials $f_1, \ldots, f_n$ on $\kk^n$ is bounded by \eqref{wt-multi-homogeneous-bound}. This bound is exact if and only if the following holds: for each nonempty subset $J$ of $[s]$, there is no common zero of $\ld_{\Omega, J}(f_1), \ldots, \ld_{\Omega, J}(f_n)$ on $\prod_{j \in J} (\Kii{I_j} \setminus \{\origin\}) \times \prod_{j \not\in J} \Kii{I_j}$ (upon identification of $\kk^n$ with $\prod_{j=1}^s \Kii{I_j}$).
\end{thm}

\begin{proof}
We have already proved that \eqref{wt-multi-homogeneous-bound} holds. Let $\tilde f_i$ be the {\em $\Omega$-homogenization} of $f_i$ defined as in \cref{multi-compactisection}. If the bound is not exact, then it follows by the same arguments as in the proof of \cref{wt-bezout-1} that $\bigcap_{j=1}^s V(\tilde f_j) \setminus \kk^n \neq \emptyset$, and then \cref{multi-closure-infty} implies that the condition in the second assertion of \cref{multi-bezout} is violated. Now assume there is a nonempty subset $J$ of $[s]$ such that $\ld_{\Omega, J}(f_1), \ldots, \ld_{\Omega, J}(f_n)$ have a common zero $a = (a_1, \ldots, a_n) \in \prod_{j \in J} (\Kii{I_j} \setminus \{\origin\}) \times \prod_{j \not\in J} \Kii{I_j}$. Replacing each $f_i$ by $f_i^m$ for some appropriate positive integer $m$, we may assume that the vertices of each $\scrP_i$ have integer coordinates. Let $I := \{i: a_i \neq 0\} \subseteq [n]$; in other words, $I$ is the smallest subset of $[n]$ such that $a \in \Kstari$. Since $J$ is nonempty, it follows that $I$ is also nonempty; in fact $I_j \cap I \neq \emptyset$ for each $j \in J$. For each $j \in J$, let $\scrT_j$ be the facet of $\scrS_j$ determined by $\langle \omega_j, \alpha \rangle = l_j$. Fix $i \in [n]$. Define
\begin{align*}
\scrQ_i := \ri \cap \left (\prod_{j \in J}(d_{i,j}/l_j)\scrT_j \times \prod_{j \not\in J} (d_{i,j}/l_j)\scrS_j \right)
\end{align*}
Then $\scrQ_i$ is a proper (nonempty) face of $\scrP_i$. Let $f_{i, \scrQ_i}$ be the component of $f_i$ supported at $\scrQ_i$, i.e.\ if $f_i = \sum_\alpha c_\alpha x^\alpha$, then $f_{i, \scrQ_i} = \sum_{\alpha \in \scrQ_i} c_\alpha x^\alpha$. It is straightforward to check that $f_{i,\scrQ_i}(a) = 0$ for each $i$. Now choose an integral element $\nu \in \rnstar$ such that $\scrQ_i = \In_\nu(\scrP_i)$ for each $i$ [you have to check that such $\nu$ exists!]. As in \cref{bkk-necessary-section}, pick BKK non-degenerate $g_1, \ldots, g_n$ with a common zero $b \in \kk^n$ such that for each $i$, $f_i(b) \neq 0$, $\np(g_i) = \scrP_i$ and $\In_\nu(g_i)(a) \neq 0$. Define a rational curve $C$ on $\kk^n$ via the parametrization from \eqref{trick-curve} and define $h_1, \ldots, h_n$ as in \eqref{trick-deformation} with $m_i = \nu(g_i)$ for each $i$. Then $t^{-m_j}f_j(c(t))$ and $t^{-m_j}g_j(c(t))$ are polynomials in $t$ and the same arguments as in \cref{bkk-necessary-section} show that $\multfkniso < \multgkniso$, as required.
\end{proof}

\section{Notes}
A version of the weighted B\'ezout theorem appears in \cite{damon}. I.\ Shafarevich gave the bound for the ``multi-homogeneous'' case  (i.e.\ the case that all weights are $1$) of the weighted multi-homogeneous B\'ezout theorem in the first edition of \cite{shaf1} in the 1970s; see also \cite{morgan-sommese}. We could not locate any past reference for the non-degeneracy condition or the estimate for the general weighted multi-homogeneous version.


\chapter{Intersection multiplicity at the origin} \label{multiplicity-chapter}
\section{Introduction}
In this chapter we consider the ``local'' version of the affine B\'ezout problem, i.e.\ the problem of estimating the intersection multiplicity of generic hypersurfaces at the origin. This computation is a crucial ingredient of the extension in \cref{affine-chapter} of Bernstein's theorem to the affine space. 
Recall that the {\em support} of a power series $f = \sum_\alpha c_\alpha x^\alpha \in \kk[[x_1, \ldots, x_n]]$ is $\supp(f):= \{\alpha: c_\alpha \neq 0\}$ and we say that $f$ is {\em supported} at $\scrA\subset \zz^n$ if $\supp(f) \subset \scrA$. Now let $\scrA_1, \ldots, \scrA_n$ be (possibly infinite) subsets of $\zzeroo{n}$.  In the case that $\scrA_j$ are finite, we saw in \cref{bkk-chapter} that within all $f_j$ supported at $\scrA_j$, $j = 1, \ldots, n$, $\multfntorusiso$ takes the maximum value when $f_1, \ldots, f_n$ are generic. It is possible to talk about ``generic'' power series supported at $\scrA_j$ even if $\scrA_j$ is infinite, and it turns out that the intersection multiplicity $\multfzero$ of $f_1, \ldots, f_n$ at the origin takes the {\em minimum} value when $f_j$ are generic power series supported at $\scrA_j$, $j = 1, \ldots, n$ (see \cref{non-degeneracy-0-thm} for the precise statement); in this chapter we compute this minimum and give a Bernstein-Kushnirenko type characterization of the systems which attain the minimum.

\section{Generic intersection multiplicity}
\index{Generic intersection multiplicity}
\index{Intersection multiplicity!generic}
Let $\scrA_j$ be a (possibly infinite) subset of $\zzeroo{n}$, $j = 1, \ldots, n$. Define
\begin{align}
\multAzero &:= \min\{\multfzero:\forall j,\ f_j \in \kk[[x_1, \ldots, x_n]],\ \supp(f_j) \subseteq \scrA_j \}
	\label{multAzero}
\end{align}
In this section we motivate, state and illustrate the formula for $\multAzero$. Its proof is given in \cref{mult-bound-section}.

\begin{center}
\begin{figure}[h]
\begin{tikzpicture}[scale=0.45]
\def\shiftone{9}
\def\opazero{0.5}
\def\tx{2.5}
\def\ty{2.5}
\def\gridx{4.5}
\def\gridy{4.5}
\def\colorzero{blue}
\def\colorone{blue}

\draw [gray,  line width=0pt] (-0.5,-0.5) grid (\gridx,\gridy);
\draw [<->] (0, \gridy) |- (\gridx, 0);
\draw[ultra thick, \colorone] (2,\gridy) --  (0,0) -- (\gridx,1);
\fill[\colorzero, opacity=\opazero ] (0,0) -- (2,\gridy) -- (\gridx,\gridy) -- (\gridx,1) -- cycle;
\draw (\tx,\ty) node {\picfontsize $\scrA_1$};

\begin{scope}[shift={(\shiftone,0)}]
	\draw [gray,  line width=0pt] (-0.5,-0.5) grid (\gridx,\gridy);
	\draw [<->] (0, \gridy) |- (\gridx, 0);
	\draw[ultra thick, \colorone]  (\gridx,1) -- (2,1) -- (1,2) -- (0,4) -- (0, \gridy);
	\fill[\colorzero, opacity=\opazero ] (\gridx,1) -- (2,1) -- (1,2) -- (0,4) -- (0,\gridy) -- (\gridx,\gridy) -- cycle;
	\draw (\tx,\ty) node {\picfontsize $\scrB_1$};
\end{scope}

\begin{scope}[shift={(2*\shiftone,0)}]
	\draw [gray,  line width=0pt] (-0.5,-0.5) grid (\gridx,\gridy);
	\draw [<->] (0, \gridy) |- (\gridx, 0);
	\draw[ultra thick, \colorone]  (4,2) -- (3,1) -- (0,2)-- (2,4) -- (4,4) -- cycle;
	\fill[\colorzero, opacity=\opazero ] (4,2) -- (3,1) -- (0,2)-- (2,4) -- (4,4) -- cycle;
	\draw (\tx,\ty) node {\picfontsize $\scrB_2$};
\end{scope}

\end{tikzpicture}

\caption{$\multpzeronodots{\scrA_1, \scrA_2} = 0$ for any $\scrA_2$ since a generic $f_1$ supported at $\scrA_1$ has a nonzero constant term; $\multpzeronodots{\scrB_1, \scrB_2} = \infty$ since every $g_j$ supported at $\scrB_j$, $j = 1, 2$, identically vanishes on the $x$-axis.} \label{fig:zeroinftyatzero}
\end{figure}
\end{center}

\subsection{Motivation for the formula} \label{mult-motivation}
In this section we study informally the case $\kk = \cc$ and try to motivate the formula for $\multAzero$. It is not hard to understand precisely when $\multAzero$ is zero or infinity - see \cref{fig:zeroinftyatzero,mult-support}. So consider the case that $0 < \multAzero < \infty$, and pick $f_j$ supported at $\scrA_j$, $j = 1, \ldots, n$, such that $\multfzero = \multAzero$. Then the origin is an isolated point of $f_1 = \cdots = f_n = 0$. Therefore \cref{order-curve,int-mult-curve} imply that near the origin $f_2 = \cdots = f_n = 0$ defines a curve $C$, and
\begin{align}
\multfzero = \sum_B \ord_\origin(f_1|_B)
\label{multfzero:sum-of-order}
\end{align}
where the sum is over all ``branches'' $B$ of $C$ at the origin. We now try to compute this sum. So fix a branch $B$ of $C$ at the origin and an analytic parametrization $\gamma = (\gamma_1, \ldots, \gamma_n): U \to B$ of $B$, where $U$ is a neighborhood of the origin on $\cc$. Let $H_B$ be the {\em smallest} coordinate subspace of $\cc^n$ containing $B$.

\subsubsection{Case 1: $H_B = \cc^n$.} \label{H_B=C^n}
 In this case no $\gamma_i$ is identically zero, so that each $\gamma_i$ can be expressed as
\begin{align*}
  \gamma_i = a_1t^{\nu_i} + \cdots
\end{align*}
where $a_i \in \cc^*$, $t$ is an analytic coordinate on $U$, and $\nu_i$ is the order (in $t$) of $\gamma_i$. Note that each $\nu_i$ is {\em positive}, since the center of $B$ is at the origin. Let $\nu_B$ be the weighted order on $\kk[x_1, \ldots, x_n]$ such that $\nu_B(x_i) = \nu_i$, $i = 1, \ldots, n$. Then for each $j$,
\begin{align*}
f_j(\gamma(t)) = \In_{\nu_B}(f_j)(a_1, \ldots, a_n) t^{\nu(f_j)} + \cdots
\end{align*}
where $\In_{\nu_B}(\cdot)$ denotes the {\em initial form} with respect to $\nu_B$, and the omitted terms have higher order in $t$. Since $f_j(\gamma(t)) \equiv 0$ for $j = 2, \ldots, n$, we have a system of $(n-1)$ {\em weighted homogeneous} polynomial equations:
\begin{align}
\In_{\nu_B}(f_j)(a_1, \ldots, a_n) = 0,\ j = 2, \ldots, n.
\label{initial-equation}
\end{align}
\Cref{bkk-bound-thm} implies that the number of solutions of \eqref{initial-equation} is the $(n-1)$-dimensional mixed volume of the Newton polytopes of $\In_{\nu_B}(f_j)$, $j = 2, \ldots, n$. It is in fact ``reasonable'' to guess that in the generic case each solution of \eqref{initial-equation} corresponds to a distinct branch $B$ of $C$ and for each such $B$, the order of $f_1|_B$ at the origin is $\nu_B(f_1)$. In other words, for each weighted order $\nu$,
\begin{align}
\sum_{\substack{H_B = \cc^n \\ \nu_B = \nu}} \ord_\origin(f_1|_B) = \nu(f_1) \mv(\In_\nu(f_2), \ldots, \In_\nu(f_n))
\label{mult-motivate-equation-0}
\end{align}

\subsubsection{Case 2: $H_B \subsetneq \cc^n$.} In this case we may assume \woutlog\ that $H_B$ is the coordinate subspace spanned by $x_1, \ldots, x_k$, $k < n$; in other words $H_B = \Ci$ in the notation of \eqref{ki} from \cref{wt-bezout-proof-section}, where $I := \{1, \ldots, k\}$. It follows that
\begin{align*}
\gamma_i
    &=
    \begin{cases}
    a_i t^{\nu_i} + \cdots & \text{if}\ i = 1, \ldots, k, \\
    0 & \text{if}\ i = k + 1, \ldots, n,
    \end{cases}
\end{align*}
where $(a_1, \ldots, a_k) \in \nctoruss{k}$. Since $\multfzero < \infty$ and since the $f_j$ are generic among those supported at $\scrA_j$, it follows that there are precisely $(n - k)$ elements among $f_2, \ldots, f_n$ which vanish identically on $\Ci$; after reordering the $f_j$ we may assume these are $f_{k+1}, \ldots, f_n$. Since $f_j(\gamma(t)) \equiv 0$ for each $j = 2, \ldots, k$, it follows that
\begin{align}
\In_{\nu_B}(f_j|_{\Ci})(a_1, \ldots, a_k) = 0,\ j = 2, \ldots, k.
\label{initial-equation-proper}
\end{align}
where $\nu_B$ is the weighted order on $\kk[x_1, \ldots, x_k]$ corresponding to weights $\nu_i$ for $x_i$, $i = 1, \ldots, k$. As in the preceding case, the number of solutions of \eqref{initial-equation-proper} in the generic situation is the $(k-1)$-dimensional mixed volume of the Newton polytopes of $\In_{\nu_B}(f_j|_{\Ci})$, $j = 2, \ldots, k$, and each solution corresponds to a distinct branch $B$ of $C$. However, each branch should be counted with proper multiplicity, and therefore this mixed volume should be multiplied by the ``intersection multiplicity of $f_{k+1}, \ldots, f_{n}$ along $\Ci$.'' It turns out (see \cref{mult-chain-1}) that for generic $f_{k+1}, \ldots, f_{n}$, this is precisely the intersection multiplicity at the origin of $f_{k+1}|_{(x_1, \ldots, x_k) = (\epsilon_1,\ldots, \epsilon_k)}, \ldots, f_n|_{(x_1, \ldots, x_k) = (\epsilon_1,\ldots, \epsilon_k)}$, where $\epsilon_1, \ldots, \epsilon_k$ are generic elements from $\cc^*$. If $\pi: \rr^n \to \rr^{n-k}$ is the projection onto the last $(n-k)$-coordinates, then the genericness of the $\epsilon_i$ imply that the support of $f_{k+j}|_{(x_1, \ldots, x_k) = (\epsilon_1,\ldots, \epsilon_k)}$ is precisely $\pi(\scrA_{k+j})$ for each $j$, and it is reasonable to guess that if the $f_{k+j}$ and $\epsilon_i$ are generic, then
\begin{align*}
\multzero{f_{k+1}|_{(x_1, \ldots, x_k) = (\epsilon_1,\ldots, \epsilon_k)}}{f_n|_{(x_1, \ldots, x_k) = (\epsilon_1,\ldots, \epsilon_k)}}
	&= \multzero{\pi(\scrA_{k+1})}{\pi(\scrA_{n})}
\end{align*}
It then follows as in the first case that for each weighted order $\nu$ on $\kk[x_1, \ldots, x_k]$,
\begin{align}
\begin{split}
\sum_{H_B = \Ci, \nu_B = \nu}
    & \ord_\origin(f_1|_B) \\
	&= \nu(f_1|_{H_B}) \mv(\In_\nu(f_2|_{H_B}), \ldots, \In_\nu(f_k|_{H_B}))  \\
    & \qquad  \times
	   \multzero{f_{k+1}|_{(x_1, \ldots, x_k) = (\epsilon_1,\ldots, \epsilon_k)}}{f_n|_{(x_1, \ldots, x_k) = (\epsilon_1,\ldots, \epsilon_k)}} \\
	&=  \min_{\scrA_1 \cap \ri} (\nu) \mv(\In_\nu(\scrA_2 \cap  \ri), \ldots, \In_\nu(\scrA_k \cap \ri))
        \multzero{\pi(\scrA_{k+1})}{\pi(\scrA_{n})}
\end{split}
\label{mult-motivation-equation}
\end{align}
Therefore $\multAzero$ should be the sum of the right hand side of identity \eqref{mult-motivation-equation} over all appropriate $I$ and $\nu$. \Cref{multiplicity-thm} states that this precisely the case, and the proof of \cref{multiplicity-thm} in \cref{multiplicity-proof-section} simply makes the preceding arguments rigorous.

\begin{center}
\begin{figure}[h]
\def\colorone{blue}
\def\colortwo{red}
\begin{tikzpicture}[
    scale=0.45,
	dot/.style = {
      draw,
      fill,
      circle,
      inner sep = 0pt,
      minimum size = 3pt,
    }
    ]
\def\shiftone{9}
\def\opazero{0.5}
\def\tx{2.5}
\def\ty{2.5}
\def\gridx{4.5}
\def\gridy{4.5}

\draw [gray,  line width=0pt] (-0.5,-0.5) grid (\gridx,\gridy);
\draw [<->] (0, \gridy) |- (\gridx, 0);
\draw[\colorone] (2,\gridy) --  (0,0) -- (\gridx,1);
\fill[\colorzero, opacity=\opazero ] (0,0) -- (2,\gridy) -- (\gridx,\gridy) -- (\gridx,1) -- cycle;
\draw (0,0) node[dot, \colortwo] {};
\draw (\tx,\ty) node {\picfontsize $\scrA_1$};

\begin{scope}[shift={(\shiftone,0)}]
	\draw [gray,  line width=0pt] (-0.5,-0.5) grid (\gridx,\gridy);
	\draw [<->] (0, \gridy) |- (\gridx, 0);
	\draw[\colorone]  (\gridx,1) -- (2,1) -- (1,2) -- (0,4) -- (0, \gridy);
	\fill[\colorzero, opacity=\opazero ] (\gridx,1) -- (2,1) -- (1,2) -- (0,4) -- (0,\gridy) -- (\gridx,\gridy) -- cycle;
	\draw[ultra thick, \colortwo]  (2,1) -- (1,2) -- (0,4);
	\draw (\tx,\ty) node {\picfontsize $\scrB_1$};
\end{scope}

\begin{scope}[shift={(2*\shiftone,0)}]
	\draw [gray,  line width=0pt] (-0.5,-0.5) grid (\gridx,\gridy);
	\draw [<->] (0, \gridy) |- (\gridx, 0);
	\draw[\colorone]  (4,2) -- (3,1) -- (0,2)-- (2,4) -- (4,4) -- cycle;
	\fill[\colorzero, opacity=\opazero ] (4,2) -- (3,1) -- (0,2)-- (2,4) -- (4,4) -- cycle;
	\draw[ultra thick, \colortwo]  (3,1) -- (0,2);
	\draw (\tx,\ty) node {\picfontsize $\scrB_2$};
\end{scope}

\end{tikzpicture}

\caption{Newton diagrams (in \colortwo) of the sets from \cref{fig:zeroinftyatzero}} \label{fig:nd}
\end{figure}
\end{center}

\subsection{Precise formulation} \label{mult-bound-statement-section}
Let $\scrA$ be a (possibly infinite) subset of $\zzeroo{n}$. The convex hull of $\scrA + \rzeroo{n}$ in $\rr^n$ is a convex polyhedron (\cref{corner-convex}); the \index{Newton!diagram}{\em Newton diagram} $\nd(\scrA)$ of $\scrA$ is the union of the compact faces of this polyhedron (\cref{fig:nd}). The Newton diagram of a power series $f$ in $(x_1, \ldots, x_n)$, denoted $\nd(f)$, is the Newton diagram of $\supp(f)$; it is the local analogue of the Newton polytope of a polynomial. 
Given diagrams $\Gamma_1, \ldots, \Gamma_n$ in $\rr^n$, define
\begin{align}
\multGammazero &:= \min\{\multfzero: \forall j,\ f_j \in \kk[[x_1, \ldots, x_n]],\ \nd(f_j) + \rzeroo{n} \subseteq \Gamma_j  + \rzeroo{n} \}
	\label{multGammazero}
\end{align}
We will see in \cref{multiplicity-thm} below that $\multGammazero$ can be expressed in terms of certain mixed volumes of the faces of $\Gamma_j$, and if $\Gamma_j$ are Newton diagrams of $\scrA_j \subseteq \znzero$, then $\multAzero = \multGammazero$. First we need to introduce some notation. We write $[n] := \{1, \ldots, n\}$. If $I \subseteq [n]$ and $k$ is a field, recall that $\ki$ is the $|I|$-dimensional coordinate subspace $\{(x_1, \ldots, x_n) \in k^n: x_i =  0\ \text{if}\ i \not\in I\}$, and $\kstari = \{(x_1, \ldots, x_n) \in \ki: x_i \neq 0\ \text{if}\ i \in I\}$. For $\scrS \subset \rr^n$, we write $\scrS^I := \scrS \cap \ri$. We denote by $\pi_I: k^n \to \ki$ the projection in the coordinates indexed by $I$, i.e.\
\begin{align}
\text{the $j$-th coordinate of $\pi_I(x_1, \ldots, x_n)$}&:=
			\begin{cases}
			x_j &\text{if}\ j \in I\\
			0	 &\text{if}\ j \not\in I.
			\end{cases} \label{pi_I}
\end{align}


Let $\nu$ be a weighted order on $\kk[x_1, \ldots, x_n]$ corresponding to weights $\nu_j$ for $x_j$, $j = 1, \ldots, n$. We identify $\nu$ with the element in $\rnstar$ with coordinates $(\nu_1, \ldots, \nu_n)$ with respect to the dual basis. We say that $\nu$ is \index{Weighted!order!centered at the origin}{\em centered at the origin} if each $\nu_i$ is positive and that $\nu$ is \index{Primitive!weighted order}{\em primitive} if it is nonzero and the greatest common divisor of $\nu_1, \ldots, \nu_n$ is $1$. If $\nu$ is centered at the origin, then it also extends to a weighted order on the ring of power series in $(x_1, \ldots, x_n)$. We write $\Vzero$ for the set of weighted orders centered at the origin and $\Vzeroprime$ for the primitive elements in $\Vzero$. Given polytopes $\Gamma_1, \ldots, \Gamma_n$ in $\rr^n$, define
\begin{align}
\multGammazerostar
	&:=  \sum_{\nu \in \Vzeroprime} \min_{\Gamma_1}(\nu) ~
		\mv'_\nu(\In_\nu(\Gamma_2), \ldots, \In_\nu(\Gamma_n))
	\label{mult-star}
\end{align}
where $\mv'_\nu(\cdot, \ldots, \cdot)$ is defined as in \eqref{mv'}.

\begin{thm}[{\cite{toricstein}}]\label{multiplicity-thm}
Let $\scrA := (\scrA_1, \ldots, \scrA_n)$ be a collection of subsets of $\zzeroo{n}$ and $\Gamma_j$ be the Newton diagram of $\scrA_j$, $j = 1, \ldots, n$. For each $I \subset [n]$, let $\tiA := \{j: \Ai_j \neq \emptyset\}$ be the set of all indices $j$ such that $\scrA_j$ touches the coordinate subspace $\ri$ of $\rr^n$. Define
\begin{align}
\tAone:= \{I \subseteq [n]: I \neq \emptyset,\ |\tiA| = |I|,\ 1 \in \tiA\} \label{I1-list}
\end{align}
Then
\begin{enumerate}
\item If $0 \not\in \bigcup_j \Gamma_j$ and there is $I \subset [n]$ such that $|\tiA| < |I|$, then $$\multAzero = \multGammazero = \infty$$
\item Otherwise
\begin{align}
\multAzero
	= \multGammazero
	= \sum_{I \in \tAone}
							\multzerostar{\Gamma^{I}_1, \Gamma^{I}_{j_2}}{\Gamma^{I}_{j_{|I|}}}
							\times
							\multzero{\pi_{[n]\setminus I}(\Gamma_{j'_1})}{\pi_{[n]\setminus I}(\Gamma_{j'_{n-|I|}})}
	\label{mult-formula}
\end{align}
where for each $I \in \tAone$, $j_1 = 1, j_2, \ldots, j_{|I|}$ are elements of $\tiA$, and $j'_1, \ldots, j'_{n-|I|}$ are elements of $[n]\setminus \tiA$.
\end{enumerate}
\end{thm}

\begin{rem} \label{formula-convention}
The product of $0$ and $\infty$, when/if it occurs in \eqref{mult-formula}, is defined to be $0$. Also empty intersection products and mixed volumes are defined as $1$. In particular, when $n = 1$, the term $\mv(\In_\nu(\Gamma_2), \ldots, \In_\nu(\Gamma_n))$ from \eqref{mult-star} is defined to be $1$.
\end{rem}

\begin{rem}[Generic intersection multiplicity is monotonic] \label{monotonic-remark-0}
The formulae for $\multAzero$ from \cref{multiplicity-thm} do not change if the $\scrA_j$ are replaced by $\scrA_j + \rzeroo{n}$. This immediately implies that $\multzero{\cdot}{\cdot}$ is {\em monotonic}, i.e.\ if  $\scrA'_j \subseteq \scrA_j + \rzeroo{n}$, $j = 1, \ldots, n$, then $\multzero{\scrA'_1}{\scrA'_n} \geq \multAzero$. Precise characterization of the cases for which $\multzero{\scrA'_1}{\scrA'_n} > \multAzero$ is given in \cref{strictly-monotone-0}.
\end{rem}

It is not obvious from the outset that the term computed by \eqref{mult-formula} is invariant under the permutations of the $\scrA_j$. Some formulae which are invariant under permutations of the $\scrA_j$ are given in \cref{other-0-section}. We now present an example to illustrate this invariance.

\def\shiftone{7.5}
\def\colorzero{blue}
\def\colorone{red}
\def\colortwo{yellow}
\def\colorfour{green}
\def\opazero{0.5}
\def\viewx{75}
\def\viewy{30}
\def\titlex{1}
\def\titley{-1}

\begin{figure}[h]
\begin{center}
\begin{tikzpicture}[scale=0.6]
\pgfplotsset{every axis title/.append style={at={(0,-0.2)}}, view={\viewx}{\viewy}, axis lines=middle, enlargelimits={upper}}

\begin{scope}
\begin{axis}
	\addplot3[fill=\colorzero,opacity=\opazero] coordinates{(1,0,0) (0,1,0) (0,0,1)};
	\addplot3 [ultra thick, \colorfour, ->] coordinates{(0.3,0.3,0.4) (0.8, 0.8, 0.9)};
	\draw (axis cs:0.8,0.8,0.9) node [above] {(1,1,1)};
\end{axis}
\draw (\titlex,\titley) node {$\Gamma_1$};
\end{scope}

\begin{scope}[shift={(\shiftone,0)}]
\begin{axis}
	\addplot3[fill=\colorone,opacity=\opazero] coordinates{(3,0,0) (0,3,0) (0,1,2) (1,0,2) (3,0,0)};
	\addplot3 [ultra thick, \colorfour, ->] coordinates{(0.9,0.9,1.2) (2.4, 2.4, 2.7)};
	\draw (axis cs:2.4, 2.4, 2.7) node [above] {(1,1,1)};
\end{axis}
\draw (\titlex,\titley) node {$\Gamma_2$};
\end{scope}

\begin{scope}[shift={(2*\shiftone,0)}]
\begin{axis}
	\addplot3[fill=\colortwo,opacity=\opazero] coordinates{(2,0,0) (0,2,0) (0,1,2) (1,0,2) (2,0,0)};
	\addplot3 [ultra thick, \colorfour, ->] coordinates{(0.8,0.8,0.8) (1.8, 1.8, 1.3)};
	\draw (axis cs:1.8, 1.8, 1.3) node [above] {(2,2,1)};
\end{axis}
\draw (\titlex,\titley) node {$\Gamma_3$};
\end{scope}
\end{tikzpicture}
\end{center}
\caption{Newton diagrams of polynomials from \cref{ex0} and inner normals to their facets} \label{fig-ex0}
\end{figure}

\begin{example} \label{ex0}
Consider the polynomials $f_1, f_2, f_3$ from \cref{ex-mv'}. If $\scrA := (\supp(f_1), \supp(f_2), \supp(f_3))$, then it is straightforward to check (see \cref{fig-ex0}) that $\tAone= \{\{1,2,3\},\{3\}\}$, so that \eqref{mult-formula} implies that
\begin{align*}
[f_1,f_2,f_3]_0
	&= [\Gamma_1, \Gamma_2, \Gamma_3]^*_0 +
		 [\pi_{\{1,2\}}(\Gamma_2), \pi_{\{1,2\}}(\Gamma_3)]_0 ~ [\Gamma^{\{3\}}_1]^*_0 
\end{align*}
The Newton diagrams of $\pi_{\{1,2\}}(\Gamma_2)$ and $\pi_{\{1,2\}}(\Gamma_3)$ are the same diagram consisting of the line segment from $(1,0)$ to $(0, 1)$, which is the Newton diagram of linear polynomials with no constant terms. It follows that $[\pi_{\{1,2\}}(\Gamma_2), \pi_{\{1,2\}}(\Gamma_3)]_0 = 1$. The Newton diagram of $\Gamma_2 + \Gamma_3$ has two facets with inner normals in $(\rr_{> 0})^3$, and these inner normals are $\nu_1 := (1,1,1)$ and $\nu_2 := (2,2,1)$ (see \cref{fig:ex0-1}). Then it follows from \cref{fig:ex0-2} and identity \eqref{vol-to-mixed-2} that
\def\scalefactor{0.3}
\begin{align*}
[f_1,f_2,f_3]_0
	&= \min_{\Gamma_1}(\nu_1) \mv'_{\nu_1}(\In_{\nu_1}(\Gamma_2), \In_{\nu_1}(\Gamma_3))
				+ \min_{\Gamma_1}(\nu_2) \mv'_{\nu_2}(\In_{\nu_2}(\Gamma_2), \In_{\nu_2}(\Gamma_3))
				+ 1 \cdot \ord_z(f_1|_{x=y=0}) \\
		&= \min_{\Gamma_1}(1,1,1) \cdot \area(
					\begin{tikzpicture}[scale=\scalefactor]
							\draw[fill=\colorfour, opacity=\opazero] (3,0) -- (5,0) -- (3,2) -- (1,2) -- cycle;
					\end{tikzpicture}
					)
				+  \min_{\Gamma_1}(2,2,1) \cdot \area(
					\begin{tikzpicture}[scale=\scalefactor]
							\draw[fill=\colorfour, opacity=\opazero] (2,0) -- (3,0) -- (2,1) -- (1,1) -- cycle;
					\end{tikzpicture}
					)
				+ 1 \cdot 1
	= 1 \cdot 4 + 1 \cdot 1 + 1 = 6.
\end{align*}

\begin{figure}[h]
\begin{center}
\begin{tikzpicture}[scale=0.6]
\pgfplotsset{every axis title/.append style={at={(0,-0.2)}}, view={\viewx}{\viewy}, axis lines=middle, enlargelimits={upper}}

\begin{scope}
\begin{axis}
	\addplot3[fill=\colorzero,opacity=\opazero] coordinates{(4,0,0) (0,4,0) (0,1,3) (1, 0,3) (4, 0, 0)};
\end{axis}
\draw (\titlex,\titley) node {$\nd(\Gamma_1 + \Gamma_2)$};
\end{scope}

\begin{scope}[shift={(\shiftone,0)}]
\begin{axis}
	\addplot3[fill=\colorone,opacity=\opazero] coordinates{(5,0,0) (0,5,0) (0,3,2) (3,0,2) (5,0,0)};
	\addplot3[fill=\colortwo,opacity=\opazero] coordinates{(3,0,2) (0,3,2) (0,2,4) (2,0,4) (3,0,2)};
\end{axis}
\draw (\titlex,\titley) node {$\nd(\Gamma_2 + \Gamma_3)$};
\end{scope}

\begin{scope}[shift={(2*\shiftone,0)}]
\begin{axis}
	\addplot3[fill=\colorzero,opacity=\opazero] coordinates{(3,0,0) (0,3,0) (0,2,1) (2,0,1) (3,0,0)};
	\addplot3[fill=\colortwo,opacity=\opazero] coordinates{(2,0,1) (0,2,1) (0,1,3) (1,0,3) (2,0,1)};
\end{axis}
\draw (\titlex,\titley) node {$\nd(\Gamma_3 + \Gamma_1)$};
\end{scope}
\end{tikzpicture}
\end{center}
\caption{Sum of the Newton diagrams of polynomials from \cref{ex0}} \label{fig:ex0-1}
\end{figure}

On the other hand with $\scrA' := (\supp(f_3), \supp(f_1), \supp(f_2))$, one has $\tAprimeone = \{\{1,2,3\}\}$. Since the Newton diagram of $\Gamma_1  + \Gamma_2$ has only one facet and that the primitive inner normal to that facet is $\nu_1$, we have from \cref{fig:ex0-2} and identity \eqref{vol-to-mixed-2} that
\begin{align*}
[f_1,f_2,f_3]_0
	&= [\Gamma_3, \Gamma_1, \Gamma_2]^*_0
	= \min_{\Gamma_3}(\nu_1) \mv'_{\nu_1}(\In_{\nu_1}(\Gamma_1), \In_{\nu_1}(\Gamma_2)) \\
	&
	= \min_{\Gamma_3}(1,1,1) \cdot \area(
			\begin{tikzpicture}[scale=\scalefactor]
					\draw[fill=\colorfour, opacity=\opazero] (3,0) -- (4,0) -- (1,3) -- (0,3) -- cycle;
			\end{tikzpicture}
			)
	= 2 \cdot 3
	= 6.
\end{align*}
Similarly, with $\scrA'' := (\supp(f_2), \supp(f_3), \supp(f_1))$, one has $\tAAone{\scrA''} = \{\{1,2,3\}\}$. The Newton diagram of $\Gamma_3  + \Gamma_1$ have two facets, with inner normals $\nu_1$ and $\nu_2$, we have from \cref{fig:ex0-2} and identity \eqref{vol-to-mixed-2} that
\begin{align*}
[f_1,f_2,f_3]_0
	&= [\Gamma_2, \Gamma_3, \Gamma_1]^*_0 \\
	&= \min_{\Gamma_2}(\nu_1)\mv'_{\nu_1}(\In_{\nu_1}(\Gamma_3), \In_{\nu_1}(\Gamma_1))
		+
		\min_{\Gamma_2}(\nu_2)\mv'_{\nu_2}(\In_{\nu_2}(\Gamma_3), \In_{\nu_2}(\Gamma_1))  \\
	&= \min_{\Gamma_2}(1,1,1) \cdot \area(
						\begin{tikzpicture}[scale=\scalefactor]
								\draw[fill=\colorfour, opacity=\opazero] (1,0) -- (3,0) -- (2,1) -- (0,1) -- cycle;
						\end{tikzpicture}
						)
					+  \min_{\Gamma_2}(2,2,1) \cdot \area(\emptyset)
	= 3 \cdot 2 + 4 \cdot 0 = 6.
\end{align*}
\end{example}

\begin{figure}[h]
\begin{center}
\begin{tikzpicture}[scale=0.4]
\def\shifttwo{6}
\def\nux{-0.5}
\def\nuy{4.5}
\def\opazero{0.5}
\def\tx{-0.5}
\def\ty{-0.5}
\def\gridx{4.5}
\def\gridy{3.5}
\def\bigshiftdver{2}

\draw (\nux,\nuy) node [right] {\picfontsize $\nu_1 = (1,1,1)$};

\draw [gray,  line width=0pt] (-0.5,-0.5) grid (\gridx,\gridy);
\draw[ultra thick, fill=\colorzero, opacity=\opazero ] (0,0) -- (1,0) -- (0,1) -- cycle;
\draw (\tx,\ty) node [below right] {\picfontsize $\In_{\nu_1}(\Gamma_1)$};    	

\begin{scope}[shift={(\shifttwo,0)}]
	\draw [gray,  line width=0pt] (-0.5,-0.5) grid (\gridx,\gridy);
	\draw[ultra thick, fill=\colorone, opacity=\opazero ] (0,0) -- (3,0) -- (1,2) -- (0,2) -- cycle;
	\draw (\tx,\ty) node [below right] {\picfontsize $\In_{\nu_1}(\Gamma_2)$};    	
\end{scope}

\begin{scope}[shift={(2*\shifttwo,0)}]
	\draw [gray,  line width=0pt] (-0.5,-0.5) grid (\gridx,\gridy);
	\draw[ultra thick, \colortwo ] (0,0) -- (2,0);
	\draw (\tx,\ty) node [below right] {\picfontsize $\In_{\nu_1}(\Gamma_3)$};    	
\end{scope}

\begin{scope}[shift={(3*\shifttwo,0)}]
	\draw [gray,  line width=0pt] (-0.5,-0.5) grid (\gridx,\gridy);
	\draw[ultra thick] (0,0) -- (4,0) -- (1,3) --  (0,3) --  cycle;
	\draw[fill=\colorzero, opacity=\opazero, dashed, ultra thick] (0,2) -- (1,2) -- (0,3) -- cycle;
	\draw[fill=\colorone, opacity=\opazero, dashed, ultra thick ] (0,0) -- (3,0) -- (1,2) -- (0,2) -- cycle;
	\draw[fill=\colorfour, opacity=\opazero] (3,0) -- (4,0) -- (1,3) -- (0,3) -- cycle;
	\draw (\tx,\ty) node [below right] {\picfontsize $\In_{\nu_1}(\Gamma_1 + \Gamma_2)$};    	
\end{scope}

\begin{scope}[shift={(4*\shifttwo,0)}]
	\draw [gray,  line width=0pt] (-0.5,-0.5) grid (\gridx,\gridy);
	\draw[ultra thick] (0,0) -- (5,0) -- (3,2) --  (0,2) --  cycle;
	\draw[fill=\colorone, opacity=\opazero, dashed, ultra thick ] (0,0) -- (3,0) -- (1,2) -- (0,2) -- cycle;
	\draw[fill=\colorfour, opacity=\opazero] (3,0) -- (5,0) -- (3,2) -- (1,2) -- cycle;
	\draw (\tx,\ty) node [below right] {\picfontsize $\In_{\nu_1}(\Gamma_2 + \Gamma_3)$};    	
\end{scope}
	
\begin{scope}[shift={(5*\shifttwo,0)}]
	\draw [gray,  line width=0pt] (-0.5,-0.5) grid (\gridx,\gridy);
	\draw[ultra thick] (0,0) -- (3,0) -- (2,1) --  (0,1) --  cycle;
	\draw[fill=\colorzero, opacity=\opazero, dashed, ultra thick] (0,0) -- (1,0) -- (0,1) -- cycle;
	\draw[fill=\colorfour, opacity=\opazero] (1,0) -- (3,0) -- (2,1) -- (0,1) -- cycle;
	\draw (\tx,\ty) node [below right] {\picfontsize $\In_{\nu_1}(\Gamma_3 + \Gamma_1)$};    	
\end{scope}

\begin{scope}[shift={(0,-\nuy-\bigshiftdver)}]
	\def\gridy{1.5}
	\def\nuy{2.5}

	\draw (\nux,\nuy) node [right] {\picfontsize $\nu_2 = (2,2,1)$};
	
	\draw [gray,  line width=0pt] (-0.5,-0.5) grid (\gridx,\gridy);
	\tikzstyle{dot} = [blue, circle, minimum size=3pt, inner sep = 0pt, fill]
	\node[dot] at (0,0) {};
	\draw (\tx,\ty) node [below right] {\picfontsize $\In_{\nu_2}(\Gamma_1)$};    	
	
	\begin{scope}[shift={(\shifttwo,0)}]
		\draw [gray,  line width=0pt] (-0.5,-0.5) grid (\gridx,\gridy);
		\draw[ultra thick, \colorone ] (0,0) -- (1,0);
		\draw (\tx,\ty) node [below right] {\picfontsize $\In_{\nu_2}(\Gamma_2)$};    	
	\end{scope}
	
	\begin{scope}[shift={(2*\shifttwo,0)}]
		\draw [gray,  line width=0pt] (-0.5,-0.5) grid (\gridx,\gridy);
		\draw[ultra thick, fill=\colortwo, opacity=\opazero ] (0,0) -- (2,0) -- (1,1) -- (0,1) -- cycle;
		\draw (\tx,\ty) node [below right] {\picfontsize $\In_{\nu_2}(\Gamma_3)$};    	
	\end{scope}
	
	\begin{scope}[shift={(3*\shifttwo,0)}]
		\draw [gray,  line width=0pt] (-0.5,-0.5) grid (\gridx,\gridy);
		\draw[ultra thick, \colorone ] (0,0) -- (1,0);
		\draw (\tx,\ty) node [below right] {\picfontsize $\In_{\nu_2}(\Gamma_1 + \Gamma_2)$};    	
	\end{scope}
	
	\begin{scope}[shift={(4*\shifttwo,0)}]
		\draw [gray,  line width=0pt] (-0.5,-0.5) grid (\gridx,\gridy);
		\draw[ultra thick] (0,0) -- (3,0) -- (2,1) --  (0,1) --  cycle;
		\draw[ultra thick, fill=\colortwo, opacity=\opazero, dashed, ultra thick ] (0,0) -- (2,0) -- (1,1) -- (0,1) -- cycle;
		\draw[fill=\colorfour, opacity=\opazero] (2,0) -- (3,0) -- (2,1) -- (1,1) -- cycle;
		\draw (\tx,\ty) node [below right] {\picfontsize $\In_{\nu_2}(\Gamma_2 + \Gamma_3)$};    	
	\end{scope}
		
	\begin{scope}[shift={(5*\shifttwo,0)}]
		\draw [gray,  line width=0pt] (-0.5,-0.5) grid (\gridx,\gridy);
		\draw[ultra thick, fill=\colortwo, opacity=\opazero, dashed, ultra thick ] (0,0) -- (2,0) -- (1,1) -- (0,1) -- cycle;
		\draw (\tx,\ty) node [below right] {\picfontsize $\In_{\nu_2}(\Gamma_3 + \Gamma_1)$};    	
	\end{scope}
\end{scope}
\end{tikzpicture}
\caption{Normalized faces of the diagrams of \cref{ex0}}  \label{fig:ex0-2}
\end{center}
\end{figure}

\section{Characterization of minimal multiplicity systems} \label{non-degeneracy-0-section}
Given a collection $\scrA = (\scrA_1, \ldots, \scrA_m)$ of (possibly infinite) subsets of $\znzero$, we write $\scrL_0(\scrA_j)$ for the space of all power series in $(x_1, \ldots, x_n)$ supported at $\scrA_j$, $j = 1, \ldots, n$, and $\scrL_0(\scrA) :=  \prod_{j= 1}^n \scrL_0(\scrA_j)$. For the case that $m = n$, in this section we characterize the systems $(f_1, \ldots, f_n) \in \scrL_0(\scrA)$, which achieve the minimum possible intersection multiplicity at the origin. The proofs of the results of this section are given in \cref{mult-non-degeneracy-section,efficient-section,strictly-less-generic-0}.

\subsection{Non-degeneracy at the origin}
As in the case of Bernstein's theorem, we try to guess the correct non-degeneracy condition by considering the case $\kk = \cc$. Also assume for simplicity that each $\scrA_j$ is {\em finite}, so that every power series supported at $\scrA_j$ is in fact a {\em polynomial}. Now pick $f_j, g_j \in \kk[x_1, \ldots, x_n]$ over $\cc$ such that $\scrA_j = \supp(f_j) \supseteq \supp(g_j)$ for each $j$, and $\multfzero > \multgzero$. Write $h_j := (1-t)f_j + tg_j$, $j = 1, \ldots, n$. Then it seems reasonable to expect that there is a curve $C(t)$ on $\cc^n$ such that $h_j (C(t)) =0$ and $\lim_{t \to 0}C(t) = \origin$ (see \cref{fig:h-t-0}). Pick a parametrization $U \to C(t)$, where $U$ is a neighborhood of the origin on $\cc$ the form $\gamma: t \mapsto (\gamma_1(t), \ldots, \gamma_n(t))$, where each $\gamma_i$ is a power series in $t$. Fix $j = 1, \ldots, n$. As in \cref{sec:bkk-non-degeneracy-statement}, we examine the initial part of the expansion of $h_j(\gamma(t))$.

\subsubsection{Base case.} At first we consider the case that the image of $\gamma$ intersects $\ntorus$, i.e.\ no $\gamma_i$ is identically zero. Let $\gamma_i = a_i t^{\nu_i} + \cdots$, where $a_i \in \cc^*$, and $\nu_i := \ord_t(\gamma_i)$. Let $\nu$ be the weighted order on $\kk[x_1, \ldots, x_n]$ corresponding to weights $\nu_i$ for $x_i$, $i = 1, \ldots, n$. Then it follows exactly as in \cref{sec:bkk-non-degeneracy-statement} that $\In_{\scrA_j, \nu}(f_j)(a) = 0$ for each $j = 1, \ldots, n$, where $\In_{\scrA_j, \nu}(\cdot)$ are defined as in \eqref{insnu}. Since $\lim_{t \to 0} \gamma(t) = \origin$, it follows in addition that each $\nu_j$ is positive, i.e.\ $\nu$ is {\em centered at the origin}. This leads to the following notion.

\begin{defn} \label{defn:b-non-degeneracy-0}
Let $\scrA := (\scrA_1, \ldots, \scrA_m)$ be a collection of (possibly infinite) subsets of $\znzero$ and $(f_1, \ldots, f_m) \in \scrL_0(\scrA)$.  We say that $f_1, \ldots, f_m$ are {\em $(\scrA, *)$-non-degenerate at the origin} if they satisfy the following condition:
\begin{align}
\parbox{0.54\textwidth}{for each weighted order $\nu$ centered at the origin, there is no common root of $\In_{\scrA_j, \nu}(f_j)$, $j = 1, \ldots, m$, on $\nktorus$.}
\tag{$\text{ND}^*_0$}
\label{b-non-degeneracy-0}
\end{align}
We say that $f_1, \ldots, f_m$ are {\em $*$-non-degenerate at the origin} if they are $(\scrB, *)$-non-degenerate at the origin with $\scrB := (\supp(f_1), \ldots, \supp(f_m))$.
\end{defn}

\begin{figure}[h]
\begin{center}

\begin{tikzpicture}[scale=0.33]

\def\a{6}
\def\b{1}
\def\c{0}
\def\xshift{-6}
\def\tlabelx{0}
\def\tlabely{-0.3}
\def\nsamples{103}
\def\tmin{-3}
\def\tmax{3}

\def\xmin{-1.2}
\def\xmax{1}
\def\ymin{-1.3}
\def\ymax{1.3}
\def\initialshift{3}

\begin{scope}[shift={(0,0)}]
\begin{axis}[
xmin = \xmin, xmax=\xmax, ymin = \ymin, ymax= \ymax,
axis equal=true, axis equal image=true, hide axis
]
\addplot[blue, thick, domain=\tmin:\tmax, samples=\nsamples] ({x^2 - 1} ,{x^3 - x});
\addplot[red, thick, domain=\xmin:\xmax,samples=2](x,x);
\end{axis}
\draw (\tlabelx,\tlabely) node [below right] {\picfontsize $t = 0$};    	
\end{scope}

\begin{scope}[shift={(3*\xshift,0)}]
\begin{axis}[
xmin = \xmin, xmax=\xmax, ymin = \ymin, ymax= \ymax,
axis equal=true, axis equal image=true, hide axis
]
\addplot[blue, thick, domain=\tmin:\tmax, samples=\nsamples] ({x^2 - 1} ,{x^3 - x});
\addplot[red, thick, domain=\xmin:\xmax,samples=2](x,{(\b*x + \c)/\a});
\end{axis}
\draw (\tlabelx,\tlabely) node [below right] {\picfontsize $t = 1$};    	
\end{scope}

\begin{scope}[shift={( 2*\xshift,0)}]
\def\t{0.4}
\begin{axis}[
xmin = \xmin, xmax=\xmax, ymin = \ymin, ymax= \ymax,
axis equal=true, axis equal image=true, hide axis
]
\addplot[blue, thick, domain=\tmin:\tmax, samples=\nsamples] ({x^2 - 1} ,{x^3 - x});
\addplot[red, thick, domain=\xmin:\xmax,samples=2](x,{(x*((1-\t) + \t*\b) +\t*\c)/(1-\t + \t*\a)});
\end{axis}
\draw (\tlabelx,\tlabely) node [below right] {\picfontsize $t = \t$};    	
\end{scope}

\begin{scope}[shift={(\xshift,0)}]
\def\t{0.2}
\begin{axis}[
xmin = \xmin, xmax=\xmax, ymin = \ymin, ymax= \ymax,
axis equal=true, axis equal image=true, hide axis
]
\addplot[blue, thick, domain=\tmin:\tmax, samples=\nsamples] ({x^2 - 1} ,{x^3 - x});
\addplot[red, thick, domain=\xmin:\xmax,samples=2](x,{(x*((1-\t) + \t*\b) +\t*\c)/(1-\t + \t*\a)});
\end{axis}
\draw (\tlabelx,\tlabely) node [below right] {\picfontsize $t = \t$};    	
\end{scope}

\def\textshiftx{1}
\def\textshifty{5}
\def\polyshiftone{1}
\def\polyshifttwo{9}

\begin{scope}[shift={(-\xshift + \textshiftx,0)}]	
\node [
right,
text width= 5.4cm,align=justify] at (0,\textshifty) {
\picfontsize
$(f_1, f_2) = (y-x, y^2 - x^2 - x^3)$ \\
$(g_1, g_2) = (6y - x, y^2 - x^2 - x^3)$
};

	\begin{scope}[shift={(\polyshiftone,0)}]					
	\draw [gray,  line width=0pt] (-0.5,-0.5) grid (2.5,2.5);
	\draw [<->] (0,2.5)  |- (2.5,0);
	\draw[thick, blue ] (1,0) -- (0,1);
	\draw (-0.5, \tlabely) node [below right] {
	\picfontsize
	$\nd(f_1) = \nd(g_1)$ \newline
	};
	\end{scope}

	\begin{scope}[shift={(\polyshiftone+\polyshifttwo,0)}]					
	\draw [gray,  line width=0pt] (-0.5,-0.5) grid (2.5,2.5);
	\draw [<->] (0,2.5)  |- (2.5,0);
	\draw[thick, blue ] (2,0) -- (0,2);
	\draw (-0.5, \tlabely) node [below right] {
	\picfontsize
	$\nd(f_2) = \nd(g_2)$ \newline
	};
	\end{scope}
\end{scope}

\end{tikzpicture}
\caption{\mbox{A common (non-fixed) root of $(1-t)f_j + tg_j = 0$, $j = 1,2$, approaches the origin as $t \to 0$}}
\label{fig:h-t-0}
\end{center}

\end{figure}

\subsubsection{General case.}
So far we ignored the possibility that some of the $\gamma_j$ can be identically zero. This happens if $\gamma(t)$ belongs to a proper coordinate subspace of $\cc^n$. Incorporating this possibility and running the same arguments as in the first case leads to the following notion.

\begin{defn} \label{defn:non-degeneracy-0}
Let $\scrA := (\scrA_1, \ldots, \scrA_m)$ be a collection of (possibly infinite) subsets of $\znzero$ and $(f_1, \ldots, f_m) \in \scrL_0(\scrA)$. For each $I \subseteq [n]$ and each $j$, write $f_j|_{\Ki}$ for the power series obtained from $f_j$ by substituting $0$ for each $x_k$ such that $k \not\in I$, and write $\Ai := (\Ai_1, \ldots, \Ai_m) = (\scrA_1 \cap \ri, \ldots, \scrA_m \cap \ri)$. We say that $f_1, \ldots, f_m$ are \index{Non-degeneracy!centered at the origin}{\em $\scrA$-non-degenerate at the origin} if they satisfy the following condition:
\begin{align}
\parbox{0.45\textwidth}{$f_1|_{\Ki}, \ldots, f_m|_{\Ki}$ are $(\Ai , *)$-non-degenerate at the origin for each nonempty subset $I$ of $[n]$.}
\tag{ND$_0$}
\label{non-degeneracy-0}
\end{align}
We say that $f_1, \ldots, f_m$ are {\em non-degenerate at the origin} if they are $\scrB$-non-degenerate at the origin with $\scrB := (\supp(f_1), \ldots, \supp(f_m))$.
\end{defn}

The preceding discussion suggests that for $\kk = \cc$, $\scrA$-non-degeneracy at the origin is a sufficient condition for minimum intersection multiplicity at the origin. We will see that it is in fact necessary and sufficient for all (algebraically closed) $\kk$.

\subsection{The results} \label{sec:non-degenerate-results-at-origin}
Let $\scrA := (\scrA_1, \ldots, \scrA_n)$ be a collection of (possibly infinite) subsets of $\zzeroo{n}$. \Cref{non-degeneracy-0-thm} below states the necessary and sufficient condition for the minimality of $\multfzero$ for $(f_1, \ldots, f_n) \in \scrL_0(\scrA)$. It also states that $\multfzero$ is minimal for ``generic'' $(f_1, \ldots, f_n) \in \scrL_0(\scrA)$. We have to be careful about the notion of ``genericness'' though, since the spaces $\scrL_0(\scrA_j)$ and $\scrL_0(\scrA)$ are in general infinite dimensional vector spaces over $\kk$, and therefore they are {\em not} algebraic varieties. Let $\scrA' := (\scrA_1 \cap \nd(\scrA_1), \ldots, \scrA_n \cap \nd(\scrA_n))$. Then $\scrL_0(\scrA')$ is an algebraic variety isomorphic to $\kk^{\sum_j |\scrA_j \cap \nd(\scrA_j)|}$. Let $\pi: \scrL_0(\scrA) \to \scrL_0(\scrA')$ be the natural projection which ``ignores'' the coefficients corresponding to exponents not in $\nd(\scrA_j)$, $j = 1, \ldots, n$. Write $\scrM_0(\scrA)$ (respectively, $\scrM_0(\scrA')$) for the set of all $(f_1, \ldots, f_n)$ in $\scrL_0(\scrA)$ (respectively, $\scrL_0(\scrA')$) with the minimum possible value for $\multfzero$. We will show that $\scrM_0(\scrA')$ is a nonempty Zariski open (and therefore Zariski dense) subset of $\scrL_0(\scrA')$, and $\scrM_0(\scrA) = \pi^{-1}(\scrM'_0)$.

\begin{rem} \label{ind-remark}
An \index{Ind-variety}{\em ind-variety} over a field $k$ is a set X along with a chain of subsets $X_0 \subset X_1 \subset \cdots$ such that
\begin{enumerate}
\item $X = \bigcup_i X_i$,
\item Each $X_i$ is an algebraic variety over $k$, and
\item The inclusions $X_i \into X_{i+1}$ are closed embeddings of algebraic varieties.
\end{enumerate}
It is not hard to see, taking arbitrary sequences of finite subsets $\scrA_{j,0} \subset \scrA_{j,1} \subset \cdots $ such that $\scrA_j = \bigcup_i \scrA_{j,i}$, that $\scrL_0(\scrA)$ is an ind-variety. The notion of Zariski topology has a natural extension to the case of ind-varieties. \Cref{non-degeneracy-0-thm} implies in particular that $\scrM_0(\scrA)$ is a nonempty dense open subset of $\scrL_0(\scrA)$ in this topology.
\end{rem}

\begin{thm}[{\cite{toricstein}}] \label{non-degeneracy-0-thm}
$\scrM_0(\scrA) = \pi^{-1}(\scrM'_0)$ and $\scrM_0(\scrA')$ is a nonempty Zariski open subset of $\scrL_0(\scrA')$. If $\multAzero = \infty$, then $\scrM_0(\scrA) = \scrL_0(\scrA)$. Otherwise the following are equivalent:
\begin{enumerate}
\item \label{condition:M0}  $(f_1, \ldots, f_n) \in \scrM_0(\scrA)$
\item \label{condition:non-degenerate-0} $f_1, \ldots, f_n$ are $\scrA$-non-degenerate at the origin.
\end{enumerate}
\end{thm}

\Cref{non-degeneracy-0-thm} is proven in \cref{mult-non-degeneracy-section}. To check $\scrA$-non-degeneracy at the origin, one needs to check $(\Ai, *)$-non-degeneracy for {\em all} (nonempty) subsets $I$ of $[n]$. The following theorem, which we prove in \cref{efficient-section}, often significantly limits the number of test cases.

\begin{thm}[{\cite{toricstein}}] \label{non-degeneracy-0'-thm}
Let $\scrA := (\scrA_1, \ldots, \scrA_m)$ be a collection\footnote{Note that the number of subsets is $m$, which may be distinct from $n$.} of (possibly infinite) subsets of $\zzeroo{n}$. For each $I \subseteq [n]$, let $\emptyiA := \{j: \Ai_j = \emptyset\}$ and $\mscrIA := \{I \subseteq [n]:\  I \neq \emptyset,\ |\emptyiA| \geq n - |I|\}$. Then for $(f_1, \ldots, f_m) \in \scrL_0(\scrA)$ the following are equivalent:
\begin{enumerate}
\item $f_1, \ldots, f_m$ are $\scrA$-non-degenerate at the origin.
\item $f_1|_{\Ki}, \ldots, f_m|_{\Ki}$ are $(\Ai, *)$-non-degenerate at the origin for every $I \in \mscrIA$.
\end{enumerate}
\end{thm}

\begin{rem}\label{non-degeneracy-0'-remark}
If $m = n$ and $0 < \multAzero < \infty$, then (due to \cref{mult-support} below)
\begin{align*}
\mscrIA = \{I \subseteq [n]:\  I \neq \emptyset,\ |\emptyiA| = n- |I|\} = \{I \subseteq [n]:\  I \neq \emptyset,\ |\tiA| = |I|\}
\end{align*}
where $\tiA := \{j: \Ai_j \neq \emptyset\}$.
\end{rem}

Now we go back to the $m = n$ case, i.e.\ $\scrA := (\scrA_1, \ldots, \scrA_n)$ is a collection of subsets of $\zz^n$. Define $\mscrIA$ as in \cref{non-degeneracy-0'-thm}. Similar to the characterization of strict monotonicity of mixed volume in \cref{strictly-mixed-monotone}, we give (in \cref{strictly-monotone-0}) a combinatorial characterization of strict monotonicity of $\multAzero$. As a corollary in \cref{strictly-less-generic-0} we prove the following result, which says that in the same way as in the case of $\nktorus$ (\cref{alternate-bernstein}), $\scrA$-non-degeneracy of a system of power series at the origin is equivalent to a combinatorial condition plus non-degeneracy at the origin with respect to their supports.

\begin{cor} \label{alternate-non-degeneracy-0}
If $0 < \multAzero < \infty$, then the following are equivalent for $(f_1, \ldots, f_n) \in \scrL_0(\scrA)$:
\begin{enumerate}
\item \label{assn:non-deg-0} $f_1, \ldots, f_n$ are $\scrA$-non-degenerate at the origin.
\item \label{assn:non-deg-0-combinatorial}
\begin{enumerate}
\item  for each nonempty subset $I$ of $[n]$, and each $\nu \in \rnstar$ which is centered at the origin, the collection $\{\In_\nu(\nd(f_j) \cap \ri): \In_\nu(\Ai_j) \cap \supp(f_j) \neq \emptyset\}$ of polytopes is dependent, and
\item $f_1, \ldots, f_n$ are non-degenerate at the origin.
\end{enumerate}
\item \label{assn:non-deg-0-efftorial}
\begin{enumerate}
\item for each $I \in \mscrIA$ and each $\nu \in \rnstar$ which is centered at the origin, the collection $\{\In_\nu(\nd(f_j) \cap \ri): \In_\nu(\Ai_j) \cap \supp(f_j) \neq \emptyset\}$ of polytopes is dependent, and
\item $f_1, \ldots, f_n$ are non-degenerate at the origin.
\end{enumerate}
\end{enumerate}
\end{cor}

\section{Proof of the non-degeneracy condition} \label{mult-non-degeneracy-section}
In this section we prove \cref{non-degeneracy-0-thm}. Let $\scrA:= (\scrA_1, \ldots, \scrA_m)$, $m \geq 1$, be a collection of subsets of $\zzeroo{n}$. Let $\scrL_0(\scrA), \scrA'$ be as in \cref{sec:non-degenerate-results-at-origin}. Let $I \subseteq [n]$; define $\tiA := \{j: \Ai_j \neq \emptyset\}$ as in \cref{multiplicity-thm}. Note that $\tiA = \tiAprime$.

\begin{lemma}\label{isolated-prop}
Assume $0 \not\in \bigcup_j \scrA_j$. Then
\begin{enumerate}
\item If $|\tiA| < |I|$, then $\dim_\kk(\kk[[x_i : i \in I]]/\langle f_1|_{\Ki}, \ldots, f_m|_{\Ki} \rangle) = \infty$ for all $(f_1, \ldots, f_m) \in \scrL_0(\scrA)$.
\item If $|\tiA| \geq |I|$, then $V(f_1, \ldots, f_m) \cap \Kstari$ is isolated for generic $f_1, \ldots, f_m \in \scrL_0(\scrA')$.
\end{enumerate}
\end{lemma}

\begin{proof}
Due to \cref{finite-mult-determinacy} it suffices to prove the first assertion for $m$-tuple $(f_1, \ldots, f_m)$ of {\em polynomials} in $\scrL_0(\scrA)$. If $|\tiA| < |I|$ then the number of $f_j$ such that $f_j|_{\Ki}$ is nonzero is less than $|I|$. Since $0 \not\in \scrA_j$ for any $j$, each $f_j|_{\Ki}$ is in the maximal ideal of $R_I := \kk[[x_i : i \in I]]/\langle f_1|_{\Ki}, \ldots, f_m|_{\Ki} \rangle$. \Cref{thm:pure-dimension} implies that the transcendence degree of $R_I$ over $\kk$ is positive, so that $\dim_\kk(R_I) = \infty$. The second assertion follows from Bernstein's theorem.
\end{proof}

\begin{cor}[cf.\ {\cite[Lemma 2]{rojas-toric}, \cite[Proposition 5]{herrero}}]
\label{mult-support}
Assume $m = n$.
\begin{enumerate}
\item $\multAzero = 0$ if and only if $0 \in \bigcup_{i=1}^n \scrA_i$.
\item $\multAzero = \infty$ if and only if $0 \not\in \bigcup_{i=1}^n \scrA_i$ and there is $I \subseteq [n]$ such that $|\tiA| < |I|$. \qed
\end{enumerate}
\end{cor}

Let $\scrN_0(\scrA)$ be the set of all $(f_1, \ldots, f_m) \in \scrL_0(\scrA)$ such that $f_1, \ldots, f_m$ are $\scrA$-non-degenerate at the origin. Note that $\scrN_0(\scrA) = \pi^{-1}(\scrN_0(\scrA'))$, where $\pi: \scrL_0(\scrA) \to \scrL_0(\scrA')$ is the natural projection.

\begin{prop} \label{non-degeneracy-0-existence}
$\scrN_0(\scrA')$ is a Zariski open subset of $\scrL_0(\scrA')$. If either $0 \in \bigcup_{i=1}^m \scrA_i$ or $|\tiA| \geq |I|$ for all $I \subseteq [n]$, then $\scrN_0(\scrA')$ is nonempty.
\end{prop}

\begin{proof}
For each $m$-tuple $\scrB = (\scrB_1, \ldots, \scrB_m)$ of subsets of $\rr^n$ and for each $\nu \in \rnstar$, we write $\In_\nu(\scrB) := (\In_\nu(\scrB_1), \ldots, \In_\nu(\scrB_m))$. If $\scrB' = \In_\nu(\scrB)$ for some $\nu \in \rnstar$, we say that $\scrB'$ is a {\em face} of $\scrB$ and write that $\scrB' \preceq \scrB$; if in addition $\nu$ is centered at the origin, we write that $\scrB' \preceq_0 \scrB$.

\begin{proclaim} \label{transitive-claim}
If $\scrB' \preceq \scrB \preceq_0 \scrA'$, then $\scrB' \preceq_0 \scrA'$.
\end{proclaim}

\begin{proof}
By assumption there is $\nu \in \rnstar$ centered at the origin such that $\scrB = \In_\nu( \scrA')$. Pick $\nu' \in \rnstar$ such that $\scrB' = \In_{\nu'}(\scrB)$. If $k$ is a sufficiently large positive integer, then each of the coordinates of $k\nu + \nu'$ is {\em positive} with respect to the basis dual to the standard basis on $\rr^n$, and $\In_{k\nu + \nu'}(\scrA') = \scrB'$, so that $\scrB' \preceq_0 \scrA'$.
\end{proof}

Given $\scrB  \preceq \scrA'$, and $f = (f_1, \ldots, f_m) \in \scrL_0(\scrA')$, define $f_{j,\scrB_j}$ as in \cref{bkk-existential-section}. Let $\scrD_\scrB(\scrA')$ be the set of all $f \in \scrL_0(\scrA')$ such that there is a common root of $f_{1,\scrB_1}, \ldots, f_{m, \scrB_m}$ on $\nktorus$. Let $\scrD_0(\scrA') := \bigcup_{\scrB \preceq_0 \scrA'} \scrD_\scrB(\scrA')$. \Cref{transitive-claim} implies that $\scrD_0(\scrA') = \bigcup_{\scrB \preceq_0 \scrA'}\bigcup_{\scrB' \preceq \scrB} \scrD^0_{\scrB'}(\scrA')$, so that \cref{DJBbar-closed} implies that $\scrD_0(\scrA')$ is a Zariski closed subset of $\scrL_0(\scrA')$. Let $I \subseteq [n]$. Replacing $\scrA'$ by $\scrA'^I := (\scrA'_1 \cap \ri, \ldots, \scrA'_m \cap \ri)$, it follows that $\scrD_0(\scrA'^I)$ is a Zariski closed subset of $\scrL_0(\scrA'^I)$. Let $\bar \pi_{0,I}: \scrL_0(\scrA') \to \scrL_0(\scrA'^I)$ be the natural projection. Then $\scrN_0(\scrA') = \scrL_0(\scrA') \setminus \bigcup_{I \subseteq [n]}\bar \pi_{0,I}^{-1}(\scrD_0(\scrA'^I))$ is Zariski open in $\scrL_0(\scrA')$. It now remains to prove the second assertion of \cref{non-degeneracy-0-existence}. If $0 \in \scrA_i$ for some $i$, then any polynomial supported at $\scrA'_i$ with a nonzero constant term would lead to an element in $\scrN_0(\scrA')$. On the other hand, if $|\tiA| \geq |I|$ for every $I \subseteq [n]$, then \cref{DJBbar-proper} implies that $\scrD_0(\scrA'^I)$ is a proper Zariski closed subset of $\scrL_0(\scrA'^I)$ for every $I \subseteq [n]$, so that $\scrN_0(\scrA')$ is nonempty, as required.
\end{proof}

We now explore the relation between non-degeneracy at the origin and the intersection multiplicity at the origin. At first we need to extend the notion of weighted orders and ``initial coefficients'' corresponding to branches on $\nktorus$ to the case of branches on $\kk^n$.

\begin{defn} \label{IB-defn}
Let $B := (Z,z)$ be a branch of a curve $C \subset \kk^n$. Identify $Z^* := Z \setminus z$ with its image on $C$ and let $I_B :=  \{i: x_i|_{Z^*} \not\equiv 0\}$. Note that $\Kii{I_B}$ is the smallest coordinate subspace of $\kk^n$ which contains $Z^*$. We write $\nu_B$ for the weighted order on $\kk[x_i, x_i^{-1}: i \in I_B]$ corresponding to the weight $\ord_z(x_i|_Z)$ for each $i \in I_B$. Fix a parameter $\rho_B$ of $B$ and define
\begin{align}
\begin{split}
\In_B(x_j) &:=
	\begin{cases}
		0	&\text{if}\ j \not\in I_B \\
		\left. \frac{x_j}{(\rho_B)^{\nu_B(x_j)}} \right|_z
			&\text{if}\ j \in I_B.
	\end{cases}\\
\In(B) &:= (\In_B(x_1), \ldots, \In_B(x_n)) \in \Kstarii{I_B}
\end{split}
\end{align}
Compare this definition with the case of branches on $\nktorus$ defined in \eqref{in-B-torus} in \cref{toric-center-section}. The following result is immediate from the definition.
\end{defn}

\begin{lemma} \label{center-lemma}
If the center of $B$ is on $\Kstari$, then $I \subset I_B$.   \qed
\end{lemma}

\begin{lemma} \label{branch-lemma-IB}
Let $(f_1, \ldots, f_m) \in \scrL_0(\scrA) \cap \kk[x_1, \ldots, x_n]$ and $B$ be a branch of a curve contained in $V(f_1, \ldots, f_m)$. Then $\In(B) \in V(\Innub(f_1|_{\Kib}), \ldots, \Innub(f_m|_{\Kib})) \cap \Kstarib$. If in addition $B$ is a branch at the origin, then $f_1|_{\Kib}, \ldots, f_m|_{\Kib}$ violate condition \eqref{b-non-degeneracy-0} (from \cref{defn:b-non-degeneracy-0}) with $\nu = \nu_B$; in particular, $f_1, \ldots, f_m$ violate \eqref{non-degeneracy-0} (from \cref{defn:non-degeneracy-0}) with $I = I_B$.
\end{lemma}

\begin{proof}
The first assertion is a direct corollary of \cref{branch-lemma-1}. If $B$ is a branch at the origin, then $\nu_B$ is centered at the origin, so that the second assertion follows from the first one.
\end{proof}

\begin{cor} \label{finite-0-cor}
If $f_1, \ldots, f_m \in \scrL_0(\scrA) \cap \kk[x_1, \ldots, x_n]$ are $\scrA$-non-degenerate at the origin, then the origin can not be a non-isolated point of $V(f_1, \ldots, f_m) \subset \kk^n$. \qed
\end{cor}

\subsection{Proof of \cref{non-degeneracy-0-thm}}
Below sometimes we work with $\kk^{n+1}$ with coordinates $(x_1, \ldots, x_n, t)$. In those cases we usually denote the coordinates of elements of $\kk^{n+1}$ as pairs, with the last component of the pair denoting the $t$-coordinate. In particular, the origin of $\kk^{n+1}$ is denoted as $(0,0)$. Take $m = n$ and define $\scrM_0 (\scrA)$ as in \cref{sec:non-degenerate-results-at-origin}.


\begin{claim}  \label{sufficient-0-lemma}
$\scrM_0(\scrA) \supseteq \scrN_0(\scrA)$.
\end{claim}

\begin{proof}
Let $(f_1, \ldots, f_n) \in \scrL_0(\scrA) \setminus \scrM_0(\scrA)$. We will show that $f_1, \ldots, f_n$ are {\em $\scrA$-degenerate} at the origin. By our assumption there is $(g_1, \ldots, g_n)  \in \scrL_0(\scrA)$ such that $\multgzero < \multfzero$. Due to \cref{finite-mult-determinacy} we may assume all $g_i$ and $f_j$ are {\em polynomials} in $x_1, \ldots, x_n$. Let $t$ be a new indeterminate. An application of \cref{mult-deformation-0} with $h_j := (1-t)f_j + tg_j$, $j= 1, \ldots, n$, $X = \kk^n$, and $(b_0, \epsilon_0) = (0,1)$ implies that
\begin{prooflist}
\item \label{0-non-isolated} either the origin is a non-isolated zero of $f_1, \ldots, f_n$,
\item \label{0-extra-component} or there is an irreducible component $V$ of $V(h_1, \ldots, h_n)$ in $\kk^n \times \kk$ containing $(0,0)$ such that $V$ is different from $\{0\} \times \kk$.
\end{prooflist}
In case \ref{0-non-isolated} $f_1, \ldots, f_n$ are $\scrA$-degenerate at the origin (\cref{finite-0-cor}), so consider that we are in case \ref{0-extra-component}. Choose a branch $B$ at the origin of a curve contained in $V(h_1, \ldots, h_n) \subset \kk^{n+1}$ which is different from $\{0\} \times \kk$. Since $B \not\subset \{0\} \times \kk$, it follows that $I := I_B \cap [n] \neq \emptyset$. Let $\nu$ be the restriction of $\nu_B$ to $\kk[x_i: i \in I]$. Then it follows as in the proof of \cref{bkk-sufficiency} that for each $j = 1, \ldots, n$,

\begin{align}
\In_{\Ai_j, \nu}(f_j|_{\Ki})
	&=
	\begin{cases}
	\In_\nu(f_j|_{\Ki}) = \In_{\nu_B}(h_j|_{\Kib})& \text{if}\ \supp(f_j|_{\Ki}) \cap \In_\nu(\Ai_j) \neq \emptyset,\\
	0 & \text{otherwise.}
	\end{cases}
\label{in-fj-Ki}
\end{align}

\Cref{branch-lemma-IB} implies that $h_1|_{\Kib}, \ldots, h_n|_{\Kib}$ violate \eqref{b-non-degeneracy-0} with $\nu = \nu_B$, and therefore identity \eqref{in-fj-Ki} implies that $f_1, \ldots, f_n$ are $\scrA$-degenerate at the origin, as desired.
\end{proof}

\begin{claim}  \label{necessary-0-lemma}
Assume $\multAzero < \infty$. Then $\scrM_0(\scrA) \subseteq \scrN_0(\scrA)$.
\end{claim}

\begin{proof}
If $\multAzero = 0$, then $(f_1, \ldots, f_n) \in \scrM_0(\scrA)$  if and only if one of the $f_j$ has a nonzero constant term, which immediately implies that $f_1, \ldots, f_n$ are $\scrA$-non-degenerate. So assume $0 < \multAzero < \infty$. Pick $(f_1, \ldots, f_n) \in \scrL_0(\scrA) \setminus \scrN_0(\scrA)$. We will show that $(f_1, \ldots, f_n) \not\in \scrM_0(\scrA)$ following (the adapted version of) Bernstein's trick from \cref{bkk-necessary-section}. \Woutlog\ we may assume that $\multfzero < \infty$, and due to \cref{finite-mult-determinacy} we may assume in addition that each $f_j$ is a polynomial in $(x_1, \ldots, x_n)$.  Since $f_1, \ldots, f_n$ violate \eqref{non-degeneracy-0}, there is a nonempty subset $I$ of $[n]$ and a weighted order $\nu$ centered at the origin on $\kk[x_i: i \in I]$ such that $\In_{\Ai_1,\nu}(f_1|_{\Ki}), \ldots, \In_{\Ai_n, \nu}(f_n|_{\Ki})$ have a common zero $a = (a_1, \ldots, a_n) \in \Kstari$. Let $\tiA = \{j: \Ai_j \neq \emptyset\}$. Since $\multAzero < \infty$, \cref{mult-support,non-degeneracy-0-existence,finite-mult-determinacy} imply that there is a system $(g_1, \ldots, g_n) \in \scrL_0(\scrA)$ of {\em polynomials} in $(x_1, \ldots, x_n)$ such that $\In_{\Ai_j, \nu}(g_j|_{\Ki})(a) \neq 0$ for each $j \in \tiA$. Define $c(t):= (c_1(t), \ldots, c_n(t)) :\kk \to \Ki$ as follows:
\begin{align*}
c_i(t)
	:=
	\begin{cases}
	a_it^{\nu_i} &\text{if}\ i \in I,\\
	0	&\text{otherwise.}
	\end{cases}
\end{align*}
For each $j \in \tiA$, let $m_j := \min_{\Ai_j}(\nu)$. Define
\begin{align*}
h_j &:=
	\begin{cases}
	t^{-m_j}g_j(c(t))f_j - t^{-m_j}f_j(c(t))g_j &\text{if}\ j \in \tiA,\\
	f_j & \text{otherwise.}
	\end{cases}
\end{align*}
Note that each $h_j$ is a polynomial in $(x_1, \ldots, x_n,t)$. Since $\In_{\Ai_j, \nu}(g_j|_{\Ki})(a) \neq 0 = \In_{\Ai_j, \nu}(f_j|_{\Ki})(a)$ for each $j \in \tiA$, it follows as in \cref{modified-proof} that $h_j(x,0)$ is a nonzero constant multiple of $f_j$ for each $j$. By our assumption the origin is an isolated zero of $f_1, \ldots, f_n$. Since $h_1, \ldots, h_n$ vanish on the curve $\{(c(t),t): t \in \kk\}$, \cref{mult-deformation-0} implies that $\multfzero > \multzero{h_1(x,\epsilon)}{h_n(x,\epsilon)}$ for generic $\epsilon \in \kk$. Since $\supp(h_j(x,\epsilon)) \subset \scrA_j$ for each $\epsilon$, it follows that $\multfzero > \multAzero$, as required.
\end{proof}

\Cref{non-degeneracy-0-thm} now follows from \cref{mult-support,non-degeneracy-0-existence,sufficient-0-lemma,necessary-0-lemma}. \qed

\section{Proof of the bound} \label{mult-bound-section}
In this section we prove \cref{multiplicity-thm}. The computation of intersection multiplicity becomes easier if a generic system satisfies a property which is stronger than \eqref{non-degeneracy-0}; at first we prove that such systems exist. The proof of \cref{multiplicity-thm} is then given in \cref{multiplicity-proof-section}. We start with a notation: if $g$ is a polynomial in $(x_1, \ldots, x_n)$, $a = (a_1, \ldots, a_n) \in \kk^n$ and $I \subseteq [n]$, we write $g^{I}_a$ for the polynomial in $(x_i: i \in I)$ obtained from substituting $a_{i'}$ for $x_{i'}$ for each $i'$ {\em not} in $I$.

\begin{center}
\newcommand\conecolor{blue}
\newcommand\cylcolor{olive}
\newcommand\zxiscolor{magenta}
\newcommand\intersectcolor{magenta}
\tikzstyle{sdot} = [yellow, circle, minimum size = 3pt, inner sep = 0pt, fill]
\pgfmathsetmacro\scalefactor{1.5}
\pgfmathsetmacro\axiscalefactor{1}
\pgfmathsetmacro\samplenum{5}
\pgfmathsetmacro\sampley{45}
\pgfmathsetmacro\opacity{0.3}
\pgfmathsetmacro\fillopacity{0.3}

\begin{figure}[h]
\begin{tikzpicture}[scale = \scalefactor]

\pgfmathsetmacro\b{0}
\pgfmathsetmacro\c{1}
\pgfmathsetmacro\a{1}

\pgfmathsetmacro\conerpluslimit{1.5}
\pgfmathsetmacro\conerminuslimit{1.5}
\pgfmathsetmacro\cylzlowlim{-\conerminuslimit/\c + \b}
\pgfmathsetmacro\cylzuplim{\conerpluslimit/\c + \b}

\begin{scope}
\begin{axis}[
    axis lines = middle,
    axis equal,
    view/h = 45, view/v = 15,
    scale = \axiscalefactor,
    xticklabels = {,,},
    yticklabels = {,,},
    zticklabels = {,,},
]
\addplot3 [
    surf,
    color = \conecolor,
    opacity = \opacity,
    fill opacity = \fillopacity,
    faceted color = \conecolor,
    samples = \samplenum,
    samples y = \sampley,
    domain = 0 : \conerpluslimit,
    y domain = 0 : 360
] (
	{x*cos(y)},
    {x*(1-sin(y))/\c},
   	{x*sin(y)/\c + \b}
);

\addplot3 [
    surf,
    color = \conecolor,
    opacity = \opacity,
    fill opacity = \fillopacity,
    faceted color = \conecolor,
    samples = \samplenum,
    samples y = \sampley,
    domain = 0 : \conerpluslimit,
    y domain = 0 : 360
] (
	{x*cos(y)},
    {x*(-1-sin(y))/\c},
   	{x*sin(y)/\c + \b}
);

\addplot3 [
    surf,
    color = \cylcolor,
    opacity = \opacity,
    fill opacity = \fillopacity,
    faceted color = \cylcolor,
    samples = 2,
    samples y = \sampley,
    domain = \cylzlowlim : \cylzuplim,
    y domain = 0 : 360
] (
	{\a*cos(y)},
	{\a*(1 + sin(y))},
    {x}
);

\addplot3 [
    samples = 2,
    samples y = 0,
    domain = \cylzlowlim : \cylzuplim,
    \zxiscolor, dashed, thick
] (
	{0},
	{0},
    {x}
);

\addplot3 [
    samples = \sampley,
    samples y = 0,
    domain = 0 : 360,
    \intersectcolor, dashed, thick
] (
	{\a*sin(x)},
	{\a*(1 - cos(x))},
    {\a*((cos(x/2))^2/(\c^2) - (sin(x/2))^2) + \b}
);

\node[sdot] at (axis cs:0,0,{\a/(\c*\c) + \b}) {};

\end{axis}
\end{scope}

\end{tikzpicture}
\caption{Illustration of \cref{degenerate-containment}: $f := x^2 + (y-1)^2 - 1 = 0$ and $g := z^2 + x^2 - (y+z)^2 = 0$ intersect along the $z$-axis and an ellipse which intersect at $(0,0,1)$. $f|_{z=1} = x^2 - 2y + y^2$ and $g|_{z = 1} = x^2 - 2y - y^2$ are degenerate at the origin, since with $\nu = (1,2) \in \rnstar$, $\In_\nu(f|_{z=1}) = \In_\nu(g|_{z=1}) = x^2 - 2y$.} \label{fig:degenerate-specialization-0}
\end{figure}
\end{center}

\begin{lemma} \label{degenerate-containment}
Let $I \subseteq [n]$ and $f_1, \ldots, f_k \in \kk[x_1, \ldots, x_n]$ such that $f_j|_{\Ki} = 0$ for each $j =1, \ldots, k$. Assume $V(f_1, \ldots, f_k)$ has an irreducible component $V$ which is not contained in $\Ki$. Then $f^{[n]\setminus I}_{1,a}, \ldots, f^{[n]\setminus I}_{k,a}$ are {\em degenerate} at the origin for each $a \in \Kstari \cap V$ (see \cref{fig:degenerate-specialization-0}).
\end{lemma}

\begin{proof}
Let $a \in \Kstari \cap V$ and $B$ be a branch centered at $a$ of a curve contained in $V$ such that $B \not\subset \Ki$. Then $I_B \supsetneqq I$ (\cref{center-lemma}) and $\pi_I(\In(B)) = a$, where $\pi_I: \kk^n \to \Ki$ is the natural projection. Let $I' := I_B \setminus I $ and $\nu'$ be the restriction of $\nu_B$ to $\kk[x_{i'}: i' \in I']$. Since the center of $B$ is on $\Kstari$, it follows that $\nu_B(x_i) = 0$ for each $i \in I$. For each $j = 1, \ldots, k$, since $f_j|_{\Ki} \equiv 0$, it then follows that $\In_{\nu_B}(f_j|_{\Kib})(\In(B)) = \In_{\nu'}(f^{[n]\setminus I}_{j,a}|_{\Kii{I'}})(a')$, where $a' := \pi_{I'}(\In(B)) \in \Kstarii{I'}$. The result now follows from \cref{branch-lemma-IB}.
\end{proof}

Let $f_1, \ldots, f_m$ be polynomials in $(x_1, \ldots, x_n)$ and $I \subseteq [n]$. Define
\begin{align}
D^I_0(f_1, \ldots, f_m):= \{a \in \Kstari: f^{[n]\setminus I}_{1,a}, \ldots, f^{[n]\setminus I}_{m,a}\ \text{are degenerate at the origin}\}
\label{DI0}
\end{align}
The following result is immediate from \cref{thm:pure-dimension,degenerate-containment}.

\begin{cor} \label{strong-cor}
Let $J, J' \subseteq [m]$ and $V'$ be an irreducible component of $V(f_{j'}: j' \in J')$ such that $V' \not\subseteq \Ki$. Assume $\dim(D^I_0(f_{j'}: j' \in J') \cap V(f_j: j \in J)) < |I| - |J|$. Then $V'$ does {\em not} contain any irreducible component of $V(f_j: j \in J) \cap \Kstari$. \qed
\end{cor}

\begin{defn} \label{strong-defn}
Given a collection $\scrA = (\scrA_1, \ldots, \scrA_m)$, $m \geq 1$, of {\em finite} subsets of $\zzeroo{n}$ and $I \subseteq [n]$, define
 \begin{align*}
 \tiA &:= \{j: \scrA_j \cap \ri \neq \emptyset\} \subseteq [m] \\
 \ziA &:= \{J = \{j_1, \ldots, j_{n - |I|}\}\subseteq [m] \setminus \tiA: |J| = n - |I|,\ \multzero{\pi_{[n]\setminus I}(\scrA_{j_1})}{\pi_{n \setminus I}(\scrA_{j_{[n] - |I|}})} < \infty\} \\
 \Ai &:= (\scrA_1 \cap \ri, \ldots, \scrA_m \cap \ri)
 \end{align*}
 where $\pi_{[n] \setminus I}: \rr^n \to \rii{[n]\setminus I}$ is defined as in \eqref{pi_I}. Due to the finiteness of the $\scrA_j$ we can identify $\scrL_0(\scrA)$ with the collection $\scrL(\scrA)$ of {\em polynomials} $(f_1, \ldots, f_m)$ such that $\supp(f_j) \subseteq \scrA_j$ for each $j$. Given $(f_1, \ldots, f_m) \in \scrL(\scrA)$, we say that $f_1, \ldots, f_m$ are \index{Non-degeneracy!strong}{\em strongly $\scrA$-non-degenerate} if for all $I \subseteq [n]$,
\begin{defnlist}
\item \label{strongly-non-degenerate} $f_1|_{\Ki}, \ldots, f_m|_{\Ki}$ are {\em properly $\Ai$-non-degenerate} (see \cref{properly-non-degenerate-section}),
\item \label{recursively-non-degenerate} for all $J' \in \ziA$, $f^{[n]\setminus I}_{j',(1, \ldots, 1)}$, $j' \in J'$, are {\em non-degenerate at the origin} and $\np(f^{[n]\setminus I}_{j',(1, \ldots, 1)}) = \conv(\pi_{[n]\setminus I}(\scrA_{j'}))$ for each $j' \in J'$,
\item \label{intersectionally-non-degenerate} for all $J\subseteq \tiA$ and $J' \in \ziA$,
\begin{align}
\dim(D^I_0(f_{j'}: j' \in J') \cap V(f_j: j \in J)) < |I| - |J|,
\label{dim-degen-IJ}
\end{align}
where $D^I_0(\cdot)$ are as in \eqref{DI0}.
\end{defnlist}
Note that property \ref{recursively-non-degenerate} with $I = \emptyset$ in particular implies that
\begin{defnlist}[resume]
\item \label{strong-np} $\np(f_j) = \conv(\scrA_j)$ for each $j = 1, \ldots, m$.
\end{defnlist}
We write $\Nstrong(\scrA)$ be the collection of all $(f_1, \ldots, f_m) \in \scrL(\scrA)$ which are strongly $\scrA$-non-degenerate. Recall that $\scrN_0(\scrA)$ stands for the collection of systems which are $\scrA$-non-degenerate at the origin. For \cref{A-non-degeneracy,A-non-degeneracy-existence} below we assume each $\scrA_j$ is finite.
\end{defn}

\begin{prop} \label{A-non-degeneracy}
Assume either $0 \in \bigcup_{i=1}^m \scrA_i$ or $|\tiA| \geq |I|$ for all $I \subseteq [n]$. Then $\Nstrong(\scrA) \subseteq \scrN_0(\scrA)$. In particular, if $m = n$ and $\multAzero < \infty$, then $\Nstrong(\scrA) \subseteq \scrN_0(\scrA)$.
\end{prop}

\begin{proof}
Follows from \cref{properly-non-degenerate-claim,mult-support}.
\end{proof}

\begin{prop} \label{A-non-degeneracy-existence}
$\Nstrong(\scrA)$ is {\em constructible} and it contains a nonempty Zariski open subset of $\scrL(\scrA)$.
\end{prop}

\begin{proof}
By \cref{properly-non-degenerate-claim,non-degeneracy-0-existence} the collection of systems that satisfy properties \ref{strongly-non-degenerate} and \ref{recursively-non-degenerate} of \cref{strong-defn} is a nonempty Zariski open subset of $\scrL(\scrA)$. Therefore we can concentrate only on property \ref{intersectionally-non-degenerate}. For a subset $J \subseteq [m]$, write $\scrA_J := (\scrA_j: j \in J)$. Let $I \subseteq [n]$, $J\subseteq \tiA$ and $J' \in \ziA$. Write $\scrN^I(\scrA_J, \scrA_{J'})$ for the subset of $\scrL(\scrA_J) \times \scrL(\scrA_{J'})$ consisting of all $((f_j: j \in J), (f_{j'}: j' \in J'))$ which satisfy property \eqref{dim-degen-IJ}. Consider the set of maps from \cref{fig:strong-projections}, where
\begin{itemize}
\item $\pi_{J'}, \pi_{J,J'}, \pi_{J, I}, \pi_{J',I}$ are natural projections,
\item $\sigma$ is the ``substitution map'' which maps $((f_{j'}: j'\in J'), a) \in \scrL(\scrA_{J'}) \times \Kstari$ to $(f^{[n]\setminus I}_{j',a} : j' \in J')$.
\item $\scrD_0(\cdot)$ denotes the collection of systems which are {\em degenerate} at the origin, and
\item $\scrV_{J,I} := \{((f_j: j \in J), a) \in \scrL(\scrA_J) \times \Kstari: f_j(a) = 0$ for each $j \in J\}$
\end{itemize}

\begin{center}
\begin{figure}[h]
\begin{tikzcd}[column sep=tiny]
&
&
&
\scrL(\scrA_J) \times \scrL(\scrA_{J'}) \times \Kstari
\arrow[lld, "\pi_{J,J'}"]
\arrow[d, "\pi_{J,I}"]
\arrow[rd, "\pi_{J',I}"]
&
\\
\scrN^I(\scrA_J, \scrA_{J'})
\arrow[r, hook]
&
\scrL(\scrA_J) \times \scrL(\scrA_{J'})
\arrow[d, "\pi_{J'}"]
&
\scrV_{J,I}
\arrow[r, hook]
&
\scrL(\scrA_J) \times \Kstari
&
\scrL(\scrA_{J'}) \times \Kstari
\arrow[d, "\sigma = substitution"]
\\
&
\scrL(\scrA_{J'})
&
&
\scrD_0(\pi_{[n]\setminus I}(\scrA_{J'}))
\arrow[r, hook]
&
\scrL(\pi_{[n]\setminus I}(\scrA_{J'}))
\end{tikzcd}
\caption{
	Maps from the proof of \cref{A-non-degeneracy-existence}
} \label{fig:strong-projections}
\end{figure}
\end{center}
Let $\scrZ$ be the subset of $\scrL(\scrA_J) \times \scrL(\scrA_{J'}) \times \Kstari$ consisting of all $((f_j: j \in J), (f_{j'}: j' \in J'), a)$ such that $a \in V(f_j: j \in J)$ and $(f^{[n]\setminus I}_{j',a} : j' \in J')$ are $\pi_{[n]\setminus I}(\scrA_{J'})$-degenerate at the origin. Then $\scrZ = \pi_{J,I}^{-1}(\scrV_{J,I}) \cap (\sigma \circ \pi_{J',I})^{-1}(\scrD_0(\pi_{[n]\setminus I}(\scrA_{J'})))$. Since both $\scrV_{J,I}$ and $\scrD_0(\pi_{[n]\setminus I}(\scrA_{J'}))$ are Zariski closed (the closedness of $\scrD_0(\cdot)$ follows from \cref{non-degeneracy-0-existence}), it follows that $\scrZ$ is also Zariski closed. Since $\scrN^I(J,J')$ is the set of all elements in $\scrL(\scrA_J) \times \scrL(\scrA_{J'})$ whose pre-image under $\pi_{J,J'}|_\scrZ$ has dimension less than $|I| - |J|$, \cref{constructible-fiber-dimension} implies that $\scrN^I(J,J')$ is constructible. We now show that it contains a nonempty Zariski open subset of $\scrL(\scrA_J) \times \scrL(\scrA_{J'})$. Fix any $a_0 \in \Kstari$, and let $\sigma_0$ be the composition
\begin{align*}
\scrL(\scrA_{J'}) \into \scrL(\scrA_{J'}) \times \{a_0\} \xrightarrow{\sigma} \scrL(\pi_{[n]\setminus I}(\scrA_{J'}))
\end{align*}
where the first map simply takes $(f_{j'}: j' \in J') \mapsto ((f_{j'}: j; \in J'), a_0)$. Since $J' \in \ziA$, \cref{non-degeneracy-0-existence} implies that $\scrY' := \sigma_0^{-1}(\scrN_0(\pi_{[n] \setminus I} (\scrA_{J'}))$ is a nonempty Zariski open subset of $\scrL(\scrA_{J'})$. Pick an arbitrary system $(f_{j'}: j' \in J') \in \scrY'$, and let $\sigma'$ be the composition
\begin{align*}
\Kstari \into \{(f_{j'}: j' \in J')\} \times \Kstari \xrightarrow{\sigma} \scrL(\pi_{[n]\setminus I}(\scrA_{J'}))
\end{align*}
where the first map simply takes $a \mapsto ((f_{j'}: j' \in J'), a)$. Since $\sigma'(a_0) = \sigma_0(f_{j'}: j' \in J') \in \scrN_0(\pi_{[n] \setminus I} (\scrA_{J'}))$, \cref{non-degeneracy-0-existence} implies that $D^I_0(f_{j'}: j' \in J') = (\sigma')^{-1}(\scrD_0(\pi_{[n] \setminus I} (\scrA_{J'}))$ is a proper Zariski closed subset of $\Kstari$; in particular, $\dim(D^I_0(f_{j'}: j' \in J')) < |I|$. Since $J \subset \tiA$, \cref{intersection-lemma} below implies that there is a nonempty open subset $\scrW$ of $\scrL(\scrA_J)$ such that $((f_j: j \in J), (f_{j'}: j' \in J')) \in \scrN^I(J,J')$ for each $(f_j: j \in J) \in \scrW$. Since $(f_{j'}: j' \in J')$ was an arbitrary element from $Y'$, \cref{exercise:dimstructibly-open} implies that $\scrN^I(J,J')$ contains a nonempty Zariski open subset of $\scrL(\scrA_J) \times \scrL(\scrA_{J'})$, as required.
\end{proof}

\begin{lemma} \label{intersection-lemma}
Let $W$ be an irreducible subvariety $\nktorus$ and $\scrB = (\scrB_1, \ldots, \scrB_k)$ be a collection of finite nonempty subsets of $\zz^n$. Let $\scrW := \{(f_1, \ldots, f_k) \in \scrL(\scrB): \dim(V(f_1, \ldots, f_k) \cap W) \leq \dim(W) - k\}$. Then $\scrW$ is a constructible subset of $\scrL(\scrB)$ and it contains a nonempty Zariski open subset of $\scrL(\scrB)$.
\end{lemma}

\begin{proof}
Let  $\scrW' := \{((f_1, \ldots, f_k), (x_1, \ldots, x_n)) \in \scrL(\scrB) \times W: f_j(x_1, \ldots, x_n) = 0$ for each $j\}$ and $\pi_W: \scrW' \to W$ be the natural projection. For each $w \in W$, $\pi_W^{-1}(w)$ is a linear subspace of $\scrL(\scrB)$ defined by $k$ linearly independent linear equations, so that $\dim(\pi_W^{-1}(w)) = \dim(\scrL(\scrB) - k =  \sum_j |\scrB_j| - k$. \Cref{fiber-dimension} then implies that $\dim(\scrW') = \sum_j |\scrB_j|- k + \dim(W) = \dim(\scrL(\scrB)) + (\dim(W) - k)$. Now the result follows from applying \cref{fiber-dimension,constructible-fiber-dimension} to $\pi_{\scrB}|_{\scrW'} : \scrW' \to \scrL(\scrB)$, where $\pi_\scrB: \scrL(\scrB) \times \nktorus \to \scrL(\scrB)$ is the natural projection.
\end{proof}

\newcommand{\Bzeroi}{\scrB^I_\origin}
\newcommand{\Bzeroinu}{\scrB^I_{\origin, \nu}}
\newcommand{\Bzeroijnu}{\scrB^I_{\origin ,j, \nu}}
\newcommand{\Vzeroi}{\scrV^I_\origin}
\newcommand{\Vzeroii}[1]{\scrV^{#1}_\origin}
\newcommand{\Vzeroiprime}{\scrV'^I_\origin}
\newcommand{\Zzero}{\scrZ_\origin}
\newcommand{\Zzeroi}{\scrZ_\origin^I}

\subsection{Proof of \cref{multiplicity-thm}} \label{multiplicity-proof-section}
\Cref{mult-support} implies that \cref{multiplicity-thm} holds when $\multAzero = 0$ or $\infty$. So assume $0 < \multAzero < \infty$. Let $\scrA' := (\scrA_1 \cap \Gamma_1, \ldots, \scrA_n \cap \Gamma_n)$. Pick strongly $\scrA'$-non-degenerate $(f_1, \ldots, f_n) \in \scrL_0(\scrA')$. \Cref{non-degeneracy-0-thm,A-non-degeneracy} imply that
\begin{align*}
\multfzero = \multAzero 
\end{align*}
Therefore it suffices to show that
\begin{align}
\multfzero
	= \sum_{I \in \tAone}
							\multzerostar{\Gamma^{I}_1, \Gamma^{I}_{j_2}}{\Gamma^{I}_{j_{|I|}}}
							\times
							\multzero{\pi_{[n]\setminus I}(\Gamma_{j'_1})}{\pi_{[n]\setminus I}(\Gamma_{j'_{n-|I|}})}
	\label{mult-formula'}
\end{align}
where for each $I \in \tAone$, $j_1 = 1, j_2, \ldots, j_{|I|}$ are elements of $\tiA$, and $j'_1, \ldots, j'_{n-|I|}$ are elements of $[n]\setminus \tiA$. We proceed by induction on $n$. It is true for $n = 1$ (see \cref{formula-convention}), so assume it is true for all dimensions smaller than $n$. Since $0 < \multfzero < \infty$, \cref{int-mult-curve} implies that on a sufficiently small Zariski open neighborhood $U$ of  the origin in $\kk^n$, the closed subscheme of $U$ defined by $f_2, \ldots, f_n$ is a possibly non-reduced curve $C$. For each $I \subseteq [n]$, let $\{C^I_j\}_j$ be the set of irreducible components of $C$ such that $\Ki$ is the smallest coordinate subspace of $\kk^n$ containing each $C^I_j$.

\begin{claim} \label{I-A-1-claim}
Let $I \subseteq [n]$, $\tAone$ be as in \cref{multiplicity-thm} and $\ziA$ be as in \cref{strong-defn}.
\begin{enumerate}
\item If $\{C^I_j\}_j$ is nonempty, then $I \in \tAone$.
\item If $I \in \tAone$, then $[n] \setminus \tiA \in \ziA$.
\end{enumerate}
\end{claim}

\begin{proof}
For the first assertion, pick $I \subseteq [n]$ such that $\{C^I_j\}_j$ is nonempty. Since $0 <\multAzero < \infty$, \cref{mult-support} implies that $|\tiA| \geq |I|$. On the other hand, if $|\tiA \setminus \{1\}| \geq |I|$, then the proper non-degeneracy of $f_1|_{\Ki}, \ldots, f_n|_{\Ki}$ (property \ref{strongly-non-degenerate} of strong $\scrA$-non-degeneracy) and \cref{branch-lemma-IB} implies that $\{C^I_J\}_j$ is empty, which is a contradiction. Accordingly $|\tiA \setminus \{1\}| = |I| - 1$ and $|\tiA| = |I|$, which imply that $I \in \tAone$, as required. For the second assertion, pick $I \in \tAone$ and set $J' := [n] \setminus \tiA$. Since $|J'| = n - |I|$, we have to show that $\multzero{\pi_{[n]\setminus I}(\scrA_{j'_1})}{\pi_{n \setminus I}(\scrA_{j'_{[n] - |I|}})} < \infty$, where $j'_1, \ldots, j'_{n-|I|}$ are elements of $J'$. Indeed, otherwise \cref{mult-support} would imply that $|\tiiAA{I'}{\scrA}| < |I'|$ for some $I' \supsetneq I$, which would in turn imply (since by assumption $0 \not\in \bigcup_j \scrA_j$) that $\multAzero = \infty$, which is a contradiction.
\end{proof}

Pick $I \in \tAone$. Let $j_1 = 1, j_2, \ldots, j_{|I|}$ be the elements of $\tiA$ and $j'_1, \ldots, j'_{n- |I|}$ be the elements of $J' := [n] \setminus \tiA$. Since $J' \in \ziA$ (\cref{I-A-1-claim}), property \ref{recursively-non-degenerate} of strong $\scrA$-non-degeneracy and \cref{degenerate-containment} imply that $\Ki$ is an irreducible component of $V(f_{j'_1}, \ldots, f_{j'_{n-|I|}})$. On the other hand, applying property \ref{intersectionally-non-degenerate} of strong $\scrA$-non-degeneracy with $J = \{j_2, \ldots, j_{|I|}\}$, and then using \cref{strong-cor} shows that no irreducible component of $V(f_{j'_1}, \ldots, f_{j'_{n-|I|}})$ other than $\Ki$ contains any irreducible component of $V(f_{j_2}, \ldots, f_{j_{|I|}}) \cap \Kstari$. Therefore \cref{I-A-1-claim,order-curve,int-mult-curve,mult-chain-0} imply that
\begin{align}
\multfzero
	= \sum_{I \in \tAone} \multzero{f^{[n]\setminus I}_{j'_1, \epsilon}}{f^{[n]\setminus I}_{j'_{n-|I|}, \epsilon}}
	\sum_j  \ord_\origin(f_1|_{C^I_j}) \multp{f_{j_2}|_{\Kstari}}{f_{j_{|I|}}|_{\Kstari}}{C^I_j}
	\label{mult-formula'''}
\end{align}
where $\epsilon$ is a generic element of $\nktorus$ and $f^{[n]\setminus I}_{\cdot, \epsilon}$ are as in \cref{degenerate-containment}. Since $J' \in \ziA$, property \ref{recursively-non-degenerate} of strong $\scrA$-non-degeneracy and \cref{non-degeneracy-0-thm} imply that
\begin{align*}
\multzero{f^{[n]\setminus I}_{j'_1, \epsilon}}{f^{[n]\setminus I}_{j'_{n-|I|}, \epsilon}}
	&= \multzero{\pi_{[n]\setminus I}(\scrA_{j'_1})}{\pi_{[n]\setminus I}(\scrA_{j'_{n-|I|}})}
	= \multzero{\pi_{[n]\setminus I}(\Gamma_{j'_1})}{\pi_{[n]\setminus I}(\Gamma_{j'_{n-|I|}})}
\end{align*}
It remains to compute the inner sum of the right hand side of \eqref{mult-formula'''}. Let $I \subseteq [n]$. Write $R^I := \kk[x_i : i \in I]$. Let $\Vzeroi$ be the set of weighted orders on $R^I$ which are centered at the origin and $\Vzeroiprime$ be the set of primitive elements in $\Vzeroi$. For each $\nu \in \Vzeroiprime$, let $\Bzeroijnu$ be the set of all branches at the origin of $C^I_j$ such that $\nu_B$ is proportional to $\nu$. \Cref{order-curve} implies that for each $I,j$,
\begin{align*}
\ord_\origin(f_1|_{C^I_j}) = \sum_{\nu \in \Vzeroiprime} \sum_{(Z,z) \in \Bzeroijnu} \ord_z(f_1|_{C^I_j})
\end{align*}
Therefore it suffices to show that
\begin{align}
\sum_j \sum_{(Z,z) \in \Bzeroijnu} \ord_z(f_1|_{C^I_j}) \multp{f_{j_2}|_{\Kstari}}{f_{j_{|I|}}|_{\Kstari}}{C^I_j}
	&=  \min_{\Gamma^I_1}(\nu)
			\mv'_\nu(\In_\nu(\Gamma^I_{j_2}), \ldots, \In_\nu(\Gamma^I_{j_{|I|}}))
			\label{order-0-sum}
\end{align}
To see it, apply \cref{order-nu-cor} (with $n = |I|$) to $f_1|_{\Kstari}, f_{j_2}|_{\Kstari}, \ldots, f_{j_{I}}|_{\Kstari}$. Property \ref{strongly-non-degenerate} of strong $\scrA$-non-degeneracy implies that all the assumptions of \cref{order-nu-prop,order-nu-cor} are satisfied. Part \ref{nonempty-implication} of \cref{order-nu-prop} implies that each irreducible component of the resulting curve $C'$ comes from an irreducible component of $V(f_{j_2}|_{\Kstari}, \ldots, f_{j_{I}}|_{\Kstari}) \subset \Kstari$ and therefore the collections $\scrB'_{j,\nu}$ from \cref{order-nu-cor} are precisely the collections $\Bzeroijnu$. \Cref{order-nu-cor} then implies identity \eqref{order-0-sum} and completes the proof of \cref{multiplicity-thm}.

\section{The efficient version of the non-degeneracy condition} \label{efficient-section}
In this section we prove \cref{non-degeneracy-0'-thm}. Given $I \subseteq [n]$, we write $\Vi$ for the set of weighted orders on $\kk[x_i :i \in I]$. Given $I \subseteq \tilde I$, we say that $\nu \in \Vi$ and $\tilde \nu \in \Vitilde$ are \index{Compatible weighted order}\index{Weighted!order!compatible}{\em compatible} if $(\nu(x_i): i \in I)$ and $(\tilde \nu(x_i): i \in I)$ are proportional, with a {\em positive} constant of proportionality, and $\tilde \nu(x_{\tilde i}) > 0$ for each $\tilde i \in \tilde I \setminus I$. \Cref{non-degeneracy-0'-thm} follows directly from \cref{reduction-lemma} below.

\begin{lemma} \label{reduction-lemma}
Let $\scrA := (\scrA_1, \ldots, \scrA_m)$ be a collection of (possibly infinite) subsets of $\zzeroo{n}$, $I$ be a nonempty subset of $[n]$, and $\nu \in \Vi$ be such that
\begin{align}
\parbox{0.63\textwidth}{
$\In_{\tilde \nu}(\scrA_j)$ is finite for each $\tilde \nu \in \rnstar$ which is compatible with $\nu$.
} \label{initially-finite}
\end{align}
Let $f_1, \ldots, f_m \in \kk[[x_1, \ldots, x_n]]$ such that $\supp(f_j) \subseteq \scrA_j$ for each $j$. Assume
\begin{enumerate}
\item \label{less-than-assumtion} $|\emptyiA| <  n - |I|$, where $\emptyiA := \{j: \Ai_j = \emptyset\}$, and
\item \label{initial-assumption} $\In_{\Ai_j, \nu}(f_j|_{\Ki})$, $j = 1, \ldots, m$, have a common zero $u \in \nktorus$.
\end{enumerate}
Then there exists $\tilde I \supsetneqq I$ and $\tilde \nu \in \Vitilde$ such that
\begin{enumerate}
\setcounter{enumi}{2}
\item \label{compatible-reduction} $\tilde \nu$ is compatible with $\nu$.
\item \label{zero-reduction} $\In_{\Atildei_j, \tilde \nu}(f_j|_{\Kitilde})$, $j = 1, \ldots, m$, have a common zero $\tilde u \in \nktorus$ such that $\pi_I(\tilde u) = \pi_I(u)$, where $\pi_I:\nktorus \to \Kstari$ is defined as in \eqref{pi_I}.
\end{enumerate}
\end{lemma}

\begin{proof}
Due to \eqref{initially-finite} we may assume without any loss of generality that the support of each $f_j$ is finite, i.e.\ the $f_j$ are {\em polynomials} in $(x_1, \ldots, x_n)$. We may also assume that $I = \{1, \ldots, k\}$, $1 \leq k \leq n$. Let $a := \pi_I(u) \in \Kstari$ and $(a_1, \ldots, a_n)$ be the coordinates of $a$. At first consider the case that $\nu(x_i) = 0$ for each $i \in I$. Assumption \eqref{initial-assumption} then says that $a$ is a common zero of $f_1, \ldots, f_m$ on $\Kstari$. Let $y_j := x_j - a_j$, $j = 1, \ldots, n$, so that $a$ is the origin of $\kk^n$ with respect to $(y_1, \ldots, y_n)$ coordinates. Choose any integral $\nu' \in \rnstar$ with positive coordinates with respect to the dual basis, and let $\pi: \blnuprimekn \to \kk^n$ be the $\nu'$-weighted blow up of $\kk^n$ with respect to $(y_1, \ldots, y_n)$ coordinates (see \cref{weighted-blow-up-section}). Let $\Enuprime$ be the exceptional divisor of $\pi$, and $W'$ be the strict transform of $\Kstari$ on $\blnuprimekn$. Since $|\emptyiA| < n-|I|$, there is an irreducible component $V$ of $V(f_j: j \in \emptyiA) \cap \nktorus$ properly containing $\Kstari$. Then the strict transform $V'$ of (the closure of) $V$ properly contains $W'$. Pick $a' \in \Enuprime \cap W'$, and choose an irreducible curve $C' \subset V'$ such that $a' \in C' \not\subseteq W' \cup \Enuprime$, and a branch $B' = (Z',z')$ of $C'$ centered at $a'$. Let $\tilde I := I_{B'}$ and $\tilde \nu := \nu_{B'} \in \Vitilde$ (\cref{IB-defn}). Since $\pi$ is centered at $a \in \Kstari$ and since $\pi(B') \not\subset \Ki$, it follows that $I \subsetneqq \tilde I$, $\In_{B'}(x_i) = a_i$ and $\tilde \nu(x_i) = 0$ for each $i \in I$, and $\tilde \nu(x_{\tilde i})$ is positive for each $\tilde i \in \tilde I \setminus I$. Fix $j \in [m]$. If $j \not\in \emptyiA$, it follows that $\min_{\Atildei_j}(\tilde \nu) = 0$ and $\In_{\Atildei_j,\tilde \nu}(f_j|_{\Kitilde}) = f_j|_{\Ki} = \In_{\Ai_j, \nu}(f_j|_{\Ki})$. This implies that $\In_{\Atildei_j, \tilde \nu}(f_j|_{\Kitilde})(\In(B')) = f_j|_{\Ki}(a) = 0$. On the other hand, if $j \in \emptyiA$, then $\In_{\Atildei_j, \tilde \nu}(f_j|_{\Kitilde})(\In(B')) = 0$ due to \cref{branch-lemma-IB}. The lemma is therefore true in the case that $\nu$ is the trivial weighted order. \\

Now assume $\nu$ is not the trivial weighted order. Identify $\nu$ with the element in $\rnstar$ with coordinates $(\nu(x_1), \ldots, \nu(x_k))$ with respect to the basis dual to the standard basis of $\rr^n$. Choose a basis $\alpha_1, \ldots, \alpha_k$ of $\zz^k$ such that $\langle \nu, \alpha_j \rangle = 0$ for $j = 1, \ldots, k-1$, and $\langle \nu, \alpha_k \rangle = 1$. Then $(x^{\alpha_1}, \ldots, x^{\alpha_k}, x_{k+1}, \ldots, x_n)$ are coordinates on $X := \Kstarii{k} \times \kk^{n-k}$. Define
\begin{align*}
y_j
	&:=
	\begin{cases}
		x^{\alpha_j} - a^{\alpha_j} & \text{if}\ 1 \leq j \leq k-1,\\
		x^{\alpha_k} & \text{if}\ j = k,\\
		x_j	& \text{if}\ k+1 \leq j \leq n.
	\end{cases}
\end{align*}
Write $Y$ for the affine space $\kk^n$ with coordinates $(y_1, \ldots, y_n)$. Choose positive integers $\nu'_1, \ldots, \nu'_n$ such that $\nu'_k = 1$ and $\nu'_j \gg 1$ for $j = k+1, \ldots, n$. Let $\nu'$ be the element in $\rnstar$ with coordinates $(\nu'_1, \ldots, \nu'_n)$ with respect to the basis dual to the standard basis, and $\pi: Y' \to Y$ be the $\nu'$-weighted blow up of $Y$ with respect to $(y_1, \ldots, y_n)$ coordinates, $E$ be the exceptional divisor of $\pi$, and $W'$ be the strict transform on $Y'$ of $W := V(y_{k+1}, \ldots, y_n) \subset Y$. \Cref{blow-up-coordinates} implies that there is an affine open subset $U$ of $Y'$ such that
\begin{prooflist}
\item \label{U-coordinates} $U \cong \kk \times \nktoruss{k-1} \times \kk^{n-k}$ with respect to coordinates $(z_1, \ldots, z_n)$ where $z_1, \ldots, z_k$ are monomials in $(y_1, \ldots, y_k)$, $z_j = y_j/z_1^{\nu'_j}$ for $j = k+1, \ldots, n$, $\nu'(z_1) = 1$ and $\nu'(z_j) =0$ for $j = 2, \ldots, n$,
\item \label{U-E} $U \cap E = V(z_1) \cong \nktoruss{k-1} \times \kk^{n-k}$, and
\item \label{U-W'} $U \cap W' = V(z_{k+1}, \ldots, z_n) \cong  \kk \times \nktoruss{k-1}$.
\end{prooflist}
We treat $X$ as an open subset of $Y$ via the natural embedding. There is an irreducible component $V$ of $V(f_j: j \in \emptyiA) \cap X$ such that its closure $\bar V$ in $Y$ properly contains $W$. The strict transform $V'$ of $\bar V$ on $Y'$ properly contains $W'$. Pick $a' \in U \cap W' \cap E$. Choose an irreducible curve $C' \subset V'$ such that $a' \in C'$, and $C' \not\subset E \cup W'$, and $C' \cap \pi^{-1}(X) \neq \emptyset$. Pick a branch $B' = (Z',z')$ of $C'$ centered at $a'$. Since $\pi(B') \cap X \neq \emptyset$, we may treat $B'$ as a branch (possibly at infinity) of a curve on $X$. Define $\nu_{B'}$ and $I_{B'}$ as in \cref{IB-defn}. Since each of $x_1, \ldots, x_k$ is everywhere nonzero on $X$, it follows that $I_{B'} \supset  \{1, \ldots, k\} = I$. On the other hand, since $\pi(B') \not\subset W$, it follows that there exists $j > k$ such that $x_j|_{B'} \not\equiv 0$. It follows that $I_{B'} \supsetneq I$. We show that properties \eqref{compatible-reduction} and \eqref{zero-reduction} are true with $\tilde I:= I_{B'}$ and $\tilde \nu := \nu_{B'}$. Indeed, since $a' \in E$, for each $j = 1, \ldots, n$, either $y_j|_{B'} \equiv 0$, or $\ord_{z'}(y_j|_{B'}) > 0$. Therefore, for each $j = 1, \ldots, k-1$,
\begin{prooflist}[resume]
\item \label{nu'-alpha-<k} $\nu_{B'}(x^{\alpha_j}) = \ord_{z'}((a^{\alpha_j} + y_j)|_{B'}) = 0$, since $\ord_{z'}(y_j|_{B'}) > 0$.
\end{prooflist}
Since $\nu_{B'}(x^{\alpha_k}) = \ord_{z'}(y_k|_{B'}) > 0$, it follows that $\nu_{B'}$ and $\nu$ are proportional on $\kk[x_i : i \in I]$ with a positive constant of proportionality. Pick $j \in I_{B'} \setminus I$. Then $j > k$, and $\nu_{B'}(x_j) = \ord_{z'}(y_j|_{B'}) > 0$. It follows that $\nu_{B'}$ and $\nu$ are compatible. It remains to exhibit property \eqref{zero-reduction}. Since the center of $B$ is on $U \cap E$, properties \ref{U-coordinates} and \ref{U-E} of $U$ imply that $(\ord_{z'}(y_1|_{B'}), \ldots, \ord_{z'}(y_k|_{B'}))$ is proportional to $(\nu'_1, \ldots, \nu'_k)$. Since $\nu'(z_1) = 1$, it follows that the constant of proportionality is $q :=\ord_{z'}(z_1|_{B'})$. Therefore $\nu_{B'}(x^{\alpha_k}) = \ord_{z'}(y_k|_{B'}) = \nu'_k q = q = q\nu(x^{\alpha_k})$. Since $\nu_{B'}$ is compatible with $\nu$, it follows that
\begin{prooflist}[resume]
\item \label{nu'-1-k}  $\nu_{B'}(x_j) = q\nu(x_j)$ for $j = 1, \ldots, k$.
\end{prooflist}
On the other hand, since $a' \in U \cap W' \cap E$, properties \ref{U-E} and \ref{U-W'} imply that $\ord_{z'}(z_j|_{B'}) \geq 1$ if $j > k$. It follows that
\begin{prooflist}[resume]
\item \label{nu'->k}   for each $j \in I_{B'} \setminus I$,  $\nu_{B'}(x_j) = \ord_{z'}(z_j|_{B'}) + \nu'_j \ord_{z'} (z_1|_{B'})  >  q\nu'_j$.
\end{prooflist}
Let $u' := \In(B') \in \kstarii{I_{B'}}$. Observation \ref{nu'-alpha-<k} implies that $u'^{\alpha_j} = a^{\alpha_j}$ for $j = 1, \ldots, k-1$. \Cref{torus-image} then implies that there is $t \in \kk^*$ such that $(a_1, \ldots, a_k) = (t^{\nu(x_1)}u'_1, \ldots, t^{\nu(x_k)}u'_k)$. Choose a $q$-th root $t'$ of $t$ in $\kk$ and let $\tilde u = (\tilde u_1, \ldots, \tilde u_n)$ be an element with coordinates
\begin{align*}
\tilde u_j
    & :=
        \begin{cases}
            t'^{\nu_{B'}(x_j)}u'_j & \text{if}\ j \in I_{B'}, \\
            \text{arbitrary element in}\ \kk^* & \text{otherwise}.
        \end{cases} \\
    &=
        \begin{cases}
            t^{\nu(x_j)}u'_j = a_j & \text{if}\ j \in I, \\
            t'^{\nu_{B'}(x_j)}u'_j & \text{if}\ j \in I_{B'} \setminus I, \\
            \text{arbitrary element in}\ \kk^* & \text{otherwise}.
        \end{cases}
\end{align*}
Note that $\pi_I(\tilde u) = a = \pi_I(u)$. Fix $j \in [m]$. If $j \not\in \emptyiA$, then \ref{nu'-1-k} and \ref{nu'->k} imply that choosing $\nu'_{k+1}, \ldots, \nu'_n$ sufficiently large we can ensure that $\In_{\Aibprime_j, \nu_{B'}}(f_j|_{\kk^{I_{B'}}}) = \In_{\Ai_j, \nu}(f_j|_{\Ki})$, which would imply that $\In_{\Aibprime_j, \nu_{B'}}(f_j|_{\kk^{I_{B'}}})(\tilde u) = \In_{\Ai_j, \nu}(f_j|_{\Ki})(a) = 0$. On the other hand, if $j \in \emptyiA$, then 
\begin{align*}
\In_{\Aibprime_j, \nu_{B'}}(f_j|_{\kk^{I_{B'}}})(\tilde u) 
	&= t'^{\min_{\Aibprime_j}(\nu_{B'})}\In_{\Aibprime_j, \nu_{B'}}(f_j|_{\kk^{I_{B'}}})(\In(B')) 
	= 0
\end{align*}
due to \cref{branch-lemma-IB}. This completes the proof of property \eqref{zero-reduction}.
\end{proof}

\section{Other formulae for generic intersection multiplicity} \label{other-0-section}

\subsection{The formula of Huber-Sturmfels and Rojas}
Let $t$ be a new indeterminate. Fix positive integers $k_1, \ldots, k_n$. Note that for each $f_1, \ldots, f_n \in \kk[[x_1, \ldots, x_n]]$,
\begin{align*}
\multfzero
	&= \multzero{t, f_1 + c_1t^{k_1}}{f_n + c_nt^{k_n}}
\end{align*}
for any $c_1, \ldots, c_n \in \kk$. It follows that, for each collection of subsets $\scrA_1, \ldots, \scrA_n$ of $\znzero$,
\begin{align*}
\multAzero
	&= \multzero{\hat \scrA_0}{\hat \scrA_n}
\end{align*}
where $\hat \scrA_0 := \{(1, 0, \ldots, 0)\} \subset \zzeroo{n+1}$ and $\hat \scrA_j :=  \{(k_j, 0, \ldots, 0)\} \cup (\{0\} \times \scrA_j) \subset \zzeroo{n+1}$ for $j = 1, \ldots, n$. Let $\hat \scrA := (\hat \scrA_0, \ldots, \hat \scrA_n)$. It follows from \eqref{I1-list} that $\tAAone{\hat \scrA} = \{[n+1]\}$ and therefore, if $\multAzero < \infty$, then \cref{multiplicity-thm} implies that
\begin{align}
\multAzero
	= \multzerostar{\hat \Gamma_0}{\hat \Gamma_n}
	= \sum_{\hat \nu \in \hatVzeroprime}	\hat \nu_0 \mv'_{\hat \nu}(\In_{\hat \nu}(\hat \Gamma_1), \ldots, \In_{\hat \nu} (\hat \Gamma_n))
	\label{mult-formula-rojas}
\end{align}
where $\hat \Gamma_j$ are the Newton diagrams of $\hat \scrA_j$, and $\hat \nu$ ranges over the primitive weighted orders on $\kk[t, x_1, \ldots, x_n]$ which are centered at the origin, and $\hat \nu_0 := \hat \nu(t)$. Note that $ \mv'_{\hat \nu}(\In_{\hat \nu}(\hat \Gamma_1), \ldots, \In_{\hat \nu} (\hat \Gamma_n))$ is positive only if $\hat \nu'$ is the inner normal to a ``lower'' facet of $\hat \Gamma_1 + \cdots + \hat \Gamma_n$ (the designation ``lower'' comes from the fact that $\hat \nu'$ points ``upward'' along the $t$-coordinate). B.\ Huber and B.\ Sturmfels presented in \cite{hurmfels-bern} the idea of ``lifting'' subsets of $\zz^n$ to one extra dimension and summing the mixed volumes of faces corresponding to certain lower facets of the sum of the lifted bodies. J.\ M.\ Rojas \cite{rojas-toric} observed that the expression in the right hand side of \eqref{mult-formula-rojas} gives the generic intersection multiplicity at the origin. Note that unlike the formula \eqref{mult-formula} from \cref{multiplicity-thm}, the expression in \eqref{mult-formula-rojas} is symmetric in $\scrA_1, \ldots, \scrA_n$ (provided the $k_j$ are chosen to be equal).

\newcommand\diagramcolor{blue}
\newcommand\subdiagramcolor{green}

\begin{center}
\begin{figure}[h]
\begin{tikzpicture}[
    scale=0.45,
	dot/.style = {
      draw,
      fill,
      circle,
      inner sep = 0pt,
      minimum size = 3pt,
    }
    ]

\def\shiftone{12}
\def\opazero{0.5}
\def\tx{2.5}
\def\ty{2.5}
\def\gridx{4.5}
\def\gridy{4.5}
\def\colorone{blue}
\def\colortwo{red}

\draw [gray,  line width=0pt] (-0.5,-0.5) grid (\gridx,\gridy);
\draw [<->] (0, \gridy) |- (\gridx, 0);
\fill[\diagramcolor, opacity=\opazero ] (\gridx,1) -- (3,0) -- (1,1) -- (0,3) -- (3,4) -- cycle;
\fill[\subdiagramcolor, opacity=\opazero ] (3,0) -- (1,1) -- (0,3) -- (0,0) -- cycle;
\draw[thick] (\gridx,1) -- (3,0) -- (1,1) -- (0,3) -- (3,4) -- cycle;

\draw (\tx,\ty) node {\picfontsize $\scrA_1$};

\begin{scope}[shift={(\shiftone,0)}]
	\draw [gray,  line width=0pt] (-0.5,-0.5) grid (\gridx,\gridy);
	\draw [<->] (0, \gridy) |- (\gridx, 0);
	\fill[\diagramcolor, opacity=\opazero ] (\gridx,1) -- (2,1) -- (1,2) -- (0,4) -- (0,\gridy) -- (\gridx,\gridy) -- cycle;
	\fill[\subdiagramcolor, opacity=\opazero ] (2,1) -- (1,2) -- (0,4) -- (0,0) -- cycle;
	\draw[thick] (\gridx,1) -- (2,1) -- (1,2) -- (0,4) -- (0,\gridy);
	\draw (\tx,\ty) node {\picfontsize $\scrA_2$};
\end{scope}

\end{tikzpicture}

\caption{$\scrA_1$ is convenient, whereas $\scrA_2$ is not. The subdiagram volume of $\scrA_j$ is the area of the region shaded in \subdiagramcolor.} \label{fig:convenient}
\end{figure}
\end{center}

\subsection{Convenient Newton diagrams} \label{convenient-section}
\index{Convenient subset of $\rzeroo{n}$}
We say that a subset of $\rzeroo{n}$ is {\em convenient} if it contains a point on each coordinate axis. The \index{Subdiagram volume}{\em subdiagram volume} $\volsubn(\scrA)$ of a subset $\scrA$ of $\rzeroo{n}$ is the $n$-dimensional volume of the ``cone'' whose base is the Newton diagram of $\scrA$ and apex is at the origin; in other words, $\volsubn(\scrA)$ is the $n$-dimensional volume of the union of all line segments from the origin to $\nd(\scrA)$ (\cref{fig:convenient}).

\begin{prop} \label{convenient-prop}
Let $\scrA_1, \ldots, \scrA_n$ be subsets of $\znzero$. Let $\Gamma_j := \nd(\scrA_j)$, $j = 1, \ldots, n$.
\begin{enumerate}
\item \label{mult-convenient} If $\Gamma_2, \ldots, \Gamma_n$ are convenient, then
\begin{align}
\multAzero
	= \sum_{\nu \in \Vzeroprime} \min_{\Gamma_1}(\nu) \mv'_{\nu}(\In_{\nu}(\Gamma_2), \ldots, \In_{\nu} (\Gamma_n))
	\label{mult-formula-convenient}
\end{align}
\item \label{mult-convenient-equal} (Kushnirenko \cite[Theorem 22.8]{aizenberg-yuzhakov}) If $\Gamma$ is a convenient Newton diagram, and if $\Gamma_j = \Gamma$ for each $j$, then
\begin{align}
\multAzero
	= n!\volsubn(\Gamma)
	\label{mult-formula-convenient-equal}
\end{align}
\item \label{mult-convenient-alternate} (Ajzenberg and Yuzhakov \cite[Theorem 22.10]{aizenberg-yuzhakov}) If $\Gamma_1, \ldots, \Gamma_n$ are convenient, then
\begin{align}
\multAzero
	&= \sum_{\substack{I \subseteq [n]\\ I \neq \emptyset}}(-1)^{n - |I|}\volsubn  (\sum_{i \in I} \Gamma_i )
	\label{mult-formula-convenient-alternate}
\end{align}
\end{enumerate}
\end{prop}

\begin{proof}
If $\Gamma_2, \ldots, \Gamma_n$ are convenient, then $\tAone = \{[n]\}$, and \eqref{mult-formula-convenient} follows from \eqref{mult-formula}. Now we prove assertion \eqref{mult-convenient-equal}. Let $\{\scrQ_j\}_j$ be the facets of $\Gamma_2 + \cdots + \Gamma_n$ with inner normals in $\zpluss{n}$. Then \eqref{mult-formula-convenient} implies that
\begin{align*}
\multAzero
	= \sum_j \min_{\Gamma}(\nu_j) \mv'_{\nu_j}(\In_{\nu_j}(\Gamma), \ldots, \In_{\nu_j} (\Gamma))
	= (n-1)! \sum_j \min_{\Gamma}(\nu_j) \vol'_{\nu_j}(\In_{\nu_j}(\Gamma))
\end{align*}
where $\nu'_j$ are the inner normals to $\scrQ_j$ and $\vol'_{\nu_j}$ are as in \cref{volume-facet-rational}. Now fix $j$, and let $\scrR_j := \conv(\scrQ_j \cup \{\origin\})$. Then $\scrQ_j$ is a facet of $\scrR_j$ with {\em outer} primitive normal $\nu_j$, and all other facets of $\scrR_j$ passes through the origin. Since $\max_{\scrR_j}(\nu_j) = \min_\Gamma(\nu_j)$, \cref{volume-facet-rational} implies that $\vol_n (\scrR_j) = (1/n)\vol'_{\nu_j}(\In_{\nu_j}(\Gamma))$. Since $\volsubn(\Gamma) = \sum_j \vol_n(\scrR_j)$, identity \eqref{mult-formula-convenient-equal} follows. Since $\multAzero$ is multi-additive and symmetric in the $\scrA_j$, assertion \eqref{mult-convenient-alternate} then follows from \cref{tom-lemma}.
\end{proof}

The following is a more precise version of assertion \eqref{mult-convenient-equal} of \cref{convenient-prop}.

\begin{prop} \label{kushniplicity}
Let $\scrA_1, \ldots, \scrA_n$ be subsets of $\znzero$. Let $\Gamma_j := \nd(\scrA_j)$, $j = 1, \ldots, n$, and $\Gamma := \nd(\bigcup_{j =1}^n \scrA_j)$. For each $I \subseteq [n]$, let $\tiA$, where $\scrA:= (\scrA_1, \ldots, \scrA_n)$, be as in \cref{multiplicity-thm}. Then
\begin{enumerate}
\item \label{kushniplicity-bound} (Kushnirenko \cite[Theorem 22.8]{aizenberg-yuzhakov}) $\multAzero \geq n! \volsubn(\Gamma)$.
\item \label{kushniplicity-equal} $\multAzero = n! \volsubn(\Gamma)$ if and only if for each nonempty $I \subseteq [n]$, $|\tiA| \geq |I|$ and for each weighted order $\nu$ centered at the origin on $\kk[x_i: i \in I]$, the collection $\{\In_\nu(\Gamma_j \cap \ri):j \in \tiA,\ \Gamma_j \cap \In_\nu(\Gamma \cap \ri) \neq \emptyset\}$ of polytopes is dependent.
\end{enumerate}
\end{prop}

\begin{proof}
If $\multAzero = \infty$ then both assertions of the proposition are satisfied (for the second assertion one needs to use \cref{mult-support}). So assume $\multAzero < \infty$. Then $\Gamma$ is convenient, so that assertion \eqref{kushniplicity-bound} follows from assertion \eqref{mult-convenient-equal} of \cref{convenient-prop} and the definition of generic intersection multiplicity. Regarding the second assertion, \cref{non-degeneracy-0-thm} implies that $\multAzero = \multzero{\Gamma}{\Gamma}$ if and only if generic $(f_1, \ldots, f_n) \in \scrL_0(\scrA)$ are $\scrB$-non-degenerate at the origin, where $\scrB := (\bigcup_{j=1}^n \scrA_j, \ldots, \bigcup_{j=1}^n \scrA_j)$. The second assertion now follows from \cref{positively-mixed}.
\end{proof}

\subsection{Making $\scrA_j$ convenient without changing $\multAzero$.} \label{convenient-reduction}
If $\multfzero < \infty$, then the ideal generated by $f_1, \ldots, f_n$ in $\kk[[x_1, \ldots, x_n]]$ contains all sufficiently large powers of the maximal ideal of $\kk[[x_1, \ldots, x_n]]$. It follows that if we replace $f_j$ by $f_j + \sum_j c_{i,j} x_j^{d_{i,j}}$, then $\multfzero$ does not change for sufficiently large $d_{i,j}$. Since the Newton diagrams of the $f_j$ become convenient after these replacements, it follows that given any set of subsets $\scrA_1, \ldots, \scrA_n$ of $\znzero$ such that $\multAzero < \infty$, we may use \eqref{mult-formula-convenient} or \eqref{mult-formula-convenient-alternate} to compute $\multAzero$ after adding to each $\Gamma_j$ appropriate vertices on the coordinate axes. In this section we derive a ``sharp'' explicit bound on the placement of these vertices which guarantees that the intersection multiplicity at the origin remains unchanged. A.\ Khovanskii told the author in 2017 that he also had obtained, but never published, such a bound.

\begin{figure}[h]
\def\xmin{-0.5}
\def\xmax{16.5}
\def\ymin{-0.5}
\def\ymax{8.5}
\def\tx{0}
\def\ty{5}
\def\tw{5cm}
\def\numultiplier{1}
\def\colornu{purple}
\def\colorone{olive}

\tikzstyle{dot} = [\colordot, circle, minimum size=4pt, inner sep = 0pt, fill]

\begin{center}

\begin{tikzpicture}[scale=\scalefactor]
\draw [<->] (0, \ymax) |- (\xmax, 0);

\coordinate (A) at (2.5,2);
\node[dot] (B) at (13,1) {};
\node[dot] (C) at (7,8) {};
\node[dot] (D) at (8,5) {};
\coordinate (nu) at (2,1);
\coordinate (a1) at (0,7);
\coordinate (a2) at (3.5,0);

\draw[thick, \colorone] (a1) --  (a2);

\node[dot] at (A) {};

\node[anchor = east] at (A) {\picfontsize $A$};
\node[anchor = north] at (B) {\picfontsize $B$};
\node[anchor = south] at (C) {\picfontsize $C$};
\node[anchor = east] at (D) {\picfontsize $D$};
\node[anchor = east] at (a1) {\picfontsize $m_2(\scrS, \nu)$};
\node[anchor = north] at (a2) {\picfontsize $m_1(\scrS, \nu)$};

\coordinate (O) at (0.9, 5);
\node[anchor = west] at (O) {\picfontsize $H(\scrS, \nu)$};

\coordinate (nutip) at ($(A) + \numultiplier*(nu)$);
\draw [thick, \colornu, ->] (A) -- (nutip);
\node[anchor = west] at (nutip) {\picfontsize $\nu$};
\end{tikzpicture}

\caption{
	$m_i(\scrS, \nu)$ for $\scrS = \{A, B, C, D\}$
} \label{fig:minu}
\end{center}
\end{figure}

Let $\scrS$ be a compact subset of $\rr^n$ and $\nu$ be an element of $\rnstar$ centered at the origin. Let $H(\scrS, \nu) := \{ \alpha \in \rr^n: \langle \nu, \alpha \rangle = \min_\scrS(\nu)\}$ be the hyperplane perpendicular to $\nu$ which contains the ``face'' $\In_\nu(\scrS)$ of $\scrS$ corresponding to $\nu$. We write $m_i(\scrS, \nu)$ for the $i$-th coordinate of the point of the intersection of $H(\scrS,\nu)$ and the $i$-th coordinate axis (see \cref{fig:minu}). Note that
\begin{align*}
m_i(\scrS, \nu)
    & = \frac{\min_{\scrS}(\nu)}{\nu_i}
\end{align*}
where $\nu = (\nu_1, \ldots, \nu_n)$ with respect to the coordinate dual to the standard basis of $\rr^n$. Given a collection $\scrA = (\scrA_1, \ldots, \scrA_n)$ of subsets of $\znzero$, pick $I \in \tAone$, where $\tAone$ is as in \eqref{I1-list}. Let $j_1 = 1, j_2, \ldots, j_{|I|}$ be the elements of $\tiA$. For each $j$, let $\Gamma^I_j$ be the Newton diagram of $\scrA^I_j := \scrA_j \cap \ri$. Let $\scrV'^I_{0,1}(\scrA)$ be the set of primitive weighted orders centered at the origin on $\kk[x_i : i \in I]$ such that the $(|I|-1)$-dimensional mixed volume of $\In_\nu(\Gamma^I_{j_2}), \ldots, \In_\nu(\Gamma^I_{j_{|I|}})$ is nonzero; recall that faces with nonzero mixed volume can be detected combinatorially (\cref{positively-mixed}). Let $e_1, \ldots, e_n$ be the standard unit vectors in $\rr^n$. Define
\begin{align}
m^I_{i,1}(\scrA)
    & :=
       \left\{
            \begin{array}{@{}lll}
            k_i & \text{if}\ I = \{i\},\ \Gamma^I_1 = \{k_ie_i\}, \text{and}\ f_j|_{\kk^I} = 0\ \text{for each}\ j > 1
                & \text{(case 1)}\\
            1   & \text{else if}\ \scrV'^I_{0,1}(\scrA) = \emptyset
                & \text{(case 2)}\\
            \max_{\nu \in \scrV'^I_{0,1}(\scrA)} m_i(\Gamma^I_1, \nu)
                & \text{otherwise}
                & \text{(case 3)}
            \end{array}
        \right.
\label{mI}
\end{align}

\begin{center}
\begin{figure}[h]
\def\shiftone{7}
\def\opazero{0.5}
\def\tx{2.5}
\def\ty{2.5}
\def\gridx{4.5}
\def\gridy{4.5}
\def\colorzero{blue}
\def\colorone{blue}

\tikzstyle{sdot} = [red, circle, minimum size=3pt, inner sep = 0pt, fill]
\tikzstyle{bsdot} = [blue, circle, minimum size=3pt, inner sep = 0pt, fill]

\begin{subfigure}[b]{0.45\textwidth}
\begin{tikzpicture}[scale=0.45]
\draw [gray,  line width=0pt] (-0.5,-0.5) grid (\gridx,\gridy);
\draw [<->] (0, \gridy) |- (\gridx, 0);
\draw[ultra thick, \colorone]  (2,1) --  (0,3);
\fill[\colorzero, opacity=\opazero] (2,1) --  (0,3) -- (0,\gridy) -- (\gridx,\gridy) -- (\gridx,1) -- cycle;
\draw (\tx,\ty) node {\picfontsize $\scrA_1$};
\draw[thick, dashed]  (0,3) -- (1.5,0);

\node[sdot] at (0,3) {};
\node[anchor = east] at (0,3) {\picfontsize $m^{\{2\}}_{2,1}(\scrA)$};
\node[anchor = east] at (0,2) {\picfontsize (Case 1)};

\node[sdot] at (1.5,0) {};
\node[anchor = north] at (1.5,0) {\picfontsize $m^{\{1,2\}}_{1,1}(\scrA)$};
\node[anchor = north] at (1.5,-1) {\picfontsize (Case 3)};
\begin{scope}[shift={(\shiftone,0)}]
	\draw [gray,  line width=0pt] (-0.5,-0.5) grid (\gridx,\gridy);
	\draw [<->] (0, \gridy) |- (\gridx, 0);
	\draw[ultra thick, \colorone]  (2,0) -- (1,2);
	\fill[\colorzero, opacity=\opazero] (2,0) -- (1,2)--  (1,\gridy) -- (\gridx,\gridy) -- (\gridx,0) -- cycle;
	\draw (\tx,\ty) node {\picfontsize $\scrA_2$};
\end{scope}
\end{tikzpicture}
\caption{Cases 1 and 3}
\end{subfigure}
\hfill
\begin{subfigure}[b]{0.45\textwidth}
\begin{tikzpicture}[scale=0.45]
\draw [gray,  line width=0pt] (-0.5,-0.5) grid (\gridx,\gridy);
\draw [<->] (0, \gridy) |- (\gridx, 0);
\draw[ultra thick, \colorone]  (2,1) --  (0,3);
\fill[\colorzero, opacity=\opazero] (2,1) --  (0,3) -- (0,\gridy) -- (\gridx,\gridy) -- (\gridx,1) -- cycle;
\draw (\tx,\ty) node {\picfontsize $\scrA_1$};

\node[sdot] at (0,3) {};
\node[anchor = east] at (0,3) {\picfontsize $m^{\{2\}}_{2,1}(\scrA)$};
\node[anchor = east] at (0,2) {\picfontsize (Case 1)};

\node[sdot] at (0,1) {};
\node[anchor = east] at (0,1) {\picfontsize $m^{\{1,2\}}_{2,1}(\scrA)$};
\node[anchor = east] at (0,0) {\picfontsize (Case 2)};
\node[anchor = north] at (0,-1) {\phantom{(Case 3)}};
\begin{scope}[shift={(\shiftone,0)}]
	\draw [gray,  line width=0pt] (-0.5,-0.5) grid (\gridx,\gridy);
	\draw [<->] (0, \gridy) |- (\gridx, 0);
	\fill[\colorzero, opacity=\opazero] (1,0)--  (1,\gridy) -- (\gridx,\gridy) -- (\gridx,0) -- cycle;
	\node[bsdot] at (1,0) {};
	\draw (\tx,\ty) node {\picfontsize $\scrA_2$};
\end{scope}
\end{tikzpicture}
\caption{Cases 1 and 2}
\end{subfigure}
\caption{Different cases of \eqref{mI}}  \label{fig:mI}
\end{figure}
\end{center}
See \cref{fig:mI} for an illustration of different cases in the definition of $m^I_{i, 1}$. Define
\begin{align}
m_{i,1}(\scrA)
    &:= \max_{\substack {I \in \tAone \\ i \in I}} m^I_{i,1}
\label{convenient-reduction-bound}
\end{align}
Since $[n] \in \tAone$, it follows that $m^{[n]}_{i,1}(\scrA)$, and therefore $m_{i,1}(\scrA)$ is a well defined nonnegative rational number for each $i$. Let $\scrA'_1 := \scrA_1 \cup \{m'_1e_1, \ldots, m'_ne_n\}$, where $m'_i$ are arbitrary integers greater than or equal to $m_{i,1}(\scrA)$. Note that $\scrA'_1$ is convenient. Let $\scrA' := (\scrA'_1, \scrA_2, \ldots, \scrA_n)$.
\begin{prop}
Assume $\multAzero < \infty$. Then $\multzero{\scrA'_1,\scrA_2}{\scrA_n} = \multAzero$. If in addition $\multAzero > 0$, then the transformation $\scrA_1 \mapsto \scrA'_1$ is sharp in the following sense: if $\scrA''_1 \supset \scrA_1 \cup \{m''_ie_i\}$ for any $i$ and any nonnegative integer $m''_i$ such that $m''_i <  m_{i,1}(\scrA)$, then $\multzero{\scrA''_1,\scrA_2}{\scrA_n} < \multAzero$.
\end{prop}

\begin{proof}
If $\multAzero = 0$, then $0 \in \bigcup_j \scrA_j \subset \scrA'_1 \cup \bigcup_{j \geq 2}\scrA_j$, so that $\multzero{\scrA'_1,\scrA_2}{\scrA_n} = 0$ as well. Now assume $0 < \multAzero < \infty$. Then $\origin \not\in \scrA_1$, so that $m_i(\Gamma^I_1, \nu) > 0$ for each $i \in [n]$, $I \subseteq [n]$ and each nontrivial weighted order centered at the origin on $\kk[x_j : j \in I]$. It follows that $m_{i,1}(\scrA) > 0$ for each $i$, and therefore \cref{mult-support} implies that $0 < \multzero{\scrA'_1,\scrA_2}{\scrA_n}  < \infty$ and $\tAone= \tAprimeone$. Now pick $I \in \tAprimeone$. Let $j_1 = j,j_2, \ldots, j_{|I|}$ be the elements of $\tiAprime$. If $\Gamma'^I_1$ is the Newton diagram of $\scrA'_1 \cap \ri$, it is straightforward to see using the definition of $\multzerostar{\cdot}{\cdot}$ that $\multzerostar{\Gamma'^I_1, \Gamma^I_{j_2}}{\Gamma^I_{j_{|I|}}} =  \multzerostar{\Gamma^I_1, \Gamma^I_{j_2}}{\Gamma^I_{j_{|I|}}}$, and therefore \eqref{mult-formula} implies that $\multzero{\scrA'_1,\scrA_2}{\scrA_n} = \multAzero$. On the other hand, if for some $i \in I$, $\Gamma'^I_1$ contains an element on the $i$-th axis with coordinates $m''_ie_i$ such that $m''_i < m_{i,1}(\Gamma^I_1,\nu)$ for some $\nu \in  \scrV'^I_{0,1}(\scrA)$, then $\min_{\Gamma'^I_1}(\nu) < \min_{\Gamma^I_1}(\nu)$, and it would follow that $\multzerostar{\Gamma'^I_1, \Gamma^I_{j_2}}{\Gamma^I_{j_{|I|}}}
	<  \multzerostar{\Gamma^I_1, \Gamma^I_{j_2}}{\Gamma^I_{j_{|I|}}}$, which implies the last assertion.
\end{proof}

It is clear that given $\scrA_1, \ldots, \scrA_n$ such that $\multAzero < \infty$, repeating the above process $n$ times would yield a collection $\scrA'_1, \ldots, \scrA'_n$ of convenient subsets of $\znzero$ such that $\multzero{\scrA'_1}{\scrA'_n} = \multAzero$, as required. However, as \cref{fig:non-unique-convenientize} shows, the outcome of the process is in general {\em not} unique: different ordering of the $\scrA_j$ might result in different $\scrA'_1, \ldots, \scrA'_n$.

\begin{center}
\begin{figure}[h]
\def\shiftone{7}
\def\opazero{0.5}
\def\gridx{4.5}
\def\gridy{4.5}
\def\colorzero{green}
\def\colorone{red}

\tikzstyle{sdot} = [red, circle, minimum size=3pt, inner sep = 0pt, fill]
\tikzstyle{bsdot} = [blue, circle, minimum size=3pt, inner sep = 0pt, fill]

\subcaptionbox{
    $\Gamma'_1 = \Gamma'_2 =$ the line segment from $(2,0)$ to $(0,4)$
    \label{fig:non-unique-convenientize-a}
}[0.45\textwidth]{
    \begin{tikzpicture}[scale=0.45]
    \draw [gray,  line width=0pt] (-0.5,-0.5) grid (\gridx,\gridy);
    \draw [<->] (0, \gridy) |- (\gridx, 0);
    \draw[ultra thick, \colorzero]  (2,2) --  (0,4);
    \draw (1.5,3.5) node {\picfontsize $\Gamma_1$};
    \draw[ultra thick, loosely dashed, \colorone]  (0,4) -- (2,0);
    \draw (2,1) node {\picfontsize $\Gamma'_1$};

    \begin{scope}[shift={(\shiftone,0)}]
    	\draw [gray,  line width=0pt] (-0.5,-0.5) grid (\gridx,\gridy);
    	\draw [<->] (0, \gridy) |- (\gridx, 0);
    	\draw[ultra thick, \colorzero]  (2,0) -- (1,2);
    	\draw (2,1) node {\picfontsize $\Gamma_2$};
        \draw[ultra thick, loosely dashed, \colorone]  (0,4) -- (2,0);
        \draw (1,3) node {\picfontsize $\Gamma'_2$};
    \end{scope}
    \end{tikzpicture}
}
\hfill
\subcaptionbox{
    Reversing the order of the $\scrA_j$ from \cref{fig:non-unique-convenientize-a} changes $\Gamma'_1$ (respectively $\Gamma'_2$) to the line segment from $(2,0)$ to $(0,2)$ (respectively, from $(4,0)$ to $(0,4)$).
}[0.45\textwidth]{
    \begin{tikzpicture}[scale=0.45]
    \draw [gray,  line width=0pt] (-0.5,-0.5) grid (\gridx,\gridy);
    \draw [<->] (0, \gridy) |- (\gridx, 0);
    \draw[ultra thick, \colorzero]  (2,0) -- (1,2);
    \draw (2,1) node {\picfontsize $\Gamma_1$};
    \draw[ultra thick, loosely dashed, \colorone]  (0,2) -- (2,0);
    \draw (0.5,0.5) node {\picfontsize $\Gamma'_1$};
    \begin{scope}[shift={(\shiftone,0)}]
        \draw [gray,  line width=0pt] (-0.5,-0.5) grid (\gridx,\gridy);
        \draw [<->] (0, \gridy) |- (\gridx, 0);
        \draw[ultra thick, \colorzero]  (2,2) --  (0,4);
        \draw (1.5,3.5) node {\picfontsize $\Gamma_2$};
        \draw[ultra thick, loosely dashed, \colorone]  (0,4) -- (4,0);
        \draw (3.5,1.5) node {\picfontsize $\Gamma'_2$};
    \end{scope}
    \end{tikzpicture}
}
\caption{Dependence of $\{\scrA'_1, \scrA'_2\}$ on the ordering of $\scrA_j$: here $\Gamma_j, \Gamma'_i$ denote respectively the Newton diagram of $\scrA_j, \scrA'_i$.}  \label{fig:non-unique-convenientize}
\end{figure}
\end{center}

%

\section{Monotonicity of generic intersection multiplicity} \label{strictly-less-generic-0}
In \cref{monotonic-remark-0} we saw that $\multzero{\cdot}{\cdot}$ is ``monotonic'' as a function on $n$-tuples of subsets of $\znzero$. In this section we characterize in \cref{strictly-monotone-0} the conditions under which it is ``strictly monotonic,'' and as an application we prove the alternate formulation of non-degeneracy at the origin in \cref{alternate-non-degeneracy-0} (the counterparts of these results in the toric case are \cref{strictly-mixed-monotone,alternate-bernstein}). We also state a curious implication (\cref{initial-cor}) of \cref{non-degeneracy-0-thm} that in the case the monotonicity is not strict, the intersection multiplicity is determined by the Newton diagram of the intersection.

\begin{thm} \label{strictly-monotone-0}
Let $\scrB_j \subseteq \conv(\scrA_j) + \rzeroo{n}$, $j = 1, \ldots, n$. Then $\multAzero \leq \multBzero$. If $0 < \multAzero < \infty$, then the following are equivalent:
\begin{enumerate}
\item $\multAzero = \multBzero$,
\item for each nonempty subset $I$ of $[n]$, and each $\nu \in \rnstar$ which is centered at the origin, the collection $\{\In_\nu(\nd(\Bi_j)): \In_\nu(\nd(\Ai_j)) \cap \supp(\scrB_j) \neq \emptyset\}$ of polytopes is dependent,
\item for each $I \in \mscrIA := \{I \subseteq [n]:\  I \neq \emptyset,\ |\emptyiA| \geq n - |I|\}$ and each $\nu \in \rnstar$ which is centered at the origin, the collection $\{\In_\nu(\nd(\Bi_j)): \In_\nu(\nd(\Ai_j)) \cap \supp(\scrB_j) \neq \emptyset\}$ of polytopes is dependent.
\end{enumerate}
\end{thm}

\begin{proof}
This follows exactly in the same way as the proof of \cref{strictly-mixed-monotone} by considering generic $(f_1, \ldots, f_n) \in \scrL_0(\scrB)$, then observing that $\multfzero = \multAzero$ if and only if $f_1, \ldots, f_n$ are $\scrA$-non-degenerate at the origin, and finally applying \cref{bkk-bound-thm,positively-mixed} together with the genericness of $f_1, \ldots, f_n$.
\end{proof}

\begin{cor} \label{alternate-corollary-0}
\Cref{alternate-non-degeneracy-0} holds.
\end{cor}

\begin{proof}
The equivalence of assertions \eqref{assn:non-deg-0} and \eqref{assn:non-deg-0-combinatorial} follows exactly as in the proof of \cref{alternate-corollary}. \Cref{non-degeneracy-0'-thm} implies that assertions \eqref{assn:non-deg-0-combinatorial} and \eqref{assn:non-deg-0-efftorial} are equivalent.
\end{proof}

In \cref{initial-cor} below we use the following notation: given $f = \sum_\alpha c_\alpha x^\alpha \in \kk[[x_1, \ldots, x_n]]$ and $\scrS \subseteq \rr^n$, we write $f_\scrS := \sum_{\alpha \in \scrS} c_\alpha x^\alpha$.

\begin{prop} \label{initial-cor}
Let $(f_1, \ldots, f_n) \in \scrL_0(\scrA)$ be such that $\multfzero = \multAzero$. Then
\begin{align*}
\multzero{f_{1, \nd(\scrA_1)}}{f_{n,\nd(\scrA_n)}} = \multfzero = \multAzero
\end{align*}
In particular, if $\scrB_j$ are subsets of $\scrA_j + \rzeroo{n}$ such that $\multBzero = \multAzero$, then
\begin{align*}
\multzero{\scrB_1 \cap \nd(\scrA_1)}{\scrB_n \cap \nd(\scrA_n)} = \multBzero = \multAzero
\end{align*}
\end{prop}

\begin{proof}
Follows immediately from \cref{non-degeneracy-0-thm}, since $f_1, \ldots, f_n$ are non-degenerate at the origin if and only if $f_{1, \nd(\scrA_1)}, \ldots, f_{n,\nd(\scrA_n)}$ are non-degenerate at the origin.
\end{proof}

The requirement that $\scrB_j \subseteq \scrA_j + \rzeroo{n}$ for each $j$ is necessary for \cref{initial-cor}. Indeed, if $\scrA_j, \scrB_j$, $j = 1, 2$, are from \cref{initial-figure}, then $\multpzeronodots{\scrA_1, \scrA_2} = \multpzeronodots{\scrB_1, \scrB_2} = 9$, but $\multpzeronodots{\scrB_1 \cap \nd(\scrA_1), \scrB_2 \cap \nd(\scrA_2)} = \infty$.

\begin{figure}[h]
\begin{center}

\begin{tikzpicture}[scale=0.45]
\def\shiftone{7}
\def\opazero{0.5}
\def\tx{2.5}
\def\ty{2.5}
\def\gridx{4.5}
\def\gridy{4.5}
\def\colorone{red}

\draw [gray,  line width=0pt] (-0.5,-0.5) grid (\gridx,\gridy);
\draw [<->] (0, \gridy) |- (\gridx, 0);
\draw[ultra thick, \colorone]  (3,0) --  (0,3);
\fill[\colorzero, opacity=\opazero ] (3,0) --  (0,3) -- (0,\gridy) -- (\gridx,\gridy) -- (\gridx,0) -- cycle;
\draw (\tx,\ty) node {\picfontsize $\scrA_1$};

\begin{scope}[shift={(\shiftone,0)}]
	\draw [gray,  line width=0pt] (-0.5,-0.5) grid (\gridx,\gridy);
	\draw [<->] (0, \gridy) |- (\gridx, 0);
	\draw[ultra thick, \colorone]  (4,0) -- (2,1) -- (1,2)--  (0,4);
	\fill[\colorzero, opacity=\opazero ] (4,0) -- (2,1) -- (1,2)--  (0,4)-- (0,\gridy) -- (\gridx,\gridy) -- (\gridx,0) -- cycle;
	\draw (\tx,\ty) node {\picfontsize $\scrA_2$};
\end{scope}

\begin{scope}[shift={(2*\shiftone,0)}]
	\draw [gray,  line width=0pt] (-0.5,-0.5) grid (\gridx,\gridy);
	\draw [<->] (0, \gridy) |- (\gridx, 0);
	\draw[ultra thick, \colorone]  (4,0) -- (2,1) -- (1,2)--  (0,4);
	\fill[\colorzero, opacity=\opazero ] (4,0) -- (2,1) -- (1,2)--  (0,4)-- (0,\gridy) -- (\gridx,\gridy) -- (\gridx,0) -- cycle;
	\draw (\tx,\ty) node {\picfontsize $\scrB_1$};
\end{scope}

\begin{scope}[shift={(3*\shiftone,0)}]
	\draw [gray,  line width=0pt] (-0.5,-0.5) grid (\gridx,\gridy);
	\draw [<->] (0, \gridy) |- (\gridx, 0);
	\draw[ultra thick, \colorone]  (3,0) --  (0,3);
	\fill[\colorzero, opacity=\opazero ] (3,0) --  (0,3) -- (0,\gridy) -- (\gridx,\gridy) -- (\gridx,0) -- cycle;
	\draw (\tx,\ty) node {\picfontsize $\scrB_2$};
\end{scope}
\end{tikzpicture}

\end{center}
\caption{Failure of \cref{initial-cor} when there is $j$ such that $\scrB_j \not\subseteq \scrA_j + \rzeroo{n}$} \label{initial-figure}
\end{figure}

\section{Notes}
For convenient Newton diagrams there is a formula for generic intersection multiplicity in terms of integer lattice points in the region bounded by the diagram and the coordinate hyperplanes (see e.g.\ \cite[Theorem 5]{esterov}). A.\ Khovanskii informed the author that he had obtained (but did not publish) a bound equivalent to \eqref{convenient-reduction-bound} which reduces the computation of generic intersection multiplicity to the convenient case. Recently M.\ Herrero, G.\ Jeronimo and J.\ Sabia \cite{herrero-jeronimo-sabia-multiplicity} gave some other formulae for generic intersection mulitplicity in the general case.

\chapter{Number of zeroes on the affine space II: the general case} \label{affine-chapter}
\chaptermark{Number of zeroes on the affine space II}

\newcommand{\DBI}{\scrD^I_{\scrB}}
\newcommand{\DBprimeI}{\scrD^I_{\scrB'}}

\newcommand{\dUA}{\mscrD(U,\scrA)}
\newcommand{\dUUA}[1]{\mscrD(#1,\scrA)}
\newcommand{\dUAstar}{\mscrD^*(U,\scrA)}
\newcommand{\dUB}{\mscrD(U,\scrB)}
\newcommand{\dUBstar}{\mscrD^*(U,\scrB)}
\newcommand{\DUI}{\scrD^I_U}
\newcommand{\iUA}{\mscrI(U,\scrA)}
\newcommand{\iUUA}[1]{\mscrI(#1,\scrA)}
\newcommand{\iUAstar}{\mscrI^*(U,\scrA)}
\newcommand{\iUB}{\mscrI(U,\scrB)}
\newcommand{\iUBstar}{\mscrI^*(U,\scrB)}

\newcommand{\exc}{\mscrE}
\newcommand{\excbarA}{\bar \exc(\scrA)}
\newcommand{\excbarAA}[1]{\bar \exc(#1)}
\newcommand{\excA}{\exc(\scrA)}
\newcommand{\excU}{\exc(U)}
\newcommand{\excUA}{\exc(U,\scrA)}
\newcommand{\excUAI}{\exc^I(U,\scrA)}
\newcommand{\excUAprimeI}{\exc^I(U,\scrA')}

\newcommand{\multPstarS}{\multpstar{\scrP_1}{\scrP_n}{\mscrS}}

\newcommand{\VstarI}{(V^*)^I}
\newcommand{\VstarIf}{(V^*)^I(f)}

\section{Introduction}

In this chapter we compute the number of solutions on $\kk^n$ (or more generally, on any given Zariski open subset of $\kk^n$) of generic systems of polynomials with given supports, and give explicit BKK-type characterizations of genericness in terms of initial forms of the polynomials. As a special case we derive generalizations of weighted (multi-homogeneous)-B\'ezout theorems involving arbitrary weighted degrees (i.e.\ weighted degrees with possibly negative or zero weights). 
\section{The bound}

\subsection{Khovanskii's formula} \label{first-section}
For polynomials $f_1, \ldots, f_n$, and any Zariski open subset $U$ of $\kk^n$, as in \cref{bkk-section} let $\multfiso{U}$ be the sum of intersection multiplicities of $f_1, \ldots, f_n$ at all the isolated points of $V(f_1, \ldots, f_n) \cap U$. Given a collection $\scrA := (\scrA_1, \ldots, \scrA_n)$ of $n$ finite subsets of $\znzero$, define
\begin{align*}
\multAiso{U} &:= \max\{\multfiso{U}:  \supp(f_j) \subset \scrA_j,\ j = 1, \ldots, n\}
\end{align*}
In this section we give a formula for $\multAiso{U}$ in terms of (mixed volumes of) convex hulls of $\scrA_j$. For $I \subseteq [n]$, let $\tiA := \{j: \Ai_j \neq \emptyset\}$ as in \cref{multiplicity-thm}, and let
\begin{align}
\excA := \{I \subseteq [n]:\ \text{there is $\tilde I \supseteq I$ such that $|\titildeA| < |\tilde I|$}\} \label{excA}
\end{align}

The following result, which follows immediately from \cref{thm:pure-dimension}, implies that $\multAiso{U} = \multAiso{\ua}$, where $\ua := U\setminus \bigcup_{i \in \excA} \Ki$.

\begin{lemma} \label{non-isolated-lemma}
Let $f_1, \ldots, f_m \in \kk[x_1, \ldots, x_n]$ and let $V := V(f_1, \ldots, f_m) \subset \kk^n$. Given $I \subseteq [n]$, if there exists $\tilde I \supseteq I$ such that $|\{j: f_j|_{\Kitilde} \not\equiv 0\}| < |\tilde I|$, then no point of $V \cap \Ki$ is isolated in $V$. \qed
\end{lemma}

Define
\begin{align}
\excU  &:= \{I \subseteq [n]: \Ki \cap U = \emptyset\}, \label{excU}\\
\touchUA & := \{I \subseteq [n]: I \notin \excU \cup \excA,\ |\tiA| = |I|\}. \label{touchUA}
\end{align}

The following result is also a straightforward consequence of \cref{thm:pure-dimension}.

\begin{lemma} \label{isolated-lemma}
If $I \not\in \touchUA$, then $\Kstari \cap V(f_1, \ldots, f_n) \cap U = \emptyset$ for generic $f_1, \ldots, f_n$ such that $\supp(f_j) \subseteq \scrA_j$, $j =1 , \ldots, n$. \qed
\end{lemma}

\begin{rem} \label{empty-touch}
It is possible that $\emptyset$ (i.e.\ the empty set) is in $\touchUA$; this is the case if and only if the origin is in $U$ and $0 < \multAzero < \infty$ (see \cref{ex-bkk-0}).
\end{rem}

\begin{thm}[Khovanskii\footnote{A.\ Khovanskii described this unpublished formula to the author during the Askoldfest, 2017}] \label{extended-bkk-bound-0}
Let $\scrP_j$ be the convex hull of $\scrA_j$, $j = 1, \ldots, n$. For $I \subseteq [n]$, let $\scrP^I_j := \scrP_j \cap \rr^I$, and let $\pi_I: \rr^n \to \ri$ be the natural projection (as in \eqref{pi_I}). Then
\begin{align}
\multAiso{U}
	= \sum_{I \in \touchUA}
							\mv(\scrP^I_{j_1}, \ldots, \scrP^I_{j_{|I|}})
							\times
							\multzero{\pi_{[n]\setminus I}(\scrP_{j'_1})}{\pi_{[n]\setminus I}(\scrP_{j'_{n-|I|}})}
	\label{extended-bkk-formula-0}
\end{align}
where for each $I \in \touchUA$, $j_1 , \ldots, j_{|I|}$ are the elements of $\tiA$, and $j'_1, \ldots, j'_{n-|I|}$ are the elements of $[n]\setminus \tiA$, and $\multzero{\cdot}{\cdot}$ is defined as in \eqref{multAzero}.
\end{thm}

The interpretation of the right hand side of \eqref{extended-bkk-formula-0} is straightforward - for each $I \in \touchUA$, the corresponding summand counts with multiplicity the number of solutions on $\Kstari \cap U$ of generic systems supported at $\scrA_1 ,\ldots, \scrA_n$.
In the next section we present another formula which sometimes is more efficient, since it involves summing over elements from a proper subset of $\touchUA$.

\begin{example} \label{ex-bkk-0}
Let $\scrA_j$ be the support of $f_j$ from \cref{ex-mv',ex0}, and $U$ be a nonempty Zariski open subset of $\kk^3$. Then $\excA = \emptyset$, and
\begin{align*}
\touchUA &=
	\begin{cases}
	\{\{1,2,3\},\{3\}, \emptyset\} &\text{if $\origin \in U$ (Case 1),}\\
	\{\{1,2,3\},\{3\}\} &\text{if $\origin \not\in U$, but $U$ contains a point on the $z$-axis (Case 2),}\\
	\{\{1,2,3\}\} &\text{otherwise (Case 3).}
	\end{cases}
\end{align*}
In Case 3, identity \eqref{extended-bkk-formula-0} and \cref{ex-mv'} imply that
\begin{align*}
\multpisonodots{\scrA_1,\scrA_2,\scrA_3}{U}  = \mv(\scrP_1, \scrP_2, \scrP_3) = 5
\end{align*}
Since the projections of $\scrP_2$ and $\scrP_3$ onto the $(x,y)$-plane have nontrivial linear part, and $\scrP^{\{3\}}_1$ has integer length $1$ (see \cref{fig-ex-bkk-0}), identity \eqref{extended-bkk-formula-0} implies that in Case 2,
\begin{align*}
\multpisonodots{\scrA_1,\scrA_2,\scrA_3}{U}
	&= 5+ \mv(\scrP^{\{3\}}_1) \times \multpnodots{\pi_{\{1,2\}}(\scrP_2),\pi_{\{1,2\}}(\scrP_3)}{\origin}
	= 5+ 1\cdot 1
	= 6
\end{align*}
 Finally, in Case 1, identity \eqref{extended-bkk-formula-0}  and the computation from \cref{ex0} imply that
\begin{align*}
\multpisonodots{\scrA_1,\scrA_2,\scrA_3}{U}
	&= 6 + \multpnodots{\scrP_1,\scrP_2,\scrP_3}{\origin}
	= 6 + 6
	= 12
\end{align*}

\begin{figure}[h]
\begin{center}

\begin{tikzpicture}[scale=0.6]
\pgfplotsset{every axis/.append style = {view={\viewx}{\viewy}, axis lines=middle, enlargelimits={upper}}}
\begin{scope}
\begin{axis}
	\addplot3 [fill=blue,opacity=0] coordinates{(3,0,0) (0,3,0) (0,1,0) (1,0,0) (3,0,0)};
	\addplot3 [blue, line width=4pt] coordinates{(0,0,1) (0,0,2)};
\end{axis}
\draw (1,4.5) node {$\scrP^{\{3\}}_1$};
\end{scope}

\begin{scope}[shift={(\shiftone,0)}]
\begin{axis}
	\addplot3 [black] coordinates{(0,0,1) (0,0,2)};
	\addplot3 [fill=blue,opacity=\opazero, thick] coordinates{(3,0,0) (0,3,0) (0,1,0) (1,0,0) (3,0,0)};
\end{axis}
\draw (4,0.5) node {$\pi_{\{1,2\}}(\scrP_2)$};
\end{scope}

\begin{scope}[shift={(2*\shiftone,0)}]
\begin{axis}
	\addplot3 [black] coordinates{(0,0,1) (0,0,2)};
	\addplot3 [fill=blue,opacity=\opazero, thick] coordinates{(2,0,0) (0,2,0) (0,1,0) (1,0,0) (2,0,0)};
\end{axis}
\draw (4,0.5) node {$\pi_{\{1,2\}}(\scrP_3)$};
\end{scope}
\end{tikzpicture}
\end{center}
\caption{Computing $\mv(\scrP^{\{3\}}_1) \times \multpnodots{\pi_{\{1,2\}}(\scrP_2),\pi_{\{1,2\}}(\scrP_3)}{\origin}$}
\label{fig-ex-bkk-0}
\end{figure}

\end{example}

\subsection{A formula in the same spirit as the formula for generic intersection multiplicity} \label{second-section}
If $\nu$ is a weighted order on $\kk[x_1, \ldots, x_n]$, we say that $\nu$ is \index{Weighted!order!centered at infinity}{\em centered at infinity} if $\nu(x_i) < 0$ for some $i \in [n]$. Given $I \subseteq [n]$, we say that $\nu$ is {\em centered at $\Ki$} if $\nu(x_i) \geq 0$ for each $i \in I$ and $\nu(x_{i'}) > 0$ for each $i' \in [n]\setminus I$. Given a collection $\mscrS$ of subsets of $[n]$, we denote by $\Kns$ the complement in $\kk^n$ of the coordinate subspaces $\Ki$ for all $I \in \mscrS$, i.e.\
\begin{align}
\Kns &:= \kk^n \setminus \bigcup_{I \in \mscrS} \Ki \label{Kns}
\end{align}
We write $\Vs$ for the union, over all $I \in \scrS$, of the sets of weighted orders centered at $\Ki$, and $\Vinfty$ for the set of weighted orders centered at infinity; the collection of primitive elements in $\Vs$ and $\Vinfty$ are denoted respectively as $\Vprimes$ and $\Vprimeinfty$.

\begin{example} \label{Kns-remark}
Taking $\mscrS = \emptyset$ gives $\Kns =  \kk^n$ and $\Vs= \emptyset$. If we take $\mscrS=\{\emptyset\}$, then $\Kns= \kk^n \setminus \{0\}$ and $\Vs$ is the set $\Vzero$ of weighted orders {\em centered at the origin} (see \cref{mult-bound-statement-section}). If $\mscrS$ is the set of all subsets of $[n]$ consisting of $n-1$ elements, then $\Kns = \nktorus$ and $\Vs$ is the set of all nonzero weighted orders which are not centered at infinity.
\end{example}

Let $\scrA := (\scrA_1, \ldots, \scrA_n)$ be a collection of subsets of $\zzeroo{n}$ and $\scrP_j$ be the convex hull in $\rr^n$ of $\scrA_j$, $j = 1, \ldots, n$. Given a collection $\mscrS$ of subsets of $[n]$, define
\begin{align}
\multPstarS
	&:=  -\sum_{\nu \in \Vprimes \cup \Vprimeinfty} \min_{\scrP_1}(\nu) ~
		\mv'_\nu(\In_\nu(\scrP_2), \ldots, \In_\nu(\scrP_n))
	\label{multPstar}
\end{align}
where $\mv'_\nu(\cdot, \ldots, \cdot)$ is defined as in \eqref{mv'}.

\begin{thm}[{\cite{toricstein}}]\label{extended-bkk-bound}
Let $U$ be a Zariski open subset of $\kk^n$. We continue to use the notation of \cref{extended-bkk-bound-0}. Define $\touchUAone  := \{I \in \touchUA:  1 \in \tiA\}$, and for each $I \subseteq [n]$, set $\excUAI := \{J \subseteq I: J \in \excU \cup \excA\}$. Then
\begin{align}
\multAiso{U}
	= \sum_{I \in \touchUAone}
							\multpstar{\scrP^{I}_1, \scrP^{I}_{j_2}}{\scrP^{I}_{j_{|I|}}}{\excUAI}
							\times
							\multzero{\pi_{[n]\setminus I}(\scrP_{j'_1})}{\pi_{[n]\setminus I}(\scrP_{j'_{n-|I|}})}
	\label{extended-bkk-formula}
\end{align}
where for each $I \in \touchUAone$, $j_1 = 1, j_2, \ldots, j_{|I|}$ are the elements of $\tiA$, and $j'_1, \ldots, j'_{n-|I|}$ are the elements of $[n]\setminus \tiA$.
\end{thm}

There is an obvious analogy between formula \eqref{extended-bkk-formula} and the formula \eqref{mult-formula} for intersection multiplicity at the origin. The interpretation of the terms on the right hand side of \eqref{extended-bkk-formula} is also analogous to the interpretation of the corresponding terms of \eqref{mult-formula} described in \cref{mult-motivation}; in particular, for each $I \in \touchUAone$, the corresponding summand on the right hand side of \eqref{extended-bkk-formula} is the sum, for generic $f_1, \ldots, f_n$ supported respectively at $\scrA_1, \ldots, \scrA_n$, of the {\em negative} of orders of $f_1$ along the branches (counted with appropriate multiplicities) of the curve determined by $f_2 = \cdots = f_n = 0$ which lie on $\Ki$ and are centered either at infinity or at $\Kiprime$ for some $I' \in \excUAI$. As in the case of \eqref{mult-formula}, and unlike \eqref{extended-bkk-formula-0}, the symmetry of the right hand side of \eqref{extended-bkk-formula} with respect to permutations of the $\scrP_j$ is not at all obvious. We use \eqref{extended-bkk-formula} to derive the symmetric formula of Huber and Sturmfels \cite{hurmfels-bern} and Rojas \cite{rojas-toric} in \cref{hsr-affine}.

\def\shiftx{6}
\def\colorzero{blue}
\def\opazero{0.5}
\def\viewx{70}
\def\viewy{5}
\def\titlex{1}
\def\titley{-1}
\def\scalefactor{0.6}
\def\colorzero{blue}
\def\colorone{green}

\begin{figure}[h]
\begin{center}
\begin{tikzpicture}[scale=0.6]
\pgfplotsset{every axis/.append style = {view={\viewx}{\viewy}, axis lines=middle, enlargelimits={upper}}}

\begin{scope}
\begin{axis}
	\addplot3 [fill=blue,opacity=\opazero, thick] coordinates{(3,0,0) (0,3,0) (0,1,4) (1,0,4) (3,0,0)};
	\addplot3 [fill=red,opacity=\opazero, thick] coordinates{(3,0,0) (1,0,4) (1,0,3) (2,0,1) (3,0,0)};
	\addplot3 [dashed, black, thick] coordinates{(1,0,3) (0,1,3) (0,1,4) (1,0,4) (1,0,3)};
	\addplot3 [dashed, black, thick] coordinates{(1,0,3) (0,1,3) (0,2,1) (2,0,1) (1,0,3)};
	\addplot3 [dashed, black, thick] coordinates{(0,2,1) (2,0,1) (3,0,0) (0,3,0) (0,2,1)};
\end{axis}
\draw (\titlex,\titley) node {$\scrP_1 + \scrP_3$};
\end{scope}

\def\factorone{4.5}
\def\factortwo{4.5}
\def\factorthree{1}
\begin{scope}[shift={(\shiftx,0)}]
\begin{axis}[hide axis, enlargelimits={0.15}]
	\addplot3 [ultra thick, ->] coordinates{(0,0,0) (2,2,1)};
	\addplot3 [ultra thick, blue, ->] coordinates{(0,0,0) (-2*\factorthree,-2*\factorthree,-\factorthree)};
	\addplot3 [ultra thick, ->] coordinates{(0,0,0) (1,1,1)};
	\addplot3 [ultra thick, ->] coordinates{(0,0,0) (\factortwo,\factortwo,0)};	
	\addplot3 [ultra thick, red, ->] coordinates{(0,0,0) (0,\factorone,0)};
	\addplot3 [ultra thick, ->] coordinates{(0,0,0) (\factortwo,0,0)};
	\draw (axis cs:2,2,1) node [right] {(2,2,1)};
	\draw (axis cs:1,1,1) node [above] {(1,1,1)};
	\draw (axis cs:0,\factorone,0) node [right] {(0,1,0)};
	\draw (axis cs: \factortwo,0,0) node [below] {(1,0,0)};
	\draw (axis cs:\factortwo,\factortwo,0) node [below] {(1,1,0)};
	\draw (axis cs:-2*\factorthree,-2*\factorthree,-\factorthree) node [below] {(-2,-2,-1)};
\end{axis}
\draw (\titlex,\titley) node [right, text width=3.5cm] {inner normals to facets of $\scrP_1 + \scrP_3$};
\end{scope}

\end{tikzpicture}
\end{center}
\caption{$\scrP_1+ \scrP_3$ and the inner normals to its facets} \label{fig-ex-bkk-1-1}
\end{figure}

\begin{example}\label{ex-bkk-1}
We continue with $\scrA_1, \scrA_2, \scrA_3$ from \cref{ex-bkk-0}, and compute $\multpisonodots{\scrA_1,\scrA_2,\scrA_3}{U} $ using \cref{extended-bkk-bound} for nonempty Zariski open subsets $U$ of $\kk^3$. It is straightforward to check that
\begin{align*}
\touchUAone
	=
	\begin{cases}
	\{\{1,2,3\},\{3\}\} &\text{if $U$ contains a point on the $z$-axis (Case 1 or 2 of \cref{ex-bkk-0}),}\\
	\{\{1,2,3\}\} &\text{otherwise (Case 3 of \cref{ex-bkk-0}).}
	\end{cases}
\end{align*}
On the other hand, if we change the order of the $\scrA_j$, or equivalently, consider the collection $\scrA' := (\scrA_2, \scrA_1, \scrA_3)$, then one checks that $\touchUAprimeone =\{ \{1,2,3\}\}$ for any (nonempty) $U$, and we apply identity \eqref{extended-bkk-formula} to $\scrA'$ to reduce computation. In particular, we have
\begin{align*}
\multpisonodots{\scrA_1,\scrA_2,\scrA_3}{U}
	&= \multpisonodots{\scrA_2,\scrA_1\scrA_3}{U}
	= \multpstarnodots{\scrP_2, \scrP_1, \scrP_3}{\excUA}
	= -\sum_{\nu \in \Vprimess{\excUA} \cup \Vprimeinfty} \min_{\scrP_2}(\nu) ~
			\mv'_\nu(\In_\nu(\scrP_1), \In_\nu(\scrP_3))
\end{align*}

\begin{figure}[h]
\begin{center}
\begin{tikzpicture}[scale=0.4]

\def\shiftone{4}
\def\shifttwo{5}
\def\nux{-0.5}
\def\nuy{3.5}
\def\plusy{1.5}
\def\opazero{0.5}
\def\tx{-0.5}
\def\ty{-0.5}
\def\gridx{3.5}
\def\gridxx{3.5}
\def\gridy{2.5}
\def\bigshiftdhor{3}
\def\bigshiftdver{3}

\draw (\nux,\nuy) node [right] {\picfontsize $\nu = (-2,-2,-1)$};
\draw (\shiftone,\plusy) node {$+$};    	
\draw (\shiftone+\shifttwo,\plusy) node {$=$};    	

\draw [gray,  line width=0pt] (-0.5,-0.5) grid (\gridx,\gridy);

\draw[ultra thick, fill=\colorzero, opacity=\opazero ] (0,0) -- (1,0) -- (0,1) -- cycle;
\draw (\tx,\ty) node [below right] {\picfontsize $\In_\nu(\scrP_1)$};    	

\begin{scope}[shift={(\shifttwo,0)}]
	\draw [gray,  line width=0pt] (-0.5,-0.5) grid (\gridx,\gridy);
	
	\draw[ultra thick, fill=\colorzero, opacity=\opazero] (0,0) -- (2,0) -- (1,1) -- (0,1) --cycle;
	\draw (\tx,\ty) node [below right] {\picfontsize $\In_\nu(\scrP_3)$};    	
\end{scope}

\begin{scope}[shift={(2*\shifttwo,0)}]
	\draw [gray,  line width=0pt] (-0.5,-0.5) grid (\gridxx,\gridy);
	
	\draw[ultra thick] (0,0) -- (3,0) -- (1,2) --  (0,2) --  cycle;
	\draw[fill=\colorone, opacity=\opazero] (2,0) -- (1,1) -- (2,1) -- (3,0) -- cycle;
	\draw[fill=\colorone, opacity=\opazero] (0,1) -- (1,1) -- (1,2) -- (0,2) -- cycle;
	\draw[fill=\colorzero, opacity=\opazero, dashed, ultra thick] (0,0) -- (2,0) -- (1,1) -- (0,1) --cycle;
	\draw[fill=\colorzero, opacity=\opazero, dashed, ultra thick ] (1,1) -- (2,1) -- (1,2) -- cycle;
	
	\draw (\tx,\ty) node [below right] {\picfontsize $\In_\nu(\scrP_1 + \scrP_3)$};    	
\end{scope}

\begin{scope}[shift={(2*\shifttwo+\gridxx+\bigshiftdhor,0)}]
	\draw (\nux,\nuy) node [right] {\picfontsize $\nu = (2,2,1)$};
	\draw (\shiftone,\plusy) node {$+$};    	
	\draw (\shiftone+\shifttwo,\plusy) node {$=$};    	
	
	\draw [gray,  line width=0pt] (-0.5,-0.5) grid (\gridx,\gridy);
	
   \fill[blue] (0,0) circle (6pt);
	\draw (\tx,\ty) node [below right] {\picfontsize $\In_\nu(\scrP_1)$};    	
	
	\begin{scope}[shift={(\shifttwo,0)}]
		\draw [gray,  line width=0pt] (-0.5,-0.5) grid (\gridx,\gridy);
		
		\draw[ultra thick, fill=\colorzero, opacity=\opazero] (0,0) -- (2,0) -- (1,1) -- (0,1) --cycle;
		\draw (\tx,\ty) node [below right] {\picfontsize $\In_\nu(\scrP_3)$};    	
	\end{scope}
	
	\begin{scope}[shift={(2*\shifttwo,0)}]
		\draw [gray,  line width=0pt] (-0.5,-0.5) grid (\gridxx,\gridy);
		
		\draw[fill=\colorzero, opacity=\opazero, dashed, ultra thick] (0,0) -- (2,0) -- (1,1) -- (0,1) --cycle;
		
		\draw (\tx,\ty) node [below right] {\picfontsize $\In_\nu(\scrP_1 + \scrP_3)$};    	
	\end{scope}
\end{scope}		

\begin{scope}[shift={(0,-\nuy-\bigshiftdver)}]
\draw (\nux,\nuy) node [right] {\picfontsize $\nu = (1,1,1)$};
\draw (\shiftone,\plusy) node {$+$};    	
\draw (\shiftone+\shifttwo,\plusy) node {$=$};    	

\draw [gray,  line width=0pt] (-0.5,-0.5) grid (\gridx,\gridy);

\draw[ultra thick, fill=\colorzero, opacity=\opazero ] (0,0) -- (1,0) -- (0,1) -- cycle;
\draw (\tx,\ty) node [below right] {\picfontsize $\In_\nu(\scrP_1)$};    	

\begin{scope}[shift={(\shifttwo,0)}]
	\draw [gray,  line width=0pt] (-0.5,-0.5) grid (\gridx,\gridy);
	
	\draw[ultra thick, blue] (0,0) -- (2,0);
	\draw (\tx,\ty) node [below right] {\picfontsize $\In_\nu(\scrP_3)$};    	
\end{scope}

\begin{scope}[shift={(2*\shifttwo,0)}]
	\draw [gray,  line width=0pt] (-0.5,-0.5) grid (\gridxx,\gridy);
	
	\draw[ultra thick] (0,0) -- (3,0) -- (2,1) --  (0,1) --  cycle;
	\draw[fill=\colorone, opacity=\opazero] (1,0) -- (3,0) -- (2,1) -- (0,1) -- cycle;
	\draw[fill=\colorzero, opacity=\opazero, dashed, ultra thick]  (0,0) -- (1,0) -- (0,1) -- cycle;

	\draw (\tx,\ty) node [below right] {\picfontsize $\In_\nu(\scrP_1 + \scrP_3)$};    	
\end{scope}

\begin{scope}[shift={(2*\shifttwo+\gridxx+\bigshiftdhor,0)}]
	\draw (\nux,\nuy) node [right] {\picfontsize $\nu = (1,1,0)$};
	\draw (\shiftone,\plusy) node {$+$};    	
	\draw (\shiftone+\shifttwo,\plusy) node {$=$};    	
	
	\draw [gray,  line width=0pt] (-0.5,-0.5) grid (\gridx,\gridy);
	
	\draw[ultra thick, blue] (0,0) -- (0,1);
	\draw (\tx,\ty) node [below right] {\picfontsize $\In_\nu(\scrP_1)$};    	
	
	\begin{scope}[shift={(\shifttwo,0)}]
		\draw [gray,  line width=0pt] (-0.5,-0.5) grid (\gridx,\gridy);
		
		\draw[ultra thick, blue] (0,0) -- (1,0);
		\draw (\tx,\ty) node [below right] {\picfontsize $\In_\nu(\scrP_3)$};    	
	\end{scope}
	
	\begin{scope}[shift={(2*\shifttwo,0)}]
		\draw [gray,  line width=0pt] (-0.5,-0.5) grid (\gridxx,\gridy);
		
		\draw[fill=\colorone, opacity=\opazero, ultra thick] (0,0) -- (1,0) -- (1,1) -- (0,1) --cycle;
		
		\draw (\tx,\ty) node [below right] {\picfontsize $\In_\nu(\scrP_1 + \scrP_3)$};    	
	\end{scope}
\end{scope}		
\end{scope}
\end{tikzpicture}
\caption{The image under $\psi_\nu$ of initial faces of $\scrP_1, \scrP_3$ and $\scrP_1 + \scrP_3$}  \label{fig-ex-bkk-1-2}
\end{center}
\end{figure}

The inner normals to the facets of $\scrP_1 + \scrP_2$ are listed in \cref{fig-ex-bkk-1-1}. The only element in $\Vprimeinfty$ is $(-2,-2,-1)$. If the origin is in $U$, then $\excU = \emptyset$, and therefore it follows from \cref{mixed-example,fig-ex-bkk-1-2} that
\def\scalefactor{0.3}
\def\colorzero{blue}
\def\colorone{green}
\begin{align*}
\multpisonodots{\scrA_1,\scrA_2,\scrA_3}{U}
	&= - \min_{\scrP_2}(-2,-2,-1) \cdot  	
	\area(
		\begin{tikzpicture}[scale=\scalefactor]
			\draw[fill=\colorone, opacity=\opazero] (2,0) -- (1,1) -- (2,1) -- (3,0) -- cycle;
			\draw[fill=\colorone, opacity=\opazero] (0,1) -- (1,1) -- (1,2) -- (0,2) -- cycle;
		\end{tikzpicture}
		)
	= 6 \cdot 2
	=12
\end{align*}
If the origin is not in $U$, but $U$ contains some other points of the $z$-axis, then for $I = \{1,2,3\}$, the set $\excUAprimeI$ contains the emptyset, but not $\{3\}$. It follows that $\Vprimess{\excUAprimeI}$ does not contain $(1,1,0)$, but it contains each (primitive) weighted order centered at the origin (see \cref{Kns-remark}); in particular, it contains $(2,2,1)$ and $(1,1,1)$. Since $\min_{\scrP_2}(\nu) = 0$ when $\nu = (1,0,0)$ or $(0,1,0)$, it does not matter if these two elements are in $\Vprimess{\excUAprimeI}$. It then follows from \cref{fig-ex-bkk-1-2} that
\begin{align*}
\multpisonodots{\scrA_1,\scrA_2,\scrA_3}{U}
	&= 12
		-  \min_{\scrP_2}(2,2,1) \cdot  	
			\area(\emptyset	)
	 	- \min_{\scrP_2}(1,1,1) \cdot  	
			\area(
				\begin{tikzpicture}[scale=\scalefactor]
					\draw[fill=\colorone, opacity=\opazero] (1,0) -- (0,1) -- (2,1) -- (3,0) -- cycle;
				\end{tikzpicture}
				)
	= 12 - 4 \cdot 0 - 3 \cdot 2
	= 6
\end{align*}
If $U$ does not intersect the $z$-axis, then $(1,1,0)$ is also an element of $\Vprimess{\excUAprimeI}$, and it follows that
\begin{align*}
\multpisonodots{\scrA_1,\scrA_2,\scrA_3}{U}
	&= 6
	 	- \min_{\scrP_2}(1,1,0) \cdot  	
			\area(
				\begin{tikzpicture}[scale=\scalefactor]
					\draw[fill=\colorone, opacity=\opazero] (0,0) -- (1,0) -- (1,1) -- (0,1) -- cycle;
				\end{tikzpicture}
				)
	= 6  -  1 \cdot 1
	= 5
\end{align*}
The computations therefore agree with those from \cref{ex-bkk-0}. Note that formulae \eqref{extended-bkk-formula-0} and \eqref{extended-bkk-formula} resolve the cases in the opposite order!
\end{example}

\section{Derivation of the formuale for the bound} \label{extended-bkk-bound-section}
In this section we prove \cref{extended-bkk-bound-0,extended-bkk-bound}. Throughout this section $\scrA := (\scrA_1, \ldots, \scrA_n)$ denotes a collection of finite subsets of $\znzero$, and $\scrP_j$ denotes the convex hull of $\scrA_j$ in $\rr^n$, $j  = 1, \ldots, n$. As in the preceding chapters, we write $\scrL(\scrA)$ for the space of $n$-tuples of polynomials supported at $\scrA$, and as in \cref{strong-defn}, we write $\Nstrong(\scrA)$ be the collection of all strongly $\scrA$-non-degenerate $(f_1, \ldots, f_n) \in \scrL(\scrA)$. Given a Zariski open subset $U$ of $\kk^n$, write $\Nstrong(U,\scrA)$ for all $(f_1, \ldots, f_n) \in \Nstrong(\scrA)$ such that for all $I \not\in \excU \cup \excA$,
\begin{align}
(V(f_1, \ldots, f_n) \cap \Kstari) \setminus U = \emptyset \label{strong-U}
\end{align}
\Cref{A-non-degeneracy-existence,intersection-lemma} imply that $\Nstrong(U,\scrA)$ contains a nonempty Zariski open subset of $\scrL(\scrA)$. Assertion \eqref{openly-max} of \cref{mult-deformation-global} therefore implies that in order to prove \cref{extended-bkk-bound-0,extended-bkk-bound}, it suffices to show that $\multf{U}$ equals the quantities from the right hand sides of \eqref{extended-bkk-formula-0} and \eqref{extended-bkk-formula} for all $(f_1, \ldots, f_n) \in \Nstrong(U,\scrA)$.

\subsection{Proof of \cref{extended-bkk-bound-0}} \label{extended-bkk-bound-0-proof}

\Cref{extended-bkk-bound-0} follows from the following result.

\begin{prop}
Let $(f_1, \ldots, f_n) \in \Nstrong(U, \scrA)$ and $I \subseteq [n] \setminus (\excU \cup \excA)$.
\begin{enumerate}
\item \label{empty-I-assertion} If $I \not\in \touchUA$, then $V(f_1, \ldots, f_n) \cap \Kstari = \emptyset$.
\item \label{nonempty-I-assertion} If $I \in \touchUA$, then
\begin{align}
\sum_{a \in \Kstari} \multf{a}
	&= \mv(\scrP^I_{j_1}, \ldots, \scrP^I_{j_{|I|}})
							\times
							\multzero{\pi_{[n]\setminus I}(\scrP_{j'_1})}{\pi_{[n]\setminus I}(\scrP_{j'_{n-|I|}})}
	\label{bkk-formula-0-I}
\end{align}
where $j_1 , \ldots, j_{|I|}$ are elements of $\tiA$, and $j'_1, \ldots, j'_{n-|I|}$ are elements of $[n]\setminus \tiA$.
\end{enumerate}
\end{prop}

\begin{proof}
If $I \not\in \touchUA \cup \excA$, then $|\tiA| > |I|$, and property \ref{strongly-non-degenerate} of strongly $\scrA$-non-degenerate systems imply that $V(f_1, \ldots, f_n) \cap \Kstari = \emptyset$, which proves assertion \eqref{empty-I-assertion}. Now pick $I \in \touchUA$ and $a \in V(f_1, \ldots, f_n) \cap \Kstari$. Since proper non-degeneracy implies BKK non-degeneracy when dimension of the ambient affine space is equal to the number of polynomials, properties \ref{strongly-non-degenerate} and \ref{strong-np} of strongly $\scrA$-non-degenerate systems imply that $(f_j|_{\Ki}: j \in \tiA)$ are $(\Ai_j: j \in \tiA)$-non-degenerate, and \cref{finite-cor} implies that $ V(f_1, \ldots, f_n) \cap \Kstari$ is finite. If $j \in \tiA$, \cref{thm:pure-dimension} implies that $C := V(f_i: i\neq j) \cap \Kstari$ is purely one dimensional near $a$. Property \ref{intersectionally-non-degenerate} of strongly $\scrA$-non-degenerate systems and \cref{strong-cor} imply that $\Ki$ is an irreducible component of $V' := V(f_{j'}: j' \not\in \tiA)$, and no irreducible component of $V'$ other than $\Ki$ contains any irreducible component of $C$. \Cref{mult-chain-1} then implies that
\begin{align}
\multp{f_1}{f_n}{a} = \multp{f_{j_1}|_{\Ki}}{f_{j_{|I|}}|_{\Ki}}{a} \times \multzero{f_{j'_1, \epsilon}}
 	{f_{j'_{n-|I|}, \epsilon}} \label{bkk-formula-0-I-a}
\end{align}
for generic $\epsilon = (\epsilon_1, \ldots, \epsilon_n) \in \Kstari$, where $f_{j'_k, \epsilon}$ are formed from $f_{j'_k}$ by substituting $\epsilon_i$ for $x_i$ for all $i \in I$. Due to property \eqref{strong-U}, assertion \eqref{nonempty-I-assertion} follows by summing \eqref{bkk-formula-0-I-a} over all $a \in \Kstari$ due to \cref{bkk-non-degenerate-thm,non-degeneracy-0-thm} (after using $(\Ai_{j_1}, \ldots, \Ai_{j_{|I|}})$-non-degeneracy of $f_{j_1}|_{\Ki}, \ldots,  f_{j_{|I|}}|_{\Ki}$ and non-degeneracy at the origin of $f_{j'_1, \epsilon}, \ldots, f_{j'_{n-|I|}, \epsilon}$ for generic $\epsilon$).
\end{proof}

\subsection{Proof of \cref{extended-bkk-bound}}
Let $(f_1, \ldots, f_n) \in \Nstrong(U, \scrA)$ and $C$ be the closed subscheme of $U' := U \setminus \bigcup_{I\in \excA}\Ki$. defined by $f_2, \ldots, f_n$  For each $I \subseteq [n]$ such that $U' \cap \Ki \neq \emptyset$, let $\{C^I_j\}_j$ be the set of irreducible components of $C$ such that $\Ki$ is the smallest coordinate subspace of $\kk^n$ containing each $C^I_j$.

\begin{claim} \label{C^I_j}
Let $\touchUAone$ be as in \cref{extended-bkk-bound}, $I \subseteq [n]$ such that $U' \cap \Ki \neq \emptyset$, and $\ziA$ be as in \cref{strong-defn}.
\begin{enumerate}
\item If $\{C^I_j\}_j$ is nonempty, then $I \in \touchUAone$.
\item If $I \in \touchUAone$, then $[n] \setminus \tiA \in \ziA$.
\end{enumerate}
\end{claim}

\begin{proof}
For the first assertion, pick $I \subseteq [n]$ such that $\{C^I_j\}_j$ is nonempty. Since $I \not\in \excA$, it follows that $|\tiA| \geq |I|$. On the other hand, if $|\tiA \setminus \{1\}| \geq |I|$, then property \ref{strongly-non-degenerate} of strong $\scrA$-non-degeneracy and \cref{branch-lemma-IB} implies that each $\{C^I_J\}_j$ is a point, which contradicts \cref{thm:pure-dimension}. Accordingly $|\tiA \setminus \{1\}| = |I| - 1$ and $|\tiA| = |I|$, which imply that $I \in \touchUAone$, as required. The second assertion follows from the definition of $\excA$ and \cref{mult-support}.
\end{proof}

For each $I \in \touchUAone$, property \ref{strongly-non-degenerate} of strong $\scrA$-non-degeneracy and assertion \eqref{2-n-dim<n-1} of \cref{order-sum-nktorus} imply that either $\{C^I_j\}_j$ is empty, or each $C^I_j$ has dimension one. It then follows from \cref{C^I_j} that $C$ is a curve, and therefore \cref{cor:Macaulay,int-mult-curve,order-curve} imply that
\begin{align*}
\multf{U}
	&= \sum_{a \in C} \ord_a(f_1|_C) = \sum_{I,j,a} \multp{f_2}{f_n}{C^I_j}\ord_a(f_1|_{C^I_j})
\end{align*}
where $\multp{f_2}{f_n}{C^I_j}$ are defined as in \cref{complete-mult-section}. Pick $I \in \touchUAone$. Let $j_1 = 1, j_2, \ldots, j_{|I|}$ be the elements of $\tiA$, and $j'_1, \ldots, j'_{n- |I|}$ be the elements of $ [n] \setminus \tiA$. Property \ref{intersectionally-non-degenerate} of strongly $\scrA$-non-degenerate systems and \cref{strong-cor} imply that $\Ki$ is an irreducible component of $V' := V(f_{j'_1}, \ldots, f_{j'_{n-|I|}})$, and if $\{C^I_j\}_j$ is nonempty, then no irreducible component of $V'$ other than $\Ki$ contains any irreducible component of $C$. \Cref{mult-chain-0} then implies that
\begin{align*}
\multf{U}
	&= \sum_{I,j,a} \ord_a(f_1|_{C^I_j}) \multp{f_{j_2}|_{\Ki}}{f_{j_{|I|}}|_{\Ki}}{C^I_j} \times \multzero{f_{j'_1, \epsilon}}
		 	{f_{j'_{n-|I|}, \epsilon}}
\end{align*}
for generic $\epsilon = (\epsilon_1, \ldots, \epsilon_n) \in \Kstari$, where $f_{j'_k, \epsilon}$ are formed from $f_{j'_k}$ by substituting $\epsilon_i$ for $x_i$ for all $i \in I$. Properties \ref{recursively-non-degenerate}, \ref{strong-np} of strongly $\scrA$-non-degenerate systems and \cref{non-degeneracy-0-thm} imply that $\multzero{f_{j'_1, \epsilon}}{f_{j'_{n-|I|}, \epsilon}}
	= \multzero{\pi_{[n]\setminus I}(\scrP_{j'_1})}{\pi_{[n]\setminus I}(\scrP_{j'_{n-|I|}})}$ for generic $\epsilon \in \Ki$. In order to prove \cref{extended-bkk-bound} therefore it suffices to show that for each $I \in \touchUA$,
\begin{align*}
\sum_{j,a} \ord_a(f_1|_{C^I_j}) \multp{f_{j_2}|_{\Ki}}{f_{j_{|I|}}|_{\Ki}}{C^I_j}
	= 	\multpstar{\scrP^{I}_1, \scrP^{I}_{j_2}}{\scrP^{I}_{j_{|I|}}}{\excUAI}
\end{align*}
Identify $\Ki$ with $\kk^k$, where $k := |I|$. \Cref{order-nu-prop,C-Q} imply that we can find a $k$-dimensional convex rational polytope $\scrP$ such that
\begin{itemize}
\item $\Ki \into \xp$,
\item $V(f_{j_2}|_{\Ki}, \ldots, f_{j_k}|_{\Ki})$ extends to a complete curve $\bar C^I$ on $\xp$, and
\item $f_1$ restricts to a nonzero rational function on $\bar C^I$ which is representable near every point of $\bar C^I$ as a quotient of non zero-divisors.
\end{itemize}
\Cref{zero-sum} implies that
\begin{align*}
\sum_j \sum_{a \in U'} \ord_a(f_1|_{C^I_j}) \multp{f_{j_2}|_{\Ki}}{f_{j_k}|_{\Ki}}{C^I_j}
	= 	- \sum_j \sum_{a \in \bar C^I_j \setminus U'}  \ord_a(f_1|_{C^I_j}) \multp{f_{j_2}|_{\Ki}}{f_{j_k}|_{\Ki}}{C^I_j}
\end{align*}
where $\bar C^I_j$ are the closures of $C^I_j$ in $\xp$. Property \eqref{strong-U} of the $f_j$ implies that $a \in \bar C^I_j \setminus U'$ if and only if either $a \in \Kii{I'}$ for some $I' \in \excU \cup \excA$, or (the germ at) $a$ is a branch at infinity of $C^I_j$. \Cref{order-nu-cor} then implies that
\begin{align*}
\sum_j \sum_{a \in U'} \ord_a(f_1|_{C^I_j}) \multp{f_{j_2}|_{\Ki}}{f_{j_k}|_{\Ki}}{C^I_j}
	= 	- \sum_{\nu \in \Viprimess{\excUAI} \cup \Viprimeinfty} \min_{\scrP^I_1}(\nu) ~
			\mv'_\nu(\In_\nu(\scrP^I_{j_2}), \ldots, \In_\nu(\scrP^I_{j_k}))
\end{align*}
where $\Viprimeinfty$ (respectively $\Viprimess{\excUAI}$) is the set of primitive weighted orders on $\kk[x_i: i \in I]$ which are centered at infinity (respectively at $\Kii{I'}$ for some $I' \in \excUAI$). Since the right hand side of the preceding identity is precisely $\multpstar{\scrP^{I}_1, \scrP^{I}_{j_2}}{\scrP^{I}_{j_{|I|}}}{\excUAI}$, this completes the proof of \cref{extended-bkk-bound}.

\section{Other formulae for the bound} \label{other-section}
Throughout this section we continue to use $\scrA $ to denote a collection $(\scrA_1, \ldots, \scrA_n)$ of finite subsets of $\znzero$ and $\scrP_j$ to denote the convex hull of $\scrA_j$, $j = 1, \ldots, n$.

\subsection{The formula of Huber-Sturmfels and Rojas} \label{hsr-affine}
Let $t$ be a new indeterminate. Fix positive integers $k_1, \ldots, k_n$. Note that for each $f_1, \ldots, f_n \in \kk[x_1, \ldots, x_n]$, and each Zariski open subset $U$ of $\kk^n$,
\begin{align}
\multfiso{U}
	&= \multpisonodots{t, f_1 + c_1t^{k_1}, \ldots, f_n + c_nt^{k_n}}{U \times \kk} \label{lifted-system}
\end{align}
for any $c_1, \ldots, c_n \in \kk$. It follows that
\begin{align*}
\multAiso{U}
	&= \multpiso{\hat \scrA_0}{\hat \scrA_n}{\hat U}
\end{align*}
where $\hat U := U \times \kk$, $\hat \scrA_0 := \{(1, 0, \ldots, 0)\} \subset \zzeroo{n+1}$ and $\hat \scrA_j :=  \{(k_j, 0, \ldots, 0)\} \cup (\{0\} \times \scrA_j) \subset \zzeroo{n+1}$  for $j = 1, \ldots, n$. Let $\hat \scrA := (\hat \scrA_0, \ldots, \hat \scrA_n)$. It is straightforward to check that $\touchUUAAone{\hat U}{\hat \scrA} = \{[n+1]\}$, so that \cref{extended-bkk-bound} implies that
\begin{align*}
\multAiso{U}
	= \multpstar{\hat \scrP_0}{\hat\scrP_n}{\exc( \hat U)\cup \exc(\hat \scrA)}
	= - \sum_{\hat \nu \in \Vhatprimess{\exc(\hat U) \cup \exc(\hat \scrA)} \cup \Vhatprimeinfty}	\hat \nu_0 \mv'_{\hat \nu}(\In_{\hat \nu}(\hat \scrP_1), \ldots, \In_{\hat \nu} (\hat \scrP_n))
\end{align*}
where $\hat \scrP_j$ are the convex hulls of $\hat \scrA_j$, and $\hat \nu$ ranges over the collection $\Vhatprimess{\exc(\hat U) \cup \exc(\hat \scrA)} \cup \Vhatprimeinfty$ of all primitive weighted orders on $\kk[t, x_1, \ldots, x_n]$ which are either centered at infinity or centered at $\Ki$ for some $I \in \exc( \hat U)\cup \exc(\hat \scrA)$, and $\hat \nu_0 := \hat \nu(t)$. Now, since $\dim(\hat P_0) = 0$, either \cref{bkk-bound-thm} or \cref{positively-mixed} implies that $\mv(\hat \scrP_0, \ldots, \hat \scrP_n) = 0$. Therefore assertion \eqref{rational-mv'} of \cref{rational-mv'-prop} implies that
\begin{align}
\multAiso{U}
&= \sum_{\hat \nu \not\in \Vhatprimess{\exc(\hat U) \cup \exc(\hat \scrA)} \cup \Vhatprimeinfty}	
	\hat \nu_0 \mv'_{\hat \nu}(\In_{\hat \nu}(\hat \scrP_1), \ldots, \In_{\hat \nu} (\hat \scrP_n))
\label{extended-bkk-formula-rojas-0}
\end{align}
If $\nu$ is a primitive weighted order on $\kk[t, x_1, \ldots, x_n]$, then $\nu \not\in \Vhatprimess{\exc(\hat U) \cup \exc(\hat \scrA)} \cup \Vhatprimeinfty$ if and only if all the following hold:
\begin{prooflist}
\item \label{non-neg} $\nu$ is nonnegative, and
\item \label{exc'} for each $I \in \excU \cup \excA$, there is $i' \not\in I$ such that $\nu(x_i) = 0$.
\end{prooflist}
Let $\Vhatprime(U,\scrA)$ be the set of all primitive weighted orders $\hat \nu$ on $\kk[t, x_1, \ldots, x_n]$ which satisfy properties \ref{non-neg}, \ref{exc'} and in addition satisfy the following:
\begin{prooflist}[resume]
\item  $\hat \nu_0 := \hat \nu(t)$ is positive.
\end{prooflist}
Since a summand on the right hand side of \eqref{extended-bkk-formula-rojas-0} has a nonzero contribution only if $\hat \nu_0$ is positive, it follows that for any Zariski open subset $U$ of $\kk^n$,
\begin{align}
\multAiso{U}
&= \sum_{\hat \nu \in \Vhatprime(U,\scrA)}	
	\hat \nu_0 \mv'_{\hat \nu}(\In_{\hat \nu}(\hat \scrP_1), \ldots, \In_{\hat \nu} (\hat \scrP_n))
\label{extended-bkk-formula-rojas}
\end{align}
In the case that $\kk = \cc$, $\excA = \emptyset$, and $U = \kk^n \setminus V( \prod_{j \in J} x_j) \cong \nktoruss{|J|} \times \kk^{n-|J|}$ for some $J \subseteq [n]$, formula \eqref{extended-bkk-formula-rojas} appeared in \cite{hurmfels-bern}. In this case $\Vhatprime(U, \scrA)$ consists of all primitive nonnegative weighted orders $\hat \nu$ on $\kk[t, x_1, \ldots, x_n]$ such that $\hat \nu_0$ is positive, and $\hat \nu(x_j) = 0$ for each $j \in J$. The sum on the right hand side of \eqref{extended-bkk-formula-rojas} in this case was termed in \cite{hurmfels-bern} as the {\em $I$-stable mixed volume} (where $I := [n]\setminus J$) of $\scrA_1, \ldots, \scrA_n$. J.\ M.\ Rojas \cite{rojas-toric} observed that the formula of \cite{hurmfels-bern} works over all algebraically closed fields.

\subsection{Estimates in terms of single mixed volumes} \label{li-wang-section}
If $U$ is nonempty, identity \eqref{extended-bkk-formula-0} implies that $\multAiso{U} \geq \mv(\scrP_1, \ldots, \scrP_n)$. On the other hand, since $\scrA_j \subseteq \scrA'_j := \scrA_j \cup \{\origin\}$, it trivially follows that $\multAiso{U} \leq \multAprimeiso{U} = \mv(\scrP'_1, \ldots, \scrP'_n)$, where $\scrP'_j$ are the convex hull of $\scrA'_j$, and the last equality follows from \eqref{extended-bkk-formula-0}. It follows that for nonempty Zariski open subsets $U$ of $\kk^n$,
\begin{align}
\mv(\conv(\scrA_1), \ldots, \conv(\scrA_n))
	\leq \multAiso{U}
	\leq \mv(\conv(\scrA_1 \cup \{\origin\}) \ldots, \conv(\scrA_n \cup \{\origin\}))
	\label{mv-estimate}
\end{align}
The upper bound in \eqref{mv-estimate} is due to T.\ Y.\ Li and X.\ Wang \cite{li-wang-bkk}. We now examine when these bounds are exact. The lower bound is easier to handle; the following result follows directly from \cref{positively-mixed,extended-bkk-bound-0}.

\begin{prop} \label{lower-prop}
Let $U$ be a nonempty Zariski open subset of $\kk^n$. Then the following are equivalent:
\begin{enumerate}
\item  The first inequality in \eqref{mv-estimate} holds with equality.
\item \label{lower-condition} For each $I \in \touchUA\setminus \{[n]\}$, $\scrP_{j_1} \cap \ri, \ldots, \scrP_{j_{|I|}} \cap \ri$ are {\em dependent}, where $j_1, \ldots, j_{|I|}$ are the elements of $\tiA$. \qed
\end{enumerate}
\end{prop}

\begin{rem}
Since an empty collection of convex polytopes is by definition {\em independent}, condition \eqref{lower-condition} of \cref{lower-prop} implies in particular that $\emptyset \not\in \touchUA$, which in turn implies that either $U$ does not contain the origin, or $\multAzero$ is zero or $\infty$ (\cref{empty-touch}).
\end{rem}

Following A.\ Khovanskii \cite{khovanus}, We say that $\scrA = (\scrA_1, \ldots, \scrA_n)$ is {\em regularly attached to the coordinate cross} if for each proper subset $I$ of $[n]$, the set of nonempty elements of $\{\scrP_j \cap \ri: j = 1, \ldots, n\}$ is dependent; in particular this implies (taking $I = \emptyset$) that the origin belongs to at least one of the $\scrA_j$. The following is immediate from \cref{lower-prop}.

\begin{cor}[Khovanskii {\cite{khovanus}}]
If $U$ is nonempty and $\scrA$ is regularly attached to the coordinate cross, then $\multAiso{U} = \mv(\scrP_1, \ldots, \scrP_n)$. \qed
\end{cor}

Now let $M := \mv(\conv(\scrA_1 \cup \{\origin\}) \ldots, \conv(\scrA_n \cup \{\origin\}))$ be the upper bound from \eqref{mv-estimate}. Consider as in \eqref{lifted-system} the system
\begin{align}
f_j  + t = 0,\ j = 1, \ldots, n, \label{lifted-system'}
\end{align}
where $f_j$ are generic polynomials supported at $\scrA_j$, with $\np(f_j) = \conv(\scrA_j) = \scrP_j$. For generic $t \neq 0$, the corresponding system has precisely $M$ isolated solutions, all of which are on $\nktorus$. Therefore, the number of solutions of the system at $t = 0$ is also $M$ if and only if there is no curve of solutions of the system \eqref{lifted-system'} on $\nktorus$ that escapes $U$ or becomes non-isolated as $t$ approaches $0$. \Cref{bkk-bound-thm,order-nu-prop,positively-mixed} imply that such a curve exists if and only if there is a weighted order $\hat \nu$ on $\kk[x_1, \ldots, x_n, t]$ such that
\begin{prooflist}
\item $\hat \nu(t) > 0$,
\item the restriction of $\hat \nu$ to $\kk[x_1, \ldots, x_n]$ is either centered at infinity or at $\Ki$ for some $I \in \excU \cup \excA$, and
\item $\In_{\hat \nu}(\hat \scrP_1), \ldots,  \In_{\hat \nu}(\hat \scrP_n)$ are {\em independent}, where $\hat \scrP_j := \np(f_j + t) \subset \rzeroo{n+1}$.
\end{prooflist}
Let $\nu$ be the restriction of $\hat \nu$ to $\kk[x_1, \ldots, x_n]$. Let $m_\nu := \max\{\min_{\scrP_j}(\nu): j = 1, \ldots, n\}$. Since $\dim( \In_{\hat \nu}(\sum_j \hat \scrP_j)) \leq \dim( \In_{\nu}(\sum_j \scrP_j)) + 1$, and since $\hat \nu(t) > 0$, it follows from the definition of dependence of polytopes that $\In_{\hat \nu}(\hat \scrP_1), \ldots,  \In_{\hat \nu}(\hat \scrP_n)$ are independent if and only if
\begin{prooflist}[resume]
\item  \label{m_nu-condition}$m_\nu > 0$, and
\item \label{dimension-condition} for all $J \subseteq [n]$,
\begin{align}
\dim(\In_\nu(\sum_{j \in J} \scrP_j))
	\geq
	\begin{cases}
	|J| &\text{if $\min_{\scrP_j}(\nu) < m_\nu$ for each $j \in J$,}\\
	|J| - 1 & \text{if there is $j \in J$ such that $\min_{\scrP_j}(\nu) = m_\nu$.}
	\end{cases}
	\label{li-wang-dimension-constraint}
\end{align}
\end{prooflist}
Taking $J = [n]$ in \eqref{li-wang-dimension-constraint} implies in particular that $\dim(\In_\nu(\sum_j \scrP_j)) = n-1$. These observations imply the following result.

\begin{prop} \label{li-wang-condition}
Let $\scrP := \sum_j \scrP_j$. The second inequality in \eqref{mv-estimate} holds with equality if and only if
\begin{enumerate}
\item \label{leq-n-2} either $\dim(\scrP) \leq n -2$, or
\item $\dim(\scrP) \geq n -1$, and for each face of dimension $(n-1)$ of $\scrP$ such that its primitive inner normal\footnote{If $\dim(\scrP) = n-1$, then both of the primitive normals to $\scrP$ are considered to be inner.} $\nu$ is centered at infinity or at $\Ki$ for some $I \in \excU \cup \excA$, at least one of the conditions \ref{m_nu-condition} and \ref{dimension-condition} fails for $\nu$. \qed
\end{enumerate}
\end{prop}

As a corollary we get a situation where {\em both} of the bounds of \eqref{mv-estimate} are exact.

\begin{cor} \label{conveniently-li-wang}
Assume each $\scrA_j$ is {\em convenient} and $U$ is a nonempty Zariski open subset of $\kk^n$. Assume in addition that at least one of the following holds:
\begin{enumerate}
\item either $U$ does not contain the origin,
\item or at least one of the $\scrA_j$ contains the origin.
\end{enumerate}
Then
\begin{align*}
 \multAiso{U}
 	= \mv(\conv(\scrA_1), \ldots, \conv(\scrA_n))
	=  \mv(\conv(\scrA_1 \cup \{\origin\}) \ldots, \conv(\scrA_n \cup \{\origin\}))
\end{align*}
\end{cor}

\begin{proof}
$\scrA_j$ satisfy the hypotheses of both \cref{lower-prop,li-wang-condition}.
\end{proof}

%

\section{Examples motivating the non-degeneracy conditions} \label{non-deg-example-section}


In this section we give some examples to motivate the necessary and sufficient conditions for the equality $\multfkniso = \multAkniso$, where $\scrA := (\scrA_1, \ldots, \scrA_n)$ is a given collection of finite subsets of $\znzero$, and $(f_1, \ldots, f_n) \in \scrL(\scrA)$. We consider the case $U = \kk^n$; in a sense this is the most important case, and it already captures the essence of the general case. We also assume for simplicity that $\excA = \emptyset$, which ensures in particular that for generic $(f_1, \ldots, f_n) \in \scrL(\scrA)$, all points in the set $V(f_1, \ldots, f_n)$ of common zeroes of $f_1, \ldots, f_n$ in $\kk^n$ are isolated. In this scenario, if we apply the intuitive reasoning from \cref{sec:bkk-non-degeneracy-statement} that motivated the non-degeneracy condition for Bernstein's theorem, we are led to the following condition:\index{Non-degeneracy!centered at infinity}
\begin{align}
\parbox{0.69\textwidth}{$f_1|_{\Ki}, \ldots, f_n|_{\Ki}$ are {\em $\Ai$-non-degenerate at infinity} for each $I \subseteq [n]$.}
\tag{ND$_\infty$}
\label{non-deg-infinity}
\end{align}
where $\Ai := (\Ai_1, \ldots, \Ai_n) = (\scrA_1 \cap \ri, \ldots, \scrA_n \cap \ri)$, and {\em non-degeneracy at infinity} is defined as follows: given a collection $\scrB := (\scrB_1, \ldots, \scrB_m)$ of finite subsets of $\znzero$ and $g_j \in \kk[x_1, \ldots, x_n]$, such that $\supp(g_j) \subseteq \scrB_j$, $j = 1, \ldots, m$, we say that $g_1, \ldots, g_m$ are $\scrB$- non-degenerate at infinity if for each weighted order $\nu$ {\em centered at infinity} (see \cref{second-section}), there is no common root of $\In_{\scrB_j,\nu}(g_j)$, $j = 1, \ldots, m$, on $\nktorus$. We now present a series of examples which illustrate how condition \eqref{non-deg-infinity} falls short of characterizing the correct non-degeneracy condition, and which also suggest the ways to amend it. In all these examples $\scrP_j$ would denote the convex hull of $\scrA_j$, $j = 1, \ldots, n$.

\def\shiftx{9}
\def\colorzero{blue}
\def\opazero{0.5}
\def\viewx{70}
\def\viewy{5}
\def\titlex{3}
\def\titley{-1}
\def\scalefactor{0.6}
\def\colorzero{blue}
\def\colorone{green}

\begin{figure}[h]
\begin{center}
\begin{tikzpicture}[scale=0.6]
\pgfplotsset{every axis/.append style = {view={\viewx}{\viewy}, axis lines=middle, enlargelimits={upper}, axis equal}}

\begin{scope}
\begin{axis}
	\addplot3 [fill=blue,opacity=\opazero, thick] coordinates{(0,1,1) (1,1,0) (1,0,1)};
	\addplot3 [fill=red,opacity=\opazero, thick] coordinates{(0,0,0) (1,1,0) (1,0,1)};
	\addplot3 [thick, dashed] coordinates{(0,0,0) (0,1,1)};
\end{axis}
\draw (\titlex,\titley) node {$\scrP_1 = \scrP_2$};
\end{scope}

\begin{scope}[shift={(\shiftx,0)}]
\begin{axis}
	\addplot3 [fill=blue,opacity=\opazero, thick] coordinates{(0,1,2) (1,1,1) (1,0,2)};
	\addplot3 [fill=red,opacity=\opazero, thick] coordinates{(0,0,1) (1,1,1) (1,0,2)};
	\addplot3 [thick, dashed] coordinates{(0,0,1) (0,1,2)};
\end{axis}
\draw (\titlex,\titley) node {$\scrP_3 = \scrP_1 + \{(0,0,1)\}$};
\end{scope}


\end{tikzpicture}
\end{center}
\caption{Newton polytopes of \cref{ex-aff-non-deg-1}} \label{fig-ex-aff-non-deg-1}
\end{figure}

\begin{example}\label{ex-aff-non-deg-1}
Let $\scrA_1 = \scrA_2 = \{(0,0,0), (1,1,0), (0, 1, 1), (1,0,1)\} \subset \zzeroo{3}$ and $\scrA_3 = \scrA_1 + \{(0,0,1)\}$.  Then $\scrP_1 = \scrP_2$ is a tetrahedron and $\scrP_3$ is a translation of $\scrP_1$, and therefore
\begin{align*}
\mv(\scrP_1, \scrP_2, \scrP_3) = 3! \vol(\scrP_1) = 2
\end{align*}
If $f_j$ are polynomials such that $\np(f_j) = \scrP_j$, then
\begin{align*}
f_1 &= a_1 + b_1x_1x_2 + c_1x_2x_3 + d_1x_3x_1\\
f_2 &= a_2 + b_2x_1x_2+ c_2x_2x_3  + d_2x_3x_1 \\
f_3 &= x_3(a_3 + b_3x_1x_2 + c_3x_2x_3  + d_3x_3x_1)
\end{align*}
where $a_j, b_j, c_j, d_j \in \ktorus$. We write $V$ for the set of common zeroes of $f_1, \ldots, f_3$ on $\kk^3$, and $V^*$ for $V \cap \nktoruss{3}$.
\begin{defnlist}
\item \label{aff-non-deg-1-0} If all $a_j, b_j, c_j, d_j$ are generic, then it is straightforward to check directly that $V = V^*$. Therefore \cref{bkk-bound-thm} implies that $\multpisonodots{\scrA_1, \scrA_2, \scrA_3}{\kk^3}  = \multpisonodots{f_1, f_2, f_3}{\kk^3}  = \mv(\scrP_1, \scrP_2, \scrP_3) = 2$.
\item \label{aff-non-deg-1-1} If $a_1 = a_2$, $b_1 = b_2$, and the remaining coefficients are generic, then $V = V^* \cup C$, where $C := \{x_3= a_1 + b_1x_1x_2 = 0\}$ is a positive dimensional component of $V(f_1, f_2, f_3)$. However, $f_1, f_2, f_3$ still satisfy \eqref{b-non-degeneracy}, and \cref{bkk-non-degenerate-thm} implies that $\multpisonodots{f_1, f_2, f_3}{\kk^3} = 2 = \multpisonodots{\scrA_1, \scrA_2, \scrA_3}{\kk^3}$.
\item \label{aff-non-deg-1-2} If $a_1 = a_2 = a_3$, $b_1 = b_2 = b_3$, and the rest of the coefficients are generic, then again $V = V^* \cup C$. However, \eqref{b-non-degeneracy} fails for the weighted order $\nu$ corresponding to weights $(-1,1,2)$ for $(x,y,z)$, and \cref{bkk-non-degenerate-thm} implies that $\multpisonodots{f_1, f_2, f_3}{\kk^3} < 2 = \multpisonodots{\scrA_1, \scrA_2, \scrA_3}{\kk^3}$. (It is straightforward to verify directly that in this case $V^* = \emptyset$ and $\multpisonodots{f_1, f_2, f_3}{\kk^3} = 0$.)
\end{defnlist}
\end{example}

Part \ref{aff-non-deg-1-1} of \cref{ex-aff-non-deg-1} shows that it is possible that $V(f_1, \ldots, f_n)$ has a positive dimensional component on $\kk^n$, but still $\multfkniso = \multAkniso$ (where $\scrA_j = \supp(f_j)$). (This does not happen in the case of $\nktorus$, see \cref{positive-dimensional-prop}.) Moreover, since the intersection of the curve $C$ with the ``torus'' of the $(x_1,x_2)$-plane is nonempty, in part \ref{aff-non-deg-1-1} of \cref{ex-aff-non-deg-1}, condition \eqref{non-deg-infinity} is violated for $I = \{1,2\}$. However, note that the intersections of $\scrP_1$ and $\scrP_2$ with the $(x_1,x_2)$-plane (in $\rr^n$) are {\em dependent} in the terminology of \cref{dependent-defn}. Moreover, in part \ref{aff-non-deg-1-2} of \cref{ex-aff-non-deg-1}, where $\multpisonodots{f_1, f_2, f_3}{\kk^3} < \multpisonodots{\scrA_1, \scrA_2, \scrA_3}{\kk^3}$, condition \eqref{non-deg-infinity} is violated with $I = \{1, 2, 3\}$, and the corresponding polytopes are {\em independent}. This motivates the following definition.

\begin{defn} \label{ri-dependent-defn}
An ordered collection $\scrB = (\scrB_1, \ldots, \scrB_m)$, $m \geq 1$, of collections of finite subsets of $\rr^n$ is called {\em $\ri$-dependent} if there is a nonempty subset $J$ of $[m]$ such that $\scrB^I_j := \scrB_j \cap \ri$ is nonempty for each $j \in J$, and the collection $\{\conv(\scrB^I_j): j \in J\}$ of convex polytopes is \index{Dependence!of finite subsets of $\rr^n$}dependent (\cref{dependent-defn}); otherwise we say that $\scrB$ is \index{Independence!of finite subsets of $\rr^n$}{\em $\ri$-independent}.
\end{defn}

\Cref{ex-aff-non-deg-1} suggests that
\begin{prooflist}
\item \label{non-deg-ind-infinity} Condition \eqref{non-deg-infinity} should be checked only for those $I \subseteq [n]$ such that $\scrA$ is $\ri$-independent.
\end{prooflist}
This, however, is not enough, as the next example shows.

\def\shiftx{9}
\def\colorzero{blue}
\def\opazero{0.5}
\def\viewx{70}
\def\viewy{5}
\def\titlex{2}
\def\titley{-1}
\def\scalefactor{0.6}
\def\colorzero{blue}
\def\colorone{green}

\begin{figure}[h]
\begin{center}
\begin{tikzpicture}[scale=0.6]
\pgfplotsset{every axis/.append style = {view={\viewx}{\viewy}, axis lines=middle, enlargelimits={upper}}}

\begin{scope}
\begin{axis}
	\addplot3 [fill=blue,opacity=\opazero, thick] coordinates{(0,1,1) (1,0,0) (1,0,1)};
	\addplot3 [fill=red,opacity=\opazero, thick] coordinates{(0,0,0) (1,0,0) (1,0,1)};
	\addplot3 [thick, dashed] coordinates{(0,0,0) (0,1,1)};
\end{axis}
\draw (\titlex,\titley) node {$\scrP_1 = \scrP_2$};
\end{scope}

\begin{scope}[shift={(\shiftx,0)}]
\begin{axis}
	\addplot3 [fill=red,opacity=\opazero, thick] coordinates{(0,0,1) (0,1,1) (1,0,2)};
	\addplot3 [fill=blue,opacity=\opazero, thick] coordinates{(0,1,1) (1,0,2) (0,1,2)};
	\addplot3 [thick, dashed] coordinates{(0,0,1) (0,1,2)};
\end{axis}
\draw (\titlex,\titley) node {$\scrP_3$};
\end{scope}

\def\factorone{1}
\def\factortwo{1}
\def\factorthree{1}
\def\factorfour{1}
\begin{scope}[shift={(2*\shiftone,0)}]
\begin{axis}[hide axis, enlargelimits={0.15}]
	\addplot3 [ultra thick, red, ->] coordinates{(0,0,0) (0,-\factorone,0)};
	\addplot3 [ultra thick, blue, ->] coordinates{(0,0,0) (\factortwo, \factortwo,0)};
	\addplot3 [ultra thick, ->] coordinates{(0,0,0) (-\factorthree,-\factorthree,\factorthree)};
	\addplot3 [ultra thick, ->] coordinates{(0,0,0) (0,\factorfour,-\factorfour)};

	\draw (axis cs:0,-\factorone,0) node [below] {(0,-1,0)};
	\draw (axis cs: \factortwo,\factortwo,0) node [above] {(1,1,0)};
	\draw (axis cs:-\factorthree,-\factorthree,\factorthree) node [right] {(-1,-1,1)};
	\draw (axis cs: 0,\factorfour,-\factorfour) node [below] {(0,1,-1)};
\end{axis}
\draw (\titlex,\titley) node [right, text width=3.5cm] {outer normals to facets of $\scrP_1$};
\end{scope}
\end{tikzpicture}
\end{center}
\caption{Newton polytopes of \cref{ex-aff-non-deg-2}} \label{fig-ex-aff-non-deg-2}
\end{figure}

\begin{example}\label{ex-aff-non-deg-2}
Consider the following system of polynomials:
\begin{align*}
f_1 &= a_1 + b_1x_1 + c_1x_2x_3 + d_1x_3x_1\\
f_2 &= a_2 + b_2x_1 + c_2x_2x_3  + d_2x_3x_1 \\
f_3 &= x_3(a_3 + b_3x_2 + c_3x_2x_3  + d_3x_3x_1)
\end{align*}
where $a_j,b_j, c_j, d_j \in \ktorus$, and $\scrA_j  := \supp(f_j)$, $j = 1, 2, 3$. We continue to write $V := V(f_1, f_2, f_3)$ and $V^* := V \cap \nktoruss{3}$. When all the coefficients are generic, then it is straightforward to check that $V = V^*$, so that \cref{bkk-bound-thm} implies that
\begin{align*}
\multpisonodots{f_1, f_2, f_3}{\kk^3}
	&= \multpisonodots{\scrA_1, \scrA_2, \scrA_3}{\kk^3}
	= \mv(\scrP_1, \scrP_2, \scrP_3)
\end{align*}
We compute $\mv(\scrP_1, \scrP_2, \scrP_3)$ using \cref{rational-mv'-prop}. If $\nu$ is the primitive outer normal to any of the facets of $\scrP_1 = \scrP_2$, it is straightforward to check that the image of the corresponding facet under the map $\psi_\nu$ from the definition of $\mv'_\nu(\cdot, \ldots, \cdot)$ is (up to a translation) the triangle $\scrT$ with vertices $(0,0),(0,1),(1,0)$, and therefore $\mv'_\nu(\ld_\nu(\scrP_1), \ld_\nu(\scrP_2)) = 2! \area(\scrT) = 1$.
It follows that
\begin{align*}
\multpisonodots{\scrA_1, \scrA_2, \scrA_3}{\kk^3}
	&= \mv(\scrP_1, \scrP_2, \scrP_3) \\
	&= \max_{\scrP_3}(0,-1,0) + \max_{\scrP_3}(0,1,-1) + \max_{\scrP_3}(1,1,0) + \max_{\scrP_3}(-1,-1,1) \\
	&= 0 + 0 + 1+ 1 \\
	& = 2
\end{align*}
If $a_1 = a_2$ and $b_1 = b_2$ and the other coefficients are generic, then it is straightforward to check directly that $V^* = \emptyset$ and $V$ is the curve $\{x_3 = a_1+b_1x_1 = 0\}$, so that $\multpisonodots{f_1, f_2, f_3}{\kk^3} = 0$. However, the only $I \subseteq \{1, 2,3\}$ such that $\scrA$ is $\ri$-independent is $I = \{1, 2, 3\}$, and it is straightforward to check that $f_j|_{\Ki} = f_j$ are in fact $\scrA$-non-degenerate at infinity. In particular, $f_1, f_2, f_3$ satisfy condition \ref{non-deg-ind-infinity}, but $\multpisonodots{f_1, f_2, f_3}{\kk^3} < \multpisonodots{\scrA_1, \scrA_2, \scrA_3}{\kk^3} $.
\end{example}

Given $I' \subseteq I \subseteq [n]$ and a weighted order $\nu$ on $\kk[x_i: i \in I]$, we say that $\nu$ is \index{Weighted!order!centered at $\Kstari$}centered at $\Kstariprime$ if $\nu(x_i) = 0$ for all $i \in I'$ and $\nu(x_i) > 0$ for all $i \in I\setminus I'$. It is straightforward to check that in \cref{ex-aff-non-deg-2}, the only nonzero weighted orders $\nu$ such that $\In_\nu(f_j)$, $j = 1, 2, 3$, have a common zero on $\nktoruss{3}$ are of the form $(0, 0, \epsilon)$ for $\epsilon > 0$, i.e.\ they are centered at $\Kstariprime$ with $I' := \{1,2\}$. It turns out that $\scrA$ is {\em hereditarily $\riprime$-dependent} (see \cref{non-deg-section}) and therefore \cref{ex-aff-non-deg-2} suggests that for $\multfkniso$ to be equal to $\multAkniso$, the following condition needs to be satisfied in addition to \ref{non-deg-ind-infinity}:
\begin{prooflist}[resume]
\item \label{non-deg-ind-finity} For each $I \subseteq [n]$ such that $\scrA$ is {\em not} hereditarily $\ri$-dependent, $\In_{\Ai_j, \nu}(f_j|_{\Ki})$ do not have any common zero on $\nktorus$ for all weighted orders $\nu$ on $\kk[x_i: i \in I]$ which are centered at $\Kstariprime$ for some $I' \subsetneq I$ such that $\scrA$ is hereditarily $\riprime$-dependent.
\end{prooflist}

\begin{example}\label{ex-aff-non-deg-3}
Let $f_1, \ldots, f_4$ be polynomials in $(x_1, \ldots, x_4)$ such that $f_1$ and $f_2$ are polynomials in $(x_1, x_2)$ with nonzero constant terms, the (two dimensional) mixed volume of $\np(f_1)$ and $\np(f_2)$ is nonzero, and
\begin{align*}
f_j &= x_3f_{j,1}(x_1, x_2) + x_4f_{j,2}(x_1,x_2) \in \kk[x_1, x_2, x_3, x_4],\ j = 3, 4,
\end{align*}
where $f_{3,1}, f_{3,2}, f_{4,1}, f_{4,2}$ are nonzero polynomials in $(x,y)$ such that $\np(f_{3,1})$ and $\np(f_{4,1})$ are positive dimensional. Let $\scrA_j$ be the support of $f_j$ and $\scrP_j$ be the Newton polytope of $f_j$, $j = 1, \ldots, 4$. The only $I \subseteq \{1, 2, 3, 4\}$ such that the supports of $f_j$ are $\ri$-independent is $I = \{1,2\}$. It follows that for generic coefficients, all the common zeroes of $f_1, \ldots, f_4$ are isolated and contained in the $(x_1,x_2)$-plane. Moreover, \cref{extended-bkk-bound-0} implies that
\begin{align*}
\multpisonodots{\scrA_1, \scrA_2, \scrA_3, \scrA_4}{\kk^4}
	= 	\mv(\scrP^{\{1,2\}}_1, \scrP^{\{1,2\}}_2) \times
		\multpnodots{\pi_{\{3,4\}}(\scrP_3),\pi_{\{3,4\}}(\scrP_4)}{\origin}
	= \mv(\scrP_1, \scrP_2)
\end{align*}
Now fix BKK non-degenerate $f_1, f_2$, and a common zero $z = (z_1, z_2)$ of $f_1,f_2$ on $\nktoruss{2}$. Take $f_{3,1}$ and $f_{4,1}$ such that $f_{3,1}(z) = f_{4,1}(z) = 0$. Then $\{(z_1,z_2, t, 0): t \in \kk\} \subseteq  V(f_1, f_2, f_3, f_4)$, so that $(z_1, z_2, 0,0)$ is no longer an isolated point of $V(f_1, f_2, f_3, f_4)$ (even though it is an isolated point of $V(f_1, f_2, f_3, f_4) \cap \Kii{\{1,2\}}$). It follows that $\multpisonodots{f_1, f_2,f_3, f_4}{\kk^4} <  \mv(\scrP_1, \scrP_2) = \multpisonodots{\scrA_1, \scrA_2, \scrA_3, \scrA_4}{\kk^4}$. However, $f_1, \ldots, f_4$ satisfy both conditions \ref{non-deg-ind-infinity} and \ref{non-deg-ind-finity}.
\end{example}

\Cref{ex-aff-non-deg-3} leads us to another condition that needs to be satisfied in addition to \ref{non-deg-ind-infinity} and \ref{non-deg-ind-finity} for $\multfkniso$ to be equal to $\multAkniso$:
\begin{prooflist}[resume]
\item \label{non-deg-dep} For each $I \subseteq [n]$ such that $\scrA$ is hereditarily $\ri$-dependent, $\In_{\Ai_j, \nu}(f_j|_{\Ki})$ do not have any common zero on $\nktorus$ for all weighted orders $\nu$ on $\kk[x_i: i \in I]$ which are centered at $\Kstariprime$ for some $I' \subsetneq I$ such that $\scrA$ is not hereditarily $\riprime$-dependent.
\end{prooflist}
Note that the coordinate subspaces in condition \ref{non-deg-dep} are in a sense ``dual'' to those in condition \ref{non-deg-ind-finity}.

\section{Non-degeneracy conditions} \label{non-deg-section}
Let $\scrA := (\scrA_1, \ldots, \scrA_n)$ be a collection of finite subsets of $\znzero$ and $I \subseteq [n]$. We say that $\scrA$ is \index{Hereditary dependence}\index{Dependence!hereditary}{\em hereditarily $\ri$-dependent} if $\scrA$ is $\ri$-dependent (see \cref{ri-dependent-defn}) and there is $I' \supseteq I$ such that
\begin{defnlist}
\item $\scrA$ is $\riprime$-dependent,
\item $|\tiprimeA| = |I'|$,
\item $|\titildeA| > |\tilde I|$ for each $\tilde I$ such that $I \subseteq \tilde I \subsetneq I'$.
\end{defnlist}

\begin{rem} \label{empty-remark}
If $I = \emptyset$, then $\ri$ is the origin, and $\scrA$ is $\ri$-dependent if and only if the origin belongs to some $\scrA_j$. However, even if $\scrA$ is $\ri$-dependent, it might not be hereditarily $\ri$-dependent, see \cref{empty-figure}.
\end{rem}

\begin{figure}[h]
\begin{center}

\begin{tikzpicture}[scale=0.45]
\def\shiftone{7.5}
\def\opazero{0.5}
\def\gridx{4.5}
\def\gridy{4.5}
\def\colorone{red}

\draw [gray,  line width=0pt] (-0.5,-0.5) grid (\gridx,\gridy);
\draw [<->] (0, \gridy) |- (\gridx, 0);
\draw[ultra thick, \colorone]  (3,0) -- (0,3) -- (0,0) -- cycle;
\fill[\colorzero, opacity=\opazero ] (3,0) -- (0,3) -- (0,0) -- cycle;
\def\tx{0.5}
\def\ty{0.5}
\draw (\tx,\ty) node {\picfontsize $\scrA_1$};

\begin{scope}[shift={(\shiftone,0)}]
	\draw [gray,  line width=0pt] (-0.5,-0.5) grid (\gridx,\gridy);
	\draw [<->] (0, \gridy) |- (\gridx, 0);
	\draw[ultra thick, \colorone]  (1,1) --  (4,1) -- (4,4) -- (1,4) -- cycle;
	\fill[\colorzero, opacity=\opazero ] (1,1) --  (4,1) -- (4,4) -- (1,4) -- cycle;
	\def\tx{2.5}
	\def\ty{2.5}
	\draw (\tx,\ty) node {\picfontsize $\scrA_2$};
\end{scope}

\begin{scope}[shift={(2*\shiftone,0)}]
	\draw [gray,  line width=0pt] (-0.5,-0.5) grid (\gridx,\gridy);
	\draw [<->] (0, \gridy) |- (\gridx, 0);
	\draw[ultra thick, \colorone]  (4,0) -- (1,4) -- (0,0) -- cycle;
	\fill[\colorzero, opacity=\opazero ] (4,0) -- (1,4) -- (0,0) -- cycle;
	\def\tx{1.5}
	\def\ty{1.5}
	\draw (\tx,\ty) node {\picfontsize $\scrA_3$};
\end{scope}
\end{tikzpicture}

\end{center}
\caption{Both $(\scrA_1, \scrA_2)$ and $(\scrA_2, \scrA_3)$ are $\rr^I$-dependent for $I = \emptyset$. However, only the latter pair is hereditarily $\rr^I$-dependent (with $I' = \{2\}$).} \label{empty-figure}
\end{figure}

The relevance of ``hereditary dependence'' to affine B\'ezout problem is given by \cref{hereditary-prop} below: it states that if $\scrA$ is hereditarily $\ri$-dependent, then for all $(f_1, \ldots, f_n) \in \scrL(\scrA)$, either $V^I := V(f_1, \ldots, f_n) \cap \Kstari$ is empty, or all points of $V^I$ are non-isolated in (the possibly larger set) $V(f_1, \ldots, f_n)$ (however, points of $V^I$ might be isolated in $V^I$ itself!). This is not necessarily true if $\scrA$ is simply $\ri$-dependent, e.g.\ the system $(x + y - 1, 2x - y - 2)$ (over a field of characteristic not equal to two) has an isolated zero on the coordinate subspace $y = 0$.

\begin{prop} \label{hereditary-prop}
If $\scrA$ is hereditarily $\ri$-dependent, then $V(f_1, \ldots, f_n)$ has no isolated point on $\Kstari$ for each $(f_1, \ldots, f_n) \in \scrL(\scrA)$.
\end{prop}

We prove \cref{hereditary-prop} in \cref{reducsection}. Now we introduce the correct non-degeneracy condition for (arbitrary open subsets of) the affine space. Let $U$ be a Zariski open subset of $\kk^n$. Let $\excA$ and $\excU$ be as in \cref{first-section}. Define
\begin{align}
\dUA &:= \{I \subseteq [n]: I \not \in \excU \cup \excA,\ \scrA\ \text{is hereditarily $\ri$-dependent}\}  \label{dUA}\\
\iUA &:= \{I \subseteq [n]: I \not \in \excU \cup \excA,\ \scrA\ \text{is not hereditarily $\ri$-dependent}\}  \label{iUA}
\end{align}
We say that polynomials $f_1, \ldots, f_n$ are \index{Non-degeneracy!with respect to an open subset of $\kk^n$ and tuples of finite subsets of $\znzero$}{\em $(U,\scrA)$-non-degenerate} if the following conditions are true:
\begin{defnlist}
\item \label{U-property} For each nonempty $I \in \iUA$, there is no common zero of $f_1|_{\Ki}, \ldots, f_n|_{\Ki}$ on $\Kstari \setminus U$ (note that this condition is vacuously true when $U = \kk^n$, or more generally, when $U = \Kns$ (see \eqref{Kns}) for some collection $\mscrS$ of subsets of $[n]$),
\item \label{iUA-property} For each nonempty $I \in \iUA$ and for each weighted order $\nu$ on $\kk[x_i: i \in I]$ such that $\nu$ is centered at infinity or at $\Kstariprime$ for some $I' \in \excU \cup \excA \cup \dUA$, there is no common zero of $\In_{\Ai_1, \nu}(f_1|_{\Ki}), \ldots, \In_{\Ai_n,\nu}(f_n|_{\Ki})$ on $\nktorus$.
\item \label{dUA-property}  For each nonempty $I \in \dUA$ and for each weighted order $\nu$ on $\kk[x_i: i \in I]$ such that $\nu$ is centered at $\Kstariprime$ for some $I' \in \iUA$, there is no common zero of $\In_{\Ai_1, \nu}(f_1|_{\Ki}), \ldots, \In_{\Ai_n,\nu}(f_n|_{\Ki})$ on $\nktorus$.	
\end{defnlist}

We prove the following result in \cref{non-deg-1-proof-section}.

\begin{thm}[{\cite{toricstein}}] \label{non-deg-1}
Assume $\multAiso{U} > 0$. Then for each $(f_1, \ldots, f_n) \in \scrL(\scrA)$, the following are equivalent:
\begin{enumerate}
\item \label{non-deg-1-0} $\multfiso{U} = \multAiso{U}$,
\item \label{non-deg-1-1}  $f_1, \ldots, f_n$ are $(U,\scrA)$-non-degenerate.
\end{enumerate}
The collection $\scrN(U, \scrA)$ of all $(U, \scrA)$-non-degenerate $(f_1, \ldots, f_n) \in \scrL(\scrA)$ is a nonempty Zariski open subset of $\scrL(\scrA)$.
\end{thm}

In both \cref{ex-aff-non-deg-1,ex-aff-non-deg-2}, $\dUUA{\kk^3}$ is the set $\{\emptyset, \{1\}, \{2\}, \{1,2\}\}$ of all subsets of $\{1,2\}$, and $\iUUA{\kk^3}$ is the set of remaining subsets of $\{1,2,3\}$. It is straightforward to check that in cases \ref{aff-non-deg-1-0} and \ref{aff-non-deg-1-1} of \cref{ex-aff-non-deg-1}, $f_1, f_2, f_3$ are $(\kk^3,\scrA)$-non-degenerate, but in case \ref{aff-non-deg-1-2}, condition \ref{iUA-property} of $(\kk^3,\scrA)$-non-degeneracy fails for $I = \{1, 2, 3\}$ and $\nu = (-1,1,2)$, which is centered at infinity. In \cref{ex-aff-non-deg-2}, condition \ref{iUA-property} of $(\kk^3,\scrA)$-non-degeneracy fails for $I = \{1, 2, 3\}$ and $\nu = (0,0,1)$, which is centered at $\Kstariprime$, where $I' := \{1,2\}$. In the scenario of \cref{ex-aff-non-deg-3}, $\iUUA{\kk^4}$ is the collection of all subsets of $\{1,2\}$, and $\dUUA{\kk^4}$ is the set of remaining subsets of $\{1,2,3,4\}$, and condition \ref{dUA-property} of $(\kk^4,\scrA)$-non-degeneracy fails for $I = \{1, 2, 3,4\}$ and $\nu = (0,0,1,2)$, which is centered at $\Kstariprime$, where $I' := \{1,2\}$.

\subsection{A more efficient formulation of $(U,\scrA)$-non-degeneracy}
In this section we describe a criterion equivalent to $(U,\scrA)$-non-degeneracy, but which involves checking fewer conditions. Recall that for $I \subseteq [n]$, we write $\tiA := \{j \in [n]:  \scrA_j \cap \ri \neq \emptyset\}$. Define
\begin{align}
\begin{split}
\dUAstar
	&:= \{I\in \dUA: I \neq \emptyset,\ |\tiA| = |I|\} \\
	&= \{I \subseteq [n]:
		I \neq \emptyset,\ |\tiA| = |I|,\ \text{$\scrA$ is $\ri$-dependent} \}
\end{split} \label{dUAstar} \\
\begin{split}
\iUAstar
	&:= \{I \in \iUA:  I \neq\emptyset,\ |\tiA| = |I|\} \\
	&= \{I \subseteq [n]:
		I \neq \emptyset,\ |\tiA| = |I|,\ \text{$\scrA$ is $\ri$-independent} \}
\end{split} \label{iUAstar}
\end{align}

\begin{prop}[{\cite{toricstein}}] \label{non-deg-2}
For polynomials $f_1, \ldots, f_n$ in $(x_1, \ldots, x_n)$, the following are equivalent:
\begin{enumerate}
\item $f_i$'s are $(U, \scrA)$-non-degenerate.
\item
\begin{enumerate}[label=(\roman{enumii})]
\item property \ref{U-property} of $(U, \scrA)$-non-degeneracy holds,
\item property \ref{iUA-property} of $(U, \scrA)$-non-degeneracy holds with $\iUA$ replaced by $\iUAstar$,
\item property \ref{dUA-property} of $(U, \scrA)$-non-degeneracy holds with $\dUA$ replaced by $\dUAstar$.
\end{enumerate}
\end{enumerate}
\end{prop}

We prove \cref{non-deg-2} in \cref{reducsection}.

\begin{example}[Warning!]
In property \ref{iUA-property} of $(U, \scrA)$-non-degeneracy $\dUA$ can {\em not} be replaced by $\dUAstar$, and in property \ref{dUA-property} of $(U, \scrA)$-non-degeneracy $\iUA$ can {\em not} be replaced by $\iUAstar$. Indeed, at first consider the system
\begin{align*}
f_1 &= 1 + x_1\\
f_2 &= 1 + x_1 + x_2\\
f_3 &= 1 + x_1 + 2x_2 + ax_3x_4\\
f_4 &= x_4(1 +  x_1 + bx_2 + cx_3 + dx_4)
\end{align*}
for generic $a, b, c, d \in \ktorus$. Let $U = \kk^4$ and $\scrA_j =  \supp(f_j)$, $j = 1, \ldots, 4$. Then $\dUA$ is the collection of all subsets of $\{1,2,3\}$, $\iUA$ is the collection of all subsets of $\{1, 2, 3, 4\}$ containing $4$, $\dUAstar = \{\{1,2,3\}\}$, and $\iUAstar = \{\{1,2,3,4\}\}$. Condition \ref{iUA-property} of $(U, \scrA)$-non-degeneracy fails with $I = \{1, 2,3,4\}$ and $\nu = (0,1,1,1)$, so that the center of $\nu$ is $\Kstarii{I'}$, where $I' := \{1\}\in \dUA$. However, it is straightforward to check that condition \ref{iUA-property} would not have been violated had $\dUA$ been replaced by $\dUAstar$. Now consider the system
\begin{align*}
g_1 &= 1 + x_1 + x_2 + x_3 \\
g_2 &= 1 + x_1 + 2x_2 + 3x_3 \\
g_3 &= x_2(ax_1 + bx_2 + cx_3)+ x_3(a'x_1+b'x_2+c'x_3) \\
g_4 &= x_4(1+x_1)
\end{align*}
where $a,b,c,a',b',c'$ are generic elements in $\ktorus$. Let $\scrB := (\supp(g_1), \ldots, \supp(g_4))$, and $U := \kk^4$. Then $\dUB$ is the collection of all subsets of $\{1,2,3,4\}$ containing $4$, $\iUB$ is the collection of all subsets of $\{1,2,3,4\}$ not containing $4$, $\dUBstar = \{\{1,2,3,4\}\}$, and $\iUBstar = \{\{1,2,3\}\}$. It is straightforward to check that $g_1, g_2, g_3, g_4$ violate condition \ref{dUA-property} of $(U, \scrB)$-non-degeneracy with $I = \{1, 2,3,4\}$ and $\nu = (0,1,1,1)$ (so that the center of $\nu$ is $\Kstarii{I'}$, where $I' := \{1\}\in \iUB$). It is also straightforward to check that condition \ref{dUA-property} would not have been violated had $\iUB$ been replaced by $\iUBstar$.
\end{example}


The following combinatorial description of strict monotonicity of $\multpiso{\cdot}{\cdot}{U}$ follows from \cref{non-deg-1,non-deg-2} exactly as in the proof of \cref{strictly-monotone-0}.

\begin{cor}
Let $\scrB_j \subseteq \conv(\scrA_j)$, $j = 1, \ldots, n$. Then $\multAiso{U} \geq \multBiso{U}$. If $\multAiso{U} > 0$, then the following are equivalent:
\begin{enumerate}
\item $\multAiso{U} = \multBiso{U}$,
\item
\begin{enumerate}
\item For each nonempty $I \in \iUA$ and for each weighted order $\nu$ on $\kk[x_i: i \in I]$ such that $\nu$ is centered at infinity or at $\Kstariprime$ for some $I' \in \excU \cup \excA \cup \dUA$, the collection $\{\In_\nu(\nd(\Bi_j)): \In_\nu(\nd(\Ai_j)) \cap \supp(\scrB_j) \neq \emptyset\}$ of polytopes is dependent, and
\item for each nonempty $I \in \dUA$ and for each weighted order $\nu$ on $\kk[x_i: i \in I]$ such that $\nu$ is centered at  $\Kstariprime$ for some $I' \in \iUA$, the collection $\{\In_\nu(\nd(\Bi_j)): \In_\nu(\nd(\Ai_j)) \cap \supp(\scrB_j) \neq \emptyset\}$ of polytopes is dependent.
\end{enumerate}
\item
\begin{enumerate}
\item For each nonempty $I \in \iUAstar$ and for each weighted order $\nu$ on $\kk[x_i: i \in I]$ such that $\nu$ is centered at infinity or at $\Kstariprime$ for some $I' \in \excU \cup \excA \cup \dUA$, the collection $\{\In_\nu(\nd(\Bi_j)): \In_\nu(\nd(\Ai_j)) \cap \supp(\scrB_j) \neq \emptyset\}$ of polytopes is dependent, and
\item for each nonempty $I \in \dUAstar$ and for each weighted order $\nu$ on $\kk[x_i: i \in I]$ such that $\nu$ is centered at  $\Kstariprime$ for some $I' \in \iUA$, the collection $\{\In_\nu(\nd(\Bi_j)): \In_\nu(\nd(\Ai_j)) \cap \supp(\scrB_j) \neq \emptyset\}$ of polytopes is dependent. \qed
\end{enumerate}
\end{enumerate}
\end{cor}

\section{Proof of the non-degeneracy conditions} \label{non-deg-proof-section}

\subsection{Reduction to the more efficient non-degeneracy criterion} \label{reducsection}
In this section we prove \cref{hereditary-prop,non-deg-2}.

\begin{proof}[Proof of \cref{hereditary-prop}]
Let $I' \supseteq I$ be as in the definition of hereditary dependence. By restricting all $f_j$'s to $\Kiprime$, we may assume that $I' = [n]$. Let $Z^I$ be the set of isolated points $V(f_1, \ldots, f_n)$ which are on $\Kstari$. Assume to the contrary of the claim that $Z^I \neq \emptyset$. Let $(g_1, \ldots, g_n) \in \scrL(\scrA)$ be a system such that
\begin{enumerate}
\item $\np(g_j) = \conv(\scrA_j)$, $j = 1, \ldots, n$, and
\item $g_1|_{\Kitilde}, \ldots, g_n|_{\Kitilde}$ are {\em properly $\Atildei$-non-degenerate} (see \cref{properly-non-degenerate-section}) for each $\tilde I \subseteq [n]$, where $\Atildei := (\scrA_1 \cap \rii{\tilde I}, \ldots, \scrA_n \cap \rii{\tilde n})$ (\cref{properly-non-degenerate-claim} implies that the set of such systems is a nonempty Zariski open subset of $\scrL(\scrA)$).
\end{enumerate}
Let $t$ be a new indeterminate, and $h_j := tg_j + (1-t)f_j \in \kk[x_1, \ldots, x_n, t]$, $j = 1, \ldots, n$. Pick $z \in Z^I$. Then $(z,0)$ is an isolated point of $V(h_1, \ldots, h_n, t)$ on $\kk^{n+1}$, and therefore \cref{thm:pure-dimension} implies that there is a Zariski open neighborhood $U$ of $(z,0)$ in $\kk^{n+1}$ such that $C := V(h_1, \ldots, h_n) \cap U$ is a curve. Let $C'$ be an irreducible component of $C$ containing $(z,0)$, and $\bar C'$ be the closure of $C'$ in $\pp^n \times \pp^1$. Then $C' \not\subseteq \kk^n \times \{\epsilon\}$ for any $\epsilon \in \kk$, and \cref{linear-intersection} implies that $\bar C'$ intersects $\pp^n \times \{1\}$. Pick a branch $B$ of $\bar C'$ centered at a point $(z',1) \in \bar C' \cap (\pp^n \times \{1\})$. Define $I_B \subseteq [n+1]$ and the weighted order $\nu_B$ on the coordinate ring of $\Kib$ as in \cref{IB-defn}. Let $\tilde I := I_B \cap [n]$ and $\tilde \nu$ be the restriction of $\nu_B$ to $\kk[x_i : i \in \tilde I]$. Since $\np(f_j) \subseteq \np(g_j)$ for each $j$, it follows that $\tilde \nu(f_j|_{\Kitilde}) \geq \tilde \nu(g_j|_{\Kitilde})$. Since $\nu_B(t-1) > 0$, it follows that $\nu_B(t) = 0$, and therefore $\In_{\nu_B}(h_j|_{\Kib}) = t\In_{\tilde \nu}(g_j|_{\Kitilde})$ for each $j$. \Cref{branch-lemma-IB} then implies that $\In_{\tilde \nu}(g_j|_{\Kitilde})$ have a common zero in $\Kstaritilde$. On the other hand, since $\bar C'$ contains the point $(z,0) \in \Kstari$, it follows that $I \subseteq \tilde I$, and therefore it follows from the definition of hereditary dependence (applied with $I' = [n]$) that either $|\titildeA| > |\tilde I|$, or $|\titildeA| =|\tilde I|$ and $\scrA$ is $\ri$-dependent. In any event, there is a nonempty subset $J$ of $\titildeA$ such that $\dim(\sum_{j \in J}\np(g_j|_{\Kitilde})) < |J|$. Since $g_1|_{\Kitilde}, \ldots, g_n|_{\Kitilde}$ are properly $\Atildei$-non-degenerate, it then follows by definition of proper non-degeneracy that  there is no common zero of $\In_{\tilde \nu}(g_j|_{\Kitilde})$ on $\Kstaritilde$. This contradiction completes the proof.
\end{proof}

\begin{proof}[Proof of \cref{non-deg-2}]
Since $\dUAstar \subseteq \dUA$ and $\iUAstar \subseteq \iUA$, it suffices to show the following:
\begin{enumerate}
\item property \ref{iUA-property} of $(U, \scrA)$-non-degeneracy holds if it holds with $\iUA$ replaced by $\iUAstar$,
\item property \ref{dUA-property} of $(U, \scrA)$-non-degeneracy holds if it holds with $\dUA$ replaced by $\dUAstar$.
\end{enumerate}
It follows from the definition of hereditary dependence that for every $I \in \iUA \setminus \iUAstar$, there exists $I' \in \iUAstar$ such that $I \subsetneq I'$ and $|\tiprimeA| > |I'|$ for each $I'$ such that $I \subseteq I' \subsetneq \tilde I$. The same statement also holds if we replace $\iUA$ and $\iUAstar$ respectively by $\dUA$ and $\dUAstar$. Since restricting all $f_j$ to $\Kitilde$ yields a system with the same number of nonzero polynomials as the number of variables, it suffices to prove the following claim: ``if there is $I \subseteq [n]$ such that $|\tiA| > |I|$ and $\In_{\Ai_1,\nu}(f_1|_{\Ki}),  \ldots, \In_{\Ai_n, \nu}(f_n|_{\Ki})$ have a common zero on $\nktorus$ for some weighted order $\nu$ on $\kk[x_i: i \in I]$, then there is $\tilde I \supsetneq I$ and a weighted order $\tilde \nu$ on $\kk[x_{\tilde i}: \tilde i \in \tilde I]$ such that $\tilde \nu$ is compatible with $\nu$, and $\In_{\Atildei_1, \tilde \nu}(f_1|_{\Kitilde}), \ldots, \In_{\Atildei_n, \tilde \nu}(f_n|_{\Kitilde})$ have a common zero in $\tilde u \in \nktorus$.'' This follows from \cref{reduction-lemma}.
\end{proof}

\subsection{Understanding condition \ref{dUA-property} of $(U, \scrA)$-non-degeneracy}
In this section we show that if condition \ref{dUA-property} of $(U, \scrA)$-non-degeneracy is violated for $(f_1, \ldots, f_n) \in \scrL(\scrA)$ with $I$, $\nu$ and $I'$, then the common zero on $\nktorus$ of $\In_{\Ai_1, \nu}(f_1|_{\Ki}), \ldots, \In_{\Ai_n, \nu}(f_1|_{\Ki})$ corresponds to a non-isolated point of $V(f_1, \ldots, f_n)$. Recall the definition of $\pi_I: \kk^n \to \Ki$ from \eqref{pi_I}.

\begin{lemma} \label{trivial-lemma}
Let $(f_1, \ldots, f_n) \in \scrL(\scrA)$, and $\nu$ be a weighted order on $\kk[x_1, \ldots, x_n]$ centered at $\Kstari$, $I \subseteq [n]$.
Assume
\begin{enumerate}
\item \label{trivially-dependent} $\conv(\scrA_1), \ldots, \conv(\scrA_n)$ are dependent polytopes in $\rr^n$, and
\item \label{trivially-common} $\In_{\scrA_1, \nu}(f_1),  \ldots,\In_{\scrA_n, \nu}(f_n)$ have a common zero $a \in \nktorus$.
\end{enumerate}
Then $\pi_{I}(a) \in \Kstari$ is a non-isolated point of the set $V(f_1, \ldots, f_n)$ of common zeroes of $f_1, \ldots, f_n$ on $\kk^n$.
\end{lemma}

Note that both conditions \eqref{trivially-dependent} and \eqref{trivially-common} are necessary for the conclusion of \cref{trivial-lemma} to hold. For example, $f_1 = 1 + x_1, f_2 = 1 + x_1 + x_2$ satisfy condition \eqref{trivially-common}, but not \eqref{trivially-dependent}, with $I = \{1\}$, $\nu = (0,1)$, $\scrA_j = \supp(f_j)$, $j = 1,2$, and $a = (-1,c)$ for some arbitrary $c \in \kk^*$. In this case $\pi_I(a) = (-1,0)$ is an {\em isolated} point of $V(f_1, f_2)$. On the other hand, $f_1 = 1 + x_1, f_2 = x_2$, and $\scrA_j = \supp(f_j)$, $j = 1,2$, satisfy, with the same $I$, $\nu$ and $a$, condition \eqref{trivially-dependent}, but not \eqref{trivially-common}, and $\pi_I(a) = (-1,0)$ is again an isolated point of $V(f_1, f_2)$.

\begin{proof}[Proof of \cref{trivial-lemma}]
It is immediate to check that $\pi_{I}(a)$ is in $V(f_1, \ldots, f_n)$; we only have to show that it is non-isolated in there. Since $\scrP_j := \conv(\scrA_j)$, $j = 1, \ldots, n$, are dependent, it follows that there is $J \subseteq [n]$ such that  $p:= \dim(\sum_{j \in J} \scrP_j) <|J|$. Let $\Pi$ be the (unique) $p$-dimensional linear subspace of $\rr^n$ such that $\sum_{j \in J} \scrP_j$ is contained in a translate of $\Pi$. Let $\nu_j := \nu(x_j)$, $j = 1, \ldots, n$. We identify $\nu$ with the element of $\rnstar$ with coordinates $(\nu_1, \ldots, \nu_n)$ with respect to the basis dual to the standard basis of $\rr^n$. Let $\Pi_0 := \{\alpha \in \Pi: \langle \nu, \alpha \rangle = 0\}$ and $r := \dim(\Pi_0)$. Choose a basis $\alpha_1, \ldots, \alpha_r$ of $\Pi_0 \cap \zz^n$. Let $c_i := a^{\alpha_i}$, $i = 1, \ldots, n$. Let $Y$ be the subvariety of $\nktorus$ determined by $x^{\alpha_i} - c_i$, $i = 1, \ldots, r$, and $\bar Y$ be the closure of $Y$ in $\kk^n$. Then $Y$, and therefore $\bar Y$, is an irreducible variety of codimension $r$ in $\kk^n$.

\begin{proclaim} \label{pi_I(a)-claim}
$\pi_I(a) \in \bar Y$. Moreover, if $g$ is any Laurent polynomial in $(x_1, \ldots, x_n)$ such that $\supp(g) \subseteq \Pi_0$, then $g$ restricts to a constant function on $Y$ with value $g(a)$.
\end{proclaim}

\begin{proof}
The second assertion is obvious, so we prove the first assertion. If $\nu$ is the trivial weighted order, then $\pi_I$ is the identity and therefore $\pi_I(a) = a \in Y$. Otherwise let $C$ be the curve on $\kk^n$ parametrized by $t \mapsto c(t) := (a_1t^{\nu_1}, \ldots, a_nt^{\nu_n})$. Then $C \cap \nktorus \subseteq Y$, so that $C \subseteq \bar Y$. Since $\pi_I(a) = c(0) \in C$, the claim is proved.
\end{proof}

Note that $r$ equals either $p$ or $p -1$. If $r = p$, then $\Pi_0 = \Pi$ and for each $j \in J$, $f_j$ is $\nu$-homogeneous and is of the form $f_j = x^{\beta_j}g_j$ for some $\beta_j \in \zz^n$ and Laurent polynomial $g_j$ such that $\supp(g_j) \subseteq \Pi_0$. \Cref{pi_I(a)-claim} implies that $f_j|_Y \equiv 0$ for each $j \in J$, so that $\pi_I(a) \in \bar Y \subseteq V(f_j :j \in J)$. It follows that at least one of the irreducible components of $V(f_j :j \in J)$ containing $\pi_I(a)$ has codimension smaller than $|J|$ in $\kk^n$. The lemma then follows due to \cref{thm:pure-dimension}. It remains to consider the case that $r = p-1$. \Cref{basis-lemma} implies that $\alpha_1, \ldots, \alpha_r$ can be extended to a basis $\alpha_1, \ldots, \alpha_n$ of $\zz^n$ such that $\alpha_1, \ldots, \alpha_{r+1}$ is a basis of $\Pi$, and $\langle \nu, \alpha_j \rangle \geq 0$ for each $j = 1, \ldots, n$; this in particular implies that $\langle \nu, \alpha_{r+1} \rangle > 0$. Let $y_i := x^{\alpha_i}$, $i = 1, \ldots, n$. Then the $y_i$ form a system of coordinates on $\nktorus$ and the projection onto $(y_{r+1}, \ldots, y_n)$ restricts to an isomorphism $Y \cong \nktoruss{n-r}$. Pick $\beta_j := (\beta_{j,1}, \ldots, \beta_{j,n}) \in \zz^n$ such that $x_j = \prod_{i=1}^n y_i^{\beta_{j,i}}$, $j = 1, \ldots, n$. Then $\bar Y$ is the closure of the image of the map $\psi: \nktoruss{n-r} \to \kk^n$ given by
\begin{align*}
\psi (y_{r+1}, \ldots, y_n) :=  ( c'_1 \prod_{i=r+1}^n y_i^{\beta_{1,i}}, \ldots, 	 c'_n \prod_{i=r+1}^n y_i^{\beta_{n,i}}),
\end{align*}
where $c'_j := \prod_{i=1}^r c_i^{\beta_{j,i}}$, $j = 1, \ldots, n$. Let $\bar Y'$ be the closure of the image of the map $\psi' :\nktoruss{n-r} \to \kk^n$ given by
\begin{align*}
\psi'(y_{r+1}, \ldots, y_n) :=  (\prod_{i=r+1}^n y_i^{\beta_{1,i}}, \ldots, \prod_{i=r+1}^n y_i^{\beta_{n,i}})					
\end{align*}
Then $\bar Y'$ is isomorphic to $\bar Y$ via the map $\rho: (x_1, \ldots, x_n) \mapsto (c'_1x_1, \ldots, c'_nx_n)$.
Let $\beta'_j := (\beta_{j,r+1}, \ldots, \beta_{j,n})$, $j = 1, \ldots, n$, and $\scrB' := \{\beta'_0, \ldots, \beta'_n\}$, where $\beta'_0$ is the origin in $\zz^{n-r}$. Let $\xbprime$ be the corresponding toric variety. In the notation of \cref{xa-thm} $\bar Y'$ is isomorphic to the affine open subset $\xbprime \cap U_{\beta'_0}$ of $\xbprime$. Since $\dim(\bar Y') = n-r$, it follows that the convex hull $\scrP'$ of $\scrB'$ in $\rr^{n-r}$ has dimension $n-r$. Let $\nu'_j := \langle \nu, \alpha_j \rangle$, $j = r+1, \ldots, n$, and $\nu' \in \rnnstar{n-r}$ be the element with coordinates $(\nu'_{r+1}, \ldots, \nu'_n)$. Then for each $j = 1, \ldots, n$, $\langle \nu', \beta'_j \rangle = \sum_{i = r+1}^n \beta_{j,i} \langle \nu, \alpha_i \rangle = \langle \nu, \sum_{i=1}^n \beta_{j,i} \alpha_i \rangle = \nu_j \geq 0$. It follows that $\min_{\scrP'}(\nu') = 0$ (since $\langle \nu', \beta'_0 \rangle = 0$). Since $\nu'_{r+1} > 0$, there is a facet $\scrQ'$ of $\scrP'$ containing $\In_{\nu'}(\scrP')$ such that the first coordinate of the inner normal with respect to the dual basis is positive (this follows e.g.\ from \cref{prop:edge-hull,normal-faces,stritctly-convex-fan}). Let $\scrC' := \scrB' \cap \scrQ'$, and $\corbit{\scrC'}$ be the corresponding torus invariant codimension one subvariety of $\xbprime$ (see \cref{xa-thm}). Let $Z' := \bar Y' \cap \corbit{\scrC'}$, and $Z := \rho(Z') \subseteq \bar Y$.

\begin{proclaim} \label{pi_I(a)-claim-2}
$\pi_I(a) \in Z \subseteq V(f_j: j \in J)$.
\end{proclaim}

\begin{proof}
Let $\eta \in \rnnstar{n-r}$ be the primitive inner normal to $\scrQ'$, and $(\eta_{r+1}, \ldots, \eta_n)$ be the coordinates of $\eta$ in the dual basis of $\rnnstar{n-r}$. For each $b = (b_{r+1}, \ldots, b_n) \in \nktoruss{n-r}$, consider the rational curve on $\bar Y'$ parametrized by
$$t \mapsto c'_b(t) := \psi'(b_{r+1}t^{\eta_{r+1}}, \ldots, b_nt^{\eta_n}) = (b^{\beta'_1}t^{\eta'_1}, \ldots, b^{\beta'_n}t^{\eta'_n})$$
where $\eta'_j := \langle \eta, \beta'_j \rangle$ for each $j = 1, \ldots, n$. Since $0 = \beta_0 \in \In_{\nu'}(\scrP') \subset \scrQ'$, it follows that $\min_{\scrP'}(\eta) = 0$, so that each $\eta'_j$ is nonnegative, and it is zero if and only if $\beta'_j \in \scrQ'$. In particular, $c'_b(0)$ is well defined, and is an element in $\kk^n$. Moreover, \cref{xa-thm} implies that the set $\{c'_b(0): b \in \nktoruss{n-r}\}$ is precisely $\bar Y' \cap \orbit{\scrC'}$; in particular it is a dense Zariski open subset of $Z'$. Therefore, in order to prove $Z \subseteq V(f_j: j \in J)$, it suffices to show that $f_j (\rho( c'_b(0))) = 0$ for each $j \in J$. Fix $j \in J$. Write $f_j$ as $f_j = f_{j,0} + f_{j,1} + \cdots $, where $f_{j,k}$ are $\nu$-homogeneous polynomials with $\nu(f_{j,0}) < \nu(f_{j,1}) < \cdots$. We will show that
\begin{align}
f_{j,k} (\rho(c'_b(0)))= 0 \label{f-j-k-0}
\end{align}
for each $k$. Let $m_j := \min_\nu(\scrA_j)$. Note that $\nu(f_{j,k}) \geq m_j$ for each $k$. At first consider the case that $\nu(f_{j,k}) = m_j$. This implies that $k = 0$, $\In_{\nu}(\scrA_j) \cap \supp(f_j)$ is nonempty, and $f_{j,k} = \In_{\scrA_j, \nu}(f_j)$. Since $\In_{\scrA_j, \nu}(f_j)(a) = 0$, in this case \cref{pi_I(a)-claim} implies that $f_{j,k}$ is identically zero on $\rho(c'_b(t))$, $t \in \kk$; in particular, \eqref{f-j-k-0} holds. Now fix $k$ such that $\nu(f_{j,k}) > m_j$. It suffices to show that
\begin{align}
\ord_t(f_{j,k} \circ \rho \circ  c'_b(t)) > 0 \label{ord-t-positive}
\end{align}
Pick $\alpha_0 \in \In_\nu(\scrA_j)$ and $\alpha \in \supp(f_{j,k})$. Then $\alpha - \alpha_0 \in \Pi$, and $\langle \nu, \alpha - \alpha_0 \rangle = \nu(f_{j,k}) - m_j > 0$. It follows that $\alpha - \alpha_0 = \sum_{j=1}^{r+1} m_j \alpha_j$ with integers $m_j$ such that $m_{r+1} > 0$. Since $x^{\alpha_j}|_Y \equiv c_j$ for $j = 1, \ldots, r$, and $x^{\alpha_{r+1}}|_Y \equiv  y_{r+1}|_Y$, it follows that
\begin{align*}
\ord_t(x^\alpha \circ \rho \circ c'_b(t))
	&= \ord_t(x^{\alpha_0} \circ \rho \circ c'_b(t)) + m_{r+1}\eta_{r+1}
\end{align*}
%
Since $\eta_{r+1} > 0$ (due to our choice of $\scrQ'$) and $\alpha_0 \in \znzero$, \eqref{ord-t-positive} follows. It remains to prove that $\pi_I(a) \in Z$. Consider the rational curve $C'$ on $\bar Y'$ parametrized by
\begin{align*}
t &\mapsto \psi'(c_{r+1}t^{\nu'_{r+1}}, \ldots, c_n t^{\nu'_n})
	= ((\prod_{i=r+1}^n c_i^{\beta_{1,i}}) t^{\nu_1}, \ldots, (\prod_{i=r+1}^n c_i^{\beta_{n,i}}) t^{\nu_n})
\end{align*}
Since $\In_{\nu'}(\scrP')$ is a face of $\scrQ'$, \cref{curve-to-orbit} implies that the center of the branch of $C'$ at $t = 0$ is contained in $Z'$. Note that $\rho(C')$ is precisely the curve $C$ from \cref{pi_I(a)-claim}. It follows from \cref{pi_I(a)-claim} that $\pi_I(a)$ is the center of the branch of $C$ at $t = 0$, and therefore it is on $Z$, as required.
\end{proof}
Since $\dim(Z) = n- r - 1 > n-|J|$, the lemma follows from \cref{pi_I(a)-claim-2,thm:pure-dimension}.
\end{proof}

\begin{cor}  \label{trivial-cor}
Let $I \in \dUA \setminus \{\emptyset\}$ and $\nu$ be a weighted order on $\kk[x_i: i \in I]$ such that $\nu$ is centered at $\Kstariprime$ for some $I' \subseteq I$, and $\In_{\Ai_1, \nu}(f_1|_{\Ki}), \ldots, \In_{\Ai_n,\nu}(f_n|_{\Ki})$ have a common zero $a \in \nktorus$. Then $\pi_{I'}(a) \in \Kstariprime$ is a non-isolated point of the zero-set $V(f_1, \ldots, f_n)$ of $f_1, \ldots, f_n$ on $\kk^n$.
\end{cor}

\begin{proof}
It is straightforward to check that $\pi_{I'}(a) \in V(f_1, \ldots, f_n)$. There is $J \supseteq I$ such that $\scrA$ is $\rj$-dependent, $|\tjA| = |J|$, and $|\titildeA| > |\tilde I|$ for each $\tilde I$ such that $I \subseteq \tilde I \subsetneq J$. Replacing the $f_j$ by $f_j|_{\Kii{J}}$ and applying \cref{reduction-lemma} reduces the corollary to the case that $I = [n]$. Then it follows from \cref{trivial-lemma}.
\end{proof}

\subsection{Proof of \cref{non-deg-1}} \label{non-deg-1-proof-section}
We divide the proof of \cref{non-deg-1} in three parts:
\subsubsection{Proof of the implication \eqref{non-deg-1-1} $\im$ \eqref{non-deg-1-0} from \cref{non-deg-1}}
We proceed as in the proof of \cref{bkk-sufficiency}. Pick $(f_1, \ldots, f_n), (g_1, \ldots, g_n) \in \scrL(\scrA)$ such that $\multgiso{U}> \multfiso{U}$. We will show that $f_1, \ldots, f_n$ are $(U,\scrA)$-degenerate. Define $h_j := (1-t)f_j + tg_j$, where $t$ is a new indeterminate. The set of isolated zeroes of $g_1, \ldots, g_n$ on $U$ is nonempty. \Cref{mult-deformation-global} implies that we can find an irreducible curve $C$ contained in the set of zeroes of $h_1, \ldots, h_n$ on $U \times \kk$ such that
\begin{prooflist}
\item \label{intersection-1} $C$ intersects $Z \times \{1\}$,
\item \label{generically-isolated} for generic $\epsilon \in \kk$, the set $C_\epsilon := C \cap (U \times \{\epsilon\})$ is nonempty, and each point of $C_\epsilon$ is an isolated zero of $h_1|_{t = \epsilon}, \ldots, h_n|_{t= \epsilon}$ on $U$; and
\item \label{property-at-zero}
\begin{prooflist}
\item \label{non-isolated-case} either there is $(z,0) \in C$ such that $z$ is a non-isolated zero of $f_1, \ldots, f_n$ on $U$,
\item \label{infinite-case} or $C$ has a ``point at infinity with respect to $U$'' at $t = 0$, i.e.\ if $\bar C$ is the closure of $C$ in $\pp^n \times \kk$, then there is $(z,0) \in \bar C$ such that $z \not\in U$.
\end{prooflist}
\end{prooflist}

\begin{claim} \label{tilde-I-claim}
Let $B$ be a branch of $\bar C \subset \pp^n \times \kk$. Let $\tilde I := I_B \cap [n]$ (where $I_B \subseteq [n+1]$ is defined as in \cref{IB-defn}) and $\tilde \nu$ be the restriction of $\nu_B$ to $\kk[x_i: i \in \tilde I]$.
\begin{enumerate}
\item \label{tilde-I-generic} $C_\epsilon \subseteq \Kstaritilde \times \{\epsilon\}$ for generic $\epsilon \in \kk$.
\item \label{tilde-I-UA} $\tilde I \in \iUA$.
\item \label{tilde-I-t0} Assume the center of $B$ is $(z,0)$ where $z \in \pp^n$. Then
\begin{enumerate}
\item \label{tilde-I-t0-zero} $\In_{\Atildei_j, \tilde \nu}(f_j|_{\Kitilde})$, $j = 1, \ldots, n$, have a common zero on $\nktorus$.
\item \label{tilde-I-t0-center-infinity} If $z \in \pp^n\setminus \kk^n$, then $\tilde I \neq \emptyset$ and $\tilde \nu$ is centered at infinity. In particular, $f_1, \ldots, f_n$ violate condition \ref{iUA-property} of $(U,\scrA)$-non-degeneracy with $I = \tilde I$ and $\nu = \tilde \nu$.
\item \label{tilde-I-t0-center-finity} Otherwise let $J$ be the (unique) subset of $[n]$ such that $z \in \Kstarii{J}$. Then either $\tilde I = J = \emptyset$, or $\tilde I \neq \emptyset$ and $\tilde \nu$ is centered at $\Kstarii{J}$. In particular, if $J \in \excU \cup \excA \cup \dUA$, then $f_1, \ldots, f_n$ violate property \ref{iUA-property} of $(U,\scrA)$-non-degeneracy with $I = \tilde I$, $I' = J$ and $\nu = \tilde \nu$.

\end{enumerate}
\end{enumerate}
\end{claim}

\begin{proof}
Assertion \eqref{tilde-I-generic} follows from the definition of $I_B$. Since $(h_1|_{t = \epsilon}, \ldots, h_n|_{t= \epsilon}) \in \scrL(\scrA)$ for all $\epsilon$, assertion \eqref{tilde-I-generic}, property \ref{generically-isolated} and \cref{non-isolated-lemma} imply that $\tilde I \not\in \excU \cup \excA$, and then \cref{hereditary-prop} implies that $\tilde I \in \iUA$. This proves assertion \eqref{tilde-I-UA}. 
Now assume we are in the situation of assertion \eqref{tilde-I-t0}. Fix $j$, $1 \leq j \leq n$. Since $\nu_B(t|_{\Kib}) > 0$, it follows that
\begin{align*}
\In_{\Atildei_j, \tilde \nu}(f_j|_{\Kitilde})
	&=
	\begin{cases}
	\In_{\tilde \nu}(f_j|_{\Kitilde}) = \In_{\nu_B}(h_j|_{\Kib})& \text{if}\ \supp(f_j) \cap \In_{\tilde \nu}(\Atildei_j) \neq \emptyset,\\
	0 & \text{otherwise.}
	\end{cases}
\end{align*}
Assertion \eqref{tilde-I-t0-zero} now follows from \cref{branch-lemma-IB}. If $z \in \pp^n \setminus \kk^n$, then there is at least one $j$ such that $1/x_j$ is a regular function near $z$ which vanishes at $z$. This $j$ has to be in $\tilde I$ and $\tilde \nu(x_j)$ has to be negative, which proves the first statement of part \eqref{tilde-I-t0-center-infinity}. The second statement then follows from assertions \eqref{tilde-I-UA} and \eqref{tilde-I-t0-zero}. The first statement of part \eqref{tilde-I-t0-center-finity} is obvious. If $J \in \excU \cup \excA \cup \dUA$, then assertion \eqref{tilde-I-UA} and the first statement of part \eqref{tilde-I-t0-center-finity} imply that $\tilde I \neq \emptyset$, and then the second statement follows from assertions \eqref{tilde-I-UA} and \eqref{tilde-I-t0-zero}.
\end{proof}

Now we resume the proof of the implication \eqref{non-deg-1-1} $\im$ \eqref{non-deg-1-0} from \cref{non-deg-1}. At first assume \ref{non-isolated-case} holds. Pick a branch $B$ of $C$ centered at $(z,0)$ and let $\tilde I, J, \tilde \nu$ be as in part \eqref{tilde-I-t0-center-finity} of \cref{tilde-I-claim}. Part \eqref{tilde-I-t0-center-finity} of \cref{tilde-I-claim} implies that $f_1, \ldots, f_n$ are $(U, \scrA)$-degenerate if $J \in \excU \cup \excA \cup \dUA$. So assume $J \in \iUA$. Pick an irreducible curve $C'$ on $U$ such that $z \in C' \subseteq V(f_1, \ldots, f_n)$. Let $J'$ be the smallest subset of $[n]$ such that $C' \subseteq \Kii{J'}$. Since $J' \supseteq J$ and $J \not\in \excU \cup \excA$, it follows that $J' \not\in \excU \cup \excA$ as well. If $J' \in \iUA$, then \cref{branch-lemma-IB} implies that $f_1, \ldots, f_n$ violate condition \ref{iUA-property} of $(U,\scrA)$-non-degeneracy with $I = J'$ and some weighted order $\nu$ on $\kk[x_j: j \in J']$ centered at infinity (take a branch $B'$ of $C'$ centered at infinity, and set $\nu = \nu_{B'}$). On the other hand, if $J' \in \dUA$, then $f_1, \ldots, f_n$ violate condition \ref{dUA-property} of $(U,\scrA)$-non-degeneracy with $I = J'$ and $I' = J$ (take a branch $B'$ of $C'$ centered at $z$ and take $\nu := \nu_{B'}$). This completes the proof of \eqref{non-deg-1-1} $\im$ \eqref{non-deg-1-0} in the case that \ref{non-isolated-case} holds. Now assume we are in case \ref{infinite-case}. Pick $z \in \pp^n\setminus U$ such that $(z,0) \in \bar C$. If $z \in \pp^n \setminus \kk^n$, then part \eqref{tilde-I-t0-center-infinity} of \cref{tilde-I-claim} implies that $f_1, \ldots, f_n$ are $(U,\scrA)$-degenerate. So assume $z \in \kk^n \setminus U$. Define $J$ as in part \eqref{tilde-I-t0-center-finity} of \cref{tilde-I-claim}. If $J \in \excU \cup \excA \cup \dUA$, then part \eqref{tilde-I-t0-center-finity} of \cref{tilde-I-claim} implies that $f_1, \ldots, f_n$ are $(U, \scrA)$-degenerate. So assume $J \in \iUA$. But then $f_1, \ldots, f_n$ violate condition \ref{U-property} of $(U,\scrA)$-non-degeneracy. This completes the proof of the implication \eqref{non-deg-1-1} $\im$ \eqref{non-deg-1-0}. \\

\subsubsection{Proof of the implication \eqref{non-deg-1-0} $\im$ \eqref{non-deg-1-1} from \cref{non-deg-1}}
Assume $\multAiso{U} > 0$ and pick an $(U,\scrA)$-degenerate system $f_1, \ldots, f_n \in \scrL(\scrA)$. We will show that $\multAiso{U} > \multfiso{U}$. Recall the definition of $\iUAstar$ from \eqref{iUAstar}.

\begin{claim} \label{case-degenerate}
There is $I \in \iUAstar$ such that one of the following holds:
\begin{enumerate}
\item \label{case-nonempty} Either $I$ is nonempty, and there is a weighted order $\nu$ on $\kk[x_i: i \in I]$ such that $\In_{\Ai_j,\nu}(f_j|_{\Ki})$, $j = 1, \ldots, n$, have a common zero $a \in \nktorus$, and one of the following holds:
\begin{enumerate}
\item \label{case-infinity} $\nu$ is centered at infinity,
\item \label{case-U} or $\nu$ is centered at $\Kstariprime$ for some $I' \subseteq I$, and $\pi_{I'}(a) \not\in U$,
\item \label{case-non-isolated} or $\nu$ is centered at $\Kstariprime$ for some $I' \subseteq I$ and $\pi_{I'}(a)$ is a non-isolated zero of $f_1, \ldots, f_n$.
\end{enumerate}
\item \label{case-isolated-but-not} Or there is an isolated point $a$ of $V(f_1, \ldots, f_n) \cap \Ki \cap U$ which is {\em not} isolated in $V(f_1, \ldots, f_n) \subset \kk^n$.
\end{enumerate}
\end{claim}

\begin{proof}
If $f_1, \ldots, f_n$ violate property \ref{U-property} of $(U, \scrA)$-non-degeneracy, then there is a common zero $a'$ of $f_1|_{\Kiprime}, \ldots, f_n|_{\Kiprime}$ on $\Kstariprime \setminus U$ for some $I' \in \iUA$. \Cref{reduction-lemma} then implies that the claim holds with case \eqref{case-U} holds. If property \ref{iUA-property} of $(U, \scrA)$-non-degeneracy fails with $I \in \iUAstar$, then either the claim holds with case \eqref{case-infinity}, or there is a weighted order $\nu$ on $\kk[x_i: i \in I]$ such that $\In_{\Ai_j, \nu}(f_j|_{\Ki})$ have a common zero $a \in \nktorus$, and $\nu$ is centered at $\Kstariprime$ for some $I' \in \excU \cup \excA \cup \dUA$. It is straightforward to check that $\pi_{I'}(a)$ is in the set $V$ of common zeroes of $f_1, \ldots, f_n$ on $\kk^n$. If $I' \in \excU$, then we are in case \eqref{case-U}, since $\pi_{I'}(a) \in \Kstariprime$ and $\Kstariprime \cap U = \emptyset$. If $I' \in \excA \cup \dUA$, then $\pi_{I'}(a)$ has to be a non-isolated point of $V$ due to \cref{non-isolated-lemma,hereditary-prop}, which is case \eqref{case-non-isolated}. Due to \cref{non-deg-2} the only case left to consider is that of $f_1, \ldots, f_n$ violating property \ref{dUA-property} of $(U,\scrA)$-non-degeneracy. Then there is $J \in \dUA$ and a weighted order $\eta$ on $\kk[x_j: j \in J]$ centered at $\Kstarii{J'}$ for some $J' \in \iUA$ such that $\In_{\Aj_1, \eta}(f_1|_{\Kii{J}}), \ldots, \In_{\Aj_n, \eta}(f_n|_{\Kii{J}})$ have a common zero $b$ on $\nktorus$. \Cref{trivial-cor} implies that $\pi_{J'}(b)$ is a non-isolated point of $V(f_1, \ldots, f_n)$. Pick the smallest subset $I$ of $[n]$ containing $J'$ such that $|\tiA| = |I|$. Since $J' \not\in \excU \cup \excA$, it follows that $I \not\in \excU \cup \excA$, and since $\scrA$ is not hereditarily $\rii{J'}$-dependent, it follows that $\scrA$ is $\ri$-independent; in particular, $I \in \iUAstar$. If $\pi_{J'}(b)$ is an isolated point of $V(f_1, \ldots, f_n) \cap \Ki$ (which is e.g.\ the case if $I = J' = \emptyset$), then case \eqref{case-isolated-but-not} holds with $a = \pi_{J'}(b)$. Otherwise picking a branch at $\pi_{J'}(b)$ of a curve contained in $V(f_1, \ldots, f_n) \cap \Ki$ and applying \cref{branch-lemma-IB,reduction-lemma} shows that case \eqref{case-non-isolated} holds.
\end{proof}

\begin{claim} \label{case-degenerate-2}
Assume case \eqref{case-isolated-but-not} of \cref{case-degenerate} holds. Then there is $(g_1, \ldots, g_n) \in \scrL(\scrA)$ such that $g_j|_{\Ki} = f_j|_{\Ki}$ for each $j$, and $a$ is an isolated point of $V(g_1, \ldots, g_n) \subset \kk^n$.
\end{claim}

\begin{proof}
Fix $\tilde I \supsetneq I$. Since $|\tiA| = |I|$ and $I \not\in \excA$, it follows that $\titildeA \setminus \tiA$ contains at least $|\tilde I| - |I|$ elements. For each weighted order $\nu$ on $\kk[x_i: i \in \tilde I]$ centered at $\Kstari$, choosing generically coefficients of $f_j$, $j \in \titildeA \setminus \tiA$, it can be ensured that $\In_{\Atildei_j, \nu}(f_j)$, $j \in \titildeA \setminus \tiA$, have no common zero $\tilde a$ on $\nktorus$ such that $\pi_I(\tilde a) = a$. The claim now follows due to \cref{branch-lemma-IB}.
\end{proof}

At first consider case \eqref{case-isolated-but-not} of \cref{case-degenerate}. Pick $(g_1, \ldots, g_n)$ as in \cref{case-degenerate-2}. Apply \cref{mult-deformation-global} to $X = U$ and $h_j = (1-t)f_j + tg_j$, $j = 1, \ldots, n$. It is straightforward to see that in this case $\{a\} \times \kk$ is an irreducible component of the curve $C$ from \cref{mult-deformation-global}, so that
\begin{align*}
\multfiso{U} < \multpiso{h_1|_{t=\epsilon}}{h_n|_{t=\epsilon}}{U}
\end{align*}
for generic $\epsilon \in \kk$. It follows that $\multfiso{U} < \multAiso{U}$, as required. Now assume there are $I$, $\nu$ and $a$ as in case \eqref{case-nonempty} of \cref{case-degenerate}. We may assume \woutlog\ that $I = \tiA = \{1, \ldots, k\}$ for some $k$, $1 \leq k \leq n$.

\begin{claim} \label{case-degenerate-1}
There is $(g_1, \ldots, g_n) \in \scrL(\scrA)$ such that
\begin{enumerate}
\item \label{g-a-non-zero} $\In_{\Ai_j, \nu}(g_j|_{\Ki})(a) \neq 0$ for each $j = 1, \ldots, k$.
\item \label{g-b-zero} there is a common zero $b$ of $g_1, \ldots, g_k$ on $\Kstari \cap U$ such that
\begin{enumerate}
\item \label{g-b-isolated} $b$ is an isolated point of $V(g_1, \ldots, g_n) \subset \kk^n$, and
\item \label{f-b-non-zero} $f_j(b) \neq 0$ for each $j = 1, \ldots, k$.
\end{enumerate}
\end{enumerate}
\end{claim}

\begin{proof}
Let $\scrL'(\scrA)$ be the collection of all $(g_1, \ldots, g_n) \in \scrL(\scrA)$ such that
\begin{prooflist}
\item \label{g-bkk-non-deg} $g_1|_{\Kii{J}}, \ldots, g_n|_{\Kii{J}}$ are $\Aj$-non-degenerate (in the sense of \cref{defn:b-non-deg}) for all $J \supset I$, and
\item \label{g-U-non-zero} there is no common zero of $g_1, \ldots, g_k$ on $\Kstari \setminus U$.
\end{prooflist}
Since $I \not\in \excU \cup \excA$, \cref{bkk-0,intersection-lemma} imply that $\scrL'(\scrA)$ contains a nonempty Zariski open subset of $\scrL(\scrA)$. It follows that $\scrL'_a(\scrA) := \{(g_1, \ldots, g_n) \in \scrL'(\scrA): \In_{\Ai_j, \nu}(g_j|_{\Ki})(a) \neq 0$ for each $j = 1, \ldots, k\}$ also contains a nonempty Zariski open subset of $\scrL(\scrA)$. Since $\scrA$ is $\ri$-independent, it follows that the ($k$-dimensional) mixed volume of $\conv(\scrA_j) \cap \ri$, $j = 1, \ldots, k$, is nonzero. Due to \ref{g-U-non-zero}, the arguments of \cref{claim:dominant} then imply that we can find $(g_1, \ldots, g_n) \in \scrL'_a(\scrA)$ and a common zero $b$ of $g_1, \ldots, g_n$ on $\Kstari \cap U$ such that $f_j(b) \neq 0$ for each $j = 1, \ldots, k$. Property \ref{g-bkk-non-deg} together with \cref{branch-lemma-IB} then imply that $b$ must be an isolated zero of $g_1, \ldots, g_n$ on $U$.
\end{proof}

Now we follow the process from \cref{modified-proof}. Fix integers $\nu'_j > \nu_j := \nu(x_j)$. Let $C$ be the rational curve on $\kk^n$ parametrized by $c(t):= (c_1(t), \ldots, c_n(t)) :\kk \to \kk^n$ given by
\begin{align*}
c_j(t) :=
\begin{cases}
 a_jt^{\nu_j} + (b_j - a_j)t^{\nu'_j}
 	& \text{if}\ 1 \leq j \leq k,\\
 0
 	& \text{otherwise.}
\end{cases}
\end{align*}
Let $m_j := \min_{\Ai_j}(\nu)$, $1 \leq j \leq k$. Define
\begin{align*}
h_j &:=
\begin{cases}
 t^{-m_j}f_j(c(t))g_j - t^{-m_j}g_j(c(t))f_j
 	& \text{if}\ 1 \leq j \leq k,\\
 (1-t)f_j + tg_j
 	& \text{otherwise.}
\end{cases}
\end{align*}
Assertion \eqref{g-a-non-zero} of \cref{case-degenerate-1} implies that each $h_j|_{t=0}$ is a nonzero constant times $f_j$, and assertion \eqref{g-b-zero} of \cref{case-degenerate-1} implies that each $h_j|_{t=1}$ is a nonzero constant times $g_j$ and $c(1) = b$ is an isolated zero of $g_1, \ldots, g_n$ on $U$. The assumptions of case \eqref{case-nonempty} of \cref{case-degenerate} implies that the center of $C$ at $t = 0$ is either out of $U$, or it is a non-isolated zero of $f_1, \ldots, f_n$. Since each $h_j$ vanishes on the curve $C' := \{(c(t),t): t \in \kk\} \subset \kk^{n+1}$, assertion \eqref{non-max-condition} of \cref{mult-deformation-global} implies that $\multfiso{U} < \multntorusiso{h_1|_{t=\epsilon}}{h_n|_{t=\epsilon}}$ for generic $\epsilon \in \kk$. It follows that $\multfiso{U} < \multAiso{U}$, as required.\\

\subsubsection{Proof that $\scrN(U, \scrA)$ is a nonempty Zariski open subset of $\scrL(\scrA)$}
As in \cref{bkk-existential-section} we call $\scrB = (\scrB_1, \ldots, \scrB_n)$ a {\em face} of $\scrA$ and write $\scrB \preceq \scrA$, if there is $\nu \in \rnstar$ such that $\scrB_j = \In_\nu(\scrA_j)$ for each $j$; in that case we also write $\scrB := \In_\nu(\scrA)$ and we say that $\scrB$ is {\em centered at infinity} (respectively, {\em centered at $\Kstari$} for some $I \subseteq [n]$) if $\nu$ is centered at infinity (respectively, centered at $\Kstari$). Moreover, if $g$ is a polynomial supported at $\scrA_j$, we write $g_{\scrB_j}$ for $\In_{\scrA_j, \nu}(g)$. Consider the systems of polynomials that violate either property \ref{iUA-property} or property \ref{dUA-property} of $(U, \scrA)$-non-degeneracy: pick $I \subseteq [n]$ and a weighted order $\nu$ on $\kk[x_i: i \in I]$. Let $\scrB := \In_\nu(\Ai)$ and $\DBI$ be the set of all $(g_1, \ldots, g_n) \in \scrL(\scrA)$ such that there is a common root of $(g_j|_{\Ki})_{\scrB_j}$, $j = 1, \ldots, n$, on $\nktorus$. If $(f_1, \ldots, f_n) \in \scrL(\scrA)$ is in the closure of $\DBI$, then \cref{DJBbar-closed} implies that $(f_1, \ldots, f_n) \in \DBprimeI$ for some $\scrB' \preceq \scrB$. Note that
\begin{prooflist}
\item \label{B-at-infinity} If $\scrB$ is centered at infinity, then $\scrB'$ is also centered at infinity.
\item If $\scrB$ is centered at $\Kstariprime$ for some $I' \in  \excU \cup \excA$, then $\scrB'$ is also centered at $\Kstarii{I''}$ for some $I'' \in  \excU \cup \excA$.
\item \label{B-d-B'-I} If $\scrB$ is centered at $\Kstariprime$ for some $I' \subseteq [n]$, then for each $j = 1, \ldots, n$,
\begin{align*}
f_j|_{\Kiprime}
    &=
    \begin{cases}
    f_{j, \scrB_j}
        &\text{if}\ \scrB_j \cap \riprime \neq \emptyset,\\
    0 &\text{otherwise.}
    \end{cases}
\end{align*}
It follows that
\begin{prooflist}
\item if $I' \in \dUA$ and $\scrB'$ is centered at $\Kstarii{I''}$ for some $I'' \in \iUA$, then $(f_1, \ldots, f_n)$ violate property \ref{dUA-property} of $(U, \scrA)$-non-degeneracy with $I, I'$ replaced respectively by $I',I''$;
\item if $I' \in \iUA$ and $\scrB'$ is centered at infinity, then $(f_1, \ldots, f_n)$ violate property \ref{iUA-property} of $(U, \scrA)$-non-degeneracy with $I, I'$ replaced respectively by $I',I''$;
\item if $I' \in \iUA$ and $\scrB'$ is centered at $\Kstarii{I''}$ for some $I'' \in \dUA$, then $(f_1, \ldots, f_n)$ violate property \ref{iUA-property} of $(U, \scrA)$-non-degeneracy with $I, I'$ replaced respectively by $I',I''$;
\end{prooflist}
\end{prooflist}
It follows from these observations that the set of systems which violate at least one of the properties \ref{iUA-property} and \ref{dUA-property} of $(U, \scrA)$-non-degeneracy is Zariski closed in $\scrL(\scrA)$. Now we tackle property \ref{U-property}. Pick $I \in \iUA$ and let $\DUI$ be the set of all $(g_1, \ldots, g_n) \in \scrL(\scrA)$ such that there is a common root of $g_1, \ldots, g_n$ on $\Kstari \setminus U$. If $(f_1, \ldots, f_n) \in \scrL(\scrA)$ is in the closure of $\DUI$, then the arguments from the proof of \cref{DJBbar-closed} imply that there is a weighted order $\nu$ on $\kk[x_i: i \in I]$ and a common zero $a$ of $\In_{\Ai_j, \nu}(f_j|_{\Ki})$, $j = 1, \ldots, n$, on $\nktorus$, such that
\begin{prooflist}[resume]
\item \label{U-nu-infty} either $\nu$ is centered at infinity, in which case $f_1, \ldots, f_n$ violate property \ref{iUA-property} of $(U, \scrA)$-non-degeneracy,
\item \label{U-nu-I'} or $\nu$ is centered at $\Kstariprime$ for some $I' \subseteq I$, and $\pi_{I'}(a) \in \Kstariprime \setminus U$; in this case $f_1, \ldots, f_n$ violate property \ref{U-property} (with $I$ replaced by $I'$) of $(U, \scrA)$-non-degeneracy if $I' \in \iUA$, and property \ref{iUA-property} of $(U, \scrA)$-non-degeneracy if $I' \in \dUA$.
\end{prooflist}
It follows that the set of $(U, \scrA)$-degenerate systems is Zariski closed in $\scrL(\scrA)$, as required. \qed

\section{Weighted B\'ezout theorem: general version} \label{weighted-application}
\subsection{Weighted B\'ezout theorem II: all \bm{$\omega(f_j)\geq 0$}}%
For weights $\omega = (\omega_1, \ldots, \omega_n)$ to be applicable in weighted B\'ezout theorem (\cref{wt-bezout-1}), each $\omega_i$ has to be positive. In this section we replace this condition by a weaker one - that $\omega(f_j)$ has to be nonnegative for each $j$. This opens up a new possibility: if some $\omega_i$ is nonpositive, then the set of polynomials $f$ with $\omega(f)$ bounded above by a given integer will be an infinite dimensional vector space over $\kk$, and the number of (isolated) solutions of $n$ such polynomials can be arbitrarily large. Therefore to estimate number of solutions one has to bound the degree in each $x_i$ such that $\omega_i \leq 0$. It is natural then to consider, given an integer $d$, and a nonnegative integer $m_i$ for each $i$ such that $\omega_i \leq 0$, the set of polynomials supported at the polytope
\begin{align*}
\scrP(\omega, d, \vec m)
	&:= \{
		\alpha \in \rr^n: \langle \omega, \alpha \rangle \leq d,\
		\alpha_i \geq 0,\ i = 1, \ldots, n,\
		\alpha_i \leq m_i,\ i \in I_0 \cup I_-
		\}
\end{align*}
where $I_- := \{i: \omega_i < 0\}$ and $I_0 := \{i : \omega_i = 0\}$ (see \cref{fig:pomegadm}). The reason for our restriction to the case of nonnegative $\omega(f_j)$ is the following observation:

\begin{figure}[h]
\def\colorone{red}
\def\colortwo{blue}
\def\colorg{green}

\def\shiftx{7.5}
\def\colorzero{green}
\def\colorone{orange}
\def\colortwo{blue}
\def\colorthree{red}
\def\opazero{0.5}
\def\viewx{60}
\def\viewy{30}
\def\titlex{2}
\def\titley{-1}
\def\scalefactor{0.6}
\begin{center}
\begin{subfigure}[b]{0.3\textwidth}
\begin{tikzpicture}[scale=\scalefactor]
\pgfplotsset{every axis title/.append style={at={(0,-0.2)}}, view={\viewx}{\viewy}, axis lines=middle, enlargelimits={upper}}

\begin{axis}
\addplot3[ thick, draw, fill=\colorone,opacity=\opazero] coordinates{(0,1,0) (2,1,0) (2,1,1) (0,1,1)};
\addplot3[ thick, draw, fill=\colortwo,opacity=\opazero] coordinates{(1,0,0) (2,1,0) (2,1,1) (1,0,1)};
\addplot3[ thick, draw, fill=\colorthree,opacity=\opazero] coordinates{(0,0,0) (1,0,0) (1,0,1) (0,0,1)};
\addplot3[ thick, draw, fill=\colorg,opacity=\opazero] coordinates{(1,0,1) (0,0,1) (0,1,1) (2,1,1)};
\end{axis}
\end{tikzpicture}
\caption{$\omega = (1,-1,0)$,\\ $d = 1$, $m_2 = m_3 = 1$}  \label{fig:pomega-p0n}
\end{subfigure}
\begin{subfigure}[b]{0.3\textwidth}
\begin{tikzpicture}[scale=\scalefactor]
\pgfplotsset{every axis title/.append style={at={(0,-0.2)}}, view={\viewx}{\viewy}, axis lines=middle, enlargelimits={upper}}

\begin{axis}
\addplot3[ thick, draw, fill=\colorone,opacity=\opazero] coordinates{(3,1,1) (2,1,0) (0,1,0) (0,1,1)};
\addplot3[ thick, draw, fill=\colortwo,opacity=\opazero] coordinates{(1,0,0) (2,0,1) (3,1,1) (2,1,0)};
\addplot3[ thick, draw, fill=\colorg,opacity=\opazero] coordinates{(3,1,1) (2,0,1) (0,0,1) (0,1,1)};
\addplot3[ thick, draw, fill=\colorthree,opacity=\opazero] coordinates{(0,0,0) (1,0,0) (2,0,1) (0,0,1)};
\end{axis}
\end{tikzpicture}
\caption{$\omega = (1,-1,-1)$,\\ $d = 1$, $m_2 = m_3 = 1$}  \label{fig:pomega-pnn}
\end{subfigure}
\begin{subfigure}[b]{0.3\textwidth}
\begin{tikzpicture}[scale=\scalefactor]
\pgfplotsset{every axis title/.append style={at={(0,-0.2)}}, view={\viewx}{\viewy}, axis lines=middle, enlargelimits={upper}}

\begin{axis}
	\addplot3 [draw, ultra thick, fill=\colortwo,opacity=\opazero] coordinates{(0,0,1) (0,1,0) (1,0,0) (0,0,1)};
	\addplot3 [draw, ultra thick, fill=\colorthree,opacity=\opazero] coordinates{(0,0,0) (1,0,0) (0,0,1) (0,0,0)};
\end{axis}
\end{tikzpicture}
\caption{$\omega = (1,1,1)$,\\ $d = 1$}  \label{fig:pomega-ppp}
\end{subfigure}
\caption{$P(\omega,d,\vec m)$ for different $\omega,d,\vec m$}  \label{fig:pomegadm}
\end{center}
\end{figure}

\begin{prop} \label{sum-pomega}
Let $d_1, \ldots, d_k \in \zz$ and $m_{j,i} \in \zzero$, $j = 1, \ldots, k$, $i \in I_0 \cup I_-$. If $d_j \geq 0$ for each $j$, then
\begin{align}
\sum_j \scrP(\omega, d_j, \vec m_j)
	= \scrP(\omega, \sum_j d_j, \sum_j \vec m_j)
\label{sumpomega}
\end{align}
\end{prop}

\begin{proof}
$\scrP(\omega, d_j, \vec m_j)$ is the product of an $(n-|I_0|)$-dimensional polytope of the same form with the $I_0$-dimensional box formed by the product over all $i \in I_0$ of the closed intervals from $0$ to $m_{j,i}$ on $x_i$-axis. It is then straightforward to see that to prove \eqref{sumpomega} it suffices to prove it under the additional condition that $I_0 = \emptyset$. So assume $I_0 = \emptyset$. Let $m := \max_{i,j}m_{j,i}$ and define
\begin{align*}
x'_i	&:=
	\begin{cases}
		x_i & \text{if}\ i \not\in I_- \\
		m-x_i & \text{if}\ i \in I_-
	\end{cases}\\
\omega'_i &:= |\omega_i|,\ i = 1, \ldots, n.%
\end{align*}
It is straightforward to check that in $(x'_1, \ldots, x'_n)$-coordinates, in the notation of \cref{exercise:omegadm}, up to a reordering of the $x_i$ if necessary,
\begin{align}
\scrP(\omega, d_j, \vec m_j)
	&= \scrQ(\vec \omega', d'_j, \vec m_j) + (m'_{j,1}, \ldots, m'_{j,n})
\label{pomega=qomega'}
\end{align}
where
\begin{align*}
d'_j
	&:= d_j + \sum_{i \in I_-}|\omega_i|m_i,\ \text{and}\\
m'_{j,i}
	&:=
	\begin{cases}
		0 & \text{if}\ i \not\in I_- \\
		m - m_{j,i} & \text{if}\ i \in I_-
	\end{cases}
\end{align*}
It then follows from \cref{exercise:omegadm} that
\begin{align*}
\sum_j \scrP(\omega, d_j, \vec m_j)
	= \scrQ(\vec \omega', \sum_j d'_j, \sum_j \vec m_j) + \sum_j (m'_{j,1}, \ldots, m'_{j,n})
	=  \scrP(\omega, \sum_j d_j, \sum_j \vec m_j)
\end{align*}
as required.
\end{proof}

\begin{rem}
\eqref{sumpomega} may fail to hold if $d_j < 0$ for some $j$. Indeed, it follows from \cref{exercise:omegadm} and \eqref{pomega=qomega'} that with $\omega = (-1, -1, 1)$, $m_{1,1} = m_{1,2} = 1$ and $m_{2,1} = m_{2,2} = 3$, $\scrP(\omega, 1, \vec m_1) + \scrP(\omega, -3, \vec m_1) \subsetneq \scrP(\omega, -2, \vec m_1 + \vec m_2)$.
\end{rem}
We will now compute the mixed volume of $\scrP(\omega, d_j, \vec m_j)$, $j = 1, \ldots, n$. \Cref{exercise:vol-omegadm} and \eqref{pomega=qomega'} imply that
\begin{align}
\vol(\scrP(\omega, d, \vec m))
	&=  \frac{\prod_{i \in I_0}m_i}{(n-|I_0|)!\prod_{i \not\in I_0} |\omega_i|} 	
	   \sum_{\substack{I \subseteq I_- \\ d + \sum_{i \in I}|\omega_i| m_i > 0}} (-1)^{|I_-| - |I|}
	   (d + \sum_{i \in I}|\omega_i| m_i)^{n - |I_0|}
\label{volume-pomega}
\end{align}
In particular, if $d \geq 0$, then
\begin{align*}
\vol(\scrP(\omega, d, \vec m))
	&=  \frac{\prod_{i \in I_0}m_i}{(n-|I_0|)!\prod_{i \not\in I_0} |\omega_i|} 	
	   \sum_{I \subseteq I_- } (-1)^{|I_-| - |I|}
	   (d + \sum_{i \in I}|\omega_i| m_i)^{n - |I_0|}
\end{align*}
Therefore \cref{sum-pomega} implies that for each $d_1, \ldots, d_n, \lambda_1, \ldots, \lambda_n \geq 0$,
\begin{align*}
\vol(\sum_{j=1}^n \lambda_j\scrP(\omega, d_j, \vec m_j))
	&=  \frac{\prod_{i \in I_0}( \sum_{j=1}^n \lambda_j m_{j,i})}
		{(n-|I_0|)!\prod_{i \not\in I_0} |\omega_i|} 	
	   \sum_{I \subseteq I_- } (-1)^{|I_-| - |I|}
	   ( \sum_{j=1}^n \lambda_j (d_j + \sum_{i \in I}|\omega_i| m_{j,i}))^{n - |I_0|}
\end{align*}
The mixed volume of $\scrP(\omega, d_1, \vec m_1), \ldots, \scrP(\omega, d_n, \vec m_n)$ is the coefficient of $\lambda_1\cdots \lambda_n$ in the right hand side of the above expression. For each $J_0 \subseteq [n]$ such that $|J_0| = |I_0|$, the coefficient of $\prod_{j \in J_0} \lambda_j$ in $\prod_{i \in I_0}(\sum_{j=1}^n \lambda_j m_{j,i})$ is the {\em permanent} of the $|I_0| \times |I_0|$ matrix
\begin{align}
D_{I_0,J_0}
	&:=
	  \begin{pmatrix}
	  m_{j_1,i_1} & \cdots & m_{j_1, i_k}\\
	  \vdots & & \vdots \\
	  m_{j_k, i_1} & \cdots & m_{j_k,i_k}
	  \end{pmatrix}
\label{DI0J0}
\end{align}
where $k := |I_0| = |J_0|$ and $i_1, \ldots, i_k$ (respectively, $j_1, \ldots, j_k$) are elements of $I_0$ (respectively, $J_0$). (Note that $\perm(D_{I_0,J_0})$ does not depend on the ordering of the elements of $I_0$ or $J_0$. If $I_0 = J_0 = \emptyset$, then $D_{I_0,J_0}$ is the empty matrix, and its permanent is by convention $1$.) On the other hand, the coefficient of $\prod_{j \not\in J_0} \lambda_j$ in $(\sum_{j=1}^n \lambda_j (d_j + \sum_{i \in I}|\omega_i| m_{j,i}))^{n - k}$ is $(n-k)!\prod_{j \not\in J_0}(d_j + \sum_{i \in I}|\omega_i| m_{j,i})$. Combining all these together yields:

\begin{prop}\label{prop:mvpomega}
Let $d_1, \ldots, d_n \in \zz$ and $m_{j,i} \in \zzero$, $j = 1, \ldots, k$, $i \in I_0 \cup I_-$. If $d_j \geq 0$ for each $j$, then the mixed volume of $\scrP(\omega, d_1, \vec m_1), \ldots, \scrP(\omega, d_n, \vec m_n)$ is
\begin{align*}
\sum_{\substack{J_0 \subseteq [n]\\|J_0| = |I_0|}} \perm(D_{I_0,J_0})
	\sum_{I \subseteq I_-} (-1)^{|I_-| - |I|}
		\cfrac{
			\prod_{j \in [n] \setminus J_0} (d_j + \sum_{i \in I_-}|\omega_i|m_{j,i})
		}{
			\prod_{i \in [n] \setminus I_0}|\omega_i|}
\end{align*}
\end{prop}

The following result describes the faces of $\scrP(\omega, d, \vec m)$. Its proof is left as an exercise. Let $\eta \in \rnstar$ with coordinates $(\eta_1, \ldots, \eta_n)$ with respect to the basis dual to the standard basis of $\rr^n$. 

\begin{prop} \label{pomega-face}
Define $M := \sup\{\eta_i/\omega_i: \omega_i > 0\}$. Assume $d \geq 0$. Then
\begin{enumerate}
\item If $M \leq 0$ (which is the case if in particular $\omega_i \leq 0$ for each $i$), then $\ld_\eta(\scrP(\omega, d, \vec m))$ is the set of all $(\alpha_1, \ldots, \alpha_n) \in \scrP$ such that
\begin{align*}
\alpha_i
	&=
	\begin{cases}
	0   &\text{if}\ \eta_i < 0,\\
    m_i &\text{if}\ \eta_i > 0.
    \end{cases}
\end{align*}
\item If $M > 0$, then $\ld_\eta(\scrP(\omega, d, \vec m))$ is the set of all $(\alpha_1, \ldots, \alpha_n) \in \scrP$ such that $\sum_i \alpha_i \omega_i = d$ and
\begin{align*}
\alpha_i
	&=
	\begin{cases}
	0   &\text{if}\ (\omega_i = 0,\ \eta_i < 0)\
            \text{or}\ (\omega_i > 0,\ \cfrac{\eta_i}{\omega_i} < M)\
            \text{or}\ (\omega_i < 0,\ \cfrac{\eta_i}{\omega_i} > M), \\
    m_i &\text{if}\ (\omega_i = 0,\ \eta_i > 0)\
            \text{or}\ (\omega_i < 0,\ \cfrac{\eta_i}{\omega_i} < M). \qed
    \end{cases}
\end{align*}
\end{enumerate}
\end{prop}

We are now ready to prove the second version of the weighted B\'ezout theorem. Let $f_j$ be polynomials supported at $\scrP(\omega, d_j, \vec m_j)$, $j = 1, \ldots, n$. For each $I \subseteq I_0 \cup I_-$, define
\begin{align*}
\scrL_{\vec m_j, I}
	&:=
	\{
	\alpha = (\alpha_1, \ldots, \alpha_n): \alpha_i = m_{j,i}\ \text{for each}\ i \in I
	\} \\
\scrL_{\omega, d_j, \vec m_j, I}
	&:= \{\alpha \in \scrL_{\vec m_j, I}: \sum_{i=1}^n \omega_i \alpha_i = d_j \}
\end{align*}
Denote the coefficient of $x^\alpha$ in $f_j$ by $c_{j\alpha}$, i.e.\ $f_j = \sum_{\alpha} c_{j\alpha}x^\alpha$. Write
\begin{align*}
\ld_{\vec m_j, I}(f_j)
	&:= \sum_{\alpha \in \scrL_{\vec m_j, I}} c_{j\alpha}x^\alpha,\\
\ld_{\omega, d_j, \vec m_j, I}(f_j)
	&:= \sum_{\alpha \in \scrL_{\omega, d_j, \vec m_j, I}} c_{j\alpha}x^\alpha
\end{align*}

\begin{thm}[Weighted B\'ezout theorem II] \label{wt-bezout-21}
\index{Weighted!B\'ezout theorem!version II}
\index{B\'ezout's theorem!weighted!version II}
Given polynomials $f_1, \ldots, f_n$ in $(x_1, \ldots, x_n)$ and any weighted degree $\omega$ on $\kk[x_1, \ldots, x_n]$, the number $N$ of isolated solutions of $f_1, \ldots, f_n$ on $\kk^n$ satisfies:
\begin{align}
N
	&\leq
		\sum_{\substack{J_0 \subseteq [n]\\|J_0| = |I_0|}} \perm(D_{I_0,J_0})
		\sum_{I \subseteq I_-} (-1)^{|I_-| - |I|}
		\cfrac{
			\prod_{j \in [n] \setminus J_0} \left(\max\{\omega(f_j),0\} + \sum_{i \in I_-}|\omega_i|\deg_{x_i}(f_j) \right)
		}{
			\prod_{i \in [n] \setminus I_0}|\omega_i|}
\label{wt-II}
\end{align}
where $D_{I_0,J_0}$ is defined as in \eqref{DI0J0} with $m_{j,i} := \deg_{x_i}(f_j)$. More generally, if $d_1, \ldots, d_n$ and $m_{j,i}$ are {\em nonnegative} integers such that each $d_j \geq \omega(f_j)$ and $m_{j,i} \geq \deg_{x_i}(f_j)$, then
\begin{align}
N
	&\leq
		\sum_{\substack{J_0 \subseteq [n]\\|J_0| = |I_0|}} \perm(D_{I_0,J_0})
		\sum_{I \subseteq I_-} (-1)^{|I_-| - |I|}
		\cfrac{
			\prod_{j \in [n] \setminus J_0} \left(d_j + \sum_{i \in I_-}|\omega_i|m_{j,i} \right)
		}{
			\prod_{i \in [n] \setminus I_0}|\omega_i|}
\label{wt-II-1}
\end{align}
If the right hand side of \eqref{wt-II-1} is nonzero, i.e.\ $\scrP(\omega, d_1, \vec m_1), \ldots, \scrP(\omega, d_n, \vec m_n)$ are independent (which is the case e.g.\ if each $d_j$ and each $m_{j,i}$ is positive), then \eqref{wt-II-1} holds with an equality if and only if both of the following are true:
\begin{enumerate}
\item \label{I-condition} For each {\em nonempty} subset $I$ of $I_0 \cup I_-$, $\ld_{\vec m_j, I}(f_j)$, $j = 1, \ldots, n$, have no common zero on $\kk^n \setminus \bigcup_{i \in I} V(x_i)$.
\item \label{omegaI-condition} For each (possibly empty) subset $I$ of $I_0 \cup I_-$, $\ld_{\omega, d_j, \vec m_j, I}(f_j)$, $j = 1, \ldots, n$, have no common zero on $\kk^n \setminus \left(\cup_{i \in I} V(x_i)) \bigcup (\cap_{i \in I_+}V(x_i)) \right)$ (where $I_+ := \{i: \omega_i > 0\}$).
\end{enumerate}
\end{thm}

\begin{proof}
Inequalities \eqref{wt-II} and \eqref{wt-II-1} follow directly from \cref{extended-bkk-bound-0,prop:mvpomega}. To find the non-degeneracy conditions we apply \cref{non-deg-1,non-deg-2} with $\scrA := (\scrP(\omega, d_1, \vec m_1) \cap \zz^n, \ldots, \scrP(\omega, d_n, \vec m_n) \cap \zz^n)$ and $U := \kk^n$. It is straightforward to check that $\excA = \excU = \emptyset$, and if $\scrP(\omega, d_j, \vec m_j)$, $j = 1, \ldots, n$, are independent then $\dUA = \emptyset$ and $\iUAstar = [n]$. Therefore \cref{non-deg-1,non-deg-2} imply that if $\scrP(\omega, d_j, \vec m_j)$, $j = 1, \ldots, n$, are independent (which due to \cref{positively-mixed,prop:mvpomega} is equivalent to the right hand side of \eqref{wt-II-1} being nonzero), then \eqref{wt-II-1} holds with an equality if and only if the following holds:
\begin{align}
\parbox{0.69\textwidth}{
for each weighted order $\nu$ centered at infinity on $\kk[x_1, \ldots, x_n]$, there is no common zero of $\In_{\scrA_1, \nu}(f_1), \ldots, \In_{\scrA_n,\nu}(f_n)$ on $\nktorus$.
}\label{wt-non-deg-20}
\end{align}
Pick a weighted order $\nu$ centered at infinity on $\kk[x_1, \ldots, x_n]$. Now apply \cref{pomega-face} to $\eta := -\nu$. If $M := \sup\{\eta_i/\omega_i: \omega_i > 0\} \leq 0$, then $I := \{i \in I_0 \cup I_-: \eta_i > 0\}$ is nonempty and \cref{pomega-face} implies that
\begin{align*}
\In_{\scrA_j, \nu}(f_j)
    &= \ld_{\vec m_j, I}(f_j)|_{\{x_{i'} = 0: i' \in I'\}}
\end{align*}
where $I' := \{i \in [n]: \eta_i < 0\}$. Therefore \eqref{wt-non-deg-20} is equivalent to the condition that $\ld_{\vec m_j, I}(f_j)$, $j = 1, \ldots, n$, have no common zero on $\left( \kk^n \setminus \bigcup_{i \in [n] \setminus I'} V(x_i)\right) \cap \bigcap_{i' \in I'} V(x_{i'})$. Since it is possible for $I'$ to be any subset of $[n]\setminus I$, taking into account all such possibilities leads to condition \eqref{I-condition}. Now consider the case that $M > 0$. In this case define $I := \{i \in I_0: \eta_i > 0\} \cup \{i \in I_-: \eta_i/\omega_i < M\}$ and $I' := \{i \in I_0: \eta_i < 0\} \cup \{i \in I_-: \eta_i/\omega_i > M\} \cup \{i \in I_+: \eta_i/\omega_i < M\}$. \Cref{pomega-face} implies that
\begin{align*}
\In_{\scrA_j, \nu}(f_j)
    &= \ld_{\omega, d_j, \vec m_j, I}(f_j)|_{\{x_{i'} = 0: i' \in I'\}}
\end{align*}
It then follows as in the preceding case that to satisfy \eqref{wt-non-deg-20} for all choices of $I'$ is equivalent to condition \eqref{omegaI-condition}, and this completes the proof.
\end{proof}

\begin{rem}
One does not really need the main results of this chapter to establish the bounds \eqref{wt-II} and \eqref{wt-II-1} of \cref{wt-bezout-21} - these can be established in the same way as the proof of the classical weighted B\'ezout bound from \cref{wt-bezout-proof-section} once the mixed volume computation from \cref{prop:mvpomega} is available. However, we used \cref{non-deg-1,non-deg-2} in an essential way to establish the nondegeneracy conditions of \cref{wt-bezout-21}.
\end{rem}

\subsection{Weighted B\'ezout theorem III: the general case}%
In this section we consider weighted B\'ezout theorem without any restriction on the $\omega_i$ or $\omega(f_j)$. In this generality we do not know of any compact expressions for either the mixed volume or the faces of $\scrP(\omega, d_j, \vec m_j)$. Therefore our version of general weighted B\'ezout theorem below is little more than direct application of the extension of Bernstein's theorem to the affine space.

\begin{thm}[Weighted B\'ezout theorem III] \label{wt-bezout-3}
\index{Weighted!B\'ezout theorem!version III}
\index{B\'ezout's theorem!weighted!version III}
Let $d_j\in \zz$ and $m_{j,i} \in \zzero$ be such that $d_j \geq \omega(f_j)$ for each $j$ and each $m_{j,i} \geq \deg_{x_i}(f_j)$. Given $I \subseteq I_-$, write $\badj(I) := \{j: \sum_{i \in I_-\setminus I} m_{j,i}\omega_i > d_j\}$. Let $\goodI := \{I \subseteq I_-: |\badj(\tilde I)| \leq |\tilde I|$ for each $\tilde I \subseteq I\}$. Then the number $N$ of isolated solutions of $f_1, \ldots, f_n$ on $\kk^n$ is
\begin{align}
N \leq
	\sum_{I \in \goodI, |\badj(I)| = |I|}
	\mv_{n-k}(\scrP_{j'_1} \cap \riprime, \ldots, \scrP_{j'_{n-k}} \cap \riprime)
	\times \multzero{\pi_I(\scrP_{j_1})}{\pi_I(\scrP_{j_k})}
\label{general-wt-formula}
\end{align}
where
\begin{itemize}
\item $\scrI' := [n]\setminus I$,
\item $\scrP_j$ are short for $\scrP(\omega, d_j, \vec m_j)$,
\item $j_1, \ldots, j_k$ (respectively, $j'_1, \ldots, j'_{n-k}$) are elements of $\badj(I)$ (respectively, $\{1, \ldots, n\} \setminus \badj(I)$), and
\item $\mv_{n-k}(\cdot, \ldots, \cdot)$ is the $(n-k)$-dimensional mixed volume.
\end{itemize}
The bound in \eqref{general-wt-formula} holds with an equality if and only if $f_1, \ldots, f_n$ are $(\kk^n, \scrA)$-non-degenerate, where $\scrA := (\scrP_1 \cap \zz^n, \ldots, \scrP_n \cap \zz^n)$.
\end{thm}

\begin{proof}
Follows immediately from \cref{extended-bkk-bound-0,non-deg-1,non-deg-2}.
\end{proof}

\section{Weighted multi-homogeneous B\'ezout theorem: general version} \label{weighted-multapplication}
In this section we generalize the weighted multi-homogeneous version of B\'ezout's theorem (\cref{multi-bezout}) by replacing the weighted degrees with positive weights by weighted degrees from \cref{wt-bezout-21}. As in the setting of \cref{multi-bezout}, let $\mscrI := (I_1, \ldots, I_s)$ be an ordered partition of $[n] := \{1, \ldots, n\}$, and for each $j = 1, \ldots, s$, let $\omega_j$ be a weighted degree on $\kk[x_k: k \in I_j]$ with weights $\omega_{j,k}$ for $x_k$, $k \in I_j$. Let $I_{j,+}$ (respectively, $I_{j,0}, I_{j,-}$) be the set of all $k \in I_j$ such that $\omega_{j,k} > 0$ (respectively, $\omega_{j, k} = 0$, $\omega_{j,k} < 0$). Given nonnegative integers $d_{i,j} \geq \omega_j(f_i)$ and $m_{i,j,k} \geq \deg_{x_k}(f_i)$ for each $k \in I_{j,0} \cup I_{j,-}$, we consider the polytope
\begin{align*}
\scrP_i
    &:= \prod_{j = 1}^s \scrP(\omega_j, d_{i,j}, \vec m_{i,j})
\end{align*}
Let $n_j := |I_j|$, $j = 1, \ldots, s$. The mixed volume of $\scrP_1, \ldots, \scrP_n$ is the coefficient of $\lambda_1 \cdots \lambda_n$ in the polynomial
\begin{align*}
\vol_n(\sum_{i=1}^n \lambda_i \scrP_i)
	&= \vol_n(\sum_{i=1}^n \lambda_i \prod_{j=1}^s \scrP(\omega_j, d_{i,j}, \vec m_{i,j}))
	= \vol_n(\prod_{j=1}^s (\sum_{i=1}^n \lambda_i \scrP(\omega_j, d_{i,j}, \vec m_{i,j}))) \\
    &= \vol_n(\prod_{j=1}^s \scrP(\omega_j, \sum_{i=1}^n \lambda_i d_{i,j}, \sum_{i=1}^n \lambda_i \vec m_{i,j}))
    = \prod_{j=1}^s \vol_{n_j}(\scrP(\omega_j, \sum_{i=1}^n \lambda_i d_{i,j}, \sum_{i=1}^n \lambda_i \vec m_{i,j}))
\end{align*}
where the third equality follows from \cref{sum-pomega}. After a refinement of $\mscrI$ if necessary, we may, and will, assume that $\omega_{j,k} \neq 0$ for each $j,k$, i.e.\ $I_j = I_{j,+} \cup I_{j,-}$ for each $j$. Then \eqref{volume-pomega} implies that
\begin{align*}
\vol_n(\sum_{i=1}^n \lambda_i \scrP_i)
    &= \prod_{j=1}^s \left(
            \frac{1}{n_j!\prod_{k \in I_j} |\omega_{j,k}|} 	
        	   \sum_{I \subseteq I_{j,-} } (-1)^{|I_{j,-}| - |I|}
        	   ( \sum_{i=1}^n \lambda_i (d_{i,j} + \sum_{k \in I}|\omega_{j,k}| m_{i,j,k}))^{n_j}
            \right)
\end{align*}
Let $I_- := \bigcup_j I_{j,-}$. For each $I \subseteq I_-$, write
\begin{align}
d_{I,i,j}
    &:= d_{i,j} + \sum_{k \in I \cap I_j}|\omega_{j,k}| m_{i,j,k}
\end{align}
Let $\Omega := \{\omega_1, \ldots, \omega_s\}$, and $D(\Omega, \vec d, \vec m, I)$ be the following $n \times n$ matrix:
\vspace{2mm}

\begin{align*}
\begin{matrix}
D(\Omega, \vec d, \vec m, I)
	&:=
	  \begin{pmatrix}
	  \bovermat{$n_1$ times}{ d_{I,1,1}  & \cdots & d_{I,1,1}}
	  	&
	  \overmat{$\mkern-3.5mu\cdots$}{\cdots & \cdots}
	  	&
	  \bovermat{$n_s$ times}{ d_{I,1,s}  & \cdots & d_{I,1,s}}
	   	\\[0.5em]
	   	\vdots & & \vdots & & & \vdots & & \vdots \\
	   	d_{I,n, 1}  & \cdots & d_{I,n,1} & \cdots & \cdots & d_{I,n,s}  & \cdots & d_{I,n,s}
	   \end{pmatrix}
\end{matrix}
\end{align*}
The preceding discussion together with \cref{extended-bkk-bound-0} or \cref{equal-multpn} imply that
\begin{align}
\multfkniso
    \leq \multPkniso
    = \mv(\scrP_1, \ldots, \scrP_n)
	= \sum_{I \subseteq I_-} (-1)^{|I_{-}| - |I|}
        \frac{ \perm(D(\Omega, \vec d, \vec m, I))}
		 {(\prod_j n_j!) (\prod_{j,k}\omega_{j,k})}
		 \label{wt-multi-homogeneous-bound-II}
\end{align}
It is straightforward to check that the bound \eqref{wt-II-1} from weighted B\'ezout theorem II corresponds to the special case of \eqref{wt-multi-homogeneous-bound-II} in which $n_j = 1$ for all but possibly one $j \in \{1, \ldots, s\}$. We will now see that the criterion for attainment of this bound is an amalgam of the non-degeneracy criteria of weighted multi-homogeneous B\'ezout theorem (\cref{multi-bezout}) and weighted B\'ezout theorem II (\cref{wt-bezout-21}). Given $I \subseteq I_-$, $J \subseteq [s]$, and $l \in [n]$, let $\scrL_{\Omega, \vec d, \vec m, I, J, i}$ be the set of all $\alpha = (\alpha_1, \ldots, \alpha_n) \in \znzero$ such that
\begin{itemize}
\item $\alpha_k = m_{i,j,k}$ for each $j \in [s]$ and $k \in I_{j,-} \cap I$, and
\item $\sum_{k \in I_j} \omega_{j,k}\alpha_k = d_{i,j}$ for each $j \in J$.
\end{itemize}
Given $g = \sum_\alpha c_\alpha x^\alpha \in \kk[x_1, \ldots, x_n]$ supported at $\scrP_i$, define
\begin{align*}
\ld_{\Omega, \vec d, \vec m, I, J, i}(g)
    &:= \sum_{\alpha \in \scrL_{\Omega, \vec d, \vec m, I, J, i}} c_\alpha x^\alpha
\end{align*}
We are now ready to prove version II of the weighted multi-homogeneous B\'ezout theorem. But at first we recall the assumptions:
\begin{defnlist}
\item $\omega_{j,k} \neq 0$ for each $j = 1, \ldots, s$, and $k \in I_j$.
\item $d_{i,j} \geq \max\{\omega_j(f_i), 0\}$ for each $i = 1, \ldots, n$, and $j = 1, \ldots, s$.
\item $m_{i,j,k} \geq \deg_{x_k}(f_i)$ for each $i = 1, \ldots, n$, and $j = 1, \ldots, s$, and $k \in I_{j,-}$.
\end{defnlist}

\begin{thm}[Weighted multi-homogeneous B\'ezout theorem II] \label{multi-bezout-II}
\index{Weighted!multi-homogeneous B\'ezout theorem!version II}
\index{B\'ezout's theorem!weighted multi-homogeneous!version II}
Under the above assumptions the number of isolated solutions of polynomials $f_1, \ldots, f_n$ on $\kk^n$ is bounded by \eqref{wt-multi-homogeneous-bound-II}. This bound is exact if and only if the following holds: for each pair $I,J$ such that $I \subseteq I_-$, $J \subseteq [s]$, and at least one of $I$ and $J$ is nonempty, there is no common zero of $\ld_{\Omega, \vec d, \vec m, I, J, 1}(f_1), \ldots, \ld_{\Omega, \vec d, \vec m, I, J, n}(f_n)$ on $\kk^n \setminus \left(\cup_{i \in I} V(x_i)) \bigcup (\cup_{j \in J}\cap_{k \in I_{j,+}}V(x_k)) \right)$.
\end{thm}

\begin{proof}
This follows from \cref{extended-bkk-bound-0,non-deg-1,non-deg-2} via arguments similar to those in the proof of \cref{wt-bezout-21}.
\end{proof}

\section{Open problems}
\subsection{Systems with isolated zeroes on a given coordinate subspace}
Given finite subsets $\scrA_1, \ldots, \scrA_n$ of $\znzero$ and a coordinate subspace $\Ki$ of $\kk^n$, it is straightforward to identify if for {\em generic} $f_j$ supported at $\scrA_j$, there is any common zero of $f_1, \ldots, f_n$ on $\Ki$, and in case there are such points, if they are isolated in $V(f_1, \ldots, f_n)$ or not. This is the content of the next result, which is straightforward to prove from \cref{positively-mixed,non-isolated-lemma}.

\begin{prop}[Zeroes of generic systems] \label{generic-characterization}
Let $\scrA := (\scrA_1, \ldots, \scrA_n)$. For all $f = (f_1, \ldots, f_n) \in \scrL(\scrA)$, write $V(f) := V(f_1, \ldots, f_n) \subset \kk^n$ and $\VstarIf := V(f_1, \ldots, f_n) \cap \Kstari$.
\begin{enumerate}
\item The following are equivalent:
\begin{enumerate}
\item $\VstarIf= \emptyset$ for generic $f \in \scrL(\scrA)$.
\item $\scrA$ is $\ri$-dependent.
\end{enumerate}
\item The following are equivalent:
\begin{enumerate}
\item $\VstarIf \neq \emptyset$ and all points of $\VstarIf$ are isolated in $V(f)$ for generic $f \in \scrL(\scrA)$.
\item $\scrA$ is $\ri$-independent and $I \not\in \excA$.
\end{enumerate}
\item The following are equivalent:
\begin{enumerate}
\item $\VstarIf \neq \emptyset$ and all points of $\VstarIf$ are non-isolated in $V(f)$ for generic $f \in \scrL(\scrA)$.
\item $\scrA$ is $\ri$-independent and $I \in \excA$. \qed
\end{enumerate}	
\end{enumerate}
\end{prop}

\begin{problem}[Existence of systems with given support and isolated zeroes on given coordinate subspace] \label{isolated-problem}
Characterize those $I \subseteq [n]$ for which there are $f_1, \ldots, f_n$ such that $\supp(f_j) = \scrA_j$, $j = 1, \ldots, n$, $\VstarI(f_1, \ldots, f_n) \neq \emptyset$ and all points (or some points) of $\VstarI(f_1, \ldots, f_n)$ are isolated in $V(f_1, \ldots, f_n)$.
\end{problem}

If $I$ is as in \cref{isolated-problem}, then
\begin{defnlist}
\item \cref{non-isolated-lemma} implies that $I \not\in \excA$;
\item \cref{bkk-non-degenerate-thm,positively-mixed} imply that
\begin{defnlist}
\item either $I$ is $\ri$-independent and $|\tiA| = |I|$ (i.e.\ $I \in \iUAstar$),
\item or $|\tiA| > |I|$ (in particular, $\scrA$ is $\ri$-dependent);
\end{defnlist}
\item \cref{hereditary-prop} implies that $\scrA$ is not hereditarily $\ri$-dependent.
\end{defnlist}
In the context of these observations \cref{isolated-problem} boils down to the following problem.

\addtocounter{thm}{-1}
\let\oldthethm\thethm
\renewcommand{\thethm}{\oldthethm$'$}
\begin{problem} \label{isolated-problem'}
If $I \not\in \excA$ and $\scrA$ is not hereditarily $\ri$-dependent and $|\tiA| > |I|$, does there exist $f_1, \ldots, f_n$ such that $\supp(f_j) = \scrA_j$ for each $j$, $\VstarI(f_1, \ldots f_n) \neq \emptyset$ and all points (or some points) of $\VstarI(f_1, \ldots, f_n)$ are isolated in $V(f_1, \ldots, f_n)$? If not, then characterize those $I \subseteq [n]$ which satisfy the hypothesis of the preceding question but fail the conclusion.
\end{problem}
\let\thethm\oldthethm

\subsection{Non-isolated zeroes and non-degeneracy}
In contrast to the case of $\nktorus$, \cref{ex-aff-non-deg-1} shows that for $f = (f_1, \ldots, f_n) \in \scrL(\scrA)$, the existence of non-isolated solutions of $f_1, \ldots, f_n$ does not automatically mean that $f_1, \ldots, f_n$ are $(\kk^n, \scrA)$-non-degenerate. More precisely, part \ref{aff-non-deg-1-1} of \cref{ex-aff-non-deg-1} shows that if $I \in \dUUA{\kk^n}$, then it is possible for $f_1, \ldots, f_n$ to be $(\kk^n,\scrA)$-non-degenerate even if $\VstarIf$ has non-isolated points. On the other hand, condition \ref{iUA-property} of $(U, \scrA)$-non-degeneracy implies that if $I \in \iUUA{\kk^n}$ and $\VstarIf$ has non-isolated points, then $f_1, \ldots, f_n$ are $(\kk^n,\scrA)$-degenerate. The question is if it is sufficient.

\begin{problem} \label{non-isolated-problem}
If $I \in \dUUA{\kk^n}$, does there exist $f = (f_1, \ldots, f_n)$ such that $\supp(f_j) = \scrA_j$ for each $j$, and $(V^*)^I(f)$ is nonempty (and due to \cref{hereditary-prop} necessarily positive dimensional), but $f_1, \ldots, f_n$ are $(\kk^n,\scrA)$-non-degenerate? If not, then characterize those $I \in \dUUA{\kk^n}$ for which there is no such $f \in \scrL(\scrA)$.
\end{problem}

\subsection{Simple criteria for equality of Li and Wang's bound}
Since the upper bound of Li and Wang from \eqref{mv-estimate} is so simple, it would be interesting to find simple criteria under which it holds with equality. \Cref{li-wang-condition} gives a characterization of all such scenarios, so the question is if it can be made ``more explicit'' in any sense, or if there are simpler criteria (e.g.\ as in \cref{conveniently-li-wang} or assertion \eqref{leq-n-2} of \cref{li-wang-condition}) which are sufficient. One possible criterion was proposed in \cite{rojas-wang} and \cite{rojas-toric}, but that turns out to be incorrect. Indeed, both \cite[Theorem 1]{rojas-wang} and \cite[Affine Point Theorem II]{rojas-toric} imply the following: if $\excA = \emptyset$ and the intersection of each $\scrA_j$ with each of the $n$ coordinate hyperplanes is nonempty, then Li and Wang's bound holds with equality. This is indeed the case for $n \leq 2$, but as the following example shows, it is false in higher dimensions.

\begin{example}
Let $f_1 := ax +by + cx^2$, $f_2 := a'x + b'y +c'x^2$, $f_3 := pz^kx + q$, where $a,b,c,a',b',c',p,q$ are generic elements in $\kk$ and $k \geq 1$. It is straightforward to check directly that on $\kk^n$ there are precisely $k$ solutions for $f_1 = f_2 = f_3 = 0$ and $2k$ solutions for $f_1 = f_2 = f_3 = t$ for generic $t \neq 0$, so that Li and Wang's bound fails for $\scrA_j := \supp(f_j)$, $j = 1, 2, 3$. Note that both conditions \ref{m_nu-condition} and \ref{dimension-condition} from \cref{li-wang-section} hold with the weighted degree $\nu$ on $\kk[x,y,z]$ corresponding to weights $x \mapsto k$, $y \mapsto k$, $z \mapsto -1$.
\end{example}

\subsection{``Compact'' formulae for general weighted and weighted multi-homogeneous versions of B\'ezout's theorem}
There are scenarios not covered in weighted B\'ezout theorem II (\cref{wt-bezout-21}) for which very similar bound exists, e.g.\ in the case that $|I_-| = 1$ (and no restriction that the $d_j$ have to be nonnegative). This motivates the question: is it possible to find a formula that is more explicit than the one from weighted B\'ezout theorem III (\cref{wt-bezout-3}), and which is more general than version II? That would also lead to a more general version of weighted multi-homogeneous B\'ezout theorem II (\cref{multi-bezout-II}). 


\chapter{Milnor number of a hypersurface at the origin} \label{milnor-chapter}
\newcommand{\milnor}{\mu_0}
\newcommand{\newton}{\nu}

\renewcommand{\partialxone}[1]{\partial_1 #1}
\renewcommand{\partialxtwo}[1]{\partial_2 #1}
\renewcommand{\partialxi}[1]{\partial_i #1}
\newcommand{\partialxioriginal}[1]{\partial #1/\partial x_i}
\renewcommand{\partialxifrac}[1]{\partial_i #1}
\renewcommand{\partialxii}[2]{\partial_{#2} #1}
\renewcommand{\partialxiifrac}[2]{\partial_{#2} #1}
\renewcommand{\partialxj}[1]{\partial_j #1}
\renewcommand{\partialxk}[1]{\partial_k #1}
\renewcommand{\partialxn}[1]{\partial_n #1}

\section{Introduction}
The modern theory of applications of Newton polyhedra to affine B\'ezout problem started from A.\ Kushnirenko's work aimed at answering V.\ I.\ Arnold's question on {\em Milnor numbers} of generic singularities. In \cite{kush-poly-milnor} Kushnirenko gave a beautiful formula for a lower bound of the Milnor number at the origin in terms of volumes of the region bounded by the Newton diagram, and showed that the bound is attained in the case that the singularity is {\em Newton non-degenerate}. In this chapter we 
show that the notion of {\em non-degeneracy at the origin} introduced in \cref{multiplicity-chapter} can be used to derive (and generalize) Kushnirenko's result on Milnor numbers. In particular, based on non-degeneracy at the origin we introduce a non-degeneracy criterion which generalizes Newton non-degeneracy and {\em inner Newton non-degeneracy}\footnote{This terminology is taken from \cite{boubakri-greuel-markwig}.}, the latter introduced by C.\ T.\ C.\ Wall \cite{wall}. We show that in zero characteristic the new criterion is necessary and sufficient for the Milnor number to be the minimum, and the minimum Milnor number can be obtained by Kushnirenko's bound. In positive characteristic this criterion is sufficient, but not necessary.

\section{Milnor number} \label{milnor-defn}

The \index{Milnor number}{\em Milnor number} $\milnor(f)$ of a power series $f$ in $(x_1, \ldots, x_n)$ is the dimension over $\kk$ of the quotient of $\kk[[x_1, \ldots, x_n]]$ by the ideal generated by the partial derivatives of $f$, i.e.\ $\milnor(f) = \multzero{\partialxone{f}}{\partialxn{f}}$ (where $\partialxi(\cdot)$ is short for $\partialxioriginal{(\cdot)}$). It is a fundamental measure of complexity of the singularity of $V(f)$ at the origin (in the case that $f$ is the Taylor series of a rational function, or in the case that $\kk = \cc$ and $f$ is a analytic at the origin).

\begin{prop} \label{milnor-prop}
Let $f \in \kk[x_1, \ldots, x_n]$ such that $f(\origin) = 0$.
\begin{enumerate}
\item \label{milnor-0}$\milnor(f) = 0$ if and only the origin is a nonsingular point of $V(f)$.
\item \label{milnor-infty} If $V(f)$ has a non-isolated singularity at the origin, then $\milnor(f) = \infty$. The converse holds if $\character(\kk) = 0$.
\end{enumerate}
\end{prop}

\begin{proof}
The first assertion is clear, so we prove the second assertion. If $V(f)$ has a non-isolated singularity at the origin, then the origin is a non-isolated point of $V(\partialxone{f}, \ldots, \partialxn{f})$, so that $\milnor(f) = \infty$ (\cref{int-mult-curve}). Now assume $\milnor(f) = \infty$. Then the origin is a non-isolated point of $V(\partialxone{f}, \ldots, \partialxn{f})$ (\cref{int-mult-curve}) and therefore there is an irreducible curve $C$ containing the origin such that $C \subseteq V(\partialxone{f}, \ldots, \partialxn{f})$ (\cref{closure-curve-lemma}). It suffices to show that $f|_C \equiv 0$ if $\character(\kk) = 0$. Pick a nonsingular point $a \in C$. Then there is an isomorphism $\phi: \kk[[t]] \cong \hatlocal{C}{a}$, and $d(f \circ \phi)/dt = \sum_{j=1}^n (\partialxj{f} ) (d(\phi_j)/dt) \equiv 0 \in \hatlocal{C}{a}$. Since $\character(\kk) = 0$, it follows that $f$ is constant on $C$. Since $\origin \in C$ and $f(\origin) = 0$, it follows that $C \subseteq V(f)$, as required.
\end{proof}


\begin{example}
The converse to assertion \eqref{milnor-infty} of \cref{milnor-prop} may not be true if $p := \character(\kk)$ is positive; consider e.g.\ the case that $n = 1$ and $f(x) = x^p$, or $n = 2$ and $f(x,y) = x^p + y^q$, where $q \geq 2$ is relatively prime to $p$. The latter example in particular shows that Milnor number can be infinite even for {\em isolated} singular points.
\end{example}

In the case that $\kk= \cc$ and the origin is an isolated singular point of $V(f)$, Milnor originally defined $\milnor(f)$ in \cite[Chapter 7]{milnor} in the following way: let $S^{2n-1}_\epsilon$ be the sphere of radius $\epsilon$ centered at the origin of $\cc^n \cong \rr^{2n}$ and $S^{2n-1} := S^{2n-1}_1$ be the unit sphere of $\cc^n$. Given a morphism $g: \cc^n \to \cc^n$ such that the origin is an isolated zero of $g^{-1}(\origin)$, the {\em multiplicity} of $g$ at the origin is the {\em degree of the mapping}\footnote{The degree of a differentiable map $\phi: M \to N$ between oriented differentiable manifolds of the same dimension, where $M$ is compact and $N$ is connected, is the sum of {\em sign} of $df_x$ over all $x \in \phi^{-1}(y)$ for a generic $y \in N$, where $df_x$ is the derivative map from the tangent space of $M$ at $x$ to the tangent space of $N$ at $y$, and the sign of $df_x$ is either $1$ or $-1$ depending on whether $df_x$ preserves or reverses orientation.} $S^{2n-1}_\epsilon \mapsto S^{2n-1}$ given by $z \mapsto g(z)/\norm{g(z)}$ (where $\norm{\cdot}$ is the Euclidean distance). Milnor defined $\milnor(f)$ as the multiplicity at the origin of the map $z \mapsto (\partialxone{f}, \ldots, \partialxn{f})$. Milnor showed that for all sufficiently small $\epsilon$, if $\phi: S^{2n-1}_\epsilon \setminus V(f) \to S^1$ is the map given by $z \mapsto f(z)/\norm{f(z)}$, then each fiber of $\phi$ is a smooth $(2n-2)$-dimensional real manifold with homotopy type of a ``bouquet'' $S^{n-1} \vee \cdots \vee S^{n-1}$ of spheres, and $\milnor(f)$ is precisely the number of spheres in the bouquet. The fact that the multiplicity of a map $g: \cc^n \to \cc^n$ at the origin equals $\multgzero$ was left in \cite[Appendix B]{milnor} as an exercise; a proof can be found in \cite[Chapter I.5]{sing1}.

\section{Generic Milnor number}
Let $\scrA$ be a (possibly infinite) subset of $\znzero$. We write $\scrL_0(\scrA)$ for the set of all power series in $(x_1, \ldots, x_n)$ supported at $\scrA$. For each $j = 1, \ldots, n$, define
\begin{align}
\partialxj{\scrA} := \{\alpha - e_j :\alpha \in \scrA,\ \alpha - e_j \in \znzero,\ p\ \text{does not divide}\ \alpha_j\}
\label{partial-A}
\end{align}
where $e_j$ is the $j$-th standard unit vector in $\zz^n$ and $p := \character(\kk)$. Note that $\partialxj{\scrA}$ is the support of $\partialxj{g}$ for {\em generic}\footnote{``Generic'' refers to elements of a nonempty Zariski open (dense) subset of $\scrL_0(\scrA)$ in the Zariski topology mentioned in \cref{ind-remark}.} $g \in \scrL_0(\scrA)$. Define
\begin{align*}
\milnor(\scrA)
	&:= \min\{\milnor(f): f\in \scrL_0(\scrA)\} 
\end{align*}
In \cref{milnor-thm} below we estimate $\milnor(\scrA)$ in terms of the intersection multiplicity $\multGammazero$ at the origin of the {\em Newton diagrams} $\Gamma_j$ of $\partialxj{\scrA}$. Given $f \in \scrL_0(\scrA)$, we say that $f$ is \index{Non-degeneracy!partial}{\em partially $\scrA$-non-degenerate at the origin} if the partial derivatives of $f$ are $(\partialxone{\scrA}, \ldots, \partialxn{\scrA})$-non-degenerate at the origin in the sense of \cref{defn:non-degeneracy-0}, i.e.\ if the following property holds:
\begin{align}
\parbox{0.84\textwidth}{
for each nonempty subset $I$ of $[n]$ and each weighted order $\nu$ centered at the origin on $\kk[x_i: i \in I]$, there is no common zero of $\In_{\Ai_j, \nu}((\partialxj{f})|_{\Ki})$, $j = 1, \ldots, n$, on $\nktorus$
} \label{partial-A-non-degeneracy}
\end{align}
where $\Ai_j := \partialxj{\scrA} \cap \ri$; also recall that $\In_{\scrS, \eta}(g)$, where $\eta$ is a weighted degree, $\scrS \subseteq \rr^n$ and $g$ is a Laurent polynomial supported at $\scrS$, is defined as follows:
\begin{align*}
\In_{\scrS, \eta}(g)
	&:= \sum_{\alpha \in \In_\eta(\scrS)} c_\alpha x^\alpha
    =
	\begin{cases}
	\In_\eta(g) & \text{if}\ \supp(g) \cap \In_\eta(\scrS) \neq \emptyset,\\
	0 & \text{otherwise.}
	\end{cases}
\end{align*}
If $f$ is partially $\scrA$-non-degenerate at the origin for $\scrA = \supp(f)$, we simply say that $f$ is {\em partially non-degenerate at the origin}. In other words, $f$ is partially non-degenerate at the origin if it satisfies the following property:
\begin{align}
\parbox{0.84\textwidth}{
for each nonempty subset $I$ of $[n]$ and each weighted order $\nu$ centered at the origin on $\kk[x_i: i \in I]$, there is no common zero of $\In_\nu((\partialxj{f})|_{\Ki})$, $j = 1, \ldots, n$, on $\nktorus$.
} \label{partial-non-degeneracy}
\end{align}
Recall that one does {\em not} have to check this condition for {\em all} nonempty subsets of $[n]$ - see \cref{non-degeneracy-0'-thm,non-degeneracy-0'-remark}.
%


\begin{thm}[{\cite{toricstein}}] \label{milnor-thm}
Let $\Gamma_j := \nd(\partialxj{\scrA})$, $j = 1, \ldots, n$. Assume $\origin \not\in \scrA$. Then
\begin{align}
\milnor(\scrA)  \geq \multGammazero \label{milnor-bound}
\end{align}
Moreover,
\begin{enumerate}
\item \label{milnor-eq}The following are equivalent for all $f \in \scrL_0(\scrA)$:
\begin{enumerate}
\item $\milnor(f) = \multGammazero < \infty$,
\item $f$ is partially $\scrA$-non-degenerate at the origin.
\end{enumerate}
\item\label{milnor-open}  Let $\scrM'_0(\scrA)$ be the set of all $f \in \scrL_0(\scrA)$ which are partially $\scrA$-non-degenerate at the origin. Let $\scrA'$ be any {\em finite} subset of $\scrA$ such that $\partialxj{\scrA'} \supseteq \Gamma_j \cap \partialxj{\scrA}$, $j = 1, \ldots, n$. Then $\scrM'_0(\scrA')$ is a Zariski open subset of $\scrL_0(\scrA')$, and $\scrM'_0(\scrA) = \pi^{-1}(\scrM'_0(\scrA'))$, where $\pi: \scrL_0(\scrA) \to \scrL_0(\scrA')$ is the natural projection.
\item \label{milnor-A-eq} If $\multGammazero = \infty$, then $\scrM'_0(\scrA) = \emptyset$. If $\multGammazero < \infty$, then $\milnor(\scrA) = \multGammazero$ if and only if $\scrM'_0(\scrA)$ is nonempty.
\item \label{milnor-existence} If $\character(\kk) = 0$ and $\multGammazero < \infty$, then $\scrM'_0(\scrA)$ is nonempty and $\milnor(\scrA) = \multGammazero$.
\end{enumerate}
\end{thm}

\begin{example} \label{milnor-nonexample}
Assertion \eqref{milnor-existence} of \cref{milnor-thm} may not be true, i.e.\ the bound in \eqref{milnor-bound} may be strict, in the case that $p := \character(\kk) > 0$. E.g.\ let $\scrA := \{(p+1,1), (1, p+1)\} \subset \zzeroo{2}$. Then $\partialxone{\scrA} = \{(p,1), (0, p+1)\}$ and $\partialxii{\scrA}{2} = \{(p+1, 0), (1, p)\}$. If $f_1 = a_{1,1}x_1^p x_2 + a_{1,2}x_2^{p+1}$ and $f_2 = a_{2,1}x_1^{p+1} + a_{2,2}x_1x_2^p$ are generic polynomials supported respectively at $\partialxone{\scrA}$ and $\partialxtwo{\scrA}$, then $\multpzeronodots{f_1,f_2} = (p+1)^2$. Therefore $\multpzeronodots{\Gamma_1, \Gamma_2} = (p+1)^2$. On the other hand, if $f = ax_1^{p+1}x_2 + bx_1x_2^{p+1}$ is a generic polynomial supported at $\scrA$, then $\partialxone{f} = x_2(ax_1^p + bx_2^p)$ and $\partialxtwo{f} = x_1(ax_1^p + bx_2^p)$ so that $\milnor(\scrA) = \multpzeronodots{x_2(ax_1^p + bx_2^p), x_1(ax_1^p + bx_2^p)} = \infty > \multpzeronodots{\Gamma_1, \Gamma_2}$. It is straightforward to check that $\partialxone{f}, \partialxtwo{f}$ are $\scrA$-degenerate at the origin, i.e.\ $\scrM'_0(\scrA) = \emptyset$.
\end{example}

\begin{proof}[Proof of \cref{milnor-thm}]
Assertions \eqref{milnor-eq}, \eqref{milnor-open} and \eqref{milnor-A-eq} follow from \cref{non-degeneracy-0-thm,finite-0-cor}. Therefore it suffices to show that $\scrM'_0(\scrA)$ is nonempty in the case that $\character(\kk) = 0$. We may assume \woutlog\ that $0 < \milnor(\scrA) < \infty$. Pick any finite subset $\scrA'$ of $\scrA$ satisfying the assumptions of assertion \eqref{milnor-open}. Let $I \subseteq [n]$ and $\nu$ be a weighted order on $\kk[x_i: i \in I]$ centered at the origin. Denote by $\scrL^I_\nu$ the set of all $g \in \scrL_0(\scrA')$ such that
\begin{itemize}
\item $\supp(g) = \scrA'$, and
\item $V(\In_\nu((\partialxone{g})|_{\Ki}), \ldots, \In_\nu((\partialxn{g})|_{\Ki})) \cap \nktorus= \emptyset$
\end{itemize}
It suffices to show that $\scrL^I_\nu$ contains a nonempty Zariski open subset of $\scrL_0(\scrA')$. We may assume without loss of generality that $I = \{1, \ldots, k\}$. Take $g \in \scrL_0(\scrA')$, and express it as
\begin{align*}
g &= g_0(x_1, \ldots, x_k) + \sum_{i = k+1}^n x_ig_i(x_1, \ldots, x_k)  + \cdots
\end{align*}
where the omitted terms have quadratic or higher order in $(x_{k+1}, \ldots, x_n)$. Since $0  < \multGammazero < \infty$, \cref{mult-support} implies that $(\partialxj{g})|_{\Ki} \not\equiv 0$ for $l \geq k$ values of $j$; denote them by $1 \leq j_1 < \cdots < j_l \leq n$. If $g_0 \equiv 0$, then each $j_i > k$, and $V(\In_\nu((\partialxone{g})|_{\Ki}), \ldots, \In_\nu((\partialxn{g})|_{\Ki})) = V(\In_\nu(g_{j_1}), \ldots, \In_\nu(g_{j_l})$. Therefore $g \in \scrL^I_\nu$ if $g_{j_1}, \ldots, g_{j_l}$ are {\em BKK non-degenerate}. Since $l \geq k$, \cref{bkk-0} then implies that $\scrL^I_\nu$ contains a nonempty Zariski open subset of $\scrL_0(\scrA')$, as required. So assume $g_0 \not\equiv 0$. Let $h_0 := \In_\nu(g_0)$ and $V'(h_0) :=  V(\partialxone{h_0}, \ldots, \partialxk{h_0}) \cap \nktorus$. Since $V(\In_\nu((\partialxone{g})|_{\Ki}), \ldots, \In_\nu((\partialxn{g})|_{\Ki})) \subseteq V'(h_0)$, it suffices to show that the set of all polynomials $h \in \scrL_0(\scrA'_0)$, where $\scrA'_0 := \supp(h_0)$, such that $V'(h) = \emptyset$ contains a nonempty Zariski open subset of $\scrL_0(\scrA'_0)$. Let $Z := \{(x, h): x \in \nktorus,\ h \in \scrL_0(\scrA'_0),\ \partialxone{h}(x) = \cdots = \partialxk{h}(x) = 0\} \subset \nktorus \times \scrL_0(\scrA'_0)$, and let $\pi_1: Z \to \nktorus$ and $\pi_2:  Z \to \scrL_0(\scrA'_0)$ be the natural projections. It suffices to show that $\dim(\pi_2(Z)) < \dim \scrL_0(\scrA'_0)$. We prove this by a dimension count. Denote the elements in $\scrA'_0$ by $\alpha_i = (\alpha_{i,1},  \ldots, \alpha_{i,n})$, $i = 1, \ldots, N$, and the coefficients of $x^{\alpha_i}$ (in a polynomial supported at $\scrA'_0$) by $a_i$. Let $x \in \nktorus$. Then $\pi_1^{-1}(x)$ is the subspace of $\scrL_0(\scrA'_0)$ defined by a system of linear equations of the form
\begin{align}
\begin{pmatrix}
\alpha_{1,1}x^{\alpha_1- e_1} & \cdots  & \alpha_{N,1}x^{\alpha_N - e_1} \\
\vdots & & \vdots \\
\alpha_{1,k}x^{\alpha_1- e_k} & \cdots  & \alpha_{N,k}x^{\alpha_N - e_k}
\end{pmatrix}
\begin{pmatrix}
a_1 \\
\vdots\\
a_N
\end{pmatrix}
&= 0 \label{Z_0-eqn-0}
\end{align}
where $e_1, \ldots, e_k$ are the standard unit vectors in $\zz^k$. Since $\bm{\character(\kk) = 0}$ (note: this is the only place the assumption of zero characteristic is used), the rank (as a matrix over $\kk$) of the left-most matrix in \eqref{Z_0-eqn-0} is the same as the rank (as a matrix over $\qq$) of
\begin{align*}
B &:=
\begin{pmatrix}
\alpha_{1,1} & \cdots  & \alpha_{N,1}\\
\vdots & & \vdots \\
\alpha_{1,k} & \cdots  & \alpha_{N,k}
\end{pmatrix}
\end{align*}
Since $\nu$ has positive weights, and the $\alpha_j$ belong to a level set (corresponding to a positive value) of $\nu$, it follows that $\rank(B) = 1 + \dim( \np(h_0))$ and therefore $\dim(\pi_1^{-1}(x)) = N - 1 - \dim( \np(h_0))$. Since this is independent of $x$, it follows that $\dim(Z) = N+n - 1 - \dim(\np(h_0))$. On the other hand, if $x \in V'(h)$ for some $h \in \scrL_0(\scrA'_0)$, then since the support of each of $\partialxj{h}$ is contained in a translation of $\np(h_0)$, it is straightforward to check that $\partialxj{h}(xz^\beta) = 0$ for each $j = 1, \ldots, k$, and each $z \in \nktorus$ and each $\beta \in \zz^n$ which is normal (with respect to the ``dot product'') to $\np(h_0)$. Therefore the dimension of each fiber of $\pi_2$ is at least $n - \dim(\np(h_0))$. It follows that $\dim(\pi_2( Z)) \leq \dim( Z) - n + \dim(\np(h_0)) = N-1$, as required.
\end{proof}

\begin{cor} \label{partial-isolated}
Let $\Gamma_j := \nd(\partialxj{f})$, $j = 1, \ldots, n$. Assume $f(\origin) = 0$ and $f$ is partially non-degenerate. Then the following are equivalent:
\begin{enumerate}
\item $\milnor(f) < \infty$.
\item $\multGammazero < \infty$.
\item $|\{j : \Gamma_j \cap \ri \neq \emptyset\}| \geq |I|$ for each $I \subseteq [n]$.
\end{enumerate}
\end{cor}

\begin{proof}
Combine \cref{mult-support,milnor-thm}.
\end{proof}

\section{Classical notions of non-degeneracy} \label{classical-section}

\subsection{Newton non-degeneracy}
For a power series $f$ in $(x_1, \ldots, x_n)$, we write $\jjj(f)$ for the ideal generated by the partial derivatives of $f$. We say that $f$ is \index{Non-degeneracy!Newton}\index{Newton!non-degeneracy}{\em Newton non-degenerate} iff $\jjj(\In_\nu(f))$ has no zero on $\nktorus$ for each weighted order $\nu$ on $\kk[[x_1, \ldots, x_n]]$ centered at the origin. Newton non-degeneracy is possibly the most studied non-degeneracy property of hypersurface germs: it is a Zariski open condition (in the same sense as partial non-degeneracy at the origin) and, in the characteristic zero case, also nonempty. In this section we discuss its relationship with partial non-degeneracy at the origin. In general Newton non-degeneracy does {\em not} imply ``finite determinacy,'' i.e.\ $f$ can be Newton non-degenerate but still $V(f)$ may have a non-isolated singularity at the origin (take e.g.\ $f:= x_1 \cdots x_n$, $n \geq 2$). However, if the Newton diagram of $f$ is {\em convenient} then Newton non-degeneracy of $f$ implies that the origin is an isolated singularity of $V(f)$ (\cref{Newton-isolated}), which can be resolved by a ``toric modification.'' As a result, the invariants of the singularity can be computed combinatorially in terms of the diagram (see e.g.\ \cite{oka-non-degenerate}). We show in \cref{Newton-to-Milnor} that for isolated singularities, Newton non-degeneracy is a special case of partial non-degeneracy at the origin. However, the following example shows that even in the case of convenient diagrams partial non-degeneracy at the origin does not imply Newton non-degeneracy.

\begin{example}\label{Milnor-not-Newton}
Let $f := x_1 + (x_2 + x_3)^q$, where $q \geq 2$. Then $\nd(f)$ is convenient and $f$ is not Newton non-degenerate (take $\nu$ with weights $(q+1, 1, 1)$ for $(x,y,z)$). However, $f$ is partially non-degenerate at the origin with $\milnor(f) =0$.
\end{example}

We now show that Newton non-degeneracy implies partial non-degeneracy at the origin. The following notation is used in its proof: let $I\subseteq I' \subseteq [n]$, $\nu \in \rnnstar{I}$ and $\nu' \in \rnnstar{I'}$. We say that $\nu$ and $\nu'$ are \index{Compatible weighted order}{\em compatible}
if the weighted order on $\kk[x_i: i \in I]$ induced by $\nu$ and the weighted order on $\kk[x_{i'}: i' \in I']$ induced by $\nu'$ are compatible in the sense of \cref{efficient-section}.

\begin{prop}[{\cite{toricstein}}] \label{Newton-to-Milnor}
Let $f \in \kk[[x_1, \ldots, x_n]]$ and $\Gamma_j := \nd(\partialxj{f})$, $j = 1, \ldots, n$. If $f$ is Newton non-degenerate and $\multGammazero < \infty$, then $f$ is partially non-degenerate at the origin. In particular, if $f$ is Newton non-degenerate and $\nd(f)$ is convenient, then $f$ is partially non-degenerate at the origin.
\end{prop}

\begin{proof}
We start with a direct proof of the second assertion since it is easier to see. Assume $\Gamma := \nd(f)$ is convenient and $f$ is Newton non-degenerate. Pick a nonempty subset $I$ of $[n]$ and a weighted order $\nu$ on $\kk[x_i: i \in I]$ which is centered at the origin. Since $\Gamma$ is convenient, $\Gamma \cap \ri \neq \emptyset$. Therefore we can find a weighted order $\nu'$ on $\kk[x_1, \ldots, x_n]$ such that $\nu'$ is compatible with $\nu$ and $\In_{\nu'}(\Gamma) \subset \ri$. Then $\In_{\nu'}(f)$ depends only on $(x_i: i \in I)$. Since $f$ is Newton non-degenerate, it follows that $\partialxi{(\In_{\nu'}(f))}$, $i \in I$, do not have any common zero in $\nktorus$. But if $i \in I$ is such that $\partialxi{(\In_{\nu'}(f))}$ is not identically zero, then $\partialxi{(\In_{\nu'}(f))} = \partialxi{(\In_{\nu}(f|_{\Ki}))} =  \In_{\nu}((\partialxi{f})|_{\Ki})$. This implies that  $\In_{\nu}((\partialxi{f})|_{\Ki})$, $i \in I$, do not have any common zero on $\Kstari$, as required for partial non-degeneracy of $f$ at the origin. \\

Now we prove the first assertion. If $\multGammazero = 0$, then \cref{mult-support} implies that $f$ is partially non-degenerate at the origin. So assume $0 < \multGammazero < \infty$ and $f$ is partially {\em degenerate} at the origin. It suffices to show that $f$ is Newton degenerate. Pick $I \subseteq [n]$ and a {\em primitive} weighted order $\nu$ centered at the origin on $\kk[x_i: i \in I]$ such that $\In_\nu((\partialxi{f})|_{\Ki})$, $i \in I$, have a common zero $(a_1, \ldots, a_n) \in \nktorus$. At first consider the case that $f|_{\Ki} \not\equiv 0$. Then as in the convenient case pick a weighted order $\nu'$ on $\kk[x_1, \ldots, x_n]$ such that $\nu'$ is compatible with $\nu$ and $\In_{\nu'}(f) \in \kk[x_i: i \in I]$. Then for each $j = 1, \ldots, n$, $\partialxj{(\In_{\nu'}(f))}$ is either $\In_\nu((\partialxj{f})|_{\Ki})$ or is identically zero. It follows that $(a_1, \ldots, a_n)$ is a common zero of $\jjj(\In_{\nu'}(f))$ on $\nktorus$, so that $f$ is Newton degenerate, as required. Now assume that $f|_{\Ki} \equiv 0$. We may assume that $I = \{1, \ldots, k\}$ for some $k$, $1 \leq k \leq n$, and $(\partialxj{f})|_{\Ki} \not\equiv 0$ if and only if $i = k+1, \ldots, k+l$. Then $f$ can be expressed as
\begin{align}
f &= \sum_{j= 1}^l x_{k+j}f_j(x_1, \ldots, x_k) + \sum_{i \geq j \geq 1} x_{k+i}x_{k+j}f_{i,j}(x_1, \ldots, x_n)
\label{fexpression}
\end{align}
Let $h_j := \In_\nu(f_j)$, $j = 1, \ldots, l$, and $B$ be the $k \times l$ matrix with $(i,j)$-th entry $(\partialxi{h_j})(a_1, \ldots, a_k)$. We claim that $\rank(B) < k$. Indeed, let $\nu_i:= \nu(x_i)$, $i = 1, \ldots, k$. For each $j = 1, \ldots, l$, by assumption $(a_1, \ldots, a_k)$ is a common zero of $\In_\nu((\partialxii{f}{k+j})|_{\Ki})= h_j(x_1, \ldots, x_k)$, so that $h_j(a_1t^{\nu_1}, \ldots, a_kt^{\nu_k}) = 0$ for all $t \in \kk$. Note that
\begin{align*}
\frac{d}{dt}(h_j(a_1t^{\nu_1}, \ldots, a_kt^{\nu_k}))
	&= \sum_{i=1}^k \frac{\partial h_j}{\partial x_i}(a_1t^{\nu_1}, \ldots, a_kt^{\nu_k}) \nu_ia_i t^{\nu_i - 1}
	= \sum_{i=1}^k \frac{\partial h_j}{\partial x_i}(a_1, \ldots, a_k) \nu_ia_i t^{\nu(h_j) - 1}
\end{align*}
Setting $t = 1$ it follows that $a'B = 0$, where $a' := (\nu_1a_1, \ldots, \nu_ka_k) \in \kk^n$. Since $\nu$ is primitive, $a' \neq 0$, so that the map $\kk^k \to \kk^l$ given by multiplication by $B$ on the right is {\em not} injective. Therefore $\rank(B) < k$, as claimed. Since $0 < \multGammazero < \infty$, \cref{mult-support} implies that $l \geq k$, so that there is $b = (b_1, \ldots, b_l) \neq 0 \in \kk^l$ such that $B$ times the transpose of $b$ is zero. Let $m$ be the number of nonzero coordinates of $b$. \Woutlog\ we may assume that $b_j \neq 0$ if and only if $j = 1, \ldots, m$. Then we have that
\begin{align}
\sum_{j=1}^m b_j \frac{\partial h_j}{\partial x_i} (a_1, \ldots, a_k) = 0
\end{align}
for all $i = 1, \ldots, k$. Fix positive integers $q, q_{k+m+1}, \ldots, q_n$ and let $\nu'$ be the weighted order on $\kk[x_1, \ldots, x_n]$ such that
\begin{align*}
\nu'(x_i)
	&=
		\begin{cases}
		\nu(x_i) &\text{if}\ i = 1, \ldots, k,\\
		q - \nu(f_{i-k}) &\text{if}\ i = k+1, \ldots, k+m,\\
		q_i &\text{if}\ i = k+m+1, \ldots, n.
		\end{cases}
\end{align*}
If $q \gg 1$ and $q_i \gg q$ for $i = k+m+1, \ldots, n$, then identity \eqref{fexpression} implies that $\In_{\nu'}(f) = \sum_{j=1}^m x_{k+j}h_j$. If $b'_1, \ldots, b'_{n-k-m}$ are arbitrary elements in $\ktorus$, it follows that $(a_1, \ldots, a_k, b_1, \ldots, b_m,  b'_1, \ldots, b'_{n-k-m})$ is a zero of $\jjj(\In_{\nu'}(f))$ on $\nktorus$. Therefore $f$ is Newton degenerate, as required.
\end{proof}

\begin{cor}[{Brzostowski and Oleksik \cite[Theorem 3.1]{brzostowski-oleksik}}] \label{Newton-isolated}
Let $\Gamma_j := \nd(\partialxj{f})$, $j = 1, \ldots, n$. Assume $f(\origin) = 0$ and $f$ is Newton non-degenerate. Then the following are equivalent:
\begin{enumerate}
\item $\milnor(f) < \infty$.
\item $\multGammazero < \infty$.
\item $|\{j : \Gamma_j \cap \ri \neq \emptyset\}| \geq |I|$ for each $I \subseteq [n]$.
\end{enumerate}
\end{cor}

\begin{proof}
Combine \cref{partial-isolated,Newton-to-Milnor}.
\end{proof}

\subsection{Inner Newton non-degeneracy} \label{inner-section}
\index{Diagram}A {\em diagram} $\Gamma$ in $\rr^n$ is the Newton diagram of some subset of $\rzeroo{n}$, and a \index{Face!of a diagram}{\em face} $\Delta$ of $\Gamma$ is a compact face of $\Gamma + \rzeroo{n}$. $\Delta$ is called an \index{Inner!face of a diagram}\index{Face!inner}{\em inner face} if it is not contained in any proper coordinate subspace of $\rr^n$. If $f = \sum_{\alpha \in \znzero} c_\alpha x^\alpha$ is a power series in $(x_1, \ldots, x_n)$ and $\Delta$ is a subset of $\rr^n$, we write $f_\Delta := \sum_{\alpha \in \Delta} c_\alpha x^\alpha$.
We say that $f$ is \index{Non-degeneracy!inner Newton}\index{Inner!Newton non-degeneracy}{\em inner Newton non-degenerate} if there is a {\em convenient} diagram $\Gamma$ such that
\begin{defnlist}
\item no point of $\supp(f)$ ``lies below'' $\Gamma$, i.e.\ $\supp(f) \subseteq \Gamma + \rzeroo{n}$, and
\item for every inner face $\Delta$ of $\Gamma$ and for every nonempty subset $I$ of $[n]$,
\begin{align}
\Delta \cap \ri \neq \emptyset \im V(\jjj(f_\Delta)) \cap \Kstari = \emptyset \label{inner-property}
\end{align}
\end{defnlist}
The difference between Newton non-degeneracy and inner Newton non-degeneracy is most evident in the case of  weighted homogeneous polynomials with isolated singularities, e.g.\ consider the polynomial $f := x_1 + (x_2 + x_3)^q$, where $q \geq 2$, from \cref{Milnor-not-Newton}.  The Newton diagram of $f$ is convenient and has only one inner face, namely the two dimensional face $\Delta$ with vertices $(1, 0,0), (0,q,0), (0,0,q)$. In particular, $f_\Delta = f$ and $\partialxone{f} = 1$, which implies that $f$ is inner non-degenerate. However, as we saw in \cref{Milnor-not-Newton}, $f$ is Newton degenerate. In fact C.\ T.\ C.\ Wall introduced inner Newton non-degeneracy in \cite{wall} in order to find a condition which (in the case of power series with convenient Newton diagrams) is weaker than Newton non-degeneracy, but still wide enough to include all weighted homogeneous polynomials with isolated singularities. We now show that inner Newton non-degeneracy implies partial non-degeneracy at the origin. Given $\Delta \subseteq \rr^n$, we write $I_\Delta$ for the smallest subset of $[n]$ such that $\Delta \subseteq \rii{I_\Delta}$.

\begin{lemma} \label{inner-face-lemma}
Let $\Gamma$ be a convenient diagram, $I \subseteq [n]$ and $\Delta$ be a face of $\Gamma \cap \ri$. Pick an integral element $\nu \in \rnnstar{I}$ centered at the origin such that $\Delta = \In_\nu(\Gamma \cap \ri)$. Then there is an integral element $\nu' \in \rnstar$ such that 
\begin{enumerate}
\item \label{inner-face-lemma:center0} $\nu'$ is centered at the origin,
\item \label{inner-face-lemma:lambda} there is $\lambda > 0$ such that 
\begin{enumerate}
\item $\langle \nu', \alpha \rangle = \lambda \langle \nu, \alpha \rangle$ for all $\alpha \in \rii{I_\Delta}$ (which means $\nu'$ is compatible with $\nu|_{\rii{I_\Delta}}$), 
\item \label{inner-face-lemma:<} $\langle \nu', \alpha \rangle < \lambda \langle \nu, \alpha \rangle$ for all $\alpha \in \rzeroo{I} \setminus \rii{I_\Delta}$ (here $\rzeroo{I}$ denotes the set of all elements in $\ri$ with nonnegative coordinates). 
\end{enumerate}
\item $\Delta' := \In_{\nu'}(\Gamma)$ is an inner face of $\Gamma$. 
\item \label{inner-face-lemma:Delta} $\Delta= \Delta' \cap \rii{I_\Delta}$.
\end{enumerate}
\end{lemma}

\begin{proof}
We proceed by induction on $\delta := n - |I_\Delta|$. The lemma is obviously true if $\delta = 0$. Now assume $\delta = 1$. Then \woutlog\ we may assume $I_\Delta = \{1, \ldots, n-1\}$. Note that $I$ can be either $I_\Delta$ or $[n]$. For each pair of positive integers $q,r$, let $\nu_{q,r}$ be the element on $\rnstar$ defined as follows:
$$\langle \nu_{q,r}, (\alpha_1, \ldots, \alpha_n) \rangle
	:= q\langle \nu, (\alpha_1, \ldots, \alpha_{n-1}, 0) \rangle +r\alpha_n$$
If $r \gg q$, then $\In_{\nu_{q,r}}(\Gamma) = \Delta$. On the other hand, since $\nu$ is centered at the origin and $\Gamma$ is convenient, if $q \gg r$, then $\In_{\nu_{q,r}}(\Gamma)$ is a point on the $x_n$-axis. Therefore we can find $q',r'$ such that $\In_{\nu_{q',r'}}(\Gamma)$ contains $\Delta$ and also a point with positive $n$-th coordinate. We claim that the lemma holds with $\nu' := \nu_{q',r'}$. Indeed, if $I = I_\Delta$, then this is clear by the construction of $\nu_{q',r'}$. If $I = [n]$, then we also need to prove assertion \eqref{inner-face-lemma:<}. But this follows from the observations that if $\alpha'$ is a point on $\In_{\nu_{q',r'}}(\Gamma)$ with positive $n$-th coordinate, then $\langle \nu, \alpha' \rangle > \min_{\Gamma}(\nu)$. 
\\

Now assume $\delta \geq 2$. Pick $\tilde i \not\in I_\Delta$. The inductive hypothesis implies that there is an integral element $\tilde \nu \in \rnnstar{ [n] \setminus \{\tilde i\}}$ such that assertions \eqref{inner-face-lemma:center0} to \eqref{inner-face-lemma:Delta} of the lemma hold with $\Gamma$ replaced by $\Gamma \cap \rii{[n]\setminus \{\tilde i\}}$, $\nu$ replaced by $\nu|_{\rnnstar{I\setminus \{\tilde i\}}}$ and $\nu'$ replaced by $\tilde \nu$. Extend $\tilde \nu$ to an element $\tilde \nu' \in \rnstar$ by defining that for all $(\alpha_1, \ldots, \alpha_n) \in \rr^n$,  
\begin{align*}
\langle \tilde \nu', (\alpha_1, \ldots, \alpha_n) \rangle
	:= \langle \tilde \nu, (\alpha_1, \ldots, \alpha_{\tilde i - 1}, 0, \alpha_{\tilde i + 1}, \ldots, \alpha_n) \rangle + \epsilon \alpha_{\tilde i}
\end{align*}
where $\epsilon$ is a very small positive number. Then we can ensure that 
\begin{prooflist}
\item $\tilde \Delta' := \In_{\tilde \nu'}(\Gamma) = \In_{\tilde \nu}(\Gamma \cap \rii{[n]\setminus \{\tilde i\}})$, 
\item there is $\tilde \lambda > 0$ such that 
\begin{prooflist}
\item $\langle \tilde \nu', \alpha \rangle = \tilde \lambda \langle \nu, \alpha \rangle$ for all $\alpha \in \rii{I_\Delta}$, 
\item $\langle \tilde \nu', \alpha \rangle < \tilde \lambda \langle \nu, \alpha \rangle$ for all $\alpha \in \rzeroo{I} \setminus \rii{I_\Delta}$
\end{prooflist}
\end{prooflist}
Since $I_{\tilde \Delta'} = [n] \setminus \{\tilde i\}$, we can apply the $\delta = 1$ case of the lemma to obtain an integral element $\nu' \in \rnstar$ such that assertions \eqref{inner-face-lemma:center0} to \eqref{inner-face-lemma:Delta} of the lemma hold with $I$ replaced by $[n]$ and $\nu$ replaced by $\tilde \nu'$. It is straightforward to check that the lemma holds with $\nu'$, which completes the proof.
\end{proof}

\begin{prop}[{\cite{toricstein}}]\label{inner-to-Milnor}
If $f \in \kk[[x_1, \ldots, x_n]]$ is inner Newton non-degenerate, then it is partially non-degenerate at the origin.
\end{prop}

\begin{proof}
Assume $f$ is inner Newton non-degenerate with respect to a convenient diagram $\Gamma$. Pick $I \subseteq [n]$ and a weighted order $\nu$ on $\kk[x_i: i \in I]$. Let $\Delta := \In_\nu(\Gamma \cap \ri)$. Pick an integral element $\nu' \in \rnstar$ which satisfies all assertions of \cref{inner-face-lemma}. In particular, $\Delta' := \In_{\nu'}(\Gamma)$ is an inner face of $\Gamma$ and $\Delta = \Delta' \cap \rii{I_\Delta}$. Fix $j$, $1 \leq j \leq n$.

\begin{proclaim} \label{f-Delta'-claim}
One of the following holds:
\begin{prooflist}
\item either $\partialxj{f_{\Delta'}}$ is identically zero on $\Kii{I_\Delta}$,
\item or $(\partialxj{f_{\Delta'}})|_{\Kii{I_\Delta}} = \In_\nu((\partialxj{f})|_{\Ki})$.
\end{prooflist}
\end{proclaim}

\begin{proof}
Assume $\partialxj{f_{\Delta'}}$ is not identically zero on $\Kii{I_\Delta}$. Then it is straightforward to check that
\begin{align}
(\partialxj{f_{\Delta'}})|_{\Kii{I_\Delta}} 
	&= (\In_{\nu'}(\partialxj{f}))|_{\Kii{I_\Delta}} \label{eqn:partial-in=in-partial}
\end{align}
Pick $\alpha \in \supp((\partialxj{f})|_{\Ki}) \setminus \supp((\partialxj{f_{\Delta'}})|_{\Kii{I_\Delta}})$. It suffices to show that $\alpha \not\in \supp(\In_\nu((\partialxj{f})|_{\Ki}))$. Identity \eqref{eqn:partial-in=in-partial} implies that $\alpha \not\in \supp((\In_{\nu'}(\partialxj{f}))|_{\Kii{I_\Delta}})$. If $\alpha \in \rii{I_\Delta}$, then the compatibility of $\nu'$ and $\nu|_{\rii{I_\Delta}}$ implies that $\alpha \not\in \supp(\In_\nu((\partialxj{f})|_{\Ki}))$, as required. So assume $\alpha \in \ri \setminus \rii{I_\Delta}$. Then \cref{inner-face-lemma} implies that 
\begin{align*}
\langle \nu, \alpha \rangle > \langle \nu', \alpha \rangle / \lambda
\end{align*}
where $\lambda$ is as in assertion \eqref{inner-face-lemma:lambda} of \cref{inner-face-lemma}. Since $\langle \nu', \alpha \rangle \geq \nu'(\partialxj{f}) = \nu'((\partialxj{f_{\Delta'}})|_{\Kii{I_\Delta}}) = \lambda \nu((\partialxj{f_{\Delta'}})|_{\Kii{I_\Delta}})$, it follows that $\langle \nu, \alpha \rangle > \nu((\partialxj{f})|_{\Ki})$, as required. 
\end{proof}

\Cref{f-Delta'-claim} implies that $V(\In_\nu((\partialxone{f})|_{\Ki}), \ldots, \In_\nu((\partialxn{f})|_{\Ki})) \cap \Kstari$ is contained in the product of $V((\partialxone{f_{\Delta'}})|_{\Kii{I_\Delta}}, \ldots, (\partialxn{f_{\Delta'}})|_{\Kii{I_\Delta}}) \cap \Kstarii{I_\Delta}$ with $\Kstarii{I\setminus I_\Delta}$. The inner non-degeneracy of $f$ with respect to $\Gamma$ then implies that $V(\In_\nu(\partialxone{f}|_{\Ki}), \ldots, \In_\nu(\partialxn{f}|_{\Ki})) \cap \Kstari = \emptyset$, as required.
\end{proof}

\begin{cor}[{Wall \cite[Lemma 1.2]{wall}}]
If $f$ is inner Newton non-degenerate and $f(0) = 0$, then $\milnor(f) < \infty$. \qed
\end{cor}

If $p := \character(\kk)$ is nonzero, partial Newton non-degeneracy at the origin is strictly weaker than inner Newton non-degeneracy, e.g.\ $x^p + y^p + x^{p+1} + y^{p+1}$ is partially non-degenerate at the origin, but it is inner Newton {\em degenerate}. We do not know if this is true in zero characteristic (see \cref{non-degenerate-relations}).

\section{Newton number: Kushnirenko's formula for the generic Milnor number}
Let $\Gamma$ be a diagram in $\rr^n$. We write $\bar \Gamma$ for the region bounded by the cone with base $\Gamma$ and apex at the origin, and $\volsubk(\Gamma)$, $0 \leq k \leq n$, for the sum of $k$-dimensional Euclidean volumes of the intersections of $\bar \Gamma$ with the $k$-dimensional coordinate subspaces of $\rr^n$ (in particular, $V_0(\Gamma)$ is defined to be $1$). The \index{Newton!number}{\em Newton number} of $\Gamma$ is
\begin{align*}
\newton(\Gamma)
	&:=
	\begin{cases}
		\sum_{k=0}^n (-1)^{n-k} k! \volsubk(\Gamma)
			&\text{if $\Gamma$ is convenient,} \\
		\sup\{\newton( \Gamma \cup \{me_1, \ldots, me_n\}): m \geq 0\}
			& \text{otherwise,}
	\end{cases}
\end{align*}
where $e_j$ are the unit vectors along the (positive direction of the) axes of $\rr^n$. Let $f \in \kk[[x_1, \ldots, x_n]]$ such that $f(\origin) = 0$. A.\ Kushnirenko proved in \cite{kush-poly-milnor} that $\milnor(f) \geq \newton(\nd(f))$, and $\milnor(f) = \newton(\nd(f))$ if $f$ is Newton non-degenerate and if either $\character(\kk) = 0$ or $\nd(f)$ is convenient. C.\ T.\ C.\ Wall proved in \cite{wall} that if $\character(\kk) = 0$, then the equality $\milnor(f) = \newton(\nd(f))$ continues to hold if $f$ is inner Newton non-degenerate. In this section we prove these results and present some generalizations.

\subsection{Preliminary results} \label{kushliminary}
Let $\scrA$ be a subset of $\znzero$ not containing the origin. For each $m \geq 0$, let $\scrA_m := \scrA \cup \{(\alpha_1, \ldots, \alpha_n): \sum_{j=1}^n \alpha_j \geq m\}$. For each $j = 1, \ldots, n$, let
\begin{align*}
\scrA'_{m,j} := \{\alpha - e_j :\alpha \in \scrA_m,\ \alpha - e_j \in \znzero\})
\end{align*}
where $e_j$ are the unit vectors along the (positive direction of the) axes of $\rr^n$. Note that $\partialxj{\scrA_m} \subseteq \scrA'_{m,j}$, and the inclusion is proper if $\character(\kk)$ is positive. In any event, \cref{milnor-thm} and the monotonicity of intersection multiplicity (\cref{monotonic-remark-0}) imply that for each $m \geq 1$,
\begin{align}
\milnor(\scrA)
	&\geq \multzero{\partialxone{\scrA}}{\partialxn{\scrA}}
	\geq \multzero{\partialxone{\scrA_m}}{\partialxn{\scrA_m}}
	\geq \multzero{\scrA'_{m,1}}{\scrA'_{m,n}}
\label{kushnirenko-chain}
\end{align}

\begin{lemma} \label{newton'}
Let $\scrA \subseteq \znzero \setminus \{\origin\}$. Define
\begin{align*}
\scrA'_j := \{\alpha - e_j :\alpha \in \scrA,\ \alpha - e_j \in \znzero\}
\end{align*}
where $e_j$ are the $j$-th standard unit vector in $\rr^n$, $j = 1, \ldots, n$. Assume each $\scrA'_j$ is {\em convenient}. Then $\scrA$ is also convenient and $\multzero{\scrA'_1}{\scrA'_n} = \nu(\nd(\scrA))$.
\end{lemma}

\begin{rem}
Kushnirenko's theorem implies that the assumption ``each $\scrA'_j$ is convenient'' in \cref{newton'} is not needed for the equality $\multzero{\scrA'_1}{\scrA'_n} = \nu(\nd(\scrA))$ to hold - see \cref{newton''} below.
\end{rem}

\begin{proof}[Proof of \cref{newton'}]
Note that $\multzero{\scrA'_1}{\scrA'_n}  = \multgzero$ for all $g_1, \ldots, g_n \in \kk[[x_1, \ldots, x_n]]$ such that 
\begin{prooflist}
\item $\supp(g_j) = \scrA'_j$ for each $j$, and 
\item $g_1, \ldots, g_n$ are non-degenerate at the origin.
\end{prooflist}
We will show that there are $g_1, \ldots, g_n$ which satisfy both of the above properties and in addition satisfy $\multgzero = \nu(\nd(\scrA))$. Indeed, choose $g_1, \ldots, g_n$ which satisfy the above properties, and in addition satisfy 
\begin{prooflist}[resume]
\item \label{newton':restriction} for each pair of subsets $I, J$ of $[n]$ such that $|I| = |J|$, the restrictions $g_j|_{\Ki}$, $j \in J$ are $(\scrA'_j \cap \ri: j \in J)$-non-degenerate at the origin. 
\end{prooflist}
Note that this is possible since each $\scrA'_j$ is convenient and since ``generic'' systems are non-degenerate at the origin (\cref{non-degeneracy-0-thm}). For each $I, J \subseteq [n]$ such that $|I| + |J| = n$, we write
\begin{align*}
\multpnodots{(g_i)_{i \in I}, (x_j)_{j \in J}}{\origin}
	:= \multpnodots{g_{i_1}, \ldots, g_{i_k}, x_{j_1}, \ldots, x_{j_{n-k}}}{\origin}
\end{align*}
where $I = \{i_1, \ldots, i_k\}$ and $J = \{j_1, \ldots, j_{n-k}\}$. Since $\scrA'_j$ are convenient, \cref{mult-support} and property \ref{newton':restriction} imply that $\multpnodots{(g_i)_{i \in I}, (x_j)_{j \in [n]\setminus I}}{\origin}$ is defined for each $I \subseteq [n]$. Consequently \cref{kushnirenko-lemma} implies that
\begin{align*}
\multgzero
	&= \sum_{I \subseteq [n]} (-1)^{n- |I|} \multpnodots{(x_ig_i)_{i \in I}, (x_j)_{j \in [n]\setminus I}}{\origin}
	= \sum_{I \subseteq [n]} (-1)^{n- |I|} \multzero{x_{i_1}g_{i_1}|_{\Ki}} {x_{i_{|I|}}g_{i_{|I|}}|_{\Ki}}
\end{align*}
where $i_1, \ldots, i_{|I|}$ are elements of $I$ for each $I \subseteq [n]$. Fix $I \subseteq [n]$. Since the Newton diagram of the union of the supports of $x_{i_j}g_{i_j}|_{\Ki}$ is $\nd(\scrA) \cap \ri$ (this is where the assumption $\origin \not\in \scrA$ is used!), property \ref{newton':restriction} and \cref{kushniplicity} imply that 
\begin{align*}
\multzero{x_{i_1}g_{i_1}|_{\Ki}} {x_{i_{|I|}}g_{i_{|I|}}|_{\Ki}}
	&= \volsub{|I|}(\nd(\scrA) \cap \ri) 
\end{align*}
The result now follows from the definition of $\newton(\cdot)$. 
\end{proof}

\begin{cor} \label{newton-0}
$\multzero{\scrA'_{m,1}}{\scrA'_{m,n}} = \nu(\nd(\scrA_m))$ for each $m \geq 1$. \qed
\end{cor}

\begin{cor} [{Kushnirenko \cite[Theorem I, part (i)]{kush-poly-milnor}}] \label{milnor-geq-newton}
\index{Kushnirenko's theorem!on Milnor number}
$\milnor(f) \geq \newton(\nd(f))$ for all $f \in \kk[[x_1, \ldots, x_n]]$. In particular, if $\newton(\nd(f)) = \infty$, then $\milnor(f) = \infty$.
\end{cor}

\begin{proof}
Combine inequation \eqref{kushnirenko-chain} and \cref{newton-0}.
\end{proof}
%
%
%
%
%

\subsection{Characteristic zero case} \label{zero-section}
Continue to assume that $\scrA$ is a subset of $\znzero$ not containing the origin.

\begin{thm}\label{milnor=newton-0}
If $\character(\kk) = 0$, then $\milnor(\scrA) = \newton(\nd(\scrA))$.
\end{thm}

\begin{proof}
If $\character(\kk) = 0$, then \cref{milnor-thm} implies that $\milnor(\scrA) = \multzero{\partialxone{\scrA}}{\partialxn{\scrA}}$. At first consider the case that $\multzero{\partialxone{\scrA}}{\partialxn{\scrA}} < \infty$. For $m \gg 1$ then \cref{finite-mult-determinacy} implies that $\multzero{\partialxone{\scrA}}{\partialxn{\scrA}} = \multzero{\partialxone{\scrA_m}}{\partialxn{\scrA_m}}$. Since in zero characteristic $\scrA'_{m,j} = \partialxj{\scrA_m}$, \cref{newton-0} implies that $\milnor(\scrA) = \newton(\nd(\scrA))$. On the other hand, if $\multzero{\partialxone{\scrA}}{\partialxn{\scrA}} = \infty$, then $\sup_m \multzero{\partialxone{\scrA_m}}{\partialxn{\scrA_m}} = \infty$ (\cref{finite-mult-determinacy}), so that \cref{newton-0} implies that $\newton(\nd(\scrA)) = \infty$, as required.
\end{proof}

\begin{cor} \label{newton''}
\Cref{newton'} holds even without the assumption that each $\scrA'_j$ is convenient.
\end{cor}

\begin{proof}
In zero characteristic $\milnor(\scrA) = \multzero{\scrA'_1}{\scrA'_n}$ (\cref{milnor-thm}); now use \cref{milnor=newton-0}.
\end{proof}

\begin{cor}[{Cf.\ \cite[Characteristic zero case of Theorem I]{kush-poly-milnor}, \cite[Theorem 1.6]{wall}, \cite[Corollaries 3.10 and 3.11]{brzostowski-oleksik}}]
\index{Kushnirenko's theorem!on Milnor number}
Let $f \in \kk[[x_1, \ldots, x_n]]$ be such that $f(\origin) = 0$. Assume $\character(\kk) = 0$. Then $\milnor(f) = \newton(\nd(f))$ whenever $f$ is partially non-degenerate at the origin. In particular, if $f$ is either Newton non-degenerate or inner Newton non-degenerate, then $\milnor(f) = \newton(\nd(f))$. 
\end{cor}

\begin{proof}
Combine \cref{milnor-thm,Newton-to-Milnor,inner-to-Milnor,milnor=newton-0}.
\end{proof}

\subsection{The general case} \label{general-section}
We continue to use the notation of \cref{kushliminary}. In particular, $\scrA$ is a subset of $\znzero$ not containing the origin. Let $\scrM'_0(\scrA)$ be as in \cref{milnor-thm} the set of all $f \in \scrL_0(\scrA)$ which are partially $\scrA$-non-degenerate at the origin. \Cref{finite-mult-determinacy,milnor-thm}, inequation \eqref{kushnirenko-chain} and \cref{newton-0} imply the following result.

\begin{prop} \label{milnor=newton}
 $\milnor(\scrA) = \newton(\nd(\scrA))$ if and only if
\begin{enumerate}
\item either $\newton(\nd(\scrA)) = \infty$, or
\item
\begin{enumerate}
\item \label{partial-existence} $\scrM'_0(\scrA)$ is nonempty, and
\item \label{mult-equality}
$\multzero{\partialxone{\scrA_m}}{\partialxn{\scrA_m}} = \newton(\nd(\scrA_m))$ for all $m \gg 1$ not divisible by $p$. \qed
\end{enumerate}
\end{enumerate}
\end{prop}

\begin{rem}
\Cref{mult-support,newton''} imply that the following are equivalent:
\begin{enumerate}
\item $\newton(\nd(\scrA)) = \infty$.
\item $\multzero{\scrA'_1}{\scrA'_n} = \infty$.
\item there is a nonempty subset $I$ of $[n]$ such that $|\{j : \scrA'_j \cap \ri \neq \emptyset\}| < |I|$.
\end{enumerate}
\end{rem}

Let $p := \character(\kk)$. \Cref{milnor-nonexample} shows that if $p$ is positive, then condition \eqref{partial-existence} of \cref{milnor=newton} is nontrivial. The example below shows that condition \eqref{mult-equality} is nontrivial as well.

\begin{figure}[h]
\begin{center}
\begin{tikzpicture}[ scale=0.6]
\tikzstyle{dot} = [red, circle, minimum size=3pt, inner sep = 0pt, fill]
\def\p{2}
\pgfmathsetmacro\r{2*\p + 4};
\pgfmathsetmacro\q{2*\p + 2};
\draw [<->] (0,{\q + 0.5}) |- ({\r + 0.5},0);

\node[dot, label=below:{\small $r$}] (r) at (\r,0) {};
\node[dot, label=right:{\small $(p,p)$}] (p) at (\p,\p) {};
\node[dot, label=left:{\small $q$}] (q) at (0,\q) {};

\draw[thick, green] (r) -- (p) -- (q);
\draw[dashed] (0,0) -- (p);
\draw (\p, 1) node [right] {\small $pr/2$};
\draw (1, \p) node {\small $pq/2$};
\end{tikzpicture}
\caption{Subdivision of the area under the Newton diagram from \cref{non-newton-2}}  \label{fig:non-newton-2}
\end{center}
\end{figure}

\begin{example} \label{non-newton-2}
Assume $p$ is positive. Let $f = x_1^r + x_1^{p}x_2^p + x_2^{q}$, where $q,r$ are large positive integers not divisible by $p$. Let $\scrA := \supp(f)$ and $\Gamma := \nd(f)$. It is straightforward to see from \cref{fig:non-newton-2} that $\newton(\Gamma) =  2p(q+r)/2- (r + q) + 1 = (p-1)(q+r) + 1$. On the other hand, $\partialxone{f} = rx_1^r$ and $\partialxtwo{f} = qx_2^{q-1}$, so that $f$ is partially non-degenerate at the origin. It follows that $\scrM'_0(\scrA) \neq \emptyset$ and $\milnor(\scrA) = \milnor(f) = \multpzeronodots{rx_1^{r-1},  qx_2^{q-1}} = (r-1)(q-1)$. It follows that $\milnor(\scrA) > \newton(\nd(\scrA))$ for sufficiently large $q,r$.
\end{example}


We now state a condition which guarantees that condition \eqref{mult-equality} of \cref{milnor=newton} holds. Given a subset $\scrB$ of $\znzero$, let $J_\scrB := \{j \in [n]: \partialxj{\scrB} \neq \emptyset\} = \{j \in [n]:$ there is $(\beta_1, \ldots, \beta_n) \in \scrB$ such that $p$ does not divide $\beta_j\}$. The condition is the following:
\begin{align}
\parbox{.78\textwidth}{
for each $m \gg 1$ not divisible by $p$ and for each face $\Delta$ of $\nd(\scrA_m)$, $J_{\scrA_m \cap \Delta} \neq \emptyset$, and the convex hulls of $\partialxj{(\scrA_m \cap \Delta)}$, $j \in J_{\scrA_m \cap \Delta}$, are dependent.
} \label{milnor=newton-condition-1}
\end{align}

\begin{prop}  \label{milnor=newton-prop-1}
Assume \eqref{milnor=newton-condition-1} holds. Then condition \eqref{mult-equality} of \cref{milnor=newton} also holds. In particular, $\milnor(\scrA) = \newton(\nd(\scrA))$ if either $\newton(\nd(\scrA)) = \infty$ or $\scrM'_0(\scrA)$ is nonempty.
\end{prop}

\begin{proof}
Pick $m \gg 1$ not divisible by $p$. The arguments from the proof of \cref{newton'} show that it suffices to prove that $\multzero{x_1g_1}{x_ng_n} = n! \volsubn(\nd(\scrA_m))$ for power series $g_j$ such that $\supp(g_j) = \partialxj{\scrA_m}$ and $g_1, \ldots, g_n$ are non-degenerate at the origin. Condition \eqref{milnor=newton-condition-1} ensures that the Newton diagram of the union of the supports of $x_jg_j$ is $\nd(\scrA_m)$, and that the condition of assertion \eqref{kushniplicity-equal} of \cref{kushniplicity} is also satisfied. Therefore the result follows from \cref{kushniplicity}.
\end{proof}

Note that condition \eqref{milnor=newton-condition-1} is not necessary for \eqref{mult-equality} of \cref{milnor=newton} to hold - see \cref{newton-problem}.

\begin{cor}[{\cite[Positive characteristic case of Theorem I]{kush-poly-milnor}}] \label{kushnirenko+}
\index{Kushnirenko's theorem!on Milnor number}
Let $f \in \kk[[x_1, \ldots, x_n]]$ be such that $f(\origin) = 0$, $\nd(f)$ is convenient, and $f$ is Newton non-degenerate. Then $\milnor(f) = \newton(\nd(f))$.
\end{cor}

\begin{proof}
Due to \cref{Newton-to-Milnor,milnor=newton-prop-1} it suffices to show that $\scrA:= \supp(f)$ satisfies condition \eqref{milnor=newton-condition-1}. Since $\scrA$ is convenient, $\nd(\scrA_m) = \nd(\scrA)$ for $m \gg 1$. So pick a face $\Delta$ of $\nd(\scrA)$. As in \cref{inner-section} let $f_\Delta$ be the ``portion'' of $f$ supported at $\Delta$. Let $g_j :=\partialxj{f_\Delta}$, $j = 1, \ldots, n$. Since $f$ is Newton non-degenerate, it follows that $g_1, \ldots, g_n$ are BKK non-degenerate. Since $g_1, \ldots, g_n$ have no common zero on $\nktorus$, \cref{bkk-bound-thm,bkk-non-degenerate-thm,positively-mixed} imply that the Newton polytopes of $g_j$ are dependent, as required.
\end{proof}

\section{Open problems}

\subsection{Existence of non-degenerate polynomials}
Let $\scrA$ be a subset of $\znzero$ not containing the origin. The estimate of $\milnor(\scrA)$ from \cref{milnor-thm} is exact if and only if the set $\scrM'_0(\scrA)$ of power series which are partially $\scrA$-non-degenerate at the origin is nonempty. \Cref{milnor-thm} shows that in characteristic zero $\scrM'_0(\scrA)$ is always nonempty, and \cref{milnor-nonexample} shows that in positive characteristic there are $\scrA$ such that $\scrM'_0(\scrA)$ is empty. This motivates the following problem.

\begin{problem} \label{milnor-existence-problem}
In the case that $\character(\kk) > 0$, characterize those $\scrA$ for which $\scrM'_0(\scrA)$ is nonempty. Compute $\milnor(\scrA)$ for those $\scrA$ such that $\scrM'_0(\scrA)$ is empty.
\end{problem}

Let $\scrN'^0_0(\scrA)$ be the set of all power series $f$ supported at $\scrA$ such that $\nd(f) = \nd(A)$, $\nd(\partialxj{f}) = \nd(\partialxj{\scrA})$ for each $j$, and $f$ is Newton non-degenerate. \Cref{Newton-to-Milnor} implies that $\scrN'^0_0(\scrA) \subseteq \scrM'_0(\scrA)$ when $\milnor(\scrA) < \infty$.  The following is therefore a subproblem of \cref{milnor-existence-problem} in that case.

\begin{problem} \label{newton-existence-problem}
In the case that $\character(\kk) > 0$, characterize those $\scrA$ for which $\scrN'^0_0(\scrA)$ is nonempty.
\end{problem}

The proof of \cref{milnor-thm} gives a sufficient condition for existence of Newton non-degenerate polynomials: let $\scrB$ be a finite subset of $\znzero$ and $B$ be the $n \times |\scrB|$ matrix whose columns are the elements of $B$. Let $\rank_\kk(\cdot)$ denote the rank of a matrix over $\kk$. The condition we are interested in is the following:
\begin{align}
\rank_\kk(B) = \dim(\conv(\scrB)) + 1 \label{rank-condition-1}
\end{align}

\begin{lemma} \label{newton-existence-lemma}
Assume \eqref{rank-condition-1} holds. Then the set of polynomials $g$ supported at $\scrB$ such that $\partialxj{g}$, $j = 1, \ldots, n$, have no common zero on $\nktorus$ contains a nonempty Zariski open subset of the space of all polynomials supported at $\scrB$.
\end{lemma}

\begin{proof}
This is in fact the main content of the proof of \cref{milnor-thm} (the only place where the zero characteristic played a role in that proof was to ensure that \eqref{rank-condition-1} holds).
\end{proof}

\begin{cor} \label{newton-existence-1}
If \eqref{rank-condition-1} holds with $\scrB = \Delta \cap \scrA$ for each face $\Delta$ of $\nd(\scrA)$, then $\scrN'^0_0(\scrA) \neq \emptyset$. \qed
\end{cor}

Assertion \eqref{milnor-eq} of \cref{milnor-thm} implies that for $\scrA$ to admit polynomials which are partially non-degenerate at the origin, it is necessary that
\begin{align}
\multzero{\partialxone{\scrA}}{\partialxn{\scrA}} < \infty \label{mult-Gamma-finite}
\end{align}
\Cref{mult-support} implies that \eqref{mult-Gamma-finite} is equivalent to the condition that $|\{j : \partial_j{\scrA} \cap \ri \neq \emptyset\}| \geq |I|$ for each $I \subseteq [n]$. The following is an immediate corollary of \cref{Newton-to-Milnor,newton-existence-1}.

\begin{cor} \label{milnor-existence-1}
Assume \eqref{mult-Gamma-finite} holds and \eqref{rank-condition-1} holds with $\scrB = \Delta \cap \scrA$ for each face $\Delta$ of $\nd(\scrA)$, then $\scrM'_0(\scrA) \neq \emptyset$. \qed
\end{cor}

\begin{question}
Is the condition from \cref{newton-existence-1} necessary for $\scrN'^0_0(\scrA)$ to be nonempty?
\end{question}

\subsection{Relation among non-degeneracy conditions} \label{non-degenerate-relations}
Given a power series $f$, let us write $(N)$, $(I)$, $(P)$ to denote respectively the conditions that $f$ is Newton non-degenerate, inner Newton non-degenerate and partially non-degenerate at the origin. If $p := \character(\kk)$ is nonzero, then $(P)$ does not imply $(I)$, e.g. $x^p + y^p + x^{p+1} + y^{p+1}$ is partially non-degenerate at the origin, but it is not inner Newton degenerate. This observation together with the discussion from \cref{classical-section} implies the relations depicted in \cref{fig:relations}, where ``$(N_{\milnor < \infty})$'' denotes the condition that $f$ is Newton non-degenerate and $\milnor(f) < \infty$.

\begin{center}
\begin{figure}[h]

\tikzcdset{arrow style=tikz, diagrams={>={Computer Modern Rightarrow[width=8pt,length=3pt]}}}

\begin{subfigure}[b]{0.3\textwidth}
\begin{tikzcd}[column sep=tiny]
& (P)
\arrow[ld, double, "/" anchor=center, shift left = 0.5ex] \arrow [ld, double, shift right=0.5ex, <-]
\arrow[rd, double,  shift right = 0.5ex, swap, "?"] \arrow [rd, double, shift left= 0.5ex, <-]
& \\
(N_{\milnor < \infty})
\arrow[double, rr, "?" , shift left=.5ex] \arrow[double, rr, "/" anchor=center, <-, shift right=0.5ex]
&
&(I)
\end{tikzcd}
\caption{zero characteristic}
\label{fig:relations-0}
\end{subfigure}
\begin{subfigure}[b]{0.3\textwidth}
\begin{tikzcd}[column sep=tiny]
& (P)
\arrow[ld, double, "/" anchor=center, shift left = 0.5ex] \arrow [ld, double, shift right=0.5ex, <-]
\arrow[rd, double,  shift right = 0.5ex, "/" anchor=center] \arrow [rd, double, shift left= 0.5ex, <-]
& \\
(N_{\milnor < \infty})
\arrow[double, rr, "?" , shift left=.5ex] \arrow[double, rr, "/" anchor=center, <-, shift right=0.5ex]
&
&(I)
\end{tikzcd}
\caption{positive characteristic}
\label{fig:relations-+}
\end{subfigure}
\caption{Relation among non-degeneracy conditions}  \label{fig:relations}
\end{figure}
\end{center}

\begin{problem}
Determine if the question-marked implications from \cref{fig:relations} are valid.
\end{problem}

We now show that in zero characteristic the implication $(P) \im (I)$ does hold in dimension $\leq 3$.

\begin{prop} \label{Milnor-to-inner}
Pick $f \in \kk[x_1, \ldots, x_n]]$ such that $f(0) = 0$. Assume $f$ is partially non-degenerate at the origin. If $n \leq 3$ and $\character(\kk)$ is zero, then $f$ is also inner Newton non-degenerate.
\end{prop}

\begin{proof}
Since all the $\partialxj{f}$ can not be identically zero on any axis, $\nd(f)$ satisfies the following property:
\begin{align}
\parbox{0.6\textwidth}{
the distance from any axis to $\nd(f)$ can not be greater than $1$.
}\label{ndistance}
\end{align}
This leads to the possibilities of \cref{fig-partial-2} in the case that $n = 2$. If $\nd(f)$ is convenient, take $\Gamma = \nd(f)$, otherwise take $\Gamma$ to be the union of $\nd(f)$ and edges (with appropriate slopes) from the end points of $\nd(f)$ to some points on the axes (e.g.\ the ``dashed edges'' in \cref{fig-partial-2}). It is straightforward to check that $f$ is inner Newton non-degenerate with respect to $\Gamma$.

\begin{figure}[h]
\begin{center}

\begin{tikzpicture}[scale=0.3]
\def\shiftone{12}
\def\opazero{0.5}
\def\tx{2}
\def\ty{3}
\def\gridx{7.5}
\def\gridy{5.5}
\def\colorzero{green}
\def\colorone{red}

\draw [gray,  line width=0pt] (-0.5,-0.5) grid (\gridx,\gridy);
\draw [<->] (0, \gridy) |- (\gridx, 0);

\draw[ultra thick, \colorone]  (7,0) -- (4,1) -- (2,2) -- (1,3) -- (0,5);
\fill[\colorzero, opacity=\opazero ]  (7,0) -- (4,1) -- (2,2) -- (1,3) -- (0,5) -- (0,\gridy) -- (\gridx,\gridy) -- (\gridx,0) -- cycle;

\begin{scope}[shift={(\shiftone,0)}]
	\draw [gray,  line width=0pt] (-0.5,-0.5) grid (\gridx,\gridy);
	\draw [<->] (0, \gridy) |- (\gridx, 0);
	\draw[ultra thick, \colorone]  (7,0) -- (4,1) -- (2,2) -- (1,3);
	\fill[\colorzero, opacity=\opazero ]  (7,0) -- (4,1) -- (2,2) -- (1,3) -- (1,\gridy) -- (\gridx,\gridy) -- (\gridx,0) -- cycle;
	\draw[ thick, dashed] (1,3) -- (0,4);
\end{scope}

\begin{scope}[shift={(2*\shiftone,0)}]
	\draw [gray,  line width=0pt] (-0.5,-0.5) grid (\gridx,\gridy);
	\draw [<->] (0, \gridy) |- (\gridx, 0);
	\draw[ultra thick, \colorone]  (4,1) -- (2,2) -- (1,3) -- (0,5);
	\fill[\colorzero, opacity=\opazero ] (4,1) -- (2,2) -- (1,3) -- (0,5)-- (0,\gridy) -- (\gridx,\gridy) -- (\gridx,1) -- cycle;
	\draw[ thick, dashed] (4,1) -- (6,0);
\end{scope}

\begin{scope}[shift={(3*\shiftone,0)}]
	\draw [gray,  line width=0pt] (-0.5,-0.5) grid (\gridx,\gridy);
	\draw [<->] (0, \gridy) |- (\gridx, 0);
	\draw[ultra thick, \colorone]  (4,1) -- (2,2) -- (1,3);
	\fill[\colorzero, opacity=\opazero ]   (4,1) -- (2,2) -- (1,3) -- (1,\gridy) -- (\gridx,\gridy) -- (\gridx,1) -- cycle;
	\draw[ thick, dashed] (1,3) -- (0,4);
	\draw[ thick, dashed] (4,1) -- (6,0);
\end{scope}

\end{tikzpicture}

\end{center}
\caption{Four possibilities for $\nd(f)$ in dimension two} \label{fig-partial-2}
\end{figure}

Now assume that $n = 3$. If $\nd(f)$ does not touch (at least) one of the three coordinate hyperplanes, then $f$ would be divisible by some $x_i$. The partial non-degeneracy of $f$ at the origin would then imply that $f = x_ig$ such that $g(0) \neq 0$, i.e.\ $\nd(f) = \{e_i\}$, where $e_i$ is the $i$-th standard unit vector in $\rr^3$. In that case $f$ would be inner Newton non-degenerate with any convenient diagram with a vertex at $e_i$. Therefore we may assume that $\nd(f)$ touches every coordinate hyperplane. Let $\Gamma$ be the Newton diagram of $f + x_1^N + x_2^N + x_3^N$ for some $N \gg 1$. We claim that $f$ is inner Newton non-degenerate with respect to $\Gamma$. Indeed, let $\Delta$ be an inner face of $\Gamma$ and $I \subseteq \{1, 2, 3\}$ be such that $\Delta \cap \ri \neq \emptyset$. We will check that the inner non-degeneracy condition \eqref{inner-property} holds for $\Delta$ and $I$. Let $\Delta' := \Delta \cap \ri$. If $\Delta' \cap \nd(f) = \emptyset$, then $\Delta' = \{Ne_i\}$ for some $i$. Assume \woutlog\ that $i = 1$. Let $\alpha^{(2)}$ (respectively, $\alpha^{(3)}$) be the point on the intersection of $\nd(f)$ with the $(x_1,x_2)$-plane (respectively, $(x_1, x_3)$-plane) which is closest to $x_1$-axis. Since $\Delta$ is an inner face, it is straightforward to check that when $N$ is sufficiently large, $\Delta$ must be the triangle with vertices $Ne_1, \alpha^{(2)}, \alpha^{(3)}$. Property \eqref{ndistance} implies that either $\alpha^{(2)} = (k, 1,0)$ for some $k \geq 0$ or $\alpha^{(3)} = (k,0,1)$ for some $k \geq 0$. Then $x_1^k \in \jjj(f_\Delta)$, and therefore $V(\jjj(f_\Delta))$ does not contain any point on $x_1$-axis other than possibly the origin, as required. It remains to consider the case that $\Delta'' := \Delta' \cap \nd(f) \neq \emptyset$. If $|I| = 1$ or $|I| = 3$, it is straightforward to check that violation of the inner non-degeneracy condition \eqref{inner-property} implies violation of partial non-degeneracy of $f$. So we may assume \woutlog\ $I = \{1, 2\}$. Note that $\dim(\Delta'') \leq \dim(\Delta') \leq 1$. If $\dim(\Delta'') = 0$, then $f_\Delta = cx_1^{\alpha_1}x_2^{\alpha_2} + x_3g$ for some polynomial $g$, $(\alpha_1, \alpha_2) \in \zzeroo{2} \setminus \{\origin\}$ and $c \neq 0$. It is then clear that either $\partialxone{f_\Delta}$ or $\partialxii{f_\Delta}{2}$ does not vanish at any point on $\Kstari$. So assume $\dim(\Delta'') = 1$. It then follows that $\Delta'' = \Delta'$, and if $N$ is sufficiently large, then $\Delta$ is in fact a (two dimensional) face of $\nd(f)$ containing $\Delta'$, see \cref{fig:3inner-Delta}.

\begin{center}
\begin{figure}[h]
\def\scale{0.6}
\def\viewx{75}
\def\viewy{15}
\def\colorzero{blue}
\def\colorzerob{black}
\def\colorone{orange}
\def\coloroneb{red}
\def\coloreta{green}
\begin{subfigure}[b]{0.3\textwidth}
\begin{tikzpicture}[scale=\scale]
\pgfplotsset{every axis title/.append style={at={(0,-0.2)}}, view={\viewx}{\viewy}, axis lines=middle, enlargelimits={false}}

\begin{axis}[
	xticklabels = \empty, yticklabels = \empty, zticklabels = \empty,
	xtick = \empty, ytick = \empty, ztick = \empty
]
	\addplot3[fill=\colorzero,opacity=\opazero] coordinates{(1,6,0) (6,1,0) (4,0,3) (1,1,5) (0,3,4) (1,6,0)};
	\addplot3 [ultra thick, \colorzerob] coordinates{(1,6,0) (6,1,0)};
	\addplot3 [ultra thick, \coloreta, ->] coordinates{(3,4,0) (4,5,0)};
\end{axis}
\node[anchor = west] at (2.5,2.5) {\picfontsize $\Delta$};
\node[anchor = west] at (2.75,0) {\picfontsize $\Delta'$};
\node[anchor = west] at (5.5,0.5) {\picfontsize $\nu$};
\end{tikzpicture}
\caption{$\Delta$}  \label{fig:3inner-Delta}
\end{subfigure}
\begin{subfigure}[b]{0.3\textwidth}
\begin{tikzpicture}[scale=\scale]
\pgfplotsset{every axis title/.append style={at={(0,-0.2)}}, view={\viewx}{\viewy}, axis lines=middle, enlargelimits={false}}

\begin{axis}[
	xticklabels = \empty, yticklabels = \empty, zticklabels = \empty,
	xtick = \empty, ytick = \empty, ztick = \empty
]
	\addplot3[fill=white,opacity=0] coordinates{(1,6,0) (6,1,0) (4,0,3) (1,1,5) (0,3,4) (1,6,0)};
	\addplot3[fill=\colorzero,opacity=\opazero] coordinates{(2,4,0) (4,2,0) (4,0,2) (1,1,4) (0,3,3) (2,4,0)};
	\addplot3 [ultra thick, \coloreta, ->] coordinates{(3,3,0) (4,4,0)};
\end{axis}
\node[anchor = west] at (1,2.5) {\picfontsize $\partialxii{\Delta}{3}$};
\node[anchor = west] at (4.5,0.5) {\picfontsize $\nu$};
\end{tikzpicture}
\caption{Case 1: $\partialxii{\Delta}{3}$\\ touches $(x_1,x_2)$-plane}  \label{fig:3inner-Delta3-1}
\end{subfigure}
\begin{subfigure}[b]{0.3\textwidth}
\begin{tikzpicture}[scale=\scale]
\pgfplotsset{every axis title/.append style={at={(0,-0.2)}}, view={\viewx}{\viewy}, axis lines=middle, enlargelimits={false}}

\begin{axis}[
	xticklabels = \empty, yticklabels = \empty, zticklabels = \empty,
	xtick = \empty, ytick = \empty, ztick = \empty
]
	\addplot3[fill=white,opacity=0] coordinates{(1,6,0) (6,1,0) (4,0,3) (1,1,5) (0,3,4) (1,6,0)};
	\addplot3[fill=\colorzero,opacity=\opazero] coordinates{(4,0,2) (1,1,4) (0,3,3) (1,3,2) (4,0,2)};
	\addplot3[fill=\colorone,opacity=\opazero] coordinates{(4,0,2) (1,3,2) (2,5,0) (4,3,0) (4,0,2)};
	\addplot3 [ultra thick, \coloroneb] coordinates{(4,0,2) (1,3,2)};
	\addplot3 [ultra thick, \coloroneb] coordinates{(2,5,0) (4,3,0)};
	\addplot3 [ultra thick, \coloreta, ->] coordinates{(3,4,0) (4,5,0)};
\end{axis}
\node[anchor = west] at (1,3.25) {\picfontsize $\partialxii{\Delta}{3}$};
\node[anchor = east] at (1.1,2) {\picfontsize $\Lambda''$};
\node at (2.75,1.75) {\picfontsize $\Lambda$};
\node[anchor = east] at (3.6,0.3) {\picfontsize $\Lambda'$};
\node[anchor = west] at (5.5,0.5) {\picfontsize $\nu$};
\end{tikzpicture}
\caption{Case 2: $\partialxii{\Delta}{3}$ does not touch $(x_1,x_2)$-plane}  \label{fig:3inner-Delta3-2}
\end{subfigure}
\caption{Scenarios when $\dim(\Delta'') = 1$} \label{fig:3inner}
\end{figure}
\end{center}

If $\partialxii{\Delta}{3}$ touches the $(x_1,x_2)$-plane, then applying the partial non-degeneracy condition \eqref{partial-non-degeneracy} with $I = \{1,2\}$ and $\nu$ equal to the ``inner normal'' to $\Delta'$ in $\ri$ shows that $V(\jjj(f_\Delta)) \cap \Kstari = \emptyset$. So assume that $\partialxii{\Delta}{3}$ does not touch $(x_1,x_2)$-plane. In this case we claim that $V(\partialxone{f_{\Delta'}},\partialxii{f_{\Delta'}}{2}) \cap \Kstari = \emptyset$. Indeed, assume to the contrary that there is a common zero $a = (a_1,a_2,0)$ of $\partialxone{f_{\Delta'}},\partialxii{f_{\Delta'}}{2}$ on $\Kstari$. The partial non-degeneracy then implies that the intersection of $\partialxii{\Gamma}{3}$ with the $(x_1,x_2)$-plane is nonempty, and if $\Lambda'$ is the face of this intersection determined by $\nu$, then $(\partialxii{f}{3})_{\Lambda'}$ does {\em not} vanish at $a$. It is straightforward to check that there is a face $\Lambda$ of $\partialxii{\Gamma}{3}$ which contains $\Lambda'$ and also intersects $\partialxii{\Delta}{3}$. Let $\Lambda'' := \Lambda \cap \partialxii{\Delta}{3}$. If $(\partialxii{f}{3})_{\Lambda''}$ does not vanish after substituting $x_1 = a_1$ and $x_2 = a_2$, then one can find a zero of $(\partialxii{f}{3})_{\Lambda}$ of the form $(a_1,a_2,b)$ for some $b \in \kk^*$, and this would contradict partial non-degeneracy condition \eqref{partial-non-degeneracy} with $I = \{1,2,3\}$ and $\nu$ being an inner normal to $\Lambda$ in $\rr^3$. On the other hand, if $(\partialxii{f}{3})_{\Lambda''}$ vanishes after substituting $x_1 = a_1$ and $x_2 = a_2$, then $\Lambda''$ must be an edge parallel to $\Delta'$ and the partial non-degeneracy condition \eqref{partial-non-degeneracy} is violated with $I = \{1,2,3\}$ and $\nu$ being an inner normal to $\Lambda''$ in $\rr^3$. This proves that $V(\partialxone{f_{\Delta'}},\partialxii{f_{\Delta'}}{2}) \cap \Kstari = \emptyset$, as claimed. But then it is clear that $V(\jjj(f_\Delta)) \cap \Kstari$ is also empty, as required.
\end{proof}

\subsection{Conditions for validity of Kushnirenko's formula in positive characteristic} \label{newton-problem}
Consider the case that $p := \character(\kk) > 0$. The characterization of the conditions under which $\milnor(\scrA) = \newton(\nd(\scrA))$ established in \cref{milnor=newton}, in particular condition \eqref{mult-equality} of \cref{milnor=newton}, is {\em not} explicit. It would be interesting to find simple criteria under which this condition holds. The criterion from \cref{milnor=newton-prop-1} is sufficient, but not necessary, e.g.\ if $f = x + y^p$, then \eqref{milnor=newton-condition-1} fails (and $f$ is not Newton non-degenerate as well), but $\milnor(f) = \newton(\nd(f)) = 0$. Note that $f$ is {\em inner Newton non-degenerate}. Y.\ Boubakri, G.-M.\ Greuel and T.\ Markwig claim in \cite[Theorem 3.5]{boubakri-greuel-markwig} that $\milnor(f) = \newton(\nd(f))$ if $f$ is inner Newton non-degenerate, i.e.\ Wall's result \cite[Theorem 1.6]{wall} extends to positive characteristics; however, their proof is wrong, as explained by J.\ Stevens in \cite[Example 2.5]{stevens_conjectures_stably_newton_nondeg}. 

\begin{problem}
Determine necessary and sufficient conditions for $\milnor(f) = \newton(\nd(f))$ to hold when $\character(\kk) > 0$. In particular, does it hold when $f$ is inner Newton non-degenerate?
\end{problem}

\subsection{Monotonicity of Newton number}
\Cref{milnor-thm,milnor=newton-0} and the monotonicity of intersection multiplicity (\cref{monotonic-remark-0}) implies that the Newton number $\newton(\Gamma)$ of a diagram $\Gamma$ grows monotonically with $\Gamma$. Even though the arguments in this implication involve nontrivial algebraic geometry, the monotonicity of Newton numbers is a purely ``elementary'' convex geometric statement. It is a problem of V.\ I.\ Arnold \cite[Problem 1982-16]{arnold-problems} to find an elementary proof of this statement. S.\ K.\ Lando wrote in the commentary of that problem that he gave an elementary proof for $n = 2$. For $n = 3$ an elementary proof was given by S.\ Brzostowski, T.\ Krasi{\'n}ski and J.\ Walewska \cite{brzostowski-krasinski-walewska}. Since the monotonicity of the formula \eqref{mult-formula} for intersection multiplicity is obvious (see \cref{monotonic-remark-0}), one possible strategy would be to find an elementary proof of the following identity (which is an immediate consequence of \cref{milnor-thm,milnor=newton-0}):
\begin{align*}
\newton(\Gamma)
	&=  \sum_{I \in \tAone}
		\multzerostar{\Gamma^{I}_1, \Gamma^{I}_{j_2}}{\Gamma^{I}_{j_{|I|}}}
		\times
		\multzero{\pi_{[n]\setminus I}(\Gamma_{j'_1})}{\pi_{[n]\setminus I}(\Gamma_{j'_{n-|I|}})}
\end{align*}
where the right hand side is as in \cref{multiplicity-thm} with $\Gamma_j := \nd(\{\alpha - e_j:\alpha \in \Gamma + \rzeroo{n},\ \alpha - e_j \in \znzero\})$, $j = 1, \ldots, n$.


\chapter{Beyond this book} \label{further-chapter}

In this final chapter we outline some of the natural directions of further study for a reader of this book, and point out a few interesting recent works which are accessible to someone equipped with the knowledge of algebraic and toric varieties developed in this book.

\section{Toric varieties}
An obvious direction for further pursuit is to study toric varieties. Standard introductions to the theory of toric varieties include \cite{oda,fultoric,littlehalcox}. Unlike this book, in these expositions toric varieties are defined via {\em fans}, which reveals their combinatorial structure more cleanly. The combinatorics makes delightful appearances even in dimension two, e.g.\ continued fractions appear in resolutions of singularities \cite[Section 2.6]{fultoric} and Pick's formula appears in intersection theory \cite[Section 5.3]{fultoric}. A main attraction of toric varieties is their ``computability,'' due to which, in the words of Fulton, ``toric varieties have provided a remarkably fertile testing ground for general theories'' of algebraic geometry. As a recent example we mention the work \cite{grigoriev-milman-nash} of D.\ Grigoriev and P.\ Milman where a study of toric varieties provides an evidence that {\em Nash blow ups}\footnote{The tangent space $\tangent{X}{a}$ at a nonsingular point $a$ of a projective variety $X \subseteq \pp^n$ of dimension $d$ can be identified with a point on the {\em Grassmannian} $\grassman(d+1,n+1)$, which is the space of $(d+1)$-dimensional linear subspaces of $\kk^{n+1}$. Let $X'$ be the set of nonsingular points on $X$. The {\em Nash blow up} of $X$ is the closure in $X \times \grassman(d+1,n+1)$ of the graph of the map $X' \to \grassman(d+1, n+1)$ given by $x \mapsto \tangent{X}{a}$.}, which are in general somewhat poorly understood, might hold the key to a surprisingly efficient tool in the problem of resolution of singularities\footnote{John Nash asked in a private communication to Hironaka in early 1960s (as noted by Spivakovsky \cite{spinash}) whether singularities (in characteristic zero) can be resolved by a finite sequence of Nash blow ups. It remains an open problem for all dimensions $\geq 2$. M.\ Spivakovsky \cite{spinash} showed that in dimension two a finite sequence of {\em normalized Nash blow ups}, i.e.\ Nash  blow ups followed by normalizations, does suffice. Grigoriev and Milman \cite{grigoriev-milman-nash} show that for toric varieties the Nash blow-up corresponds to a higher dimensional analogue of the Euclidean division of integers, and the resolution of singularities via Nash blow ups becomes a problem in combinatorics; moreover for toric surfaces the resolution of singularities via normalized Nash blow ups has a {\em polynomial} complexity which was very surprising since all general algorithms, which are essentially based on the original algorithm of H.\ Hironaka, have a much higher complexity.}. The article \cite{grigoriev-milman-nash}, which only uses elementary properties of ``binomial varieties'' (which are slight generalizations of toric varities), would be very much accessible for a reader of this book after some familiarity with the notion of normalization (which is covered e.g.\ in \cite[Section 1.3]{littlehalcox}).

\section{Newton-Okounkov bodies}
The Newton-Okounkov body is a recent very fruitful generalization of the Newton polytope of a (Laurent) polynomial. It was originally introduced by A.\ Okounkov, who associated in \cite{okounkov-brunn-minkowski,okounkov-why} convex bodies to ample divisors on a smooth variety. R.\ Lazarsfeld and M.\ Mustata \cite{lazarsfeld-mustata-paper}, and independently K.\ Kaveh and A.\ Khovanskii \cite{kaveh-khovanskii-annals} made further generalizations of this construction. Since then there have been numerous articles on these structures. The construction of Kaveh and Khovanskii in particular is very elementary and leads to simple proofs of nontrivial results from intersection theory and also convex geometry. The series of articles by Kaveh and Khovanskii on Newton-Okounkov bodies, in particular \cite{khovanskii-kaveh,kaveh-khovanskii-annals}, would constitute an excellent reading material for a reader of this book. 

\section{B\'ezout problem}
Another natural extension of the material of this book would be a deeper study of the B\'ezout problem of counting numbers of solutions of systems of polynomials. A.\ Khovanskii studied higher dimensional analogues of this problem. In particular, he computes in \cite{khovanus} the Euler characteristics and in \cite{khovanskii-irred-components} the number of irreducible components of generic complete intersections on the torus. P.\ Philippon and M.\ Sombra \cite{philippon-sombra} studied a parametrized version of Bernstein's theorem over a nonsingular curve, and gave an answer in terms of an associated ``mixed integral,'' which is a generalization of mixed volume to concave functions. In a sequel to this book we describe an inductive algorithm to compute the precise number (counted with multiplicity) of solutions of any given system of $n$ polynomials in $n$ variables starting from the base estimate given by \cref{extended-bkk-bound-0,extended-bkk-bound}.

\section{Newton diagrams}
Newton diagrams continue to be deeply studied in numerous algebraic and geometric problems. G.\ Rond and B.\ Schober \cite{irred-power-series-rond-shoeber} gave a simple proof of a new criterion for irreducibility of a polynomial with power series coefficients in terms of its Newton diagram. The {\em \L ojasiewicz exponent}\footnote{If $f$ is a polynomial or an analytic function near the origin on $\cc^n$, its {\em \L ojasiewicz exponent} at the origin is the smallest $\theta > 0$ such that $|\nabla{f}| \geq C|z|^\theta$ for some $C > 0$ and all $z$ in a neighborhood of the origin.} of a complex hypersurface is an important invariant which is simple to define, but very hard to compute in many concrete situations. S.\ Brzostowski, T.\ Krasi\'nski and G.\ Oleksik \cite{brzostowski-krasinski-oleksik-lorface} determined the \L ojasiewicz exponent of a Newton non-degenerate surface singularity in $\cc^3$ in terms of the Newton diagram. In general, many invariants of a hypersurface singularity can be computed from its Newton diagram if the singularity is non-degenerate in some sense. As we have seen in \cref{milnor-chapter}, almost all singularities with a given diagram are non-degenerate. On the other hand, a given function is degenerate for most choices of coordinates, and for a degenerate function it is usually very difficult, if possible at all, to find coordinates in which it is non-degenerate. However, sometimes a degenerate function becomes non-degenerate after adding a quadratic form in new variables, and invariants computed from the Newton diagram of the new function can shed light on the original singularity. Motivated by these considerations, V.\ I.\ Arnold asked whether this is always possible \cite[Problems 1975-3, 1976-8]{arnold-problems}. This problem is little understood - see the article \cite{stevens_conjectures_stably_newton_nondeg} by J.\ Stevens for an exposition.

\section{Counting real zeroes}
In parallel to studying the relation between the topological complexity (as measured e.g.\ by number of zeroes) of systems defined over $\cc$ and their Newton polyhedra, V.\ I.\ Arnold and his students also explored the relation between the topological complexity of systems defined over $\rr$ and the number of nonzero terms (or some other algebraic/combinatorial measure of the complexity) of the systems, initiating the theory of ``Fewnomials.'' The canonical introduction to this subject is by A.\ Khovanskii \cite{khovanomials}, who proved the conjecture of A.\ Kushnirenko that the number of real zeroes is bounded by the number of nonzero terms of systems of equations. Khovanskii's bound however is not optimal, and many authors have worked on improving it, e.g.\ T-Y. Li, J.\ M.\ Rojas and X.\ Wang \cite{li-rojas-wang-real}, F.\ Bihan and F.\ Sottile \cite{bihan-sottile}. Sottile \cite{sottile-real} gave a beautiful introduction to these and other recent developments on the subject. A highly recommended read is Kushnirenko's letter\footnote{Dated February 26, 2008, available from Frank Sottile's website.} to Sottile where he describes the philosophy behind his conjectures, and early developments of the Fewnomial theory, including the contributions of K.\ Sevastyanov who tragically died young in an accident.

\appendix

\chapter{Commutative algebra results used in \cref{var-chapter} without a proof} \label{noproof-section}
The first major result we use without a proof is a special case of the ``principal ideal theorem'' of W.\ Krull, which is a fundamental tool in algebraic treatments of dimension. We use it frequently, starting from the proof that the dimension of the set of common zeroes of $k$ polynomials in $n$-variables has dimension at least $n-k$ (\cref{thm:pure-dimension}).

\begin{thm} [{\cite[Corollary 11.17, Theorem 11.25]{am}}] \label{thm:principal-ideal}
\index{Krull's principal ideal theorem}
\index{Principal ideal!theorem}
Let $R$ be a finitely generated integral domain over a field $k$, $f$ be a nonzero element of $R$ and $\ppp$ be an {\em isolated prime ideal}\footnote{This means $\ppp$ is minimal among the prime ideals containing the ideal $\langle f \rangle$ of $R$ generated by $f$. It follows from the results discussed in \cref{appendix:primary-decomposition} that under the hypotheses of \cref{thm:principal-ideal} the isolated prime ideals of $\langle f \rangle$ are precisely those (finitely many) prime ideals $\ppp_1, \ldots, \ppp_s$, determined by $f$ uniquely up to reordering, which satisfy the following property: $\sqrt{\langle f \rangle} = \bigcap_j \ppp_j$ and this presentation is {\em minimal} in the sense that $\sqrt{\langle f \rangle} \neq \bigcap_{j' \neq j} \ppp_{j'}$ for any $j = 1, \ldots, s$.} of the principal ideal $fR$ of $R$ generated by $f$. If $f$ is not a unit, then $\trd_k(R/\ppp) = \trd_k(R) - 1$ (where $\trd_k$ is the transcendence degree over $k$).
\end{thm}

We use a few basic properties of {\em localization} and {\em completion} of Noetherian rings. Short introductions to these notions are provided in \cref{local-ring-section,completion-section}. The following result, also due to Krull, is used in our proof that near a nonsingular point a variety is a ``locally complete intersection'' (\cref{local-equations}). 

\begin{thm}[{\cite[Corollary 10.18]{am}}] \label{thm:Krull}
Let $R$ be a Noetherian local ring with maximal ideal $\mmm$. If $\qqq$ is an ideal of $R$, then $\qqq = \bigcap_{r \geq 0} (\qqq + \mmm^r)$.
\end{thm}

\Cref{thm:exactly-complete} below is used in this book for the first time in the proof that the completion of the local ring at a point of a variety is isomorphic to the quotient of a power series ring (\cref{complete-quotient}). We use \cref{thm:Krull,thm:exactly-complete} only in the case that $R$ is of the form $\{f/g: f, g \in k[x_1, \ldots, x_n],\ g(0) \neq 0\}$, where $k$ is a field. 

\begin{thm}[{Exactness of completion \cite[Proposition 10.12]{am}}] \label{thm:exactly-complete}
Let $I \subset J$ be ideals of a Noetherian ring $R$ and let $\bar R := R/I$. Let $\hat R$ (respectively $\hat {\bar R}$) be the completion of $R$ (respectively $\bar R$) with respect to $J$ (respectively $J\bar R$). Then $\hat {\bar R} \cong \hat R/I \hat R$.
\end{thm}

\chapter{Miscellaneous commutative algebra} \label{algebra-section}
By a ``ring'' in this book we mean a commutative ring with identity. In this chapter we briefly recall several notions related to rings which are used in the book. 

\section{Integral domain, UFD, PID}
A \index{Zero-divisor}{\em zero-divisor} in a ring $R$ is a nonzero element $f$ such that $fg = 0$ for some nonzero $g \in R$. An \index{Integral domain}{\em integral domain} is a ring without zero-divisors. An \index{Irreducible!element of an integral domain}\index{Prime!element of an integral domain}{\em irreducible element} or {\em prime} of an integral domain is a nonzero element which is a non-unit\footnote{A unit of a ring is an element which has a multiplicative inverse.} and which can not be expressed as a product of two non-units. The notion of prime elements is modeled after properties of prime numbers; indeed, the prime elements in the ring of integers are precisely the prime numbers (and their negatives). A \index{Unique factorization domain}\index{UFD}{\em unique factorization domain} (in short, {\em UFD}) is an integral domain such that every non-unit can be represented as a product of irreducible elements and this representation is unique up to multiplication by a unit or reordering of the irreducible factors. The ring of integers is a UFD, and so are polynomial rings (in a finite or infinite number of variables) over the integers or over a field. The ring of integers satisfies the additional property that each of its ideals is {\em principal}, i.e.\ generated by a single element. An integral domain with the latter property is called a \index{Principal ideal!domain}\index{PID}{\em principal ideal domain} (in short, {\em PID}).  The following is usually called the ``fundamental theorem'' or the ``structure theorem'' of finitely generated modules over a PID. Its proof can be found in any standard introductory abstract algebra text, e.g.\ \cite{dummit-foote}.

\begin{thm} \label{fundamentally-over-PID}
Let $M$ be a finitely generated module over a PID $R$. Then
$$M \cong R^r \dsum R/\langle p_1^{\alpha_1} \rangle
			  \dsum R/\langle p_2^{\alpha_2} \rangle \cdots
			  \dsum R/\langle p_k^{\alpha_k} \rangle$$
for some $r, \alpha_1, \ldots, \alpha_k \geq 0$ and primes $p_1, \ldots, p_k$.
\end{thm}

\section{Prime and maximal ideals} 
Given elements $g_1, \ldots, g_k$ of a ring $R$, we write $\langle g_1, \ldots, g_k\rangle$ to denote the ideal generated by $R$. Recall that a \index{Prime!ideal}{\em prime} ideal of $R$ is a proper ideal $\ppp$ such that $fg \in \ppp$ implies either $f$ or $g$ is in $\ppp$. A \index{Maximal ideal}{\em maximal} ideal of $R$ is a proper ideal $\mmm$ such that the only proper ideal of $R$ containing $\mmm$ is $\mmm$ itself. It is straightforward to check that every maximal ideal is prime. The following are some other basic properties of prime and maximal ideals - see e.g.\ \cite[Chapter 1]{am} for their proof. 

\begin{thm} \label{thm:mp}
Let $R$ be a nonzero ring (remember that by ``ring'' we mean a commutative ring with identity). 
\begin{enumerate}
\item \label{mp:max-existence} Given any non-unit $g \in R$, there is a maximal ideal of $R$ containing $g$.
\item \label{mp:prime-avoidence} If $\qqq$ is an ideal of $R$ contained in the union of finitely many prime ideals $\ppp_1, \ldots, \ppp_k$ of $R$, then $\qqq \subseteq \ppp_j$ for some $j$. 
\item \label{mp:prime-containment} If $\ppp$ is a prime ideal of $R$ containing the intersection of finitely many ideals $\qqq_1, \ldots, \qqq_k$, then $\ppp \supseteq \qqq_j$ for some $j$. If $\ppp = \bigcap_{j=1}^k \qqq_j$, then $\ppp = \qqq_j$ for some $j$.
\item \label{mp:nilradical} The intersection of all prime ideals of $R$ is precisely its \index{Nilradical}{\em nilradical}, i.e.\ the ideal consisting of all \index{Nilpotent}{\em nilpotent}\footnote{$g \in R$ is called {\em nilpotent} if there is $n \geq 1$ such that $g^n = 0$.} elements of $R$. 
\end{enumerate}
\end{thm}

Assertion \eqref{mp:prime-avoidence} of \cref{thm:mp} is sometimes referred to as \index{Prime!avoidance}``prime avoidance''; see \cite[Lemma 3.3]{eisenview} for a stronger variation and geometric interpretation. 

\section{Noetherian rings, Hilbert's basis theorem, annihilators}
A ring $R$ is called \index{Noetherian ring}{\em Noetherian}\footnote{in honour of Emmy Noether, who was one of the founders of modern algebra. A very influential and highly respected mathematician, she was denied a permanent academic position for most of her career on account of being a woman.} if every nonempty set $\scrI$ of ideals in $R$ has a maximal element with respect to $\subseteq$, i.e.\ there is an ideal $\qqq \in \scrI$ such that the only element of $\scrI$ containing $\qqq$ is $\qqq$ itself. The following is a basic characterization of Noetherian rings (see e.g.\ \cite[Proposition 6.2]{am}). 

\begin{prop} \label{prop:characterize-noetherian}
Given a ring $R$, the following are equivalent: 
\begin{enumerate}
\item $R$ is Noetherian. 
\item Every ideal of $R$ is {\em finitely generated} (i.e.\ generated by finitely many elements). 
\end{enumerate}
\end{prop}

 A field is Noetherian, since the only proper ideal of a field is the zero ideal. It is straightforward to check that the ring of integers, or more generally, any PID is Noetherian. Hilbert's basis theorem (\cref{thm:Hilbert-basis}), which we prove next, implies that (quotients of) polynomial rings in finitely many variables over a field are also Noetherian. Most of the rings that appear in this book are Noetherian. The ring $R[x_1, x_2,\cdots]$ of polynomials in infinitely many variables $x_j$, $j \geq 1$, over a ring $R$ is {\em not} Noetherian, since e.g.\ the sequence of ideals $\langle x_1, \ldots, x_n \rangle$, $n \geq 1$, does not have any maximal element. 
 
 \begin{thm}[Hilbert's basis theorem] \label{thm:Hilbert-basis}
 \index{Hilbert's theorem!Basis theorem}
 The ring of polynomials in finitely many variables over a Noetherian ring is also Noetherian.
 \end{thm}
 
 \begin{proof}
 It suffices to show that $R[x]$ is Noetherian if $R$ is a Noetherian ring and $x$ is an indeterminate over $R$. Let $\qqq$ be an ideal of $R[x]$. We will find a finite set of generators of $R$. Given $f = a_0 +a_1 x + \cdots + a_d x^d \in R[x ]$, $a_d \neq 0$, we call $a_d$ the {\em leading coefficient} of $f$ . Let $\rrr \subseteq R$ be the set of leading coefficients of all polynomials in $\qqq$. It is easy to check that $\rrr$ is an ideal in $R$. Since $R$ is Noetherian, there are polynomials $f_1, \ldots, f_r \in \qqq$ whose leading coefficients generate $\rrr$. Let $M := \max\{\deg(f_i): i = 1, \ldots, r\}$. For each $m < M$, let $\rrr_m$ be the ideal in $R$ generated by the leading coefficients of all polynomials $f \in \qqq$ such that $\deg(f) = m$. By Noetherianity of $R$, there is a finite set $\{f_{mj}\}_j$ of polynomials in $\qqq$ of degree $m$ whose leading coefficients generate $\rrr_m$. We claim that $\qqq$ is generated by all the $f_i$ and all the $f_{mj}$. Indeed, let $f \in \qqq$. If $d := \deg(f) \geq M$, then since the leading coefficient of $f$ can be expressed as a sum of $b_i \in R$ times the leading coefficient of $f_i$, it follows that the degree of $f - \sum_i b_i x^{d - \deg(f_i)}f_i$ is smaller than $d$. Similarly, if $d < M$, then the degree of $f - \sum_j c_j f_{dj}$ is smaller than $d$ for some $\{c_j\}_j \subseteq R$. In any event, the degree of $f$ can be reduced by subtracting an element from the ideal $\qqq'$ generated by the $f_i$ and the $f_{mj}$. Repeating this process finitely many times yields that $f \in \qqq'$, i.e.\ $\qqq = \qqq'$, which proves the theorem.
 \end{proof}
 
 Given $g \in R$, one writes $(0:g)$ for the set of all $f \in R$ such that $fg = 0$; in other words $(0:g)$ is the \index{Annihilator}{\em annihilator} of $g$. It is clear that the set of zero-divisors is precisely the union of the annihilators of the nonzero elements of $R$. 

\begin{prop} \label{prop:zero-union-prime}
The set of zero-divisors in a Noetherian ring $R$ is a union of prime ideals which are of the form $(0:g)$ for some $g \in R$. 
\end{prop} 

\begin{proof}
Let $\scrI$ be the set of all ideals of a Noetherian ring $R$ of the form $(0:g)$, where $g$ varies over nonzero elements of $R$. By Noetherianity of $R$, each $\qqq \in \scrI$ is contained in a maximal element $\ppp \in \scrI$. We will show that $\ppp$ is prime. Indeed, pick $g \in R$ such that $\ppp = (0:g)$. If $\ppp$ is not prime, then there are $f_1, f_2 \in R \setminus \ppp$ such that $f_1f_2 \in \ppp$, i.e.\ $f_1f_2g = 0$. But then $(0:f_2g)$ properly contains $\ppp$. This contradiction (to the maximality of $\ppp$) completes the proof. 
\end{proof}

In \cref{prop:zero-union-minimally-prime} below we prove a stronger version of \cref{prop:zero-union-prime}. 

\section{(Algebraic) Field extensions} \label{algebraic-extension-section}
A \index{Field extension}{\em field extension} $K/F$ is simply a pair consisting of a field $K$ and a subfield $F$ of $K$. The \index{Degree!of a field extension}{\em degree} $[K: F]$ of the extension $K/F$ is the dimension of $K$ as a vector space over $F$; we say that $K/F$ is a {\em finite} extension if $[K:F] < \infty$. 

\begin{example}
$\cc$ is a finite extension of $\rr$ with $[\cc: \rr] = 2$. $\rr$ is {\em not} a finite extension of $\qq$ (since e.g.\ $\qq$ is countable, but $\rr$ is not). 
\end{example}

An element $\alpha \in K$ is \index{Algebraic!element}{\em algebraic} over $F$ if there is a nonzero polynomial $f \in F[t]$ (where $t$ is an indeterminate) such that $f(\alpha) = 0 \in K$. A \index{Transcendental element}{\em transcendental} element over $F$ is an element which is not algebraic over $F$. The extension $K/F$ is called \index{Algebraic!extension}\index{Field extension!algebraic}{\em algebraic} if all elements of $K$ are algebraic over $F$.

\begin{example}
All elements of $F$ are trivially algebraic over $F$. All finite extensions of $F$ are algebraic. Indeed, if $d := [K:F] < \infty$, then $1, \alpha, \alpha^2, \ldots, \alpha^d$ are linearly dependent over $F$ for each $\alpha \in K$, which leads to an algebraic equation of $\alpha$ over $F$. If $x$ is an indeterminate over $F$, then no element of $F(x) \setminus F$ is algebraic over $F$ (where $F(x)$ is the field of rational functions in $x$). For $\qq \subseteq \rr$, one has $\sqrt{2}$ is algebraic over $\qq$, but $\pi$ is not. All elements of $\cc$ are algebraic over $\rr$. 
\end{example}

$K$ is said to be an \index{Algebraic!closure}{\em algebraic closure} of $F$ if $K$ is algebraic over $F$ and if every polynomial $f \in F[x]$ has a root in $K$. We say that $K$ is \index{Algebraically closed field}{\em algebraically closed} if $K$ is an algebraic closure of itself. The existence of algebraic closures and algebraically closed fields are by now standard results of algebra (see e.g.\ \cite[Section 13.4]{dummit-foote}). The {\em fundamental theorem of algebra}\footnote{See \cite{mathoverflow-fundamental-thm-algebra} for a delightful collection of many proofs of the fundamental theorem of algebra.} states that $\cc$ is algebraically closed, and in particular, $\cc$ is the algebraic closure of $\rr$. An algebraically closed field must be infinite, since given finitely many elements $\alpha_1, \ldots, \alpha_d$ of a field $K$, no $\alpha_j$ is a root of $1 + \prod_{j=1}^d (x - \alpha_j)$. 

\section{Hilbert's Nullstellensatz}
Let $K$ be a field. Given a subset $X$ of $K^n$ we write $I(X)$ for the set of all polynomials in $K[x_1, \ldots, x_n]$ which are zero at all points of $X$; it is straightforward to check that $I(X)$ is an ideal of $K[x_1, \ldots, x_n]$. Conversely, if $\qqq$ is an ideal of $K[x_1, \ldots, x_n]$, we write $V(\qqq)$ for the set of points on $K^n$ on which each element of $\qqq$ is zero; it is straightforward to check that $I(V(\qqq)) \supseteq \qqq$ and $V(I(X)) \supseteq X$. The following result of David Hilbert describes the basic correspondence between $I(\cdot)$ and $V(\cdot)$ in the case that $K$ is algebraically closed.

\begin{thm}[Hilbert's Nullstellensatz] \label{thm:Hilbert-nulls}
\index{Hilbert's theorem!Nullstellensatz}
\index{Nullstellensatz}
Assume $K$ is algebraically closed. Then for each ideal $\qqq$ of $K[x_1, \ldots, x_n]$, $I(V(\qqq))$ is the radical $\sqrt{\qqq}$ of $\qqq$. In particular, the maximal ideals of $K[x_1, \ldots, x_n]$ are of the form $I(a)$ for $a \in K^n$.
\end{thm}

In this section we prove the Nullstellensatz following an argument of Enrique Arrondo \cite{arrondo2006}. %
We start with a simple version of the ``normalization lemma'' of Emmy Noether. 

\begin{lemma} \label{lemma:noether-norrmalization-0}
Let $K$ be an infinite field, and $f \in K[x_1, \ldots, x_n]$ be a polynomial with $d := \deg(f) \geq 1$. Then there are $\lambda_1, \ldots, \lambda_{n-1} \in K$ such that the coefficient of $x_n^d$ in $f(x_1 + \lambda_1 x_n, \ldots, x_{n-1} + \lambda_{n-1} x_n, x_n)$ is nonzero. 
\end{lemma}

\begin{proof}
Let $f_d$ be the homogeneous component of $f$ of degree $d$. The coefficient of $x_n^d$ in $f(x_1 + \lambda_1 x_n, \ldots, x_{n-1} + \lambda_{n-1} x_n, x_n)$ is $f_d(\lambda_1, \ldots, \lambda_{n-1}, 1)$. Note that $f_d(x_1, \ldots, x_{n-1}, 1)$ is a nonzero polynomial in $(x_1, \ldots, x_{n-1})$. It is then straightforward to show, e.g.\ by induction on $n$, that there are $\lambda_1, \ldots, \lambda_{n-1} \in K$ such that $f_d(x_1, \ldots, x_{n-1}, 1) \neq 0$ (this is \cref{exercise:pol-zero-on-Kn}). 
\end{proof}

\begin{thm}[Weak Nullstellensatz] \label{thm:Hilbert-nulls-weak}
Let $K$ be an algebraically closed field and $\qqq$ be a proper ideal of $K[x_1, \ldots, x_n]$. Then $V(\qqq) \neq \emptyset$.  
\end{thm}

\begin{proof}
\Woutlog\ we may assume $\qqq \neq 0$. We proceed by induction on $n$. If $n = 1$, then $\qqq$ must be a principal ideal generated by a nonconstant polynomial $f$. Since $K$ is algebraically closed, $f$ has a root $a \in K$. Then $V(\qqq)$ contains $a$ and therefore is nonempty. Now assume $n > 1$. Due to \cref{lemma:noether-norrmalization-0} after a change of coordinates if necessary we may assume that $\qqq$ contains a nonzero polynomial $g$ of the form $g = g_0 + g_1x_n + \cdots + g_{e-1}x_n^{e-1} + x_n^e$ where $e \geq 1$ and $g_0, \ldots, g_{e-1} \in K[x_1, \ldots, x_{n-1}]$. Let $\qqq' := \qqq \cap K[x_1, \ldots, x_{n-1}]$. By the inductive hypothesis there is $(a_1, \ldots, a_{n-1}) \in V(\qqq')$. 

\begin{proclaim}
$\rrr := \{f(a_1, \ldots, a_{n-1}, x_n): f \in \qqq\}$ is a proper ideal of $K[x_n]$. 
\end{proclaim}

\begin{proof}
It is clear that $\rrr$ is an ideal of $K[x_n]$. Assume to the contrary that $\rrr = K[x_n]$. Then there is $f \in \qqq$ such that $f(a_1, \ldots, a_{n-1}, x_n) = 1$. Write $f = f_0 + f_1x_n + \cdots + f_dx_n^d$ with $f_0, \ldots, f_d \in K[x_1, \ldots, x_{n-1}]$. Then $f_0(a_1, \ldots, a_{n-1}) = 1$ and $f_1(a_1, \ldots, a_{n-1}) = \cdots = f_d(a_1, \ldots, a_{n-1}) = 0$. Let $R \in K[x_1, \ldots, x_{n-1}]$ be the \index{Resultant}{\em resultant} of $f$ and $g$ with respect to $x_n$, i.e.\ $R$ is the determinant of the following $(d+e) \times (d+e)$ matrix: 
\begin{align*}
	  \begin{pmatrix}
	  f_0  	& f_1 		& \cdots & \cdots 	& \cdots  & f_d		 & 0 		& \cdots   & 0 \\
	  0 	 & f_0		 & f_1	   & \cdots   & \cdots  & \cdots  & f_d		& 0			 & 0 \\
	  		 &  			& \ddots & 			  &			   &			&			& \ddots  &	\\
	  0		 & \cdots   & 0  	   & f_0	   & f_1      & \cdots  & \cdots & \cdots  & f_d \\
	  g_0  	& g_1 		& \cdots & g_{e-1}& 1  		  & 0		  & 0 		& \cdots  	& 0 \\
	  0 	 & g_0		 & g_1	   & \cdots  & g_{e-1}& 1		  & 0		& 0			  & 0 \\
		  	  &  			& \ddots & 			 &		   	  &			   & \ddots&			 & \\
		  	  &  			&			& \ddots &		   	  &			   &		  &	\ddots	 & \\
	  0		 & \cdots   & \cdots &	0		 & g_0      & g_1	  & \cdots & g_{e-1} & 1 
   \end{pmatrix}
\end{align*}
where the nonzero entries of the first $e$ rows are translations of $(f_0, \ldots, f_d)$, and the next $d$ rows are translations of $(g_0, \ldots, g_{e-1}, 1)$. If the first column of the above matrix is replaced by the first column plus $x_n$ times the second column plus $x_n^2$ times the third column and so on, then the first colum turns into the column vector (which is the transpose of) $(f, x_nf, \cdots, x_n^{e-1}f, g, x_ng, \ldots, x_n^ {d-1} g)$. Expanding the resulting determinant along the first column then shows that $R$ is in the ideal generated by $f$ and $g$. It follows that $R \in \qqq'$. On the other hand, evaluating the entries of the above matrix at $(a_1, \ldots, a_{n-1})$ converts it into a lower triangular matrix whose diagonal entries are all $1$; this implies that $R(a_1, \ldots, a_{n-1}) = 1$, which contradicts the fact that $(a_1, \ldots, a_{n-1}) \in V(\qqq')$ and proves the claim.
\end{proof}
The preceding claim and the $n = 1$ case of the theorem (which we already proved) then show that there is $a_n \in K$ such that $f(a_1, \ldots, a_n) = 0$ for all $f \in \qqq$, which completes the proof of the theorem.
\end{proof}

Now we prove the Nullstellensatz (\cref{thm:Hilbert-nulls}) via a classical argument of Rabinowitsch\footnote{Although Rabinowitsch's trick, which first appeared in \cite{rabinowtrick}, is widely known in algebraic geometry, the man Rabinowitsch is not. The most widely accepted account is that he is the mathematical physicist George Yuri Rainich, who had published several articles under his original name ``Rabinowitsch'' before immigrating to USA and changing his name. There are however some doubt about this claim - see the comments to the {\em MathOverflow} answer \cite{mathoverflow-rabinowitsch}.
}. Let $\qqq$ be an ideal of $K[x_1, \ldots, x_n]$. The containment $\sqrt{\qqq} \subseteq I(V(\qqq))$ is straightforward to verify. We will show that $\sqrt{\qqq} \supseteq I(V(\qqq))$. Fix a set of generators $f_1, \ldots, f_m$ of $\qqq$. Given $f \in I(V(\qqq))$, consider the ideal $\rrr$ of $K[x_1, \ldots, x_{n+1}]$ generated by $f_1, \ldots, f_m, x_{n+1}f - 1$. Then $V(\rrr) = \emptyset$, since whenever $f_1 = \cdots = f_m = 0$, then $f = 0$ and therefore $x_{n+1}f - 1 = -1$. \Cref{thm:Hilbert-nulls-weak} then implies that $1 \in \rrr$, i.e.\ there is an equation of the form $1 = \sum_i f_i h_i + (x_{n+1}f - 1)h$. Substituting $x_{n+1} = 1/y$ in the preceding equation and clearing the denominator via multiplying by an appropriate power of $y$ yields an equation of the form $y^N = \sum_i f_i \tilde h_i + (f - y)\tilde h$ for some $\tilde h_1, \ldots, \tilde h_m, \tilde h \in K[x_1, \ldots, x_n, y]$. Substituting $y = f$ in the preceding equation then shows that $f^N \in \qqq$, as required. \qed

\section{Nakayama's lemma}
In this section we present a brief discussion of {\em Nakayama's lemma} following \cite{am}. Let $R$ be a ring, and $\jjj$ be its \index{Jacobson ideal}{\em Jacobson ideal}, i.e.\ the intersection of all the maximal ideals of $R$.

\begin{lemma} \label{jacobson-unit}
$1 + f$ is a unit in $R$ for each $f \in \jjj$.
\end{lemma}

\begin{proof}[Proof of \cref{jacobson-unit}]
Pick $f \in \jjj$. If $1 + f$ is not a unit, then there is a maximal ideal $\mmm$ of $R$ that contains $1 + f$ (\cref{thm:mp}, assertion \eqref{mp:max-existence}). Since $\mmm$ contains both $f$ and $1 + f$, it contains $1$, which is a contradiction. 
\end{proof}

\begin{lemma}[Nakayama's lemma] \label{Nakayama}
\index{Nakayama's lemma}
Let $M$ be a finitely generated $R$-module and $\qqq$ be an ideal of $R$ contained in the Jacobson ideal of $R$. If $\qqq M = M$, then $M = 0$.
\end{lemma}

\begin{proof}
Assume $M \neq 0$ and let $u_1, \ldots, u_k$ be a minimal set of generators of $M$. Since $u_k \in \qqq M$, there is an equation of the form $u_k = \sum_{j=1}^k f_ju_j$, with each $f_j \in \qqq$. It follows that $(1-f_k)u_k = \sum_{j=1}^{k-1}u_jf_j$. Since $1 - f_k$ is a unit in $R$ (\cref{jacobson-unit}), it follows that $u_k$ is an $R$-linear combination of $u_1, \ldots, u_{k-1}$. This means that $M$ can be generated by $u_1, \ldots, u_{k-1}$, a contradiction.
\end{proof}
Let $R$ be a local ring with (the unique) maximal ideal $\mmm$, and $\kk = R/\mmm$ be its residue field. Let $M$ be a finitely generated $R$-module. Since $M/\mmm M$ is annihilated by $\mmm$, it is a finitely generated $R/\mmm$-module, i.e.\ a finite dimensional vector space over $\kk$.

\begin{cor} \label{Nakayama-generation}
Let $m_1, \ldots, m_k \in M$ be such that their images span $M/\mmm M$ as a vector space over $\kk$. Then $m_1, \ldots, m_k$ generate $M$ as a module over $R$.
\end{cor}

\begin{proof}
Let $N$ be the submodule of $M$ generated by $m_1, \ldots, m_k$. Since the image of $N$ in $M/\mmm M$ is all of $M/\mmm M$, it follows that $M = N + \mmm M$. But then $\mmm (M/N) = (N + \mmm M)/N = M/N$. \Cref{Nakayama} then implies that $M = N$, as required.
\end{proof}

\section{Localization, local rings} \label{local-ring-section}
If $S$ is a multiplicative closed subset of a ring $R$, then the \index{Localization}{\em localization} $R_S$ of $R$ with respect to $S$ is the equivalence class of ``fractions'' $\{f/g: f \in R,\ g \in S\}$ under the equivalence relation that $f/g \sim f'/g'$ provided $(fg' - f'g)h = 0$ for some $h \in S$. It is easy to check that $R_S$ is a commutative ring with respect to the usual rules of addition and multiplication of fractions. Two cases of localizations are especially relevant for our purpose:
\begin{itemize}
\item Case 1: $S = \{f^k\}_k$ for some $f \in R$. In this case we denote $R_S$ by $R_f$, and say that $R_f$ is the {\em localization of $R$ at $f$.}
\item Case 2: $S = R \setminus \ppp$ for some prime ideal $\ppp$ of $R$. In this case we denote $R_S$ by $R_\ppp$, and say that $R_\ppp$ is the {\em localization of $R$ at $\ppp$.}
\end{itemize}

The following proposition compiles some of the basic properties of localizations. We refer to \cite[Chapter 3]{am} for a lucid discussion of these and other basic properties of localizations of commutative rings.

\begin{prop} \label{prop:localization-basic}
$R_S$ is the zero ring if and only if $0 \in S$. If $R_S$ is not the zero ring, then (the equivalence class of) every element of $S$ is invertible in $R_S$, and every ideal $I$ of of $R_S$ is generated by $R \cap I$. In particular, if $R$ is Noetherian, then so is $R_S$. 
\end{prop}

A \index{Local ring}{\em local ring} is a ring with a unique maximal ideal. \Cref{prop:localization-basic} implies that the localization $R_\ppp$ of a ring $R$ at a prime ideal $\ppp$ is a local ring, and the maximal ideal of $R_\ppp$ is generated by $\ppp$. 

\begin{example} \label{example:local-ring-0}
Given $a = (a_1, \ldots, a_n) \in \kk^n$, where $\kk$ is a field, the set $S_a$ of all polynomials $f \in R := \kk[x_1, \ldots, x_n]$ such that $f(a) \neq 0$ is multiplicatively closed. In fact $\mmm_a := R \setminus S_a = \{f \in R: f(a) = 0\}$ is a maximal ideal\footnote{It is straightforward to check that $\mmm_a$ is an ideal. It is maximal since if $f \not\in \mmm_a$, then $f - f(a) \in \mmm_a$ and $1 = \frac{1}{f(a)}(f - (f - f(a))$.}. It follows that $R_{S_a} = R_{\mmm_a}$ is a local ring and its maximal ideal is generated by $\mmm_a$. 
\end{example}

\section{Discrete valuation rings}\label{discrete-valuection}

A \index{Discrete valuation}{\em discrete valuation} on a field $K$ is a surjective map $\nu$ from $K$ onto $\zz \cup \{\infty\}$ such that $\nu(0) = \infty$, $\nu(fg) = \nu(f) + \nu(g)$ and $\nu(f+g) \geq \min\{\nu(f), \nu(g)\}$ for each $f,g \in K$. It is straightforward to check that the set of all $f \in K$ such that $\nu(f) \geq 0$ is a subring of $K$; it is called the \index{Valuation ring}{\em valuation ring} of $\nu$.

\begin{example}
A basic example of a discrete valuation is the \index{Order!of a rational function in one variable}{\em order} of rational functions in one variable over a field $\kk$. Recall that the order of $f = \sum_j a_j t^j \in \kk[t]$, where $t$ is an indeterminate, is $\ord(f) := \min\{j: a_j \neq 0\}$, and the order can be extended to the field $\kk(t)$ of rational functions by defining $\ord(f/g) := \ord(f) - \ord(g)$. The valuation ring of $\ord$ is precisely the localization of $\kk[t]$ at the maximal ideal generated by $t$, i.e.\ the subring $\{f/g: f, g \in \kk[t],\ g(0) \neq 0\}$ of $\kk(t)$. 
\end{example}

A \index{Discrete valuation!ring}{\em discrete valuation ring} is an integral domain which is the valuation ring of a discrete valuation on its field of fractions. We now record some properties of discrete valuation rings: their verification is straightforward, and is left as an exercise. 

\begin{prop} \label{prop:discrete-properties}
Let $\nu$ be a discrete valuation on a field $K$, and $R := \{f \in K: \nu(f) \geq 0\}$ be the valuation ring of $\nu$. 
\begin{enumerate}
\item The units of $R$ are precisely the elements $f$ in $K$ with $\nu(f) = 0$.
\item $\mmm := \{f \in K: \nu(f) > 0\}$ is the unique maximal ideal of $R$; in particular, $R$ is a {\em local ring}.
\item If $f$ is any element in $\mmm$ such that $\nu(f) = 1$, then $\mmm$ is the (principal) ideal of $R$ generated by $f$. More generally, every proper ideal of $R$ is a principal ideal generated by $f^k$ for some $k \geq 1$.
\item \label{exercise:discrete-properties-dim} If $R$ contains a field $\kk$ which is isomorphic to $R/\mmm$, then for each $g \in R$, $R/gR$ is a vector space over $\kk$ of dimension $\nu(g)$.
\item \label{uniquely-dvr} $R$ uniquely determines $\nu$, i.e.\ if $\nu'$ is a discrete valuation on $K$ such that $R$ is also the valuation ring of $\nu'$, then $\nu' = \nu$. \qed
\end{enumerate}
\end{prop}

\Cref{prop:discrete-properties} in particular implies that a discrete valuation ring $R$ is a local ring whose maximal ideal is principal. A \index{Parameter!of a discrete valuation ring}{\em parameter} of $R$ is a generator of its maximal ideal.

\section{Krull dimension} \label{Krull-section}
A {\em chain} of prime ideals of length $n \geq 0$ in a ring $R$ is a finite sequence $\ppp_0 \subsetneq \ppp_1 \subsetneq \cdots \subsetneq \ppp_n$ of prime ideals of $R$. The \index{Krull dimension}{\em Krull dimension} of $R$ is the supremum of the lengths of all chains of prime ideals in $R$. 

\begin{example} \label{Krull-examples}
Since the only maximal ideal of a nontrivial field $F$ is the zero ideal, the Krull dimension of $F$ is zero. Each of the following rings has Krull dimension one: the ring of integers, the ring of polynomials in one variable over a field, a discrete valuation ring.
\end{example}

\section{Primary decomposition} \label{appendix:primary-decomposition}
A proper ideal $\qqq$ of a ring $R$ is \index{Primary!ideal}{\em primary} if $fg \in \qqq$ implies either $f \in \qqq$ or $g^k \in \qqq$ for some $k \geq 1$.

\begin{prop} \label{prop:radical-of-primaries}
Let $\qqq$ be an ideal of a ring $R$ and $\sqrt{\qqq}$ be the \index{Radical!of an ideal}radical\footnote{The radical of $\qqq$ is the ideal consisting of all $g \in R$ such that $g^k \in \qqq$ for some $k \geq 1$.} of $\qqq$. 
\begin{enumerate}
\item If $\qqq$ is primary, then $\sqrt{\qqq}$ is prime. 
\item \label{radical-of-primaries:maximal} If $\sqrt{\qqq}$ is {\em maximal}, then $\qqq$ is primary. 
\item \label{radical-of-primaries:intersection} Write $\ppp := \sqrt{\qqq}$. Assume $\qqq$ and $\qqq'$ are primary ideals of $R$ such that $\sqrt{\qqq'} = \sqrt{\qqq} = \ppp$. Then $\qqq \cap \qqq'$ is also primary and $\sqrt{\qqq \cap \qqq'} = \ppp$. 
\end{enumerate}   
\end{prop}

\begin{proof}
The proof of the first assertion is straightforward from the definition of a primary ideal - it is left as an exercise. For the second assertion, assume $\sqrt{\qqq}$ is maximal. 

\begin{proclaim} \label{claim:radmax-unit}
If $g \not\in \sqrt{\qqq}$, then $g$ is a unit in $R/\qqq$.
\end{proclaim}

\begin{proof}
Indeed, there is $f \in \sqrt{\qqq}$ and $u, v \in R$ such that $gu + fv = 1$. If $k \geq 1$ is such that $f^k \in \qqq$, then the relation $(gu + fv)^k = 1$ reduces to $gu' = 1$ in $R/\qqq$ for an appropriate $u' \in R$.
\end{proof}

Pick $f, g \in R$ such that $fg \in \qqq$. If $g \not\in \sqrt{\qqq}$, then \cref{claim:radmax-unit} implies that $f \in \qqq$. Therefore $\qqq$ is primary. It remains to prove the third assertion. It is straightforward to check that $\sqrt{\qqq \cap \qqq'} = \sqrt{\qqq} \cap \sqrt{\qqq'} = \ppp$. Now pick $fg \in \qqq \cap \qqq'$ such that $f \not\in \qqq \cap \qqq'$. Then $f$ is not in either $\qqq$ or $\qqq'$, and since both are primary with radical $\ppp$, it follows that $g \in \ppp = \sqrt{\qqq \cap \qqq'}$, so that $\qqq \cap \qqq'$ is also primary.
\end{proof}

In the examples below $\kk$ denotes a field. 

\begin{example} \label{example:xy-y^2}
The ideal $\qqq$ of $\kk[x,y]$ generated by $x^2$ and $xy$ is {\em not} primary, since $xy \in \qqq$, but $x \not\in \qqq$ and $y \not\in \sqrt{\qqq} = \langle x \rangle$. In particular, assertion \eqref{radical-of-primaries:maximal} of \cref{prop:radical-of-primaries} does not hold if ``maximal'' is replaced by ``prime.''
\end{example}

\begin{example}
The primary ideals of $\kk[t]$ are the zero ideal and the ideals generated by $f^n$, $n \geq 1$, for irreducible polynomials $f$. More generally, the nonzero primary ideals of a PID are precisely the ideals generated by powers of irreducible elements. 
\end{example}

\begin{example}
The ideal $\qqq$ of $\kk[x,y]$ generated by $x$ and $y^2$ is primary (\cref{prop:radical-of-primaries}, assertion \eqref{radical-of-primaries:maximal}). The radical of $\qqq$ is the maximal ideal $\mmm$ generated by $x, y$. Note that $\mmm \supsetneq \qqq \supsetneq \mmm^2$, so that $\qqq$ is {\em not} a prime-power. 
\end{example}

A \index{Primary!decomposition}{\em primary decomposition} of an ideal $\qqq$ is an expression 
\begin{align}
\qqq = \bigcap_{j = 1}^k \qqq_j \label{eqn:minimal-decomposition}
\end{align}
of $\qqq$ as a finite intersection of primary ideals $\qqq_j$. Given such a primary decomposition, excluding any redundant ideal from the right hand side and then grouping the ideals with the same radical, it can be ensured that
\begin{itemize}
\item $\qqq$ can {\em not} be expressed as an intersection of less than $k$ of the $\qqq_j$, and
\item the radicals $\sqrt{\qqq_j}$ of $\qqq_j$ are distinct (due to assertion \eqref{radical-of-primaries:intersection} of \cref{prop:radical-of-primaries}); 
\end{itemize}
in that case we say that \eqref{eqn:minimal-decomposition} is a {\em minimal} primary decomposition, and that the prime ideals $\sqrt{\qqq_j}$ are \index{Prime!ideal!associated with an ideal}{\em associated with} $\qqq$ (that $\sqrt{\qqq_j}$ are prime follows from \cref{prop:radical-of-primaries}). Note that
\begin{align*}
\sqrt{\qqq} = \bigcap_{j=1}^k \sqrt{\qqq_j}
\end{align*}
is a (not necessarily minimal) primary decomposition of the radical $\sqrt{\qqq}$ of $\qqq$; in particular, $\sqrt{\qqq}$ is a finite intersection of {\em prime} ideals associated with $\qqq$. 

\begin{example}
Let $\qqq := \langle x^2, xy \rangle \subset \kk[x,y]$. We saw in \cref{example:xy-y^2} that $\qqq$ is not primary. Each of the following is a minimal primary decomposition of $\qqq$ (since $\langle x \rangle$ is prime and the radical of the other ideal in the decomposition is maximal)
\begin{align*}
\qqq = \langle x \rangle \cap \langle x^2, y \rangle 
	= \langle x \rangle \cap \langle x^2, xy, y^k \rangle\ (k \geq 1)
\end{align*}
In particular, minimal primary decompositions are in general {\em not} unique\footnote{However, the {\em isolated} primary ideals (i.e.\ primary ideals whose radicals are minimal among the radicals of all primary ideals appearing in the decomposition) in a primary decomposition are in fact unique (e.g.\ the ideal $\langle x \rangle$ will appear in every primary decomposition of $\langle x^2, xy \rangle$) - see \cite[Corollary 4.11]{am}.}.
\end{example}

\begin{example} \label{example:local-primary}
Let $S$ is a multiplicatively closed subset of $R$ and $\qqq$ be an ideal of $R$ such that $\qqq \cap S = \emptyset$. If $\qqq$ is primary, then it is straightforward to check that the ideal $\qqq_S$ generated by $\qqq$ in the localization $R_S$ of $R$ is also primary, and moreover, $\sqrt{\qqq_S} = (\sqrt{\qqq})_S$. It follows that if $\qqq = \bigcap_{j=1}^k \qqq_j$ is a minimal primary decomposition of $\qqq$, then the following is a minimal primary decomposition of $\qqq_S$: 
\begin{align*}
\qqq_S = \bigcap_i (\qqq_{j_i})_S
\end{align*}
where $\{j_i\}$ is the subset of $\{1, \ldots, k\}$ consisting of those $j$ such that $\qqq_j \cap S = \emptyset$.   
\end{example}

\begin{prop} \label{prop:minimal-prime}
Let $\qqq$ be an ideal in $R$ which has a primary decomposition. Then every prime ideal in $R$ containing $\qqq$ contains one of the prime ideals associated with $\qqq$. The minimal ideals among those associated with $\qqq$ are precisely the minimal elements in the set of all ideals in $R$ containing $\qqq$.
\end{prop}

\begin{proof}
Follows immediately from assertion \eqref{mp:prime-containment} of \cref{thm:mp}.  
\end{proof}

The following is a fundamental property of Noetherian rings; it is a combination of \cite[Theorem 4.5, Proposition 7.17]{am}. 

\begin{thm} \label{thm:primary-decomposition}
Assume $R$ is Noetherian. Then every proper ideal has a primary decomposition. In particular, every \index{Radical!ideal}radical ideal\footnote{An ideal is {\em radical} if it equals its own radical.} is a finite intersection of prime ideals. The prime ideals associated with an ideal $\qqq$ are precisely those prime ideals of $R$ which occur in the set of ideals $(\qqq:f) := \{g \in R: fg \in \qqq\}$ as $f$ varies over $R$, and hence are uniquely determined by $\qqq$. 
\end{thm}

\begin{cor} \label{prop:zero-union-minimally-prime}
The set of zero-divisors in a Noetherian ring $R$ is the union of the prime ideals associated with the zero ideal. Every prime ideal of $R$ contains a prime ideal associated with the zero ideal. The minimal ideals among the prime ideals associated with the zero ideal are precisely the minimal elements of the set of {\em minimal}\footnote{A {\em minimal} prime ideal $\ppp$ of a ring $R$ is a prime ideal $\ppp$ such that the only prime ideal of $R$ contained in $\ppp$ is $\ppp$ itself.} prime ideals of $R$. In particular, every element of a minimal prime ideal of $R$ is a zero-divisor. 
\end{cor}

\begin{proof}
The first assertion follows from combining \cref{prop:zero-union-prime,thm:primary-decomposition}. The remaining assertions then follow from \cref{prop:minimal-prime}.
\end{proof}

\begin{cor}\label{prop:non-zero-divisor-restriction}
Let $f$ be a non zero-divisor in a ring $R$ and $\ppp$ be a minimal prime ideal of $R$. If $R$ is Noetherian, then (the image of) $f$ remains a non zero-divisor in $R/\ppp$. 
\end{cor}

\begin{proof}
\Cref{prop:zero-union-minimally-prime} implies that $f \not\in \ppp$. Therefore, if $g \not\in \ppp$, then $fg \not\in \ppp$. Therefore $f$ is not a zero-divisor in $R/\ppp$, as required. 
\end{proof}

\section{Length of modules} \label{length-section}
Let $M$ be a module over a ring $R$. A \index{Composition series}{\em composition series} of $M$ of length $n \geq 0$ is a sequence
\begin{align}
M = M_0 \supsetneq M_1 \supsetneq \cdots \supsetneq M_n = 0 \label{composition-series}
\end{align}
of $R$-submodules which is ``maximal,'' or equivalently, each quotient $M_{i-1}/M_i$, $1 \leq i \leq n$, is {\em simple}, (that is, has no nonzero proper submodule). Not every module has a composition series. In fact the following are equivalent (\cite[Proposition 6.8]{am}):
\begin{enumerate}
\item $M$ has a composition series,
\item $M$ satisfies {\em both} ascending and descending chain conditions\footnote{$M$ satisfies {\em ascending} (respectively, {\em descending}) {\em chain condition} if for every chain $M_0 \subseteq M_1 \subseteq \cdots$ (respectively, $M_0 \supseteq M_1 \supseteq \cdots$) of submodules of $M$, there is $k$ such that $M_j = M_k$ for all $j \geq k$.}. 
\end{enumerate}

If $M$ has a composition series, then all composition series of $M$ have the same length (\cite[Proposition 6.7]{am}). The \index{Length of a module}{\em length} $l(M)$ of $M$ is defined to be infinite if it has no composition series; otherwise $l(M)$ is the length of any composition series of $M$. If \eqref{composition-series} is a composition series of $M$, then each $M_{i-1}/M_i$ is isomorphic to $R/\mmm$ for some maximal ideal $\mmm$ of $R$. 

\begin{example}
$\zz$ does not have a composition series as a module over itself. Indeed, $\zz$-modules of $\zz$ are simply ideals of $\zz$. Thereore if $1 = n_0 < n_1 < \cdots$ is an infinite sequence of positive integers, then there is an infinite sequence $\langle n_0 \rangle \supsetneq \langle n_0n_1 \rangle \supsetneq \langle n_0n_1n_2 \rangle \supsetneq \cdots$ of ideals of $\zz$; consequently $l(\zz) = \infty$. On the other hand, if $R = \zz/12\zz$, then both of the following are composition series of $R$ (as a module over $R$ or over $\zz$): 
\begin{align*}
&\zz/12\zz  \supsetneq 2\zz/12\zz \supsetneq 4\zz/12\zz \supsetneq 0 \\
&\zz/12\zz  \supsetneq 3\zz/12\zz \supsetneq 6\zz/12\zz \supsetneq 0 
\end{align*}
It follows that $l(\zz/12\zz) = 3$. Note that for each of these composition series the successive quotients are isomorphic to either $\zz/2\zz$ or $\zz/3\zz$. In general, given primes $p_j$ and nonnegative integers $n_j$, $j = 1, \ldots, k$, it is straightforward to check that the length of $\zz/(\prod_j p_j^{n_j})\zz$ as a module over $\zz$ is $\sum_j n_j$ and each successive quotient of each of its composition series is isomoprhic to $\zz/p_j\zz$ for some $j$.  
\end{example}

\begin{prop} \label{exercise:composition-quotient}
Let $M$ be a module over a ring $R$ with a composition-series $M = M_0 \supsetneq M_1 \supsetneq \cdots \supsetneq M_n = 0$. Fix $i$, $1 \leq i \leq n$. 
\begin{enumerate}
\item \label{composition-quotient:max} Each $M_{i-1}/M_i$ is isomorphic to $R/\mmm$ for some maximal ideal $\mmm$ of $R$. 
\item \label{composition-quotient:k} Let $\mmm$ be as in assertion \eqref{composition-quotient:max}. If there is a subfield $k$ of $R$ which maps isomorphically onto $R/\mmm$, then $M_{i-1}/M_i$ is a one dimensional vector space over $k$. 
\end{enumerate}
\end{prop}

\begin{proof}
Pick any nonzero element $\bar m \in M_{i-1}/M_i$. Since $M_{i-1}/M_i$ is simple, the map $R \to M_{i-1}/M_i$ given by $f \mapsto fm$ is a surjective map and its kernel must be a maximal ideal of $R$. This implies the first assertion. The second assertion follows immediately from the first. 
\end{proof}

\begin{cor} \label{prop:local-field-finite-length}
Assume $R$ is a local ring with maximal ideal $\mmm$, and there is a field $k \subseteq R$ which maps isomorphically onto $R/\mmm$. Then for any $R$-module $M$, the length of $M$ as an $R$-module is equal to the length of $M$ as a $k$-module, which in turn is equal to the dimension of $M$ as a $k$-vector space. \qed
\end{cor}

\begin{prop} \label{prop:minimal-localization-finite-length}
Let $\ppp$ be a \index{Minimal!prime ideal}{\em minimal} prime ideal of a Noetherian ring $R$. Then the localization $R_\ppp$ of $R$ at $\ppp$ has a finite length as an $R_\ppp$-module. 
\end{prop}

\begin{proof}
Since $\ppp$ is a minimal prime ideal of $R$, the ideal generated by $\ppp$ is the unique prime ideal of $R_\ppp$, i.e.\ the Krull dimension of $R_\ppp$ is zero. Since $R_\ppp$ is also Noetherian, it follows that it is \index{Artinian ring}{\em Artinian}\footnote{A ring $R$ is {\em Artinian} if it satisfies the descending chain condition as a module over itself, i.e.\ if for every descending chain of ideals $I_0 \supseteq I_1 \supseteq \cdots$, there is $k$ such that $I_j = I_k$ for all $j \geq k$.} (\cite[Theorem 8.5]{am}) as well. This implies that the length $R_\ppp$ as a module over itself is finite (\cite[Proposition 6.8]{am}), as required.
\end{proof}

\section{(In)Separable field extensions} \label{separable-section}
Let $K/F$ be a finite extension of a field $F$ and $\alpha \in K$ be algebraic over $F$. Pick $f \in F[t]$ such that $f$ is \index{Monic polynomial}{\em monic}\footnote{A polynomial in a single variable $t$ is monic if the coefficient of its highest degree term is $1$.}, $f(\alpha) = 0$, and $f$ has the smallest possible degree among all nonzero polynomials $g \in F[t]$ such that $g(\alpha) = 0$. It is straightforward to see, since $F[t]$ is a PID, that these properties uniquely determine $f$, and $f$ is {\em irreducible} in $F[t]$; we say that $f$ is the \index{Minimal!polynomial}{\em minimal polynomial} of $\alpha$ over $F$. We say that $\alpha$ is \index{Separable!element}{\em separable} over $F$ if $f$ has distinct roots in the algebraic closure $\bar F$ of $F$. Recall that a polynomial $g$ in one variable (over any field) has distinct roots if and only if the ``greatest common divisor'' $\gcd(g, g')$ of $g$ and its derivative $g'$ is a nonzero constant. Since $\gcd(f,f')$ is an element of $F[t]$ and since $f$ is irreducible in $F[t]$, this immediately implies the following:

\begin{prop} \label{prop:separable-condition}
$\alpha$ is separable over $F$ if and only if the derivative of its minimal polynomial is a nonzero polynomial in $F[t]$. 
\end{prop} 

We say that $K/F$ is a \index{Separable!extension}\index{Field extension!separable}{\em separable} extension if every element of $K$ is separable (and in particular, algebraic) over $F$.

\begin{example} \label{example:separable-0}
If $p := \character(F) = 0$, then the derivative of any nonconstant polynomial over $F$ is a nonzero polynomial. It follows that every element which is algebraic over $F$ is also separable, i.e.\ every algebraic extension of $F$ is separable. Now assume $p > 0$. Let $x, y$ be indeterminates over $F$. Then $y^p - x$ is irreducible in $F(x)[y]$, so that $R := F(x)[y]/\langle y^p - x \rangle$ is an integral domain. Let $K$ be the field of fractions of $R$. Then $K/F$ is a algebraic, but {\em not} separable (since the image of $y$ in $K$ not separable over $F(x)$). 
\end{example}

We now prove the \index{Primitive!element theorem}{\em primitive element theorem} in the case that the base field is infinite. It is also true for finite fields (see e.g.\ \cite[Chapter 14, Proposition 17]{dummit-foote}), but we do not use that case in this book. The following proof is taken from \cite[Problem 6-31]{fulturve}. 

\begin{thm} \label{thm:primitive-element-infinite}
Let $K/F$ be a finite separable extension of fields.  Assume $|F| = \infty$. Then there is $\alpha \in K$ such that $K = F(\alpha)$. Moreover, given $\alpha_1, \ldots, \alpha_n \in K$ such that $K = F(\alpha_1, \ldots, \alpha_n)$, one can ensure that $\alpha = \sum_j \lambda_j \alpha_j$ for some $\lambda_1, \ldots, \lambda_n \in F$.
\end{thm}

\begin{proof}
A straightforward induction on $n$ shows that it suffices to prove the second assertion of the theorem for $n = 2$. So assume $K = F(\alpha_1, \alpha_2)$. Let $f_i \in F[t]$ be the minimal polynomial of $\alpha_i$. Since $K/F$ is separable, over the algebraic closure $\bar K$ of $K$ there are factorizations of the form $f_1= \prod_{j = 1}^{d_j} (t - \alpha_{1,j})$, where $\alpha_1 = \alpha_{1,1} \neq \alpha_{1, j}$ for any $j > 1$. Since $F$ is infinite, there is $\lambda \in F$ such that $\lambda \alpha_1 + \alpha_2 \neq \lambda \alpha_{1, j_1} + \alpha'_{2}$ for all $j_1 > 1$ and all roots $\alpha'_2$ of $f_2(t)$ in $\bar K$. Let $\alpha := \lambda \alpha_1 + \alpha_2$ and $f(t) := f_2(\alpha - \lambda t) \in F'(t)$, where $F' := F(\alpha)$. Then it is straightforward to check that $f(\alpha_1) = 0$, and for each $j_1 > 1$, $f(\alpha_{1,j_1}) \neq 0$. It follows that the greatest common divisor of $f$ and $f_1$ in $F'(t)$ is $t - \alpha_1$, which means the ideal generated by $f$ and $f_1$ in $F'(t)$ is the same as the ideal generated by $t - \alpha_1$. In particular, $t - \alpha_1 = g(t)f(t) + g_1(t)f_1(t)$ for some $g, g_1 \in F'(t)$. But then $\alpha_1 = -(g(0)f(0) + g_1(0)f_1(0)) \in F'$. It follows that $\alpha_2 := \alpha - \lambda \alpha_1 \in F'$ as well, so that $F' = K$, as required. 
\end{proof}

\begin{rem} \label{rem:primitive-element-infinite}
The arguments of the proof of \cref{thm:primitive-element-infinite} does not use the separability of $\alpha_2$. Therefore the following generalization of \cref{thm:primitive-element-infinite} holds: ``Let $K = F(\alpha_1, \ldots, \alpha_n)$ be a finite extension of $F$ such that $\alpha_1, \ldots, \alpha_{n-1}$ are separable over $F$. If $|F| = \infty$, then $K = F(\sum_j \lambda_j \alpha_j)$ for some $\lambda_1, \ldots, \lambda_n \in F$.''
\end{rem}

A field $F$ of characteristic $p$ is called {\em perfect} if either $p = 0$ or if for every $\alpha \in F$, there is $\beta \in F$ such that $\beta^p = \alpha$. Following \cite[Appendix 5, Proposition 1]{shaf1} we now prove a result of \index{Schmidt's theorem}F.\ K.\ Schmidt for infinite perfect fields (note that it is also true in the case when the (perfect) field is finite, see e.g.\ \cite[Chapter II, Theorem 31]{zsI}).

\begin{cor} \label{thm:Schmidt-infinite}
Let $F$ be a perfect field and $K/F$ is a finitely generated field extension. If $|F| = \infty$, then there are $\alpha_1, \ldots, \alpha_{d+1} \in K$ such that $\alpha_1, \ldots, \alpha_d$ are algebraically independent over $F$, $\alpha_{d+1}$ is separable over $F(\alpha_1, \ldots, \alpha_d)$, and $K = F(\alpha_1, \ldots, \alpha_{d+1})$. 
\end{cor}

\begin{proof}
Pick $\beta_1, \ldots, \beta_n \in K$ such that $K = F(\beta_1, \ldots, \beta_n)$. Let $d$ be the maximal number of the $\beta_j$ which are algebraically independent over $F$. \Woutlog\ we may assume that $\beta_1, \ldots, \beta_d$ are algebraically independent over $F$. 

\begin{proclaim} \label{claim:Schmidt-infinite}
For each $j = d, \ldots, n$, reordering $\beta_1, \ldots, \beta_j$ if necessary, we may ensure that $\beta_1, \ldots, \beta_d$ are algebraically independent over $F$ and $F(\beta_1, \ldots, \beta_j) = F(\beta_1, \ldots, \beta_d, \gamma_j)$ for some $\gamma_j \in K$. 
\end{proclaim}

\begin{proof}
The claim is clearly true for $j = d$. We proceed by induction and assume it is true for $j$, $d \leq j \leq n - 1$. By definition of $d$, $\beta_{j+1}$ is algebraic over $F(\beta_1, \ldots, \beta_d)$. Let $f$ be a nonzero irreducible polynomial in $F[t_1, \ldots, t_d, t_{j+1}]$ (where the $t_i$ are indeterminates) such that $f(\beta_1, \ldots, \beta_d, \beta_{j+1}) = 0$. We claim that there is $i$ such that $\partial{f}/\partial{t_i}$ is a nonzero polynomial. Indeed, otherwise  $p := \character(F)$ is nonzero and each $t_i$ occurs in $f$ in powers which are multiples of $p$. Then $f$ is of the form $f = \sum b_{i_1, \ldots, i_d, i_{j+1}}t_1^{pi_1} \cdots t_d^{pi_d}t_{j+1}^{pi_{j+1}}$. Choose $a_{i_1, \ldots, i_d, i_{j+1}} \in F$ such that $a_{i_1, \ldots, i_d, i_{j+1}}^p = b_{i_1, \ldots, i_d, i_{j+1}}$. Then $f = (\sum a_{i_1, \ldots, i_d, i_{j+1}} t_1^{i_1} \cdots t_d^{i_d}t_{j+1}^{i_{j+1}})^p$, which contradicts the irreducibility of $f$. Reorder $\beta_1, \ldots, \beta_d, \beta_{j+1}$ if necessary so that $i = j+1$. Then $\beta_1, \ldots, \beta_d$ still remain algebraically independent over $F$, and $\beta_{j+1}$ is {\em separable} over $F(\beta_1, \ldots, \beta_d)$ (\cref{prop:separable-condition}). It then follows from the inductive hypothesis and \cref{rem:primitive-element-infinite} that $F(\beta_1, \ldots, \beta_{j+1}) = F(\beta_1, \ldots, \beta_d, \gamma_j, \beta_{j+1}) =  F(\beta_1, \ldots, \beta_d, \gamma_{j+1})$ for some $\gamma_{j+1} \in K$, as required.
\end{proof}

Let $\gamma_n$ be as in \cref{claim:Schmidt-infinite}. The arguments of the proof of \cref{claim:Schmidt-infinite} show that reordering $\beta_1, \ldots, \beta_d, \gamma_n$ if necessary, we may ensure that $\gamma_n$ is separable over $F(\beta_1, \ldots, \beta_d)$, which proves the corollary. 
\end{proof}

Given an algebraic extension $K/F$, it is a standard result from algebra (see e.g.\ \cite[Section 14.9]{dummit-foote}) that the set $\sepclosure{F}$ of all elements in $K$ which are separable over $F$ is a field; we say that $\sepclosure{F}$ is the \index{Separable!closure}{\em separable closure} of $F$ in $K$. The \index{Separable!degree!of a field extension}\index{Degree!separable}{\em separable degree} (respectively \index{Inseparable!degree!of a field extension}\index{Degree!inseparable}{\em inseparable degree}) of $K/F$ is $[\sepclosure{F}: F]$ (respectively $[K: \sepclosure{F}]$). The following result compiles a few basic properties regarding these notions (see e.g.\ \cite[Section 14.9]{dummit-foote}):

\begin{prop} \label{prop:sep-insep-deg}
Assume $p := \character(K)$ is nonzero\footnote{In zero characteristic $\sepclosure{F} = K$ (\cref{example:separable-0}), so that the conclusions of \cref{prop:sep-insep-deg} are trivially true.}. Then the extension $K/\sepclosure{F}$ is {\em purely inseparable}, i.e.\ for each $\alpha \in K$, there is a nonzero integer $m$ such that the minimal polynomial of $\alpha$ over $\sepclosure{F}$ is of the form $t^{p^m} - \alpha$. Moreover, the degree of $K/F$ is the product of its separable degree and the inseparable degree, i.e.\ 
\begin{align}
[K: F] = [K :\sepclosure{F}] [\sepclosure{F}: F]
\label{eqn:sep-insep-deg}
\end{align}
\end{prop}

\begin{example} \label{example:separable-1}
Assume $p := \character(F) > 0$. Choose positive integers $k, q$, where $q$ is relatively prime to $p$. If $x, y$ are indeterminates, it follows as in \cref{example:separable-0} that $R := F(x)[y]/\langle y^{qp^k} - x \rangle$ is an integral domain and the quotient field $K$ of $R$ is not a separable extension over $F(x)$. Note that (the image of) $y^{p^k}$ in $K$ is separable over $F(x)$, since its minimal polynomial in $F(x)[t]$ is $t ^q - x$, which is separable. It is not hard to see $\sepclosure{F(x)}$ is the field $F(x)(y^{p^k})$ generated over $F(x)$ by $y^{p^k}$. It follows that $[\sepclosure{F(x)}: F(x)] = q$, $[K: \sepclosure{F(x)}] = p^k$, so that $[K: \sepclosure{F(x)}][\sepclosure{F(x)}: F(x)] = qp^k = [K:F(x)]$, as implied by \eqref{eqn:sep-insep-deg}. 
\end{example}

\section{Rings of formal power series over a field} \label{formal-power-section}
A \index{Formal power series}\index{Power series}{\em (formal)  power series} over a field $\kk$ in variables $(x_1, \ldots, x_n)$ is a formal expansion of the form $f = \sum_{\alpha \in \zzeroo{n}} c_\alpha x^\alpha$, where $c_\alpha \in \kk$, and $x^\alpha$ is a shorthand for the monomial $x_1^{\alpha_1} \cdots x_n^{\alpha_n}$. The \index{Power series!ring}\index{Formal power series!ring}{\em power series ring} $\kk[[x_1, \ldots, x_n]]$ is the set of all such power series; it has the structure of a ring induced by the usual multiplication and product of power series. In this section we write $\hat R$ for $\kk[[x_1, \ldots, x_n]]$. It is straightforward to see that $\hat R$ is an algebra over $\kk$, and also an {\em integral domain}.  

\begin{prop} \label{prop:power-local}
$\hat R$ is a {\em local ring}, i.e.\ it has a unique maximal ideal. The maximal ideal $\hat \mmm$ of $\hat R$ consists of all power series with zero constant term. Every element in $\hat R \setminus \hat m$ is a unit in $\hat R$. 
\end{prop}

\begin{proof}
That $\hat \mmm$ is a maximal ideal follows from the isomorphism $\hat R/\hat R \cong k$. If $f \in \hat R \setminus \hat \mmm$, then $f = c(1 + g)$ for some $c \in \kk \setminus \{0\}$ and $g \in \hat \mmm$. It is straightforward to check that \begin{align}
\tilde f := \frac{1}{c}(1 + \sum_{d \geq 1}(-1)^dg^d) \label{eqn:power-inverse}
\end{align}
is a well defined element in $\hat R$ and $\tilde f f = 1$. Therefore $f$ is a unit. It follows that $\hat \mmm$ is the unique maximal ideal of $\hat R$.
\end{proof}

We now describe all the $\kk$-algebra automorphisms\footnote{A $\kk$-algebra automorphism of $\hat R$ is a $\kk$-algebra isomorphism from $\hat R$ to itself.} of $\hat R$. Given $f = \sum_\alpha c_\alpha x^\alpha \in \hat R$, Recall that the \index{Homogeneous!component!of a formal power series}{\em homogeneous component of $f$ of degre $d$} is $f_d := \sum_{|\alpha| = d} c_\alpha x^\alpha$, where $|\alpha| := \sum_j \alpha_j$; note that $f_d$ is a {\em polynomial} for each $d \geq 0$. The \index{Order!of a formal power series}{\em order} $\ord(f)$ of a power series $f$ is the smallest $m$ such that $c_\alpha \neq 0$ for some $\alpha$ with $|\alpha| = m$. 

\begin{thm}\label{thm:linear-power}
Let $f_1, \ldots, f_n \in \hat \mmm$. 
\begin{enumerate}
\item \label{linear-power:hom} The map $x_j \mapsto f_j$, $j = 1, \ldots, n$, induces a well-defined $\kk$-algebra homomorphism from $\Phi:\hat R \to \hat R$. All $\kk$-algebra homomorphisms from $\hat R$ to itself is of this form. 
\item \label{linear-power:iso} $\Phi$ is an isomorphism if and only if the linear parts of the $f_j$ are linearly independent over $\kk$. 
\end{enumerate}
\end{thm}

\begin{proof}
Assertion \eqref{linear-power:hom} is straightforward. The ($\Leftarrow$) implication of assertion \eqref{linear-power:iso} is also straightforward in the case that all $f_j$ are linear. Moreover, the $(\Rightarrow)$ implication of assertion  \eqref{linear-power:iso} also follows from the linear case of the ($\Leftarrow$) implication. Now we prove the general case of the ($\Leftarrow$) implication of assertion \eqref{linear-power:iso}. Due to the linear case, after an automorphism of $\hat R$ if necessary we may assume that the linear part of each $f_j$ is precisely $x_j$, $j = 1, \ldots, n$. It is then straightforward to check that for any $f \in \hat R$, 
\begin{align}
\parbox{0.66\textwidth}{
$\ord(\Phi(f)) = \ord(f)$. Moreover, if $m = \ord(f)$, then the homogeneous component of degree $m$ of $\Phi(f)$ is the same as that of $f$.  
}\label{linear-power:order}
\end{align}
Property \eqref{linear-power:order} immediately implies that $\Phi$ is injective. To see that $\Phi$ is surjective, pick $f \in \hat R$. Let $m_0 := \ord(f)$ and $g_0$ be the degree-$m_0$ homogeneous component of $f$. Property \eqref{linear-power:order} implies that the order $m_1$ of $f_1 := f - \Phi(g_0)$ is greater than $ m_0$; let $g_1$ be the degree-$m_1$ homogeneous component of $f_1$ and set $f_2 := f_1 - \Phi(g_1)$. Continuing this way we get a power series $g = g_0 + g_1 + \cdots  \in \hat R$ such that $\Phi(g) = f$, as required.
\end{proof}

\begin{cor} \label{cor:hat-prime}
If $f_1, \ldots, f_r \in \hat \mmm$ has linearly independent (over $\kk$) linear parts, then they generate a prime ideal in $\hat R$. \qed
\end{cor}
%

Let $R := \{g/f: f, g \in \kk[x_1, \ldots, x_n],\ f(0) \neq 0\}$; in other words, $R$ is the {\em localization} of the polynomial ring $k[x_1, \ldots, x_n]$ at the ideal generated by polynomials which vanish at the origin (\cref{example:local-ring-0}). There is a natural map $R \to \hat R$ which is identity on polynomials, and for each polynomial $f$ such that $f(\origin) \neq 0$, maps $1/f$ to a power series as in \eqref{eqn:power-inverse}. The following is straightforward to verify:

\begin{prop} \label{prop:R-into-R-hat}
The map $R \to \hat R$ is injective. 
\end{prop}

From now on we will treat $R$ as a subring of $\hat R$. Recall that $R$ is also a local ring and $\mmm := \hat \mmm \cap R$ is the unique maximal ideal of $R$ (\cref{example:local-ring-0}). 

\begin{prop} \label{prop:m-hat-m}
$(\hat \mmm)^d \cap R = \mmm^d$ for each $d \geq 0$. 
\end{prop}

\begin{proof}
It is clear that $(\hat \mmm)^d \cap R \supseteq \mmm^d$ for each $d$. for the opposite inclusion, assume $g/f = h \in \hat \mmm^d$ for some $f, g \in \kk[x_1, \ldots, x_n]$ such that $f(0) \neq 0$. But then the identity $g = fh \in \hat R$ implies that all monomials in $g$ must have order $\geq d$. This means $g \in \mmm^d$, as required. 
\end{proof}

In the next section we use ``monomial orders'' on $\hat R$ to derive some of its basic properties in a simple way, e.g.\ we show that $\hat R$ is Noetherian (\cref{noetherian-power-series}) and that the dimensions of the quotients of $\hat R$ are ``finitely determined'' (\cref{finite-determinacy}). 

\section{Monomial orders on rings of formal power series} \label{power-order-section}

\subsection{Monomial orders on $\zzeroo{n}$} 
A \index{Monomial order!on $\zzeroo{n}$}{\em monomial order} on $\zzeroo{n}$ is a binary relation $\preceq$ on $\zzeroo{n}$ such that 
\begin{defnlist}
\item \label{total-prop} $\preceq$ is a {\em total order}\footnote{A {\em total order} on a set $S$ is a binary relation $\preceq$ on $S$ which is reflexive (i.e.\ $x \preceq x$ for all $x \in S$), transitive (i.e.\ if $x \preceq y$ and $y \preceq z$ then $x \preceq z$), antisymmetric (i.e.\ if $x \preceq y$ and $y \preceq x$ then $x = y$), and totally comparable (i.e.\ either $x \preceq y$ or $y \preceq x$ for each $x, y \in S$).},
\item \label{compatible-prop} $\preceq$ is compatible with the addition on $\zzeroo{n}$, i.e.\ if $\alpha \preceq \beta$, then $\alpha + \gamma \preceq \beta + \gamma$ for all $\gamma \in \zzeroo{n}$, and
\item \label{zero-prop} $0 \preceq \alpha$ for each $\alpha \in \zzeroo{n}$.
\end{defnlist}
We show below in \cref{corner-well} that every monomial order $\preceq$ is also a {\em well order} on $\zzeroo{n}$, i.e.\ for every nonempty subset $S$ of $\zzeroo{n}$, there is a unique $\alpha \in S$ such that $\alpha \preceq \alpha'$ for all $\alpha' \in S$. 

\begin{example} \label{lexicographic-example}
The \index{Lexicographic order}{\em lexicographic order} $\preceq_{lex}$ on $\zzeroo{n}$ is defined as follows: if $\alpha, \beta \in \zzeroo{n}$, then $\alpha \preceq_{lex} \beta$ if either $\alpha = \beta$, or $\alpha \neq \beta$ and the first nonzero coordinate from the left of $\alpha - \beta$ is negative. Replacing ``left'' to ``right'' in the preceding definition leads to \index{Reverse lexicographic order}{\em reverse lexicographic order} $\preceq_{rlex}$. It is straightforward to check that $\preceq_{lex}$ and $\preceq_{rlex}$ are monomial orders. 
\end{example}

A {\em corner point} of a subset $S$ of $\zzeroo{n}$ is an element $\alpha \in S$ such that there is no $\alpha' \in S$, $\alpha' \neq \alpha$, such that $\alpha = \alpha' + \beta$ for some $\beta \in \zzeroo{n}$. 

\begin{lemma} \label{corner-lemma}
Let $S$ be a nonempty subset of $\zzeroo{n}$. The set $\scrC_S$ of corner points of $S$ is finite and nonempty. Moreover, $S + \zzeroo{n} = \scrC_S + \zzeroo{n}$.
\end{lemma}  

\begin{proof}
For any $\alpha = (\alpha_1, \ldots, \alpha_n) \in S$, define $S_{\leq \alpha} := \{\beta \in S: \alpha - \beta \in \zzeroo{n}\}$. Since $S_{\leq \alpha}$ is a finite set, it has a corner point $\beta$. It is clear that $\beta$ is also a corner point of $S$ and $\alpha \in \beta + \zzeroo{n}$. This proves the second assertion, and in addition shows that $\scrC_S$ is nonempty. It remains to show that $\scrC_S$ is finite. Assume to the contrary that it is infinite. Let $\alpha^0 = (\alpha^0_1, \ldots, \alpha^0_n)$ be an arbitrary element of $S$. For each $\alpha = (\alpha_1, \ldots, \alpha_n) \in \scrC_S \setminus \{\alpha^0\}$, there is $j$ such that $\alpha_j < \alpha^0_j$. Fix $j_1$ such that there are infinitely many $\alpha \in \scrC_S$ with $\alpha_{j_1} < \alpha^0_{j_1}$. Since there are finitely many choices for $\alpha_{j_1}$, it follows that there is $a_{j_1} < \alpha^0_{j_1}$ such that $\scrC^1_S:= \{\alpha \in \scrC_S: \alpha_{j_1} = a_{j_1}\}$ is infinite. Now fix $\alpha^1 \in \scrC^1_S$. Replacing $\alpha^0$ by $\alpha^1$ and $\scrC_S$ by $\scrC^1_S$ and running the above procedure yields $j_2 \neq j_1$ and $a_{j_2} < \alpha^1_{j_2}$ such that $\scrC^2_S := \{\alpha \in \scrC^1_S: \alpha_{j_2} = a_{j_2}\}$ is infinite. Continuing this process we will end up with an infinite set $\scrC^n_S$. But this is absurd, since $\scrC^n_S$ will consist of a single element $(a_1, \ldots, a_n)$ by construction. This contradiction implies that $\scrC_S$ was finite to begin with, which completes the proof.
\end{proof}

\begin{cor} \label{corner-well}
Every monomial order on $\zzeroo{n}$ is also a well order on $\zzeroo{n}$. 
\end{cor}

\begin{proof}
Let $\preceq$ be a monomial order on $\zzeroo{n}$ and $S$ be a nonempty subset of $\zzeroo{n}$. Let $\scrC_S$ be the set of corner points of $S$. Since $\scrC_S$ is finite and nonempty, it has a unique minimal element $\beta_0$ with respect to $\preceq$. For every $\alpha \in S$, \cref{corner-lemma} implies that there is $\beta \in \scrC_S$ such that $\alpha - \beta \in \zzeroo{n}$, so that properties \ref{compatible-prop} and \ref{zero-prop} of monomial orders imply that $\beta \preceq \alpha$. It follows that $\beta_0$ is the minimal element of $\scrS$ with respect to $\preceq$. 
\end{proof}

\begin{cor} \label{corner-convex}
If $S \subseteq \zzeroo{n}$, then the convex hull of $S + \rzeroo{n}$ is a (convex) polyhedron. 
\end{cor}

\begin{proof}
\Cref{corner-lemma} implies that $S + \rzeroo{n} = \scrC_S + \rzeroo{n}$, where $\scrC_S$ is the {\em finite} set of corner points of $\scrS$. It is straightforward to see that the convex hull of the Minkowski addition of any finite set with $\rzeroo{n}$ is a convex polyhedron.
\end{proof}

\begin{example}
If $S = \{(0,0)\} \cup \{(-n, 1): n \geq 0\} \subseteq \zz^2$, then the convex hull of $S + \rzeroo{2}$ is the upper half-plane excluding the negative $x$-axis $\{(a, 0): a < 0\}$. In particular, \cref{corner-convex} may not hold if $\zzeroo{n}$ is replaced by $\zz^n$.
\end{example}

We say that a monomial order $\preceq$ on $\zzeroo{n}$ has \index{Finite!depth}{\em finite depth} if for every $\alpha \in \zzeroo{n}$, the set $\znzeropreceqalpha := \{\beta \in \zzeroo{n}: \beta \preceq \alpha\}$ is finite. 

\begin{example} \label{grlex}
The lexicographic order $\preceq_{lex}$ from \cref{lexicographic-example} does {\em not} have finite depth. The \index{Graded lexicographic order}{\em graded lexicographic order} $\preceq_{grlex}$ on $\zzeroo{n}$ is defined as follows: if $\alpha = (\alpha_1, \ldots, \alpha_n)$ and $\beta = (\beta_1, \ldots, \beta_n) \in \zzeroo{n}$, then $\alpha \preceq_{grlex} \beta$ if either $\sum_j \alpha_j < \sum_j \beta_j$, or if $\sum_j \alpha_j = \sum_j \beta_j$ and $\alpha \preceq_{lex} \beta$. It is straightforward to check that $\preceq_{grlex}$ is a monomial order on $\zzeroo{n}$ of finite depth.
\end{example}

\subsection{Monomial orders on rings of formal power series over a field} \label{monomial-power-section}
Following \cref{formal-power-section} we write $\hat R$ for the ring $\kk[[x_1, \ldots, x_n]]$ of formal power series in indeterminates $(x_1, \ldots, x_n)$ over a field $\kk$ and $\hat m$ for the (unique) maximal ring of $\hat R$. A \index{Monomial order!on the ring of formal power series}{\em monomial order} on $\hat R$ is simply a monomial order $\preceq$ on $\zzeroo{n}$, which induces an ordering on the set of monomials in $(x_1, \ldots, x_n)$ by the relation: $x^\alpha \preceq x^\beta$ if and only if $\alpha \preceq \beta$. In this section we use monomial orders to deduce some of the basic properties of $\hat R$.  Fix a monomial order $\preceq$ on $\hat R$, or equivalently, on $\zzeroo{n}$. For each nonempty subset $S$ of $\zzeroo{n}$, we write $\In_{\preceq}(S)$ for the minimal element of $S$ with respect to $\preceq$. For $f = \sum_\alpha c_\alpha x^\alpha \in \hat R$, the \index{Support!of a formal power series}{\em support} $\supp(f)$ of $f$ is the set $\{\alpha: c_\alpha \neq 0\} \subseteq \zzeroo{n}$; if $f \neq 0$ and $\alpha := \In_\preceq(\supp(f))$, we say that $\alpha$ is the \index{Initial!exponent with respect to a monomial order}{\em initial exponent} of $f$ denoted by $\exp_\preceq(f)$ and $c_\alpha x^\alpha$ is the \index{Initial!form!with respect to a monomial order}{\em initial form} $\In_\preceq(f)$ of $f$. For each subset $Q$ of $\hat R$, we write $\exp_\preceq(Q) := \{\exp_\preceq(f): f \in Q,\ f \neq 0\} \subseteq \znzero$ for the set of initial exponents of nonzero elements in $Q$.

\begin{example} \label{example:exp-example}
Consider $f_1 := x_1^2 + x_1^3 - x_2$ and $f_2 := x_1^3 + x_1x_2 + x_2^2 \in \kk[[x_1, x_2]]$. Then $\In_{\preceq_{lex}}(f_1) = x_2$ and $\In_{\preceq_{lex}}(f_2) = x_2^2$, so that $\exp_{\preceq_{lex}}(f_1) =  (0, 1)$ and $\exp_{\preceq_{lex}}(f_2) =  (0, 2)$. Substituting $x_2 = x_1^2 + x_1^3 - f_1$ in the expression of $f_2$ shows that $g :=2x_1^3 + x_1^4 + (x_1^2 + x_1^3)^2$ is in the ideal $I$ of $\kk[[x_1, x_2]]$ generated by $f_1, f_2$. If $\character(\kk) \neq 2$, then $\In_{\preceq_{lex}}(g) = 2x_1^3$, so that $\exp_{\preceq_{lex}}(g) =  (3, 0)$, and it is not hard to see that $\exp_{\preceq_{lex}}(I)$ is in fact the set of all elements on $\zzeroo{2}$ on and ``above'' the line joining $(0, 1)$ and $(3,0)$ (see \cref{fig:exp-example-lex}). 
\end{example}

\begin{center}
\def\figwidth{0.2}
\def\figsep{0.05}
\tikzstyle{dot} = [\colordot, circle, minimum size=4pt, inner sep = 0pt, fill]
\def\gridx{4.5}
\def\gridy{3.5}
\def\opazero{0.5}
\def\colordot{red}
\def\colorone{red}
\def\colorzero{blue}
\def\scale{0.5}
\begin{figure}[htb]
\begin{subfigure}[b]{\figwidth\textwidth}

\begin{tikzpicture}[scale=\scale]
\draw [gray,  line width=0pt] (-0.5,-0.5) grid (\gridx,\gridy);
\draw [<->] (0, \gridy) |- (\gridx, 0);
\node[dot, \colordot] at (2, 0) {};
\node[dot, \colordot] at (3, 0) {};
\node[dot, \colordot] at (0, 1) {};
\end{tikzpicture}

\caption{$\supp(f_1)$}
\label{fig:exp-example-f1}
\end{subfigure}\hspace{\figsep\textwidth}
\begin{subfigure}[b]{\figwidth\textwidth}

\begin{tikzpicture}[scale=\scale]
\draw [gray,  line width=0pt] (-0.5,-0.5) grid (\gridx,\gridy);
\draw [<->] (0, \gridy) |- (\gridx, 0);
\node[dot, \colordot] at (3, 0) {};
\node[dot, \colordot] at (1, 1) {};
\node[dot, \colordot] at (0, 2) {};
\end{tikzpicture}

\caption{$\supp(f_2)$}
\label{fig:exp-example-f2}
\end{subfigure}\hspace{\figsep\textwidth}
\begin{subfigure}[b]{\figwidth\textwidth}

\begin{tikzpicture}[scale=\scale]
\draw [gray,  line width=0pt] (-0.5,-0.5) grid (\gridx,\gridy);
\draw [<->] (0, \gridy) |- (\gridx, 0);
\draw[ultra thick, \colorone]  (3,0) --  (0,1);
\fill[\colorzero, opacity=\opazero ] (3,0) --  (0,1) -- (0,\gridy) -- (\gridx,\gridy) -- (\gridx,0) -- cycle;
\end{tikzpicture}

\caption{$\exp_{\preceq_{lex}}(\langle f_1, f_2 \rangle)$}
\label{fig:exp-example-lex}
\end{subfigure}\hspace{\figsep\textwidth}
\begin{subfigure}[b]{\figwidth\textwidth}

\begin{tikzpicture}[scale=\scale]
\draw [gray,  line width=0pt] (-0.5,-0.5) grid (\gridx,\gridy);
\draw [<->] (0, \gridy) |- (\gridx, 0);
\draw[ultra thick, \colorone]  (2,0) --  (0,2);
\fill[\colorzero, opacity=\opazero ] (2,0) --  (0,2) -- (0,\gridy) -- (\gridx,\gridy) -- (\gridx,0) -- cycle;
\end{tikzpicture}

\caption{$\exp_{\preceq_{rlex}}(\langle f_1, f_2 \rangle)$}
\label{fig:exp-example-rlex}
\end{subfigure}
\caption{Supports of $f_1 := x_1^2 + x_1^3 - x_2$ and $f_2 := x_1^3 + x_1x_2 + x_2^2$, and exponents of $\langle f_1, f_2  \rangle$ when $\character(\kk) \neq 2$}
\label{fig:exp-example}
\end{figure}
\end{center}

\begin{example} \label{example:exp-example-2}
To run the same computations as in \cref{example:exp-example} with $\preceq_{rlex}$, note that $\In_{\preceq_{rlex}}(f_1) = x_1^2$ and $\In_{\preceq_{rlex}}(f_2) = x_1^3$, so that  $\exp_{\preceq_{rlex}}(f_1) =  (2, 0)$ and $\exp_{\preceq_{rlex}}(f_2) =  (3, 0)$. Now $h_1 := f_2 - x_1(f_1 - f_2) = 2x_1x_2 + x_1^2x_2 + x_2^2 + x_1x_2^2$ and $h_2 := x_1h_1 - x_2(2f_1 - f_2) = 2x_2^2 + 2x_1x_2^2 + x_1^2x_2^2 + x_2^3$, so that if $\character(\kk) \neq 2$, then
\begin{alignat*}{4}
& \In_{\preceq_{rlex}}(h_1) &&= 2x_1x_2, \quad && \exp_{\preceq_{rlex}}(h_1) &&=  (1, 1) \\
& \In_{\preceq_{rlex}}(f_1) &&= 2x_2^2, \quad && \exp_{\preceq_{rlex}}(h_2) &&=  (0, 2) 
\end{alignat*}
It is then not hard to see that $\exp_{\preceq_{rlex}}(I)$ is the set of all elements on $\zzeroo{2}$ on and above the line joining $(0, 2)$ and $(2,0)$ (see \cref{fig:exp-example-rlex}). 
\end{example}

%
%
%

Using \cref{corner-lemma} it is straightforward to check that every ideal of $\hat R$ which is generated by monomials is finitely generated. We now prove the more general fact that every ideal of $\hat R$ is finitely generated. 

%
%
%


\begin{thm} \label{power-basis}
Let $\preceq$ be a monomial order of finite depth on $\hat R$. Let $I$ be an ideal of $\hat R$ and  $\scrC_I = \{\alpha_1, \ldots, \alpha_s\}$ be the set of corner points of $\exp_\preceq(I)$. For each $i = 1, \ldots, s$, pick $f_i\in I$ such that $\In_\preceq(f_i) = x^{\alpha_i}$.
\begin{enumerate}
\item Each $g \in \hat R$ can be expressed as $g = \sum_i f_i h_i + g'$ for some $h_1, \ldots, h_s, g'  \in \hat R$ such that either $g' = 0$, or $\exp_\preceq(g) \preceq \exp_\preceq(g')$ and $\exp_\preceq(g') \not \in \exp_\preceq(I)$.  
\item $I$ is generated by $f_1, \ldots, f_s$. 
\item If $g \in \hat R\setminus \{0\}$ is such that $\supp(g) \subset \zzeroo{n} \setminus\exp_\preceq(I)$, then $g \not\in I$. In particular, $\dim_\kk( \hat R/I) = | \zzeroo{n} \setminus \exp_\preceq(I)|$ and if $\dim_\kk( \hat R/I) < \infty$, then $\{x^\alpha: \alpha \in \zzeroo{n} \setminus \exp_\preceq(I)\}$ form a $\kk$-basis of $ \hat R/I$. 
\end{enumerate}
\end{thm}  

\begin{proof}
At first we prove the first assertion. Pick $g \in \hat R$. If $\alpha := \exp_\preceq(g) \not\in \exp_\preceq(I)$, there is nothing to do. Otherwise pick the smallest $i_1$, $1 \leq i_1 \leq s$, such that $\alpha = \alpha_{i_1} + \beta_1$ for some $\beta_1 \in \zzeroo{n}$. Then $\In_\preceq(g) = c_1x^{\alpha_{i_1} + \beta_1}$ for some $c_1 \in \kk$. Write $g_0 := g$ and $g_1 := g - c_1x^{\beta_1}f_{i_1}$. Continuing with $g_1$ and repeating this procedure, yields a sequence of elements $(g_k)_{k \geq 0}$ such that $\exp_\preceq(g_k) \precneqq \exp_\preceq(g_{k+1})$ for each $k$. Either this sequence is infinite, in which case the finite depth of $\preceq$ ensures that $f$ is a $\hat R$-linear combination of the $f_j$, or it stops at a stage $k$, in which case $\exp_\preceq(g_k) \not\in \exp_\preceq(I)$. This implies the first assertion. The second assertion follows from the first. The third assertion follows from the first assertion and finite depth of $\preceq$.
\end{proof}

\begin{cor}\label{noetherian-power-series}
$\hat R$ is Noetherian. \qed
\end{cor}

\begin{rem}
The requirement that $\preceq$ in \cref{power-basis} have finite depth is only a technical convenience which makes the proof shorter. All the assertions of \cref{power-basis} are true without this requirement. E.g.\ in \cref{example:exp-example,example:exp-example-2} we have applied two different choices of $\preceq$ (none of which is of finite depth) to the same ideal $I$, and have seen that $|\zzeroo{2} \setminus \exp_{\preceq}(I)|$ remains the same even though $\exp_{\preceq}(I)$ were different. In particular, $\dim_\kk(\kk[[x_1, x_2]]/I) = 3$, and taking the monomials with exponents in  $\zzeroo{2} \setminus \exp_{\preceq}(I)$ we get two bases of $\kk[[x_1, x_2]]/I$, namely $(1, x_1, x_1^2)$ and $(1, x_1, x_2)$. 
\end{rem}

Given formal power series $f_1, \ldots, f_s \in \hat R$, we now show that the dimension of $\hat R/\langle f_1, \ldots, f_s \rangle$ as a vector space over $\kk$ can be \index{Finite!determinacy}{\em finitely determined}, i.e.\ it can be determined by {\em polynomials} $g_1, \ldots, g_s$ provided the $g_j$ are ``sufficiently close'' to to $f_j$, where the ``closeness'' of elements of $R$ will be measured by monomial orders of finite depth. 
%
%
%
Let $\preceq$ be a monomial order on $\hat R$. For each $f = \sum_\alpha c_\alpha x^\alpha \in \hat R$ and each $\beta \in \zzeroo{n}$ write $[f]_{\preceq \beta} := \sum_{\alpha \preceq \beta} c_\alpha x^\alpha$.

\begin{thm} \label{finite-determinacy}
Assume $\preceq$ has finite depth. Let $I$ be an ideal of $\hat R$ generated by $f_1, \ldots, f_s$. For each $\beta \in \zzeroo{n}$, let $[I]_{\preceq \beta}$ be the ideal of $\hat R$ generated by $[f_j]_{\preceq \beta}$, $j = 1, \ldots, s$. 
\begin{enumerate}
\item \label{finite-determinacy:finite} If $\dim_\kk(\hat R/I) < \infty$, then there is $\beta \in \zzeroo{n}$ such that $\exp_\preceq(I) = \exp_{\preceq}([I]_{\preceq \beta'})$ and $\dim_\kk(\hat R/I) = \dim_\kk(\hat R/[I]_{\preceq \beta'})$ for all $\beta' \succeq \beta$. 
\item If $\dim_\kk(\hat R/I) = \infty$, then for each $N \geq 0$, there is $\beta \in \znzero$ such that $\dim_\kk(\hat R/[I]_{\preceq \beta'}) \geq N$ for all $\beta' \succeq \beta$.
\end{enumerate}
\end{thm}

\begin{proof}
At first assume $\dim_\kk(\hat R/I) < \infty$. Due to \cref{power-basis} to prove assertion \cref{finite-determinacy:finite} it suffices to show that $\exp_\preceq([I]_{\preceq {\beta'}}) = \exp_\preceq(I)$ if $\beta'$ is sufficiently ``high'' with respect to $\preceq$. \Cref{power-basis} implies that $\zzeroo{n} \setminus \exp_\preceq(I)$ is finite. For every finite subset $S$ of $\zzeroo{n}$, write $\ld_\preceq(S)$ for the maximal element of $S$ with respect to $\preceq$. Let $\beta_1 := \ld_\preceq (\zzeroo{n} \setminus \exp_\preceq(I))$. Then for each $\beta' \succeq \beta_1$ and each $g_1, \ldots, g_s \in \hat R$, we have 
\begin{align*}
\sum_j g_j[f_j]_{\preceq \beta'} 
	&= \sum_j g_j([f_j]_{\preceq \beta'} -f_j) + \sum_j g_jf_j
\end{align*}
Let $h := \sum_j g_j[f_j]_{\preceq \beta'}$, $h_1 := \sum_j g_j([f_j]_{\preceq \beta'} -f_j)$ and $h_2 := \sum_j g_jf_j$. Since $\exp_\preceq(h_2) \in \exp_\preceq(I)$ and $\exp_\preceq(h_1) \succneqq \ld_\preceq (\zzeroo{n} \setminus \exp_\preceq(I))$, it is straightforward to see that if $h = h_1 + h_2\neq 0$, then $\exp_\preceq(h ) \in \exp_\preceq(I)$. It follows that $\exp_\preceq ([I]_{\preceq \beta'}) \subseteq \exp_\preceq(I)$. Now let $\scrC_I$ be the set of corner points of $\exp_\preceq(I)$ and $\beta_2 := \ld_\preceq (\exp_\preceq(\scrC_I))$. Fix $\alpha \in \scrC_I$. There is $g_1, \ldots, g_s \in \hat R$ such that $\exp(\sum_j g_jf_j) = \alpha$. Pick $\beta' \succeq \beta_2$, and define $h, h_1,h_2$ as above. Then $\exp_\preceq (h_2) \precneqq \exp_\preceq(h_1)$, so that $\exp_\preceq(h) = \exp_\preceq(h_2) = \alpha$. Therefore $\scrC_I\subset \exp_\preceq ([I]_{\preceq \beta'})$ and consequently \cref{corner-lemma} implies that $\exp_\preceq (I) \subset \exp_\preceq ([I]_{\preceq \beta'})$. It follows that $\exp_\preceq (I) = \exp_\preceq ([I]_{\preceq \beta'})$ whenever $\beta' \succeq \beta := \ld_\preceq\{\beta_1, \beta_2\}$, which proves the first assertion. For the second assertion, fix $N \geq 0$. Take an arbitrary finite subset $S$ of $\znzero \setminus \exp_{\preceq}(I)$ such that $|S| \geq N$ and let $\beta := \ld_\preceq(S)$. The same argument as in the beginning of the proof suggests that for each $\beta' \succeq \beta$ and each $g_1, \ldots, g_s \in \hat R$, $\exp_\preceq(\sum_j g_j[f_j]_{\preceq \beta'}) \not\in S$, so that \cref{power-basis} implies that $\dim_\kk(\hat R/[I]_{\preceq \beta'}) \geq |S| \geq N$, as required. 
\end{proof}

\section{Primitive elements of $\zz^n$}
The results we prove in this section follow almost directly from the fundamental theorem of finitely generated modules over a PID (\cref{fundamentally-over-PID}); here we give more elementary proofs. An element in $\zz^n$ is {\em primitive} if it is nonzero and the greatest common divisor of its nonzero coordinates is one. Every member of a basis of $\zz^n$ is primitive. The first result below shows that the converse is also true. In this section we use $\langle \cdot, \cdot \rangle$ to denote the standard coupling of elements of $\zz^n$ and $\znstar$. 

\begin{lemma} \label{basis-lemma}
Let $n$ be a positive integer and $G$ be a subgroup of $\zz^n$.
\begin{enumerate}
\item If $\alpha$ is a primitive  element in $\zz^n$, then there is a basis of $\zz^n$ containing $\alpha$.
\item $G \cong \zz^m$ for some $m \leq n$, and there is a basis $(\alpha_1, \ldots, \alpha_n)$ of $\zz^n$ and positive integers $k_1, \ldots, k_m$ such that $(k_1 \alpha_1, \ldots, k_m \alpha_m)$ is a basis of $G$.
\end{enumerate}
\end{lemma}

\begin{proof}
For the first assertion, let $\alpha_1 :=\alpha$. Since $\alpha_1$ is primitive, there is $\beta_1 \in \znstar$ such that $\langle \beta_1, \alpha_1 \rangle = 1$. Let $H_1 := \beta_1^\perp := \{ \gamma \in \zz^n:  \langle \beta_1, \gamma \rangle= 0\}$. If $H_1 \neq 0$, pick a primitive element $\alpha_2 \in H_1$ and $\beta_2 \in \znstar$ such that $\langle \beta_2, \alpha_2 \rangle = 1$, and set $H_2 := H_1 \cap \beta_2^\perp$. Continue in this way up to the $n$-th step. It is straightforward to see that $\beta_1, \ldots, \beta_n$ are linearly independent (over $\rr$), so that $H_n = 0$. 

\begin{proclaim}
$(\alpha_1, \ldots, \alpha_n)$ is a basis of $\zz^n$. 
\end{proclaim}

\begin{proof}
Indeed, given any $\delta \in \zz^n$, let $d_1 := \langle \beta_1, \delta \rangle$. Then $\delta_1 := \delta - d_1\alpha_1 \in H_1$. Let $d_2 :=   \langle \beta_2, \delta_1 \rangle$. Then $\delta_2 := \delta_1 - d_2\alpha_2 \in H_2$. In this way we get that $\delta_n = \delta - \sum_j d_j \alpha_j \in H_n = 0$, as required.
\end{proof}

The above claim proves the first asertion. For the second assertion, we may assume $G \neq \{0\}$. For each nonzero $\alpha \in G$, let $d_\alpha$ be the greatest common divisor of the nonzero coordinates of $\alpha$. Pick $\alpha \in G$ with the smallest possible $d_\alpha$. Due to the first assertion we may assume \woutlog\ that $\alpha = d_\alpha(1, 0, \ldots, 0)$. The minimality of $d_\alpha$ implies that for all $\beta = (\beta_1, \ldots, \beta_n) \in G$, the first coordinate $\beta_1$ of $\beta$ is divisible by $d_\alpha$. It follows that $G = \zz \alpha + G'$, where $G' := G \cap (\{0\} \times \zz^{n-1})$. The second assertion holds for $G'$ by induction on $n$, which in turn implies that it holds for $G$.
\end{proof}

\begin{cor} \label{change-of-basis-lemma}
Let $\phi: \zz^n \to \zz^m$ be a homomorphism of abelian groups and $r$ be the rank (over $\qq$) of the matrix of $\phi$. Then the matrix of $\phi$ with respect to appropriate bases of $\zz^n$ and $\zz^m$ is of the form
\begin{align*}
\mat{\phi}
	&=
		\left[
			\begin{array}{c|c}
			D & 0\\
			\hline
			0 & 0
			\end{array}
		\right]
\end{align*}
where $D$ is an $r\times r$-diagonal matrix whose diagonal entries are positive integers.
\end{cor}

\begin{proof}
Let $e^m_i$ and $e^n_j$ denote respectively the $i$-th standard unit element in $\zz^m$ and the $j$-th standard unit element in $\zz^n$. Let $\alpha_1, \ldots, \alpha_m$ be the rows of the matrix of $\phi$ with respect to the standard bases of $\zz^n$ and $\zz^m$. Let $G$ be the subgroup of $\zz^n$ generated by $\alpha_1, \ldots, \alpha_m$. \Cref{basis-lemma} implies that after a change of basis of $\zz^m$ we may assume that $G$ is generated by $k_1e^n_1, \ldots, k_re^n_r$, where each $k_j$ is a positive integer. Pick $\gamma_1 = (\gamma_{1,1}, \ldots, \gamma_{1,m}) \in \zz^m$ such that $\sum_j \gamma_{1,j}\alpha_j = k_1e^n_1$. Then $\gamma_1$ must be a primitive element of $\zz^m$, so after a change of basis in $\zz^m$, we may assume that $\gamma_1 = e^m_1$, which in turn implies that $\alpha_1 = k_1e^n_1$. Note that for each $j = 2, \ldots, m$, $\alpha_j $ is of the form $d_jk_1e^n_1 + \alpha'_j$, with $\alpha'_j$ in the subgroup $G'$ of $\zz^n$ generated by $k_2e^n_2, \ldots, k_re^n_r$. Therefore, after a change of basis of $\zz^m$ the form $e^m_j \mapsto e^m_j +  d_jk_1e^m_1$, $j = 2, \ldots, m$, the matrix of $\phi$ is of the form
\begin{align*}
		\left[
			\begin{array}{c|c}
			k_1 & 0\\
			\hline
			0 & M
			\end{array}
		\right]
\end{align*}
for some $(m-1) \times (n-1)$ matrix $M$. Now apply induction (say, on $n$) to the homomorphism from $\zz^{n-1} \to \zz^{m-1}$ induced by $M$.
\end{proof}

\section{Symmetric multiadditive functions on a commutative semigroup}
Throughout this section $(\scrK, +)$ is a commutative semigroup and $n$ is a positive integer. Let $\rho$ be a mapping from $\scrK^n$ to $\rr$. We say that $\rho$ is \index{Symmetric function}{\em symmetric} if $\rho(f_1, \ldots, f_n) = \rho(f_{\sigma_1}, \ldots, f_{\sigma_n})$ for each permutation $\sigma$ of $(1, 2, \ldots, n)$. We say that $\rho$ is \index{Multiadditive function}{\em multiadditive} if
\begin{align}
\begin{split}
\rho(h_1, \ldots, h_{j-1}, qf + rg, h_{j+1} \ldots, h_n) 
	&= q\rho(h_1, \ldots, h_{j-1}, f, h_{j+1} \ldots, h_n) \\
	&\qquad + r \rho(h_1, \ldots, h_{j-1}, g, h_{j+1} \ldots, h_n)  
\end{split}
\label{multiadditivity}
\end{align}
for each $q, r \in \zzero$, $j \in \{1, \ldots, n\}$ and $f,g, h_1, \ldots, h_{j-1}, h_{j+1}, \ldots, h_n \in \scrK$. Throughout this section we use $\rho$ to denote a symmetric multiadditive function from $\scrK^n$ to $\rr$. 

\subsection{Existence from polynomial functions} 
For nonnegative integers $j_1, \ldots, j_k, n$ such that $n = j_1 + \cdots + j_k$, we write $\binom{n}{j_1, \ldots, j_k}$ for the {\em multinomial coefficient} $n!/(j_1! \cdots j_k!)$. Recall that $\binom{n}{j_1, \ldots, j_k}$ is precisely the coefficient of $x_1^{j_1} \cdots x_k^{j_k}$ in $(x_1 + \cdots + x_k)^n$.

\begin{lemma} \label{mixed-lemma}
Let $\nu: \scrK \to \rr$ be a function which satisfies the following property: for each $s \geq 1$ and $f_1, \ldots, f_s \in \scrK$, there are $\nu_\alpha(f_1, \ldots, f_s) \in \rr$ for all $\alpha \in \scrE_s:= \{(\alpha_1, \ldots, \alpha_s) \in (\zzeroo{s}: \alpha_1 + \cdots  + \alpha_s = n\}$ such that for all $\lambda_1, \ldots, \lambda_s \in \zzero$,
\begin{align}
\nu(\lambda_1f_1 + \cdots + \lambda_s f_s) = \sum_{\alpha \in \scrE_s} \nu_\alpha(f_1, \ldots, f_s) \lambda_1^{\alpha_1} \cdots \lambda_s^{\alpha_s}
\label{rho-pol}
\end{align}
Then $\rho(f_1, \ldots, f_n) := \frac{1}{n!}\nu_{(1, \ldots, 1)}(f_1, \ldots, f_n)$ is a symmetric multiadditive function from $\scrK^n$ to $\rr$ such that $\rho(f, \ldots, f) = \nu(f)$. 
\end{lemma}

\begin{proof}
Fix $f\in \scrK$. Applying \eqref{rho-pol} with $s = 1$ shows that $\nu(\lambda f) = \nu_n(f)\lambda^n$ for each $\lambda \in \zzero$.  Setting $\lambda = 1$, we have $\nu_n(f) = \nu(f)$, and therefore $\nu(\lambda f) = \lambda^n \nu(f)$ for all $\lambda \in \zzero$. It follows that 
\begin{align*}
\nu(\lambda_1f+ \cdots + \lambda_n f) 
	&=  (\lambda_1 + \cdots + \lambda_n)^n \nu(f) 
	= \sum_\alpha \binom{n}{\alpha_1, \ldots, \alpha_n} \lambda_1^{\alpha_1} \cdots \lambda_n^{\alpha_n} \nu(f)
\end{align*}
equating the coefficients of $\lambda_1 \cdots \lambda_n$ of the middle and the rightmost polynomial in $(\lambda_1, \ldots, \lambda_n)$ yields that
\begin{align*}
\rho(f, \ldots, f) = \frac{1}{n!}\binom{n}{1, \ldots, 1} \nu(f) =  \nu(f)
\end{align*}
It is clear that $\rho$ is symmetric in its arguments. For multiadditivity, write elements of $\zz^{n+1}$ as $(\alpha, \beta, \ggamma) := (\alpha, \beta, \gamma_1, \ldots, \gamma_{n-1})$ and note that 
\begin{align*}
\nu(\lambda f + \mu g + \tau_1h_1 + \cdots + \tau_{n-1} h_n) 
	&= \sum_{(\alpha, \beta, \ggamma) \in \scrE_{n+1}} \nu_{(\alpha, \beta, \ggamma)} (f, g, h_1, \ldots, h_{n-1}) \lambda^\alpha \mu^\beta\tau_1^{\gamma_1} \cdots \tau_{n-1}^{\gamma_{n-1}} 
\end{align*}
so that 
\begin{align*}
\nu(\lambda (f + g) + \tau_1h_1 + \cdots + \tau_{n-1} h_n) 
	&= \sum_{(\alpha, \beta, \ggamma) \in \scrE_{n+1}} \nu_{(\alpha, \beta, \ggamma)} (f, g, h_1, \ldots, h_{n-1}) \lambda^{\alpha +\beta} \tau_1^{\gamma_1} \cdots \tau_{n-1}^{\gamma_{n-1}} 
\end{align*}
It follows that 
\begin{align*}
\nu_{(\delta, \ggamma)}(f+g, h_1, \ldots, h_{n-1})
	&= \sum_{\alpha + \beta = \delta}\nu_{(\alpha, \beta, \ggamma)}(f, g, h_1, \ldots, h_{n-1})
\end{align*}
In particular, 
\begin{align*}
\rho(f + g, h_1, \ldots, h_{n-1})
	&= \frac{1}{n!}\nu_{(1, \ldots, 1)}(f + g, h_1, \ldots, h_{n-1}) \\
	&= \frac{1}{n!}(\nu_{(1, 0, 1, \ldots, 1)}(f, g, h_1, \ldots, h_{n-1})
		+ \nu_{(0, 1, \ldots, 1)}(f , g, h_1, \ldots, h_{n-1})) \\
	&= \rho(f, h_1, \ldots, h_{n-1}) + \rho(g, h_1, \ldots, h_{n-1})
\end{align*}
which implies that $\rho$ is multiadditive, and completes the proof.
\end{proof}

\subsection{Identities from homogeneous polynomials}
The properties of symmetric multiadditive functions imply that its summands can be ``expanded like polynomials,'' i.e.\ given $f_1, \ldots, f_N \in \scrK$ and $\lambda_{ij} \in \zzero$, 
\begin{align*}
 \rho(\sum_{j=1}^N \lambda_{1j}f_j, \ldots, \sum_{j=1}^N \lambda_{nj}f_j) 
 	&= \sum_{1 \leq j_1, \ldots, j_n \leq N} ( \prod_{i=1}^n \lambda_{i, j_i}) \rho(f_{j_1}, \ldots, f_{j_n})
 \end{align*}
This is an analogue of the polynomial identity 
\begin{align*}
\prod_{i=1}^n (\sum_{j=1}^N \lambda_{ij}x_j)
 	&=  \sum_{1 \leq j_1, \ldots, j_n \leq N}  \prod_{i=1}^n \lambda_{i,j_i} \prod_{i=1}^n x_{j_i}
 \end{align*}
 where $x_1, \ldots, x_N$ are indeterminates. This observation immediately implies the principle that symmetric multiadditive functions ``respect'' polynomial identities:

\begin{prop} \label{prop:polynomial-identity}
Let $x_1, \ldots, x_N$ be indeterminates and $f_1, \ldots, f_N \in \scrK$. Assume the following identiy holds in $\rr[x_1, \ldots, x_N]$: 
\begin{align*}
\sum_i r_i \prod_{j = 1}^N (\sum_{k = 1}^N \lambda_{ijk} x_k)^{\alpha_{ij}} = 0
\end{align*}
where $r_i \in \rr$ and $\lambda_{ijk}, \alpha_{ij} \in \zzero$, and $\sum_{j = 1}^N \alpha_{ij} = n$. Then 
\begin{align*}
\sum_i r_i \rho(\sum_{k = 1}^N \lambda_{i1k} f_k, \ldots, \sum_{k = 1}^N \lambda_{i1k} f_k, \ldots, 
		\sum_{k = 1}^N \lambda_{iNk} f_k, \ldots, \sum_{k = 1}^N \lambda_{iNk} f_k) = 0
\end{align*}
where $\sum_{k = 1}^N \lambda_{ijk} f_k$ is repeated $\alpha_{ij}$-times for each $i,j$. \qed
\end{prop}

We now use this observation to show that a symmetric multiadditive function is uniquely determined by its diagonal part. 

\begin{lemma} \label{r-s-I-lemma}
Let $x_1, \ldots, x_n$ be indeterminates and $I \subseteq [n]$. Let $I := \{i_1, \ldots, i_k\}$, where $k = |I|$. Write 
\begin{align*}
s_I & := (x_{i_1} + \cdots + x_{i_k})^n \\
r_I &:= \sum_{\substack{j_1 + \cdots + j_k = n \\j_l \geq 1,\ l = 1, \ldots, k}} \binom{n}{j_1, \ldots, j_k} x_{i_1}^{j_1} \cdots x_{i_k}^{j_k}
\end{align*}
Then 
\begin{align}
r_I 
	&= s_I -	\sum_{\substack{J \subset I \\ |J| = k -1}}s_J 		
			+ \sum_{\substack{J \subset I \\ |J| = k -2}}s_J + \cdots 
						+ (-1)^{k-1} \sum_{\substack{J \subset I \\ |J| = 1}}s_J \label{r-s-I}
\end{align}
In particular, 
\begin{align}
n!x_1 \cdots x_n 
	&= \sum_{\substack{I \subseteq [n]\\ I \neq \emptyset}}(-1)^{n - |I|} (\sum_{i \in I} x_i )^n\ \label{r-s-[n]}
\end{align}
\end{lemma}

\begin{proof}
It suffices to prove identity \eqref{r-s-I}, since identity \eqref{r-s-[n]} follows from \eqref{r-s-I} by setting $I = [n]$. Straightforward algebra shows that
\begin{align*}
s_I &= \sum_{j_1 + \cdots + j_k = n } \binom{n}{j_1, \ldots, j_k} x_{i_1}^{j_1} \cdots x_{i_k}^{j_k} 
	= r_I + \sum_{J \subsetneq I}r_J 
	= r_I + \sum_{\substack{J \subsetneq I \\ |J| = k -1}}r_J + \sum_{\substack{J \subsetneq I \\ |J| \leq  k -2}}r_J \\
	&= r_I + \sum_{\substack{J \subsetneq I \\ |J| = k -1}}(s_J -  \sum_{J' \subsetneq J} r_{J'}) + \sum_{\substack{J \subsetneq I \\ |J| \leq  k -2}}r_J \\
	&= r_I + 	\sum_{\substack{J \subsetneq I \\ |J| = k -1}}s_J + \sum_{\substack{J \subsetneq I \\ |J| \leq  k -2}}(1 - k + |J|)r_J
	= r_I + 	\sum_{\substack{J \subsetneq I \\ |J| = k -1}}s_J + \sum_{l=1}^{k-2}\sum_{\substack{J \subset I \\ |J| =  l}}(1 - k + l)r_J \\
	&= r_I + 	\sum_{\substack{J \subsetneq I \\ |J| = k -1}}s_J 
		 - \sum_{\substack{J \subset I \\ |J| = k -2}}(s_J -  \sum_{J' \subsetneq J} r_{J'}) + \sum_{l=1}^{k-3}\sum_{\substack{J \subset I \\ |J| =  l}}(1 - k + l)r_J \\
	&= r_I + 	\sum_{\substack{J \subsetneq I \\ |J| = k -1}}s_J 
			 - \sum_{\substack{J \subset I \\ |J| = k -2}}s_J 
			 + \sum_{l=1}^{k-3}\sum_{\substack{J \subset I \\ |J| =  l}}(1 - (k - l) + \binom{k-l}{2})r_J 
\end{align*}
In this way, at every step writing $r_J = s_J - \sum_{J' \subsetneq J} r_{J'}$ and rearranging terms, and observing that $r_J = s_J$ whenever $|J| = 1$, we will have
\begin{align*}
s_I
	&= r_I + 	\sum_{\substack{J \subsetneq I \\ |J| = k -1}}s_J 
				 - \sum_{\substack{J \subset I \\ |J| = k -2}}s_J + \cdots 
				 + (-1)^{k-3}\sum_{\substack{J \subset I \\ |J| = 2}}s_J \\
	& \qquad 			 
				 + \sum_{\substack{J \subset I \\ |J| =  1}}(1 - (k - 1) + \binom{k-1}{2} + \cdots + (-1)^{k-2} \binom{k-1}{k-2})s_J 
\end{align*}
Since $\sum_{j = 0}^{k-2}(-1)^j \binom{k-1}{j} = (1-1)^{k-1} - (-1)^{k-1} = -(-1)^{k-1}$, the lemma follows. 
\end{proof}

\begin{cor}[See e.g.\ \cite{mathoverflow-tom-leinster}] \label{tom-lemma}
$\rho$ is uniquely determined by its diagonal part, i.e.\ the map from $\scrK$ to $\rr$ which sends $f \mapsto \rho(f, \ldots, f)$. More precisely,   
\begin{align}
\rho(f_1, \ldots, f_n)
	&= \frac{1}{n!}\sum_{\substack{I \subseteq [n]\\ I \neq \emptyset}}(-1)^{n - |I|}\rho  (\sum_{i \in I} f_i, \ldots, \sum_{i \in I} f_i )
\end{align}
\end{cor}

\begin{proof}
Combine \cref{prop:polynomial-identity} and identity \eqref{r-s-[n]}. 
\end{proof}

\begin{lemma}
Let $x_1, \ldots, x_k, y_1, \ldots, y_k$ be indeterminates. Then
\begin{align}
x_1 \cdots x_k 
	&= \sum_{I \subseteq [k]} (-1)^{k- |I|} \prod_{i \in I}(x_i + y_i) \prod_{j \in [k]\setminus I} y_j \label{kushnirenko-pol-identity}
\end{align}
\end{lemma}

\begin{proof}
We prove this by induction on $k$. For $k = 1$ it boils down to the identity $x_1 = (x_1 +y_1) - y_1$. In the general case, write
\begin{align*}
x_1 \cdots x_k 
	&= \prod_{i=1}^k (x_i + y_i) - \sum_{j = 1}^k (\prod_{i=1}^{j-1}x_i)y_j (\prod_{i=j+1}^k(x_i + y_i))
\end{align*}
Applying induction to the expression $\prod_{i=1}^{j-1}x_i$ yields that 
\begin{align*}
x_1 \cdots x_k 
	&= \prod_{i=1}^k (x_i + y_i) - \sum_{j = 1}^k 
		(\sum_{I \subseteq [j-1]} (-1)^{j-1- |I|} \prod_{i \in I}(x_i + y_i) \prod_{i' \in [j-1]\setminus I} y_{i'})
		y_j (\prod_{i=j+1}^k(x_i + y_i)) \\
	&= \prod_{i=1}^k (x_i + y_i) - \sum_{j = 1}^k 
			\sum_{\substack{[k] \supseteq I \supseteq [k]\setminus [j]   \\ I \not\ni j}}
			(-1)^{k - |I| -1} \prod_{i \in I}(x_i + y_i) \prod_{i' \in [k]\setminus I} y_{i'}\\
	&= 	\sum_{I \subseteq [k]}
				(-1)^{k - |I|} \prod_{i \in I}(x_i + y_i) \prod_{i' \in [k]\setminus I} y_{i'}
\end{align*}
This completes the proof.
\end{proof}

For the next result we assume $\rho: \scrK^n \to \rr$ is a symmetric multiadditive ``rational'' map, i.e.\ $\rho$ may not be defined everywhere on $\scrK^n$, but if $\rho(f_1, \ldots, f_n)$ is defined, then $\rho(f_{\sigma_1}, \ldots, f_{\sigma_n})$ are defined for all permutations $\sigma$ of $(1, \ldots, n)$, and all of them take the same value in $\rr$; and identity \eqref{multiadditivity} holds whenever $\rho$ is defined on at least two of the three elements of $\scrK^n$ that appear on \eqref{multiadditivity}. Given $h_1, \ldots, h_n \in \scrK$ and $I, J \subseteq [n]$ such that $|I| + |J| = n$, we write $\rho((h_i)_{i \in I}, (h_j)_{j \in J})$ for $\rho(h_{i_1}, \ldots, h_{i_k}, h_{j_1}, \ldots, h_{j_{n-k}})$ (provided it is defined), where $I = \{i_1, \ldots, i_k\}$ and $J = \{j_1, \ldots, j_{n-k}\}$. The following is an immediate consequence of the $k = n$ case of identity \eqref{kushnirenko-pol-identity} and \cref{prop:polynomial-identity}. 

\begin{cor} \label{kushnirenko-lemma}
Let $f_1, \ldots, f_n, g_1, \ldots, g_n \in \scrK$ be such that $\rho((f_i)_{i \in I}, (g_j)_{j \in [n] \setminus I})$ is defined for each $I \subseteq [n]$. Then 
\begin{align}
\rho(f_1, \ldots, f_n)
	&= \sum_{I \subseteq [n]} (-1)^{n- |I|} \rho((f_i + g_i)_{i \in I}, (g_j)_{j \in [n] \setminus I})  \label{kushnirenko-semigroup-identity} 
	\qed
\end{align}
\end{cor}
%
%

\chapter{Some results related to schemes} \label{proof-scheme-section}
\section{Macaulay's Unmixedness Theorem}
In this section we prove the {\em unmixedness theorem} of F.\ S.\ Macaulay for polynomial rings $R := k[x_1, \ldots, x_n]$ over an arbitrary field $k$ following \cite[Section VII.8]{zsII}. This theorem is a crucial ingredient of intersection theory on nonsingular varieties, and was used in a fundamental way (via \cref{cor:Macaulay}) in many places in this book. The statement of this theorem requires the notion of \index{Dimension!of an ideal}{\em dimension}\footnote{It is in fact the same as the dimension of the algebraic variety determined by $\qqq$ as defined in \cref{dimension-section}.} $\dim(\qqq)$ of an ideal $\qqq$ of $R$. If $\qqq$ is prime, then $\dim(\qqq)$ is the transcendence degree (of the field of fractions) of $R/\qqq$ over $k$. In general $\dim(\qqq)$ is the maximum of the dimensions of the prime ideals {\em associated with} $\qqq$. Since the minimal elements of the set of prime ideals of $R$ containing $\qqq$ are precisely the minimal elements of the set of ideals associated with $\qqq$ (\cref{prop:minimal-prime}), it follows that $\dim(\qqq)$ is the maximum of the dimensions of the {\em minimal} prime ideals containing $\qqq$. One of the most fundamental properties of dimensions is given by Krull's principal ideal theorem (\cref{thm:principal-ideal}) - the following result is its straightforward consequence. 

\begin{prop} \label{prop:Macaulay-dim} 
Let $k$ be an arbitrary field, and $\qqq$ be a {\em proper} ideal of $k[x_1, \ldots, x_n]$ generated by $f_1, \ldots, f_m$, $m \leq n$. Then $\dim(\qqq) \geq n - m$. If in addition $\dim(\qqq) = n - m$, then $\dim(\langle f_1, \ldots, f_j \rangle) = n - j$, for each $j = 1, \ldots, m$. \qed
\end{prop}

The following is a special case of Macaulay's unmixedness theorem (\cref{thm:Macaulay}). 

\begin{prop} \label{Macaulay:1}
Let $g = g_1^{q_1}\cdots g_r^{q_r} \in R := k[x_1, \ldots, x_n]$, where $g_j$ are irreducible (non-constant) polynomials and $q_j$ are positive integers. Then $h \in R$ is a zero-divisor in $R/\langle g \rangle$ if and only if $h$ is divisible by some $g_j$.  
\end{prop}

\begin{proof}
This immediately follows from the fact that $R$ is a unique factorization domain.
\end{proof}

%

\begin{thm}[%
Macaulay's unmixedness\footnote{An ideal is {\em unmixed} if all its associated prime ideals have the same dimension. \Cref{thm:Macaulay} is called the ``unmixedness theorem'' due to this equivalent formulation: ``Let $\qqq$ be an ideal of $k[x_1, \ldots, x_n]$ generated by $m$ polynomials. If $\dim(\qqq) = n - m$, then $\qqq$ is unmixed.'' The equivalence of this statement with \cref{thm:Macaulay} follows in straightforward manner from \cref{prop:zero-union-minimally-prime} and Krull's principal ideal theorem (\cref{thm:principal-ideal}).} theorem {\cite[Theorem VII.26]{zsII}
}] \label{thm:Macaulay} 
\index{Macaulay's unmixedness theorem}
\index{Unmixedness theorem}
Let $k$ be an arbitrary field, and $f_1, \ldots, f_m \in k[x_1, \ldots, x_n]$, $m \leq n$, be such that $\dim(\langle f_1, \ldots, f_m\rangle) = n-m$. Then for each $j = 1, \ldots, m$, $f_j$ is a non zero-divisor in $k[x_1, \ldots, x_n]/\langle f_1, \ldots, f_{j-1} \rangle$. 
\end{thm}

\begin{proof}
Write $R := k[x_1, \ldots, x_n]$. We proceed by induction on $m$. The case of $m = 1$ is obvious, and the case of $m = 2$ follows from \cref{Macaulay:1}. So assume $m \geq 3$. If $f_m$ is a zero-divisor in $R/\langle f_1, \ldots, f_{m-1} \rangle$, then it is contained in a prime ideal $\ppp$ associated with $\qqq := \langle f_1, \ldots, f_{m-1} \rangle$. 

\begin{proclaim} \label{Macaulay:claim-p}
$\dim(\ppp) < n - m + 1$. 
\end{proclaim}
\begin{proof}
By assumption $\ppp$ is {\em not} a minimal prime ideal containing $\qqq$, i.e.\ $\ppp \supsetneq \ppp' \supseteq \qqq$ for some prime ideal $\ppp'$ of $R$ with $\dim(\ppp') = n - m + 1$. The claim follows from applying Krull's principal ideal theorem (\cref{thm:principal-ideal}) to $R/\ppp'$ (or from \cref{thm:dim-smaller}). 
\end{proof}

Let $d := \dim(\ppp)$. Then there are $x_1, \ldots, x_d$ which are algebraically independent over $k$ in $R/\ppp$. This means $\kk[x_1, \ldots, x_d] \cap \ppp = 0$. Let $\tilde R := k(x_1, \ldots, x_d)[x_{d+1}, \ldots, x_n]$. Given an ideal $\rrr$ of $R$, we write $\tilde \rrr$ for the ideal of $\tilde R$ generated by $\rrr$. 

\begin{proclaim} \label{Macaulay:claim-p-2}
Let $\rrr$ be a prime ideal of $R$. 
\begin{enumerate}
\item If $\rrr \cap k[x_1, \ldots, x_d] = \{0\}$, then $\dim(\tilde \rrr) = \dim(\rrr) - d$. 
\item The prime ideals associated to $\tilde \qqq$ are precisely the ideals $\tilde \rrr$ corresponding to those prime ideals $\rrr$ associated with $\qqq$ such that $\rrr \cap k[x_1, \ldots, x_d] = \{0\}$. 
\end{enumerate}
\end{proclaim}

\begin{proof}
The first assertion follows immediately from choosing a transcendence basis (over $k$) of $R/\rrr$ containing $x_1, \ldots, x_d$. The second assertion follows from \cref{example:local-primary}.
\end{proof}

Since $\dim(\qqq) = n - (m -1)$ (\cref{prop:Macaulay-dim} ), the first assertion of \cref{Macaulay:claim-p-2} implies that $\tilde \qqq$ is a proper ideal of $\tilde R$ with dimension $\leq n - d - (m-1)$. Since $\tilde \qqq$ is generated by $m - 1$ elements, \cref{prop:Macaulay-dim} then implies that $\dim(\tilde \qqq)$ is precisely $n - d - (m - 1)$. On the other hand, applying the first assertion of \cref{Macaulay:claim-p-2} with $\rrr = \ppp$ shows that $\tilde \ppp$ is a zero dimensional ideal of $\tilde R$ which is also associated with $\tilde \qqq$. But then \cref{Macaulay:maximal-case} below gives a contradiction and proves \cref{thm:Macaulay}.

\begin{proclaim} \label{Macaulay:maximal-case}
Let $\bbb \subset \aaa$ be ideals in $R := k[x_1, \ldots, x_r]$ such that 
\begin{enumerate}
\item $\aaa$ is prime, 
\item $\dim(\aaa) = 0$, 
\item $\bbb$ is generated by $m - 1$ elements, and 
\item $\dim(\bbb) = r - (m - 1)$.
\end{enumerate}
If $r \geq m$, then $\aaa$ is {\em not} associated with $\bbb$. 
\end{proclaim}

To complete the proof of \cref{thm:Macaulay} it remains to prove \cref{Macaulay:maximal-case}. Pick a set $g_1, \ldots, g_{m-1}$ of generators of $\bbb$. Let $\ccc$ be the ideal of $R$ generated by $g_1, \ldots, g_{m-2}$ (recall that $m \geq 3$). \Cref{prop:Macaulay-dim} implies that $\dim(\ccc) = r - (m - 2)$. If $h \in R$ is not in any of the minimal prime ideal associated with $\ccc$, the induction hypothesis then applies to $g_1, \ldots, g_{m-2}, h$ and implies that $h$ is {\em not} a zero-divisor in $R/\ccc$. \Cref{thm:primary-decomposition,prop:Macaulay-dim} then imply that all prime ideals associated to $\ccc$ have dimension $r - m + 2$; enumerate these ideals as $\ppp_1, \ldots, \ppp_s$. Let $\ppp_{s+1}, \ldots, \ppp_{s'}$ be the {\em minimal} (i.e.\ $(r - m + 1)$-dimensional) prime ideals associated to $\bbb$. Since $r \geq m$, each $\ppp_j$ has dimension greater than zero. We will use the following result to construct a special type of polynomial which does not belong to any $\ppp_j$. 

\begin{prolemma} \label{Macaulay:transcendental-choice}
Given finitely many positive dimensional prime ideals $\ppp_1, \ldots, \ppp_{s'}$ in $R := k[x_1, \ldots, x_r]$, there is $e$, $1 \leq e \leq s'$, and a polynomial $h(x_1, \ldots, x_{e-1})$ such that the image of $y_e := x_e + h(x_1, \ldots, x_{e-1})$ is transcendental over $k$ in $R/\ppp_j$ for each $j$. 
\end{prolemma}

\begin{proof}
Pick the smallest integer $i_1$ such that $x_{i_1}$ is transcendental over $k$ in $R/\ppp_j$ for some $j$. Reorder the $\ppp_j$ in a way that $x_{i_1}$ is transcendental over $k$ in $R/\ppp_1, \ldots, R/\ppp_{s_1}$ and algebraic over $k$ in $R/\ppp_j$ for $j > s_1$. If $s_1 < s'$, then pick the smallest integer $i_2 > i_1$ such that $x_{i_2}$ is transcendental over $k$ in $R/\ppp_j$ for some $j > s_1$. Then reorder $\ppp_{s_1 + 1}, \ldots, \ppp_{s'}$ in a way that $x_{i_2}$ is transcendental over $k$ in $R/\ppp_{s_1+1}, \ldots, R/\ppp_{s_2}$ and algebraic over $k$ in $R/\ppp_j$ for $j > s_2$. Since each $\ppp_j$ is positive dimensional, we can continue in this way until there is an integer $t$ such that $s_t = s'$. In particular, 
\begin{defnlist} 
\item $x_{i_t}$ is transcendental over $k$ in $R/\ppp_j$ for each $j = s_{t-1} + 1, \ldots, s_t = s'$.
\end{defnlist}
Now fix $j$, $s_{t-2} + 1 \leq j \leq s_{t-1}$. We claim that there is at most one integer $q$ such that $x_{i_t} + x_{i_{t-1}}^q$ is algebraic over $k$ in $R/\ppp_j$. Indeed, otherwise there would be $q_1 \neq q_2$ such that $x_{i_{t-1}}^{q_1} - x_{i_{t-1}}^{q_2}$ is algebraic over $k$ in $R/\ppp_j$; but this is impossible since by construction $x_{i_{t-1}}$ is transcendental over $k$ in $R/\ppp_j$. Therefore 
\begin{defnlist}[resume]
\item there is a positive integer $q_{t-1}$ such that $x_{i_t} + x_{i_{t-1}}^{q_{t-1}}$ is transcendental over $k$ in $R/\ppp_j$ for each $j = s_{t-2} + 1, \ldots, s_{t-1}$. 
\end{defnlist}
Continuing in this way one can choose integers $q_{t-1}, q_{t-2}, \ldots, q_1$ such that for each $t' = 1, \ldots, t - 1$,
\begin{defnlist}[resume]
\item $x_{i_t} + x_{i_{t-1}}^{q_{t-1}} + \cdots + x_{i_{t-t'}}^{q_{t-t'}}$ is transcendental over $k$ in $R/\ppp_j$ for each $j = s_{t-t'-1} + 1, \ldots, s_{t-t'}$. 
\end{defnlist}
Since for each $t' = 1, \ldots, t$ and each $i < i_t'$, $x_i$ is algebraic over $k$ in $R/\ppp_j$ for each $j > s_{t'-1}$, it follows that $x_{i_t} + x_{i_{t-1}}^{q_{t-1}} + \cdots + x_{i_1}^{q_1}$ satisfies the conclusion of the lemma with $e := i_t$. 
\end{proof}

We go back to the proof of \cref{Macaulay:maximal-case}. Let $y_e$ be the polynomial from \cref{Macaulay:transcendental-choice}. Note that $R = k[x_1, \ldots, x_{e-1}, y_e, x_{e+1}, \ldots, x_r]$. Since $\dim(\aaa) = 0$, the image of $y_e$ in $R/\aaa$ is algebraic over $k$. Let $f := \phi(y_e) \in k[y_e]$ be the polynomial with the minimum degree in $y_e$ which is zero in $R/\aaa$. Note that 
\begin{prooflist}
\item \label{Macaulay:f-in-a} $f \in \aaa$, and
\item \label{Macaulay:f-not-in-pj} $f$ is {\em not} in any $\ppp_j$, since the image of $y_e$ is transcendental over $k$ in $R/\ppp_j$. 
\end{prooflist}

\begin{proclaim} \label{Macaulay:induction}
$g_{m-1}$ is not a zero-divisor in $R/\langle \ccc, f \rangle$. 
\end{proclaim}

\begin{proof}
Note that $R/\langle f \rangle \cong k(\bar y_e)[x_1, \ldots, x_{e-1}, x_{e+1}, \ldots, x_r]$, where $\bar y_e$ is the image of $y_e$ in $R/\aaa$. Let $\bar \bbb$ be the ideal of $R/ \langle f \rangle$ generated by $\bbb$. Krull's principal ideal theorem (\cref{thm:principal-ideal}) and observation \ref{Macaulay:f-not-in-pj} above imply that $\dim(\bar \bbb) = \dim(\bbb) - 1 = r - (m - 1) - 1 =  r -1 - (m-1)$. Since $\bar \bbb$ is generated by $m - 1$ elements (namely the images $\bar g_i$ of $g_i$,\ $1 \leq i \leq m-1$) in the polynomial ring $R/\langle f \rangle$ over $k(\bar y_e)$, the induction hypothesis implies that $\bar g_{m-1}$ is not a zero-divisor in $(R/\langle f \rangle)/ \langle \bar g_1, \ldots, \bar g_{m-2} \rangle \cong R/\langle f, \ccc \rangle$, as required. 
\end{proof}

Returning to the proof of \cref{Macaulay:maximal-case}, assume to the contrary that $\aaa$ is associated with $\bbb$. Then there is $g \in R$ such that $\aaa = (\bbb: g)$ (\cref{thm:primary-decomposition}); in particular $gf \in \bbb = \langle  \ccc, g_{m-1} \rangle$. Therefore there is $h \in R$ such that $gf - g_{m-1}h \in \ccc$, i.e.\ $g_{m-1}h \in \langle \ccc, f \rangle$. \Cref{Macaulay:induction} then implies that $h \in \langle \ccc, f \rangle$, i.e.\ $h - af \in \ccc$ for some $a \in R$. It follows that $gf - g_{m-1}af = (g - g_{m-1}a)f \in \ccc$. Since $f$ is not in any prime ideal associated with $\ccc$ (observation \ref{Macaulay:f-not-in-pj}), $f$ is {\em not} a zero divisor in $R/\ccc$ (\cref{thm:primary-decomposition,prop:zero-union-minimally-prime}). It follows that $g - g_{m-1}a \in \ccc$, so that $g \in \langle \ccc, g_{m-1} \rangle = \bbb$. But then $(\bbb : g) = R \not\subseteq \aaa$. This contradiction finishes the proof of \cref{Macaulay:maximal-case} and \cref{thm:Macaulay}. 
\end{proof}

\section{Properties of order at a point on a possibly non-reduced curve} \label{order-app}
In \cref{nrcurve-section} we defined the notion of ``order at a point on a possibly non-reduced curve,'' and stated some of its properties without proof. In this section we prove these results, namely \cref{order-properties,order-curve}, in respectively \cref{order-properties-app,order-curve-app}. The proof uses somewhat more involved commutative algebra than the rest of the book. Recall that for a point $a$ on a possibly non-reduced curve $C$ and $f \in \local{C}{a}$, the order $\ord_a(f)$ of $f$ at $a$ is the dimension of $\local{C}{a}/f\local{C}{a}$ as a vector space over $\kk$.

\begin{prop}[\cref{order-properties}] \label{order-properties-app}
Let $a$ be a point on a possibly non-reduced curve $C$ and $f \in \local{C}{a}$. Let $C' := \supp(C)$. Recall that $\local{C'}{a}$ is a quotient of $\local{C}{a}$. 
\begin{enumerate}
\item If the image of $f$ is a non zero-divisor in $\local{C'}{a}$, then $\ord_a(f) < \infty$. In particular, if $f$ is a non zero-divisor in $\local{C}{a}$, then $\ord_a(f) < \infty$. 
\item \label{order-properties:finite} If $f$ is a non zero-divisor in $\local{C}{a}$, then $\ord_a(f) < \infty$.
\item $\ord_a(f) = 0$ if and only if $f$ is invertible in $\local{C}{a}$.
\item If $f$ is a non zero-divisor in $\local{C}{a}$ and $g \in \local{C}{a}$, then $\ord_a(fg) = \ord_a(f) + \ord_a(g)$.
\end{enumerate}
\end{prop}

\begin{proof}
If the image of $f$ is a non zero-divisor in $\local{C'}{a}$, then a version of Krull's principal ideal theorem \cite[Corollary 11.18]{am} implies that the quotient $\local{C}{a}/f\local{C}{a}$ of $\local{C}{a}$ by the ideal generated by $f$ is a Noetherian local ring of Krull dimension zero. Therefore it is also Artinian \cite[Theorem 8.5]{am} and a finite dimensional vector space over $\kk$ \cite[Exercise 8.3]{am}, which proves the first assertion. The second assertion follows from the first, and the third assertion is straightforward to prove. For the last assertion, \woutlog\ we may assume that $\ord_a(g) < \infty$. Let $(f_1, \ldots, f_l)$ be a basis of $\local{C}{a}/f\local{C}{a}$ and $(g_1, \ldots, g_m)$ be a basis of $\local{C}{a}/g\local{C}{a}$ over $\kk$. Let $h \in \local{C}{a}$. Write $K_f$ and $K_g$ for respectively the $\kk$-linear span of the $f_i$ and of the $g_j$. Then $\local{C}{a} = K_f + f\local{C}{a} = K_f + fK_g + fg\local{C}{a}$, which implies that $f_1, \ldots, f_l, fg_1, \ldots, fg_m$ spans $\local{C}{a}$ modulo the ideal generated by $fg$. We claim that they are also linearly independent over $\kk$. Indeed, pick $c_i, d_j \in \kk$ such that $\sum_ic_if_i + \sum_j d_jfg_j = fgh$, $h\in \local{C}{a}$. Then $\sum_ic_if_i  \in f\local{C}{a}$, so that $c_1 = \cdots = c_l = 0$. It follows that $f\sum_jd_jg_j = fgh$. Since $f$ is a non zero-divisor, $\sum_j d_j g_j = gh \in g\local{C}{a}$. It follows that $d_1 = \cdots = d_m = 0$ as well. 
\end{proof}

The proof of \cref{order-curve} will be long. We start with a few auxiliary results. 

\begin{lemma} \label{lemma:k[[t]]-module}
Let $R$ be a ring containing a field $k$, and $t \in R$ be such that $t$ is transcendental over $k$, and $R$ contains the ring $k[[t]]$ of power series in $t$ over $k$. Assume $R/tR$ is a finite dimensional vector space over $k$ generated by (the images in $R/tR$ of) $f_1, \ldots, f_m \in R$. Then $f_1, \ldots, f_m$ generate $R$ as a module over $k[[t]]$. 
\end{lemma}

\begin{proof}
Given $g \in R$, it can be successively expressed as $g	= \sum_j f_jh_{0,j} + tg_1 = \sum_j f_j(h_{0,j} + th_{1,j}) + t^2g_2$, and so on, resulting in an expression of the form $g = \sum_j f_jh_j$ with $h_j \in \kk[[t]]$. 
\end{proof}

\begin{prop} \label{irreducible-order}
Let $a$ be a point on an possibly non-reduced affine curve $C$. Assume $C' := \supp(C)$ is irreducible and nonsingular at $a$. Let $t$ be a regular function on $C$ such that $t|_{C'}$ is a parameter of $\local{C'}{a}$. Assume $t$ is not a zero-divisor in $\local{C}{a}$. Then
\begin{enumerate}
\item $\ord_a(t) = \mu_{C'}(C)$, where $\mu_{C'}(C)$ is the multiplicity of $C'$ in $C$ (defined in \cref{nrcurve-section} in the paragraph preceding \cref{order-curve}), 
\item  $\ord_a(f) = \ord_a(f|_{C'})\ord_a(t)$ for each regular function $f$ on $C$. 
\end{enumerate}
\end{prop}

\begin{proof}
We start with the first assertion. \Woutlog\ we may assume that $C$ is the closed subscheme of $\kk^{n+1}$ determined by an ideal $J$ in $S := \kk[x_0, \ldots, x_n]$, and $a$ corresponds to the point $x_0 = \cdots = x_n = 0$, and $t$ is the restriction of $x_0$. Then $\hatlocal{C}{a} \cong \hat S/J \hat S$ where $\hat S := \kk[[x_0, \ldots, x_n]]$ (\cref{thm:exactly-complete,smooth-completion}). Let $r := \dim_\kk(\hat S/(x_0\hat S + J\hat S)) =  \dim_\kk(\hatlocal{C}{a}/t\hatlocal{C}{a})$. Assertion \eqref{order-properties:finite} of \cref{order-properties-app} and \cref{exercise:finite-quotient} imply that $r = \dim_\kk(\local{C}{a}/t\local{C}{a}) = \ord_a(t)$. In particular, $r < \infty$, so that \cref{lemma:k[[t]]-module} implies that $\hat S/J \hat S$ is a finitely generated $\kk[[x_0]]$-module. Since $t$ is a non zero-divisor in $\local{C}{a}$, it follows that $x_0$ is a non zero-divisor in $\hatlocal{C}{a} = \hat S/J \hat S$ (this is due to the general version of exactness of completions, see e.g.\ \cite[Exercise 10.4.]{am}). The fundamental theorem of finitely generated modules over a PID (\cref{fundamentally-over-PID}) then implies that $\hat S/J\hat S \cong \kk[[x_0]]^{s}$ as a module over $\kk[[x_0]]$ for some $s \geq 0$. But then $s = \dim_\kk(\hat S/(x_0\hat S + J \hat S)) = r$. Since $t|_{C'}$ is a parameter of $\local{C'}{a}$, for each $i =1, \ldots, n$, there is $u_i \in \local{C}{a}$ and a positive integer $m_i$ such that $x_i - u_ix_0^{m_i}$ is nilpotent in $\local{C}{a}$. Pick a representative $\psi_i \in \hat S$ of $u_i$ and set $y_i := x_i - x_0^{m_i}\psi_i$. Pick $\phi_1, \ldots, \phi_r \in \hat S$ such that $\hat S/J\hat S = \dsum_{i=1}^r \phi_i \kk[[x_0]]$ as a $\kk[[x_0]]$-module. Since the linear parts of $(x_0, y_1, \ldots, y_n)$ are linearly independent, \cref{thm:linear-power} implies that $\hat S = \kk[[x_0, y_1, \ldots, y_n]]$. Moreover, since each $y_j$ is nilpotent modulo $J\hat S$, each $\phi_i$ can be expressed as a $\kk[[x_0]]$-linear combination of finitely many monomials in $(y_1, \ldots, y_n)$. Let $y^{\alpha_j}$, $j = 1, \ldots, N$, be a minimal collection of such monomials such that every other monomial in $(y_1, \ldots, y_n)$ is their $\kk((x_0))$-linear combination modulo $J\hat S$. After a Gauss-Jordan elimination process and re-orderings of the $\alpha_j$ if necessary we may assume that $\phi_i= y^{\alpha_i} + \sum_{j = r+1}^N \phi_{i,j}(x_0)y^{\alpha_j}$, $i = 1, \ldots, r$, where $\phi_{i,j}(x_0) \in \kk((x_0))$. Now express $y^{\alpha_1}$ as a $\kk[[x_0]]$-linear combination of the $\phi_i$ modulo $J \hat S$:
\begin{align*}
y^{\alpha_1} 
	= \sum_{i=1}^r \rho_i(x) \phi_i 
	= \sum_{i=1}^r \rho_i(x)y^{\alpha_i} + \sum_{i=1}^r\sum_{j=r+1}^N \rho_i(x) \phi_{i,j}(x)y^{\alpha_j} 
\end{align*}
The minimality assumption on the $y^{\alpha_i}$ implies that $\rho_1(x) = 1$, and $\rho_i(x)\phi_{i,j}(x) = 0$ if either $i > 1$ or $j > 1$. It follows that $y^{\alpha_1} = \phi_1$. The same arguments inductively show that $\phi_i = y^{\alpha_i}$ for each $i = 1, \ldots, r$. In particular, $N = r$. 

\begin{proclaim} \label{claim:1/t}
$\local{C}{C'}$ (which was defined in \cref{scheme:local-section}) can be identified with a $\kk$-subalgebra of $\local{C}{a}[1/t]$.
\end{proclaim} 

\begin{proof}
Elements of $\local{C}{C'}$ are of the form $f/g$ where $f,g \in \local{C}{a}$ such that $g|_{C'} \not\equiv 0$. Let $m := \ord_a(g|_{C'})$. \Cref{reduced-invertible} implies that there is invertible $u \in \local{C}{a}$ such that $g = ut^m - h$ where $h$ is a nilpotent in $\local{C}{a}$. Pick $k$ such that $h^{k+1} = 0$. Write $f' := f/u$ and $h' := h/u$. Then $f/g = f'/(t^m - h') = t^{-m(k+1)}f'\sum_{j=0}^k t^{mj}h'^{k-j}$ which is an element of $\local{C}{a}[1/t]$. Since $t$ is a non zero-divisor in $\local{C}{a}$, it follows that the map $\local{C}{C'} \to \local{C}{a}[1/t]$ is injective. 
\end{proof}

\Cref{claim:1/t,prop:complete-inclusion} imply that there are $\kk$-subalgebra homomorphisms 
\begin{align}
\local{C}{C'} \into \hatlocal{C}{a}[1/t] \cong \hat S_{x_0}/ J\hat S_{x_0}
\label{eqn:1/t}
\end{align}
where the first map is injective and $\hat S_{x_0} := \hat S[1/x_0] = \kk((x_0))[[y_1, \ldots, y_n]]$. Given any finite collection of elements $\beta_1, \ldots, \beta_s \in \znzero$, let $\hat J$ be the ideal of $\hat S_{x_0}/J \hat S_{x_0}$ generated by $y^{\beta_1}, \ldots, y^{\beta_s}$ and $M := \hat J \cap \local{C}{C'}$. 

\begin{proclaim} \label{claim:M_i-generators}
$y^{\beta_1}, \ldots, y^{\beta_s}$ generate $M$ as an ideal of $\local{C}{C'}$.
\end{proclaim} 

\begin{proof}
Pick $h \in M$. Due to \eqref{eqn:1/t} there is $k \geq 0$ such that $ht^k$ is represented by an element in $\local{C}{a}$ and $hx_0^k$ is represented by a $\kk[[x_0]]$-linear combination of the $y^{\beta_1}, \ldots, y^{\beta_s}$ in $\hat S/J \hat S$. Write $M'$ for the ideal of $\local{C}{a}$ generated by $y^{\beta_1}, \ldots, y^{\beta_s}$ and $N'$ for the ideal of $\local{C}{a}$ generated by $M'$ and $ht^k$. Let $L':= N'/M'$ be the quotient of $N'/M'$ as a module over $\local{C}{a}$. Let $\hat L', \hat M', \hat N'$ be the completion (with respect to the maximal ideal of $\local{C}{a}$) of respectively $L', M', N'$. \Cref{thm:exactly-complete} implies that $\hat L' \cong \hat N'/\hat M'$, and then \cite[Proposition 10.13]{am} implies that $\hat L' \cong (\hatlocal{C}{a} \otimes_{\local{C}{a}} N')/(\hatlocal{C}{a} \otimes_{\local{C}{a}} M') = 0$. A theorem of Krull \cite[Theorem 10.17]{am} then implies that $L' = 0$, so that $ht^k \in M'$. Since $t$ is invertible in $\local{C}{C'}$, it follows that $h \in M'$, as required. 
\end{proof}

Recall that the only condition satisfied by $\alpha_1, \ldots, \alpha_r$ is that 
\begin{prooflist}
\item $\{y^{\alpha_i}: i = 1, \ldots, r\}$ is a minimal collection of monomials in $(y_1, \ldots, y_n)$ such that every other monomial is a $\kk((x_0))$-linear combination modulo $J\hat S$. 
\end{prooflist}
While choosing such $\{\alpha_i\}$, we can start from monomials with {\em maximum} possible degree and then sequentially adjoin monomials with smaller degrees to ensure that for each $i \geq 2$, 
\begin{prooflist}[resume]
\item $y^{\alpha_i + e_j}$ is a $\kk((x_0))$-linear combination modulo $J\hat S$ of $y^{\alpha_1}, \ldots, y^{\alpha_{i-1}}$ for each $j = 1, \ldots, n$ (where $e_j$ is the $j$-th standard unit vector in $\zz^n$). 
\end{prooflist}
For each $i = 0, \ldots, r$, let $\hat J_i$ be the ideal of $\hat S_{x_0}/J \hat S_{x_0}$ generated by $y^{\alpha_1}, \ldots, y^{\alpha_{r-i}}$ and $M_i := \hat J_i \cap \local{C}{C'}$. Fix $i$, $0 \leq i \leq r-1$. \Cref{claim:M_i-generators} implies that $M_i/M_{i+1}$ is generated by $y^{\alpha_{r-i}}$ as an $\local{C}{C'}$-module. Since $y^{\alpha_{r - i}} \in M_i \setminus M_{i+1}$, it follows that $M_i/M_{i+1} \neq 0$. On the other hand, if $h$ is a nilpotent element in $\local{C}{C'}$, then as an element of $\hat S_{x_0} / \hat J\hat S_{x_0}$, $h$ is in the ideal generated by $y_1, \ldots, y_n$, and the choice of the $\alpha_j$ ensures that $hM_i \subseteq M_{i+1}$. Combining these observations we see that as an $\local{C}{C'}$-module, $M_i/M_{i+1} \cong \local{C}{C'}/\nnn$, where $\nnn$ is the (maximal) ideal of nilpotent elements of $\local{C}{C'}$. It follows that $\local{C}{C'} = M_0 \supset M_1 \supset \cdots \supset M_r = 0$ is a composition series of $\local{C}{C'}$, so that $\mu_{C'}(C) = r$, which proves the first assertion of \cref{irreducible-order}. \\

For the second assertion, we may assume \woutlog\ that $q := \ord_a(f|_{C'}) < \infty$. Then there is $u \in \local{C}{a}$ such that $u|_{C'}$ is invertible and $(t^q - uf)|_{C'} = 0$. \Cref{reduced-invertible} implies that $u$ is invertible in $\local{C}{a}$. Note that $h := t^q - uf$ is nilpotent in $\local{C}{a}$. By replacing $C$ by a smaller neighborhood\footnote{Since $C$ is a {\em possibly non-reduced} curve, by a ``neighborhood of $a$ in $C$'' we mean an open subscheme of $C$ containing $a$.} of $a$ in $C$ if necessary we may assume $h$ is nilpotent in $S/J$. Let $z$ be a new indeterminate which we think of as the last coordinate of $\kk^{n+1} \times \kk = \kk^{n+2}$, and let $D$ be the closed subscheme of $\kk^{n+2}$ determined by the ideal $K$ in $T := \kk[x_0, \ldots, x_n, z]$ generated by $J$ and $z^q - h$. In $T/K$ we have $uf = t^q - z^q = \prod_{i=1}^q (t - \zeta_i z)$, where the $\zeta_i$ are the $q$-th roots of unity in $\kk$. Note that $D' := \supp(D)$ is isomorphic to $C'$; in particular $a$ can be naturally identified with a point on $D$, which by an abuse of notation we also denote by $a$. Since $z$ is nilpotent in $T/K$, it follows that $(t- \zeta_i z)|_{D'} = t|_{D'}$ is a parameter of $\local{D'}{a}$ for each $i$. Since $T/K = \sum_{i=0}^{q-1} z^i S/J \cong (S/J)^q$ as a module over $S/J$, it follows that 
\begin{prooflist}
\item $t$ is not a zero-divisor in $\local{D}{a}$, which implies that $t - \zeta_i z$ is not a zero-divisor in $\local{D}{a}$ for any $i$ (since it is easy to check that if $g_1$ is a zero-divisor and $g_2$ is nilpotent, then $g_1 + g_2$ is a zero-divisor), and
\item $\ord_a(g|_D) = q\ord_a(g|_C)$ for each for each $g \in S/J$. 
\end{prooflist}
These observations together with the first assertion and \cref{order-properties-app} imply that $q\ord_a(f|_C) = \ord_a((uf)|_D) = \sum_{i=0}^{q-1}\ord_a((t - \zeta_i z)|_D) = q\ord_a(t|_D) = q^2 \ord_a(t|_C)$, so that $\ord_a(f|_C) = q\ord_a(t|_C)$, as required. 
\end{proof}

\begin{lemma} \label{finite-order}
Let $D$ be a possibly non-reduced affine curve over $\kk$ defined by an ideal $I$ of $\kk[x_1, \ldots, x_n]$. If $f \in R :=\kk[x_1, \ldots, x_n]/I$ is such that $r := \dim_\kk(R/fR) < \infty$. Then $r = \sum_{a \in D} \ord_a(f)$. 
\end{lemma}

\begin{proof}
Since $r < \infty$, there are finitely many zeroes of $f$ on $D$. Denote them by $a_1, \ldots, a_k$; let $\mmm_j$ be the maximal ideal of $a_j$ in $R$ and $\iota_j: R \to R_{\mmm_j}$ be the natural map.

\begin{proclaim} \label{chinese-claim-1}
$fR = \bigcap_j \iota_j^*(fR_{\mmm_j})$.
\end{proclaim}

\begin{proof}[Proof (\`a la Mumford {\cite[Proposition 1.11]{mummetry}})]
We only need to show the ``$\supset$'' inclusion. Let $h \in \iota_j^*(fR_{\mmm_j})$ for each $j$. Then for each maximal ideal $\mmm$ of $R$, there exists $u \not\in \mmm$ such that $uh \in fR$. It follows that the ideal $(fR : h) := \{u\in R: uh \in fR\}$ of $R$ is not contained in any maximal ideal of $R$. The Nullstellensatz then implies that $1 \in (fR: h)$, as required. 
\end{proof}

\begin{proclaim} \label{chinese-claim-2}
If $j \neq j'$, then $\iota_j^*(fR_{\mmm_j}) + \iota_{j'}^*(fR_{\mmm_{j'}}) = R$. 
\end{proclaim}

\begin{proof}
Since $R_{\mmm_j}/fR_{\mmm_j}$ and  $R_{\mmm_{j'}}/fR_{\mmm_{j'}}$ are Artinian local rings, their maximal ideals are nilpotent \cite[Proposition 8.6]{am}. Therefore there exists $q$ such that $(\mmm_j)^qR_{\mmm_j}\subseteq fR_{\mmm_j}$ and $(\mmm_{j'})^qR_{\mmm_{j'}}\subseteq fR_{\mmm_{j'}}$. Since $1 \in (\mmm_j)^q + (\mmm_{j'})^q$ [why?], the claim follows. 
\end{proof}

\Cref{chinese-claim-1,chinese-claim-2} and the Chinese remainder theorem \cite[Proposition 1.10]{am} imply that 
\begin{align}
r = \dim_\kk(R/fR) = \sum_{j=1}^k \dim_\kk(R/\iota_j^*(fR_{\mmm_j})) \label{r-1}
\end{align}  
\begin{proclaim} \label{m_j-claim}
$R/\iota_j^*(fR_{\mmm_j}) \cong R_{\mmm_j}/fR_{\mmm_j} $ for each $j$. 
\end{proclaim}

\begin{proof}
It is straightforward to see that the natural map $R/\iota_j^*(fR_{\mmm_j}) \to R_{\mmm_j}/fR_{\mmm_j}$ is injective. For surjectivity, note that if $h = g_1/g_2$ is an element of $(\mmm_j)^qR_{\mmm_j}$ for some $q \geq 0$, where $g_1\in R$ and $g_2 \in R\setminus \mmm_j$, then $h - cg_1 \in (\mmm_j)^{q+1}R_{\mmm_j}$, where $c \in \kk$ is the image of $g_2^{-1}$ in $\kk \cong R_{\mmm_j}/\mmm_jR_{\mmm_j}$. It follows by an induction on $q$ that for each $h \in R_{\mmm_j}$ and $q \geq 1$, there exists $h' \in R$ such that $h - h' \in (\mmm_j)^{q+1}R_{\mmm_j}$. Choosing $q$ such that $(\mmm_j)^{q+1}R_{\mmm_j} \subseteq fR_{\mmm_j}$ (which is possible due to the arguments in the proof of \cref{chinese-claim-2}) yields the required result.
\end{proof}

The result follows from \eqref{r-1} and \cref{m_j-claim}.
\end{proof}

\begin{thm}[\cref{order-curve}] \label{order-curve-app}
Let $a$ be a point on a possibly non-reduced curve $C$. Let $C_1, \ldots, C_s$ be the irreducible components of $\supp(C)$ containing $a$ and $\pi_i:\tilde C_i \to C_i$ be the desingularizations of $C_i$. If $f$ is a non zero-divisor in $\local{C}{a}$, then
\begin{align}
\ord_a(f) 
	&= \sum_i \mu_{C_i}(C) \ord_{ a}(f|_{C_i})
	= \sum_i \mu_{C_i}(C) \sum_{\tilde a \in \pi_i^{-1}(a)} \ord_{\tilde a}(\pi_i^*(f|_{C_i})) 
	\label{order-curve-identity-app}
\end{align}
\end{thm}

\begin{proof}
Write $C' := \supp(C)$. \Woutlog\ we may assume that 
\begin{prooflist}
\item $C$ is affine, i.e.\ the closed subscheme of $\kk^n$ determined by an ideal $\aaa$ of $\kk[x_1, \ldots, x_n]$; in particular $\kk[C'] = \kk[x_1, \ldots, x_n]/\sqrt{\aaa}$;
\item $C' \setminus\{a\}$ is nonsingular;
\item $f$ is in the maximal ideal of $\local{C}{a}$; 
\item $f$ restricts to a regular function on $C'$, and 
\item the restriction of $f$ to every irreducible component of $C'$ is non-constant.
\end{prooflist}
Let $\bar C'$ be the (unique) compactification of $C'$ such that $\bar C' \setminus \{a\}$ is nonsingular (given a closed embedding of $C'$ into an affine space, $\bar C'$ can be explicitly constructed by taking the closure of $C'$ in an projective completion of the affine space and then resolving the singularities at infinity of the closure of $C'$). The restriction $f' := f|_{C'}$ of $f$ to $C'$ induces a morphism $\bar C' \to \pp^1$ (\cref{prop:curve-nonsingular-morphextension}); we use $f'$ to denote this morphism as well. Note that $f'(a) = 0 \in \kk \subset \pp^1$ (since $f$ is in the maximal ideal of $\local{C}{a}$). It follows that
\begin{prooflist}[resume]
\item There is a finite set $S$ of $\kk \setminus \{0\}$ such that $f'^{-1}(\kk \setminus S)$ is an {\em affine} curve\footnote{In fact $f'^{-1}(\kk)$ is also affine; we did not use this fact since all the proofs we know of it use the {\em Riemann-Roch} theorem, which we did not cover in this book.} (\cref{exercise:affine-inverse}). 
\end{prooflist}
Let $D' := f'^{-1}(\kk \setminus S)$ and $p' := \prod_{s \in S}(f'-s) \in \kk[C'] \cap \kk[D']$. Then $D' \supseteq C' \setminus V(p')$. Let $q' \in \kk[D']$ such that $q'(a) \neq 0$ and $D' \setminus V(q') \subseteq C'$. There is $N \geq 0$ such that $p'^Nq' \in \kk[C']$ (since $\kk[D'] \subseteq \kk[C']_{p'}$). Let $g' := p'^{N+1}q' \in \kk[f']$. It is straightforward to check that $g'$ is regular on $C' \cup D'$, $g'(a) \neq 0$, $D' \setminus V(g') = C' \setminus V(g')$, and $\kk[D']_{g'} = \kk[C']_{g'}$. In other words, if $g \in \kk[x_1, \ldots, x_n]$ is such that $g' = g|_{C'}$, and we write $X$ for the open subscheme $C \setminus V(g)$ of $C$, and $X'$ for the support of $X$, then there is a sequence of morphisms as in \cref{fig:XX'DD'-0}.  

\newcommand{\subfigwidth}{0.4\textwidth}
\begin{center}
\begin{figure}[h]


\begin{subfigure}[b]{\subfigwidth}
\begin{tikzcd}
X 
\arrow[d, "\supp"]
&
& \\
X'
\arrow[r, "\cong"]
&
D'\setminus V(g')
\arrow[hook, r]
&D'
\end{tikzcd}
\caption{}
\label{fig:XX'DD'-0}
\end{subfigure}
\begin{subfigure}[b]{\subfigwidth}
\begin{tikzcd}
X 
\arrow[d, "\supp"]
\arrow[r, "\cong"]
&
D \setminus V(g)
\arrow[r, hook]
\arrow[d, "\supp"]
& 
D
\arrow[d, "\supp"]
\\
X'
\arrow[r, "\cong"]
&
D'\setminus V(g')
\arrow[r, hook]
&D'
\end{tikzcd}
\caption{}
\label{fig:XX'DD'-1}
\end{subfigure}
\caption{Compactification of $C\setminus V(g)$}  \label{fig:XX'DD'}
\end{figure}
\end{center}

Up to an embedded isomorphism, we can treat $X$ as an ``affine scheme,'' i.e.\ a closed subscheme of an affine space (\cref{example:open-embedded-affine}). We now construct an affine scheme $D$ with support $D'$ such that the diagram in \cref{fig:XX'DD'-1} commutes, where the ``$\cong$'' on the top left denotes an embedded isomorphism. Indeed, choose $\kk$-algebra generators $f_1 = 1/g, f_2 = g, f_3, \ldots, f_k$ of $A:= (\kk[x_1, \ldots, x_n]/\aaa)_g$ such that 
\begin{prooflist}[resume]
\item for each $j = 2, \ldots, k$, the image $\bar f_j$ in $\kk[X']$ of $f_j$ is in the image of the natural map $\kk[D'] \to \kk[X']$, and 
\item $\kk[D'] = \kk[\bar f_2, \ldots, \bar f_k]$. 
\end{prooflist}
Let $B := \kk[f_2, \ldots, f_k] \subset A$, and $D$ be the ``closed subscheme of $\kk^{k-1}$ with coordinate ring $B$,'' i.e.\ $D$ is the closed subscheme of $\kk^{k-1}$ corresponding to the kernel of the surjective map $\kk[x_2, \ldots, x_k] \to B$ which maps each $x_j$ to $f_j$. Then it is straightforward to check that the diagram in \cref{fig:XX'DD'-1} commutes. For each $i = 1, \ldots, s$, we write $D_i$ for the irreducible component of $D$ containing $C_i$. Let $p :=  \prod_{s \in S} (f-s) \in \kk[f]$, so that $p' = p|_{C'}$. Then $p'$ is invertible in $\kk[D']$, which implies that $p$ is invertible in $B$.  

\begin{proclaim}
$B$ (respectively $\kk[D']$) is a finitely generated module over $\kk[f]_p$ (respectively $\kk[f']_{p'}$). 
\end{proclaim}

\begin{proof}
\Cref{curve-finite} implies that each $\kk[D_i]$ is a finitely generated $\kk[f']_{p'}$-module. Since the natural map $\kk[D'] \to \prod_i \kk[D_i]$ is injective, $\kk[D']$ is isomorphic to a $\kk[f']_{p'}$-submodule of $\prod_i \kk[D_i]$, and therefore also a finitely generated module over $\kk[f']_{p'}$. Pick a finite collection $g'_1, \ldots, g'_{k'}$ of $\kk[f']_{p'}$-module generators of $\kk[D']$. Pick $g_i \in B$ such that $g'_i = g_i|_{D'}$, $i = 1, \ldots, k'$. Let $h_1, \ldots, h_l$ be generators of the ideal $\nnn$ of nilpotent elements of $B$. There is $m \geq 0$ such that $\nnn^{m+1} = 0$. We claim that $B$ is generated as a $\kk[f]_p$-module by $g_ih_1^{\alpha_1} \cdots h_m^{\alpha_m}$, $1 \leq i \leq k'$, $\alpha \in \zzeroo{l}$, $\alpha_1 + \cdots + \alpha_l \leq m$. Indeed, let $u \in B$. Then there are $\phi_1, \ldots, \phi_{k'} \in \kk[f]_p$ such that $u_1 := u - \sum_{j=1}^{k'} \phi_j g_j \in \nnn$. Then $u_1 = \sum_{j=1}^l u_{1,j} h_j$. Expressing the $u_{1,j}$ as $\kk[f]_p$-linear combinations of the $g_j$ modulo $\nnn$ and continuing as above gives the claim. 
\end{proof}

The fundamental theorem of finitely generated modules over a PID (\cref{fundamentally-over-PID}) implies that as a $\kk[f]_p$-module $B$ has a decompositions of the form:
\begin{align*}
B &\cong (\kk[f]_p)^r \dsum \left(\dsum_j \kk[f]_p/\langle \phi_j(f) \rangle \right) 
\end{align*}
such that each $\phi_j$ is a polynomial in an indeterminate $t$. For each $c \in \kk$, it follows that $f-c$ is a zero-divisor in $B$ if and only if $c \not \in S$ (since $f-c$ is invertible in $\kk[f]_p$) and $t-c$ divides some $\phi_j(t)$ in $\kk[t]$; in particular, there are only finitely many such $c \in \kk$. It is then straightforward to see that $r = \dim_\kk(B/\langle f-c \rangle)$ for all $c \in \kk$ such that $f-c$ is a non zero-divisor in $B$. \Cref{finite-order} then implies that 
\begin{align*}
r 	= \sum_{\substack{b \in D\\ f(b) = c}} \ord_b(f-c) 
\end{align*}
for all $c \in \kk$ such that $f-c$ is a non zero-divisor in $\kk[D]$. Now pick an arbitrary point $a^*  \in C\setminus \{a\}$. Applying the above construction with $a^*$ and $C^* := C \setminus\{a\}$ respectively in place of $a$ and $C$ yields an affine possibly non-reduced curve $D^*$ such that $\supp(D^*)$ containing an open neighborhood of $a^*$, and the same arguments show that there is an integer $r^*$ such that 
\begin{align*}
r^*	= \sum_{\substack{b \in D^*\\ f(b) = c}} \ord_b((f-c)|_{D^*}) 
\end{align*}
for all $c \in \kk$ such that $f-c$ is a non zero-divisor in the ``coordinate ring'' of $D^*$. Since $D$ and $D^*$ are ``birational''\footnote{i.e.\ there are open subschmes $U, U^*$ respectively of $D, D^*$ such that $\supp(U)$ (respectively, $\supp(U^*)$) intersects each irreducible component of $D$ (respectively, $D^*$), and there is an embedded isomorphism between $U$ and $U'$.} and since our construction guarantees that $f$ is a non zero-divisor in coordinate rings of both $D$ and $D^*$, it follows that $r = r^*$ and 
\begin{align}
\sum_{\substack{b \in D\\ f(b) = 0}} \ord_b(f)
	= \sum_{\substack{b \in D^*\\ f(b) = 0}} \ord_b(f|_{D^*}) 
\label{eqn:order-D=D*}
\end{align}
Let $a_1 = a, a_2, \ldots, a_k$ be the points of $f^{-1}(0)$ on $D$. For each $j > 1$, the construction of $D^*$ shows that $a_j \in D^*$ and $a_j$ has a neighborhood on $D$ which is isomorphic to $D^*$. Let $S^* := (f|_{\supp(D^*)})^{-1}(0) \setminus \{a_2, \ldots, ,a_k\} \subset D^*$. Identity \eqref{eqn:order-D=D*} implies that 
\begin{align*}
\ord_a(f)
	= \sum_{b \in S^*} \ord_b(f|_{D^*}) 
\end{align*}
It follows from the construction of $D^*$ that $\supp(D^*)$ is nonsingular at every point on $S$. For each $b \in S^*$, let $D^*_{i_b}$ be the (unique) irreducible component of $D^*$ containing $b$. \Cref{irreducible-order} then implies that 
\begin{align*}
\ord_a(f)
	= \sum_{b \in S^*} \mu_{D^*_{i_b}}(D^*)\ord_b(f|_{D^*_{i_b}}) 
\end{align*}
Since there $D^*_{i_b}$ and $C_i$ have isomorphic nonempty Zariski open subsets such that the corresponding open subschemes of $D^*$ and $C$ are isomorphic, it follows that $\mu_{D^*_{i_b}}(D^*) = \mu_{C_{i_b}}(C)$, and since the desingularization $\tilde C_i$ of $C_i$ is isomorphic to $D^*_{i_b}$ near $b$, it follows that $\ord_b(f|_{D^*_{i_b}})  = \ord_b(\pi_{i_b}^*(f|_{C_{i_b}}))$, which proves the theorem. 
\end{proof}

\chapter{Notation}
\begingroup
\allowdisplaybreaks
\begin{alignat*}{2}
&:=
	&& \text{is defined as} \\
&\cong
	&& \text{isomorphic} \\
&\norm{\cdot}
	&& \text{Euclidean length} \\
&\coprod
	&& \text{disjoint union} \\
&\langle \nu, \cdot \rangle,\ \nu \in \rnstar
	&& \text{the function induced by $\nu$ on $\rr^n$}\\
&\langle \beta, \cdot \rangle,\ \beta \in \rr^n
	&& \text{dot product with $\beta$}\\
&\eta^\perp
	&& \{\alpha \in \rr^n: \langle \eta, \alpha \rangle = 0\}\\
&\langle f, g, \ldots, \rangle
	&& \text{ideal generated by}\ f,g, \ldots \\
&\kk[[x]]
	&& \text{ring of formal power series in $x$ with coefficients in $\kk$}\\
&\kk((x))
	&& \text{field of Laurent series in $x$ with coefficients in $\kk$}\\
&\kk[X]
	&& \text{ring of rational functions on the algebraic variety $X$ defined over $\kk$}\\
&{[n]},\ n \in \zz
	&& \{1, 2, \ldots, n\} \\
&\mat{T},\ T: \nktorus \to \nktoruss{m}
	&& \text{$n \times m$ matrix of exponents of coordinates of $T$}\\
&x^\alpha,\ \alpha = (\alpha_1, \ldots, \alpha_n)
	&& x_1^{\alpha_1} \cdots x_n^{\alpha_n}\\
&\id_n
	&& \text{identity matrix of size $n \times n$}\\
&\aff(\cdot)
	&& \text{affine hull}\\
&\cone(S)
	&& \text{cone generated by $S$} \\
&\conv(\cdot)
	&& \text{convex hull}\\
&\partialxi(\cdot)
 	&& \partialxioriginal{(\cdot)} \\
&\fund(H)
	&& \text{volume of any fundamental lattice parallelotope of $H$}\\
&\gcd(a, b, \cdots)
	&& \text{greatest common positive divisor of the nonzero elements from $a, b, \cdots$}\\
&\In_\nu(\cdot)
	&& \text{initial form of a (Laurent) polynomial or minimizing face of a polyhedron} \\
&\jjj(\cdot) 
	&& \text{ideal generated by the partial derivatives of a polynomial} \\
&\lcm(a, b, \cdots)
	&& \text{lowest common positive multiple of the nonzero elements from $a, b, \cdots$}\\
&\ld
	&& \text{leading form of a (Laurent) polynomial or maximizing face of a polyhedron} \\
&\max_\scrP(\nu)
	&& \max\{\langle \nu, \alpha \rangle: \alpha \in \scrP\} \\
&\min_\scrP(\nu)
	&& \min\{\langle \nu, \alpha \rangle: \alpha \in \scrP\} \\
&\mv
	&& \text{mixed volume} \\
&\mv'_\nu
	&& \text{normalized mixed volume} \\
&\nnn(\cdot)
	&& \text{nilradical}\\
&\nd(\cdot)
	&& \text{Newton diagram}\\
&\np(\cdot)
	&& \text{Newton polytope}\\
&\local{X}{Z}
	&& \text{local ring of the variety $X$ at its subvariety $Z$}\\	
&\qplus
	&& \{q \in \qq: q > 0\} \\
&\qzero
	&& \{q \in \qq: q \geq  0\} \\
&\rplus
	&& \{r \in \rr: r > 0\} \\
&\rzero
	&& \{r \in \rr: r \geq  0\} \\
&\relint(\cdot)
	&& \text{relative interior} \\
&\supp(f)
	&& \text{support of $f$}\\
&\trd_k(\cdot)
	&& \text{transcendence degree over (a field) $k$}\\
&V(f, g, \ldots)
	&& \text{the set of zeroes of (or depending on the context, the closed subscheme defined by) $f, g, \ldots$} \\
&\vol'_{H}
	&& \text{normalized lattice vlume on an affine subspace $H$} \\
&\vol'_\nu
	&& \text{$\vol'_H$, where $H := \{\alpha \in \rr^n: \langle \nu, \alpha \rangle = 0\}$} \\
&\zplus
	&&  \{0,1, 2, \ldots, \} \\
&\zzero
	&&  \{1, 2, \ldots, \}  \\
\end{alignat*}
\endgroup 


\bibliographystyle{alpha}
\bibliography{../../../../../utilities/bibi}

\printindex
\end{document}